\documentclass[a4paper,12pt,oneside,leqno]{book}
\usepackage{amssymb,amsthm,amsmath,amscd,longtable,verbatim}
\input{cyrfonts.def}
\author{D.~A.~Timashev}
\title{Homogeneous spaces and \\ equivariant embeddings}
\date{January~10, 2006}

\newcommand{\adheres}{\prec}
\newcommand{\adhereseq}{\preceq}
\newcommand{\ann}{\perp}
\newcommand{\biject}{\leftrightarrow}
\newcommand{\by}{/\!\!/}
\newcommand{\dbloverline}[1]{%
  \overline{\raisebox{0.3ex}{\vphantom{\ensuremath{#1}}}%
  \smash[t]{\overline{#1}}}}
\newcommand{\divby}{\smash{\,\vdots\,}}
\newcommand{\embeds}{\hookrightarrow}
\newcommand{\embof}{\hookleftarrow}
\newcommand{\iso}{\simeq}
\newcommand{\isoto}{\stackrel{\sim}{\to}}
\newcommand{\itimes}[1]{*_{#1}}
\newcommand{\normin}{\vartriangleleft}
\newcommand{\onto}{\twoheadrightarrow}
\newcommand{\semitimes}{\rightthreetimes}
\newcommand{\timessemi}{\leftthreetimes}
\newcommand{\sort}{\angle}
\newcommand{\tran}{\top}
\newcommand{\stran}{\dag}

\renewcommand{\[}{\mathopen{[\mspace{-1mu}[}}
\renewcommand{\]}{\mathclose{]\mspace{-1mu}]}}
\renewcommand{\(}{\mathopen{(\mspace{-3mu}(}}
\renewcommand{\)}{\mathclose{)\mspace{-3mu})}}

\newcommand{\Aa}{\mathcal A}
\newcommand{\AAA}{\mathbb A}
\newcommand{\AAa}{\mathbf A}
\newcommand{\ab}{\mathfrak a}
\newcommand{\B}{\mathfrak B}
\newcommand{\br}{\mathfrak b}
\newcommand{\Bb}{\mathcal B}
\newcommand{\BBb}{\mathbf B}
\newcommand{\CC}{\mathbb C}
\newcommand{\Cc}{\mathcal C}
\newcommand{\CCc}{\mathbf C}
\newcommand{\Dd}{\mathcal D}
\newcommand{\DDd}{\mathbf D}
\newcommand{\e}{\mathfrak e}
\newcommand{\Ee}{\mathcal E}
\newcommand{\EEe}{\mathbf E}
\newcommand{\eps}{\varepsilon}
\newcommand{\f}{\mathfrak f}
\newcommand{\Ff}{\mathcal F}
\newcommand{\FF}{\mathbb F}
\newcommand{\FFf}{\mathbf F}
\newcommand{\Gg}{\mathcal G}
\newcommand{\GGg}{\mathbf G}
\newcommand{\g}{\mathfrak g}
\newcommand{\Hh}{\mathcal H}
\newcommand{\h}{\mathfrak h}
\newcommand{\Ii}{\mathcal I}
\newcommand{\II}{\mathbb I}
\newcommand{\ka}{\mathfrak k}
\newcommand{\Kk}{\mathcal K}
\newcommand{\kk}{\Bbbk}
\newcommand{\lv}{\mathfrak l}
\newcommand{\Ll}{\mathcal L}
\newcommand{\M}{\mathfrak M}
\newcommand{\m}{\mathfrak m}
\newcommand{\Mm}{\mathcal M}
\newcommand{\N}{\mathfrak N}
\newcommand{\n}{\mathfrak n}
\newcommand{\NN}{\mathbb N}
\newcommand{\Nn}{\mathcal N}
\newcommand{\Oo}{\mathcal O}
\newcommand{\p}{\mathfrak p}
\newcommand{\PP}{\mathbb P}
\newcommand{\Pp}{\mathcal P}
\newcommand{\RR}{\mathbb R}
\newcommand{\Rr}{\mathcal R}
\newcommand{\q}{\mathfrak q}
\newcommand{\QQ}{\mathbb Q}

\newcommand{\s}{\mathfrak s}
\newcommand{\Ss}{\mathcal S}
\newcommand{\tr}{\mathfrak t}
\newcommand{\Tt}{\mathcal T}
\newcommand{\un}{\mathfrak u}
\newcommand{\Uu}{\mathcal U}
\newcommand{\VV}{\mathbb V}
\newcommand{\Vv}{\mathcal V}
\newcommand{\vv}{\mathbf v}
\newcommand{\Ww}{\mathcal W}
\newcommand{\XX}{\mathbb X}
\newcommand{\Xx}{\mathcal X}
\newcommand{\XXx}{\mathbf X}
\newcommand{\Yy}{\mathcal Y}
\newcommand{\YYy}{\mathbf Y}
\newcommand{\ZZ}{\mathbb Z}
\newcommand{\Zz}{\mathcal Z}
\newcommand{\ZZz}{\mathbf Z}
\newcommand{\z}{\mathfrak z}

\newcommand{\LO}{\mathrm{L}}
\newcommand{\GL}{\mathrm{GL}}
\newcommand{\gl}{\mathfrak{gl}}
\newcommand{\SG}{\mathrm{S}}
\newcommand{\SL}{\mathrm{SL}}
\newcommand{\sgl}{\mathfrak{sl}}
\newcommand{\Or}{\mathrm{O}}
\newcommand{\SO}{\mathrm{SO}}
\newcommand{\Sp}{\mathrm{Sp}}
\newcommand{\Spin}{\mathrm{Spin}}
\newcommand{\PSL}{\mathrm{PSL}}
\newcommand{\PO}{\mathrm{PO}}
\newcommand{\PSO}{\mathrm{PSO}}
\newcommand{\PL}{\mathrm{PL}}

\newcommand{\1}{e}
\newcommand{\Alg}{\Aa}
\newcommand{\Aut}{\mathop{\mathrm{Aut}}\nolimits}
\newcommand{\Ad}{\mathop{\mathrm{Ad}}}
\newcommand{\ad}{\mathop{\mathrm{ad}}}
\newcommand{\Ann}{\mathop{\mathrm{Ann}}}
\newcommand{\As}{\mathop{\mathrm{As}}}
\newcommand{\average}[1]{#1^{\natural}}
\newcommand{\bp}{o}
\newcommand{\CA}{\mathop{\mathrm{Cent}}}
\newcommand{\CaDiv}{\mathop{\mathrm{CaDiv}}}
\newcommand{\Cl}{\mathop{\mathrm{Cl}}}
\newcommand{\Cd}{\mathbf{CD}}
\newcommand{\Ch}{\mathfrak{X}}
\newcommand{\CoCh}{\Ch^{*}}
\newcommand{\ch}{\mathop{\mathrm{char}}}
\newcommand{\codim}{\mathop{\mathrm{codim}}}
\newcommand{\const}{\text{\upshape{const}}}
\newcommand{\conv}{\mathop{\mathrm{conv}}}
\newcommand{\cork}{\mathop{\mathrm{cork}}}
\newcommand{\Cox}{\mathop{\mathrm{Cox}}}
\newcommand{\df}{\mathop{\mathrm{def}}}
\newcommand{\diag}{\mathop{\mathrm{diag}}}

\newcommand{\divr}{\mathop{\mathrm{div}}}
\newcommand{\DWC}{\mathbf{C}}
\newcommand{\E}{\textstyle\bigwedge}
\newcommand{\ef}[1]{\mathbf{f}_{#1}}
\newcommand{\ES}{\Ee}
\newcommand{\EHS}[1]{\ES_{#1,+}}
\newcommand{\Env}{\mathop{\mathrm{Env}}}
\newcommand{\Fl}{\mathop{\mathrm{Fl}}}
\newcommand{\Gr}{\mathop{\mathrm{Gr}}}
\newcommand{\gr}{\mathop{\mathrm{gr}}}
\newcommand{\Hilb}{\mathrm{Hilb}}
\newcommand{\Ho}{\mathrm{H}}
\newcommand{\Hom}{\mathop{\mathrm{Hom}}\nolimits}
\newcommand{\Homsh}{\mathop{\mathcal{H}om}\nolimits}
\newcommand{\HS}{\Oo}
\newcommand{\Hyp}{\breve{\ES}}
\newcommand{\Ideal}[1]{\II(#1)}
\newcommand{\id}{\mathbf{1}}
\renewcommand{\Im}{\mathop{\mathrm{Im}}}
\newcommand{\Ind}{\mathop{\mathrm{Ind}}\nolimits}
\newcommand{\ind}[2][\relax]{\Ll_{#1}(#2)}
\newcommand{\idx}{\mathop{\mathrm{ind}}}
\newcommand{\intr}{\mathop{\mathrm{int}}}
\newcommand{\Isom}{\mathop{\mathrm{Isom}}}
\newcommand{\Ker}{\mathop{\mathrm{Ker}}}
\newcommand{\Lin}[2][\relax]{\Oo_{#1}(#2)}
\newcommand{\Mat}{\mathrm{M}}
\newcommand{\Mor}{\mathop{\mathrm{Mor}}\nolimits}
\newcommand{\modl}{\mathop{\mathrm{mod}}}
\newcommand{\nil}{\mathrm n}
\newcommand{\norm}{\text{\upshape{norm}}}
\newcommand{\NE}{\mathop{\mathrm{NE}}}
\newcommand{\NS}{\mathop{\mathrm{NS}}}
\newcommand{\oo}[1]{\mathaccent"7017{#1}}
\newcommand{\ord}{\mathop{\mathrm{ord}}\nolimits}
\newcommand{\Pic}{\mathop{\mathrm{Pic}}\nolimits}
\newcommand{\pr}{\text{\upshape{pr}}}
\newcommand{\PrDiv}{\mathop{\mathrm{PrDiv}}}
\newcommand{\pt}{\text{pt}}
\newcommand{\Quot}{\mathop{\mathrm{Quot}}}
\newcommand{\Rad}[1]{\mathrm{R}(#1)}
\newcommand{\Radu}[1]{\mathrm{R_u}(#1)}
\newcommand{\Rf}{\mathrm{R}}
\newcommand{\red}{\text{\upshape{red}}}
\newcommand{\reg}{\text{\upshape{reg}}}
\newcommand{\res}{\rho}
\newcommand{\Res}{\mathop{\mathrm{Res}}\nolimits}
\newcommand{\RG}{\Lambda}
\newcommand{\rk}{\mathop{\mathrm{rk}}}
\newcommand{\Ru}[1]{#1_{\mathrm{u}}}
\newcommand{\sms}{\mathrm s}
\newcommand{\Sym}{\mathrm{S}}
\newcommand{\Spec}{\mathop{\mathrm{Spec}}\nolimits}
\newcommand{\Supp}{\mathop{\mathrm{Supp}}}
\newcommand{\trdeg}{\mathop{\mathrm{tr.\,deg}}}
\newcommand{\Uni}{\mathrm{U}}
\newcommand{\Vect}{\mathop{\mathrm{Vect}}}
\newcommand{\vol}{\mathop{\mathrm{vol}}}
\newcommand{\X}{\oo{X}}
\newcommand{\Y}{\oo{Y}}
\newcommand{\Zeros}[1]{\VV(#1)}

\newtheorem{theorem}{Theorem}[section]
\newtheorem{proposition}{Proposition}[section]
\newtheorem{lemma}{Lemma}[section]
\newtheorem{corollary}{Corollary}[section]
\newtheorem*{Conjecture}{Conjecture}
\theoremstyle{definition}
\newtheorem{definition}{Definition}[section]
\newtheorem{example}{Example}[section]
\theoremstyle{remark}
\newtheorem{remark}{Remark}[section]
\newtheoremstyle{claim}{\topsep}{\topsep}{\itshape}{}{\bfseries}
                       {.}{.5em}{\thmnote{#3}}
\theoremstyle{claim}
\newtheorem*{Claim}{}

\newcommand{\reftag}[1]{\textup{(#1)}}
\newcommand{\typeA}{A}
\newcommand{\typeB}{B}
\renewcommand{\theenumi}{\textup{(\arabic{enumi})}}
\renewcommand{\labelenumi}{\theenumi}
\newenvironment{roster}{%
  \begin{list}{\labelenumi}{\usecounter{enumi}
    \setlength{\labelwidth}{0pt}
    \setlength{\leftmargin}{0pt}
    \setlength{\itemindent}{\labelsep}}}
 {\end{list}}

\newcounter{tablenum}[table]
\newcommand{\tabitem}{\refstepcounter{tablenum}\thetablenum}
\newcounter{num}
\newcommand{\newnum}{\setcounter{num}{0}}
\newcommand{\num}{\refstepcounter{num}\thenum}

\makeatletter

\newcommand{\Dots}[2]
  {\ifx#10\smash{\vdots}\else
   \ifx#20\smash{\cdots}\else
   \ifx#1#2\smash{\mathinner{%
                  \mkern1mu\raise\p@\hbox{.}
                  \mkern2mu\raise4\p@\hbox{.}\mkern2mu
                  \raise7\p@\vbox{\kern7\p@\hbox{.}}\mkern1mu}}
   \else\smash{\ddots}\fi\fi\fi}

\newcommand{\No}{No.}
\newcommand{\Nos}{Nos.}
\@addtoreset{equation}{section}

\let\c@chapter\c@part
\let\cl@chapter\cl@part
\renewcommand{\thesection}{\arabic{section}}
\renewcommand{\p@section}{\S}

\def\@thefnmark{\textbf{\checkmark}}
\newenvironment{question}
  {\comment}{\endcomment}

\newcommand{\ordinal}[1]
  {\ifcat0#1\ifcase#1\or1-st\or2-nd\or3-rd\else#1-th\fi\else#1-th\fi}

\makeatother

\setlongtables

\begin{document}
\maketitle
\tableofcontents

\chapter*{Introduction}
\addcontentsline{toc}{chapter}{Introduction}
\markboth{\uppercase{Introduction}}{\uppercase{Introduction}}

Groups entered mathematics as transformation groups. From the
works of Caley and Klein it became clear that any geometric theory
studies the properties of geometric objects that are invariant
under the respective transformation group. This viewpoint
culminated in the celebrated Erlangen program~\cite{Erlangen}. An
important feature of each one of the classical
geometries---affine, projective, euclidean, spherical, and
hyperbolic---is that the respective transformation group is
transitive on the underlying space. Another feature of these
examples is that the transformation groups are linear algebraic
and their action is regular. In this way algebraic homogeneous
spaces arise in geometry.

Another source for algebraic homogeneous spaces are varieties of
geometric figures or tensors of certain type. Examples are
provided by Grassmannians, flag varieties, varieties of quadrics,
of triangles, of matrices with fixed rank etc. These homogeneous
spaces are of great importance in algebraic geometry. They were
explored intensively, starting with the works of Chasles,
Schubert, Zeuthen et~al, which gave rise to the enumerative
geometry and intersection theory.

Homogeneous spaces play an important role in representation
theory, since representations of algebraic groups are often
realized in spaces of sections or cohomologies of line (or vector)
bundles over homogeneous spaces. The geometry of a homogeneous
space can be used to study representations of the respective
group, and conversely. A bright example is the Borel--Weil--Bott
theorem~\cite{BWB} and Demazure's proof of the Weyl character
formula~\cite{W-char}.

In the study of an algebraic homogeneous space $G/H$, it is often
useful by standard reasons of algebraic geometry to pass to a
$G$-equivariant completion, or more generally, to an
\emph{embedding}, i.e., a $G$-variety $X$ containing a dense open
orbit isomorphic to~$G/H$.

An example is provided by the following classical problem of
enumerative algebraic geometry: compute the number of plane
quadrics tangent to 5~given ones. Equivalently, one has to compute
the intersection number of certain 5~divisors on the space of
quadrics $\PSL_3/\PSO_3$, which is an open orbit in
$\PP^5=\PP(S^2\CC^3)$. To solve our enumerative problem, we pass
to a good compactification of $\PSL_3/\PSO_3$. Namely, consider
the closure~$X$ in $\PP^5\times(\PP^5)^{*}$ of the graph of a
rational map sending a quadric to the dual one. Points of $X$ are
called \emph{complete quadrics}. It happens that our 5~divisors
intersect the complement of the open orbit in~$X$ properly. Hence
the sought number is just the intersection number of the
5~divisors in~$X$, which is easier to compute, because $X$ is
compact.

Embeddings of homogeneous spaces arise naturally as orbit
closures, when one studies arbitrary actions of algebraic groups.
Such questions as normality of the orbit closure, the nature of
singularities, adherence of orbits, the description of orbits in
the closure of a given orbit etc are of importance.

Embeddings of homogeneous spaces of reductive algebraic groups are
the subject of this survey. The reductivity assumption is natural
for two reasons. First, reductive groups have a good structure and
representation theory, and a deep theory of embeddings can be
developed under this restriction. Secondly, most applications to
algebraic geometry and representation theory deal with homogeneous
spaces of reductive groups. However, homogeneous spaces of
non-reductive groups and their embeddings are also considered.
They arise naturally even in the study of reductive group actions
as orbits of Borel and maximal unipotent subgroups and their
closures. (An example: Schubert varieties.)

The main topics of our survey are:
\begin{itemize}
\item The description of all embeddings of a given homogeneous
space. \item The study of geometric properties of embeddings:
affinity, (quasi)pro\-jec\-tiv\-ity, divisors and line bundles,
intersection theory, singularities etc. \item Application of
homogeneous spaces and their embeddings to algebraic geometry,
invariant theory, and representation theory. \item Determination
of a ``good'' class of homogeneous spaces, for which the above
problems have a good solution. Finding and studying natural
invariants that distinguish this class.
\end{itemize}

Now we describe briefly the content of the survey.

In Chapter~\ref{hom.spaces} we recall basic facts on algebraic
homogeneous spaces and consider basic classes of homogeneous
spaces: affine, quasiaffine, projective. We give group-theoretical
conditions that distinguish these classes. Also bundles and
fibrations over a homogeneous space $G/H$ are considered. In
particular, we compute $\Pic(G/H)$.

In Chapter~\ref{c&r} we introduce and explore two important
numerical invariants of $G/H$---the complexity and the rank. The
\emph{complexity} of $G/H$ is the codimension of a generic
$B$-orbit in~$G/H$, where $B\subseteq G$ is a Borel subgroup. The
\emph{rank} of $G/H$ is the rank of the lattice $\RG(G/H)$ of
weights of rational $B$-eigenfunctions on~$G/H$. These invariants
are of great importance in the theory of embeddings. Homogeneous
spaces of complexity $\leq1$ form a ``good'' class. It was noted
by Howe~\cite{FFT} and Panyushev~\cite{c&r} that a number of
invariant-theoretic problems admitting a nice solution have a
certain homogeneous space of complexity $\leq1$ in the background.

Complexity and rank may be defined for any action $G:X$. We prove
some semicontinuity results for complexity and rank of
$G$-sub\-va\-ri\-eties in~$X$. General methods for computing
complexity and rank of $X$ were developed by Panyushev,
see~\cite[\S\S1--2]{c&r}. We describe them in this chapter, paying
special attention to the case $X=G/H$. The formulas for complexity
and rank are given in terms of the geometry of the doubled action
$G:X\times X^{*}$ and of cotangent bundle $T^{*}X$.

The general theory of embeddings developed by Luna and Vust
\cite{LV} is the subject of Chapter~\ref{LV-theory}. The basic
idea of Luna and Vust is to patch all embeddings $X\embof G/H$
together in a huge prevariety and consider particular embeddings
as Noetherian separated open subsets determined by certain
conditions. It appears, at least for normal embeddings, that $X$
is determined by the collection of closed $G$-subvarieties
$Y\subseteq X$, and each $Y$ is determined by the collection of
$B$-stable divisors containing~$Y$. This leads to a
``combinatorial'' description of embeddings, which can be made
really combinatorial in the case of complexity $\leq1$. In this
case, embeddings are classified by certain collections of convex
polyhedral cones, as in the theory of toric
varieties~\cite{toric.intro} (which is in fact a particular case
of the Luna--Vust theory). The geometry of embeddings is also
reflected in these combinatorial data, as in the toric case. In
fact the Luna--Vust theory is developed here in more generality as
a theory of $G$-varieties in a given birational class (not
necessarily containing an open orbit).

%Affine embeddings of homogeneous spaces are considered in
%Chapter~\ref{aff.emb}. We characterize $G$-orbits on
%$G$-varieties admitting a $G$-stable open affine neighborhood.
%Homogeneous spaces having finitely many $G$-orbits in any
%affine embedding are characterized here. We classify affine
%embeddings of homogeneous spaces in small dimensions. A problem
%of normality of orbit closures is also considered here, in
%particular, for adjoint orbits, $\SL_2$-orbits in linear and
%projective spaces, and torus orbits.

$G$-invariant valuations of the function field of $G/H$ correspond
to $G$-stable divisors on embeddings of~$G/H$. They play a
fundamental role in the Luna--Vust theory as a key ingredient of
the combinatorial data used in the classification of embeddings.
In Chapter~\ref{inv.val} we explore the structure of the set of
invariant valuations, following Knop \cite{G-val}, \cite{inv.mot}.
This set can be identified with a certain collection of convex
polyhedral cones patched together along their common face. This
face consists of \emph{central} valuations---those that are zero
on $B$-invariant functions. It is a solid rational polyhedral cone
in $\RG(G/H)\otimes\QQ$ and a fundamental domain of a
crystallographic reflection group $W(G/H)$, which is called a
\emph{little Weyl group of~$G/H$}. The cone of central valuations
and the little Weyl group are linked with the geometry of the
cotangent bundle.

Spaces of complexity~$0$ form the most remarkable subclass of
homogeneous spaces. Their embeddings are called \emph{spherical
varieties}. They are studied in Chapter~\ref{spherical}.
Grassmannians, flag varieties, determinantal varieties, varieties
of quadrics, of complexes, algebraic symmetric spaces are examples
of spherical varieties. We give several characterizations of
spherical varieties from the viewpoint of algebraic transformation
groups, representation theory, and symplectic geometry. We
consider important classes of spherical varieties: symmetric
spaces, reductive algebraic monoids, horospherical varieties,
toroidal and wonderful varieties. The Luna--Vust theory is much
more developed in the spherical case by Luna, Brion, Knop, Pauer
et~al. We consider the structure of the Picard group of a
spherical variety, the intersection theory and its applications to
enumerative geometry, the cohomology of coherent sheaves, and a
powerful technique of Frobenius splitting, which leads to deep
conclusions on geometry and cohomology of spherical varieties by
reduction to positive characteristic.

%The last chapter is a link from homogeneous spaces to
%representation theory. It concerns \emph{moment polytopes},
%which arise naturally in symplectic geometry (as images of the
%moment map), in representation theory (as certain polytopes
%related to the $G$-module structure of coordinate algebras of
%affine cones), and in the theory of embeddings (for their
%relation to the combinatorial data of embeddings).

The theory of embeddings of homogeneous spaces is rather new and
far from being complete. This survey does not cover all
developments and deeper interactions with other areas. Links for
further reading may be found in the bibliography. We also
recommend the surveys \cite{LV-spher}, \cite{spher.intro},
\cite{spher.full} on spherical varieties and the
monograph~\cite{c&r} on complexity and rank in invariant theory.

A reader is supposed to be familiar with basic concepts of
commutative algebra, algebraic geometry, algebraic groups, and
invariant theory. Our basic sources in these areas are \cite{CA},
\cite{BAG} and \cite{AG}, \cite{Agr}, \cite{IT} and \cite{GIT},
respectively. More special topics are covered by Appendices.

\subsection*{Structure of the survey}

The paper is divided in chapters, chapters are subdivided in
sections, and sections are subdivided in subsections. A link~1.2
refers to Subsection (or Theorem, Lemma, Definition etc) 2 of
Section~1. We try to give sketches of the proofs, if they are not
very long or technical.

\subsection*{Notation and conventions}

We work over an algebraically closed base field~$\kk$. A part of
our results are valid over an arbitrary characteristic, but we
impose the assumption $\ch\kk=0$ whenever it simplifies
formulations and proofs. Let $p$ denote the characteristic
exponent of $\kk$ ($=\ch\kk$, or $1$ if $\ch\kk=0$).

Throughout the paper, $G$~denotes a reductive connected linear
algebraic group, unless otherwise specified. We always may assume
that $G$ is of simply connected type, i.e., a product of a torus
and a simply connected semisimple group. When we study the
geometry of a given homogeneous space $\HS$ and its embeddings, we
often fix a base point $\bp\in\HS$ and denote $H=G_{\bp}$, thus
identifying $\HS$ with~$G/H$.

Algebraic groups are denoted by capital Latin letters, and their
tangent Lie algebras by the respective lowercase Gothic letters.

Topological terms refer to the Zariski topology, unless otherwise
specified.

By a \emph{general point} of an algebraic variety we mean a point
in a certain dense open subset (depending on considered
situation), in contrast with the \emph{generic point}, which is
the dense schematic point of an irreducible algebraic variety.

We use the following general notation.
\begin{trivlist}
\item $A^{\times}$ is the unit group of an algebra~$A$.

\item $\Quot{A}$ is the field of quotients of~$A$.

\item $\kk[S]\subseteq A$ is the subalgebra generated by a subset
$S\subseteq A$.

\item $\Ch(H)$ is the character group of an algebraic group~$H$,
i.e., the group of homomorphisms $H\to\kk^{\times}$ written
additively.

\item $\CoCh(H)$ is the set of (multiplicative) one-parameter
subgroups of~$H$, i.e., homomorphisms $\kk^{\times}\to H$.

\item $H:M$ denotes an action of a group~$H$ on a set~$M$. As a
rule, it is a regular action of an algebraic group on an algebraic
variety.

\item $M^H$ is the set of fixed elements under an action $H:M$.

\item $M^{(H)}$ is the set of all (nonzero) $H$-eigenvectors in a
linear representation $H:M$.

\item $M_{\chi}=M^{(H)}_{\chi}\subseteq M$ is the subspace of
$H$-eigenvectors of the weight~$\chi\in\Ch(H)$.

\item $\kk[X]$ is the algebra of regular functions on an algebraic
variety~$X$.

\item $\kk(X)$ is the field of rational functions on~$X$.

\item $\Lin{\delta}=\Lin[X]{\delta}$ is the line bundle
corresponding to a Cartier divisor $\delta$ on~$X$ or, more
generally, the reflexive sheaf corresponding to a Weil
divisor~$\delta$.

\item $X\by H=\Spec\kk[X]^H$, where $H:X$ is an action of an
algebraic group on an affine variety, and $\kk[X]^H$ is finitely
generated.
\end{trivlist}
Other notation is gradually introduced in the text.

\chapter{Algebraic homogeneous spaces}
\label{hom.spaces}

In this chapter, $G$~denotes an \emph{arbitrary} linear algebraic
group (neither supposed to be connected nor reductive),
$H\subseteq G$ a closed subgroup. We begin in~\ref{basics} with
the definition of an algebraic homogeneous space $G/H$ as a
geometric quotient, and prove its quasiprojectivity. We also prove
some elementary facts on tangent vectors and $G$-equivariant
automorphisms of~$G/H$. In~\ref{bundles}, we describe the
structure of $G$-fibrations over $G/H$ and compute $\Pic(G/H)$.
Some related representation theory is discussed there: induction,
multiplicities, the structure of~$\kk[G]$. Basic classes of
homogeneous spaces are considered in~\ref{classes}. We prove that
$G/H$ is projective iff $H$ is parabolic, and consider criteria of
affinity of~$G/H$.  Quasiaffine $G/H$ correspond to
observable~$H$, which may be defined by several equivalent
conditions (see Theorem~\ref{observ}).

\section{Homogeneous spaces}
\label{basics}

We begin with basic definitions.
\begin{definition}\label{hom.space}
An algebraic group action $G:\HS$ is \emph{transitive} if $\forall
x,y\in\HS\ \exists g\in G:\ y=gx$. In this situation, $\HS$ is
said to be a \emph{homogeneous space}.

A \emph{pointed homogeneous space} is a pair $(\HS,\bp)$, where
$\HS$ is a homogeneous space and $\bp\in\HS$. The natural map
$\pi:G\to\HS$, $g\mapsto g\bp$, is called the \emph{orbit map}.
\end{definition}
A basic property of algebraic group actions is that each orbit is
a locally closed subvariety and thence a homogeneous space in the
sense of Definition~\ref{hom.space}. Homogeneous spaces are always
smooth and quasiprojective, by Sumihiro's Theorem~\ref{Sumihiro}.
\begin{question}
Is Sumihiro's theorem valid in $\ch\kk>0$?
\end{question}
The next definition provides a universal construction of algebraic
homogeneous spaces.
\begin{definition}
The space $G/H$ equipped with the quotient topology and a
structure sheaf $\Oo_{G/H}$ which is the direct image of the sheaf
$\Oo_G^H$ of $H$-invariant (w.r.t.\ the $H$-action on~$G$ by right
translations) regular functions on~$G$ is called the
\emph{(geometric) quotient} of $G$ modulo~$H$.
\end{definition}
\begin{theorem}
\begin{roster}
\item\label{G/H} $(G/H,\Oo_{G/H})$ is a quasiprojective
homogeneous algebraic variety. \item\label{uni} For any pointed
homogeneous space $(\HS,\bp)$ such that $G_{\bp}\supseteq H$, the
orbit map $\pi:G\to\HS$ factors through $\bar\pi:G/H\to\HS$.
\item\label{sep} $\bar\pi$ is an isomorphism iff $G_{\bp}=H$ and
$\pi$ is separable.
\end{roster}
\begin{proof}
To prove~\ref{G/H}, we use the following theorem of
Chevalley~\cite[11.2]{Agr}:
\begin{quote}
There exists a rational $G$-module $V$ and a $1$-dimensional
subspace $L\subseteq V$ such that
%*
\begin{align*}
 H &= N_G(L)=\{g\in G\mid gL=L\} \\
\h &= \n_{\g}(L)=\{\xi\in\g\mid\xi{L}\subseteq L\}
\end{align*}
%*
\end{quote}
Let $x\in\PP(V)$ correspond to~$L$; then it follows that $H=G_x$
and $\h=\Ker d_x\pi$, where $\pi:G\to Gx$ is the orbit map. By a
dimension argument, $d_x\pi$~is surjective, whence $\pi$ is
separable. Further, $Gx$~is homogeneous, whence smooth, and $\pi$
is smooth \cite[III.10.4]{AG}, whence open \cite[Ch.III,
Ex.9.1]{AG}.

Let $U\subseteq Gx$ be an open subset. We claim that each
$f\in\kk[\pi^{-1}(U)]^H$ is the pullback of some $h\in\kk[U]$.
Indeed, consider the rational map $\phi=(\pi,f):G\dasharrow
Gx\times\AAA^1$ and put $Z=\overline{\phi(G)}$. The projection
$Z\to Gx$ is separable and generically bijective, whence
birational. Therefore $f\in\phi^{*}\kk[Z]$ descends to
$h\in\kk(U)$, $f=\pi^{*}h$. If $h$ has the nonzero divisor of
poles~$D$, then $f$ has the nonzero divisor of poles~$\pi^{*}D$, a
contradiction. It follows that $Gx\iso G/H$ is a geometric
quotient.

The universal property~\ref{uni} is an obvious consequence of the
definition. Moreover, any morphism $\phi:G\to Y$ constant on
$H$-orbits factors through $\bar\phi:G/H\to Y$.

Finally, \ref{sep}~follows from the separability of the quotient
map $G\to G/H$: $\pi$~is separable iff $\bar\pi$ is so, and
$G_{\bp}=H$ means that $\bar\pi$ is bijective, whence birational
and, by equivariance and homogeneity, isomorphic.
\end{proof}
\end{theorem}
\begin{remark}\label{insep}
In~\ref{uni}, if $G_{\bp}=H$ and $\pi$ is not separable, then
$\bar\pi$ is bijective purely inseparable and finite~\cite[4.3,
4.6]{Agr}. The schematic fiber $\pi^{-1}(\bp)$ is then a
non-reduced group subscheme of $G$ containing $H$ as the reduced
part. The homogeneous space $\HS$ is uniquely determined by this
subscheme~\cite{gr.sch}.
\end{remark}
\begin{remark}
If $H\normin G$, then $G/H$ is equipped with the structure of a
linear algebraic group with usual properties of the quotient
group. Indeed, in the notation of Chevalley's theorem, we may
assume that $V=\bigoplus_{\chi\in\Ch(G)}V_{\chi}$ and consider the
natural linear action $G:\LO(V)$ by conjugation. The subspace
$E=\prod\LO(V_{\chi})$ of operators preserving each $V_{\chi}$ is
$G$-stable, and the image of $G$ in $\GL(V)$ is isomorphic
to~$G/H$. See \cite[11.5]{Agr} for details.
\end{remark}

Recall that the \emph{isotropy representation} for an action $G:X$
at $x\in X$ is the natural representation $G_x:T_xX$ by
differentials of translations. For a quotient, the isotropy
representation has a simple description:
\begin{proposition}
$T_{\1H}G/H\iso\g/\h$ as $H$-modules.
\end{proposition}
The isomorphism is given by the differential of the (separable)
quotient map $G\to G/H$. The right-hand representation of $H$ is
the quotient of the adjoint representation of $H$ in~$\g$.

Now we describe the group $\Aut_G(G/H)$ of $G$-equivariant
automorphisms of~$G/H$.
\begin{question}
Extend to any $\HS$ (cf.~the proof of
Proposition~\ref{cent.aut(G/H)}). Prove
$\Vect^GG/H=(\g/\h)^H=\n/\h= \mathop{\mathrm{Lie}}\Aut_GG/H$,
$N=N(H)$.
\end{question}
\begin{proposition}
$\Aut_G(G/H)\iso N(H)/H$ is an algebraic group acting on $G/H$
regularly and freely. The action $N(H)/H:G/H$ is induced by the
action $N(H):G$ by right translations: $(nH)(gH)=gn^{-1}H$,
$\forall g\in G,\ h\in N(H)$.
\end{proposition}
\begin{proof}
The regularity of the action $N(H)/H:G/H$ is a consequence of the
universal property of quotients. Clearly, this action is free.
Conversely, if $\phi\in\Aut_G(G/H)$, then $\phi(\1H)=nH$, and
$n\in N(H)$, because the $\phi$-action preserves stabilizers.
Finally, $\phi(gH)=g\phi(\1H)=gnH$, $\forall g\in G$.
\end{proof}

\section{Fibrations, bundles, and representations}
\label{bundles}

The concept of associated bundle is fundamental in topology. We
consider its counterpart in algebraic geometry in a particular
case.

Let $Z$ be an $H$-variety. Then $H$-acts on $G\times Z$ by
$h(g,z)=(gh^{-1},hz)$.
\begin{definition}
The quotient space $G\itimes{H}Z=(G\times Z)/H$ equipped with the
quotient topology and a structure sheaf which is the direct image
of the sheaf of $H$-invariant regular functions is called the
\emph{homogeneous fiber bundle} over $G/H$ associated with~$Z$.
\end{definition}
The $G$-action on $G\times Z$ by left translations of the first
factor commutes with the $H$-action and factors to a $G$-action on
$G\itimes{H}Z$. We denote by $g*z$ the image of $(g,z)$ in
$G\itimes{H}Z$ and identify $\1*z$ with~$z$. The embedding
$Z\embeds G\itimes{H}Z$, $z\mapsto\1*z$, solves the universal
problem for $H$-equivariant morphisms of $Z$ into $G$-spaces.

The homogeneous bundle $G\itimes{H}Z$ is $G$-equivariantly fibered
over $G/H$ with fibers $gZ$, $g\in G$. The fiber map is $g*z\to
gH$. This explains the terminology.

\begin{theorem}[{\cite{ind.act}, \cite[4.8]{IT}}]
If $Z$ is covered by $H$-stable quasiprojective open subsets, then
$G\itimes{H}Z$ is an algebraic $G$-variety, and the fiber map
$G\itimes{H}Z\to G/H$ is locally trivial in \'etale topology.
\end{theorem}
The proof is based on the fact that the fibration $G\to G/H$ is
locally trivial in \'etale topology~\cite{fib.alg}. We shall
always suppose that the assumption of the theorem is satisfied
when we consider homogeneous bundles. The assumption is satisfied,
e.g., if $Z$ is quasiprojective, or normal and $H$ is connected
(by Sumihiro's theorem). If $H$ is reductive and $Z$ is affine,
then $G\itimes{H}Z\iso(G\times Z)\by H$ is affine.

The universal property of homogeneous bundles implies that any
$G$-variety mapped onto $G/H$ is a homogeneous bundle over~$G/H$.
More precisely, a $G$-equivariant map $\phi:X\to G/H$ induces a
bijective $G$-map $G\itimes{H}Z\to X$, where $Z=\phi^{-1}(\1H)$.
If $\phi$ is separable, then $X\iso G\itimes{H}Z$. In particular,
any $G$-subvariety $Y\subseteq G\itimes{H}Z$ is $G$-isomorphic to
$G\itimes{H}(Y\cap Z)$.

Since homogeneous bundles are locally trivial in \'etale topology,
a number of local properties such as regularity, normality,
rationality of singularities etc is transferred from $Z$ to
$G\itimes{H}Z$ and back. The next lemma indicates when a
homogeneous bundle is trivial.
\begin{lemma}\label{triv}
$G\itimes{H}Z\iso G/H\times Z$ if the $H$-action on $Z$ extends to
a $G$-action.
\begin{question}
Is the converse true? No: a principal bundle with unipotent
structure group is trivial.
\end{question}
\end{lemma}
\begin{proof}
The isomorphism is given by $g*z\mapsto(gH,gz)$.
\end{proof}

If the fiber is an $H$-module, then the homogeneous bundle is
locally trivial in Zariski topology.
\begin{question}
Prove it (by taking the top exterior power).
\end{question}
By the above, any $G$-vector bundle over $G/H$ is $G$-isomorphic
to $G\itimes{H}M$ for some finite-dimensional rational
$H$-module~$M$. The respective sheaf of sections $\ind{M}$ is
described in the following way.
\begin{proposition}\label{sheaf.sect}
For any open subset $U\subseteq G/H$, we have
$\Ho^0(U,\ind{M})\iso\Mor_H(\pi^{-1}(U),M)$, where $\pi:G\to G/H$
is the quotient map.
\end{proposition}
\begin{proof}
It is easy to see that the pullback of $G\itimes{H}M\to G/H$
under~$\pi$ is a trivial vector bundle $G\times M\to G$. Hence for
$\forall s\in\Ho^0(U,\ind{M})$ we have
$\pi^{*}s\in\Mor(\pi^{-1}(U),M)$, and clearly $\pi^{*}s$ is
$H$-equivariant. Conversely, any $H$-morphism $\pi^{-1}(U)\to M$
induces a section $U\to G\itimes{H}M$ by the universal property of
the quotient.
\end{proof}
If $H:M$ is an infinite-dimensional rational module, we may
\emph{define} a quasicoherent sheaf $\ind{M}=\ind[G/H]{M}$ on
$G/H$ by the formula of Proposition~\ref{sheaf.sect}. The functor
$\ind[G/H]{\cdot}$ establishes an equivalence between the category
of rational $H$-modules and that of $G$-sheaves on~$G/H$.

Any $G$-line bundle over $G/H$ is $G$-isomorphic to
$G\itimes{H}\kk_{\chi}$, where $\kk_{\chi}=\kk$ with the
$H$-action via a character $\chi\in\Ch(H)$. This yields a
homomorphism $\Ch(H)\to\Pic G/H$,
$\chi\mapsto\ind{\chi}=\ind{\kk_{\chi}}$. Its kernel consists of
characters that correspond to different $G$-linearizations of the
trivial line bundle $G/H\times\kk$ over~$G/H$.

If $G$ is connected, then these characters are exactly the
restrictions to $H$ of characters of~$G$. Indeed, a fiberwise
linear $G$-action on $G/H\times\kk$ is a multiplication by an
algebraic cocycle $c:G\times G/H\to\kk^{\times}$,
$c(g_1g_2,x)=c(g_1,g_2x)c(g_2,x)$ for $\forall g_1,g_2\in G,\ x\in
G/H$. For connected~$G$, we have $c(g,x)=\chi(g)\lambda(x)$,
because an invertible function on a product of two irreducible
varieties is a product of invertible functions on
factors~\cite{Pic_G}. Now it is easy to deduce from the cocycle
property that $\lambda(x)\equiv1$ and $\chi\in\Ch(G)$. Conversely,
if $\chi\in\Res^G_H\Ch(G)$, then $G\itimes{H}\kk_{\chi}\iso
G/H\times\kk_{\chi}$ by Lemma~\ref{triv}.

Consider the universal cover $\widetilde{G}\to G$ (see
Appendix~\ref{rat.mod&lin}). By $\widetilde{H}$ denote the inverse
image of $H$ in $\widetilde{G}$; then
$G/H\iso\widetilde{G}/\widetilde{H}$.
\begin{question}
if $\widetilde{G}\to G$ is separable? In good characteristic?
\end{question}
Since any line bundle over $G/H$ is
$\widetilde{G}$-linearizable~(Corollary~\ref{G^-lin}), we obtain
the following theorem of Popov \cite{Pic(G/H)}, \cite{Pic_G}.
\begin{theorem}
$\Pic_G(G/H)\iso\Ch(H)$. If $G$ is connected, then $\Pic
G/H\iso\Ch(\widetilde{H})/\Res^{\widetilde{G}}_{\widetilde{H}}
\Ch(\widetilde{G})$.
\end{theorem}
(Here $\Pic_G$ denotes the group of $G$-linearized invertible
sheaves.)
\begin{example}
Let $G$ be a connected reductive group, $B\subseteq G$ a Borel
subgroup. Then $\Pic G/B$ is isomorphic to the weight lattice of
the root system of~$G$.
\end{example}

Let $X$ be a $G$-variety, and $Z\subseteq X$ an $H$-stable
subvariety. By the universal property, we have a $G$-equivariant
map $\mu:G\itimes{H}Z\to X$, $\mu(g*z)=gz$.
\begin{proposition}\label{collapse}
If $H$ is parabolic, then $\mu$ is proper and $GZ$ is closed
in~$X$.
\end{proposition}
\begin{proof}
The map $\mu$ factors as
$\mu:G\itimes{H}Z\stackrel{\iota}{\embeds}G\itimes{H}X \iso
G/H\times X$ (Lemma~\ref{triv}) ${}\stackrel{\pi}{\to}X$, where
$\iota$ is a closed embedding and $\pi$ is a projection along a
complete variety by Theorem~\ref{parab}.
\end{proof}
\begin{example}
Let $\N\subseteq\g$ be the set of nilpotent elements and
$U=\Radu{B}$, a maximal unipotent subgroup of~$G$. Then the map
$G\itimes{B}\un\to\N$ is proper and birational, see, e.g.,
\cite[5.6]{IT}. This is a well-known Springer's resolution of
singularities of~$\N$.
\end{example}

Now we discuss some representation theory related to homogeneous
spaces and to vector bundles over them.

We always deal with rational modules over algebraic groups (see
Appendix~\ref{rat.mod&lin}) and often drop the word ``rational''.
As usual in representation theories, we may define functors of
induction and restriction on categories of rational modules. Let
$H$ act on $G$ by right translations, and $M$ be an $H$-module.
\begin{definition} A $G$-module
$\Ind^G_HM=\Mor_H(G,M)\iso(\kk[G]\otimes M)^H$ is said to be
\emph{induced} from $H:M$ to~$G$. It is a rational
$G$-$\kk[G/H]$-module. By definition, we have
$\Ind^G_HM=\Ho^0(G/H,\ind{M})$.

A $G$-module $N$ considered as an $H$-module is denoted by
$\Res^G_HN$.
\end{definition}
\begin{example}\label{Ind}
$\Ind^G_H\kk=\kk[G/H]$, where $\kk$ is the trivial $H$-module.
More generally, $\Ind^G_H\kk_{\chi}=\kk[G]_{-\chi}$,
$\forall\chi\in\Ch(H)$.
\end{example}

Clearly, $\Ind^G_H$ is a left exact functor from the category of
rational $H$-modules to that of rational $G$-modules. The functor
$\Res^G_H$ is exact. We collect basic properties of induction
\begin{question}
What about the universal property, transitivity?
\end{question}
in the following
\begin{theorem}
\begin{roster}
\item\label{Ind(G:M)} If $M$ is a $G$-module, then
$\Ind^G_HM\iso\kk[G/H]\otimes M$. \item\label{Frob}(Frobenius
reciprocity) For rational modules $G:N$, $H:M$, we have
%*
\begin{equation*}
\Hom_G(N,\Ind^G_HM)\iso\Hom_H(\Res^G_HN,M)
\end{equation*}
%*
\item\label{Ind^G} For any $H$-module $M$, $(\Ind^G_HM)^G\iso
M^H$. \item\label{Ind(alg)} If $M,N$ are rational algebras, then
\ref{Ind(G:M)} and \ref{Ind^G} are isomorphisms of algebras, and
\ref{Frob} holds for equivariant algebra homomorphisms.
\end{roster}
\end{theorem}
\begin{proof}
\begin{roster}
\item[\ref{Ind(G:M)}] The isomorphism
$\iota:\Mor_H(G,M)\isoto\Mor(G/H,M)$ is given by
$\iota(m)(gH)=g\cdot m(g)$, $\forall m\in\Mor_H(G,M)$. The inverse
mapping is $\mu\mapsto m$, $m(g)=g^{-1}\mu(gH)$,
$\forall\mu\in\Mor(G/H,M)$. \item[\ref{Frob}] The isomorphism is
given by the map $\Phi\mapsto\phi$, $\forall\Phi:N\to\Mor_H(G,M)$,
where $\phi:N\to M$ is defined by $\phi(n)=\Phi(n)(\1)$, $\forall
n\in N$. The inverse map $\phi\mapsto\Phi$ is given by
$\Phi(n)(g)=\phi(g^{-1}n)$. \item[\ref{Ind^G}] Any $G$-invariant
$H$-equivariant morphism $G\to M$ is constant, and its image lies
in~$M^H$. Alternatively, one may apply the Frobenius reciprocity
to $N=\kk$. \item[\ref{Ind(alg)}] It is easy. \qedhere\end{roster}
\end{proof}
\begin{remark}\label{transfer}
The union of \ref{Ind(G:M)} and \ref{Ind^G} yields the following
assertion: if $M$ is a $G$-module, then $(\kk[G/H]\otimes M)^G\iso
M^H$. This is often called the \emph{transfer principle}, because
it allows to transfer information from $\kk[G/H]$ to~$M$. For
example, if $G$ is reductive, $\kk[G/H]$~is finitely generated,
and $M=A$ is a finitely generated $G$-algebra, then $A^H$ is
finitely generated. Other applications are discussed below.
\begin{question}
Which?
\end{question}
A good treatment of induced modules and the transfer principle can
be found in~\cite{HS}.
\end{remark}

We are interested in the $G$-module structure of $\kk[G/H]$ and of
global sections of line bundles over~$G/H$.

For any two rational $G$-modules $V,M$ ($\dim V<\infty$), put
%*
\begin{equation*}
m_V(M)=\dim\Hom_G(V,M),
\end{equation*}
%*
the \emph{multiplicity} of $V$ in~$M$. If $V$ is simple and $M$
completely reducible (e.g., $G$~is an algebraic torus or a
reductive group in characteristic zero), then $m_V(M)$ is the
number of occurrences of $V$ in a decomposition of $M$ into simple
summands.

For any $G$-variety $X$ and a $G$-line bundle $\Ll\to X$, we
abbreviate:
%*
\begin{align*}
m_V(X)&=m_V(\kk[X]),& m_V(\Ll)&=m_V(\Ho^0(X,\Ll))
\end{align*}
%*
Here is a particular case of Frobenius reciprocity:
\begin{corollary}\label{m(G/H)}
$m_V(G/H)=\dim(V^{*})^H$,
$m_V(\ind{\chi})=\dim(V^{*})^{(H)}_{-\chi}$
\end{corollary}
\begin{proof}
We have $\Ho^0(G/H,\ind{\chi})=\Ind^G_H\kk_{\chi}$, whence
%*
\begin{equation*}
\Hom_G\bigl(V,\Ho^0(G/H,\ind{\chi})\bigr)=\Hom_H(V,\kk_{\chi})=
(V^{*})^{(H)}_{-\chi}
\end{equation*}
%*
The first equality follows by taking $\chi=0$.
\end{proof}

A related problem is to describe the module structure of~$\kk[G]$.
Namely, $G$~itself is acted on by $G\times G$ via
$(g_1,g_2)g=g_1gg_2^{-1}$. Hence $\kk[G]$ is a $(G\times
G)$-algebra.

Every finite-dimensional $G$-module $V$ generates a $(G\times
G)$-stable subspace $\Mat(V)\subset\kk[G]$ spanned by matrix
entries $f_{\omega,v}(g)=\langle\omega,gv\rangle$ ($v\in V$,
$\omega\in V^{*}$) of the representation $G\to\GL(V)$. Clearly
$\Mat(V)$ is the image of a $(G\times G)$-module homomorphism
$V^{*}\otimes V\to\kk[G]$, $\omega\otimes v\mapsto f_{\omega,v}$,
and $\Mat(V)\iso V^{*}\otimes V$ is a simple $(G\times G)$-module
whenever $V$ is simple.

Matrix entries behave well w.r.t.\ algebraic operations:
%*
\begin{align}\label{mat}
\Mat(V)+\Mat(V')&=\Mat(V\oplus V'),&
\Mat(V)\cdot\Mat(V')&=\Mat(V\otimes V')
\end{align}
%*
The inversion of $G$ sends $\Mat(V)$ to~$\Mat(V^{*})$.

\begin{proposition}\label{U(Mat)}
$\kk[G]=\bigcup\Mat(V)$, where $V$ runs through all
finite-di\-men\-sional $G$-modules.
\end{proposition}
\begin{proof}
Take any finite-dimensional $G$-submodule $V\subset\kk[G]$ w.r.t.\
the $G$-action by right translations. We claim $V\subset\Mat(V)$.
Indeed, let $\omega\in V^{*}$ be defined by
$\langle\omega,v\rangle=v(\1)$, $\forall v\in V$; then $\forall
v\in V, g\in G:\ v(g)=f_{\omega,v}(g)$.
\end{proof}

\begin{theorem}\label{reg.rep}
Suppose $\ch\kk=0$ and $G$ is reductive. Then there is a $(G\times
G)$-module isomorphism
%*
\begin{equation*}
\kk[G]=\bigoplus\Mat(V)\iso\bigoplus V^{*}\otimes V
\end{equation*}
%*
where $V$ runs through all simple $G$-modules.
\end{theorem}
\begin{proof}
All the $\Mat(V)\iso V^{*}\otimes V$ are pairwise non-isomorphic
simple $(G\times G)$-modules. By Proposition~\ref{U(Mat)} and
\eqref{mat} they span the whole~$\kk[G]$.
\end{proof}
\begin{remark}
Corollary~\ref{m(G/H)} can be derived from Theorem~\ref{reg.rep}
by taking $H$-(semi)invariants from the right.
\end{remark}

The dual object to the coordinate algebra of $G$ provides a
version of the group algebra for algebraic groups.
\begin{definition}
The (\emph{algebraic}) \emph{group algebra} of $G$ is
$\Aa(G)=\kk[G]^{*}$ equipped with the multiplication law coming
from the comultiplication in~$\kk[G]$.
\end{definition}
For finite $G$ we obtain the usual group algebra. Generally,
$\Aa(G)$ can be described by finite-dimensional approximations.
The group algebra $\Aa(V)$ of a finite-dimensional $G$-module $V$
is defined as the linear span of the image of $G$ in~$\LO(V)$.
Note that $\Aa(V)$ is the $(G\times G)$-module dual to~$\Mat(V)$.
We have $\Aa(V)=\LO(V)$ whenever $V$ is simple. Given a
subquotient module $V'$ of~$V$, there is a canonical epimorphism
$\Aa(V)\onto\Aa(V')$. Therefore the algebras $\Aa(V)$ form an
inverse system over all $V$ ordered by the relation of being a
subquotient. It readily follows from Proposition~\ref{U(Mat)} that
$\Aa(G)\iso\varprojlim\Aa(V)$. One deduces that $\Aa(G)$ is a
universal ambient algebra containing both $G$ and~$\Uni\g$, the
(restricted) envelopping algebra of $\g$ \cite{gr.alg}.

\begin{definition}
The algebra $\Aa(G/H)$ of all $G$-equivariant linear endomorphisms
of $\kk[G/H]$ is called the \emph{Hecke algebra} of~$G/H$, or of
$(G,H)$.
\end{definition}

\begin{remark}
If $\ch\kk=0$ and $G$ is reductive, then $\Aa(V)=\prod\LO(V_i)$
over all simple $G$-modules $V_i$ occurring in $V$ with positive
multiplicity. Furthermore, $\Aa(G)=\prod\LO(V_i)$ and
$\Aa(G/H)=\prod\LO(V_i^H)$ over all simple~$V_i$ by
Theorem~\ref{reg.rep} and Schur's lemma.
\end{remark}

\begin{proposition}[E.~B.~Vinberg]\label{Hecke=biinv}
If $\ch\kk=0$ and $H$ is reductive, then
$\Aa(G/H)\iso\Aa(G)^{H\times H}$. In particular, the above
notation is compatible for $H=\{\1\}$.
\end{proposition}
\begin{proof}
First consider the case $H=\{\1\}$. The algebra $\Aa(V)$ acts on
$\Aa(V)^{*}=\Mat(V)$ by right translations: $af(x)=f(xa)$,
$\forall a,x\in\Aa(V),\ f\in\Mat(V)$. These actions commute with
the $G$-action by left translations and merge together into a
$G$-equivariant linear $\Aa(G)$-action on~$\kk[G]$.

Conversely, every $G$-equivariant linear map
$\phi:\kk[G]\to\kk[G]$ preserves all the spaces $\Mat(V)$. Indeed,
it follows from the proof of Proposition~\ref{U(Mat)} by applying
the inversion that $W\subseteq\Mat(W^{*})$ for any $G$-submodule
$W\subset\kk[G]$. For $W=\Mat(V)$ one easily deduces
$W=\Mat(W^{*})$ and
$\phi{W}\subseteq\Mat(\phi{W}^{*})\subseteq\Mat(W^{*})$.

The restriction of $\phi$ to $\Mat(V)$ is the right translation by
some $a_V\in\Aa(V)$. These $a_V$ give rise to $a\in\Aa(G)$
representing $\phi$ on~$\kk[G]$. Hence the group algebra coincides
with the Hecke algebra of~$G$.

In the general case, every linear $G$-endomorphism $\phi$ of
$\kk[G]^H$ extends to a unique $a\in\Aa(G)^{H\times H}$, which
annihilates the right-$H$-invariant complement of $\kk[G]^H$
in~$\kk[G]$.
\end{proof}

If $G$ is a connected reductive group, $B\subseteq G$ a Borel
subgroup, then isomorphism classes of simple $G$-modules are
indexed by $B$-\emph{dominant weights}, which form a subsemigroup
$\Ch_{+}\subseteq\Ch(B)$ (the intersection of $\Ch(B)$ with the
positive Weyl chamber). Any simple $G$-module $V$ contains a
unique, up to proportionality, $B$-eigenvector (a \emph{highest
vector}) of weight $\lambda\in\Ch_{+}$ (the \emph{highest weight})
\cite[\S31]{Agr}. The highest weight of $V^{*}$ is
$\lambda^{*}=-w_G\lambda$, where $w_G$ is the longest element of
the Weyl group.

By Corollary~\ref{m(G/H)},
%*
\begin{equation*}
m_V(\ind[G/B]{\mu})=\dim(V^{*})^{(B)}_{-\mu}=
\begin{cases}
1,&\mu=-\lambda^{*} \\
0,&\text{otherwise}
\end{cases}
\end{equation*}
%*
It follows that $V^{*}(\lambda)=\Ind^G_B\kk_{-\lambda}$ contains a
unique simple $G$-module (of highest weight~$\lambda^{*}$)
whenever $\lambda\in\Ch_{+}$, otherwise $V^{*}(\lambda)=0$. The
dual $G$-module $V(\lambda)=(\Ind^G_B\kk_{-\lambda})^{*}$ is
called a \emph{Weyl module} \cite[II.2]{Jan}.

Put $m_{\lambda}(M)=m_{V(\lambda)}(M)$ for brevity.
\begin{proposition}\label{mult(h.w)}
$m_{\lambda}(M)=\dim M^{(B)}_{\lambda}$
\end{proposition}
\begin{proof}
As $G/B$ is a projective variety (Theorem~\ref{parab}),
$V(\lambda)=\Ho^0(G/B,\ind{-\lambda})^{*}$ is finite-dimensional.
If $\dim M<\infty$, then
%*
\begin{equation*}
\Hom_G(V(\lambda),M)\iso\Hom_G(M^{*},V^{*}(\lambda))
\iso\Hom_B(M^{*},\kk_{-\lambda})\iso M^{(B)}_{\lambda}
\end{equation*}
%*
However, any rational $G$-module $M$ is a union of
finite-dimensional submodules.
\end{proof}
Thus $V(\lambda)$ can be characterized as the universal covering
$G$-module of highest weight~$\lambda$: the generating highest
vector in $V(\lambda)$ is given by evaluation of
$\Ho^0(G/B,\ind{-\lambda})$ at~$\1B$.

By Corollary~\ref{m(G/H)}, we have
%*
\begin{align}\label{mult(G/H)}
m_{\lambda}(G/H)&=\dim V^{*}(\lambda)^H,&
m_{\lambda}(\ind[G/H]{\chi})&=\dim V^{*}(\lambda)^{(H)}_{-\chi}.
\end{align}
%*

In characteristic zero, complete reducibility yields:
\begin{Claim}[Borel--Weil theorem]
If $\ch\kk=0$, then $V(\lambda)$ is a simple $G$-module of highest
weight $\lambda$ and $V^{*}(\lambda)\iso V(\lambda^{*})$.
\end{Claim}

Furthermore, Theorem~\ref{reg.rep} yields
%*
\begin{equation}\label{k[G]}
\kk[G]\iso \bigoplus_{\lambda\in\Ch_{+}}V(\lambda^{*})\otimes
V(\lambda)
\end{equation}
%*
In arbitrary characteristic, Formula~\eqref{k[G]} is no longer
true, but $\kk[G]$ possesses a ``good'' $(G\times G)$-module
filtration with factors $V^{*}(\lambda)\otimes V^{*}(\lambda^{*})$
\cite{Don}, \cite[II.4.20]{Jan}.

Notice that all the dual Weyl modules are combined in a
multigraded algebra
%*
\begin{equation}\label{k[G/U]}
\kk[G/U]=\bigoplus_{\lambda\in\Ch_{+}}\kk[G]^{(B)}_{\lambda}
\iso\bigoplus_{\lambda\in\Ch_{+}}V^{*}(\lambda)
\end{equation}
%*
called the \emph{covariant algebra} of~$G$. Here, as above,
$U=\Radu{B}$. The covariant algebra is an example of a
multiplicity free $G$-algebra, in the sense of the following
\begin{definition}
A $G$-module $M$ is said to be \emph{multiplicity-free} if
$m_{\lambda}(M)\leq1$, $\forall\lambda\in\Ch_{+}$.
\end{definition}
The multiplication in the covariant algebra has a nice property:
\begin{lemma}[{\cite[II.14.20]{Jan}}]\label{V(l)*V(m)}
$V^{*}(\lambda)\cdot V^{*}(\mu)=V^{*}(\lambda+\mu)$
\end{lemma}
The inclusion ``$\subseteq$'' in the lemma is obvious since the
$V^{*}(\lambda)$ are the homogeneous components of $\kk[G/U]$
w.r.t. an algebra grading. In characteristic zero, the reverse
inclusion stems from the fact that the $V^{*}(\lambda)$ are simple
$G$-modules and $\kk[G/U]$ is an integral domain.

\section{Classes of homogeneous spaces}
\label{classes}

We answer the following question: when is a homogeneous space
$\HS$ projective or (quasi)affine? First, we reduce the question
to a property of the pair $(G,H)$, where $H=G_{\bp}$ is the
stabilizer of a point $\bp\in\HS$.
\begin{lemma}\label{orb&quot}
$\HS$ is projective, resp.\ (quasi)affine iff $G/H$ has this
property.
\end{lemma}
\begin{proof}
We may assume $\ch\kk=p>0$. The natural map $\iota:G/H\to\HS$ is
finite bijective purely inseparable (Remark~\ref{insep}). For
completeness and affinity, we conclude by~\cite[III, ex.4.2]{AG}.

For quasiaffinity, we argue as follows. First note that
$\kk(G/H)^{p^s}\subseteq\iota^{*}\kk(\HS)$ for some $s\geq0$.
Furthermore, $\kk[G/H]^{p^s}\subseteq\iota^{*}\kk[\HS]$ (If
$f\in\kk(\HS)$, then $\iota^{*}f$ has poles on~$G/H$.) Assume
$\HS$ is an open subset of an affine variety~$Y$. Let $B$ be the
integral closure of $\iota^{*}\kk[Y]$ in~$\kk(G/H)$. Then $B$ is
finitely generated, and $X=\Spec B$ contains $G/H$ as an open
subset. Conversely, if $G/H$ is open in an affine variety~$X$,
then $A=\kk[X]\cap\iota^{*}\kk(\HS)$ is finite over
$\kk[X]^{p^s}$. Hence $A$ is finitely generated, and $X=\Spec A$
contains $\HS$ as an open subset.
\end{proof}

In the sequel, we may assume $\HS=G/H$.
\begin{lemma}\label{transit}
If $G\supseteq H\supseteq K$ and $G/H$, $H/K$ are projective,
resp.\ (quasi)\-affine, then $G/K$ is projective, resp.\
(quasi)affine.
\end{lemma}
\begin{proof}
The natural map $\phi:G/K\to G/H$ transforms after a faithfully
flat base change $G\to G/H$ to the projection
$\pi:G/K\times_{G/H}G\iso H/K\times G\to G$. If $H/K$ is
projective (resp.\ affine), then $\pi$ is proper (resp.\ affine),
whence $\phi$ is proper (resp. affine). If in addition $G/H$ is
complete (resp.\ affine), then $G/K$ is complete (resp.\ affine),
too. Another proof for projective and affine cases relies on
Theorems \ref{parab} and~\ref{exact} below. In the quasiaffine
case, lemma follows from Theorem~\ref{observ}\ref{rep}.
\end{proof}

\begin{theorem}\label{parab}
$G/H$ is projective iff $H$ is parabolic.
\end{theorem}
\begin{proof}
If $G/H$ is quasiprojective, then a Borel subgroup $B\subseteq G$
has a fixed point $gH\in G/H$, by the Borel fixed point
Theorem~\cite[21.2]{Agr}. Hence $H\supseteq g^{-1}Bg$ is
parabolic.

To prove the converse, consider an exact representation $G:V$. The
induced action of $G$ on the variety of complete flags in $V$ has
a closed orbit. Its stabilizer $B$ is solvable, and we may assume
$B\subset H$. By Lemma~\ref{orb&quot}, $G/B$ is complete, hence
$G/H$ is complete.
\end{proof}

A group-theoretical characterization of affine homogeneous spaces
is not known at the moment. We give several sufficient conditions
of affinity and a criterion for reductive~$G$.

\begin{lemma}\label{uni=>closed}
The orbits of a unipotent group $G$ on an affine variety $X$ are
closed, whence affine.
\end{lemma}
\begin{proof}
For any $x\in X$, consider closed affine subvarieties
$Y=\overline{Gx}\subseteq X$ and $Z=Y\setminus Gx$. Since
$\II(Z)\normin\kk[Y]$ is a $G$-submodule, the Lie--Kolchin theorem
implies $\exists f\in\II(Z)^G$, $f\neq0$. However $f$ is a nonzero
constant on~$Gx$, whence on~$Y$. Thus $Z=\emptyset$.
\end{proof}

\begin{theorem}\label{solv=>aff}
$G/H$ is affine if $G$ is solvable.
\end{theorem}
\begin{proof}
We may assume that $G,H$ are connected. First suppose $G$ is
unipotent. Take a representation $G:V$ such that $\exists v\in V:\
G[v]\iso G/H$. Then $H$ normalizes $\langle v\rangle$. But
$\Ch(H)=0$, whence $G_v=H$ and $Gv\iso G/H$. We conclude by
Lemma~\ref{uni=>closed}.

In the general case, $G=T\semitimes U$ and $H=S\semitimes V$,
where $U,V$ are unipotent radicals and $T,S$ are maximal tori of
$G,H$. We have $U\supset V$ and may assume that $T\supset S$. It
is easy to see that $G/H\iso T\itimes{S}U/V=(T\times U/V)\by S$ is
affine.
\end{proof}

The following notion is often useful in the theory of homogeneous
spaces.
\begin{definition}
We say that $H$ is \emph{regularly embedded} in~$G$ if
$\Radu{H}\subseteq\Radu{G}$.
\end{definition}
For example, any subgroup of a solvable group is regularly
embedded. The next theorem generalizes Theorem~\ref{solv=>aff}.
\begin{theorem}
$G/H$ is affine if $H$ is regularly embedded in~$G$.
\end{theorem}
\begin{proof}
As $\Radu{G}$ is normal in~$G$, the quotient $G/\Radu{G}$ is
affine. By Theorem~\ref{solv=>aff}, $\Radu{G}/\Radu{H}$ is affine.
Thence by Lemma~\ref{transit}, $G/\Radu{H}$ is affine. By the Main
Theorem of GIT (see Appendix~\ref{inv.th}), $G/H=(G/\Radu{H})\by
(H/\Radu{H})$ is affine, because $H/\Radu{H}$ is reductive.
\end{proof}
Weisfeiler proved~\cite{Wei} that any subgroup $H$ of a connected
group $G$ is regularly embedded in some parabolic subgroup
$P\subseteq G$. (See also~\cite[30.3]{Agr}.) Thus $G/H$ is a
fibration with the projective base $G/P$ and affine fiber $P/H$.

The following theorem is often called Matsushima's criterion. It
was proved for $\kk=\CC$ by Matsushima~\cite{Mat} and
Onishchik~\cite{Oni}, and in the general case by
Richardson~\cite{Ri}.
\begin{theorem}
$G/H$ is affine if $H$ is reductive. If $G$ is reductive, the
converse is also true.
\end{theorem}
\begin{proof}
If $H$ is reductive, then by the Main Theorem of GIT, $G/H\iso
G\by H$ is affine. A simple proof of the converse see
in~\cite{slices} ($\ch\kk=0$) or in~\cite{Ri}.
\end{proof}

The lacking of a group-theoretical criterion of affinity is
partially compensated by a cohomological criterion.
\begin{theorem}\label{exact}
$G/H$ is affine iff $\Ind^G_H$ is exact.
\end{theorem}
\begin{proof}
Recall that $\Ind^G_H(M)=\Ho^0(G/H,\ind{M})$, the sheaf $\ind{M}$
is quasicoherent, and the functor $\ind{\cdot}$ is exact. If $G/H$
is affine, then by Serre's criterion,
 $\Ind^G_H$ is exact. For a proof of the converse,
see~\cite[\S6]{HS}.
\end{proof}

The class of quasiaffine homogeneous spaces is of interest in
invariant theory and representation theory. If $G/H$ is
quasiaffine, then the subgroup $H$ is called \emph{observable}.
Observable subgroups are exactly the stabilizers of vectors in
rational $G$-modules, since any quasiaffine $G$-variety can be
equivariantly embedded in a $G$-module~\cite[1.2]{IT}.
\begin{example}
By Chevalley's theorem, $H$ is observable if $\Ch(H)=0$. In
particular a unipotent subgroup is observable.
\end{example}
\begin{example}[{\cite{observ}}]
If $\Rad{H}$ is nilpotent, then $H$ is observable.
\end{example}

It is easy to see that an intersection of observable subgroups is
again observable. Therefore for any $H\subseteq G$, there exists a
smallest observable subgroup $\widehat{H}\subseteq G$
containing~$H$. It is called the \emph{observable hull} of~$H$.
Clearly, for any rational $G$-module $M$ we have
$M^H=M^{\widehat{H}}$. This property illustrates the importance of
observable subgroups in invariant theory, see~\cite[3.7]{IT}.

We give several characterizations of observable subgroups in the
next theorem, essentially due to Bia{\l}ynicki-Birula, Hochschild,
and Mostow~\cite{observ}.
\begin{theorem}\label{observ}
The following conditions are equivalent:
\begin{enumerate}
\item\label{qaff} $G/H$ is quasiaffine. \item\label{qaff^0}
$G^0/H^0$ is quasiaffine. \item\label{quot}
$\Quot\kk[G/H]=\kk(G/H)$. \item\label{rep} Any finite-dimensional
$H$-module is embedded as an $H$-submodule in a finite-dimensional
$G$-module. \item\label{char} $\forall\chi\in\Ch(H):\
\kk[G]_{\chi}\neq0 \implies\kk[G]_{-\chi}\neq0$ (In other words,
the semigroup of weights of $H$-eigenfunctions on $G$ is actually
a group.)
\end{enumerate}
\end{theorem}
\begin{proof}
\begin{roster}
\item[\ref{qaff}$\implies$\ref{quot}] is obvious.
\item[\ref{quot}$\implies$\ref{qaff}] We have
$\kk(G/H)=\kk(G/\widehat{H})$, whence $H=\widehat{H}$.
\item[\ref{qaff}$\iff$\ref{qaff^0}] We may assume that $G$ is
connected. The map $G/H^0\to G/H$ is a Galois covering with the
Galois group $\Gamma=H/H^0$. If $G/H$ is open in an affine
variety~$X$, then $G/H^0$ is open in $Y=\Spec A$, where $A$ is the
integral closure of $\kk[X]$ in~$\kk(G/H^0)$. Conversely, if
$G/H^0$ is open in affine $Y$, then $G/H$ is open in
$X=\Spec\kk[Y]^{\Gamma}$. \item[\ref{qaff}$\implies$\ref{char}]
For a nonzero $p\in\kk[G]_{\chi}$, consider its zero set $Z\subset
G$. The quotient morphism $\pi:G\to G/H$ maps $Z$ onto a proper
closed subset of~$G/H$. Hence $\exists f\in\kk[G/H]:\
f|_{\pi(Z)}=0$. By Nullstellensatz, $\pi^{*}f^n=pq$ for some
$n\in\NN,\ q\in\kk[G]_{-\chi}$.
\item[\ref{char}$\implies$\ref{rep}] First note that a
$1$-dimensional $H$-module $W=\kk_{\chi}$ can be embedded in a
$G$-module $V$ iff $\kk[G]_{\chi}\neq0$. (Any function
$f(g)=\langle\omega,gv\rangle$, where $w\in W$, $\omega\in V^{*}$,
belongs to~$\kk[G]_{\chi}$.)

Now for any finite-dimensional $H$-module $W$, consider the
embedding $W\embeds\Mor(H,W)$ taking each $w\in W$ to the orbit
morphism $g\mapsto gw$, $g\in H$. It is $H$-equivariant w.r.t.\
the $H$-action on $\Mor(H,W)$ by right translations of an
argument. The restriction of morphisms yields a projection
$\Mor(G,W)\to\Mor(H,W)$, and we may choose a finite-dimensional
$H$-submodule $N\subset\Mor(G,W)$ mapped onto~$W$. Embed $N$ into
a finite-dimensional $G$-submodule $M\subset\Mor(G,W)$ and put
$U=\Ker(N\to W)$, $m=\dim U$. Then $\E^mU\embeds\E^mM$ and
$W\otimes\E^mU\embeds\E^{m+1}M$. By~\ref{char} and the above
remark, $\E^mU^{*}$ is embedded in a $G$-module. We conclude by
$W\iso(W\otimes\E^mU)\otimes\E^mU^{*}$.
\item[\ref{rep}$\implies$\ref{qaff}] By Chevalley's theorem, $H$
is the projective stabilizer of a vector $v$ in some
$G$-module~$V$. As an $H$-module, $\langle v\rangle\iso\kk_{\chi}$
for some $\chi\in\Ch(H)$. Then $\kk_{-\chi}$ can be embedded in a
$G$-module, i.e., $\exists G:W,\ w\in W$ such that $H$ acts on $w$
via~$-\chi$. It follows that $G_{v\otimes w}=H$.
\qedhere\end{roster}
\end{proof}

Surprisingly, quasiaffine homogeneous spaces admit a
group-theoretical characterization. Recall that a
\emph{quasiparabolic} subgroup of a connected group is the
stabilizer of a highest weight vector in an irreducible
representation.
\begin{theorem}[{\cite{subpar}}]\label{subpar}
$H\subseteq G$ is observable iff $H^0$ is regularly embedded in a
quasiparabolic subgroup of~$G^0$.
\end{theorem}

\chapter{Complexity and rank}
\label{c&r}

We retain general conventions of our survey. In particular,
$G$~denotes a reductive connected linear algebraic group. We begin
with local structure theorems, which claim that a $G$-variety may
be covered by affine open subsets stable under parabolic subgroups
of~$G$, and describe the structure of these subsets.
In~\ref{c&r(basics)}, we define two numerical invariants of a
$G$-variety related to the action of a Borel subgroup of~$G$---the
complexity and the rank. We reduce their computation to a generic
orbit on~$X$ (i.e., a homogeneous space) and prove some basic
results including the semicontinuity of complexity and rank
w.r.t.\ $G$-subvarieties. We also introduce the notion of the
weight lattice and consider the connection of complexity with the
growth of multiplicities in $\kk[X]$ for quasiaffine~$X$. The
relation of complexity and modality of an action is considered
in~\ref{modality}. In~\ref{horosph}, we introduce the class of
horospherical varieties defined by the property that all isotropy
groups contain a maximal unipotent subgroup of~$G$. The
computation of complexity and rank is fairly simple for them. On
the other hand, any $G$-variety can be contracted to a
horospherical one of the same complexity and rank.

General formulae for complexity and rank are obtained
in~\ref{cotangent} as a by-product of the study of the cotangent
action $G:T^{*}X$ and the doubled action $G:X\times X^{*}$. These
formulae involve generic stabilizers of these actions. The
particular case of a homogeneous space $X=G/H$ is considered
in~\ref{c&r(hom)}. In~\ref{c&r<=1}, we classify homogeneous spaces
of complexity and rank~$\leq1$. An application to the problem of
decomposing tensor products of representations is considered
in~\ref{double.cone}. Decomposition formulae are obtained from the
description of the $G$-module structure of coordinate algebras on
double cones of small complexity.

\section{Local structure theorems}
\label{local}

Those algebraic group actions can be effectively studied which are
more or less reduced to linear or projective representations and
their restrictions to stable subvarieties in representation
spaces. Therefore it is natural to restrict our attention to the
following class of actions.
\begin{definition}
A regular algebraic group action $G:X$ (or a $G$-variety~$X$) is
\emph{good} if $X$ can be covered by $G$-stable quasiprojective
open subsets $X_i$ such that $G:X_i$ is the restriction of the
projective representation of $G$ in an ambient projective space.
\end{definition}
\begin{example}\label{8}
Consider a rational projective curve $X$ obtained from $\PP^1$ by
identifying $0,\infty\in\PP^1$ in an ordinary double point (a
Cartesian leaf). A $\kk^{\times}$-action on $\PP^1$ with the fixed
points $0,\infty$ goes down to~$X$. This action is not good.
(Otherwise, there is a $\kk^{\times}$-stable hyperplane section of
$X$ in an ambient $\PP^n$ that does not contain the double point.
But there are no other $\kk^{\times}$-fixed points on~$X$.)
\end{example}
The reason for the action of Example~\ref{8} fails to be good is
non-normality of~$X$.
\begin{example}
If $X$ is a $G$-stable subvariety of a normal $G$-variety (e.g.,
$X$~is itself normal), then $G:X$ is good by Sumihiro's theorem
(\ref{Sumihiro}).
\end{example}
The normalization or the equivariant Chow lemma~\cite[Th.1.3]{IT}
reduce the study of arbitrary algebraic group actions to good
ones. In the sequel, only good actions are considered unless
otherwise specified.

Now let $G$ be a connected reductive group. Fix a Borel subgroup
$B\subseteq G$ and a maximal torus $T\subseteq B$. Put $U=B'$, a
maximal unipotent subgroup of~$G$.

In order to describe the local structure of good $G$-actions, we
begin with a helpful technical construction in characteristic zero
due to Brion--Luna--Vust and Grosshans.

Let $V$ be a finite-dimensional $G$-module with a lowest weight
vector~$v$, and $\omega\in V^{*}$ be the dual highest weight
vector such that $\langle v,\omega\rangle\neq0$. Let
$P=G_{\langle\omega\rangle}=L\semitimes\Ru{P}$ be the projective
stabilizer of~$\omega$, where $L\supseteq T$ is a Levi subgroup
and $\Ru{P}=\Radu{P}$. Then $P^{-}=G_{\langle
v\rangle}=L\semitimes\Ru{P^{-}}$ is the opposite parabolic to~$P$,
where $\Ru{P^{-}}=\Radu{P^{-}}$.

Put $\oo{V}=V\setminus\langle\omega\rangle^{\ann}$,
$W=\langle\Ru{\p^{-}}\omega\rangle^{\ann}$, and
$\oo{W}=W\cap\oo{V}$. (Here ${}^{\ann}$ denotes the annihilator in
the dual subspace.)
\begin{lemma}[\cite{BLV}, \cite{Gro}]\label{loc.str.V}
In characteristic zero, the action $P:\oo{V}$ gives rise to an
isomorphism
%*
\begin{equation*}
\Ru{P}\times\oo{W}=P\itimes{L}\oo{W}\isoto\oo{V}
\end{equation*}
%*
\end{lemma}
\begin{proof}
Consider level hyperplanes $V_c=\{x\in V\mid\langle
x,\omega\rangle=c\}$, $W_c=W\cap V_c$. We have
%*
\begin{equation*}
\oo{V}=\bigsqcup_{c\neq0}V_c=\kk^{\times}v+V_0,
\end{equation*}
%*
similarly for~$\oo{W}$. Note that
$W_0=\langle\g\omega\rangle^{\ann}$ is $P$-stable. Affine
hyperplanes $V_c\subset V$ and $V_c/W_0\subset V/W_0$ are
$P$-stable. It suffices to show that the induced action
$\Ru{P}:V_c/W_0$ is transitive and free whenever $c\neq0$. But
$V_0=\Ru{\p}v\oplus W_0$, whence $cv+W_0$ has the dense
$\Ru{P}$-orbit in $V_c/W_0$ with the trivial stabilizer. It
remains to note that all orbits of a unipotent group on an affine
variety are closed (Lemma~\ref{uni=>closed}).
\end{proof}
\begin{theorem}[{\cite[2.3]{W_X},
\cite[1.2]{G-val}, \cite{BLV}}]\label{loc.str} Let $X$ be a good
$G$-variety and $Y\subseteq X$ a $G$-stable subvariety. Then there
exists a unique parabolic subgroup $P=P(Y)\supseteq B$ with a Levi
decomposition $P=L\semitimes\Ru{P}$, $L\supseteq T$,
$\Ru{P}=\Radu{P}$, and a $T$-stable locally closed affine
subvariety $Z\subseteq X$ such that:
\begin{enumerate}
\item $\X=PZ$ is an affine open subset of~$X$. \item The action
$\Ru{P}:\X$ is proper and has a geometric quotient
$\X/\Ru{P}=\Spec\kk[\X]^{\Ru{P}}$. \item\label{fin.sur} A natural
map $\Ru{P}\times Z\to\X$, $(g,z)\mapsto gz$, and the quotient map
$Z\to\X/\Ru{P}$ are finite and surjective.
\begin{question}
purely inseparable?
\end{question}
\item[\ref{fin.sur}$'$] In characteristic zero, $Z$ may be chosen
to be $L$-stable and such that the action $P:\X$ gives rise to an
isomorphism
%*
\begin{equation*}
\Ru{P}\times Z=P\itimes{L}Z\isoto\X
\end{equation*}
%*
\item $\Y=Y\cap\X\neq0$, and the kernel $L_0$ of the natural
action $L=P/\Ru{P}:\Y/\Ru{P}$ contains~$L'$. Moreover,
$\Y/\Ru{P}\iso L/L_0\times C$, where the torus $L/L_0$ acts on $C$
trivially. In characteristic zero, $Y\cap Z\iso L/L_0\times C$.
\end{enumerate}
\end{theorem}
\begin{proof}
We will assume $\ch\kk=0$. (For the general case, see~\cite[1.2,
\S2]{G-val}.) Replacing $X$ by an open $G$-subvariety, we may
assume that $X$ is quasiprojective, $Y$ is closed in~$X$, and
there is a very ample $G$-line bundle $\Ll$ on~$X$. Then $X$ is
$G$-equivariantly embedded in~$\PP(V)$, $V=\Ho^0(X,\Ll)^{*}$.

Let $\overline{X},\overline{Y}$ be the closures of $X,Y$
in~$\PP(V)$. We can find a homogeneous $B$-eigenform $\omega$ in
coordinates on~$V$ that vanishes on $\overline{X}\setminus X$ and
on any given closed $B$-subvariety $D\subset Y$, but not on~$Y$.
(Take a nonzero $B$-eigenform in the ideal of
$D\cup(\overline{X}\setminus X)$ in the homogeneous coordinate
ring of $(\overline{X}\setminus X)\cup\overline{Y}$, and extend it
to $\overline{X}$ by complete reducibility of $G$-modules.)
Replacing $\Ll$ by its power, we may assume $\omega\in
\Ho^0(X,\Ll)^{(B)}$.

Now $\X=X_{\omega}$ is an affine open subset of~$X$. By
Lemma~\ref{loc.str.V},
%*
\begin{equation*}
\X\iso P\itimes{L}Z=\Ru{P}\times Z,
\end{equation*}
%*
where $P=G_{\langle\omega\rangle}$ and $Z=\PP(\oo{W})\cap X$.

If we choose for $D$ a (maybe reducible) $B$-stable divisor in~$Y$
whose stabilizer $P$ is the smallest possible one, then any
$(B\cap L)$-stable divisor of $Y\cap Z$ is $L$-stable. It follows
that each highest weight function in $\kk[Z]$ is
$L$-semiinvariant, whence $L'$ acts on $Z$ trivially. Taking $D$
sufficiently large, we may replace $Z$ by an open subset
$L$-isomorphic to $L/L_0\times C$ (with the trivial action
on~$C$).

To complete the proof, note that $P$ is uniquely determined by the
conditions of the theorem. Namely, $P=P(Y)$ equals the smallest
stabilizer of a $B$-stable divisor in~$Y$.
\end{proof}
\begin{corollary}\label{loc.str.gen}
Let $P=P(X)$ be the smallest stabilizer of a $B$-stable divisor
in~$X$. Then there exists a $T$-stable ($L$-stable if $\ch\kk=0$)
locally closed affine subset $Z\subseteq X$ such that $\X=PZ$ is
an open affine subset of~$X$, the natural maps $\Ru{P}\times
Z\to\X$, $Z\to\X/\Ru{P}$ are finite and surjective (isomorphic if
$\ch\kk=0$), and $\X/\Ru{P}\iso L/L_0\times C$, where $L\supseteq
L_0\supseteq L'$ and the $L$-action on $C$ is trivial.
\end{corollary}
\begin{example}
Let $X=\PP(S^2{\kk^n}^{*})$ be the space of quadrics in
$\PP^{n-1}$, $\ch\kk\neq2$. Then $G=\GL_n(\kk)$ acts on $X$ by
linear variable changes with the orbits $\HS_1,\dots,\HS_n$, where
$\HS_r$ is the set of quadrics of rank~$r$, and
$\overline{\HS_1}\subset\dots\subset\overline{\HS_n}=X$. Choose
the standard Borel subgroup $B\subseteq G$ of upper-triangular
matrices and the standard maximal torus $T\subseteq B$ of diagonal
matrices.
\begin{roster}
\item Put $Y=\overline{\HS_1}$, the unique closed $G$-orbit
in~$X$, which consists of double hyperplanes. In the notation of
Lemma~\ref{loc.str.V} and Theorem~\ref{loc.str}, we have
$V=S^2{\kk^n}^{*}\ni v=x_1^2$, $V^{*}=S^2\kk^n\ni\omega=e_1^2$.
(Here $e_1,\dots,e_n$ form the standard basis of~$\kk^n$ and
$x_1,\dots,x_n$ are the standard coordinates.) Then $P$ is a
standard parabolic subgroup of matrices of the form
%*
\begin{equation*}
\begin{array}{|c|@{\qquad}c@{\qquad}|@{}l}
\cline{1-2}
       L       &          \Ru{P}        & \\
\cline{1-2}
\begin{array}{@{}c@{}}
    0    \\
 \Dots0+ \\
    0    \\
\end{array} &
\begin{array}{@{}c@{}}
                                          \\
                 \smash{\text{\LARGE$L$}} \\
                                          \\
\end{array} &
\left.\vphantom{\begin{array}{@{}c@{}} \\ \\ \\ \end{array}}
\right\}\scriptstyle{n-1} \\
\cline{1-2}
\multicolumn{3}{@{}r@{}}{}\\[-4.3ex]
\multicolumn{2}{@{}r@{}}
{\underbrace{\hphantom{\qquad\smash{\text{\LARGE$L$}}\qquad}}_{n-1}}&
\end{array}
\end{equation*}
%*
(We indicate the Levi decomposition of $P$ at the figure.)
$\oo{V}$ is the set of quadratic forms $q=cx_1^2+\dots$, $c\neq0$,
$\Ru{\p}v=\{a_{12}x_1x_2+\dots+a_{1n}x_1x_n\mid a_{ij}\in\kk\}$,
and $W$ is the space of forms $q=cx_1^2+q'(x_2,\dots,x_n)$, where
$c\in\kk$, $q'$~is a quadratic form in $x_2,\dots,x_n$.

Now $\X$ is the set of quadrics given by an equation
$x_1^2+\dots=0$, $Z$~consists of quadrics with an equation
$x_1^2+q'(x_2,\dots,x_n)=0$, and $Y\cap Z=\{\langle
x_1^2\rangle\}$. Lemma~\ref{loc.str.V} or Theorem~\ref{loc.str}
say that every quadratic form with nonzero coefficient at~$x_1^2$
can be moved by~$\Ru{P}$, i.e., by a linear change of~$x_1$, to a
form containing no products $x_1x_j$, $j>1$. This is the first
step in the Lagrange method of transforming a quadric to the
normal form.

\item More generally, put $Y=\HS_r$. It is easy to see that $P(Y)$
is the group of matrices of the form
%*
\begin{equation*}
\begin{array}{r@{}|c|@{\qquad}c@{\qquad}|}
\multicolumn{2}{@{}r@{\qquad}}
{\overbrace{\hphantom{\begin{array}{c@{}c@{}c}
\smash{*} & & \smash{\text{\Large$*$}} \\[-1ex]
          & \Dots+-        &           \\[-1.5ex]
\smash{\text{\Large$0$}} & & \smash{*} \\
\end{array}}}^{r}} &
\multicolumn{1}{@{}c@{}}{} \\[-0.5ex]
\cline{2-3} \scriptstyle{r}\left\{
\begin{array}{@{}c@{}}
\\[-1ex] \\[-1.5ex] \\
\end{array}
\right. &
\begin{array}{@{}c@{}c@{}c@{}}
\smash{*} & & \smash{\text{\Large$*$}} \\[-1ex]
          & \Dots+-        &           \\[-1.5ex]
\smash{\text{\Large$0$}} & & \smash{*} \\
\end{array} & \smash{\text{\huge$*$}} \\
\cline{2-3}
&                   &                   \\
& {\text{\huge$0$}} & {\text{\huge$*$}} \\
&                   &                   \\
\cline{2-3}
\end{array}
\end{equation*}
%*
Clearly, $P(Y)=G_{\langle\omega\rangle}$, where $\omega$ is the
product of the first $r$ upper-left corner minors of the matrix of
a quadratic form. Then $\X$ is the set of quadrics, where $\omega$
does not vanish, i.e., having non-degenerate intersection with all
subspaces $\{x_k=\dots=x_n=0\}$, $k\leq r+1$. Further, $Z$
consists of quadrics with an equation
$c_1x_1^2+\dots+c_rx_r^2+q'(x_{r+1},\dots,x_n)=0$, $c_i\neq0$, and
$Y\cap Z=\{\langle c_1x_1^2+\dots+c_rx_r^2\rangle\mid c_i\neq0\}$.
The Levi subgroup $L=(\kk^{\times})^r\times\GL_{n-r}(\kk)$ acts on
$Y\cap Z$ via the first factor, and $Y\cap
Z=(\kk^{\times})^r\times\{\langle x_1^2+\dots+x_r^2\rangle\}$.
Theorem~\ref{loc.str} says that each quadric with nonzero first
$r$ upper-left corner minors transforms by a unitriangular linear
variable change to the form
$c_1x_1^2+\dots+c_rx_r^2+q'(x_{r+1},\dots,x_n)$---this is nothing
else, but the Gram--Schmidt orthogonalization method.
\end{roster}
\end{example}

A refined version of the local structure theorem was proved by
Knop in characteristic zero.

Let $X$ be a $G$-variety. We call any formal $\kk$-linear
combination of prime Cartier divisors on~$X$ a
\emph{$\kk$-divisor}. Let $D=a_1D_1+\dots+a_sD_s$ be a $B$-stable
$\kk$-divisor, and $P=P[D]$ be the stabilizer of its support
$D_1\cup\dots\cup D_s$. Replacing $G$ by a finite cover, we may
assume that the line bundles $\Lin{D_i}$ are $G$-linearized. Let
$s_i\in\Ho^0(X,\Lin{D_i})^{(B)}$ be the sections of $B$-weights
$\lambda_i$ such that $\divr s_i=D_i$, and set $\lambda_D=\sum
a_i\lambda_i$. We say that $D$ is \emph{regular} if
$\langle\lambda_D,\alpha^{\vee}\rangle\neq0$ for any root $\alpha$
such that $\g_{\alpha}\subseteq\Ru{\p[D]}$. (For example, any
effective $B$-stable Cartier divisor is regular.)

Define a morphism $\psi_D:X\setminus D\to\g^{*}$ by the formula
%*
\begin{equation*}
\langle\psi_D(x),\xi\rangle=\sum a_i\frac{\xi s_i}{s_i}(x),
\qquad\forall\xi\in\g
\end{equation*}
%*
\begin{theorem}[{\cite[2.3]{inv.mot}}]\label{loc.str.mom}
The map $\psi_D$ is a $P$-equivariant fibration over the $P$-orbit
of $\lambda_D$ considered as a linear function on a maximal torus
$\tr\subseteq\br$ and extended to $\g$ by putting
$\langle\lambda_D,\g_{\alpha}\rangle=0,\ \forall\alpha$. The
stabilizer $L=P_{\lambda_D}$ is the Levi subgroup of $P$
containing~$T$. In particular, $X\setminus D\iso P\itimes{L}Z$,
where $Z=\psi^{-1}(\lambda_D)$.
\end{theorem}

Other versions of the local structure theorem can be found
in~\cite{struct.loc} and in~\ref{B-charts},~\ref{c=0}.
\begin{question}
Rewrite~\ref{local}. In particular, relate
Theorems~\ref{loc.str},~\ref{loc.str.mom} to the moment map.
\end{question}

\section{Complexity and rank of $G$-varieties}
\label{c&r(basics)}

As before, $G$ is a connected reductive group with a fixed Borel
subgroup~$B$, a maximal unipotent subgroup $U=B'$, and a maximal
torus $T\subseteq B$. Let $X$ be an irreducible $G$-variety.
\begin{definition}
The \emph{complexity} $c_G(X)$ of the action $G:X$ is the
codimension of a generic $B$-orbit in~$X$. By the lower
semicontinuity of the function $x\mapsto\dim Bx$,
$c_G(X)=\min_{x\in X}\codim Bx$. By the Rosenlicht
theorem~\cite[2.3]{IT}, $c_G(X)=\trdeg\kk(X)^B$.

The \emph{weight lattice} $\RG(X)$ (resp.~the \emph{weight
semigroup} $\RG_{+}(X)$) is the set of weights of all rational
(regular) $B$-eigenfunctions on~$X$. It is a sublattice in the
weight lattice $\Ch(B)=\Ch(T)$ (a submonoid in the monoid
$\Ch_{+}$ of dominant weights, respectively).

The integer $r_G(X)=\rk\RG(X)$ is the \emph{rank} of $G:X$.

We usually drop the subscript $G$ in the notation of complexity
and rank.
\end{definition}
Complexity, rank, and the weight lattice are birational invariants
of an action. Replacing $X$ by a $G$-birationally equivalent
variety, we may always assume that $X$ is good, normal,
quasiprojective, or smooth, when required. These invariants are
very important in studying the geometry of the action $G:X$ and
the related representation and compactification theory. Here we
examine the most basic properties of complexity and rank.
\begin{example}\label{cr(G/P)}
Let $X$ be a projective homogeneous $G$-space. By the Bruhat
decomposition, $U$ has a dense orbit in~$X$ (a \emph{big cell}).
Hence $c(X)=r(X)=0$.
\end{example}
\begin{example}\label{cr_T}
Assume $G=T$ and let $T_0$ be the kernel of the action $T:X$. Then
$X$ contains an open $T$-stable subset $\X=T/T_0\times Z$. Hence
$c(X)=\dim Z=\dim\kk(X)^T$, $\RG(X)=\Ch(T/T_0)$, $r(X)=\dim
T/T_0$.
\end{example}
\begin{example}
Let $G$ act on $X=G$ by left translations. Then $c(X)=\dim
G/B=\dim U$ is the number of positive roots of~$G$. By
Formula~\eqref{k[G/U]}, $\RG_{+}(X)=\Ch_{+}$.
\end{example}

It is easy to prove the following
\begin{proposition}\label{c+r}
$c(X)+r(X)=\min_{x\in X}\codim Ux=\trdeg\kk(X)^U$
\end{proposition}
(Just apply Example~\ref{cr_T} to the $T$-action on the rational
quotient $X/U$.)

Complexity and weight lattice (semigroup) are monotonous by
inclusion. More precisely, we have
\begin{theorem}[{\cite[2.3]{B-orb}}]
\label{cr(Y<X)} For any closed $G$-subvariety $Y\subseteq X$,
$c(Y)\leq c(X)$, $r(Y)\leq r(X)$, and $Y=X$ iff the equalities
hold. Furthermore, $\RG(Y)\subseteq\frac1q\RG(X)$ and if $X$ is
affine, then $\RG_{+}(Y)\subseteq\frac1q\RG_{+}(X)$, where $q$ is
a sufficiently big power of the characteristic exponent of~$\kk$
($=\ch\kk$, or $1$ if $\ch\kk=0$).
\end{theorem}
The proof relies on a helpful lemma of Knop:
\begin{lemma}\label{B-eigenfun}
Let $Y\subseteq X$ be a $G$-subvariety and $p$ be the
characteristic exponent of~$\kk$. Then
%*
\begin{equation*}
\forall f\in\kk(Y)^{(B)}\ \exists\widetilde{f}\in\Oo_{X,Y}^{(B)}\
\exists q=p^N:\ f^q=\widetilde{f}|_Y
\end{equation*}
%*
\end{lemma}
\begin{proof}
Applying normalization, we may assume that $X$ is good and even
projective and $Y$ is closed in~$X$. Embed $X$ $G$-equivariantly
in a projective space and let $\widehat{X},\widehat{Y}$ be the
cones over~$X,Y$. These cones are
$\widehat{G}=G\times\kk^{\times}$-stable ($\kk^{\times}$ acts by
homotheties), and $f$ is pulled back to
$\kk(\widehat{Y})^{(\widehat{B})}$, where
$\widehat{B}=B\times\kk^{\times}$. Thence $f=F_1/F_2$, where
$F_i\in\kk[\widehat{Y}]^{(\widehat{B})}$ are homogeneous
$B$-semiinvariant polynomials. By Corollary~\ref{(M/N)^H},
$F_i^q=\widetilde{F}_i|_{\widehat{Y}}$ for some
$\widetilde{F}_i\in\kk[\widehat{X}]^{(\widehat{B})}$, $q=p^N$. Now
$\widetilde{f}=\widetilde{F}_1/\widetilde{F}_2$ is pulled down to
a rational $B$-eigenfunction on $X$ such that
$\widetilde{f}|_Y=f^q$.
\end{proof}
\begin{proof}[Proof of the theorem]
Lemma~\ref{B-eigenfun} implies that $q\RG(Y)\subseteq\RG(X)$ and
that $\kk(Y)^B$ is a purely inseparable extension of the residue
field of $\Oo_{X,Y}^{B}$, whence the inequalities and the
inclusion of weight lattices. The inclusion of weight semigroups
stems from Corollary~\ref{(M/N)^H}.

Now suppose $c(Y)=c(X)$ and $r(Y)=r(X)$. As in
Lemma~\ref{B-eigenfun}, we may assume that $X,Y$ are closed
$G$-subvarieties in a projective space and consider the cones
$\widehat{X},\widehat{Y}$ over~$X,Y$. We have
$c_{\widehat{G}}(\widehat{X})=c_G(X)$ and
$r_{\widehat{G}}(\widehat{X})=r_G(X)+1$, in view of an exact
sequence
%*
\begin{equation*}
0\longrightarrow\RG(X)\longrightarrow
\widehat{\RG}(\widehat{X})\longrightarrow\Ch(\kk^{\times})=
\ZZ\longrightarrow0,
\end{equation*}
%*
where $\widehat{\RG}$ is the weight lattice relative
to~$\widehat{B}$. Similar equalities hold for $\widehat{Y}$.

By assumption and Proposition~\ref{c+r},
$\trdeg\kk(\widehat{Y})^U=\trdeg\kk(\widehat{X})^U$. But a
rational $U$-invariant function on an affine variety is the
quotient of two $U$-invariant polynomials~\cite[Th.3.3]{IT},
whence $\kk[\widehat{Y}]^U$ and $\kk[\widehat{X}]^U$ have the same
transcendence degree. By Lemma~\ref{(A/I)^H}, $\kk[\widehat{Y}]^U$
is a purely inseparable finite extension of
$\kk[\widehat{X}]^U|_{\widehat{Y}}$, whence $\kk[\widehat{X}]^U$
restricts to $\widehat{Y}$ injectively. Therefore the ideal of
$\widehat{Y}$ contains no nonzero $U$-invariants, hence is zero.
It follows that $\widehat{Y}=\widehat{X}$, whence $Y=X$.
\end{proof}

On the other side, there is a general procedure of ``enlarging'' a
variety which preserves complexity, rank, and the weight lattice.
\begin{definition}
Let $G,G_0$ be connected reductive groups. We say that a
$G$-variety $X$ is obtained from a $G_0$-variety $X_0$ by
\emph{parabolic induction} if $X=G\itimes{Q}X_0$, where
$Q\subseteq G$ is a parabolic subgroup acting on $X_0$ via an
epimorphism $Q\onto G_0$.
\end{definition}
\begin{proposition}
$c_G(X)=c_{G_0}(X_0)$, $r_G(X)=r_{G_0}(X_0)$, $\RG(X)=\RG(X_0)$.
\end{proposition}
The proof is easy.

The weight lattice is actually an attribute of a generic $G$- (and
even $B$-) orbit.
\begin{proposition}[{\cite[2.6]{B-orb}}]\label{r(X)=r(G*)}
$\RG(X)=\RG(Gx)$ for all $x$ in an open subset of~$X$.
\end{proposition}
\begin{proof}
Replacing $X$ by the rational quotient $X/U$, we reduce the
problem to the case $G=B=T$. Let $f_1,\dots,f_r\in\kk(X)$ be
rational $T$-eigenfunctions whose weights generate~$\RG(X)$. If
all $f_i$ are defined and nonzero at~$x$, then by
Lemma~\ref{B-eigenfun}, $\RG(X)=\RG(Tx)$. (In this case, we may
assume $q=1$ in Lemma~\ref{B-eigenfun}, since $T$ is linearly
reductive.)
\end{proof}
\begin{corollary}
The function $x\to r(Gx)$ is lower semicontinuous on~$X$.
\end{corollary}

Using Lemma~\ref{B-eigenfun}, Arzhantsev proved the following
\begin{proposition}[{\cite[\S2]{SL(2):1-p}}]\label{c(Gx)}
The function $x\to c(Gx)$ is lower semicontinuous on~$X$.
\end{proposition}

In the affine case, the weight semigroup is a more subtle
invariant of an action than the weight lattice.
\begin{proposition}
For quasiaffine~$X$, $\RG(X)=\ZZ\RG_{+}(X)$.
\end{proposition}
\begin{proof}
Any rational $B$-eigenfunction on~$X$ is a quotient of two
polynomials: $f=f_1/f_2$. By the Lie-Kolchin theorem, there exists
a nonzero $B$-semiinvariant linear combination
$\sum\lambda_i(b_if_2)$, $\lambda_i\in\kk$, $b_i\in B$. Then
$f=\widetilde{f}_1/\widetilde{f}_2$, where
$\widetilde{f}_j=\sum\lambda_i(b_if_j)$ are polynomial
$B$-eigenfunctions on~$X$.
\end{proof}
\begin{proposition}
For affine~$X$, the semigroup $\RG_{+}(X)$ is finitely generated.
\end{proposition}
\begin{proof}
The semigroup $\RG_{+}(X)$ is the semigroup of weights for the
$T$-weight decomposition of~$\kk[X]^U$, the latter algebra being
finitely generated by Theorem~\ref{U-inv}\ref{A<=>A^U}.
\end{proof}

In characteristic zero, the complexity controls the growth of
multiplicities in the spaces of global sections of $G$-line
bundles on~$X$.
\begin{theorem}\label{mult}
\begin{roster}
\item\label{mult(aff)} If $X$ is affine and $\kk[X]^G=\kk$ (e.g.,
$X$~contains an open $G$-orbit), then $c(X)$ is the minimal
integer $c$ such that $m_{n\lambda}(X)=O(n^c)$ for every dominant
weight~$\lambda$. \item\label{mult(proj)} If $X$ is projective,
then $c(X)$ is the minimal integer $c$ such that
$m_{n\lambda}(\Ll^n)=O(n^c)$ for any line bundle $\Ll$ on $X$ and
any dominant weight~$\lambda$.
\begin{question}
Is this generalized to quasiprojective or homogeneous case?
Yes---to $X=G/H$.
\end{question}
\end{roster}
\end{theorem}
\begin{proof}
\begin{roster}
\item[\ref{mult(aff)}] By Proposition~\ref{mult(h.w)},
$m_{\lambda}(X)=\dim\kk[X]^U_{\lambda}$. Replacing $X$ by $X\by
U$, we may assume that $G=B=T$. Put
%*
\begin{equation*}
A(\lambda)=\bigoplus_{n\geq0}\kk[X]_{n\lambda}
\end{equation*}
%*
Then $A(\lambda)\iso(\kk[X]\otimes\kk[t])^T$, where the
indeterminate $t$ has the $T$-weight~$-\lambda$. The function
field $K=\Quot A(\lambda)$ is purely transcendental of degree~$1$
over $K^B\subseteq\kk(X)^B$, whence $\trdeg A(\lambda)\leq c+1$.

If $\lambda$ lies in the interior of the cone $\QQ_{+}\RG_{+}(X)$,
then $\trdeg A(\lambda)=c+1$. Indeed, we have $\lambda=\sum
l_i\lambda_i$, where $\lambda_i$ are the generators
of~$\RG_{+}(X)$ and $l_i$ are rational positive numbers. Any $f\in
K^B$ is expressible as $f=h_1/h_2$, where $h_j\in\kk[X]_{\mu}$ for
some~$\mu$. Now $\mu=\sum m_i\lambda_i$, and for $n$ sufficiently
large, $n_i=nl_i-m_i$ are positive integers. Then
$f=\widetilde{h}_1/\widetilde{h}_2$, where
$\widetilde{h}_j=h_j\prod f_i^{n_i}\in\kk[X]_{n\lambda}$,
$f_i\in\kk[X]_{\lambda_i}$. Thus $K^B=\kk(X)^B$, q.e.d.

By the above, $A(\lambda)$ is a finitely generated graded algebra
of Krull dimension $d\leq c+1$, and the equality holds for
general~$\lambda$. We conclude by a standard result of dimension
theory, that $\dim A(\lambda)_n$ grows as a polynomial in~$n$ of
degree~$d$, at least for $n$ sufficiently divisible.

\item[\ref{mult(proj)}] We have
$m_{n\lambda}(\Ll^n)=\dim\Ho^0(X,\Ll^n)^U_{n\lambda}$. For a
sufficiently ample $G$-line bundle~$\Mm$, the line bundle
$\Ll\otimes\Mm$ is very ample, and for nonzero
$s\in\Ho^0(X,\Mm)^U_{\mu}$, we have an inclusion
%*
\begin{equation*}
\Ho^0(X,\Ll^n)^U_{n\lambda}\embeds
\Ho^0(X,(\Ll\otimes\Mm)^n)^U_{n(\lambda+\mu)}
\end{equation*}
%*
provided by $(\cdot\otimes s^n)$. Thus it suffices to consider
very ample~$\Ll$. Consider the respective projective embedding
of~$X$, and let $\widehat{X}$ be the affine cone over~$X$. Then
$\widehat{G}=G\times\kk^{\times}$ acts on~$\widehat{X}$ (here
$\kk^{\times}$ acts by homotheties), and
$m_{n\lambda}(\Ll^n)=m_{(n\lambda,n)}(\widehat{X})$ (at least for
$n\gg0$). The assertion~\ref{mult(aff)} applied to~$\widehat{X}$
concludes the proof. \qedhere\end{roster}
\end{proof}
\begin{example}
Let $\HS=G/H$ be a quasiaffine homogeneous space such that
$\kk[\HS]$ is finitely generated. Then $X=\Spec\kk[\HS]$ contains
$\HS$ as an open orbit. By Theorem~\ref{mult}\ref{mult(aff)},
$m_{n\lambda}(X)=\dim V(n\lambda^{*})^H$ grows as $n^{c(\HS)}$ for
general~$\lambda$, and no faster for any~$\lambda$.

E.g., for $H=\{\1\}$, $c(G)$~is the number of positive roots, and
$\dim V(n\lambda^{*})$ is the polynomial in~$n$ of degree $\leq
c(G)$ by the Weyl dimension formula.
\end{example}
\begin{remark}
For $G$-varieties of complexity $\leq1$, more precise results on
multiplicities are obtained, see~\ref{c=1},~\ref{sphericity}.
\end{remark}

\section{Complexity and modality}
\label{modality}

The notion of modality was introduced in the works of Arnold on
the theory of singularities. The modality of an action is the
maximal number of parameters in a continuous family of orbits.
More precisely,
\begin{definition}
Let $H:X$ be an algebraic group action. The integer
%*
\begin{equation*}
d_H(X)=\min_{x\in X}\codim_XHx=\trdeg\kk(X)^H
\end{equation*}
%*
is called the \emph{generic modality} of the action. The
\emph{modality} of $H:X$ is the number $\modl_HX=\max_{Y\subseteq
X}d_H(Y)$, where $Y$ runs through $H$-stable irreducible
subvarieties of~$X$.
\end{definition}
Note that $c(X)=d_B(X)$.

It may happen that the modality is greater than the generic
modality of an action. For example, the natural action
$\GL_n(\kk):\LO_n(\kk)$ by left multiplication has an open orbit,
whereas its modality equals~$[n^2/4]$. Indeed, $\LO_n(\kk)$ is
covered by finitely many locally closed $\GL_n(\kk)$-stable
subsets $Y_{i_1,\dots,i_k}$, where $Y_{i_1,\dots,i_k}$ is the set
of matrices of rank~$k$ with linearly independent columns
$i_1,\dots,i_k$. Therefore an orbit in $Y_{i_1,\dots,i_k}$ depends
on $k(n-k)$ parameters, which are the coefficients of linear
expressions of the remaining $n-k$ columns by the columns
$i_1,\dots,i_k$. The maximal number of parameters is obtained for
$k=[\frac{n+1}2]$.

Replacing $\GL_n(\kk)$ by the group $B_n(\kk)$ of non-degenerate
upper-triangular matrices and $\LO_n(\kk)$ by the space
$\overline{B_n(\kk)}$ of all upper-triangular matrices shows that
the same thing may happen for a solvable group action. The action
$B_n(\kk):\overline{B_n(\kk)}$ has an open orbit, but infinitely
many orbits in its complementary.

Remarkably, for a $G$-variety $X$ and the restricted action $B:X$,
the modality equals the generic modality (=the complexity) of the
action. This result was obtained by Vinberg~\cite{c&mod} with the
aid of Popov's technique of contracting to a horospherical variety
(cf.~\ref{cotangent}). We present a proof due to
Knop~\cite{B-orb}, who developed some earlier ideas of Matsuki. A
basic tool is an action of a certain monoid on the set of
$B$-stable subvarieties.

Let $W=N_G(T)/T$ be the Weyl group of~$G$. By the Bruhat
decomposition, the only irreducible closed $B\times B$-stable
subvarieties in $G$ are the Schubert varieties
$D_w=\overline{BwB}$, $w\in W$.
\begin{definition}[{\cite[\S2]{B-orb}, \cite{RS}}]
The \emph{Richardson--Springer monoid} (\emph{RS-monoid}) of $G$
is the set of all Schubert subvarieties in $G$ with the
multiplication as of subsets in~$G$. Equivalently, RS-monoid is
the set $W$ with a new multiplication $*$ defined by
$D_{v*w}=D_vD_w$. We denote the set $W$ equipped with this product
by~$W^{*}$.
\end{definition}

Clearly, $W^{*}$ is an associative monoid with the unity~$\1$. It
is easy to describe $W^{*}$ by generators and relations. Namely,
$W$ is defined by generators $s_1,\dots,s_l$ (simple reflections)
and relations $s_i^2=\1$ and
%*
\begin{equation*}
\underbrace{s_is_js_i\dots}_{\text{$n_{ij}$ terms}}=
\underbrace{s_js_is_j\dots}_{\text{$n_{ij}$ terms}}
\qquad\text{(braid relations),}
\end{equation*}
%*
where $(n_{ij})$ is the Coxeter matrix of~$W$. The monoid $W^{*}$
has the same generators and relations $s_i^2=s_i$ and braid
relations. If $w=s_{i_1}\dots s_{i_n}$ is a reduced decomposition
of $w\in W$, then $w=s_{i_1}*\dots*s_{i_n}$ in~$W^{*}$. All these
assertions follow from standard facts on multiplication of
Schubert cells in~$G$~\cite[\S29]{Agr}.

Let $\B(X)$ be the set of all closed irreducible $B$-stable
subvarieties in~$X$. The RS-monoid acts on $\B(X)$ in a natural
way: $w*Z=D_wZ$ is $B$-stable and closed as the image of
$D_w\itimes{B}Z$ under the proper morphism $G\itimes{B}X\iso
G/B\times X\to X$.
\begin{proposition}\label{cr&RS}
$c(w*Z)\geq c(Z)$, $r(w*Z)\geq r(Z)$ for any $Z\in\B(X)$.
\end{proposition}
\begin{proof}
It suffices to consider the case of a simple reflection $w=s_i$.
In this case, $D_w=P_i$ is the respective minimal parabolic
subgroup of~$G$ If $Z$ is $P_i$-stable, there is nothing to prove.
Otherwise, the map $P_i\itimes{B}Z\to P_iZ$ is generically finite,
and we may replace $s_i*Z$ by $P_i\itimes{B}Z$ and further, by an
open subset $Bs_iB\itimes{B}Z=B\itimes{B_i}s_iZ$, where
$B_i=B\cap_iBs_i^{-1}$. Therefore the complexity (rank) of $s_iZ$
equals the complexity (resp.~rank) of $s_iZ$ w.r.t.\ the
$B_i$-action or of $Z$ w.r.t.\ the action of
$s_i^{-1}Bs_i\subseteq B$. The assertion follows.
\end{proof}
\begin{theorem}\label{c=mod}
For any $B$-stable irreducible subvariety $Y\subseteq X$, we have
$c(Y)\leq c(X)$, $r(Y)\leq r(X)$. In particular,
$\modl_H(X)=d_H(X)$, where $H=B$ or~$U$.
\end{theorem}
\begin{proof}
Follows from Proposition~\ref{cr&RS}, Theorem~\ref{cr(Y<X)}, and
Proposition~\ref{c+r}.
\end{proof}
\begin{corollary}[\cite{c&mod}, \cite{B-fin}]
\label{B-fin} Every spherical variety contains finitely many
$B$-orbits.
\end{corollary}

In the spherical case, elements of $\B(X)$ are just $B$-orbit
closures. The set of all $B$-orbits on a spherical variety,
identified with $\B(X)$, is an interesting combinatorial object.
It is finite and partially ordered by the adherence relation
$\adhereseq$ (=inclusion of orbit closures). This partial order is
compatible with the action of the RS-monoid and with the dimension
function in the following sense:
\begin{enumerate}
\item\label{monot} $\HS\adhereseq s_i*\HS$ \item
$\HS_1\adhereseq\HS_2\implies s_i*\HS_1\adhereseq s_i*\HS_2$ \item
$\HS_1\adheres\HS_2\implies\dim\HS_1<\dim\HS_2$ \item $\HS\adheres
s_i*\HS\implies\dim(s_i*\HS)=\dim\HS+1$ \item\label{1-step} (One
step property) $(s_i*\HS)_{\adhereseq}=W_i*\HS_{\adhereseq}$,
where $W_i=\{\1,s_i\}$ is a minimal standard Coxeter subgroup
in~$W$, and
$\HS_{\adhereseq}=\{\HS'\in\B(X)\mid\HS'\adhereseq\HS\}$ is the
closure of~$\HS$.
\end{enumerate}

This compatibility imposes strong restrictions on the adherence of
$B$-orbits on a spherical homogeneous space $X=G/H$.
By~\ref{1-step}, it suffices to know the closures of the minimal
orbits, i.e., such $\HS\in\B(G/H)$ that $\HS\neq w*\HS'$ for
$\forall\HS'\neq\HS,\ w\in W^{*}$. If all minimal orbits have the
same dimension then they are closed.
\begin{example}
For $H=B$, the $B$-orbits are the Schubert cells $B[w]\subset
G/B$, $w\in W$, and their closures are the Schubert subvarieties
$S_w=D_w/B$ in~$G/B$. By standard facts on the multiplication of
Bruhat cells, $B[\1]=\{[\1]\}$ is the unique minimal $B$-orbit.
Whence $S_w=\strut\overline{s_{i_1}*\dots*s_{i_n}*B[\1]}$
($w=s_{i_1}\dots s_{i_n}$ is a reduced decomposition)
$=W_{i_1}*\dots*W_{i_n}*B[\1]=P_{i_1}\dots P_{i_n}[\1]= (B\sqcup
Bs_{i_1}B)\dots(B\sqcup Bs_{i_n}B)[\1]= \bigsqcup Bs_{j_1}\dots
s_{j_k}B[\1]=\bigsqcup_{v=s_{j_1}\dots s_{j_k}} B[v]$ over all
subsequences $(j_1,\dots,j_k)$ of $(i_1,\dots,i_n)$. This is a
well-known description of the Bruhat order on~$W$.
\end{example}
\begin{example}\label{RS}
If $G/H$ is a symmetric space, i.e., $H$~is a fixed point set of
an involution, up to connected components, then $G/H$ is spherical
(Theorem~\ref{symm=>sph}) and all minimal $B$-orbits have the same
dimension~\cite{RS}. A complete description of $B$-orbits, of the
$W^{*}$-action, and of the adherence relation is obtained
in~\cite{RS} (cf.~Proposition~\ref{B&T-orb}).
\end{example}
\begin{example}
For $H=TU'$, the space $G/H$ is spherical, but the minimal
$B$-orbits have different dimensions. However, the adherence of
$B$-orbits is completely determined by the $W^{*}$-action with the
aid of properties \ref{monot}--\ref{1-step}. The set $\B(G/H)$,
the $W^{*}$-action, and the adherence relation are described
in~\cite{G/TU'}.
\end{example}
\begin{Conjecture}[\cite{G/TU'}]
For any spherical homogeneous space $G/H$, there is a unique
partial order on $\B(G/H)$ satisfying \ref{monot}--\ref{1-step}.
\end{Conjecture}

By Theorem~\ref{c=mod}, the complexity of a $G$-variety equals the
maximal number of parameters determining a continuous family of
$B$-orbits on~$X$. Generally, continuous families of $G$-orbits
depend on a lesser number of parameters. However, a result of
Akhiezer shows that the complexity of a $G$-action is the maximal
modality in the class of all actions birationally $G$-isomorphic
to the given one.
\begin{theorem}[\cite{c&mod_G}]\label{mod=c}
There exists a $G$-variety $X'$ birationally $G$-iso\-mor\-phic to
$X$ such that $\modl_GX'=c(X)$.
\end{theorem}
\begin{proof}
Let $f_1,\dots,f_c$ be a transcendence base of $\kk(X)^B/\kk$. We
may replace $X$ by a birationally $G$-isomorphic normal projective
variety. Consider an ample $G$-line bundle $\Ll$ on~$X$. Replacing
$\Ll$ by a power, we may find a section $s_0\in
\Ho^0(X,\Ll)^{(B)}$ such that $\divr{s_0}\geq\divr_{\infty}f_i$
for~$\forall i$. Put $s_i=f_is_0\in\Ho^0(X,\Ll)^{(B)}$.

Take a $G$-module $M$ generated by a highest weight vector~$m_0$
and such that there is a homomorphism $\psi_i:M\to\Ho^0(X,\Ll)$,
$\psi_i(m_0)=s_i$. Let $m_0,\dots,m_n$ be its basis of
$T$-eigenvectors with the weights $\lambda_0,\dots,\lambda_n$. Let
$E=\langle e_0,\dots,e_c\rangle$ be a trivial $G$-module of
dimension~$c+1$. A homomorphism $\psi:E\otimes M\to\Ho^0(X,\Ll)$,
$e_i\otimes m\mapsto\psi_i(m)$, gives rise to the rational
$G$-equivariant map $\phi:X\dasharrow\PP((E\otimes M)^{*})$. In
projective coordinates,
%*
\begin{equation*}
\psi(x)=[\,\cdots\,:\,\psi_i(m_j)(x)\,:\,\cdots\,]
\end{equation*}
%*

Take a one-parameter subgroup $\gamma\in\CoCh(T)$ such that
$\langle\alpha,\gamma\rangle>0$ for each positive root~$\alpha$.
If all $s_i(x)\neq0$, then
%*
\begin{xxalignat}{1}
\psi(\gamma(t)x)&=[\,\cdots\,:\,t^{-\langle\lambda_j,\gamma\rangle}
\psi_i(m_j)(x)\,:\,\cdots\,]\\
&=[\,\cdots\,:\,t^{\langle\lambda_0-\lambda_j,\gamma\rangle}
\psi_i(m_j)(x)\,:\,\cdots\,]
\longrightarrow[s_0(x):\cdots:s_c(x):0:\cdots:0]
\end{xxalignat}
%*
as $t\to0$, because $\lambda_0-\lambda_j$ is a positive linear
combination of positive roots for $\forall j>0$. Thus
%*
\begin{xxalignat}{1}
\lim_{t\to0}\gamma(t)\psi(x)&=\bigl([s_0(x):\cdots:s_c(x)],
[m_0^{*}]\bigr)\in\PP(E^{*})\times\PP(M^{*})\embeds \PP((E\otimes
M)^{*})
\end{xxalignat}
%*
(the Segre embedding), where $m_0^{*},\dots,m_n^{*}$ is the dual
basis of~$M^{*}$.

Let $X'\subseteq X\times\PP((E\otimes M)^{*})$ be the closure of
the graph of~$\phi$.  By the above,
$Y=X'\cap(X\times\PP(E^{*})\times\PP(M^{*}))$ contains points of
the form
%*
\begin{equation*}
x_0=\lim_{t\to0}\gamma(t)(x,\psi(x))=
\bigl(\lim_{t\to0}\gamma(t)x,[s_0(x):\cdots:s_c(x)],
[m_0^{*}]\bigr)
\end{equation*}
%*
The $G$-equivariant projection $Y\to\PP(E^{*})$ maps $x_0$ to
$[s_0(x):\cdots:s_c(x)]$, hence is dominant, because $f_i=s_i/s_0$
are algebraically independent on~$X$. Thence the generic modality
of any component of $Y$ dominating $\PP(E^{*})$ is greater or
equal to~$c$.
\end{proof}
\begin{corollary}[\cite{G-fin}]\label{G-fin}
A homogeneous space $\HS$ is spherical iff any $G$-variety $X$
with an open orbit isomorphic to $\HS$ has finitely many
$G$-orbits.
\end{corollary}

\section{Horospherical varieties}
\label{horosph}

There is a nice class of $G$-varieties, which is easily accessible
for study from the viewpoint of the local structure, complexity,
and rank.

\begin{definition}
A subgroup $S\subseteq G$ is \emph{horospherical} if $S$ contains
a maximal unipotent subgroup of~$G$. A $G$-variety $X$ is called
\emph{horospherical} if the stabilizer of any point on~$X$ is
horospherical. In other words, $X=GX^U$.
\end{definition}
\begin{remark}
In the definition, it suffices to require that the stabilizer of a
general point is horospherical. Indeed, this implies that $GX^U$
is dense in~$X$. On the other hand, $X^U$ is $B$-stable, whence
$GX^U$ is closed by Proposition~\ref{collapse}.
\end{remark}
\begin{example}
Consider a Lobachevsky space $L^n$ in the hyperbolic realization,
i.e., $L^n\subseteq\RR^{n+1}$ is an upper pole of a hyperboloid
$\{x\in\RR^{n+1}\mid(x,x)=1\}$ in an $(n+1)$-dimensional
pseudo-Euclidean space of signature~$(1,n)$. The group $(\Isom
L^n)^0=\SO_{1,n}^{+}$ acts transitively on the set of horospheres
in~$L^n$. Fix a horosphere $H^{n-1}\subset L^n$ and let $\langle
e_1\rangle\in\partial L^n\subseteq\RR\PP^n$ be its center lying on
the absolute. The vector $e_1\in\RR^n$ is isotropic and its
projective stabilizer $P$ is a parabolic subgroup
of~$\SO_{1,n}^{+}$. Take a line $\ell\subset L^n$ orthogonal
to~$H^{n-1}$. It intersects $\partial L^n$ in two points $\langle
e_1\rangle,\langle e_2\rangle$. The group $P$ contains a
one-parameter subgroup $A$ acting in $\langle e_1,e_2\rangle$ by
hyperbolic rotations and trivially on $\langle
e_1,e_2\rangle^{\perp}$. Then $A$ acts on $\ell$ by translations
and the complementary subgroup $S=P'$ is the stabilizer
of~$H^{n-1}$. In the matrix form,
%*
\begin{equation*}
S=\left\{\;
\begin{array}{|c|c|c|}
\hline
       1       & 0 & \text{\large\strut} u^{\tran}A \\
\hline
       0       & 1 &          0\;\Dots+0\;0         \\
\hline
       0       &   &                                \\
    \Dots0+    & u &    \smash{\text{\LARGE$A$}}    \\
       0       &   &                                \\
\hline
\end{array}
\right.\; \left|\; \vphantom{\begin{array}{|c|c|c|} \hline
       1       & 0 & \text{\large\strut} u^{\tran}A \\
\hline
       0       & 1 &         0\;\Dots+0\;0          \\
\hline
       0       &   &                                \\
    \Dots0+    & u &    \smash{\text{\LARGE$A$}}    \\
       0       &   &                                \\
\hline
\end{array}}
A\in\SO_{n-1},\ u\in\RR^{n-1} \;\right\}
\end{equation*}
%*
Recall that $H^{n-1}$ carries a Euclidean geometry, and $S=(\Isom
H^{n-1})^0$, where $\Radu{S}$ acts by translations and a Levi
subgroup of $S$ by rotations fixing an origin. Clearly, $S(\CC)$
is a horospherical subgroup of $\SO_{1,n}(\CC)$, which explains
the terminology.
\end{example}

For any parabolic $P\subseteq G$, let $P=L\semitimes\Ru{P}$ be its
Levi decomposition, and $L_0$ be any intermediate subgroup between
$L$ and~$L'$. Then a subgroup $S=L_0\semitimes\Ru{P}$ is
horospherical. Conversely,
\begin{lemma}\label{hor.sub}
Let $S\subseteq G$ be a horospherical subgroup. Then $P=N_G(S)$ is
parabolic, and for a Levi decomposition $P=L\semitimes\Ru{P}$,
$S=L_0\semitimes\Ru{P}$, where $L'\subseteq L_0\subseteq L$.
\end{lemma}
\begin{proof}
Embed $S$ regularly in a parabolic $P\subseteq G$. Since $S$ is
horospherical, $\Ru{S}=\Ru{P}$, and $S/\Ru{P}$ contains a maximal
unipotent subgroup of $P/\Ru{P}\iso L$, whence $S/\Ru{P}\iso
L_0\supseteq L'$. Now it is clear that $S=L_0\semitimes\Ru{P}$ and
$P=N_G(S)$, because $P$ normalizes~$S$ and $N_G(S)$
normalizes~$\Ru{S}$.
\end{proof}

In the sequel, assume $\ch\kk=0$ for simplicity.

Horospherical varieties can be characterized in terms of the
properties of multiplication in the algebra of regular functions.
For any $G$-module $M$ and any $\lambda\in\Ch_{+}$, let
$M_{(\lambda)}$ denote the isotypic component of type~$\lambda$
in~$M$.
\begin{proposition}[{\cite[\S4]{contr}}]\label{iso.mult}
A quasiaffine $G$-variety $X$ is horospherical iff
$\kk[X]_{(\lambda)}\cdot\kk[X]_{(\mu)}\subseteq\kk[X]_{(\lambda+\mu)}$
for $\forall\lambda,\mu\in\RG_{+}(X)$.
\end{proposition}

The local structure of a horospherical action is simple.
\begin{proposition}\label{loc.str.hor}
A horospherical $G$-variety $X$ contains an open $G$-stable subset
$\X\iso G/S\times C$, where $S\subseteq G$ is horospherical and
$G:C$ is trivial.
\end{proposition}
\begin{proof}
By a theorem of Richardson \cite[Th.7.1]{IT}, Levi subgroups of
stabilizers of general points on~$X$ are conjugate and unipotent
radicals of stabilizers form a continuous family of subgroups
in~$G$. Now it is clear from Lemma~\ref{hor.sub}, that a
horospherical subgroup may not be deformed outside its conjugacy
class, whence stabilizers of general points are all conjugate to a
certain $S\subseteq G$. Replacing $X$ by an open $G$-stable subset
yields $X\iso G\itimes{P}X^S$, where $P=N_G(S)$ \cite[2.8]{IT}.
But $P$ acts on $X^S$ via a torus $P/S$, hence $X^S$ is locally
$P$-isomorphic to $P/S\times C$, where $P$ acts on $C$ trivially
\cite[2.6]{IT}.
\end{proof}

To any $G$-variety $X$, one can relate a certain horospherical
subgroup of~$G$. Recall that by Corollary~\ref{loc.str.gen},
$X$~contains an open affine subset $\X\iso\Ru{P}\times A\times C$,
where $\Ru{P}$ is the unipotent radical of a parabolic subgroup
$P=P(X)$ and $A=L/L_0$ is a quotient torus of a Levi subgroup
$L\subseteq P$. Then $S(X)=L_0\semitimes\Ru{P}$ is the normalizer
of a generic $U$-orbit in~$X$.
\begin{definition}\label{hor.type}
The \emph{horospherical type} of~$X$ is the opposite horospherical
subgroup $S=S(X)^{-}=L_0\semitimes\Ru{P^{-}}$, where $\Ru{P^{-}}$
is the unipotent radical of the opposite parabolic subgroup
$P^{-}$ intersecting $P$ in~$L$.
\end{definition}
\begin{example}
The horospherical type of a horospherical homogeneous space $G/S$
is~$S$, because $G$ contains an open ``big cell'' $\Ru{P}\times
L\times\Ru{P^{-}}$, where $P^{-}=N_G(S)=L\semitimes\Ru{P^{-}}$.
For general horospherical varieties, the horospherical type is the
(conjugation class of) the stabilizer of general position
(Proposition~\ref{loc.str.hor}).
\end{example}

Complexity, rank and weight lattice can be read off the
horospherical type. Namely, it follows from
Corollary~\ref{loc.str.gen} that $c(X)=\dim X-\dim G+\dim S$,
$\RG(X)=\Ch(A)$, $r(X)=\dim A$, where $A=P^{-}/S$.

Every $G$-action can be deformed to a horospherical one of the
same type. This construction, called the \emph{horospherical
contraction}, was suggested by Popov~\cite{contr}. We review the
horospherical contraction in characteristic zero referring to
\cite[\S15]{HS} for arbitrary characteristic.

First consider an affine $G$-variety~$X$. Choose a one-parameter
subgroup $\gamma\in\CoCh(T)$ such that
$\langle\gamma,\lambda\rangle\geq0$ for any dominant weight and
any positive root~$\lambda$. Then
$\kk[X]^{(n)}=\bigoplus_{\langle\gamma,\lambda\rangle\leq
n}\kk[X]_{(\lambda)}$ is a $G$-stable filtration of~$\kk[X]$. The
algebra $\gr\kk[X]$ is finitely generated and has no nilpotents.
It is easy to see using Proposition~\ref{iso.mult} that
$X_0=\Spec\gr\kk[X]$ is a horospherical variety of the same type
as~$X$. Moreover, $\kk[X_0]^U\iso\kk[X]^U$ and
$\kk[X_0]\iso\kk[X]$ as $G$-modules. (Note that $S(X)$ may be
described as the common stabilizer of all $f\in\kk[X]^{(B)}$.)

Furthermore, $X_0$ may be described as the zero-fiber of a flat
family over~$\AAA^1$ with a generic fiber~$X$. Namely, let
$E=\Spec\bigoplus_{n=0}^{\infty}\kk[X]^{(n)}t^n\subseteq\kk[X][t]$.
The natural morphism $\delta:E\to\AAA^1$ is flat and
$G\times\kk^{\times}$-equivariant, where $G$ acts on $\AAA^1$
trivially and $\kk^{\times}$ acts by homotheties. Now
$\delta^{-1}(t)\iso X$ for $\forall t\neq0$, and
$\delta^{-1}(0)\iso X_0$.

If $X$ is an arbitrary $G$-variety, then we may replace it by a
birationally $G$-isomorphic projective variety, build an affine
cone $\widehat{X}$ over~$X$, and perform the above construction
for~$\widehat{X}$. Passing again to a projectivization and taking
a sufficiently small open $G$-stable subset, we obtain
\begin{proposition}\label{contr}
There exists a smooth $G\times\kk^{\times}$-variety $E$ and a
smooth $G\times\kk^{\times}$-morphism $\delta:E\to\AAA^1$ such
that $X_t=\delta^{-1}(t)$ is $G$-isomorphic to an open smooth
$G$-stable subset of~$X$ for $\forall t\neq0$, and $X_0$ is a
smooth horospherical variety of the same type as~$X$.
\end{proposition}

\section{Geometry of cotangent bundle}
\label{cotangent}

To a smooth $G$-variety~$X$, we relate a Poisson $G$-action on the
cotangent bundle $T^{*}X$ equipped with a natural symplectic
structure. Remarkably, the invariants of $G:X$ introduced
in~\ref{c&r(basics)} are closely related to the symplectic
geometry of $T^{*}X$ and to the respective moment map. In
particular, one obtains effective formulae for complexity and rank
involving symplectic invariants of $G:T^{*}X$. This theory was
developed by Knop in~\cite{W_X}. To the end of this chapter, we
assume $\ch{\kk}=0$.

Let $X$ be a smooth variety. A standard symplectic structure
on~$T^{*}X$ \cite[\S37]{SG} is given by a $2$-form
$\omega=d\mathbf{p}\wedge d\mathbf{q}=\sum dp_i\wedge dq_i$, where
$\mathbf{q}=(q_1,\dots,q_n)$ is a tuple of local coordinates
on~$X$, and $\mathbf{p}=(p_1,\dots,p_n)$ is an impulse, i.e.,
tuple of dual coordinates in a cotangent space. In a
coordinate-free form, $\omega=d\ell$, where a $1$-form $\ell$ on
$T^{*}X$ is given by $\langle\ell(\alpha),\xi\rangle=
\langle\alpha,d\pi(\xi)\rangle$, $\forall\xi\in T_{\alpha}T^{*}X$,
and $\pi:T^{*}X\to X$ is the canonical projection.

This symplectic structure defines the Poisson bracket of functions
on~$T^{*}X$. Another way to define this Poisson structure is to
consider the sheaf $\Dd_X$ of differential operators on~$X$. There
is an increasing filtration $\Dd_X=\bigcup\Dd_X^m$ by the order of
a differential operator and the isomorphism
$\gr\Dd_X\iso\Sym^{\bullet}\Tt_X=\pi_{*}\Oo_{T^{*}X}$ given by the
symbol map. Since $\gr\Dd_X$ is commutative, the commutator in
$\Dd_X$ induces the Poisson bracket on $\Oo_{T^{*}X}$ by the rule
%*
\begin{equation*}
\{\partial_1\bmod\Dd_X^{m-1},\partial_2\bmod\Dd_X^{n-1}\}=
[\partial_1,\partial_2]\bmod\Dd_X^{m+n-2},\qquad
\forall\partial_1\in\Dd_X^m,\ \partial_2\in\Dd_X^n
\end{equation*}
%*

If $X$ is a $G$-variety, then the symplectic structure on $T^{*}X$
is $G$-invariant, and for $\forall\xi\in\g$, the velocity field
$\xi*$ on $T^{*}X$ has a Hamiltonian $H_{\xi}=\xi*$, the
respective velocity field on $X$ considered as a linear function
on~$T^{*}X$ \cite[App.5]{SG}. Furthermore, the action $G:T^{*}X$
is Poisson, i.e., the map $\xi\mapsto H_{\xi}$ is a homomorphism
of $\g$ to the Poisson algebra of functions on~$T^{*}X$
\cite[App.5]{SG}. The comorphism $\Phi:T^{*}X\to\g^{*}$,
%*
\begin{equation*}
\langle\Phi(\alpha),\xi\rangle=H_{\xi}(\alpha)= \langle\alpha,\xi
x\rangle,\qquad \forall\alpha\in T_x^{*}X,\ \xi\in\g^{*},
\end{equation*}
%*
is called the \emph{moment map}. By $M_X\subseteq\g^{*}$ we denote
the closure of its image. Also set $L_X=M_X\by G$.

The moment map is $G$-equivariant \cite[App.5]{SG}. Clearly
$\langle d_{\alpha}\Phi(\nu),\xi\rangle=\omega(\nu,\xi\alpha)$ for
$\forall\nu\in T_{\alpha}T^{*}X,\ \xi\in\g$. It follows that $\Ker
d_{\alpha}\Phi=(\g\alpha)^{\sort}$, $\Im
d_{\alpha}\Phi=(\g_\alpha)^{\ann}$, where ${}^{\sort}$ and
${}^{\ann}$ denote the skew-orthocomplement and the annihilator
in~$\g^{*}$, respectively.
\begin{example}\label{moment(G/H)}
If $X=G/H$, then $T^{*}X=G\itimes{H}\h^{\ann}$ and
$\Phi(g*\alpha)=g\alpha$, $\forall g\in G,\ \alpha\in\h^{\ann}$.
Thus $M_{G/H}=\overline{G\h^{\ann}}$. In the general case, for
$\forall x\in X$, the moment map restricted to $T^{*}X|_{Gx}$
factors as $\Phi:T^{*}X|_{Gx}\onto T^{*}Gx\to M_{Gx}\subseteq
M_X$. We shall see below that for general~$x$, $M_{Gx}=M_X$.
\end{example}

The cohomomorphism $\Phi^{*}$ exists in two versions---a
commutative and a non-commutative one. Let $\Uni\g$ denote the
universal enveloping algebra of~$\g$, and $\Dd(X)$ be the algebra
of differential operators on~$X$. The action $G:X$ induces a
homomorphism $\Phi^{*}:\Uni\g\to\Dd(X)$ mapping each $\xi\in\g$ to
a $1$-order differential operator $\xi*$ on~$X$. The map
$\Phi^{*}$ is a homomorphism of filtered algebras, and the
associated graded map
%*
\begin{equation*}
\gr\Phi^{*}:\gr\Uni\g\iso\Sym^{\bullet}\g=\kk[\g^{*}]\longrightarrow
\gr\Dd(X)\subseteq\kk[T^{*}X],\qquad\xi\mapsto H_{\xi}
\end{equation*}
%*
is the commutative version of the cohomomorphism. The isomorphism
$\gr\Uni\g\iso\Sym^{\bullet}\g$ is provided by the
Poincar{\'e}--Birkhoff--Witt theorem, and the embedding
$\gr\Dd(X)\subseteq\kk[T^{*}X]$ is the symbol map.

On the sheaf level, we have the homomorphisms
$\Phi^{*}:\Oo_X\otimes\g\to\Tt_X$, $\Oo_X\otimes\Uni\g\to\Dd_X$.
Let $\Gg_X=\Phi^{*}(\Oo_X\otimes\g)$ denote the \emph{action
sheaf} (generated by velocity fields) and
$\Uu_X=\Phi^{*}(\Oo_X\otimes\Uni\g)$. Clearly,
%*
\begin{equation*}
T^{\g}X:=\Spec_{\Oo_X}\Sym^{\bullet}\Gg_X=\overline{\Im(\pi\times\Phi)}
\subseteq X\times\g^{*}
\end{equation*}
%*
The moment map factors as
%*
\begin{equation*}
\Phi:T^{*}X\longrightarrow T^{\g}X
\stackrel{\overline\Phi}{\longrightarrow}\g^{*}
\end{equation*}
%*
The (non-empty) fibers of $\pi\times\Phi$ are affine translates of
the conormal spaces to $G$-orbits. Generic fibers of $T^{\g}X\to
X$ are the cotangent spaces $\g_x^{\ann}=T^{*}_xGx$ to generic
orbits. The morphism $\overline\Phi$ is called the \emph{localized
moment map} \cite[\S2]{HC-hom}.
\begin{example}\label{comp.moment}
If $X$ is a smooth completion of a homogeneous space $\HS=G/H$,
then $T^{\g}X\supset T^{*}\HS$ and $\overline\Phi$ is a proper map
extending $\Phi:T^{*}\HS\to\g^{*}$. Thus one compactifies the
moment map of a homogeneous cotangent bundle.
\end{example}

\begin{definition}
A smooth $G$-variety $X$ is called \emph{pseudo-free} if $\Gg_X$
is locally free or $T^{\g}X$ is a vector bundle over~$X$. In other
words, the rational map
$X\dasharrow\Gr_k(\g)\iso\Gr_{n-k}(\g^{*})$,
$x\mapsto[\g_x]\mapsto[\g_x^{\ann}]$ ($n=\dim\g$, $k=\dim G_x$ for
$x\in X$ in general position) extends to~$X$, i.e., generic
isotropy subalgebras degenerate at the boundary to specific
limits.
\end{definition}
\begin{example}
For trivial reasons, $X$~is pseudo-free if the action $G:X$ is
generically free. Also, $X$~is pseudo-free if all $G$-orbits in
$X$ have the same dimension: in this case
$T^{\g}X=\bigsqcup_{Gx\subseteq X}T^{*}Gx$ \cite[2.3]{HC-hom}.
\end{example}

It is instructive to note that every $G$-variety $X$ has a
pseudo-free resolution of singularities $\check{X}\to X$: just
consider the closure of the graph of $X\dasharrow\Gr_k(\g)$ and
take for $\check{X}$ its equivariant desingularization.

It is easy to see that $\gr\Phi^{*}$ behaves well on the sheaf
level on pseudo-free varieties.
\begin{proposition}[{\cite[2.6]{HC-hom}}]
If $X$ is pseudo-free, then the filtrations on $\Uu_X$ induced
from $\Uni\g$ and $\Dd_X$ coincide and
$\gr\Phi^{*}:\Oo_X\otimes\Sym^{\bullet}\g\to\Sym^{\bullet}\Tt_X$
is surjective onto $\gr\Uu_X\iso\Sym^{\bullet}\Gg_X$.
\end{proposition}

Sometimes it is useful to replace the usual cotangent bundle by
its logarithmic version \cite[\S15]{toric.survey},
\cite[3.1]{toric.book}.
\begin{question}
Original reference?
\end{question}
Let $D\subset X$ be a divisor with normal crossings (e.g.,
smooth). The sheaf $\Omega^1_X(\log D)$ of differential 1-forms
with at most logarithmic poles along $D$ is locally generated by
$df/f$ with $f$ invertible outside~$D$. It contains $\Omega^1_X$
and is locally free. The respective vector bundle $T^{*}X(\log D)$
is said to be the \emph{logarithmic cotangent bundle}. The dual
vector bundle $TX(-\log D)$, called the \emph{logarithmic tangent
bundle}, corresponds to the subsheaf $\Tt_X(-\log D)\subset\Tt_X$
of vector fields preserving the ideal sheaf $\Ii_D\normin\Oo_X$.

If $D$ is $G$-stable, then the velocity fields of $G$ on $X$ are
tangent to~$D$, i.e., $\Gg_X\subseteq\Tt_X(-\log D)$. By duality,
we obtain the logarithmic moment map $\Phi:T^{*}X(\log D)\to
T^{\g}X\to\g^{*}$ extending the usual one on~$T^{*}(X\setminus
D)$.

\begin{remark}
We assume that $X$ is smooth in order to use the notions of
symplectic geometry. However, the moment map may be defined for
any $G$-variety~$X$.
\begin{question}
What is $T^{*}X$ for non-smooth~$X$?
\end{question}
As the definition is local and $M_X$ is a $G$-birational invariant
of~$X$, we may always pass to a smooth open subset of~$X$ and
conversely, to a (maybe singular) $G$-embedding of a smooth
$G$-variety.
\end{remark}

We are going to describe the structure of~$M_X$. We do it first
for horospherical varieties. Then we contract any $G$-variety to a
horospherical one and show that this contraction does not change
$M_X$.

Let $X$ be a horospherical variety of type~$S$. It is clear from
Proposition~\ref{loc.str.hor} that $M_X=M_{G/S}$. The moment map
factors as
%*
\begin{equation*}
\Phi_{G/S}:G\itimes{S}\s^{\ann}\stackrel{\pi_A}{\longrightarrow}
G\itimes{P^{-}}\s^{\ann}\stackrel{\overline{\Phi}}{\longrightarrow}
\g^{*},
\end{equation*}
%*
where $P^{-}=N_G(S)$, and $\pi_A$ is the quotient map modulo
$A=P^{-}/S$. By Proposition~\ref{collapse}, the map
$\overline{\Phi}$ is proper, whence $M_{G/S}=G\s^{\ann}$. If we
identify $\g^{*}$ with $\g$ via a non-degenerate $G$-invariant
inner product on~$\g$, then $\s^{\ann}$ is identified with
$\ab\oplus\Ru{\p^{-}}$ and is retracted onto $\ab$ by a certain
one-parameter subgroup of~$Z(L)$. It follows that $\ab\iso\ab^{*}$
intersects all closed $G$-orbits in~$M_{G/S}$, and
$L_{G/S}=\pi_G(\ab^{*})$, where $\pi_G:\g^{*}\to\g^{*}\by G$ is
the quotient map.

Finally, generic fibers of $\overline{\Phi}$ are finite. Indeed,
it suffices to find at least one finite fiber. But
$\overline{\Phi}^{-1}(G\Ru{\p^{-}})=G\itimes{P^{-}}\Ru{\p^{-}}$
maps onto $G\Ru{\p^{-}}$ with finite generic fibers by a theorem
of Richardson \cite[5.1]{Ad}.

We sum up in
\begin{theorem}\label{hor.moment}
For any horospherical $G$-variety $X$ of type~$S$, the natural map
$G\itimes{P^{-}}\s^{\ann}\to M_X=G\s^{\ann}$ is generically
finite, proper and surjective, and $L_X=\pi_G(\ab^{*})$.
\begin{question}
$M_X=M_{X^{*}}$ by \cite[5.5]{Ad}.
\end{question}
\end{theorem}

We have already seen that horospherical varieties, their cotangent
bundles and moment maps are easily accessible for study. A deep
result of Knop says that the closure of the image of the moment
map depends only on the horospherical type.
\begin{theorem}\label{moment}
Assume that $X$ is a $G$-variety of horospherical type~$S$. Then
$M_X=M_{G/S}$.
\end{theorem}
In the physical language, the idea of the proof is to apply
quantum technique to classical theory. We study the homomorphism
$\Phi_X^{*}:\Uni\g\to\Dd(X)$ and show that its kernel
$\Ii_X\normin \Uni\g$ depends only on the horospherical type. Then
we deduce that
$I_X=\Ker\gr\Phi_X^{*}=\Ideal{M_X}\normin\kk[\g^{*}]$ depends only
on the type of~$X$, which is the desired assertion. We retain the
notation of Proposition~\ref{contr}.

\begin{lemma}\label{ideals}
$\Ii_X=\Ii_{X_0}=\Ii_{G/S}$.
\end{lemma}
In the affine case, lemma stems from a $G$-module isomorphism
$\kk[X]\iso\kk[X_0]$. Indeed, $\Phi_X^{*}$ depends only on the
$G$-module structure of~$\kk[X]$. The general case is deduced from
the affine one~\cite[5.1]{W_X} with the help of affine cones and
of the next lemma.

Put $\Mm_X=\Im\Phi_X^{*}\subseteq\Dd(X)$. By the previous lemma,
$\Mm_X\iso\Mm_{X_0}\iso\Mm_{G/S}$.
\begin{lemma}\label{fin.mod}
$\Mm_{G/S}$ is a finite $\Uni\g$-module, and $\gr\Mm_{G/S}$ is a
finite $\kk[\g^{*}]$-module.
\end{lemma}
\begin{proof}
The first assertion follows from the second one. To prove the
second assertion, observe that $A=P^{-}/S$ acts on $G/S$ by
$G$-automorphisms. Therefore $\Mm_{G/S}\subseteq\Dd(G/S)^A$ and
$\gr\Mm_{G/S}\subseteq\kk[T^{*}G/S]^A=\kk[G\itimes{P^{-}}\s^{\ann}]$,
the latter being a finite $\kk[\g^{*}]$-module by
Theorem~\ref{hor.moment}.
\end{proof}

The restriction maps $\Mm_E\isoto\Mm_{X_t}$ agree with $\Phi^{*}$
and do not rise the order of a differential operator. We identify
$\Mm_E$ and $\Mm_{X_t}$ via this isomorphism and denote by
$\ord_E\partial$ ($\ord_{X_t}\partial$) the order of a
differential operator $\partial$ on~$E$ (resp.\ on~$X_t$).

Theorem~\ref{moment} follows from
\begin{lemma}
On $\Mm_X\iso\Mm_E\iso\Mm_{X_0}$,
$\ord_X\partial=\ord_E\partial=\ord_{X_0}\partial$
\begin{question}
Include the filtration $\ord$ induced from $\Uni\g$
(cf.~\cite[2.6]{HC-hom})? No: $\gr(\Uni\g/\Ii)$ may contain
nilpotents. Note
$\gr\Mm_{G/P^{-}}=\kk[G\itimes{P^{-}}\Ru{\p^{-}}]$
\cite[4.6]{W_X}, but
$\Spec\gr(\Uni\g/\Ii_{G/P^{-}})=G\Ru{\p^{-}}$.
\end{question}
for $\forall\partial$.
\end{lemma}
\begin{proof}
The first equality is clear, because an open subset of $E$ is
$G$-isomorphic to $X\times\kk^{\times}$. It follows from
Lemma~\ref{fin.mod} that the orders of a given differential
operator on $E$ and on~$X_0$ do not differ very much. Indeed,
$\Mm_{X_0}=\Mm_{G/S}=\sum(\Uni\g)\partial_i$ and
$\gr\Mm_{X_0}=\gr\Mm_{G/S}=\sum\kk[\g^{*}]\partial_i$ for some
$\partial_1,\dots,\partial_s\in\Mm_{X_0}$. Put
$d_i=\ord_{X_0}\partial_i$ and $d=\max_i\ord_E\partial_i$. If
$\ord_{X_0}\partial=n$, then $\partial=\sum u_i\partial_i$ for
some $u_i\in U^{n-d_i}\g$, hence $\ord_E\partial\leq n+d$.

However, if $\ord_{X_0}\partial<\ord_E\partial$, then
$\ord_{X_0}\partial^{d+1}<\ord_E\partial^{d+1}-d$, a
contradiction.
\end{proof}

The explicit description of $M_X$ in terms of the horospherical
type allows to examine invariant-theoretic properties of the
action $G:M_X$.

In the above notation, put $M=Z_G(\ab)\supseteq L$. Every
$G$-orbit in $M_X$ is of the form $Gx$,
$x\in\s^{\ann}=\ab\oplus\Ru{\p^{-}}$. Consider the Jordan
decomposition $x=x_{\sms}+x_{\nil}$. The (unique) closed orbit in
$Gx$ is~$Gx_{\sms}$. Moving $x$ by~$\Ru{P^{-}}$, we may assume
$x_{\sms}\in\ab$. If $x$ is a general point, then
$G_{x_{\sms}}=M$, $\z(x_{\sms})=\m$, thence
$x_{\nil}\in\m\cap\Ru{\p^{-}}$.

The concept of a general point can be specified as follows:
consider the principal open stratum $\ab^{\pr}\subseteq\ab$
obtained by removing all proper intersections with kernels of
roots and with $W$-translates of~$\ab$. Then $G_{\xi}=M$ for
$\forall\xi\in\ab^{\pr}$, and $G$-orbits intersect $\ab^{\pr}$ in
orbits of a finite group $W(\ab)=N_G(\ab)/M$ acting freely
on~$\ab^{\pr}$. Furthermore,
%*
\begin{equation*}
M_X^{\pr}:=\pi_G^{-1}\pi_G(\ab^{\pr})\iso
G\itimes{N(\ab)}(\ab^{\pr}+\M),\qquad\text{where }
\M=N(\ab)(\m\cap\Ru{\p^{-}})
\end{equation*}
%*
is a nilpotent cone in~$\m$.

\begin{definition}\label{symp.stab}
The Poisson action $G:T^{*}X$ is said to be \emph{symplectically
stable} if the action $G:M_X$ is stable, i.e., generic $G$-orbits
in $M_X$ are closed. By the above discussion, symplectic stability
is equivalent to $M=L$.
\end{definition}

This class of actions is wide enough.
\begin{proposition}\label{qaff=>ss}
If $X$ is quasiaffine, then $T^{*}X$ is symplectically stable.
\end{proposition}
\begin{proof}
The horospherical contraction $X_0$ and a typical orbit $G/S$
therein are quasiaffine, too. If $M\supset L$, then there is a
root $\alpha$ w.r.t\ a maximal torus $T\subseteq L$ such that
$\alpha|_{\ab}=0$ and $\g_{\alpha}\subseteq\Ru{\p^{-}}$. The
respective coroot $\alpha^{\vee}$ lies in~$\lv_0$. Let
$\s_{\alpha}=
\g_{\alpha}\oplus\langle\alpha^{\vee}\rangle\oplus\g_{-\alpha}$ be
the corresponding $\sgl_2$-subalgebra of~$\g$. Then
$\s_{\alpha}\cap\s=\langle\alpha^{\vee}\rangle\oplus\g_{-\alpha}$
is a Borel subalgebra in~$\s_{\alpha}$, and an orbit in $G/S$ of
the respective subgroup $S_{\alpha}\subseteq G$ is isomorphic
to~$\PP^1$, a contradiction with quasiaffinity.
\end{proof}

\begin{remark}\label{symp.stab=>stab}
A symplectically stable action is stable, i.e., generic orbits of
$G:T^{*}X$ are closed. Indeed, $\Phi$~is smooth along
$\Phi^{-1}(x)$ for general $x\in M_X$, whence
$(\g\alpha)^{\sort}=\Ker d_{\alpha}\Phi= T_{\alpha}\Phi^{-1}(x)$
have one and the same dimension $\dim T^{*}X-\dim M_X$ for
$\forall\alpha\in\Phi^{-1}(x)$. It follows that all orbits over
$Gx$ are closed in $\Phi^{-1}(Gx)$.
\begin{question}
Is it true for singular~$X$?
\end{question}
\end{remark}

The Poisson $G$-action on $T^{*}X$ provides two important
invariants:
\begin{definition}
The \emph{defect} $\df T^{*}X$ is the defect of the symplectic
form restricted to a generic $G$-orbit.

The \emph{corank} $\cork T^{*}X$ is the rank of the symplectic
form on the skew-orthogonal complement to the tangent space of a
generic $G$-orbit.

In other words,
%*
\begin{align*}
\df T^{*}X&=\dim(\g\alpha)^{\sort}\cap\g\alpha \\
\cork T^{*}X&=
\dim(\g\alpha)^{\sort}/(\g\alpha)^{\sort}\cap\g\alpha
\end{align*}
%*
for general $\alpha\in T^{*}X$.
\end{definition}

The cohomomorphism $\gr\Phi^{*}$ maps $\kk[\g^{*}]^G$ onto a
Poisson-commutative subalgebra $\Aa_X\subseteq\kk[T^{*}X]^G$
isomorphic to $\kk[L_X]$. Skew gradients of functions in~$\Aa_X$
commute, are $G$-stable, and both skew-orthogonal and tangent to
$G$-orbits. Indeed, for $\forall f\in\Aa_X,\ \alpha\in T^{*}X$,
$df$~is zero on~$\g\alpha$ (since $f$ is $G$-invariant) and
on~$(\g\alpha)^{\sort}=\Ker d_{\alpha}\Phi$ (because $f$ is pulled
back under~$\Phi$).

Those skew gradients generate a flow of $G$-automorphisms
preserving $G$-orbits on~$T^{*}X$, which is called a
\emph{$G$-invariant collective motion}. The restriction of this
flow to $G\alpha$ is a connected Abelian subgroup (in fact, a
torus) $A_{\alpha}\subseteq\Aut_G(G\alpha)$ with the Lie algebra
$\ab_{\alpha}\subseteq\n(\g_{\alpha})/\g_{\alpha}$.

For general~$\alpha$, $\Phi^{-1}\Phi(G\alpha)$ are level varieties
for~$\Aa_X$, because $G$-invariant regular functions separate
generic $G$-orbits in $M_X\subseteq\g^{*}$. It follows that
%*
\begin{align*}
\Ker d_{\alpha}\Aa_X&= T_{\alpha}\Phi^{-1}\Phi(G\alpha)=
(\g\alpha)+(\g\alpha)^{\sort},\qquad\text{and}\\
\ab_{\alpha}&=(\g\alpha)^{\sort}\cap(\g\alpha)=
T_{\alpha}(G\alpha\cap\Phi^{-1}\Phi(\alpha))\iso
\g_{\Phi(\alpha)}/\g_{\alpha}
\end{align*}
%*
In particular,
$\g_{\Phi(\alpha)}\supset\g_{\alpha}\supset\g_{\Phi(\alpha)}'$.

The defect of $G:T^{*}X$ is the dimension of the invariant
collective motion: $\df T^{*}X=\dim\ab_{\alpha}$ for general
$\alpha\in T^{*}X$.

The next theorem links the geometry of $X$ and the symplectic
geometry of~$T^{*}X$.
\begin{theorem}\label{c&r<=T^*}
Put $n=\dim X$, $c=c(X)$, $r=r(X)$. Then $\dim M_X=2n-2c-r$, $\df
T^{*}X=d_G(M_X)=r$, $d_G(T^{*}X)=2c+r$, $\cork T^{*}X=2c$.
\end{theorem}
\begin{proof}
In the notation of Definition~\ref{hor.type}, we have a
decomposition $\g=\Ru{\p}\oplus\lv_0\oplus\ab\oplus\Ru{\p^{-}}$,
where $\s=\lv_0\oplus\Ru{\p^{-}}$ is (the Lie algebra of) the
horospherical type of~$X$, and $\s^{\ann}$ is identified with
$\ab\oplus\Ru{\p^{-}}$ via $\g\iso\g^{*}$. By
Corollary~\ref{loc.str.gen}, $\dim\Ru{\p}=\dim\Ru{\p^{-}}=n-c-r$,
whence by Theorems~\ref{hor.moment}--\ref{moment}, $\dim M_X=\dim
G/P^{-}+\dim\s^{\ann}=\dim\Ru{\p}+\dim\Ru{\p^{-}}+\dim\ab=2n-2c-r$,
and $d_G(M_X)=\dim L_X=\dim\ab=r$. For general $\alpha\in T^{*}X$,
we have $d_G(T^{*}X)=\codim
G\alpha=\dim(\g\alpha)^{\sort}=\dim\Phi^{-1}(\alpha)=2n-\dim
M_X=2c+r$, and $\df T^{*}X=\dim
G\alpha\cap\Phi^{-1}\Phi(\alpha)=\dim
G_{\Phi(\alpha)}/G_{\alpha}=\dim G\alpha-\dim
G\Phi(\alpha)=(2n-d_G(T^{*}X))-(\dim M_X-d_G(M_X))=d_G(M_X)=r$.
Finally, $\cork T^{*}X=d_G(T^{*}X)-\df T^{*}X=2c$.
\end{proof}

Another application of the horospherical contraction and of the
moment map is the existence of the stabilizer of general position
for the $G$-action in a cotangent bundle.
\begin{theorem}[{\cite[\S8]{W_X}}]\label{sgp(T^*)}
Stabilizers in $G$ of general points in $T^{*}X$ are conjugate to
a stabilizer of the open orbit of $M\cap S$
in~$\m\cap\Ru{\p^{-}}$, in the above notation. In the
symplectically stable (e.g., quasiaffine) case, generic
stabilizers of $G:T^{*}X$ are conjugate to~$L_0$.
\end{theorem}
\begin{proof}
We prove the first assertion for horospherical~$X$. The general
case is derived form the horospherical one with the aid of the
horospherical contraction using Theorem~\ref{moment} and the
invariant collective motion, see~\cite[8.1]{W_X}, or
Remark~\ref{sgp(ss)} for the symplectically stable case.

We may assume $X=G/S$. A generic stabilizer of $G:T^{*}G/S$ equals
a generic stabilizer of $S:\s^{\ann}=\ab\oplus\Ru{\p^{-}}$. Take a
general point $x\in\s^{\ann}$ and let $x=x_{\sms}+x_{\nil}$ be the
Jordan decomposition. Moving $x$ by~$\Ru{P^{-}}$, we may assume
that $x_{\sms}$ is a general point in~$\ab$. Then
$S_{x_{\sms}}=M\cap S$, $x_{\nil}\in\m\cap\Ru{\p^{-}}$, and
$S_x=(M\cap S)_{x_{\nil}}$ is the stabilizer of a general point
in~$\m\cap\Ru{\p^{-}}$. But $M\cap S$ has the same orbits in
$\m\cap\Ru{\p^{-}}$ as a parabolic subgroup $M\cap P^{-}\subseteq
M$, because these two groups differ by a central torus in~$M$. By
\cite[Th.5.1]{Ad}, $M\cap P^{-}$ has an open orbit in the Lie
algebra of its unipotent radical $\m\cap\Ru{\p^{-}}$, which proves
the first assertion of the theorem.

If $X$ is symplectically stable, then $M=L$, whence the second
assertion.
\end{proof}

The last two theorems reduce the computation of complexity and
rank to studying generic orbits and stabilizers of a
\emph{reductive} group. Namely, it suffices to know generic
$G$-modalities of $T^{*}X$ and~$M_X$. We have a formula
%*
\begin{equation}
2c(X)+r(X)=d_G(T^{*}X)=2\dim X-\dim G+\dim G_{*},
\end{equation}
%*
where $G_{*}$ is the stabilizer of general position for
$G:T^{*}X$. For quasiaffine~$X$,
%*
\begin{equation}
r(X)=\rk G-\rk G_{*}
\end{equation}
%*
Furthermore, $\RG(X)$ is the group of characters of $T$ vanishing
on $T\cap G_{*}$, where $T$ is a maximal torus
normalizing~$G_{*}$. For homogeneous spaces, everything is reduced
even to \emph{representations} of reductive groups,
see~\ref{c&r(hom)}.

Now we explain another approach to computing complexity and rank
based on the theory of doubled actions \cite{c&r(hom)},
\cite[Ch.1]{c&r}. This approach is parallel to Knop's one and
coincides with the latter in the case of $G$-modules.

Let $\theta$ be a Weyl involution of $G$ relative to~$T$, i.e., an
involution of $G$ acting on $T$ as an inversion. Then
$\theta(P)=P^{-}$ for any parabolic $P\supseteq B$.
\begin{example}
If $G=\GL_n(\kk)$, or $\SL_n(\kk)$, and $T$ is the standard
diagonal torus, then we may put
$\theta(g)=w_G(g^{\tran})^{-1}w_G^{-1}$, where $w_G$ permutes the
standard basis of~$\kk^n$ in the reverse order:
$w_Ge_i=e_{n+1-i}$.
\end{example}
\begin{definition}
The \emph{dual} $G$-variety $X^{*}$ is a copy of $X$ equipped with
a twisted $G$-action: $gx^{*}=(\theta(g)x)^{*}$, where $x\mapsto
x^{*}$ is a fixed isomorphism $X\isoto X^{*}$.

The diagonal action $G:X\times X^{*}$ is called the \emph{doubled
action} w.r.t.\ $G:X$.
\end{definition}
\begin{remark}
If $V$ is a $G$-module, then $V^{*}$ is the dual $G$-module with a
fixed linear $G^{\theta}$-isomorphism $V\isoto V^{*}$. Similarly,
$\PP(V)^{*}\iso\PP(V^{*})$. If $X\subseteq\PP(V)$ is a
quasiprojective $G$-variety, then $X^{*}\subseteq\PP(V^{*})$.
\end{remark}
\begin{remark}
For a $G$-module $V$, the doubled $G$-variety $V\oplus V^{*}$ is
nothing else, but the cotangent bundle~$T^{*}V$.
\end{remark}

The following theorems are parallel to Theorems \ref{sgp(T^*)}
and~\ref{c&r<=T^*}.
\begin{theorem}\label{sgp(dbl)}
Stabilizers in $G$ of general points in $X\times X^{*}$ are
conjugate to~$L_0$.
\end{theorem}
\begin{theorem}\label{c&r<=dbl}
Let $G_{*}$ be the stabilizer of general position for $G:X\times
X^{*}$. Then
%*
\begin{align}
2c(X)+r(X)=d_G(X\times X^{*})&=2\dim X-\dim G+\dim G_{*}
\label{2c+r<=dbl} \\
r(X)&=\rk G-\rk G_{*} \label{r<=dbl}
\end{align}
%*
and $\RG(X)$ is the group of characters of a maximal torus $T$
normalizing~$G_{*}$ which vanish on $T\cap G_{*}$.
\end{theorem}
\begin{proof}[Proofs]
Consider an open embedding $\Ru{P}\times Z\embeds X$ from
Corollary~\ref{loc.str.gen}. Then $\Ru{P^{-}}\times Z^{*}\embeds
X^{*}$, where $Z^{*}$ is the dual $L$-variety to~$Z$. One deduces
that $G\itimes{L}(Z\times Z^{*})\embeds X\times X^{*}$ is an open
embedding,
\begin{question}
Unclear!
\end{question}
and the generic stabilizer of $L:Z\times Z^{*}$ equals~$L_0$,
whence the assertion of Theorem~\ref{sgp(dbl)}.
Theorem~\ref{c&r<=dbl} follows from Theorem~\ref{sgp(dbl)} with
the aid of Corollary~\ref{loc.str.gen}.
\end{proof}
\begin{example}\label{dbl(G/P)}
If $X=G/P$ is a projective homogeneous space, then $X^{*}=G/P^{-}$
and the stabilizer of general position for $G:X\times X^{*}$
equals~$P\cap P^{-}=L$.
\end{example}

We have a nice invariant-theoretic property of doubled actions on
affine varieties.
\begin{theorem}[{\cite[1.6]{U-res}}, {\cite[1.3.13]{c&r}}]
\label{stab(dbl)} If $X$ is affine, then generic $G$-orbits on
$X\times X^{*}$ are closed.
\end{theorem}

For a $G$-module~$V$, a stabilizer of general position for the
doubled $G$-module (or the cotangent bundle) $V\oplus V^{*}$ can
be found by an effective recursive algorithm relying on the
Brion--Luna--Vust construction (see~\ref{local}). In the notation
of Lemma~\ref{loc.str.V}, we have an isomorphism
%*
\begin{equation*}
G\itimes{L}(\oo{W}\times\oo{W}^{*})\isoto\oo{V}\times\oo{V}^{*}\subset
V\oplus V^{*}
\end{equation*}
%*
A stabilizer of general position for $G:V\oplus V^{*}$ equals that
for $L:W\oplus W^{*}$. Replacing $G:V$ by $L:W$, we apply the
Brion--Luna--Vust construction again, and so on. We obtain a
descending sequence of Levi subgroups $L_i\subseteq G$ and
$L_i$-submodules $W_i\subseteq V$. As the semisimple rank of $L_i$
decreases, the sequence terminates, and on the final \ordinal{$s$}
step, $L_s'$ acts on $W_s$ trivially. Then $G_{*}$ is just the
kernel of $L_s:W_s$.
\begin{example}
Let $G=\Spin_{10}(\kk)$ and $V=V(\omega_4)$ be one of its
half-spinor representations. In the notation of~\cite{Sem}, the
positive roots of $G$ are $\eps_i-\eps_j$, $\eps_i+\eps_j$, the
simple roots are $\alpha_i=\eps_i-\eps_{i+1}$,
$\alpha_5=\eps_4+\eps_5$, where $1\leq i<j\leq5$. The weights of
$V$ are $\frac12(\pm\eps_1\pm\dots\pm\eps_5)$, where the number of
minuses is odd, and $\omega_4^{*}=\omega_5=
\frac12(\eps_1+\dots+\eps_5)$.

Let $P_1$ be the projective stabilizer of a highest weight vector
$v_1^{*}\in V^{*}$. It's Levi subgroup $L_1$ has
$\alpha_1,\dots,\alpha_4$ as simple roots. The weights of
$\Ru{(\p_1^{-})}v_1^{*}$ are of the form
$\frac12(\pm\eps_1\pm\dots\pm\eps_5)$ with 2~minuses, and the
weights of $W_1=(\Ru{(\p_1^{-})}v_1^{*})^{\ann}$ have 1 or 5
minuses. Clearly, $L_1$ is of type $\AAa_4$, and $W_1$ is the
direct sum of a $1$-dimensional $L_1$-submodule of the weight
$\frac12(-\eps_1-\dots-\eps_5)$ and of a $5$-dimensional
$L_1$-submodule with highest weight
$\frac12(\eps_1+\dots\eps_4-\eps_5)$.

Take a highest weight vector $v_2^{*}\in W_1^{*}$ of highest
weight $\frac12(\eps_1-\eps_2-\dots-\eps_5)$, and let $P_2$ be its
projective stabilizer. The second Levi subgroup $L_2\subset P_2$
has $\alpha_2,\alpha_3,\alpha_4$ as simple roots, the weights of
$\Ru{(\p_2^{-})}v_2^{*}$ are
$\frac12(-\eps_1\pm\eps_2\pm\dots\pm\eps_5)$ with exactly 1~plus,
and $W_2=(\Ru{(\p_2^{-})}v_2^{*})^{\ann}$ has the weights
$\frac12(-\eps_1+\eps_2+\dots+\eps_5)$,
$\frac12(-\eps_1-\eps_2-\dots-\eps_5)$.

It is easy to see that $L_2'$ is exactly the kernel of $L_2:W_2$.
Thus our algorithm terminates, and we obtain $G_{*}\iso
L_2'\iso\SL_4$, $\RG(V)=\Ch(T/T\cap
L_2')=\langle\omega_1,\omega_5\rangle$ (here $\omega_1=\eps_1$),
$r(V)=2$, and $c(V)=0$. Moreover, $\kk[V]_1=V^{*}=V(\omega_5)$ and
$\kk[V]_2=S^2V(\omega_5)=V(2\omega_5)\oplus V(\omega_1)$, whence
$\RG_{+}(V)=\langle\omega_1,\omega_5\rangle$.
\end{example}

\section{Complexity and rank of homogeneous spaces}
\label{c&r(hom)}

We apply the methods developed in~\ref{cotangent} to computing
complexity and rank of homogeneous spaces.

The cotangent bundle of $G/H$ is identified
with~$G\itimes{H}\h^{\ann}$, where $\h^{\ann}\iso(\g/\h)^{*}$ is
the annihilator of~$\h$ in~$\g^{*}$. The representation
$H:\h^{\ann}$ is the \emph{coisotropy representation} at $\1H\in
G/H$. If we identify $\g$ with $\g^{*}$ via a non-degenerate
$G$-invariant inner product on~$\g$, then $\h^{\ann}$ is just the
orthogonal complement of~$\h$.

If $H$ is reductive, then $\g=\h\oplus\h^{\ann}$ as $H$-modules.
In particular, the isotropy and coisotropy representations are
isomorphic.

The following theorem is a reformulation of
Theorems~\ref{c&r<=T^*}--\ref{sgp(T^*)} for homogeneous spaces.
\begin{theorem}\label{c&r(G/H)}
Generic stabilizers of the coisotropy representation are all
conjugate to a certain subgroup $H_{*}\subseteq H$. For complexity
and rank of~$G/H$, we have the equations:
%*
\begin{align}
2c(G/H)+r(G/H)&=d_H(\h^{\ann})=\dim\h^{\ann}-\dim
H+\dim H_{*} \label{2c+r(G/H)} \\
&=\dim G-2\dim H+\dim H_{*} \notag \\
r(G/H)&=\dim G_{\alpha}-\dim H_{\alpha} \label{r(G/H)} \\
2c(G/H)&=\dim G\alpha-2\dim H\alpha \label{2c(G/H)} \\
%*
\intertext{where $\alpha\in\h^{\ann}$ is a general point
considered as an element of~$\g^{*}$. If $H$ is observable (e.g.,
reductive), then $H_{*}=L_0$ is the Levi subgroup of the
horospherical type of~$H$,}
%*
r(G/H)&=\rk G-\rk H_{*}, \label{r(qaff)}
\end{align}
%*
and $\RG(G/H)$ is the group of characters of $T$ vanishing
on~$T\cap H_{*}$.
\begin{question}
Add $c(G/H)=\dim U-\dim H+\dim B(H_{*})$
\end{question}
\end{theorem}
\begin{proof}
Generic modalities and stabilizers of the actions $G:T^{*}G/H$ and
${H:\h^{\ann}}$ coincide. This implies all assertions
except~\eqref{r(G/H)},~\eqref{2c(G/H)}. We have
$\Phi(\1*\alpha)=\alpha$ and $G_{\1*\alpha}=H_{\alpha}$, whence
the r.h.s.\ of~\eqref{r(G/H)} is the dimension of the invariant
collective motion, which yields~\eqref{r(G/H)}. Subtracting
\eqref{r(G/H)} from~\eqref{2c+r(G/H)} yields~\eqref{2c(G/H)}.
\end{proof}
Formulae \eqref{2c+r(G/H)}, \eqref{r(qaff)} are most helpful,
especially for reductive~$H$, because stabilizers of general
position for reductive group representations are known, e.g., from
Elashvili's tables~\cite{sgp(s)}, \cite{sgp(ss)}.
\begin{example}\label{Sp/Sp}
Let $G=\Sp_{2n}(\kk)$, $n\geq2$, and $H=\Sp_{2n-2}(\kk)$ be the
stabilizer of a general pair of vectors in a symplectic space, say
$e_1,e_2\in\kk^{2n}$, where a (standard) symplectic form on
$\kk^{2n}$ is $\omega=\sum x_{2i-1}\wedge x_{2i}$. Then $\Ad
G=\Sym^2\kk^{2n}$, $\Ad H=\Sym^2\kk^{2n-2}$ and $\Ad G|_H=
\Sym^2\kk^2\oplus\kk^2\otimes\kk^{2n-2}\oplus\Sym^2\kk^{2n-2}$,
where $\kk^2=\langle e_1,e_2\rangle$ and $\kk^{2n-2}=\langle
e_3,\dots,e_{2n}\rangle$. Hence the coisotropy representation of
$H$ equals $\kk^{2n-2}\oplus\kk^{2n-2}\oplus\kk^3$, where
$\kk^{2n-2}$ is the tautological and $\kk^3$ a trivial
representation of~$H$.

Clearly, $H_{*}=\Sp_{2n-4}(\kk)$ is the stabilizer of
$e_3,e_4\in\kk^{2n-2}$. It follows that $r(G/H)=2$ and
$2c(G/H)+r(G/H)=2(2n-2)+3-(n-1)(2n-1)+(n-2)(2n-3)=4$, whence
$c(G/H)=1$. Furthermore, the standard diagonal torus $T=
\{t=\diag(t_1,t_1^{-1},\dots,t_n,t_n^{-1})\}\subseteq G$
normalizes~$H_{*}$, and $\RG(G/H)=\Ch(T/T\cap H_{*})=
\langle\omega_1,\omega_2\rangle$, where $\omega_1(t)=t_1$,
$\omega_2(t)=t_1t_2$ are the first two fundamental weights of~$G$.
\end{example}

Theorem~\ref{c&r(G/H)} reduces the computation of complexity and
rank of affine homogeneous spaces to finding stabilizers of
general position for reductive group representations. The next
theorem does the same thing for arbitrary homogeneous spaces.

Consider a regular embedding $H\subseteq Q$ in a parabolic
subgroup $Q\subseteq G$. Let $K\subseteq M$ be Levi subgroups of
$H$ and~$Q$. Clearly, $K$ acts on $\Ru{Q}/\Ru{H}$ by conjugations,
and this action is isomorphic to a linear action
$K:\Ru{\q}/\Ru{\h}$ via the exponential map.

We may assume that $M\supseteq T$ and $B\subseteq Q^{-}$. Then
$B(M)=B\cap M$ is a Borel subgroup in~$M$. We may assume that
$\1K\in M/K$ is a general point, i.e., $M_{*}=K\cap\theta(K)$ is a
stabilizer of general position for $M:M/K\times(M/K)^{*}$, and
$B(M_{*})=B(M)\cap K$ is a stabilizer of general position for
$B(M):M/K$. By Theorem~\ref{sgp(dbl)}, $M_{*}$ and $B(M_{*})$ are
normalized by $T$, and $B(M_{*})^0$ is a Borel subgroup
in~$M_{*}^0$. $M_{*}$ may be non-connected, but it is a direct
product of $M_{*}^0$ and a finite subgroup of~$T$. The notions of
complexity, rank and weight lattice generalize to $M_{*}$-actions
immediately.
\begin{theorem}[{\cite[1.2]{c(nil)}, \cite[2.5.20]{c&r}}]
\label{c&r(nqaff)} With the above choice of $H$, $Q$, and $K$
among conjugates,
%*
\begin{align*}
c_G(G/H)&=c_M(M/K)+c_{M_{*}}(\Ru{Q}/\Ru{H}) \\
r_G(G/H)&=r_M(M/K)+r_{M_{*}}(\Ru{Q}/\Ru{H})
\end{align*}
%*
and there is a canonical exact sequence of weight lattices
%*
\begin{equation*}
0\longrightarrow\RG_M(M/K)\longrightarrow\RG_G(G/H)
\longrightarrow\RG_{M_{*}}(\Ru{Q}/\Ru{H})\longrightarrow0
\end{equation*}
%*
\end{theorem}
\begin{proof}
As $B\subseteq Q^{-}$, the $B$- and even $U$-orbit of $\1Q$ is
open in~$G/Q$. Hence codimensions of generic orbits and weight
lattices for the actions $B:G/H\iso G\itimes{Q}Q/H$ and $B\cap
Q=B(M):Q/H$ are equal. Further, $Q/H\iso
M\itimes{K}\Ru{Q}/\Ru{H}$. It follows with our choice of $K$ that
the codimension of a generic $B(M)$-orbit in $M/K$ is the sum of
the codimension of a generic $B(M)$-orbit in $M/K$ and of a
generic $B(M_{*})$-orbit in $\Ru{Q}/\Ru{H}$, whence the formula on
complexities.

Furthermore, stabilizers of general position of the actions
$B:G/H$, $B(M):Q/H$, $B(M_{*}):\Ru{Q}/\Ru{H}$ are all equal to
$B(L_0)=B\cap L_0$, where $L_0$ is the $T$-regular Levi subgroup
of the horospherical type of~$G/H$. The equalities
$\RG(G/H)=\Ch(B/B(L_0))$, $\RG(M/K)=\Ch(B(M)/B(M_{*}))$,
$\RG(\Ru{Q}/\Ru{H})=\Ch(B(M_{*})/B(L_0))$ imply the assertions on
ranks and weight lattices.
\end{proof}

Thus the computation of complexity and rank of $G/H$ is performed
in two steps:
\begin{enumerate}
\item Compute the group $M_{*}\subseteq K$, which is by
Theorem~\ref{sgp(T^*)} a stabilizer of general position for the
coisotropy representations $K:\ka^{\ann}$ (the orthocomplement
in~$\m$). This can be done using, e.g., Elashvili's tables. \item
Compute the stabilizer of general position for
$M_{*}:\Ru{\q}/\Ru{\h}\oplus(\Ru{\q}/\Ru{\h})^{*}$ using, e.g., an
algorithm at the end of~\ref{cotangent}.
\end{enumerate}
Complexity and rank are read off these stabilizers.
\begin{example}
Let $G=\Sp_{2n}(\kk)$, $n\geq3$, and $H$ be the stabilizer of a
general triple of vectors in a symplectic space, say
$e_1,e_2,e_3\in\kk^{2n}$, in the notation of Example~\ref{Sp/Sp}.
We choose $K=\Sp_{2n-4}(\kk)$, the stabilizer of
$e_1,e_2,e_3,e_4\in\kk^{2n}$, for a Levi subgroup of~$H$. The
unipotent radical of~$\h$ is
%*
\begin{equation*}
\Ru{\h}=\left\{\;
\begin{array}{|@{}c@{}|@{}c@{}|@{}c@{}|}
\hline
\begin{array}{cc}
0 & 0 \\
0 & 0
\end{array} &
\text{\Large$0$} & \text{\Large$0$} \\
\hline \text{\Large$0$} &
\begin{array}{cc}
0 & z \\
0 & 0
\end{array} &
\begin{array}{c}
\text{\large\strut} v^{\tran}\Omega \\
\hline 0\;\Dots+0\;0
\end{array} \\
\hline \text{\Large$0$} &
\begin{array}{c|c}
       0       &   \\
    \Dots0+    & v \\
       0       &
\end{array} &
\text{\LARGE$0$} \\
\hline
\end{array}
\right.\; \left|\;
\vphantom{\begin{array}{|@{}c@{}|@{}c@{}|@{}c@{}|} \hline
\begin{array}{cc}
0 & 0 \\
0 & 0
\end{array} &
\text{\Large$0$} & \text{\Large$0$} \\
\hline \text{\Large$0$} &
\begin{array}{cc}
0 & z \\
0 & 0
\end{array} &
\begin{array}{c}
\text{\large\strut} v^{\tran}\Omega \\
\hline 0\;\Dots+0\;0
\end{array} \\
\hline \text{\Large$0$} &
\begin{array}{c|c}
       0       &   \\
    \Dots0+    & v \\
       0       &
\end{array} &
\text{\LARGE$0$} \\
\hline
\end{array}}
z\in\kk,\ v\in\kk^{2n-4} \;\right\}
\end{equation*}
%*
where $\Omega$ is the matrix of the symplectic form on
$\kk^{2n-4}=\langle e_5,\dots,e_{2n}\rangle$. For $Q$ we take the
stabilizer of a flag $\langle e_1\rangle\subset\langle
e_1,e_3\rangle\subset\kk^{2n}$. Then
%*
\begin{align*}
M&=\left\{\;
\begin{array}{|c|c|}
\hline
\begin{array}{@{}c@{\;}c@{}c@{\;}c@{}}
t_1 & & & \\
    & t_1^{-1} & & \smash{\text{\Large$0$}} \\
 & & t_2 & \\
\smash{\text{\Large$0$}} & & & t_2^{-1} \\
\end{array} &
\text{\LARGE$0$} \\
\hline
\text{\LARGE\strut} \text{\LARGE$0$} & \text{\LARGE$A$} \\
\hline
\end{array}
\right.\; \left|\; \vphantom{\begin{array}{|c|c|} \hline
\begin{array}{@{}c@{\;}c@{}c@{\;}c@{}}
t_1 & & & \\
    & t_1^{-1} & & \smash{\text{\Large$0$}} \\
 & & t_2 & \\
\smash{\text{\Large$0$}} & & & t_2^{-1} \\
\end{array} &
\text{\LARGE$0$} \\
\hline
\text{\LARGE\strut} \text{\LARGE$0$} & \text{\LARGE$A$} \\
\hline
\end{array}}
t_1,t_2\in\kk^{\times},\ A\in\Sp_{2n-4}(\kk)
\;\right\} \\
\Ru{\q}&=\left\{\;
\begin{array}{|@{}c@{}|@{}c@{}|@{}c@{}|}
\hline
\begin{array}{cc}
0 & x \\
0 & 0
\end{array} &
\begin{array}{cc}
-y_2 & y_1 \\
  0  &  0
\end{array} &
\begin{array}{c}
\text{\large\strut} u^{\tran}\Omega \\
\hline 0\;\Dots+0\;0
\end{array} \\
\hline
\begin{array}{cc}
0 & y_1 \\
0 & y_2
\end{array} &
\begin{array}{cc}
0 & z \\
0 & 0
\end{array} &
\begin{array}{c}
\text{\large\strut} v^{\tran}\Omega \\
\hline 0\;\Dots+0\;0
\end{array} \\
\hline
\begin{array}{c|c}
       0       &   \\
    \Dots0+    & u \\
       0       &
\end{array} &
\begin{array}{c|c}
       0       &   \\
    \Dots0+    & v \\
       0       &
\end{array} &
\text{\LARGE$0$} \\
\hline
\end{array}
\right.\; \left|\;
\vphantom{\begin{array}{|@{}c@{}|@{}c@{}|@{}c@{}|} \hline
\begin{array}{cc}
0 & x \\
0 & 0
\end{array} &
\begin{array}{cc}
-y_2 & y_1 \\
  0  &  0
\end{array} &
\begin{array}{c}
\text{\large\strut} u^{\tran}\Omega \\
\hline 0\;\Dots+0\;0
\end{array} \\
\hline
\begin{array}{cc}
0 & y_1 \\
0 & y_2
\end{array} &
\begin{array}{cc}
0 & z \\
0 & 0
\end{array} &
\begin{array}{c}
\text{\large\strut} v^{\tran}\Omega \\
\hline 0\;\Dots+0\;0
\end{array} \\
\hline
\begin{array}{c|c}
       0       &   \\
    \Dots0+    & u \\
       0       &
\end{array} &
\begin{array}{c|c}
       0       &   \\
    \Dots0+    & v \\
       0       &
\end{array} &
\text{\LARGE$0$} \\
\hline
\end{array}}
x,y_1,y_2,z\in\kk,\ u,v\in\kk^{2n-4} \;\right\}
\end{align*}
%*
Clearly, $M/K\iso(\kk^{\times})^2$, and $M_{*}=K=\Sp_{2n-4}(\kk)$
acts on $\Ru{\q}$ by left multiplication of $u,v\in\kk^{2n-4}$ by
$A\in\Sp_{2n-4}(\kk)$. It follows that
$\Ru{\q}/\Ru{\h}=\kk^{2n-4}\oplus\kk^3$, a sum of the tautological
and a trivial $\Sp_{2n-4}(\kk)$-module.

We deduce that $c(M/K)=0$, $r(M/K)=2$, and $\RG(M/K)=
\langle\omega_1,\omega_2\rangle$. A generic stabilizer of
$M_{*}:\Ru{\q}/\Ru{\h}\oplus(\Ru{\q}/\Ru{\h})^{*}$ equals
$\Sp_{2n-6}(\kk)$ (=the stabilizer of $e_1,\dots,e_6$), whence
$\RG(\Ru{\q}/\Ru{\h})=\langle\overline{\omega}_3\rangle$, where
$\overline{\omega}_3$ is the first fundamental weight
of~$\Sp_{2n-6}(\kk)$, or equivalently, the restriction to the
diagonal torus in $\Sp_{2n-6}(\kk)$ of the third fundamental
weight $\omega_3(t)=t_1t_2t_3$ of~$\Sp_{2n}(\kk)$. It follows that
$r(\Ru{\q}/\Ru{\h})=1$, $2c(\Ru{\q}/\Ru{\h})+r(\Ru{\q}/\Ru{\h})=
2(2n-4+3)-(n-2)(2n-3)+(n-3)(2n-5)=7$, hence
$c(\Ru{\q}/\Ru{\h})=3$. We conclude that $c(G/H)=r(G/H)=3$, and
$\RG(G/H)=\langle\omega_1,\omega_2,\omega_3\rangle$.
\end{example}

\section{Spaces of small rank and complexity}
\label{c&r<=1}

The term ``complexity'' is justified by the fact that homogeneous
spaces of small complexity are more accessible for study. In
particular, a good compactification theory can be developed for
homogeneous spaces of complexity $\leq1$, see
Chapter~\ref{LV-theory}. On the other hand, rank and complexity
are not completely independent invariants of a homogeneous space.
In this section, we discuss the interactions between rank and
complexity, paying special attention to homogeneous spaces of
small rank and complexity. We begin with a simple
\begin{proposition}\label{r=0}
$r(G/H)=0$ iff $H$ is parabolic, i.e., $G/H$~is projective.
\end{proposition}
\begin{proof}
The ``only if'' implication follows from the Bruhat decomposition,
cf.~Example~\ref{cr(G/P)}. Conversely, if $r(G/H)=0$, then $H$
contains a maximal torus of~$G$. Replacing $H$ by a conjugate, we
may assume that $H\supseteq T$ and $B\cap H$ has the minimal
possible dimension. We claim that $H\supseteq B^{-}$. Otherwise,
there is a simple root $\alpha$ such that
$\h\not\supseteq\g_{-\alpha}$. Let
$P_{\alpha}=L_{\alpha}\semitimes N_{\alpha}$ be the respective
minimal parabolic subgroup with the $T$-regular Levi
decomposition. Then $B=B_{\alpha}\semitimes N_{\alpha}$, where
$B_{\alpha}=B\cap L_{\alpha}$ is a Borel subgroup in~$L_{\alpha}$,
and $H\cap L_{\alpha}=B_{\alpha}$ or~$T$. In both cases, we may
replace $B_{\alpha}$ by a conjugate Borel subgroup
$\widetilde{B}_{\alpha}$ in $L_{\alpha}$ so that $\dim
H\cap\widetilde{B}_{\alpha}<\dim H\cap B_{\alpha}$. Then for
$\widetilde{B}=\widetilde{B}_{\alpha}N_{\alpha}$ we have $\dim
H\cap\widetilde{B}<\dim H\cap B$, a contradiction.
\end{proof}
In particular, homogeneous spaces of rank zero have complexity
zero. This can be generalized to the following general inequality
between complexity and rank.
\begin{theorem}[{\cite[2.7]{c&r(hom)}, \cite[2.2.10]{c&r}}]
$2c(G/H)\leq\Cox G\cdot r(G/H)$, where $\Cox G$ is the maximum of
the Coxeter numbers of simple components of~$G$.
\end{theorem}
Observe that if $G$ is a simple group and $H=\{\1\}$, then the
inequality becomes an equality, since $c(G)=\dim U$, $r(G)=\rk G$,
and $\Cox G=2\dim U/\rk G$. This inequality is rather rough, and
various examples create an impression that the majority of
homogeneous spaces have either small complexity or large rank. In
particular, Panyushev proved that $r(G/H)=1$ implies
$c(G/H)\leq1$.

\begin{proposition}[\cite{r=1}]\label{r=1}
If $r(G/H)=1$, then either $c(G/H)=0$, or $G/H$ is obtained from a
homogeneous $\SL_2(\kk)$-space with finite stabilizer by parabolic
induction, whence $c(G/H)=1$.
\end{proposition}
Spherical homogeneous spaces of rank~$1$ where classified by
Akhiezer \cite{comp.div} and Brion~\cite{r(sph)=1}, see
Proposition~\ref{emb(sph,r=1)} and Table~\ref{wonder(r=1)}. The
above proposition says that, besides the spherical case, there is
only one essentially new example of rank~$1$, namely $\SL_2(\kk)$
acting on itself by left multiplications. (Factorizing by a finite
group preserves complexity and rank.) The proof (and the
classification) is based on Theorem~\ref{c&r(nqaff)}. Homogeneous
spaces of rank~$1$ are also characterized in terms of equivariant
completions, see Proposition~\ref{emb(sph,r=1)} and
Remark~\ref{emb(r=1)}.

Homogeneous spaces of small complexity are much more numerous.
Here classification results concern mainly the case, where $H$ is
reductive, i.e., $G/H$~is affine.
\begin{question}
Classification of spherical solvable $H$ (Luna), of spherical
spaces of classical type (Luna, Italians).
\end{question}
For simple $G$, affine homogeneous spaces of complexity~$0$ were
classified by Kr\"amer~\cite{sph(s)} and of complexity~$1$ by
Panyushev \cite{c=1(s)}, \cite[Ch.3]{c&r}. A complete
classification of spherical affine homogeneous spaces was obtained
by Mikityuk~\cite{sph(ss)} and Brion~\cite{class(sph)}, with a
final stroke put by Yakimova~\cite{sph.indec}.
\begin{question}
Classification of affine homogeneous spaces of complexity~$1$
(Arzhantsev--Chuvashova).
\end{question}
We expose their results in Tables~\ref{sph(s)}--\ref{sph(ss)}, see
also Table~\ref{symm}. Since the computation of complexity and
rank of a given homogeneous space represents no difficulties by
Theorems~\ref{c&r(G/H)}--\ref{c&r(nqaff)}, the main problem of
classification is to ``cut off infinity''.

Clearly, complexity and rank of $G/H$ do not change if we replace
$G$ by a finite cover and/or $H$ by a subgroup of finite index.
Thus complexity and rank depend only on the local isomorphism
class of~$G/H$, i.e., on the pair $(\g,\h)$. Therefore we may
assume that $H$ is connected.

If $G$ is not semisimple, then it decomposes in an almost direct
product $G=G'\cdot Z$, where $G'$ is its (semisimple) commutator
subgroup and $Z$ is the connected center of~$G$. It is easy to see
that complexities of $G/H$, $G/HZ$, and $G'/(HZ\cap G')$ are
equal. Therefore it suffices to solve the classification problem
for semisimple~$G$.

An initial arithmetical restriction on a subgroup $H\subseteq G$
such that $c(G/H)\leq c$ is that
%*
\begin{equation}\label{dim(H)}
\dim H\leq\dim U-c
\end{equation}
%*
A more subtle restriction is based on the notion of
$d$-decomposition~\cite{c=1(s)}. A triple of reductive groups
$(L,L_1,L_2)$ is called a \emph{$d$-decomposition} if
$d_{L_1\times L_2}(L)=d$, where $L_1\times L_2$ acts on $L$ by
left and right multiplications. Clearly, $(L,L_1,L_2)$ remains a
$d$-decomposition if one permutes $L_1,L_2$ or replaces them by
conjugates. Besides, $d_{L_1}(L/L_2)=d_{L_2}(L/L_1)=d$.
By~\cite{closed.orb}, generic orbits of each one of the actions
$L_1\times L_2:L$, $L_1:L/L_2$, $L_2:L/L_1$ are closed. In
particular, $0$-decompositions are indeed decompositions:
$L=L_1\cdot L_2$. They were classified by Onishchik~\cite{decomp}.
Some special kinds of $1$-decompositions of classical groups
occurring in the classification of homogeneous spaces of
complexity~$1$ (see below) were described by
Panyushev~\cite{c=1(s)}.

Let $H\subseteq F$ be reductive subgroups of~$G$, and $F_{*}$ be
the stabilizer of general position for the coisotropy
representation $F:\f^{\ann}$. We may assume that $\bp=\1F$ is a
general point for the $B$-action on~$G/F$, so that $\dim B\bp$ is
maximal and $B_{\bp}^0$ is a Borel subgroup in~$F_{*}^0$.
Immediately,
\begin{proposition}
We have
%*
\begin{equation*}
c_G(G/H)=c_G(G/F)+c_{F_{*}}(F/H)\geq c_G(G/F)+d_{F_{*}}(F/H)
\end{equation*}
%*
In particular, if $c(G/H)\leq c$, then $(F,F_{*},H)$ is a
$d$-decomposition for some $d\leq c$.
\end{proposition}
The latter assertion is the keystone in a method of classifying
affine spherical homogeneous spaces of simple groups suggested by
Mikityuk and extended by Panyushev to the case of complexity one.
Let us explain its core.

Let $G$ be a simple algebraic group. Maximal connected reductive
subgroups $F\subset G$ are known due to Dynkin \cite[Ch.6,
\S3]{LieG}. We choose among them those with $c(G/F)\leq1$ and
search for reductive $H\subset F$ such that still $c(G/H)\leq1$.

If $G$ is exceptional, then either $c(G/F)=0$ or $c(G/F)\geq2$.
For spherical~$F$, sorting out those $H\subset F$ which
satisfy~\eqref{dim(H)} gives only 4~new subgroups with
$c(G/H)\leq1$ (\Nos~\ref{E6/D5} of Table~\ref{sph(s)} and
\ref{E6/B4T1}--\ref{F4/D4} of Table~\ref{c=1(s)}).

If $G$ is classical, then inequality~\eqref{dim(H)} gives a finite
list of subgroups. Again $c(G/F)\neq1$ with only one exception
$G=\Sp_4(\kk)$, $F=\SL_2(\kk)$ embedded in~$\Sp_4(\kk)$ by a
$4$-dimensional irreducible representation (\No\,\ref{Sp4/SL2} of
Table~\ref{c=1(s)}). Here \eqref{dim(H)} becomes an equality, and
$F$ cannot be reduced. Sorting out $H\subset F$ with $c(G/H)\leq1$
is based on~\eqref{dim(H)} and on the fact that $(F,F_{*},H)$ is a
decomposition or $1$-decomposition. Here we find 22~new subgroups
(\Nos~\ref{SL/SLSL}--\ref{D4/G2} of Table~\ref{sph(s)} and
\ref{SL/2SL}--\ref{Sp/SL}, \ref{SO/SO}--\ref{B4/G2T1} of
Table~\ref{c=1(s)}).

If $G/H$ is a symmetric space, i.e., $H=(G^{\sigma})^0$, where
$\sigma$ is an involutive automorphism of~$G$, then $c(G/H)=0$
(Theorem~\ref{symm=>sph}). Symmetric spaces are considered in
\ref{symmetric} and classified in Table~\ref{symm}.

Up to a local isomorphism, all non-symmetric affine homogeneous
spaces of simple groups with complexity~$0$ are listed in
Table~\ref{sph(s)} and those of complexity~$1$ in
Table~\ref{c=1(s)}. In the tables, we use the following notation.

Fundamental weights of simple groups are numbered as
in~\cite{Sem}. By $\omega_i$ we denote fundamental weights of~$G$,
by $\omega'_i,\omega''_i,\dots$ those of simple components of~$H$,
and by $\eps_i$ basic characters of the central torus of~$H$. We
drop an index for a group of rank~$1$.

The column ``$H\embeds G$'' describes the embedding in terms of
the restriction to $H$ of the minimal representation of~$G$ (the
tautological representation for classical groups). We use the
multiplicative notation for representations: irreducible
representations are indicated by their highest weights expressed
in basic weights multiplicatively (i.e., products instead of sums,
powers instead of multiples, $1$~for the trivial one-dimensional
representation, etc.), and ``$+$'' stands for the sum of
representations.

The rank of $G/H$ is indicated in the column ``$r(G/H)$'', and the
column ``$\RG_{+}(G/H)$'' contains a minimal system of generators
for the weight semigroup.
\begin{question}
Verify, e.g., $\SO_{2n}/\GL_n$ ($n$ odd) and
$\SO_{10}/\Spin_7\times\SO_2$ ($k_3=0$, $0\leq
2k_1-|k_4-k_5|\divby4$), or assume $G$ simply connected, $H$
connected.
\end{question}

\begin{table}[!h]
\caption{Spherical affine homogeneous spaces of simple groups}
\label{sph(s)}
\begin{center}
\normalsize\makebox[0pt]{\begin{tabular}{|c|c|c|c|c|c|} \hline \No
& $G$ & $H$ & $H\embeds G$ & $r(G/H)$ &
$\RG_{+}(G/H)$ \\
\hline

\tabitem\label{SL/SLSL} & $\SL_n$ & $\SL_m\times\SL_{n-m}$ &
$\omega'_1+\omega''_1$ & $m+1$ &
$\omega_1+\omega_{n-1},\dots,\omega_{m-1}+\omega_{n-m+1},$ \\
&& $(m<n/2)$ &&& $\omega_m,\omega_{n-m}$ \\
\hline

\tabitem\label{SL/SpT1} & $\SL_{2n+1}$ &
$\Sp_{2n}\times\kk^{\times}$ &
$\omega'_1\eps+\eps^{-2n}$ & $2n-1$ & $\sum k_i\omega_i$ \\
&&&&& $\sum\limits_{i\text{ odd}}(2n+1-i)k_i=
\sum\limits_{i\text{ even}}ik_i$ \\
\hline

\tabitem & $\SL_{2n+1}$ & $\Sp_{2n}$ & $\omega'_1+1$ & $2n$ &
$\omega_1,\dots,\omega_{2n}$ \\
\hline

\tabitem\label{Sp/SpT1} & $\Sp_{2n}$ &
$\Sp_{2n-2}\times\kk^{\times}$ &
$\omega'_1+\eps+\eps^{-1}$ & $2$ & $2\omega_1,\omega_2$ \\
\hline

\tabitem\label{SO/GL} & $\SO_{2n+1}$ & $\GL_n$ &
$\omega'_1\eps+\omega'_{n-1}\eps^{-1}+1$ & $n$ &
$\omega_1,\dots,\omega_{n-1},2\omega_n$ \\
\hline

\tabitem & $\SO_{4n+2}$ & $\SL_{2n+1}$ & $\omega'_1+\omega'_{2n}$
& $n+1$
& $\omega_2,\omega_4,\dots,\omega_{2n},\omega_{2n+1}$ \\
\hline

\tabitem & $\SO_{10}$ & $\Spin_7\times\SO_2$ &
$\omega'_3+\eps+\eps^{-1}$ &
$4$ & $2\omega_1,\omega_2,\omega_4,\omega_5$ \\
\hline

\tabitem & $\SO_9$ & $\Spin_7$ & $\omega'_3+1$ & $2$ &
$\omega_1,\omega_4$ \\
\hline

\tabitem\label{D4/G2} & $\SO_8$ & $\GGg_2$ & $\omega'_1+1$ & $3$ &
$\omega_1,\omega_3,\omega_4$ \\
\hline

\tabitem\label{B3/G2} &
$\SO_7$ & $\GGg_2$ & $\omega'_1$ & $1$ & $\omega_3$ \\
\hline

\tabitem\label{E6/D5} & $\EEe_6$ & $\DDd_5$ &
$\omega'_1+\omega'_5+1$ & $3$ &
$\omega_1,\omega_5,\omega_6$ \\
\hline

\tabitem\label{G2/A2} & $\GGg_2$ & $\AAa_2$ &
$\omega'_1+\omega'_2+1$ & $1$ & $\omega_1$ \\
\hline

\end{tabular}}
\end{center}
\end{table}

\begin{table}[!h]
\caption{Affine homogeneous spaces of simple groups with
complexity~$1$} \label{c=1(s)}
\begin{center}
\small\makebox[0pt]{\begin{tabular}{|c|c|c|c|c|c|} \hline
\No & $G$ & $H$ & $H\embeds G$ & $r(G/H)$ & $\RG_{+}(G/H)$ \\
\hline

\tabitem\label{SL/2SL} & $\SL_{2n}$ & $\SL_n\times\SL_n$ &
$\omega'_1+\omega''_1$ & $n$ &
$\omega_1+\omega_{2n-1},\dots,\omega_{n-1}+\omega_{n+1},\omega_n$ \\
\hline

\tabitem & $\SL_n$ & $\SL_{n-2}\times(\kk^{\times})^2$ &
$\omega'_1\eps_1\eps_2+\eps_1^{2-n}+\eps_2^{2-n}$ & $3$ &
$\omega_1+\omega_{n-1},\omega_2+\omega_{n-2}$ \\
& ($n\geq3$) &&& ($n>3$) &
$2\omega_1+\omega_{n-2},\omega_2+2\omega_{n-1}$ \\
\cline{5-6}
&&&& $2$ & $\omega_1+\omega_2,3\omega_1,3\omega_2$ \\
&&&& ($n=3$) & \\
\hline

\tabitem & $\SL_n$ & $\SL_{n-2}\times\kk^{\times}$ &
$\omega'_1\eps+\eps^{d_1}+\eps^{d_2}$ & $4$ & \\
& ($n\geq5$) && $d_1\neq d_2$ && \\
&&& $d_1+d_2=2-n$ && \\
\hline

\tabitem & $\SL_6$ & $\Sp_4\times\SL_2\times\kk^{\times}$ &
$\omega'_1\eps+\omega''\eps^{-2}$ & $5$ & \\
\hline

\tabitem & $\Sp_{2n}$ & $\Sp_{2n-2}$ & $\omega'_1+2$ & $2$ &
$\omega_1,\omega_2$ \\
\hline

\tabitem & $\Sp_{2n}$ & $\Sp_{2n-4}\times\SL_2\times\SL_2$ &
$\omega'_1+\omega''+\omega'''$ & $3$ &
$\omega_1+\omega_3,\omega_2,\omega_4$ \\
& ($n\geq3$) &&& ($n>3$) & \\
\cline{5-6}
&&&& $2$ & $\omega_1+\omega_3,\omega_2$ \\
&&&& ($n=3$) & \\
\hline

\tabitem\label{Sp/SL} & $\Sp_{2n}$ & $\SL_n$ &
$\omega'_1+\omega'_{n-1}$ & $n$ &
$2\omega_1,\dots,2\omega_{n-1},\omega_n$ \\
\hline

\tabitem\label{Sp4/SL2} & $\Sp_4$ & $\SL_2$ & ${\omega'}^3$ & $2$
&
$4\omega_1,3\omega_2$ \\
&&&&& $4\omega_1+2\omega_2,6\omega_1+3\omega_2$ \\
\hline

\tabitem\label{SO/SO} & $\SO_n$ & $\SO_{n-2}$ & $\omega'_1+2$ &
$2$ & $\omega_1,\omega_2$ \\
& $(n\geq 4)$ &&&& \\
\hline

\tabitem & $\SO_{2n+1}$ & $\SL_n$ & $\omega'_1+\omega'_{n-1}+1$ &
$n$ &
$\omega_1,\dots,\omega_n$ \\
\hline

\tabitem & $\SO_{4n}$ & $\SL_{2n}$ & $\omega'_1+\omega'_{2n-1}$ &
$n$
& $\omega_2,\omega_4,\dots,\omega_{2n}$ \\
& $(n\geq 2)$ &&&& \\
\hline

\tabitem & $\SO_{11}$ & $\Spin_7\times\SO_3$ &
$\omega'_3+{\omega''}^2$ &
$5$ & $2\omega_1,2\omega_2,\omega_3,\omega_4$ \\
&&&&& $\omega_1+2\omega_5,\omega_2+2\omega_5$ \\
&&&&& $\omega_1+\omega_2+2\omega_5$ \\
\hline

\tabitem & $\SO_{10}$ & $\Spin_7$ & $\omega'_3+2$ & $4$ &
$\omega_1,\omega_2,\omega_4,\omega_5$ \\
\hline

\tabitem\label{B4/G2T1} & $\SO_9$ & $\GGg_2\times\SO_2$ &
$\omega'_1+\eps+\eps^{-1}$ &
$4$ & \\
\hline

\tabitem\label{E6/B4T1} & $\EEe_6$ & $\BBb_4\times\kk^{\times}$ &
$(\omega'_1+1)\eps^2+\omega'_4\eps^{-1}+\eps^{-4}$ &
$5$ & \\
\hline

\tabitem & $\EEe_7$ & $\EEe_6$ &
$\omega'_1+\omega'_5+2$ & $3$ & $\omega_1,\omega_2,\omega_6$ \\
\hline

\tabitem\label{F4/D4} & $\FFf_4$ & $\DDd_4$ &
$\omega'_1+\omega'_3+\omega'_4+2$ &
$2$ & $\omega_1,\omega_2$ \\
\hline

\end{tabular}}
\end{center}
\end{table}

Now we describe spherical affine homogeneous spaces of semisimple
groups.

We say that $G/H$ is \emph{decomposable} if, up to a local
isomorphism, $G=G_1\times G_2$, $H=H_1\times H_2$, and
$H_i\subseteq G_i$, $i=1,2$. Clearly, $G/H=G_1/H_1\times G_2/H_2$
is spherical iff $G_i/H_i$ are spherical. Thus it suffices to
classify indecomposable spherical spaces.

Let $H\subseteq G$ be a reductive subgroup. We say that $G/H$ is
\emph{strictly indecomposable} if $G/H'$ is still indecomposable.
All strictly indecomposable spherical affine homogeneous spaces of
semisimple (non-simple) groups are listed in Table~\ref{sph(ss)}.
\begin{question}
Verified.
\end{question}
The column ``$H\embeds G$'' describes the embedding in the
following way. White vertices of a diagram denote simple factors
of~$G$ and black vertices denote factors of~$H$. (Some factors may
vanish for small~$n$.) If a factor $H_j$ of~$H$ projects
non-trivially to a factor $G_i$ of~$G$, then the respective
vertices are joined by an edge. The product of those $H_j$ which
project to~$G_i$ is a spherical subgroup in~$G_i$, and its
embedding in~$G_i$ is described in Table~\ref{sph(s)}. It follows
from Tables~\ref{sph(s)},~\ref{sph(ss)} that $\dim Z(H)\leq1$ for
all strictly indecomposable spherical homogeneous spaces~$G/H$.

Now assume that $G/H$ is indecomposable, but not strictly. Then,
up to a local isomorphism, $G=G_1\times\dots\times G_s$ and
$H'=H_1'\times\dots\times H_s'$, where $H_i$ are the projections
of $H$ to~$G_i$, and $G_i/H_i$ are strictly indecomposable.
Furthermore, $H_i=H_i'Z_i$, where $Z_i$ is a one-dimensional
central torus, and $H=H'Z$, where $Z\subset Z_1\times\dots\times
Z_s$ is a subtorus. Since $G/H$ is indecomposable, $Z$~cannot be
decomposed as $Z'\times Z''$, where $Z',Z''$ are the projections
of $Z$ to the products of two disjoint sets of factors~$Z_i$.

If $G_i/H_i$ is spherical and $B_i\subset G_i$ is a Borel subgroup
such that $\dim B_i\cap H_i$ is minimal, then $\g_i=\br_i+\h_i'$
if $G_i/H_i'$ is spherical, or $\g_i=(\br_i+\h_i')\oplus\z_i$
otherwise. It follows that $G/H$ is spherical iff all $G_i/H_i$
are spherical and $Z$ projects onto the product of those $Z_i$ for
which $G_i/H_i'$ is not spherical. This completes the
classification.

\begin{table}[h]
\caption{Spherical affine homogeneous spaces of semisimple groups}
\label{sph(ss)}
\begin{center}
\makebox[0pt]{\begin{tabular}{|c|c||c|c|} \hline
\No & $H\embeds G$ & \No & $H\embeds G$ \\
\hline

\tabitem &
%TexCad Options
%\grade{\off}
%\emlines{\off}
%\beziermacro{\off}
%\reduce{\on}
%\snapping{\off}
%\quality{2.00}
%\graddiff{0.01}
%\snapasp{1}
%\zoom{1.00}
\unitlength 0.50ex \linethickness{0.4pt}
\begin{picture}(37.00,15.00)
\put(9.00,2.00){\circle*{2.00}} \put(24.00,2.00){\circle*{2.00}}
\put(9.00,12.00){\circle{2.00}} \put(24.00,12.00){\circle{2.00}}
\put(24.00,3.00){\line(0,1){8.00}}
\put(9.00,3.00){\line(0,1){8.00}}
\put(23.33,11.33){\line(-3,-2){13.33}}
\put(7.67,12.00){\makebox(0,0)[rc]{$\SL_n$}}
\put(7.67,2.00){\makebox(0,0)[rc]{$\SL_n$}}
\put(25.67,12.00){\makebox(0,0)[lc]{$\SL_{n+1}$}}
\put(25.67,2.00){\makebox(0,0)[lc]{$\kk^{\times}$}}
\end{picture} &

\tabitem &
%TexCad Options
%\grade{\off}
%\emlines{\off}
%\beziermacro{\off}
%\reduce{\on}
%\snapping{\off}
%\quality{2.00}
%\graddiff{0.01}
%\snapasp{1}
%\zoom{1.00}
\unitlength 0.50ex \linethickness{0.4pt}
\begin{picture}(38.00,15.00)
\put(15.00,2.00){\circle*{2.00}} \put(30.00,2.00){\circle*{2.00}}
\put(15.00,12.00){\circle{2.00}} \put(30.00,12.00){\circle{2.00}}
\put(30.00,3.00){\line(0,1){8.00}}
\put(15.00,3.00){\line(0,1){8.00}}
\put(15.67,11.33){\line(3,-2){13.33}}
\put(13.67,12.00){\makebox(0,0)[rc]{$\Sp_{2n}$}}
\put(13.67,2.00){\makebox(0,0)[rc]{$\Sp_{2n-4}$}}
\put(31.67,12.00){\makebox(0,0)[lc]{$\Sp_4$}}
\put(31.67,2.00){\makebox(0,0)[lc]{$\Sp_4$}}
\end{picture} \\
\hline

\tabitem &
%TexCad Options
%\grade{\off}
%\emlines{\off}
%\beziermacro{\off}
%\reduce{\on}
%\snapping{\off}
%\quality{2.00}
%\graddiff{0.01}
%\snapasp{1}
%\zoom{1.00}
\unitlength 0.50ex \linethickness{0.4pt}
\begin{picture}(58.00,15.00)
\put(15.00,2.00){\circle*{2.00}} \put(29.00,2.00){\circle*{2.00}}
\put(43.00,2.00){\circle*{2.00}} \put(22.00,12.00){\circle{2.00}}
\put(36.00,12.00){\circle{2.00}}
\put(13.67,2.00){\makebox(0,0)[rc]{$\GL_{n-2}$}}
\put(27.67,2.00){\makebox(0,0)[rc]{$\SL_2$}}
\put(44.33,2.00){\makebox(0,0)[lc]{$\Sp_{2m-2}$}}
\put(20.33,12.00){\makebox(0,0)[rc]{$\SL_n$}}
\put(37.67,12.00){\makebox(0,0)[lc]{$\Sp_{2m}$}}
\put(15.00,2.67){\line(3,4){6.33}}
\put(29.00,2.67){\line(-3,4){6.33}}
\put(43.00,2.67){\line(-3,4){6.33}}
\put(29.00,2.67){\line(3,4){6.33}}
\end{picture} &

\tabitem &
%TexCad Options
%\grade{\off}
%\emlines{\off}
%\beziermacro{\off}
%\reduce{\on}
%\snapping{\off}
%\quality{2.00}
%\graddiff{0.01}
%\snapasp{1}
%\zoom{1.00}
\unitlength 0.50ex \linethickness{0.4pt}
\begin{picture}(40.00,15.00)
\put(18.00,2.00){\circle*{2.00}} \put(11.00,12.00){\circle{2.00}}
\put(25.00,12.00){\circle{2.00}}
\put(16.67,2.00){\makebox(0,0)[rc]{$\SO_n$}}
\put(9.33,12.00){\makebox(0,0)[rc]{$\SO_n$}}
\put(26.67,12.00){\makebox(0,0)[lc]{$\SO_{n+1}$}}
\put(18.00,2.67){\line(-3,4){6.33}}
\put(18.00,2.67){\line(3,4){6.33}}
\end{picture} \\
\hline

\tabitem &
%TexCad Options
%\grade{\off}
%\emlines{\off}
%\beziermacro{\off}
%\reduce{\on}
%\snapping{\off}
%\quality{2.00}
%\graddiff{0.01}
%\snapasp{1}
%\zoom{1.00}
\unitlength 0.50ex \linethickness{0.4pt}
\begin{picture}(65.00,15.00)
\put(15.00,2.00){\circle*{2.00}} \put(29.00,2.00){\circle*{2.00}}
\put(43.00,2.00){\circle*{2.00}} \put(22.00,12.00){\circle{2.00}}
\put(36.00,12.00){\circle{2.00}}
\put(13.67,2.00){\makebox(0,0)[rc]{$\SL_{n-2}$}}
\put(27.67,2.00){\makebox(0,0)[rc]{$\SL_2$}}
\put(44.33,2.00){\makebox(0,0)[lc]{$\Sp_{2m-2}$}}
\put(20.33,12.00){\makebox(0,0)[rc]{$\SL_n$}}
\put(37.67,12.00){\makebox(0,0)[lc]{$\Sp_{2m}$}}
\put(15.00,2.67){\line(3,4){6.33}}
\put(29.00,2.67){\line(-3,4){6.33}}
\put(43.00,2.67){\line(-3,4){6.33}}
\put(29.00,2.67){\line(3,4){6.33}}
\put(50.00,7.00){\makebox(0,0)[lc]{($n\geq5$)}}
\end{picture} &

\tabitem &
%TexCad Options
%\grade{\off}
%\emlines{\off}
%\beziermacro{\off}
%\reduce{\on}
%\snapping{\off}
%\quality{2.00}
%\graddiff{0.01}
%\snapasp{1}
%\zoom{1.00}
\unitlength 0.50ex \linethickness{0.4pt}
\begin{picture}(68.00,27.00)
\put(15.00,2.00){\circle*{2.00}} \put(35.00,2.00){\circle*{2.00}}
\put(55.00,2.00){\circle*{2.00}} \put(15.00,12.00){\circle{2.00}}
\put(35.00,12.00){\circle{2.00}} \put(55.00,12.00){\circle{2.00}}
\put(35.00,22.00){\circle*{2.00}}
\put(13.67,2.00){\makebox(0,0)[rc]{$\Sp_{2n-2}$}}
\put(33.67,2.00){\makebox(0,0)[rc]{$\Sp_{2m-2}$}}
\put(56.33,2.00){\makebox(0,0)[lc]{$\Sp_{2l-2}$}}
\put(13.67,12.00){\makebox(0,0)[rc]{$\Sp_{2n}$}}
\put(33.67,12.00){\makebox(0,0)[rc]{$\Sp_{2m}$}}
\put(56.33,12.00){\makebox(0,0)[lc]{$\Sp_{2l}$}}
\put(33.67,22.00){\makebox(0,0)[rb]{$\Sp_2$}}
\put(15.00,3.00){\line(0,1){8.00}}
\put(15.00,11.00){\line(0,0){0.00}}
\put(35.00,3.00){\line(0,1){8.00}}
\put(55.00,3.00){\line(0,1){8.00}}
\put(35.00,13.00){\line(0,1){8.00}}
\put(34.00,21.67){\line(-2,-1){18.15}}
\put(36.00,21.67){\line(2,-1){18.15}}
\end{picture} \\
\hline

\tabitem &
%TexCad Options
%\grade{\off}
%\emlines{\off}
%\beziermacro{\off}
%\reduce{\on}
%\snapping{\off}
%\quality{2.00}
%\graddiff{0.01}
%\snapasp{1}
%\zoom{1.00}
\unitlength 0.50ex \linethickness{0.4pt}
\begin{picture}(58.00,15.00)
\put(15.00,2.00){\circle*{2.00}} \put(29.00,2.00){\circle*{2.00}}
\put(43.00,2.00){\circle*{2.00}} \put(22.00,12.00){\circle{2.00}}
\put(36.00,12.00){\circle{2.00}}
\put(13.67,2.00){\makebox(0,0)[rc]{$\Sp_{2n-2}$}}
\put(27.67,2.00){\makebox(0,0)[rc]{$\Sp_2$}}
\put(44.33,2.00){\makebox(0,0)[lc]{$\Sp_{2m-2}$}}
\put(20.33,12.00){\makebox(0,0)[rc]{$\Sp_{2n}$}}
\put(37.67,12.00){\makebox(0,0)[lc]{$\Sp_{2m}$}}
\put(15.00,2.67){\line(3,4){6.33}}
\put(29.00,2.67){\line(-3,4){6.33}}
\put(43.00,2.67){\line(-3,4){6.33}}
\put(29.00,2.67){\line(3,4){6.33}}
\end{picture} &

\tabitem &
%TexCad Options
%\grade{\off}
%\emlines{\off}
%\beziermacro{\off}
%\reduce{\on}
%\snapping{\off}
%\quality{2.00}
%\graddiff{0.01}
%\snapasp{1}
%\zoom{1.00}
\unitlength 0.50ex \linethickness{0.4pt}
\begin{picture}(72.00,15.00)
\put(15.00,2.00){\circle*{2.00}} \put(29.00,2.00){\circle*{2.00}}
\put(57.00,2.00){\circle*{2.00}} \put(22.00,12.00){\circle{2.00}}
\put(50.00,12.00){\circle{2.00}}
\put(13.67,2.00){\makebox(0,0)[rc]{$\Sp_{2n-2}$}}
\put(27.67,2.00){\makebox(0,0)[rc]{$\Sp_2$}}
\put(58.33,2.00){\makebox(0,0)[lc]{$\Sp_{2m-2}$}}
\put(20.33,12.00){\makebox(0,0)[rc]{$\Sp_{2n}$}}
\put(51.67,12.00){\makebox(0,0)[lc]{$\Sp_{2m}$}}
\put(15.00,2.67){\line(3,4){6.33}}
\put(29.00,2.67){\line(-3,4){6.33}}
\put(57.00,2.67){\line(-3,4){6.33}}
\put(29.00,2.67){\line(3,4){6.33}}
\put(43.00,2.00){\circle*{2.00}} \put(36.00,12.00){\circle{2.00}}
\put(41.67,2.00){\makebox(0,0)[rc]{$\Sp_2$}}
\put(34.33,12.00){\makebox(0,0)[rc]{$\Sp_4$}}
\put(43.00,2.67){\line(-3,4){6.33}}
\put(43.00,2.67){\line(3,4){6.33}}
\end{picture} \\
\hline

\tabitem & \multicolumn{3}{|c|}{
\begin{tabular}{c}
%TexCad Options
%\grade{\off}
%\emlines{\off}
%\beziermacro{\off}
%\reduce{\on}
%\snapping{\off}
%\quality{2.00}
%\graddiff{0.01}
%\snapasp{1}
%\zoom{1.00}
\unitlength 0.50ex \linethickness{0.4pt}
\begin{picture}(40.00,15.00)
\put(18.00,2.00){\circle*{2.00}} \put(11.00,12.00){\circle{2.00}}
\put(25.00,12.00){\circle{2.00}}
\put(16.67,2.00){\makebox(0,0)[rc]{$H$}}
\put(9.33,12.00){\makebox(0,0)[rc]{$H$}}
\put(26.67,12.00){\makebox(0,0)[lc]{$H$}}
\put(18.00,2.67){\line(-3,4){6.33}}
\put(18.00,2.67){\line(3,4){6.33}}
\end{picture}
\end{tabular}
\qquad
\begin{tabular}{c}
($H$ is any simple group)
\end{tabular}} \\
\hline

\end{tabular}}
\end{center}
\end{table}

\section{Double cones}
\label{double.cone}

The theory of complexity and rank can be applied to a fundamental
problem of representation theory: decompose a tensor product of
two simple $G$-modules into irreducibles. The idea is to realize
this tensor product as a $G$-submodule in the coordinate algebra
of a certain affine $G$-variety---a double cone---and to compute
the algebra of $U$-invariants on a double cone, which yields the
$G$-module structure of the whole coordinate algebra. In cases of
small complexity, the algebra of $U$-invariants can be effectively
computed.

We may and will assume that $G$ is a simply connected semisimple
group.

\begin{definition}
Let $\lambda\in\Ch_{+}$ be a dominant weight and
$v_{\lambda^{*}}\in V(\lambda^{*})$ be a highest weight vector. A
cone $C(\lambda)=\overline{Gv_{\lambda^{*}}}\subseteq
V(\lambda^{*})$ is called the \emph{cone of highest weight
vectors} (\emph{HV-cone}). Clearly,
$C(\lambda)=\overline{Gv_{-\lambda}}$, where $v_{-\lambda}\in
V(\lambda^{*})$ is a lowest weight vector.

The projectivization of $C(\lambda)$ is a projective homogeneous
space $G/P(\lambda^{*})\iso G/P(\lambda)^{-}$, where $P(\lambda)$
denotes the projective stabilizer of a vector of highest weight
$\lambda$ in~$V(\lambda)$.
\end{definition}

The following assertions on HV-cones are well known and easily
proved \cite[\S2]{S-var}, cf.~\ref{S-var}.
\begin{proposition}\label{HV}
\begin{enumerate}
\item $C(\lambda)=Gv_{-\lambda}\cup\{0\}$ is a normal conical
variety in~$V(\lambda^{*})$. \item $\kk[C(\lambda)]_n\iso
V(n\lambda)$ as a $G$-module. \item\label{HV.fact} $C(\lambda)$ is
factorial iff $\lambda$ is a fundamental weight.
\end{enumerate}
\end{proposition}

\begin{definition}
A variety $Z(\lambda,\mu)=C(\lambda)\times C(\mu)$ is said to be a
\emph{double cone}.
\end{definition}

The group $\widehat{G}=G\times(\kk^{\times})^2$ acts on
$Z(\lambda,\mu)$ in a natural way, where the factors
$\kk^{\times}$ act by homotheties. Thus $\kk[Z(\lambda,\mu)]$ is
bigraded and
%*
\begin{equation*}
\kk[Z(\lambda,\mu)]_{n,m}=V(n\lambda)\otimes V(m\mu)
\end{equation*}
%*
The algebra $\kk[Z(\lambda,\mu)]^U$ is finitely generated
(Theorem~\ref{U-inv}\ref{A<=>A^U}) and
$(\Ch_{+}\times\ZZ_{+}^2)$-graded, and it is clear from the above
that the knowledge of its (polyhomogeneous) generators and
syzygies provides immediately a series of decomposition rules for
$V(n\lambda)\otimes V(m\mu)$. Namely, the highest vectors of
irreducible summands in $V(n\lambda)\otimes V(m\mu)$ are (linearly
independent) products of generators of total bidegree~$(n,m)$.

The smaller is the $\widehat{G}$-complexity of~$Z(\lambda,\mu)$,
the simpler is the structure of~$\kk[Z(\lambda,\mu)]^U$. Say, if
$Z(\lambda,\mu)$ is $\widehat{G}$-spherical, then
$Z(\lambda,\mu)\by U$ is a toric $\widehat{T}$-variety, where
$\widehat{T}=T\times(\kk^{\times})^2$. Hence
$\kk[Z(\lambda,\mu)]^U$ is the semigroup algebra of the weight
semigroup of~$Z(\lambda,\mu)$ (cf.\ Example~\ref{toric}). If in
addition $Z(\lambda,\mu)$ is factorial, then $Z(\lambda,\mu)\by U$
is a factorial toric variety, hence $\kk[Z(\lambda,\mu)]^U$ is
freely generated by $\widehat{T}$-eigenfunctions of linearly
independent weights. This yields a very simple decomposition rule,
see below.

Therefore it is important to have a transparent method for
computing complexity and rank of double cones. By the theory of
doubled actions (\ref{cotangent}), the problem reduces to
computing the stabilizer of general position for the doubling
$Z\times Z^{*}$ of a double cone $Z=Z(\lambda,\mu)$. This was done
by Panyushev in~\cite{dbl.con}, see also~\cite[Ch.6]{c&r}. Here
are his results.

Let $L(\lambda)$ be the Levi subgroup of $P(\lambda)$
containing~$T$. The character $\lambda$ extends to~$L(\lambda)$.
Put $G(\lambda)=\Ker\lambda\subset L(\lambda)$. Denote by
$G(\lambda,\mu)$ the stabilizer of general position for
$G(\lambda):G/G(\mu)$ and by $L(\lambda,\mu)$ the stabilizer of
general position for $L(\lambda):G/L(\mu)$. Recall
from~\cite{closed.orb} that generic orbits of both these actions
are closed, hence the codimension of a generic orbit equals the
dimension of a categorical quotient:
%*
\begin{xalignat}{1}
\dim G(\lambda)\backslash\!\!\backslash G\by G(\mu)&= \dim G+\dim
G(\lambda,\mu)-\dim G(\lambda)-\dim G(\mu),
\label{dim(G/G*G)}\\
\dim L(\lambda)\backslash\!\!\backslash G\by L(\mu)&= \dim G+\dim
L(\lambda,\mu)-\dim L(\lambda)-\dim L(\mu), \label{dim(G/L*L)}
\end{xalignat}
%*
where $L_1\backslash\!\!\backslash L\by L_2$ denotes the
categorical quotient of the action $L_1\times L_2:L$ by left and
right multiplication. Also put
%*
\begin{equation*}
\PP(Z)=\PP(C(\lambda))\times\PP(C(\mu))\iso G/P(\lambda)^{-}\times
G/P(\mu)^{-}
\end{equation*}
%*
\begin{theorem}\label{sgp.dbl.con}
\begin{roster}
\item\label{sgp<-L1*L2:L} The stabilizers of general position for
the doubled actions $G:Z\times Z^{*}$, $G:\PP(Z)\times\PP(Z^{*})$,
$\widehat{G}:Z\times Z^{*}$ are equal to $G(\lambda,\mu)$,
$L(\lambda,\mu)$ and $\widehat{G}(\lambda,\mu)=
\{\,(g,\lambda(g),\mu(g))\mid g\in L(\lambda,\mu)\,\}$,
\begin{question}
Is it true?
\end{question}
respectively. \item\label{coiso->sgp} Put
$V(\lambda,\mu)=\Ru{\p(\lambda)}\cap\Ru{\p(\mu)}$.  Then
$G(\lambda,\mu)$ and $L(\lambda,\mu)$ are equal to the stabilizers
of general position for the doubled actions $G(\lambda)\cap
G(\mu):V(\lambda,\mu)\oplus V(\lambda,\mu)^{*}$, $L(\lambda)\cap
L(\mu):V(\lambda,\mu)\oplus V(\lambda,\mu)^{*}$, respectively.
\end{roster}
\end{theorem}
The proof of~\ref{sgp<-L1*L2:L} uses the following
\begin{lemma}\label{sgp(HV)}
The stabilizers of general position for the doubled actions
$G:C(\lambda)\times C(\lambda^{*})$,
$G:\PP(C(\lambda))\times\PP(C(\lambda^{*}))$,
$G\times\kk^{\times}:C(\lambda)\times C(\lambda^{*})$ (where
$\kk^{\times}$ acts on $C(\lambda)$ by homotheties) are equal to
$G(\lambda)$, $L(\lambda)$ and
$\widehat{L}(\lambda)=\{\,(g,\lambda(g))\mid g\in L(\lambda)\,\}$,
respectively.
\end{lemma}
\begin{proof}
Observe that $z=(v_{-\lambda},v_{\lambda})\in C(\lambda)\times
C(\lambda^{*})$ has stabilizer and projective stabilizer
$G(\lambda)$ in~$G$, whence $\codim Gz=1$, and $\overline{G\langle
z\rangle}=C(\lambda)\times C(\lambda^{*})$. It follows that
$G(\lambda)$ is the stabilizer of general position in~$G$. Other
assertions are proved similarly (cf.~Example~\ref{dbl(G/P)}).
\end{proof}
\begin{proof}[Proof of Theorem~\ref{sgp.dbl.con}]
\begin{roster}
\item[\ref{sgp<-L1*L2:L}] We have $Z\times Z^{*}=
(C(\lambda)\times C(\lambda^{*}))\times (C(\mu)\times
C(\mu^{*}))$, and similarly for $\PP(Z)\times\PP(Z^{*})$. The
stabilizer of general position for any diagonal action
$L:X_1\times X_2$ can be computed in two steps: first find the
stabilizers of general position $L_i$ for the actions $L:X_i$,
$i=1,2$, and then find the stabilizer of general position for
$L_1:L/L_2$. It remains to apply Lemma~\ref{sgp(HV)} for $L=G$ or
$\widehat{G}$ and $X_1=C(\lambda)\times C(\lambda^{*})$ or
$\PP(C(\lambda))\times\PP(C(\lambda^{*}))$, $X_2=C(\mu)\times
C(\mu^{*})$ or $\PP(C(\mu))\times\PP(C(\mu^{*}))$.
\item[\ref{coiso->sgp}] One can prove \ref{coiso->sgp} using
Luna's slice theorem, if one observes that the $L(\lambda)$-orbit
of $\1L(\mu)\in G/L(\mu)$ and the $G(\lambda)$-orbit of
$\1G(\mu)\in G/G(\mu)$ are closed, and computes the slice module.
However, the proof also stems from the theory of doubled actions.
It suffices to prove that the actions $B:Z$ (or $B:\PP(Z)$) and
$B\cap G(\lambda)\cap G(\mu):V(\lambda,\mu)$ (resp.\ $B\cap
L(\lambda)\cap L(\mu):V(\lambda,\mu)$) have the same stabilizers
of general position. (These are Borel subgroups in the generic
stabilizers of double actions.) For computing these stabilizers,
we apply the algorithm at the end of~\ref{cotangent}.

We have $Z\subseteq V=V(\lambda^{*})\oplus V(\mu^{*})$. Take a
$B$-eigenvector $\omega=(v_{\lambda},0)\in V^{*}$ and put
$\oo{Z}=Z_{\omega}$. By Lemma~\ref{loc.str.V},
$\oo{Z}\iso\Ru{P(\lambda)}\times Z'$, where $Z'\iso\langle
v_{\lambda}\rangle\times C(\mu)$ as an $L(\lambda)$-variety. Now
take $\omega'=(0,v_{\mu})$ and put $\oo{Z}'=Z'_{\omega'}$. Then
$\oo{Z}'\iso[L(\lambda)\cap\Ru{P(\mu)}]\times Z''$, where
$Z''=\langle v_{\lambda}\rangle\times\langle v_{\mu}\rangle\times
V(\lambda,\mu)$ as an $L(\lambda)\cap L(\mu)$-variety. This proves
our claim on stabilizers of general position. \qedhere\end{roster}
\end{proof}

We shall denote by $c,r,\RG$ (resp.\
$\widehat{c},\widehat{r},\widehat\RG$) the complexity, rank and
the weight lattice of a $G$- (resp.~$\widehat{G}$-) action. Since
maximal unipotent subgroups of $G$ and $\widehat{G}$ coincide, it
follows from Proposition~\ref{c+r} that
%*
\begin{equation}
c(Z)+r(Z)=\widehat{c}(Z)+\widehat{r}(Z)
\end{equation}
%*
It is also clear that $\widehat{c}(Z)\leq c(Z)\leq
\widehat{c}(Z)$. Since an open subset $Gv_{-\lambda}\times
Gv_{-\mu}\subset Z$ is a $G$-equivariant principal
$(\kk^{\times})^2$-bundle over~$\PP(Z)$, and
$\widehat\RG(Z)\subseteq\Ch(\widehat{G})=\Ch(G)\oplus\ZZ^2$
projects onto $\ZZ^2$ with the kernel~$\RG(Z)$, we have
%*
\begin{align}
\widehat{c}(Z)&=c(\PP(Z)) \label{cc(d.c)}\\
\widehat{r}(Z)&=r(\PP(Z))+2 \label{rr(d.c)}
\end{align}
%*

Theorem~\ref{sgp.dbl.con}, together with Theorem~\ref{c&r<=dbl},
yields
\begin{theorem}
The following formulae are valid:
%*
\begin{align*}
2c(Z)+r(Z)&= 2+\dim G(\lambda)\backslash\!\!\backslash G\by G(\mu)
&
r(Z)&=\rk G-\rk G(\lambda,\mu) \\
2\widehat{c}(Z)+\widehat{r}(Z)&= 2+\dim
L(\lambda)\backslash\!\!\backslash G\by L(\mu) &
\widehat{r}(Z)&=\rk\widehat{G}-\rk\widehat{G}(\lambda,\mu)
\end{align*}
%*
\end{theorem}
For the proof, just note that $\dim C(\lambda)=\frac12(\dim G-
\dim G(\lambda)+1)=\frac12(\dim\widehat G-\dim L(\lambda))$ and
recall~\eqref{dim(G/G*G)}--\eqref{dim(G/L*L)}.
\begin{corollary}
The numbers $c,r,\widehat{c},\widehat{r}$ do not change if one
transposes $\lambda$ and $\mu$ or replaces $\lambda$ (or~$\mu$) by
the dual weight~$\lambda^{*}$ (resp.~$\mu^{*}$).
\end{corollary}
Indeed, the doubled $G$-variety $Z\times Z^{*}$ and the generic
stabilizers $G(\lambda,\mu)$, $L(\lambda,\mu)$ do not change.
\begin{corollary}
For $\mu=\lambda$ or~$\lambda^{*}$,
%*
\begin{align*}
c(Z)&=c(G/G(\lambda))+1 & r(Z)&=r(G/G(\lambda)) \\
\widehat{c}(Z)&=c(G/L(\lambda)) & \widehat{r}(Z)&=r(G/L(\lambda))
\end{align*}
%*
\end{corollary}
\begin{proof}
It follows from (the proof of) Lemma~\ref{sgp.dbl.con} that a
generic orbit of $G:Z$ has codimension~$1$ and is isomorphic to
$G/G(\lambda)$, and $G:\PP(Z)$ has an open orbit isomorphic
to~$G/L(\lambda)$. Now apply~\eqref{cc(d.c)}--\eqref{rr(d.c)}.
\end{proof}

Now we restrict our attention to factorial double cones. By
Proposition~\ref{HV}\ref{HV.fact}, $Z(\lambda,\mu)$ is factorial
iff $\lambda,\mu$ are fundamental weights. We shall write $C(i)$,
$Z(i,j)$, $P(i)$, \dots\ instead of $C(\omega_i)$,
$Z(\omega_i,\omega_j)$, $P(\omega_i)$, \dots. For all simple
groups $G$ and all pairs of fundamental weights
$\omega_i,\omega_j$, complexities and ranks of $Z(i,j)$ w.r.t.\
the $G$- and $\widehat{G}$-actions were computed
in~\cite{dbl.con}. All pairs of fundamental weights
$(\omega_i,\omega_j)$ such that $\widehat{c}(Z(i,j))=0,1$ are
listed, up to the transposition of $i,j$ and an automorphism of
the Dynkin diagram, in Tables~\ref{sph(d.c)}--\ref{c=1(d.c)}.

Suppose that $\kk[Z(i,j)]^U$ is minimally generated by
bihomogeneous eigenfunctions $f_1,\dots,f_r$ of weights
$\lambda_1,\dots,\lambda_r$ and bidegrees
$(n_1,m_1),\dots,(n_r,m_r)$. We may assume that $f_1,f_2$ have the
weights $\omega_i,\omega_j$ and bidegrees $(1,0)$, $(0,1)$. The
weights of other generators and their bidegrees are indicated in
the columns ``Weights'' and ``Degrees'', respectively. Here we
assume $\omega_i=0$ whenever $i\neq1,\dots,\rk G$.

We already noted that if $\widehat{c}(Z(i,j))=0$, then
$f_1,\dots,f_r$ are algebraically independent and
$(\lambda_1,n_1,m_1),\dots,(\lambda_r,n_r,m_r)$ are linearly
independent. If $\widehat{c}(Z(i,j))=1$, then $Z(i,j)\by U$ is a
hypersurface~\cite[6.5]{dbl.con} and the (unique) syzygy between
$f_1,\dots,f_r$ is of the form $P+Q+R=0$, where $P,Q,R$ are all
monomials in $f_1,\dots,f_r$ of the same weight $\lambda_0$ and
bidegree $(n_0,m_0)$ indicated in the column ``Syzygy'' of
Table~\ref{c=1(d.c)}.
\begin{question}
Why?
\end{question}

It follows from the classification that if $i=j$, then
$\widehat{c}(Z(i,j))=0$ and $m_i=n_i=1$ for $i=3,\dots,r$. Hence
$\widehat{T}$-eigenspaces in $\kk[Z(i,i)]^U$ are one-dimensional,
and the involution transposing the factors of $Z(i,i)=C(i)\times
C(i)$ multiplies each $\widehat{T}$-eigenfunction $f$ of bidegree
$(n,n)$ by $p(f)=\pm1$. We call $p(f)$ the \emph{parity} of~$f$.
The parities of generators are given in the column ``Parity'' of
Table~\ref{sph(d.c)}. If $f=f_1^{k_1}\dots f_r^{k_r}$ ($k_1=k_2$),
then $p(f)=p(k_3,\dots,k_r):=p(f_3)^{k_3}\dots p(f_r)^{k_r}$.

\begin{longtable}{|c|c|c|c|c|}
\caption{Spherical double cones (factorial case)\label{sph(d.c)}}\\
\hline
$G$ & Pair & Weights & Degrees & Parity \\
\hline
\endfirsthead
\caption{(continued)}\\
\hline
$G$ & Pair & Weights & Degrees & Parity \\
\hline
\endhead
\hline
\endfoot

$\AAa_l$ & $(\omega_i,\omega_j)$ & $\omega_{i-k}+\omega_{j+k}$ &
$(1,1)$  & $(-1)^k$ \\*
& $i\leq j$ & $k=1,\dots,\min(i,l+1-j)$ && for $i=j$ \\
\hline

$\BBb_l$ & $(\omega_1,\omega_1)$ & $0,\omega_2$ & $(1,1)$  &
$1,-1$ \\
\cline{2-5}

& $(\omega_1,\omega_j)$ & $\omega_{j-1},\omega_{j+1}$ & $(1,1)$ &
\\*
& $2\leq j\leq l-2$ & $\omega_j$ & $(2,1)$ & \\
\cline{2-5}

& $(\omega_1,\omega_{l-1})$ & $\omega_{l-2},2\omega_l$ & $(1,1)$ &
\\*
&& $\omega_{l-1}$ & $(2,1)$ & \\
\cline{2-5}

& $(\omega_1,\omega_l)$ & $\omega_l$ & $(1,1)$ & \\*
&& $\omega_{l-1}$ & $(1,2)$ & \\
\cline{2-5}

& $(\omega_l,\omega_l)$ & $\omega_k$ & $(1,1)$ & $(-1)^{k(k+1)/2}$
\\*
&& $k=0,\dots,l-1$ && \\
\hline

$\CCc_l$ & $(\omega_1,\omega_1)$ & $0,\omega_2$ & $(1,1)$  &
$-1$ \\
\cline{2-5}

& $(\omega_1,\omega_j)$ & $\omega_{j-1},\omega_{j+1}$ & $(1,1)$ &
\\*
& $2\leq j\leq l-1$ & $\omega_j$ & $(2,1)$ & \\
\cline{2-5}

& $(\omega_1,\omega_l)$ & $\omega_{l-1}$ & $(1,1)$ & \\*
&& $\omega_l$ & $(2,1)$ & \\
\cline{2-5}

& $(\omega_l,\omega_l)$ & $0,2\omega_k$ & $(1,1)$ &
$(-1)^l,(-1)^{l-k}$ \\*
&& $k=1,\dots,l-1$ && \\
\hline

$\DDd_l$ & $(\omega_1,\omega_1)$ & $0,\omega_2$ & $(1,1)$  &
$1,-1$ \\
\cline{2-5}

& $(\omega_1,\omega_j)$ & $\omega_{j-1},\omega_{j+1}$ & $(1,1)$ &
\\*
& $2\leq j\leq l-3$ & $\omega_j$ & $(2,1)$ & \\
\cline{2-5}

& $(\omega_1,\omega_{l-2})$ & $\omega_{l-3},\omega_{l-1}+\omega_l$
& $(1,1)$ & \\*
&& $\omega_{l-2}$ & $(2,1)$ & \\
\cline{2-5}

& $(\omega_1,\omega_{l-1})$ & $\omega_l$ & $(1,1)$ & \\
\cline{2-5}

& $(\omega_l,\omega_l)$ & $\omega_{l-2k}$ & $(1,1)$ & $(-1)^l$ \\*
&& $1\leq k\leq l/2$ && \\
\cline{2-5}

& $(\omega_{l-1},\omega_l)$ & $\omega_{l-2k-1}$ & $(1,1)$ & \\
&& $1\leq k\leq(l-1)/2$ && \\
\cline{2-5}

& $(\omega_2,\omega_{l-1})$ & $\omega_{l-1},\omega_1+\omega_l$ &
$(1,1)$ & \\*
&& $\omega_{l-2}$ & $(1,2)$ & \\
\hline

$\DDd_l$ & $(\omega_3,\omega_{l-1})$ &
$\omega_1+\omega_{l-1},\omega_2+\omega_l,\omega_l$ & $(1,1)$ & \\*
$l\geq6$ && $\omega_{l-3},\omega_1+\omega_{l-2}$ & $(1,2)$ & \\*
&& $\omega_2+\omega_{l-2}$ & $(2,2)$ & \\
\hline

$\DDd_5$ & $(\omega_3,\omega_4)$ &
$\omega_1+\omega_4,\omega_2+\omega_5,\omega_5$ & $(1,1)$ & \\*
&& $\omega_2,\omega_1+\omega_3$ & $(1,2)$ & \\
\hline

$\EEe_6$ & $(\omega_1,\omega_1)$ & $\omega_2,\omega_5$ & $(1,1)$
& $-1,1$ \\
\cline{2-5}

& $(\omega_1,\omega_2)$ & $\omega_1+\omega_5,\omega_3,\omega_6$ &
$(1,1)$ & \\*
&& $\omega_2+\omega_5,\omega_4$ & $(2,1)$ & \\
\cline{2-5}

& $(\omega_1,\omega_4)$ & $\omega_2,\omega_5,\omega_5+\omega_6$ &
$(1,1)$ & \\*
&& $\omega_3,\omega_6$ & $(2,1)$ & \\
\cline{2-5}

& $(\omega_1,\omega_5)$ & $0,\omega_6$ & $(1,1)$ & \\
\cline{2-5}

& $(\omega_1,\omega_6)$ & $\omega_1,\omega_4$ & $(1,1)$ & \\*
&& $\omega_2$ & $(2,1)$ & \\
\hline

$\EEe_7$ & $(\omega_1,\omega_1)$ & $0,\omega_2,\omega_6$ &
$(1,1)$ & $-1,-1,1$ \\
\cline{2-5}

& $(\omega_1,\omega_6)$ & $\omega_1,\omega_7$ & $(1,1)$ & \\*
&& $\omega_2$ & $(2,1)$ & \\
\cline{2-5}

& $(\omega_1,\omega_7)$ & $\omega_2,\omega_5,\omega_6$ & $(1,1)$ &
\\* && $\omega_3,\omega_7$ & $(2,1)$ & \\*
&& $\omega_4$ & $(2,2)$ & \\

\end{longtable}

\begin{table}[h]
\caption{Double cones of complexity one (factorial case)}
\label{c=1(d.c)}
\begin{center}
\makebox[0pt]{\begin{tabular}{|c|c|c|c|c|} \hline
$G$ & Pair & Weights & Degrees & Syzygy \\
\hline

$\BBb_l$ & $(\omega_2,\omega_l)$ & $\omega_1+\omega_l,\omega_l$
& $(1,1)$  & $\omega_1+\omega_{l-1}+\omega_l$ \\
$l\geq4$ && $\omega_{l-2},\omega_{l-1},\omega_1+\omega_{l-1}$ &
$(1,2)$ & $(2,3)$ \\
&& $\omega_1+\omega_{l-1}$ & $(2,2)$ & \\
\hline

$\BBb_3$ & $(\omega_2,\omega_3)$ & $\omega_1+\omega_3,\omega_3$
& $(1,1)$  & $\omega_1+\omega_2+\omega_3$ \\
&& $\omega_1,\omega_2,\omega_1+\omega_2$ & $(1,2)$ & $(2,3)$ \\
\hline

$\CCc_l$ & $(\omega_2,\omega_l)$ &
$\omega_{l-2},\omega_1+\omega_{l-1}$ & $(1,1)$  &
$2\omega_1+2\omega_{l-1}+\omega_l$ \\
$l\geq4$ && $\omega_1+\omega_{l-1},2\omega_1+\omega_l,\omega_l$
& $(2,1)$ & $(4,3)$ \\
&& $2\omega_{l-1}$ & $(2,2)$ & \\
\hline

$\CCc_3$ & $(\omega_2,\omega_3)$ & $\omega_1,\omega_1+\omega_2$
& $(1,1)$  & $2\omega_1+2\omega_2+\omega_3$ \\
&& $2\omega_1+\omega_3,\omega_3$ & $(2,1)$ & $(4,3)$ \\
&& $2\omega_2$ & $(2,2)$ & \\
\hline

$\DDd_6$ & $(\omega_4,\omega_5)$ &
$\omega_2+\omega_5,\omega_5,\omega_1+\omega_6,\omega_3+\omega_6$ &
$(1,1)$  & $\omega_2+\omega_4+\omega_5$ \\
&& $\omega_4,\omega_2+\omega_4,\omega_2,\omega_1+\omega_3$ &
$(1,2)$ & $(2,3)$ \\
\hline

$\EEe_7$ & $(\omega_1,\omega_2)$ &
$\omega_1,\omega_1+\omega_6,\omega_3,\omega_7$ & $(1,1)$ &
$\omega_1+\omega_2+\omega_6$ \\
&& $\omega_2+\omega_6,\omega_2,\omega_5,\omega_6$ & $(2,1)$ &
$(3,2)$ \\
\hline

\end{tabular}}
\end{center}
\end{table}

For spherical $Z(i,j)$, the algebra $\kk[Z(i,j)]^U$ was computed
by Littelmann~\cite{sph.dbl}. He observed that a simple structure
of $\kk[Z(i,j)]^U$ leads to the following decomposition rules:
\begin{Claim}[Tensor products]
%*
\begin{alignat*}{1}
V(n\omega_i)\otimes V(m\omega_j)&=
\bigoplus\limits_{k_1(n_1,m_1)+\dots+k_r(n_r,m_r)=(n,m)}
V(k_1\lambda_1+\dots+k_r\lambda_r)
\end{alignat*}
%*
\end{Claim}
\begin{Claim}[Symmetric and exterior squares]
%*
\begin{alignat*}{1}
\Sym^2V(n\omega_i)&=\bigoplus\limits_{\substack{
k_1+k_3+\dots+k_r=n \\
p(k_3,\dots,k_r)=1    }}
V(2k_1\omega_i+k_3\lambda_3+\dots+k_r\lambda_r) \\
\E^2V(n\omega_i)&=\bigoplus\limits_{\substack{
k_1+k_3+\dots+k_r=n \\
p(k_3,\dots,k_r)=-1   }}
V(2k_1\omega_i+k_3\lambda_3+\dots+k_r\lambda_r)
\end{alignat*}
%*
\end{Claim}
\begin{Claim}[Restriction]
%*
\begin{alignat*}{1}
\Res^G_{L(i)}V_G(m\omega_j)&=
\bigoplus\limits_{k_2m_2+\dots+k_rm_r=m}
V_{L(i)}(k_2(\lambda_2-n_2\omega_i)+\dots+
k_r(\lambda_r-n_r\omega_i))
\end{alignat*}
%*
\end{Claim}
\begin{proof}[Proofs]
The first two rules stem immediately from the above discussion.
Indeed, highest weight vectors in
$\kk[Z(i,j)]^U_{n,m}=V(n\omega_i)\otimes V(m\omega_j)$ are
proportional to monomials $f=f_1^{k_1}\dots f_r^{k_r}$ with
$k_1(n_1,m_1)+\dots+k_r(n_r,m_r)=(n,m)$. The transposition of the
factors of $Z(i,i)$ transposes the factors of
$\kk[Z(i,i)]^U_{n,n}=V(n\omega_i)^{\otimes2}$, and $f$ is
(skew)symmetric iff $p(f)=1$ (resp.~$-1$).

To prove the restriction rule, observe that
$Z(i,i)_{f_1}=C(i)_{f_1}\times
C(j)=P(i)\itimes{L(i)}(\kk^{\times}v_{-\omega_i}\times
C(j))=\Ru{P(i)}\times\kk^{\times}v_{-\omega_i}\times C(j)$. Hence
$\kk[Z(i,j)]^U_{f_1}\iso \kk[\kk^{\times}v_{-\omega_i}\times
C(j)]^{U\cap L(i)}\iso \kk[f_1,f_1^{-1}]\otimes\kk[C(j)]^{U\cap
L(i)}$, and $f_2,\dots,f_r$ restrict to a free system of
generators $\bar{f}_l(y)=f_l(v_{-\omega_i},y)$ of
$\kk[C(j)]^{U\cap L(i)}$. It remains to remark that
$\kk[C(j)]_n\iso V_G(n\omega_j)$, $U\cap L(i)$ is a maximal
unipotent subgroup of $L(i)$, and $\bar{f}_l$ have
$T$-eigenweights $\lambda_l-n_l\omega_i$:
%*
\begin{alignat*}{1}
t\bar{f}_l(y)=f_l(v_{-\omega_i},t^{-1}y)=
\omega_i(t)^{-n_l}f_l(t^{-1}v_{-\omega_i},t^{-1}y)=
\lambda_l(t)\omega_i(t)^{-n_l}\bar{f}_l(y)
\end{alignat*}
%*
\end{proof}

For the cases $\widehat{c}(Z(i,j))=1$, the algebra $\kk[Z(i,j)]^U$
was computed in \cite[6.5]{c&r}. A decomposition rule
\begin{question}
Restriction rule, c.i.~$Z(i,j)\by U$.
\end{question}
for tensor products is of the form:
%*
\begin{multline*}
V(n\omega_i)\otimes V(m\omega_j)=
\bigoplus\limits_{k_1(n_1,m_1)+\dots+k_r(n_r,m_r)=(n,m)}
V(k_1\lambda_1+\dots+k_r\lambda_r) \\
-\bigoplus\limits_{l_1(n_1,m_1)+\dots+l_r(n_r,m_r)=(n-n_0,m-m_0)}
V(\lambda_0+l_1\lambda_1+\dots+l_r\lambda_r)
\end{multline*}
%*
(Here ``${-}$'' is an operation in the Grothendieck group of
$G$-modules.)

\begin{example}
Suppose $G=\SL_d(\kk)$. Consider a double cone $Z(1,1)$.
\begin{question}
Not very interesting, perhaps. Maybe to consider $Z(i,i)$ or
exceptional~$G$, to illustrate the restriction rule?
\end{question}
We have $L(1)=\GL_{d-1}(\kk)$, $V(1,1)=\kk^{d-1}$, and $L(1,1)$
consists of matrices of the form
%*
\begin{equation*}
\begin{array}{|cc|@{\qquad}c@{\qquad}|}
\hline
t & 0 &                     \\
0 & t & \smash{\text{\huge{$0$}}} \\
\hline
\rule{0em}{1em} && \\
\multicolumn{2}{|c|@{\qquad}}{\smash{\text{\huge{$0$}}}}
      & \smash{\text{\huge{$*$}}}  \\
\rule{0em}{1em} && \\
\hline
\end{array}
\qquad(t\in\kk^{\times})
\end{equation*}
%*
Its subgroup $G(1,1)$ is defined by $t=1$. Hence $r=2$,
$\widehat{r}=3$, $c=1$, $\widehat{c}=0$, and
$\RG(Z(1,1))=\langle\eps_2-\eps_1\rangle=
\langle\omega_2-2\omega_1\rangle$,
$\widehat\RG(Z(1,1))=\RG(Z(1,1))+\langle(\omega_1,1,0),
(\omega_1,0,1)\rangle=
\langle(\omega_2,1,1),(\omega_1,1,0),(\omega_1,0,1)\rangle$. (Here
$\omega_1=\eps_1$, $\omega_2=\eps_1+\eps_2$.)

Since $V(\omega_1)^{\otimes2}=(\kk^d)^{\otimes2}=
\Sym^2\kk^d\oplus\E^2\kk^d=V(2\omega_1)\oplus V(\omega_2)$, the
algebra $\kk[Z(1,1)]^U$ contains eigenfunctions of the weights
$(\omega_1,1,0),(\omega_1,0,1),(\omega_2,1,1)$, and a function of
the weight $(\omega_2,1,1)$ has parity~$-1$. Clearly, these three
functions are algebraically independent (because their weights are
linearly independent) and compose a part of a minimal generating
system of $\kk[Z(1,1)]^U$. Since $\dim Z(1,1)\by U=3$, they
generate $\kk[Z(1,1)]^U$.

As a corollary, we obtain decomposition formulae:
%*
\begin{xxalignat}{1}
V(n\omega_1)\otimes V(m\omega_1)&= \bigoplus\limits_{0\leq
k\leq\min(n,m)}
V((n+m-2k)\omega_1+k\omega_2) \\
\Sym^2V(n\omega_1)&= \bigoplus\limits_{0\leq k\leq\min(n,m)/2}
V((n+m-4k)\omega_1+2k\omega_2) \\
\E^2V(n\omega_1)&= \bigoplus\limits_{0\leq k\leq\min(n-1,m-1)/2}
V((n+m-4k-2)\omega_1+(2k+1)\omega_2)
\end{xxalignat}
%*
For $d=2$, these are well-known Clebsch--Gordan formulae.
\end{example}

All $\widehat{G}$-spherical double cones $Z=Z(\lambda,\mu)$ were
recently classified by Stembridge~\cite{mult.free}.
\begin{question}
Include the results of Moikina for complexity~$1$?
\end{question}
By~\eqref{cc(d.c)}, $\widehat{c}(Z)$~depends only on the
parabolics $P(\lambda),P(\mu)$, i.e., on the supports of
$\lambda,\mu$ w.r.t.\ fundamental weights. (The \emph{support} of
$\lambda$ is the set of fundamental weights occurring in the
decomposition of $\lambda$ with nonzero coefficients). When the
support of $\lambda$ is reduced, $P(\lambda)$~increases, whence
$\widehat{c}(Z)$ may only decrease. Therefore it suffices to find
all pairs of maximal possible supports such that
$\widehat{c}(Z)=0$ for all simple groups.

The case of one-element supports, i.e., where $\lambda,\mu$ are
multiples of fundamental weights, is already covered by
Table~\ref{sph(d.c)}.  All remaining pairs of maximal supports, up
to the transposition and an automorphism of the Dynkin diagram,
are listed in Table~\ref{sph(d.fl)}.

Note that $Z(\lambda,\mu)$ is spherical iff $V(n\lambda)\otimes
V(m\mu)$ is multiplicity-free for $\forall n,m$
(see~\ref{sphericity}).

\begin{table}[h]
\caption{Spherical double cones (non-factorial case)}
\label{sph(d.fl)}
\begin{center}
\makebox[0pt]{\begin{tabular}{|c|c|c|c|c|c|c|c|c|c|c|} \hline $G$
& \multicolumn{4}{|c|}{$\AAa_l$} &
\multicolumn{5}{|c|}{$\DDd_l$} & $\EEe_6$ \\
\hline

$\lambda$ & $\omega_1$ & $\omega_2$ & $\omega_i$ & $\omega_i$ &
$\omega_1$ & $\omega_l$ & $\omega_l$ &
$\omega_l$ & $\omega_l$ & $\omega_1$ \\
\hline

$\mu$ & $\omega_1,\dots,\omega_l$ & $\omega_i,\omega_j$ &
$\omega_1,\omega_j$ & $\omega_j,\omega_{j+1}$ &
$\omega_i,\omega_l$ & $\omega_1,\omega_{l-1}$ &
$\omega_1,\omega_l$ & $\omega_{l-1},\omega_l$ &
$\omega_1,\omega_2$ &
$\omega_1,\omega_5$ \\
\hline

\end{tabular}}
\end{center}
\end{table}

\chapter{General theory of embeddings}
\label{LV-theory}

Equivariant embeddings of homogeneous spaces are one of the main
topics of this survey. The general theory of them was developed by
Luna and Vust in a fundamental paper~\cite{LV}. However it was
noticed in \cite{LVT} that the whole theory admits a natural
exposition in a more general framework, which is discussed in this
chapter. The generically transitive case differs from the general
one by the existence of a smallest $G$-variety of a given
birational type, namely, a homogeneous space.

In~\ref{LV-gen} we discuss the general approach of Luna and Vust
based on patching all $G$-varieties of a given birational class
together in one huge prevariety and studying particular
$G$-varieties as open subsets in it. An important notion of a
$B$-chart arising in such a local study is considered
in~\ref{B-charts}. A $B$-chart is a $B$-stable affine open subset
of a $G$-variety, and any $G$-variety is covered by (finitely
many) $G$-translates of $B$-charts. $B$-charts and their
``admissible'' collections corresponding to coverings of
$G$-varieties are described in terms of \emph{colored data}
arranged of $B$-stable divisors and $G$-invariant valuations of a
given function field. This leads to a ``combinatorial''
description of $G$-varieties in terms of colored data, obtained
in~\ref{G-germs}. In the cases of complexity $\leq1$, considered
in~\ref{c=0}--\ref{c=1}, this description is indeed combinatorial,
namely, in terms of polyhedral cones, their faces, fans and other
objects of combinatorial convex geometry.

Divisors on $G$-varieties are studied in~\ref{div}. We give
criteria for a divisor to be Cartier, finitely generated and
ample, describe global sections in terms of colored data. Aspects
of the intersection theory on a $G$-variety are discussed
in~\ref{intersect}, including the role of $B$-stable cycles and a
formula for the degree of an ample divisor.

\section{The Luna--Vust theory}
\label{LV-gen}

The fundamental problem of classifying algebraic varieties has an
equivariant analogue: to describe up to a $G$-isomorphism all
varieties equipped with an action of an algebraic group $G$. A
birational classification of $G$-varieties (with a given field of
$G$-invariant functions) may be obtained in terms of Galois
cohomology~\cite[\S2]{IT}. The second, ``biregular'', part of the
problem may be formulated as follows: to describe all $G$-actions
in a given birational class. More precisely, let $K$ be a fixed
function field (i.e., a finitely generated extension of~$\kk$),
and let $G$ act on $K$ birationally. In other words, $K$ is the
function field on some $G$-variety $X$. We say that $K$ is a
\emph{$G$-field} and $X$ is a \emph{$G$-model} of~$K$. The problem
is to classify all $G$-models of $K$ in terms involving certain
invariants of~$K$ itself (such as valuations etc).
\begin{remark}
If $K^G=\kk$, or equivalently, the $G$-action on each $G$-model of
$K$ is generically transitive, then there is a minimal $G$-model
$\HS=G/H$, which is embedded as a dense orbit in any other
$G$-model of~$K$. The homogeneous space $\HS$ determines and is
determined by $K$ completely. So the problem may be thought of as
classifying all $G$-equivariant embeddings of $\HS$ in terms of
invariants of $\HS$ itself.
\end{remark}

A general approach to this problem was introduced by Luna and Vust
\cite{LV}. They considered only embeddings of homogeneous spaces.
We will follow \cite{LVT} in our more general point of view.

All models of $K$ may be glued together into one huge scheme
$\XX=\XX(K)$. By definition, points of $\XX$ are local rings that
are localizations of finitely generated $\kk$-algebras with
quotient field~$K$. Any model $X$ of $K$ (i.e., a variety with
$\kk(X)=K$) may be considered as a subset of~$\XX$, and such
subsets define the base of the \emph{Zariski topology} on $\XX$.
The structure sheaves $\Oo_X$ are patched together in a structure
sheaf $\Oo_{\XX}$. A local ring $\Oo_{\XX,Y}$ of $Y\in\XX$ in the
sense of this sheaf is exactly the ring defining $Y$ as a point
of~$\XX$.

The scheme $\XX$ is irreducible, but neither Noetherian nor
separated. It can be considered as a prevariety if we consider
only closed points $x\in\XX$ (i.e., such that the residue field
$\kk(x)=\Oo_x/\m_x$ of the respective local ring
$\Oo_x=\Oo_{\XX,x}$ equals $k$). Non-closed schematic points are
identified with closed irreducible subvarieties $Y\subseteq\XX$.

We distinguish in $\XX$ open subsets $\XX^{\reg}$,
$\XX^{\norm}$,\dots\ of smooth, normal,\dots\ points.

From this point of view, a model of $K$ is nothing but a
Noetherian separated open subset $X\subseteq\XX$.

The birational $G$-action on $K$ permutes local subrings of~$K$,
which yields an action $G:\XX$. However this is not an action in
the category of schemes or prevarieties. But the action map
$\alpha:G\times\XX\to\XX$ is rational and induces an embedding of
function fields
$\alpha^{*}:K=\kk(\XX)\embeds\kk(G\times\XX)=\kk(G)\cdot K$. (Here
$\kk(G)\cdot K=\Quot(\kk(G)\otimes_kK)$ is a free composite of
fields.) It is obvious that $G$ acts on a $G$-stable open subset
$X\subseteq\XX$ regularly iff
$\alpha^{*}(\Oo_{X,x})\subseteq\Oo_{G\times X,\1\times x}$ for
$\forall x\in X$.

Denote by $\XX_G$ the set of those $x\in\XX$ whose local rings
$\Oo_{\XX,x}$ are mapped by $\alpha^{*}$ to
$\Oo_{G\times\XX,\1\times x}$.
\begin{proposition}\label{reg.act}
$\XX_G$ is the largest open subset of $\XX$ on which $G$ acts
regularly.
\end{proposition}
\begin{proof}
We have only to prove that $\XX_G$ is open. In other words, if
$\Oo_{\XX,x}\embeds\Oo_{G\times\XX,\1\times x}$ for $x=x_0$, then
the same thing holds in a neighborhood of~$x_0$. Let $X=\Spec A$
be an affine neighborhood of~$x_0$, where $A=\kk[f_1,\dots,f_s]$
is a finitely generated algebra with quotient field~$K$. Then
$\Oo_{\XX,x}$ is a localization of $A$ at the maximal ideal
of~$x_0$. By assumption, $\alpha^{*}(f_i)$ are defined in a
neighborhood $E$ of $\1\times x_0$ (one and the same for
$i=1,\dots,s$), hence $\alpha$ restricts to a regular map $E\to
X$. The set of all $x\in X$ such that $\1\times x\in E$ is a
neighborhood of~$x_0$. In this neighborhood, we have
$\alpha(\1\times x)=x$, because this holds generically on~$\XX$.
This yields the assertion.
\end{proof}

Observe that $\g$ acts on $K$ by derivations (along velocity
fields on $\XX_G$ or on any other $G$-model).
\begin{proposition}[{\cite[1.4]{LV}}]\label{g-stab}
In characteristic zero, $x\in\XX_G$ iff $\Oo_x$ is $\g$-stable.
\end{proposition}
In particular, if $A\subset K$ is a $\g$-stable finitely generated
subalgebra, then any localization of $A$ is $\g$-stable and
consequently $X=\Spec A\subseteq\XX_G$.

\begin{example}
Let $G=\kk$ act on~$\AAA^1$ by translations. This yields a
birational action $G:K=\kk(t)$, so that $\alpha^{*}(t)=u+t$
($u$~is a coordinate on~$G$). A cuspidal curve $X\subset\AAA^2$
(the Neil parabola) defined by the equation $y^2=x^3$ becomes a
model of~$K$ if we put $t=y/x$. The local ring of the singular
point $x_0=(0,0)\in X$ consists of rational functions in $x=t^2$,
$y=t^3$ whose denominators have nonzero constant term. But
$\alpha^{*}(t^d)=u^d+du^{d-1}t+\dots$ is not defined at $0\times
x_0\in G\times\XX$ (at least when $d$ does not divide $\ch{k}$),
because $t$ is not defined at~$x_0$. Therefore $x_0\notin\XX_G$.
All other points of $X$ are in~$\XX_G$, because they are
identified via the normalization map $t\mapsto(t^2,t^3)$ with the
respective points of $\AAA^1$, where $G$ acts regularly.

The standard basic vector $\xi\in\g=\kk$ acts on $K$ as~$d/dt$,
and $\xi(t^d)=dt^{d-1}\in\Oo_{x_0}$ if $d>2$. But $\xi
x=2t\notin\Oo_{x_0}$ if $\ch{k}\neq2$, in accordance with
Proposition~\ref{g-stab}. However in characteristic~$2$, the
algebra $\kk[X]=\kk[x,y]$ is $\g$-stable and
Proposition~\ref{g-stab} is not applicable.
\end{example}
\begin{example}
Another example of this kind is the birational action of $G=\kk^n$
by translations on a blow-up $X$ of~$\AAA^n$ at~$0$. All points in
the complement to the exceptional divisor are in~$\XX_G$, since
they come from~$\AAA^n$, where $G$ acts regularly. In a
neighborhood of a point $x_0$ on the exceptional divisor, $X$~is
defined by local equations $x_i=x_1y_i$ ($1<i\leq n$) in
$\AAA^n\times\AAA^{n-1}$. We have $\alpha^{*}(x_i)=x_i+u_i$, where
$u_i$ are coordinates on~$G$, and
$\alpha^{*}(y_i)=\frac{x_i+u_i}{x_1+u_1}=
\frac{x_1y_i+u_i}{x_1+u_1}$ are not defined at~$0\times x_0$.
Hence $x_0\notin\XX_G$. On the other hand, the standard basic
vectors $\xi_1,\dots,\xi_n$ of $\g=\kk^n$ act on
$K=\kk(x_1,\dots,x_n)$ as $\partial/\partial x_1,\dots,
\partial/\partial x_n$, and
%*
\begin{equation*}
\xi_iy_j=
\begin{cases}
1/x_1, & i=j>1 \\
-y_j/x_1, & i=1<j \\
0, & \text{otherwise,}
\end{cases}
\end{equation*}
%*
so that not all $\xi_iy_j$ are in~$\Oo_{x_0}$.
\end{example}
These two examples are typical in a sense that one obtains ``bad''
birational actions if one blows up or contracts $G$-nonstable
subvarieties in a variety with a ``good'' (regular) action.

By Proposition~\ref{reg.act}, a $G$-model of $K$ is nothing but a
$G$-stable Noetherian separated open subset of $\XX_G$. The next
theorem gives a way to construct $G$-models as ``$G$-spans'',
which we use in~\ref{B-charts}.
\begin{theorem}[{\cite[1.5]{LV}}]
Assume $\X$ is an open subset in $\XX_G$. Then $X=G\X$ is
Noetherian (separated) iff $\X$ is Noetherian (separated).
\end{theorem}
\begin{proof}
If $\X$ is Noetherian, then $G\times\X$ is Noetherian, and
$X=\alpha(G\times\X)$ is Noetherian, too. If $X$ is not separated,
then the diagonal $\Delta_X$ is not closed in $X\times X$, and the
non-empty $G$-stable subset $\overline{\Delta_X}\setminus\Delta_X$
contains an orbit~$\HS$. Clearly, $\HS$~intersects the two open
subsets $\X\times X$ and $X\times\X$ of $X\times X$. Since $\HS$
is irreducible, $\HS\cap(\X\times
X)\cap(X\times\X)=\HS\cap(\X\times\X)\subseteq
\overline{\Delta_{\X}}\setminus\Delta_{\X}$ is non-empty, whence
$\X$ is not separated.

The converse implications are obvious.
\end{proof}
\begin{example}
Let $G=\kk^{\times}$ act on $\AAA^1$ by homotheties. Here
$K=\kk(t)$ and a generator $\xi$ of $\g=\kk$ acts on $K$ as
$t\frac{d}{dt}$. Put $x=\frac{t}{(1+t)^2}$, $y=\frac{t}{(1+t)^3}$.
Then $t\mapsto(x,y)$ is a birational map of $\AAA^1$ to the
Cartesian leaf $\X\subset\AAA^2$ defined by the equation
$x^3=xy-y^2$. This map provides a biregular isomorphism of
$\AAA^1\setminus\{-1\}$ onto $\X\setminus\{x_0\}$, where
$x_0=(0,0)$ is the singular point of~$\X$. Therefore
$\X\setminus\{x_0\}\subseteq\XX_G$. One can verify by direct
computation that $\alpha^{*}(x),\alpha^{*}(y)\in\Oo_{1\times
x_0}$, whence $x_0\in\XX_G$. In characteristic zero, the situation
is simpler, because $\xi x=2y-x$, $\xi y=y-3x^2$ imply that the
algebra $A=\kk[x,y]$ is $\g$-stable, hence $\X=\Spec
A\subset\XX_G$. Put $X=G\X$. Then $G$ acts on $X$ with $2$~orbits
$X\setminus\{x_0\}$ and $\{x_0\}$ (an ordinary double point),
cf.~Example~\ref{8}.
\end{example}

In the study of local geometry of a variety $X$ in a neighborhood
of its (irreducible) subvariety~$Y$, we may replace $X$ by any
open subset intersecting~$Y$, thus arriving to the notion of a
\emph{germ} of a variety in (a neighborhood of) its subvariety. If
$X$ is a model of~$K$, then a germ of $X$ in~$Y$ is essentially
the local ring $\Oo_{X,Y}$ or the respective schematic point
of~$\XX$.
\begin{definition}
A \emph{$G$-germ} (of~$K$) is a $G$-fixed schematic point
of~$\XX_G$ (or a $G$-stable irreducible subvariety of~$\XX_G$).
The set of all $G$-germs is denoted by ${}_G\XX$; a similar
notation ${}_GX$ is used for an arbitrary open subset
$X\subseteq\XX_G$. A $G$-model $X$ such that a given $G$-germ is
contained in ${}_GX$ (i.e., intersects $X$ in a $G$-stable
subvariety~$Y$) is called a \emph{geometric realization} or a
\emph{model} of the $G$-germ.
\end{definition}
Every $G$-germ admits a geometric realization: just take its
affine neighborhood $\X\subseteq\XX_G$ and put $X=G\X$. If
$X\subseteq\XX_G$ is $G$-stable, then $X$ and ${}_GX$ determine
each other. The Zariski topology is induced on ${}_G\XX$, with
${}_GX$ the open subsets. It is straightforward~\cite[6.1]{LV}
that $X$ is Noetherian iff ${}_GX$ is Noetherian.
\begin{remark}
In characteristic zero, a germ of $X$ in $Y$ is a $G$-germ iff its
local ring $\Oo_{X,Y}$ is $G$- and $\g$-stable
(cf.~Proposition~\ref{g-stab}).
\end{remark}

Germs of normal $G$-models in $G$-stable prime divisors play an
important role in the Luna--Vust theory. They are identified with
the respective $G$-invariant valuations of~$K$.
\begin{definition}
A discrete $\QQ$-valued valuation of $K/\kk$ is called
\emph{geometric}
\begin{question}
Maybe arbitrary valuations will suit?
\end{question}
if it is a multiple of the valuation corresponding to a prime
divisor in a normal model of~$K$. A \emph{$G$-valuation} is a
$G$-invariant geometric valuation. Denote by $\Vv=\Vv(K)$ the set
of all $G$-valuations of~$K/\kk$; its structure is considered in
Chapter~\ref{inv.val}, see also \ref{B-charts}--\ref{c=1} below.

The \emph{support} $\Ss_Y$ of a $G$-germ $Y$ is the set of
$v\in\Vv$ such that the valuation ring $\Oo_v$
dominates~$\Oo_{\XX,Y}$ (i.e., $v$~has center $Y$ in any geometric
realization).
\end{definition}

The support of a $G$-germ is non-empty: e.g., if $X\supseteq Y$ is
an arbitrary geometric realization of the $G$-germ, and $v$ is the
valuation corresponding to a component of the exceptional divisor
in the normalized blow-up of $X$ along~$Y$, then $v\in\Ss_Y$.

Here is a version of the valuative criterion of separation.
\begin{theorem}\label{G-sep}
A $G$-stable open subset $X\subseteq\XX_G$ is separated iff the
supports of all its $G$-germs are disjoint.
\end{theorem}
\begin{proof}
The closure $\overline{\Delta_X}$ of the diagonal
$\Delta_X\subseteq X\times X$ is a $G$-model of~$K$. The
projections of $X\times X$ to the factors induce birational
regular $G$-maps $\pi_i:\overline{\Delta_X}\to X$ ($i=1,2$). If
$X$ is not separated and
$Y\subseteq\overline{\Delta_X}\setminus\Delta_X$ is a $G$-orbit,
then the orbits $Y_i=\pi_i(Y)$ are distinct for $i=1,2$. But
$\Ss_{Y_1}\cap\Ss_{Y_2}\supseteq\Ss_Y\neq\emptyset$, a
contradiction.

The converse implication follows from the valuative criterion of
separation.
\end{proof}

Assume that $K'\subseteq K$ is a subfield containing~$\kk$. We
have a natural dominant rational map
$\phi:\XX\dasharrow\XX'=\XX(K')$. If $X\subseteq\XX$,
$X'\subseteq\XX'$ are models of $K,K'$, then $\phi:X\to X'$ is
regular iff for any $x\in X$ there exists an $x'\in X'$ such that
$\Oo_x$ dominates $\Oo_{x'}$. This $x'$ is necessarily unique
(because $X'$ is separated), and $x'=\phi(x)$.

Now assume that $K'$ is a $G$-subfield of~$K$. Suppose that $X$
and $X'$ are $G$-models of $K$ and~$K'$.
\begin{proposition}
The natural rational map $\phi:X\to X'$ is regular iff for any
$G$-germ $Y\in{}_GX$ there exists a (necessarily unique) $G$-germ
$Y'\in{}_GX'$ such that $\Oo_{X,Y}$ dominates~$\Oo_{X',Y'}$.
\end{proposition}
\begin{proof}
If $\Oo_{X,Y}$ dominates~$\Oo_{X',Y'}$, then there exist finitely
generated subalgebras $A\supseteq A'$ such that $\Oo_{X,Y}$ and
$\Oo_{X',Y'}$ are their respective localizations. Localizing $A'$
and $A$ sufficiently, we may assume that $\X=\Spec A$ and
$\X'=\Spec A'$ are open subsets of $X$ and $X'$ intersecting $Y$
and $Y'$, respectively. The regular map $\X\to\X'$ extends to the
regular map $G\X\to G\X'$. Since $Y\subseteq X$ is arbitrary,
these maps paste together in the regular map $X\to X'$.

The converse implication is obvious: just put
$Y'=\overline{\phi(Y)}$.
\end{proof}

The restriction of a $G$-valuation of $K$ to $K'$ is a
$G$-valuation, and any $G$-valuation of $K'$ can be extended to a
$G$-valuation of~$K$ (Corollary~\ref{res&ext(G-val)}). Thus the
restriction map $\phi_{*}:\Vv(K)\to\Vv(K')$ is well defined and
surjective. If $\phi:X\to X'$ is a regular map, then
$\phi_{*}(\Ss_Y)\subseteq\Ss_{Y'}$ for any $G$-germ $Y\subseteq X$
and $Y'=\overline{\phi(Y)}$.

Here is a version of the valuative criterion of properness.
\begin{theorem}[{\cite[6.4]{LV}}]\label{G-prop}
A morphism $\phi:X\to X'$ is proper iff
\begin{question}
Proof.
\end{question}
%*
\begin{equation*}
\bigcup_{Y\subseteq X}\Ss_Y=
\phi_{*}^{-1}\left(\bigcup_{Y'\subseteq X'}\Ss_{Y'}\right)
\end{equation*}
%*
\end{theorem}
\begin{corollary}
$X$ is complete iff\/ $\bigcup_{Y\subseteq X}\Ss_Y=\Vv$ (i.e.,
each $G$-valuation has center on~$X$).
\end{corollary}

\section{$B$-charts}
\label{B-charts}

From now on, $G$~is assumed to be reductive and all $G$-models to
be normal, i.e., to lie in~$\XX_G^{\norm}$. We have seen
in~\ref{LV-gen} that a $G$-model $X$ is given by a Noetherian set
${}_GX$ of $G$-germs whose supports are disjoint. The Noether
property means that $X$ is covered by $G$-spans of finitely many
``simple'', e.g., affine, open subsets $\X\subseteq X$. An
important class of such ``local charts'' is introduced in
\begin{definition}
A \emph{$B$-chart} of $X$ is a $B$-stable affine open subset
of~$X$. Generally, a \emph{$B$-chart} is a $B$-stable affine open
subset $\X\subset\XX_G^{\norm}$.
\end{definition}

It follows from the local structure theorem (\ref{local}) that any
$G$-germ admits a $B$-chart $\X\subset X$ intersecting~$Y$.
Therefore $X$ is covered by finitely many $B$-charts and their
translates. Thus it is important to obtain a compact description
for $B$-charts. We describe their coordinate algebras in terms of
their $B$-stable divisors.
\begin{definition}
Denote by $\Dd=\Dd(K)$ the set of prime divisors on $X$ that are
not $G$-stable. The valuation corresponding to a divisor $D\in\Dd$
is denoted by~$v_D$. Prime divisors that are $B$-stable but not
$G$-stable, i.e., elements of $\Dd^B$, are called
\emph{$B$-divisors}.

Let $K_B\subseteq K$ be the subalgebra of rational functions with
$B$-stable divisor of poles on $X$.
\end{definition}
\begin{remark}
The sets $\Dd$, $\Dd^B$ and $K_B$ do not depend on the choice of a
$G$-model~$X$. Indeed, a non-$G$-stable prime divisor on
$\XX_G^{\norm}$ intersects any $G$-model
$X\subseteq\XX_G^{\norm}$, and $K_B$ consists of rational
functions defined on $\XX_G^{\norm}$ everywhere outside a
non-$G$-stable divisor.
\end{remark}
Since $K_B\supseteq\kk[\X]$ for any $B$-chart~$\X$, it follows
that $\Quot{K_B}=K$.

$B$-divisors are also called \emph{colors}, and the pair
$(\Vv,\Dd^B)$ is said to be the \emph{colored equipment} (of~$K$).
It is in terms of colored equipment that $B$-charts, $G$-germs and
$G$-models are described, as we shall see below.
\begin{remark}\label{D^B(G/H)}
In the generically transitive case, $\Dd^B$~may be computed as
follows. Each $D\in\Dd$ determines a $G$-line bundle
$\ind{\chi}=\ind{\kk_{\chi}}$ over $\HS=G/H$ and a section
$\eta\in\Ho^0(\HS,\ind{\chi})=\kk[G]^{(H)}_{-\chi}$ defined
uniquely up to multiplication by an invertible function on~$\HS$,
i.e., by a scalar multiple of a character of~$G$. The section
$\eta$ may be regarded as an equation for the preimage of $D$
under the orbit map $G\to\HS$. Since $D$ is prime, $\eta$~is
indecomposable in the multiplicative semigroup
$\kk[G]^{(H)}/\kk^{\times}\Ch(G)$.

Each $f\in\kk(\HS)$ decomposes as
$f=\eta^d\eta_1^{d_1}\dots\eta_s^{d_s}$, where
$d,d_1,\dots,d_s\in\ZZ$ and
$\eta,\eta_1,\dots,\eta_s\in\kk[G]^{(H)}$ are pairwise coprime.
Then $v_D(f)=d$.

Finally, $D$ is a $B$-divisor iff $\eta$ is a $(B\times
H)$-eigenfunction. Therefore $\Dd^B$ is in bijection with the set
of generators of $\kk[G]^{(B\times H)}/\kk^{\times}\Ch(G)$.
\end{remark}

The ``dual'' object is the multiplicative group $K^{(B)}$ of
rational $B$-eigen\-func\-tions. There is an exact sequence
%*
\begin{equation}
\label{K^(B)} 1\longrightarrow(K^B)^{\times}\longrightarrow
K^{(B)} \longrightarrow\RG\longrightarrow0
\end{equation}
%*
where $\RG=\RG(K)$ is the weight lattice (of any $G$-model)
of~$K$.

In the sequel, we frequenlty use Knop's approximation
Lemma~\ref{approx}, which is crucial for reducing various
questions to $B$-eigenfunctions. In particular, it implies that
$G$-valuations are determined uniquely by their restriction to
$K^{(B)}$ (Corollary~\ref{restr}).

Let $\X$ be a $B$-chart. Then $\Alg=\kk[\X]$ is an integrally
closed finitely generated algebra, in particular, it is a Krull
ring. Therefore
%*
\begin{equation*}
\Alg=\bigcap\Oo_{\X,D} \text{ (over all prime divisors
$D\subset\X$) }= \bigcap_{w\in\Ww}\Oo_w\cap\bigcap_{D\in
\widetilde{\Rr}}\Oo_{v_D},
\end{equation*}
%*
where $\Oo_v$ is the valuation ring of~$v$, $\Ww\subseteq \Vv$,
$\Rr\subseteq\Dd^B$, $\widetilde{\Rr}=\Rr
\sqcup(\Dd\setminus\Dd^B)$. Here the $G$-valuations $w\in\Ww$ are
determined up to a rational multiple, and we shall ignore this
indeterminacy, thus passing to a ``projectivization'' of~$\Vv$. In
particular, we may assume that the group of values of every
$w\in\Ww$ is exactly $\ZZ\subset\QQ$.

The pair $(\Ww,\Rr)$ is said to be the ``colored data'' of~$\X$. A
$B$-chart is uniquely determined by its colored data. Taking
another $B$-chart changes $\Ww$ and $\Rr$ by finitely many
elements. Hence all possible $\Ww\sqcup\Rr$ lie in a certain
distinguished class $\Cd$ of equivalent subsets of
$\Vv\sqcup\Dd^B$ w.r.t.\ the equivalence relation ``differ by
finitely many elements''.

Conversely, if $\Ww\subseteq\Vv$, $\Rr\subseteq\Dd^B$, and
$\Ww\sqcup\Rr\in\Cd$, then
%*
\begin{equation*}
\Alg=\Alg(\Ww,\Rr)=\bigcap_{w\in\Ww}\Oo_w\cap
\bigcap_{D\in\widetilde{\Rr}}\Oo_{v_D}
\end{equation*}
%*
is a Krull ring. Indeed, for $\forall f\in K$ almost all
valuations from $\Ww\sqcup\Rr$ (i.e., all but finitely many)
vanish on~$f$, since it is true for colored data of $B$-charts,
hence for any subset in the class~$\Cd$.
\begin{example}
$K_B=\Alg(\emptyset,\emptyset)$
\end{example}
\begin{remark}
Here and below, we identify prime divisors and respective
valuations. Thus, for $\Vv_0\subseteq \Ww\sqcup\Rr$, we write
$\langle\Vv_0,f\rangle\geq0$ iff $v(f)\geq0$ for $\forall
v\in\Vv_0$ ($v=v_D$ for $D\in\Vv_0\cap\Rr$), and so on.
\end{remark}

\begin{proposition}
\begin{roster}
\item\label{R} All valuations from $\widetilde{\Rr}$ are essential
for~$\Alg$. \item\label{W} A valuation $w\in\Ww$ is essential for
$\Alg$ iff
%*
\begin{equation}
\exists f\in K^{(B)}:\
\langle\Ww\sqcup\Rr\setminus\{w\},f\rangle\geq0,\ w(f)<0 \tag{W}
\end{equation}
%*
\end{roster}
\end{proposition}
\begin{proof}
\begin{roster}
\item[\ref{R}] Let $X$ be a smooth $G$-model of~$K$. Consider the
$G$-line bundle $\Ll=\Lin[X]{D}$, where $D\in\widetilde{\Rr}$, and
let $\eta\in\Ho^0(X,\Ll)$ be a section with $\divr\eta=D$. Put
$f=g\eta/\eta$, where $g\in G$, $gD\neq D$. Then $v_D(f)=-1$,
$\langle\widetilde\Rr\setminus\{D\},f\rangle\geq0$, and
$\langle\Ww,f\rangle=0$ by Corollary~\ref{v(g*sec/sec)}. Thus
$v_D$ is essential for~$\Alg$. \item[\ref{W}] Assume $w\in\Ww$ is
essential for $\Alg$; then $\exists f\in K:\
\langle\Ww\sqcup\widetilde\Rr\setminus\{w\},f\rangle\geq0$,
$w(f)<0$. Applying Lemma~\ref{approx}, we replace $f$ by a
$B$-eigenfunction and obtain~\reftag{W}. The converse implication
is obvious. \qedhere\end{roster}
\end{proof}

\begin{theorem}
\begin{roster}
\item\label{C} $\Quot{\Alg}=K$ iff
%*
\begin{equation}
\forall\Vv_0\subseteq\Ww\sqcup\Rr,\ \Vv_0\text{ finite},\ \exists
f\in K^{(B)}:\ \langle\Ww\sqcup\Rr,f\rangle\geq0,\
\langle\Vv_0,f\rangle>0 \tag{C}
\end{equation}
%*
\item\label{F} $\Alg$ is finitely generated iff
%*
\begin{equation}
\Alg^U=\kk\bigl[\,f\in K^{(B)}\bigm|
\langle\Ww\sqcup\Rr,f\rangle\geq0\,\bigr] \text{ is finitely
generated} \tag{F}
\end{equation}
%*
\item\label{B-chart} Under the equivalent conditions of
\ref{C}--\ref{F}, $\X=\Spec\Alg$ is a $B$-chart.
\end{roster}
\end{theorem}
\begin{proof}
\begin{roster}
\item[\ref{C}] Assume $\Quot{\Alg}=K$. We may assume
$\Vv_0=\{v\}$; then $\exists f\in\Alg\subseteq K_B:\ v(f)>0$.
Applying Lemma~\ref{approx}, we replace $f$ by an element of
$\Alg^{(B)}$ and obtain~\reftag{C}. Conversely, assume \reftag{C}
is true and $h\in K_B$. Then we take $\Vv_0=\{v\in\Ww\sqcup\Rr\mid
v(h)\neq0\}$ and, multiplying $h$ by $f^N$ for $N\gg0$ (killing
the poles), we fall into~$\Alg$. Hence $K_B\subseteq\Quot{\Alg}$,
and this yields $K=\Quot{\Alg}$.

\item[\ref{F}] ($\ch\kk=0$) Let $X$ be a smooth $G$-model of $K$.
Take an effective divisor on $X$ with support $\Dd^B\setminus\Rr$
and consider the corresponding section $\eta\in\Ho^0(X,\Ll)^{(B)}$
of the $G$-line bundle~$\Ll$. Consider an algebra
$R=\bigoplus_{n\geq0}R_n$, where
$R_n=\{\sigma\in\Ho^0(X,\Ll^n)\mid\sigma/\eta^n\in\Alg\}$. Then
$\Alg=\bigcup\eta^{-n}R_n\subseteq\Quot{R}$. Since every
$G$-valuation of $K$ can be extended to a $G$-valuation of
$\Quot{R}$ (Corollary~\ref{res&ext(G-val)}), we see that
$R_n=\{\sigma\mid\forall w\in\Ww:\ w(\sigma)\geq nw(\eta)\}$ is
$G$-stable.

Though $\Alg$ is not a $G$-algebra, it is very close to a
$G$-algebra, so that we may apply Lemma~\ref{G-qinv} below and
reduce the problem of finite generation to $G$-algebras.

In the notation of Lemma~\ref{G-qinv}, if $\Alg=\phi(R)$ is
finitely generated, then it is easy to construct a finitely
generated graded $G$-subalgebra $S\subseteq R$ such that
$\phi(S)=\Alg$. Hence $S^U$ and also $\Alg^U=\phi(S^U)$ are
finitely generated. Conversely, if $\Alg^U$ is finitely generated,
then we construct a finitely generated graded $G$-algebra
$S=\langle GS^U\rangle$ such that $\phi(S^U)=\Alg^U$. This yields
that $\Alg=\phi(S)$
\begin{question}
Modify the proof for $\ch\kk>0$.
\end{question}
is finitely generated. \item[\ref{B-chart}] In characteristic
zero, just note that $\Alg$ is $\g$-stable, because all $\Oo_v$
($v\in\Ww\sqcup\Rr$) are. In general, since $\Alg$ is finitely
generated, it follows that $\Alg=\kk[\eta^{-n}M]$ for a
finite-dimensional $G$-submodule $M$ in some~$R_n$. Let $X'$ be
the closure of the image of the natural rational map
$X\dasharrow\PP(M^{*})$. Then $X'$ is a $G$-model and
$\X=\Spec\Alg\subseteq X'$. \qedhere\end{roster}
\end{proof}
\begin{lemma}\label{G-qinv}
Suppose $R$ is a $\ZZ_{+}$-graded $G$-algebra without zero
divisors, $S\subseteq R$ is a $G$-stable graded subalgebra. Take
$\eta\in S_1^{(B)}$ and consider the homomorphism
$\phi:R\to\Quot{R}$, $\phi(\sigma)={\sigma}/{\eta^n}$ for
$\forall\sigma\in R_n$. Then $\phi(S)^U=\phi(S^U)$,
$\phi(R)^U=\phi(R^U)$, and in characteristic zero
$\phi(S)=\phi(R)\iff\phi(S^U)=\phi(R^U)$.
\end{lemma}
\begin{proof}
Suppose $f=\phi(\sum\sigma_n)\in\phi(S)^U,\ \sigma_n\in S_n$; then
$f=\sum{\sigma_n}/{\eta^n}=
\sum{\sigma_n\eta^{n_0-n}}/{\eta^{n_0}}$ (for $n_0$ sufficiently
large) ${}=\phi(\sum\sigma_n\eta^{n_0-n})$, and $\eta\in S^U$
implies ${\sum\sigma_n\eta^{n_0-n}\in S^U}$. The same argument
shows $\phi(R)^U=\phi(R^U)$.

Now assume $\ch\kk=0$ and $\phi(S)\neq\phi(R)$. The $U$-module
$\phi(R)/\phi(S)$ contains a nonzero $U$-invariant, which is the
image of $f={\sigma}/{\eta^n},\ \sigma\in R_n$. For $\forall u\in
U\ \exists k\geq0:\ \sigma-u\sigma\in\eta^{-k}S_{n+k}$. Since the
$U$-module generated by $\sigma$ is finite-dimensional, we may
choose $k$ independent of~$u$. Replacing $\sigma$ by
$\sigma\eta^k$, we may assume that $\sigma$ determines a nonzero
element of $(R_n/S_n)^U$. By complete reducibility of $G$-modules,
we may replace $\sigma$ by an element of~$R_n^U$, without changing
$\sigma\bmod S_n$. Then $f\in\phi(R^U)\setminus\phi(S^U)$, and we
are done.
\end{proof}
\begin{corollary}\label{CFW}
A pair $(\Ww,\Rr)$ from $\Cd$ is the colored data of a $B$-chart
iff conditions \reftag{C},\reftag{F},\reftag{W} are satisfied.
\end{corollary}

Note that elements of $\Dd^B\setminus\Rr$ are exactly the
irreducible components of $X\setminus\X$, where $X=G\X$. A
$B$-chart is $G$-stable iff $\Rr=\Dd^B$.
\begin{corollary}
Affine $G$-models are in bijection with colored data $(\Ww,\Dd^B)$
satisfying \reftag{C},\reftag{F},\reftag{W}.
\end{corollary}

\begin{remark}\label{non-geom}
In this section, we never use an \emph{apriori} assumption that
$G$-invariant valuations from $\Ww$ are geometric.
\end{remark}

The local structure of $B$-charts is well understood.

The subgroup $P=N_G(\X)$ is a parabolic containing~$B$. We have
$P=P[\Dd^B\setminus\Rr]=\bigcap_{D\in\Dd^B\setminus\Rr}P[D]$,
where $P[D]$ is the stabilizer of~$D$. In the generically
transitive case, if $\eta\in\kk[G]^{(B\times H)}_{\lambda,\chi}$
is an equation of~$D$, then $P[D]=P(\lambda)$ is the parabolic
associated with~$\lambda$.

Let $P=L\Ru{P}$ be the Levi decomposition ($L\supseteq T$).

In~\ref{div}, we prove that the divisor $X\setminus\X$ is ample on
$X=G\X$ (Corollary~\ref{GX^0}). Now Lemma~\ref{loc.str.V}
\begin{question}
Rewrite the section on the local structure.
\end{question}
implies the following
\begin{proposition}[\cite{C-div}]\label{loc.str.B}
\begin{roster}
\item The action $\Ru{P}:\X$ is proper and has a geometric
quotient $\X/\Ru{P}=\Spec\kk[\X]^{\Ru{P}}$. \item\label{fin.sur.B}
There exists a $T$-stable ($L$-stable if $\ch\kk=0$) closed affine
subvariety $Z\subseteq\X$ such that $\X=PZ$ and the natural maps
$\Ru{P}\times Z\to\X$, $Z\to\X/\Ru{P}$ are finite and surjective.
\item[\ref{fin.sur.B}$'$] In characteristic zero, the $P$-action
on $\X$ induces an isomorphism
%*
\begin{equation*}
\Ru{P}\times Z=P\itimes{L}Z\isoto\X
\end{equation*}
%*
\end{roster}
\end{proposition}

\section{Classification of $G$-models}
\label{G-germs}

We begin with a description of $G$-germs in terms of colored data.

Consider a $G$-germ $Y\in{}_G\XX^{\norm}$. Let
$\Vv_Y\subseteq\Vv$, $\Dd_Y\subseteq\Dd$ be the subsets
corresponding to all $B$-stable divisors on $\XX_G^{\norm}$
containing~$Y$. The pair $(\Vv_Y,\Dd^B_Y)$ is said to be the
\emph{colored data} of the $G$-germ.

If the $G$-germ intersects a $B$-chart $\X$, then $\Y=Y\cap\X$ is
the center of any $v\in\Ss_Y$, i.e., $v|_{k[\X]}\geq0$, and the
ideal $\Ideal{\Y}\normin\kk[\X]$ is given by $v>0$. Conversely, if
a $G$-valuation $v\in\Vv$ is non-negative on $\kk[\X]$, then it
determines a $G$-germ intersecting~$\X$. If $(\Ww,\Rr)$ is the
colored data of~$\X$, then $\Vv_Y\subseteq\Ww$,
$\Dd^B_Y\subseteq\Rr$.

\begin{proposition}[{\cite[3.8]{G-val}}]\label{germ}
\begin{roster}
\item\label{c.d.} A $G$-germ is uniquely determined by its colored
data. \item\label{S} A $G$-valuation $v$ is in $\Ss_Y$ iff
%*
\begin{equation}
\begin{split}
\forall f\in K^{(B)}:\ \langle\Vv_Y\sqcup\Dd^B_Y,f\rangle\geq0
\implies v(f)\geq0\\
\text{and if $>$ occurs in the l.h.s., then $v(f)>0$}
\end{split}
\tag{S}
\end{equation}
%*
\end{roster}
\end{proposition}
\begin{proof}
Choose a geometric realization $X\supseteq Y$ and a $B$-chart
$\X\subseteq X$ intersecting~$Y$.
\begin{roster}
\item[\ref{S}] Observe that for $f\in K^{(B)}$ we have
$f\in\Oo_{X,Y}^{(B)}\iff\langle\Vv_Y\sqcup\Dd^B_Y,f\rangle\geq0$
and $f\in\m_{X,Y}^{(B)}$ iff one of these inequalities is strict.

If $\Oo_v$ dominates $\Oo_{X,Y}$, then \reftag{S} is satisfied.
Conversely, if $\exists f\in\Oo_{X,Y}$ such that $v(f)<0$, then,
applying Lemma~\ref{approx}, we replace $f$ by a $B$-eigenfunction
and see that \reftag{S} fails. Therefore
$\Oo_v\supseteq\Oo_{X,Y}\supseteq\kk[\X]$ and $v$ has center
$Y'\supseteq Y$ on~$X$. If $Y'\neq Y$, then for $\forall
v'\in\Ss_Y$ there is $f\in K$ such that $v'(f)>0$, $v(f)=0$.
Replacing $f$ by a $B$-eigenfunction again, we obtain a
contradiction with~\reftag{S}. Thus $\Oo_v$ dominates~$\Oo_{X,Y}$.
\item[\ref{c.d.}] Since $\kk[\X]\subseteq
\Alg=\Alg(\Vv_Y,\Dd^B_Y)\subseteq\Oo_{X,Y}$, the local ring
$\Oo_{X,Y}$ is the localization of $\Alg$ in the ideal
$I_Y=\Alg\cap\m_{X,Y}$. Take any $v\in\Ss_Y$; then $I_Y$ is
defined in $\Alg$ by $v>0$. But $\Ss_Y$ is determined
by~$(\Vv_Y,\Dd^B_Y)$. \qedhere\end{roster}
\end{proof}

Now we describe $G$-germs in a given $B$-chart $\X=\Spec\Alg$,
$\Alg=\Alg(\Ww,\Rr)$.

\begin{theorem}\label{germs}
\begin{roster}
\item\label{V} $v\in\Vv$ has a center on $\X$ iff
%*
\begin{equation}
\forall f\in K^{(B)}:\ \langle\Ww\sqcup\Rr,f\rangle\geq0 \implies
v(f)\geq0 \tag{V}
\end{equation}
%*
\item\label{VD'} Assume $v\in\Ss_Y$. A $G$-valuation $w\in\Ww$
belongs to $\Vv_Y$ iff
%*
\begin{equation}
\forall f\in K^{(B)}:\ \langle\Ww\sqcup\Rr,f\rangle\geq0,\
v(f)=0\implies w(f)=0 \tag{V$'$}
\end{equation}
%*
Similarly, $D\in\Rr$ belongs to $\Dd^B_Y$ iff
%*
\begin{equation}
\forall f\in K^{(B)}:\ \langle\Ww\sqcup\Rr,f\rangle\geq0,
v(f)=0\implies v_D(f)=0 \tag{D$'$}
\end{equation}
%*
\end{roster}
\end{theorem}

\begin{proof}
\begin{roster}
\item[\ref{V}] $v$ has a center iff $v|_{\Alg}\geq0$. This clearly
implies~\reftag{V}. On the other hand, if $f\in \Alg$, $v(f)<0$,
then, applying Lemma~\ref{approx}, we replace $f$ by an element of
$\Alg^{(B)}$ see that \reftag{V} is false.

\item[\ref{VD'}] Assume $w\in\Ww$ (or $D\in\Rr$) belongs to
$\Vv_Y$ (or $\Dd^B_Y$); then every function $f\in\Alg$ not
vanishing on $\Y$ (i.e., $v(f)=0$) does not vanish on the
respective $B$-stable divisor of $\X$ as well (i.e., $w(f)=0$ or
$v_D(f)=0$). For $f\in\Alg^{(B)}$, we obtain \reftag{V$'$}
(or~\reftag{D$'$}).

Conversely, assume $w\notin\Vv_Y$ (or $D\notin\Dd^B_Y$); then
there exists $f\in\Alg$ vanishing on the respective $B$-stable
divisor of~$\X$ (i.e., $w(f)>0$ or $v_D(f)>0$) but not on~$\Y$
(i.e., $v(f)=0$). Applying Lemma~\ref{approx}, we replace $f$ by
an element of $\Alg^{(B)}$ and see that \reftag{V$'$}
(or~\reftag{D$'$}) is false. \qedhere\end{roster}
\end{proof}

Summing up, we can construct every $G$-model in the following way:
\begin{enumerate}
\item Take a finite collection of colored data
$(\Ww_{\alpha},\Rr_{\alpha})$ in $\Cd$ satisfying
\reftag{C},\reftag{F}. Decrease $\Ww_{\alpha}$ if necessary so as
to satisfy~\reftag{W}. These colored data determine finitely many
$B$-charts~$\X_{\alpha}$. \item\label{charts->germs} Compute from
$(\Ww_{\alpha},\Rr_{\alpha})$ via conditions
\reftag{V},\reftag{V$'$},\reftag{D$'$} the collection of colored
data $(\Vv_Y,\Dd^B_Y)$ of $G$-germs $Y$
intersecting~$\X_{\alpha}$. \item\label{germ->supp} Compute the
supports $\Ss_Y$ from $(\Vv_Y,\Dd^B_Y)$ using~\reftag{S}.
\end{enumerate}
The $G$-models $X_{\alpha}=G\X_{\alpha}$ may be glued together in
a $G$-model $X$ iff the supports $\Ss_Y$ obtained at
Step~\ref{germ->supp} are disjoint (Theorem~\ref{G-sep}). The
collection ${}_GX$ of $G$-germs is given by
Step~\ref{charts->germs} as the collection of their colored data,
which is called the \emph{colored data} of~$X$.
\begin{remark}\label{min.chart}
We notice in addition that the collection of covering $B$-charts
$\X_{\alpha}$ is of course not uniquely determined. Furthermore,
one $G$-germ may have a lot of different $B$-charts. For example,
in the notation of Theorem~\ref{germs}, we may consider a
principal open subset $\X_f=\{x\mid f(x)\ne0\}$ in~$\X$, where
$f\in\Alg^{(B)}$, $v(f)=0$ (to avoid cutting $\Y$ off), i.e., pass
from $\Alg$ to its localization~$\Alg_f$. This corresponds to
removing from $(\Ww\setminus\Vv_Y)\sqcup(\Rr\setminus\Dd^B_Y)$ a
finite set $\Ww_0\sqcup\Rr_0$ of those valuations that are
positive on~$f$.  By~\reftag{V$'$} and~\reftag{D$'$}, this set may
contain any finite number of elements from $(\Ww\setminus\Vv_Y)
\sqcup(\Rr\setminus\Dd^B_Y)$. In particular, if
$(\Ww\setminus\Vv_Y)\sqcup(\Rr\setminus\Dd^B_Y)$ is finite, then
there exists a \emph{minimal} $B$-chart with $\Ww=\Vv_Y$,
$\Rr=\Dd^B_Y$.
\end{remark}

Parabolic induction does not change $K^{(B)}$ and~$\Vv$, while
$\Dd^B$ is extended by finitely many colors, whose valuations
vanish on~$K^B$, see Proposition~\ref{par.ind(hyp)}. The $G$-germs
of an induced variety are induced from those of the original
variety, and it is easy to prove the following result:
\begin{proposition}\label{par.ind(c.d.)}
Parabolic induction does not change the colored data of a
$G$-model.
\end{proposition}

\section{Case of complexity~0}
\label{c=0}

A practical use of the theory developed in the preceding sections
depends on whether the colored equipment of a $G$-field is
accessible for computation and operation or not. It was noted
already in~\cite{LV} that there is no hope to obtain a transparent
classification of $G$-models from the general description
in~\ref{G-germs} (maybe except particular examples) if the
complexity is $>1$. On the other hand, if the complexity is
$\leq1$, then an explicit solution to the classification problem
is obtained. An appropriate language to operate with the colored
equipment is that of convex polyhedral geometry.

We shall write $c(K)$, $r(K)$ for the complexity, resp.\ rank, of
(any $G$-model of)~$K$. If $c(K)=0$, then any $G$-model contains
an open $B$-orbit, hence an open $G$-orbit~$\HS$. Homogeneous
spaces of complexity zero (=spherical spaces) and their embeddings
are studied in details in Chapter~\ref{spherical}. Here we
classify the embeddings of a given spherical homogeneous space in
the framework of the Luna--Vust theory. This classification was
first obtained by Luna and Vust \cite[8.10]{LV}. For a modern
self-contained exposition, see \cite{LV-spher},
\cite[\S3]{spher.full}.

Let $\HS=G/H$ be a spherical homogeneous space and $K=\kk(\HS)$.
Since $K^B=\kk$, it follows from Corollary~\ref{restr} and the
exact sequence~\eqref{K^(B)} that $G$-valuations are identified by
restriction to $K^{(B)}$ with $\QQ$-linear functionals on the
lattice $\RG=\RG(\HS)$. The set $\Vv$ is a convex solid polyhedral
cone in $\ES=\Hom(\RG,\QQ)$, which is cosimplicial in
characteristic zero (see Chapter~\ref{inv.val}). The set $\Dd^B$
consists of irreducible components of the complement to the dense
$B$-orbit in~$\HS$, hence is finite. The restriction to $K^{(B)}$
yields a map $\res:\Dd^B\to\ES$, which is in general not
injective.
\begin{remark}\label{D^B(sph)}
If $\Dd^B=\{D_1,\dots,D_s\}$ and
$\eta_1,\dots,\eta_s\in\kk[G]^{(B\times H)}$ are the respective
indecomposable elements of biweights
$(\lambda_1,\chi_1),\dots,(\lambda_s,\chi_s)$, then
$(\lambda_i,\chi_i)$ are linearly independent. (Otherwise, there
is a linear dependence $\sum d_i(\lambda_i,\chi_i)=0$, and
$f=\eta_1^{d_1}\dots\eta_s^{d_s}$ is a non-constant $B$-invariant
rational function on~$\HS$.) If
$f=\eta_1^{d_1}\dots\eta_s^{d_s}\in K^{(B)}_{\lambda}$, then $\sum
d_i\lambda_i=\lambda$, $\sum d_i\chi_i=0$, and
$\langle\res(D_i),f\rangle=v_{D_i}(f)=d_i$.
\end{remark}
\begin{definition}
The space $\ES$ equipped with the cone $\Vv\subseteq\ES$ and with
the map $\res:\Dd^B\to\ES$ is the \emph{colored space} (of~$\HS$).
\end{definition}

Now we consider the structure of colored data and reorganize them
in a more convenient way. The proofs are straightforward, as soon
as we interpret $B$-eigenfunctions as linear functionals on~$\ES$.

The class $\Cd$ consists of finite sets.

Condition \reftag{C} means that $\Ww\sqcup\res(\Rr)$ generates a
strictly convex cone $\Cc=\Cc(\Ww,\Rr)$ in $\ES$ and
$\res(\Rr)\not\ni0$.

Condition \reftag{W} means that the elements of $\Ww$ are exactly
the generators of those edges of $\Cc$ that do not
intersect~$\res(\Rr)$.

Condition \reftag{F} holds automatically: $\Alg^U$ is the
semigroup algebra of $\RG\cap\Cc^{\vee}$, where $\Cc^{\vee}=
\{\lambda\in\ES^{*}\mid\langle\Cc,\lambda\rangle\geq0\}$ is the
dual cone to~$\Cc$. Since $\Cc^{\vee}$ is finitely generated, the
semigroup $\RG\cap\Cc^{\vee}$ is finitely generated by Gordan's
lemma~\cite[1.3]{toric.survey}.

Condition \reftag{V} means that $v\in\Cc$.

Conditions \reftag{V$'$} and \reftag{D$'$} say that $\Vv_Y$ and
$\Dd^B_Y$ consist of those elements of $\Ww\sqcup\Rr$ which lie in
the face $\Cc_Y=\Cc(\Vv_Y,\Dd^B_Y)\subseteq\Cc$ containing $v$ in
its (relative) interior.

Condition \reftag{S} means that $v\in\Vv\cap\intr\Cc_Y$.

Observe that every $G$-germ $Y$ has a minimal $B$-chart $\X_Y$
with $\Cc=\Cc_Y$, $\Rr=\Dd^B_Y$ (Remark~\ref{min.chart}). It
suffices to consider only such charts.

\begin{definition}
A \emph{colored cone} in $\ES$ is a pair $(\Cc,\Rr)$, where
$\Rr\subseteq\Dd^B$, $\res(\Rr)\not\ni0$, and $\Cc$ is a strictly
convex cone generated by $\res(\Rr)$ and finitely many vectors
from~$\Vv$.

A colored cone $(\Cc,\Rr)$ is \emph{supported} if
$\intr\Cc\cap\Vv\neq\emptyset$.

A \emph{face} of $(\Cc,\Rr)$ is a colored cone $(\Cc',\Rr')$,
where $\Cc'$ is a face of~$\Cc$ and $\Rr'=\Rr\cap\res^{-1}(\Cc')$.

A \emph{colored fan} is a finite set of supported colored cones
which is closed under passing to a supported face and such that
different cones intersect in faces inside~$\Vv$.
\end{definition}
\begin{theorem}\label{sph.emb}
\begin{roster}
\item\label{charts(sph)} $B$-charts are in bijection with colored
cones in~$\ES$. \item\label{germs(sph)} $G$-germs are in bijection
with supported colored cones. \item\label{col.fans} $G$-models are
in bijection with colored fans. \item\label{orb(sph)} Every
$G$-model $X$ contains finitely many $G$-orbits. If
$Y_1,Y_2\subseteq X$ are two $G$-orbits, then $Y_1\adhereseq Y_2$
iff $(\Cc_{Y_2},\Dd^B_{Y_2})$ is a face of
$(\Cc_{Y_1},\Dd^B_{Y_1})$.
\end{roster}
\end{theorem}
\begin{corollary}\label{sph.aff}
Affine $G$-models are in bijection with colored cones of the form
$(\Cc,\Dd^B)$.
\end{corollary}
\begin{corollary}\label{sph.qaff}
$\HS$~is (quasi)affine iff $\res(\Dd^B)$ can be separated from
$\Vv$ by a hyperplane (resp.\ does not contain~$0$ and spans a
strictly convex cone).
\end{corollary}
\begin{corollary}\label{sph.comp}
A $G$-model is complete iff its colored fan covers all~$\Vv$.
\end{corollary}
\begin{example}[Toric varieties]\label{toric}
Suppose $G=B=T$ is a torus. We may assume $H=\{\1\}$. Here
$\Vv=\ES$ (see~\ref{val.cones}) and there are no colors. Hence
embeddings of $T$ are in bijection with fans in~$\ES$, where a
\emph{fan} is a finite set of strictly convex polyhedral cones
which is closed under passing to a face and such that different
cones intersect in faces. Every embedding $X$ of $T$ contains
finitely many $T$-orbits, which correspond to cones in the fan.
For any orbit $Y\subseteq X$, the union $X_Y$ of all orbits
containing $Y$ in their closure is the minimal $T$-chart of $Y$
determined by~$\Cc_Y$. We have
$\kk[X_Y]=\kk[\Ch(T)\cap\Cc_Y^{\vee}]\subseteq\kk[T]$. $X$~is
complete iff its fan is the subdivision of the whole~$\ES$.

Equivariant embeddings of a torus are called \emph{toric
varieties}. Due to their nice combinatorial description, toric
varieties are a good testing site for various concepts and
problems of algebraic geometry. Their theory is well developed,
see \cite{toric.survey}, \cite{toric.book}, \cite{toric.intro}.
\end{example}

Other examples of spherical varieties are considered in
Chapter~\ref{spherical}.
\begin{question}
Non-normal (locally linear) spherical varieties: the normalization
map is bijective on $G$-orbits.
\end{question}

Now we discuss the functoriality of colored data.

Let $\overline{H}\subseteq G$ be an overgroup of~$H$. Denote by
$(\overline{\ES},\overline{\Vv},\overline{\Dd}^B,\overline{\res})$
the colored space of $\overline{\HS}=G/\overline{H}$. The
canonical map $\phi:\HS\to\overline{\HS}$ induces an embedding
$\phi^{*}:\overline{K}\embeds K$ and a linear map
$\phi_{*}:\ES\onto\overline{\ES}$. We have
$\phi_{*}(\Vv)=\overline{\Vv}$. If $\Dd^B_{\phi}$ is the set of
$B$-divisors in $\HS$ mapping dominantly to $\overline{\HS}$, then
there is a canonical surjection
$\phi_{*}:\Dd^B\setminus\Dd^B_{\phi}\onto\overline{\Dd}^B$ such
that $\overline{\res}=\res\phi_{*}$.

\begin{definition}
A colored cone $(\Cc,\Rr)$ in $\ES$ \emph{dominates} a colored
cone $(\overline{\Cc},\overline{\Rr})$ in $\overline{\ES}$ if
$\phi_{*}(\intr\Cc)\subseteq\intr\overline{\Cc}$ and
$\phi_{*}(\Rr\setminus\Dd^B_{\phi})\subseteq\overline{\Rr}$. A
colored fan $\Ff$ in $\ES$ \emph{dominates} a colored fan
$\overline{\Ff}$ in $\overline{\ES}$ if each cone from $\Ff$
dominates a cone from~$\overline{\Ff}$.

The \emph{support} of $\Ff$ is
$\Supp\Ff=\bigcup_{(\Cc,\Rr)\in\Ff}\Cc\cap\Vv$. (Observe that
$\{\Cc\cap\Vv\}$ is a polyhedral subdivision of~$\Supp\Ff$.)
\end{definition}

The next theorem is deduced from the results of~\ref{LV-gen}.
\begin{theorem}[{\cite[4.1--4.2]{LV-spher}}]\label{sph->sph}
Let $X,\overline{X}$ be the embeddings of $\HS,\overline{\HS}$
determined by fans $\Ff,\overline{\Ff}$. Then $\phi$ extends to a
morphism $X\to\overline{X}$ iff $\Ff$ dominates~$\overline{\Ff}$.
Furthermore, $\phi:X\to\overline{X}$ is proper iff
$\Supp\Ff=\phi_{*}^{-1}(\Supp\overline{\Ff})\cap\Vv$.
\end{theorem}
\begin{proof}
If $\Oo_{X,Y}$ dominates $\Oo_{\overline{X},\overline{Y}}$, then
clearly $\phi_{*}(\Dd^B_Y\setminus\Dd^B_{\phi})\subseteq
\overline{\Dd}^B_{\overline{Y}}$ and
$\phi_{*}(\Ss_Y)\subseteq\Ss_{\overline{Y}}$, or equivalently,
$(\Cc_Y,\Dd^B_Y)$ dominates
$(\Cc_{\overline{Y}},\overline{\Dd}^B_{\overline{Y}})$.
Conversely, if $(\Cc_Y,\Dd^B_Y)$ dominates
$(\Cc_{\overline{Y}},\overline{\Dd}^B_{\overline{Y}})$ for some
$Y\subseteq X$, $\overline{Y}\subseteq\overline{X}$, then
$\Alg=\Alg(\Vv_Y,\Dd^B_Y)\supseteq\overline{\Alg}=
\Alg(\overline{\Vv}_{\overline{Y}},\overline{\Dd}^B_{\overline{Y}})$
and $I_{\overline{Y}}=I_Y\cap\overline{\Alg}$, where
$I_Y=\Alg\cap\m_{X,Y}$ is defined in $\Alg$ by $v>0,\ \forall
v\in\Ss_Y$. Hence $\Oo_{X,Y}$ dominates
$\Oo_{\overline{X},\overline{Y}}$.

A criterion of properness is a reformulation of
Theorem~\ref{G-prop}.
\end{proof}

Overgroups of $H$ can be classified in terms of the colored space.
\begin{definition}
A \emph{colored subspace} of $\ES$ is a pair $(\ES_0,\Rr_0)$,
where $\Rr_0\subseteq\Dd^B$ and $\ES_0\subseteq\ES$ is a subspace
generated as a cone by $\res(\Rr_0)$ and some vectors from~$\Vv$.
\end{definition}
For example, $(\ES_{\phi},\Dd^B_{\phi})$ is a colored subspace,
where $\ES_{\phi}=\Ker\phi_{*}$ \cite[3.4]{spher.full}.
\begin{theorem}[{\cite[4.4]{LV-spher}}]
The correspondence $\overline{H}\mapsto(\ES_{\phi},\Dd^B_{\phi})$
is an order-preserving bijection between overgroups of $H$ with
$\overline{H}/H$ connected and colored subspaces of~$\ES$.
\end{theorem}
\begin{example}
If $H=B$, then $\ES=\Vv=0$ and $\Dd^B$ is the set of Schubert
divisors on $G/B$, which are in bijection with the simple roots.
Hence an overgroup of $B$ is determined by a subset of simple
roots---a well-known classification of parabolics.

More generally, parabolic overgroups $P\supseteq H$ are in
bijection with subsets $\Rr_0\subseteq\Dd^B$ such that
$\res(\Rr_0)\cup\Vv$ generates $\ES$ as a cone. Indeed, $P$ is
parabolic $\iff r(G/P)=0\iff\ES(G/P)=0\iff\ES_{\phi}=\ES$.
\end{example}

One may consider \emph{generalized} colored fans, dropping the
assumption that colored cones are strictly convex and their colors
do not map to~$0$. (These are exactly the preimages of usual
colored fans in quotients by colored subspaces.) Then there is a
bijection between dominant separable $G$-maps $\HS\to X$ to normal
$G$-varieties and generalized colored fans \cite[4.5]{LV-spher}.

Now we derive some properties of $G$-orbits (due to Brion) and
local geometry of a spherical embedding.

\begin{proposition}\label{cr(sph.orb)}
Suppose $X$ is an embedding of $\HS$ and $Y\subseteq X$ an
irreducible $G$-subvariety. Then $c(Y)=0$,
$r(Y)=r(X)-\dim\Cc_Y=\codim\Cc_Y$, and
$\RG(Y)=\Cc_Y^{\ann}\cap\RG(X)$ up to $p$-torsion.
\end{proposition}
\begin{proof}
By Theorem~\ref{cr(Y<X)}, $Y$~is spherical. By
Lemma~\ref{B-eigenfun}, for $\forall f\in\kk(Y)^{(B)}$ there is
$\widetilde{f}\in\kk(X)^{(B)}$ such that $\widetilde{f}|_Y=f^q$,
where $q$ is a sufficiently big power of~$p$. It remains to note
that $\widetilde{f}$ is defined and nonzero on~$Y$ iff
$\langle\Vv_Y\cap\Dd^B_Y,\widetilde{f}\rangle=0$, i.e., the
$B$-eigenweight of $\widetilde{f}$ lies in~$\Cc_Y^{\ann}$.
\end{proof}

The local structure of $B$-charts is given by
Proposition~\ref{loc.str.B}. For a minimal $B$-chart $\X_Y$, the
description can be refined.

Let $P=P[\Dd^B\setminus\Dd^B_Y]$ and $P=L\Ru{P}$ be its Levi
decomposition. Theorem~\ref{loc.str} yields
\begin{question}
Local structure of (toroidal) varieties in large (removing
$D\in\Dd^B\setminus\bigcup\Dd^B_Y$) \cite[2.4]{spher.full}
\end{question}
\begin{theorem}\label{loc.str.sph}
There is a $T$-stable ($L$-stable if $\ch\kk=0$) closed subvariety
$Z\subseteq\X_Y$ such that:
\begin{enumerate}
\item The natural maps $\Ru{P}\times Z\to\X_Y$ and $Z\to\X/\Ru{P}$
are finite and surjective. \item Put $\Y=Y\cap\X_Y$. Then
$\Y/\Ru{P}\iso L/L_0$, where $L_0\supseteq L'$. \item In
characteristic zero, $\X_Y\iso P\itimes{L}Z=\Ru{P}\times Z$,
$Y\cap Z\iso L/L_0$, and there exists an $L_0$-stable subvariety
$Z_0\subseteq Z$ transversal to $Y\cap Z$ at a fixed point~$y_0$
and such that $Z=L\itimes{L_0}Z_0$. The varieties $Z$ and $Z_0$
are affine and spherical, and $r(Z)=r(\HS)$, $r(Z_0)=\dim\Cc_Y$.
\end{enumerate}
\end{theorem}
The isomorphism $Z\iso L\itimes{L_0}Z_0$ stems, e.g., from Luna's
slice theorem, or is proved directly: since $L/L_0$ is a torus,
$\kk[Y\cap Z]$~is pulled back to an $L$-stable subalgebra
of~$\kk[Z]$, whence an equivariant retraction $Z\to Y\cap Z=L/L_0$
with a fiber~$Z_0$.
\begin{corollary}
$\dim Y=\codim\Cc_Y+\dim\Ru{P}$
\end{corollary}
\begin{remark}\label{c.d.(slice)}
In characteristic zero, there is a bijection $f\biject f|_Z$
between $B$-eigenfunctions on $X$ and~$Z$, which preserve the
order along a divisor. Hence $\Cc_Y=\Cc_{Y\cap Z}$ and
$\Dd^B_Y\supseteq\Dd^B_{Y\cap Z}$. However some colors on $X$ may
become $L$-stable divisors on~$Z$ (``a discoloration'').
\end{remark}

\begin{theorem}\label{norm&rat}
In characteristic zero, all $G$-subvarieties $Y\subseteq X$ are
normal and have rational singularities (in particular, they are
Cohen--Macaulay).
\end{theorem}
\begin{proof}
By the local structure theorem, we may assume that $X$ is affine.
Then $X\by U$ is an affine toric variety and $Y\by U$ its
$T$-stable subvariety. It is well known \cite[3.1,
3.5]{toric.intro} that $Y\by U$ is a normal toric variety and has
rational (in fact, Abelian quotient)
\begin{question}
This means quotient of $\AAA^n$ by a finite group? This holds only
for simplicial fans.
\end{question}
singularities. By Theorem~\ref{U-inv}\ref{X<=>X/U}, the same is
true for~$Y$.
\end{proof}

A spherical embedding defined by a fan whose colored cones have no
colors is called \emph{toroidal}. In particular, toric varieties
are toroidal. Conversely, the local structure theorem readily
implies that toroidal varieties are ``locally toric''
(Theorem~\ref{loc.str.tor}). This is the reason for most nice
geometric properties which distinguish toroidal varieties among
arbitrary spherical varieties. Toroidal varieties are discussed
in~\ref{toroidal}.

\section{Case of complexity~1}
\label{c=1}

Here we obtain the classification of $G$-models in the case
$c(K)=1$. This case splits in two subcases:
\begin{enumerate}
\item\emph{Generically transitive case}: $d_G(K)=0$. Here any
$G$-model contains a dense $G$-orbit $\HS$ of complexity~$1$.
\item\emph{One-parametric case}: $d_G(K)=1$. Here generic
$G$-orbits in any $G$-model are spherical and form a one-parameter
family. (In fact, all $G$-orbits are spherical by
Proposition~\ref{c(Gx)}.)
\end{enumerate}
We are interested mainly in the generically transitive case.
However the one-parameter case might be of interest, e.g., in
studying deformations of spherical homogeneous spaces and their
embeddings. There are differences between these two cases (e.g.,
in the description of $B$-divisors), but the description of
$G$-models is uniform~\cite{LVT}.

First we describe the colored environment.

Since $c(K)=1$, there is a (unique) non-singular projective curve
$C$ such that $K^B=\kk(C)$. Generic $B$-orbits on a $G$-model of
$K$ are parametrized by an open subset of~$C$. In the generically
transitive case, $K\subseteq\kk(G)$ is unirational, because $G$ is
a rational variety, which is proved by considering the ``big
cell'' in~$G$. Whence $C=\PP^1$ by the L\"uroth theorem.

\begin{definition}
For any $x\in C$ consider the half-space
$\EHS{x}=\QQ_{+}\times\ES$. The \emph{hyperspace} (of~$K$) is the
union $\Hyp$ of all $\EHS{x}$ glued together along their common
boundary hyperplane~$\ES$, called the \emph{center} of~$\Hyp$.
More formally,
%*
\begin{equation*}
\Hyp=\bigsqcup_{x\in C}\{x\}\times\EHS{x}\biggm/{\sim}
\end{equation*}
%*
where $(x,h,\ell)\sim(x',h',\ell')$ iff $x=x',\ h=h',\ \ell=\ell'$
or $h=h'=0,\ \ell=\ell'$.
\end{definition}

Since $\RG$ is a free Abelian group, the exact sequence
\eqref{K^(B)} splits. Fix a splitting $\ef{}:\RG\to K^{(B)}$,
$\lambda\mapsto\ef{\lambda}$.

If $v$ is a geometric valuation of~$K$, then $v|_{K^{(B)}}$ is
determined by a triple $(x,h,\ell)$, where $x\in C$, $h\in\QQ_{+}$
satisfy $v|_{K^B}=hv_x$ and
$\ell=v|_{\ef{}(\RG)}\in\ES=\Hom(\RG,\QQ)$. Therefore
$v|_{K^B}\in\Hyp$. Thus $\Vv$ is embedded in $\Hyp$, and we have a
map $\res:\Dd^B\to\Hyp$ (restriction to~$K^{(B)}$). We say that
$(\Hyp,\Vv,\Dd^B,\res)$ is the \emph{colored hyperspace}. The
valuation $v$ and the respective divisor are called \emph{central}
if $v|_{K^{(B)}}\in\ES$.

By Theorems~\ref{val.cone},~\ref{cent.cone}, and
Corollary~\ref{cosimp}, $\Vv_x=\Vv\cap\EHS{x}$ is a convex solid
polyhedral cone in~$\EHS{x}$, simplicial in characteristic~$0$,
and $\Zz=\Vv\cap\ES$ is a convex solid cone in~$\ES$.

By Corollary~\ref{D^B_v}, the set
$\Dd^B_x=\Dd^B\cap\res^{-1}(\EHS{x})$ is finite for $\forall x\in
C$. In particular, the set of central $B$-divisors is finite.

For an arbitrary $G$-model $X$, consider the rational $B$-quotient
map $\pi:X\dasharrow C$ separating generic $B$-orbits. Thus
generic $B$-orbits determine a one-parameter family of $B$-stable
prime divisors on $X$ parametrized by an open subset
$\oo{C}\subseteq C$. Decreasing $\oo{C}$ if necessary, we may
assume that these divisors do not occur in $\divr\ef{\lambda}$
($\lambda\in\RG$) and that they are pull-backs of points
$x\in\oo{C}$. Their images in $\Hyp$ are the vectors
$\eps_x=(1,0)\in\EHS{x}$. Clearly, $\{\eps_x\mid
x\in\oo{C}\}\in\Cd$.

In the generically transitive case, $\pi:\HS\dasharrow\PP^1$ is
determined by a one-dimensional linear system of $B$-divisors. In
other words, there is a $G$-line bundle $\Ll$ on $\HS$ and a
two-dimensional subspace $M$ of its $B$-eigensections which
defines this linear system. Elements of $M$ are homogeneous
coordinates on $\PP^1=\PP(M^{*})$. If $\HS=G/H$, then
$\Ll=\ind{\chi_0}$ and $M=\kk[G]^{(B\times
H)}_{(\lambda_0,-\chi_0)}$ for some $\lambda_0\in\Ch_{+}$,
$\chi_0\in\Ch(H)$. Except for finitely many lines, $M$ consists of
indecomposable elements corresponding to generic $B$-divisors.

Indecomposable elements of $M$ and the respective $B$-divisors are
called \emph{regular}. A regular $B$-divisor $D_x=\pi^{*}(x)$ is
represented in $\Hyp$ by a vector $(1,\ell)\in\EHS{x}$, and
$\ell=0$ for all but finitely many~$x$.

Besides, there is a finite set of one-dimensional subspaces
$\kk[G]^{(B\times H)}_{(\lambda_i,-\chi_i)}$, $i=1,\dots,s$,
consisting of indecomposable elements that correspond to other
$B$-divisors. If $\eta_i\in\kk[G]^{(B\times
H)}_{(\lambda_i,-\chi_i)}$ divides some $\eta\in\kk[G]^{(B\times
H)}_{(\lambda_0,-\chi_0)}$, then $D_i=\divr\eta_i$ is represented
in $\Hyp$ by $(h_i,\ell_i)\in\EHS{x}$, where $\divr\eta=D_x$ and
$h_i$ is the multiplicity of $\eta_i$ in~$\eta$ (or of $D_i$
in~$D_x$). Such $\eta_i$ and $D_i$ are called \emph{subregular}.
Other $\eta_i$ are called \emph{central} (since $D_i$ are
central).

In characteristic zero, the above description of $B$-divisors
allows to compute multiplicities in the spaces of global sections
of $G$-line bundles on~$\HS$.
\begin{proposition}
For $\forall\chi\in\Ch(H)$ and $\forall\lambda\in\Ch_{+}$, let
$k_0$ be the minimal integer such that
$(\lambda,-\chi)=\sum_{i=0}^sk_i(\lambda_i,-\chi_i)+(\mu,-\mu)$,
where $k_i\geq0$ and $\mu\in\Ch(G)$. Then
$m_{\lambda}(\ind{\chi})=k_0+1$.
\end{proposition}
\begin{proof}
Every $\eta\in\Ho^0(\HS,\ind{\chi})^{(B)}_{\lambda}=
\kk[G]^{(B\times H)}_{(\lambda,-\chi)}$ decomposes uniquely as
$\eta=\sigma_1\dots\sigma_{k_0}\eta_1^{k_1}\dots\eta_s^{k_s}$,
where $\sigma_j\in\kk[G]^{(B\times H)}_{(\lambda_0,-\chi_0)}$.
Therefore $\dim\kk[G]^{(B\times H)}_{(\lambda,-\chi)}=
\dim\kk[G]^{(B\times H)}_{(k_0\lambda_0,-k_0\chi_0)}$, and
$\kk[G]^{(B\times H)}_{(k_0\lambda_0,-k_0\chi_0)}=
\Sym^{k_0}\kk[G]^{(B\times H)}_{(\lambda_0,-\chi_0)}$ has
dimension ${k_0+1}$.
\end{proof}
\begin{corollary}[{\cite[1.2]{c=1(s)}}]
If $\Ch(H)=0$, then $m_{\lambda}(\HS)=k_0+1$, where $k_0=\max
\{k\mid\lambda-k\lambda_0\in\RG_{+}(\HS)\}$.
\end{corollary}

In the one-parameter case, generic $B$-stable divisors are
$G$-stable, whence $\eps_x\in\Vv_x$ for $x\in\oo{C}$.
\begin{lemma}
In the one-parameter case, all $B$-divisors are central.
\end{lemma}
\begin{proof}
If $D$ is a non-central $B$-divisor, then $v_D(f)>0$ for some
$f\in K^B=K^G$. Hence $D$ is $G$-stable, a contradiction.
\end{proof}

Every $B$-divisor intersects a generic $G$-orbit $\HS\subset X$
transversally, and $\ef{\lambda}$ is defined on $\HS$ for
$\forall\lambda\in\RG$. Hence $\ES(K)=\ES(\HS)$ and $\Dd^B(K)$ is
identified with $\Dd^B(\HS)$. Furthermore, it follows
from~\ref{val.cones} that $\Zz=\Vv(\HS)$ and
$\Vv_x=\Zz+\QQ_{+}\eps_x$ for $x\in\oo{C}$. Thus the colored
equipment in the one-parameter case is in the major part
determined by the colored equipment of a generic $G$-orbit: only
the structure of $\Vv_x$ for finitely many $x\in C$ depends on the
whole one-parameter family of orbits.

\begin{remark}\label{split}
Since a splitting $\ef{}:\RG\to K^{(B)}$ is not uniquely defined,
the maps $\Vv\embeds\Hyp$, $\res:\Dd^B\to\Hyp$ are not canonical.
But the change of splitting is easily controlled. If $\ef{}'$ is
another splitting, then passing from $\ef{}$ to $\ef{}'$ produces
a shift of each~$\EHS{x}$: $h'=h$ and $\ell'=\ell+h\ell_x$, where
$\langle\ell_x,\lambda\rangle=v_x(\ef{\lambda}'/\ef{\lambda})$.
The shifting vectors $\ell_x\in\RG^{*}=\Hom(\RG,\ZZ)\subset\ES$
have the property that $\sum_{x\in C}\langle\ell_x,\lambda\rangle
x$ is a principal divisor on $C$ for $\forall\lambda\in\RG$, in
particular, $\sum_{x\in C}\ell_x=0$. Conversely, any collection of
integral shifting vectors $\ell_x\in\RG^{*}$ such that $\sum_{x\in
C}\langle\ell_x,\lambda\rangle x$ is a principal divisor
for~$\forall\lambda$ defines a change of splitting. For $C=\PP^1$,
it suffices to have $\sum\ell_x=0$.
\end{remark}

Now we describe the dual object to the hyperspace.
\begin{definition}
A \emph{linear functional} on the hyperspace is a function $\phi$
on $\Hyp$ such that $\phi_x=\phi|_{\EHS{x}}$ is a $\QQ$-linear
functional for $\forall x\in C$ and $\sum_{x\in
C}\langle\eps_x,\phi_x\rangle=0$. A linear functional $\phi$ is
\emph{admissible} if $N\sum_{x\in C}\langle\eps_x,\phi_x\rangle x$
is a principal divisor on~$C$ for some $N\in\NN$. Denote by
$\Hyp^{*}$ the space of linear functionals and by $\Hyp^{*}_{\ad}$
the subspace of admissible functionals on~$\Hyp$. The set
$\Ker\phi=\bigcup_{x\in C}\Ker\phi_x$ is called the \emph{kernel}
of $\phi\in\Hyp^{*}$.
\end{definition}

If $C=\PP^1$, then any linear functional is admissible. Any
$f=f_0\ef{\lambda}\in K^{(B)}$, $f_0\in K^B$, $\lambda\in\RG$,
determines an admissible linear functional $\phi$ by means of
$\langle q,\phi_x\rangle=hv_x(f_0)+\langle\ell,\lambda\rangle$,
$\forall q=(h,\ell)\in\EHS{x}$, and $f$ is determined by $\phi$
uniquely up to a scalar multiple. Conversely, a multiple of any
admissible functional is determined by a $B$-eigenfunction.

Any collection of linear functionals $\phi_x$ on $\EHS{x}$ whose
restrictions to $\ES$ coincide can be deformed to an admissible
functional by a ``small variation''.
\begin{lemma}\label{small.var}
Let $\phi_x$ be linear functionals on $\EHS{x}$ such that
$\phi_x|_{\ES}$ does not depend on $x\in C$,
$\langle\eps_x,\phi_x\rangle=0$ for all but finitely many~$x$, and
$\sum\langle\eps_x,\phi_x\rangle<0$. Then for any finite subset
$C_0\subset C$ and $\forall\eps>0$ there exists
$\psi\in\Hyp^{*}_{\ad}$ such that $\psi_x\geq\phi_x$ on $\EHS{x}$
for $\forall x\in C$ with the equality for $x\in C_0$ (in
particular, $\psi|_{\ES}=\phi_x|_{\ES}$) and
$|\langle\eps_x,\psi_x\rangle-\langle\eps_x,\phi_x\rangle|<\eps$.
\end{lemma}
\begin{proof}
The divisor $-N\sum\langle\eps_x,\phi_x\rangle x$ is very ample on
$C$ for $N$ sufficiently large. Moving the respective hyperplane
section of $C$, we obtain an equivalent very ample divisor without
of the form $\sum n_x x$, $n_x=0,1$, $n_x=0$ whenever $x\in C_0$.
Then the divisor $\sum(N\langle\eps_x,\phi_x\rangle+n_x)x$ is
principal, and $\psi\in\Hyp^{*}_{\ad}$ defined by
$\psi|_{\ES}=\phi_x|_{\ES}$,
$\langle\eps_x,\psi_x\rangle=\langle\eps_x,\phi_x\rangle+n_x/N$,
is the desired admissible functional.
\end{proof}

For reorganizing colored data in a way similar to the spherical
case, we need some notions from the geometry of the hyperspace.

\begin{definition}
A \emph{cone} in $\Hyp$ is a cone in some~$\EHS{x}$.

A \emph{hypercone} in $\Hyp$ is a union $\Cc=\bigcup_{x\in
C}\Cc_x$ of finitely generated convex cones $\Cc_c=\Cc\cap\EHS{x}$
such that
\begin{enumerate}
\item $\Cc_x=\Kk+\QQ_{+}\eps_x$ for all but finitely many~$x$,
where $\Kk=\Cc\cap\ES$. \item
\begin{roster}
\item[Either\quad(\typeA)] $\exists x\in C:\ \Cc_x=\Kk$
\item[\hphantom{Either\quad(\typeA)}\llap{or\quad(\typeB)}]
$\Bb=\sum\Bb_x\subseteq\Kk$, where
$\eps_x+\Bb_x=\Cc_x\cap(\eps_x+\ES)$.
\end{roster}
\end{enumerate}
The hypercone is \emph{strictly convex} if all $\Cc_x$ are and
$\Bb\not\ni0$.
\end{definition}
\begin{remark}
The Minkowski sum $\sum\Bb_x$ of infinitely many polyhedral
domains $\Bb_x$ is defined as the set of all sums $\sum b_x$,
$b_x\in\Bb_x$, that make sense, i.e., $b_x=0$ for all but finitely
many~$x$. In particular, for a hypercone of type~{\typeA},
$\exists x\in C:\ \Bb_x=\emptyset$ $\implies\Bb=\emptyset$.
\end{remark}
\begin{definition}
Suppose that $Q\subseteq\Hyp$ differs from $\{\eps_x\mid
x\in\oo{C}\}$ by finitely many elements. Let $\eps_x+\Pp_x$ be the
convex hull of the intersection points of $\eps_x+\ES$ with the
rays $\QQ_{+}q$, $q\in Q$. We say that the hypercone $\Cc=\Cc(Q)$,
where $\Cc_x$ are generated by $Q\cap\EHS{x}$ and $\Pp=\sum\Pp_x$,
is \emph{generated} by~$Q$.
\end{definition}
\begin{remark}
We have $\Bb_x=\Pp_x+\Kk$ and $\Bb=\Pp+\Kk$.
\end{remark}
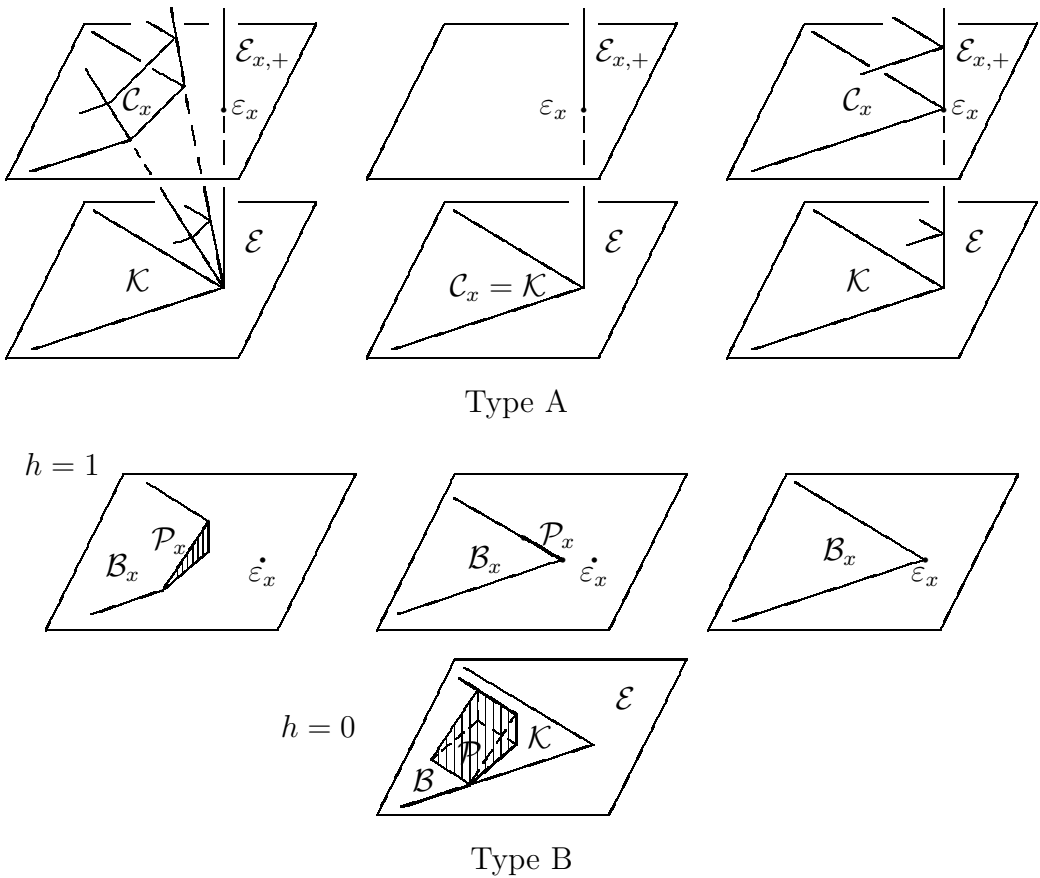
\begin{figure}[h]
\caption{Hypercones}
\begin{center}
%TexCad Options
%\grade{\off}
%\emlines{\off}
%\beziermacro{\off}
%\reduce{\on}
%\snapping{\off}
%\quality{2.00}
%\graddiff{0.01}
%\snapasp{1}
%\zoom{1.00}
\unitlength 0.25em \linethickness{0.4pt}
\begin{picture}(40.00,45.04)
%\emline(0.00,0.00)(30.00,0.00)
\put(0.00,0.00){\line(1,0){30.00}}
%\end
%\emline(30.00,0.00)(40.00,20.00)
\multiput(30.00,0.00)(0.12,0.24){84}{\line(0,1){0.24}}
%\end
%\emline(10.00,20.00)(0.00,0.00)
\multiput(10.00,20.00)(-0.12,-0.24){84}{\line(0,-1){0.24}}
%\end
%\emline(28.00,9.00)(3.00,1.00)
\multiput(28.00,9.00)(-0.37,-0.12){67}{\line(-1,0){0.37}}
%\end
%\emline(28.00,9.00)(11.00,19.00)
\multiput(28.00,9.00)(-0.20,0.12){84}{\line(-1,0){0.20}}
%\end
\put(17.00,10.00){\makebox(0,0)[cc]{$\Kk$}}
\put(32.00,15.00){\makebox(0,0)[cc]{$\ES$}}
%\emline(0.00,23.00)(30.00,23.00)
\put(0.00,23.00){\line(1,0){30.00}}
%\end
%\emline(30.00,23.00)(40.00,43.00)
\multiput(30.00,23.00)(0.12,0.24){84}{\line(0,1){0.24}}
%\end
%\emline(10.00,43.00)(0.00,23.00)
\multiput(10.00,43.00)(-0.12,-0.24){84}{\line(0,-1){0.24}}
%\end
\put(33.06,39.06){\makebox(0,0)[cc]{$\EHS{x}$}}
%\emline(28.00,32.00)(28.00,45.00)
\put(28.00,32.00){\line(0,1){13.00}}
%\end
\put(29.00,32.00){\makebox(0,0)[lc]{$\eps_x$}}
\put(17.00,33.00){\makebox(0,0)[cc]{$\Cc_x$}}
%\emline(16.00,28.00)(23.00,35.00)
\multiput(16.00,28.00)(0.12,0.12){59}{\line(0,1){0.12}}
%\end
%\emline(16.00,28.00)(3.00,24.00)
\multiput(16.00,28.00)(-0.38,-0.12){34}{\line(-1,0){0.38}}
%\end
%\emline(27.97,8.97)(19.87,21.99)
\multiput(27.97,8.97)(-0.12,0.19){68}{\line(0,1){0.19}}
%\end
%\emline(16.01,27.97)(10.22,37.13)
\multiput(16.01,27.97)(-0.12,0.19){49}{\line(0,1){0.19}}
%\end
%\emline(16.59,27.00)(17.36,25.75)
\multiput(16.59,27.00)(0.11,-0.18){7}{\line(0,-1){0.18}}
%\end
%\emline(17.94,24.88)(18.71,23.63)
\multiput(17.94,24.88)(0.11,-0.18){7}{\line(0,-1){0.18}}
%\end
%\emline(27.97,8.97)(25.46,21.99)
\multiput(27.97,8.97)(-0.12,0.62){21}{\line(0,1){0.62}}
%\end
%\emline(22.95,35.01)(21.12,45.04)
\multiput(22.95,35.01)(-0.11,0.63){16}{\line(0,1){0.63}}
%\end
%\emline(23.14,33.95)(23.63,31.44)
\multiput(23.14,33.95)(0.10,-0.50){5}{\line(0,-1){0.50}}
%\end
%\emline(23.92,29.90)(24.40,27.39)
\multiput(23.92,29.90)(0.12,-0.63){4}{\line(0,-1){0.63}}
%\end
%\emline(24.78,25.75)(25.17,23.72)
\multiput(24.78,25.75)(0.10,-0.51){4}{\line(0,-1){0.51}}
%\end
%\emline(28.03,8.96)(28.03,21.97)
\put(28.03,8.96){\line(0,1){13.01}}
%\end
%\emline(28.03,31.02)(28.03,29.00)
\put(28.03,31.02){\line(0,-1){2.02}}
%\end
%\emline(28.03,28.03)(28.03,24.95)
\put(28.03,28.03){\line(0,-1){3.08}}
%\end
%\emline(26.20,17.73)(23.89,15.41)
\multiput(26.20,17.73)(-0.12,-0.12){20}{\line(0,-1){0.12}}
%\end
%\emline(23.89,15.41)(21.48,14.64)
\multiput(23.89,15.41)(-0.34,-0.11){7}{\line(-1,0){0.34}}
%\end
%\emline(26.20,17.82)(22.93,19.75)
\multiput(26.20,17.82)(-0.19,0.11){17}{\line(-1,0){0.19}}
%\end
%\emline(13.10,32.66)(21.77,41.23)
\multiput(13.10,32.66)(0.12,0.12){72}{\line(1,0){0.12}}
%\end
%\emline(21.77,41.23)(18.21,43.35)
\multiput(21.77,41.23)(-0.20,0.12){18}{\line(-1,0){0.20}}
%\end
%\emline(13.10,32.66)(9.54,31.50)
\multiput(13.10,32.66)(-0.36,-0.12){10}{\line(-1,0){0.36}}
%\end
%\emline(23.03,34.97)(18.79,37.38)
\multiput(23.03,34.97)(-0.20,0.11){21}{\line(-1,0){0.20}}
%\end
%\emline(10.98,42.00)(17.53,38.25)
\multiput(10.98,42.00)(0.20,-0.12){32}{\line(1,0){0.20}}
%\end
\put(28.03,31.98){\circle*{0.58}}
%\emline(10.02,42.97)(15.99,42.97)
\put(10.02,42.97){\line(1,0){5.97}}
%\end
%\emline(23.03,42.97)(26.01,42.97)
\put(23.03,42.97){\line(1,0){2.98}}
%\end
%\emline(29.96,42.97)(39.98,42.97)
\put(29.96,42.97){\line(1,0){10.02}}
%\end
%\emline(10.02,20.04)(19.08,20.04)
\put(10.02,20.04){\line(1,0){9.06}}
%\end
%\emline(29.96,20.04)(39.98,20.04)
\put(29.96,20.04){\line(1,0){10.02}}
%\end
\end{picture}
\hfill
%TexCad Options
%\grade{\off}
%\emlines{\off}
%\beziermacro{\off}
%\reduce{\on}
%\snapping{\off}
%\quality{2.00}
%\graddiff{0.01}
%\snapasp{1}
%\zoom{1.00}
\unitlength 0.25em \linethickness{0.4pt}
\begin{picture}(40.00,45.00)
%\emline(0.00,0.00)(30.00,0.00)
\put(0.00,0.00){\line(1,0){30.00}}
%\end
%\emline(30.00,0.00)(40.00,20.00)
\multiput(30.00,0.00)(0.12,0.24){84}{\line(0,1){0.24}}
%\end
%\emline(10.00,20.00)(0.00,0.00)
\multiput(10.00,20.00)(-0.12,-0.24){84}{\line(0,-1){0.24}}
%\end
%\emline(28.00,9.00)(3.00,1.00)
\multiput(28.00,9.00)(-0.37,-0.12){67}{\line(-1,0){0.37}}
%\end
%\emline(28.00,9.00)(11.00,19.00)
\multiput(28.00,9.00)(-0.20,0.12){84}{\line(-1,0){0.20}}
%\end
\put(17.00,9.00){\makebox(0,0)[cc]{$\Cc_x=\Kk$}}
\put(32.00,15.00){\makebox(0,0)[cc]{$\ES$}}
%\emline(0.00,23.00)(30.00,23.00)
\put(0.00,23.00){\line(1,0){30.00}}
%\end
%\emline(30.00,23.00)(40.00,43.00)
\multiput(30.00,23.00)(0.12,0.24){84}{\line(0,1){0.24}}
%\end
%\emline(10.00,43.00)(0.00,23.00)
\multiput(10.00,43.00)(-0.12,-0.24){84}{\line(0,-1){0.24}}
%\end
\put(33.06,39.06){\makebox(0,0)[cc]{$\EHS{x}$}}
%\emline(28.00,9.00)(28.00,22.00)
\put(28.00,9.00){\line(0,1){13.00}}
%\end
%\emline(28.00,32.00)(28.00,45.00)
\put(28.00,32.00){\line(0,1){13.00}}
%\end
\put(26.00,32.00){\makebox(0,0)[rc]{$\eps_x$}}
%\emline(10.02,20.04)(26.01,20.04)
\put(10.02,20.04){\line(1,0){15.99}}
%\end
%\emline(29.96,20.04)(39.98,20.04)
\put(29.96,20.04){\line(1,0){10.02}}
%\end
%\emline(39.98,42.97)(29.96,42.97)
\put(39.98,42.97){\line(-1,0){10.02}}
%\end
%\emline(26.01,42.97)(10.02,42.97)
\put(26.01,42.97){\line(-1,0){15.99}}
%\end
%\emline(28.03,31.02)(28.03,29.00)
\put(28.03,31.02){\line(0,-1){2.02}}
%\end
\put(28.03,31.89){\circle*{0.58}}
%\emline(28.03,28.03)(28.03,24.95)
\put(28.03,28.03){\line(0,-1){3.08}}
%\end
\end{picture}
\hfill
%TexCad Options
%\grade{\off}
%\emlines{\off}
%\beziermacro{\off}
%\reduce{\on}
%\snapping{\off}
%\quality{2.00}
%\graddiff{0.01}
%\snapasp{1}
%\zoom{1.00}
\unitlength 0.25em \linethickness{0.4pt}
\begin{picture}(40.00,45.00)
%\emline(0.00,0.00)(30.00,0.00)
\put(0.00,0.00){\line(1,0){30.00}}
%\end
%\emline(30.00,0.00)(40.00,20.00)
\multiput(30.00,0.00)(0.12,0.24){84}{\line(0,1){0.24}}
%\end
%\emline(10.00,20.00)(0.00,0.00)
\multiput(10.00,20.00)(-0.12,-0.24){84}{\line(0,-1){0.24}}
%\end
%\emline(28.00,9.00)(3.00,1.00)
\multiput(28.00,9.00)(-0.37,-0.12){67}{\line(-1,0){0.37}}
%\end
%\emline(28.00,9.00)(11.00,19.00)
\multiput(28.00,9.00)(-0.20,0.12){84}{\line(-1,0){0.20}}
%\end
\put(17.00,10.00){\makebox(0,0)[cc]{$\Kk$}}
\put(32.00,15.00){\makebox(0,0)[cc]{$\ES$}}
%\emline(0.00,23.00)(30.00,23.00)
\put(0.00,23.00){\line(1,0){30.00}}
%\end
%\emline(30.00,23.00)(40.00,43.00)
\multiput(30.00,23.00)(0.12,0.24){84}{\line(0,1){0.24}}
%\end
%\emline(10.00,43.00)(0.00,23.00)
\multiput(10.00,43.00)(-0.12,-0.24){84}{\line(0,-1){0.24}}
%\end
\put(33.06,39.06){\makebox(0,0)[cc]{$\EHS{x}$}}
%\emline(28.00,9.00)(28.00,22.00)
\put(28.00,9.00){\line(0,1){13.00}}
%\end
%\emline(28.00,32.00)(28.00,45.00)
\put(28.00,32.00){\line(0,1){13.00}}
%\end
\put(29.00,32.00){\makebox(0,0)[lc]{$\eps_x$}}
%\emline(28.00,32.00)(3.00,24.00)
\multiput(28.00,32.00)(-0.37,-0.12){67}{\line(-1,0){0.37}}
%\end
\put(17.00,33.00){\makebox(0,0)[cc]{$\Cc_x$}}
%\emline(10.02,20.04)(26.01,20.04)
\put(10.02,20.04){\line(1,0){15.99}}
%\end
%\emline(29.96,20.04)(39.98,20.04)
\put(29.96,20.04){\line(1,0){10.02}}
%\end
%\emline(29.96,42.97)(39.98,42.97)
\put(29.96,42.97){\line(1,0){10.02}}
%\end
%\emline(28.03,31.02)(28.03,29.00)
\put(28.03,31.02){\line(0,-1){2.02}}
%\end
%\emline(28.03,28.03)(28.03,24.95)
\put(28.03,28.03){\line(0,-1){3.08}}
%\end
\put(28.03,31.98){\circle*{0.58}}
%\emline(28.03,15.99)(23.22,14.45)
\multiput(28.03,15.99)(-0.37,-0.12){13}{\line(-1,0){0.37}}
%\end
%\emline(28.03,15.99)(24.76,17.92)
\multiput(28.03,15.99)(-0.19,0.11){17}{\line(-1,0){0.19}}
%\end
%\emline(28.03,39.98)(17.24,36.51)
\multiput(28.03,39.98)(-0.37,-0.12){29}{\line(-1,0){0.37}}
%\end
%\emline(28.03,39.98)(21.48,43.83)
\multiput(28.03,39.98)(-0.20,0.12){33}{\line(-1,0){0.20}}
%\end
%\emline(20.33,36.61)(28.03,31.98)
\multiput(20.33,36.61)(0.20,-0.12){39}{\line(1,0){0.20}}
%\end
%\emline(10.98,42.00)(18.40,37.57)
\multiput(10.98,42.00)(0.20,-0.12){37}{\line(1,0){0.20}}
%\end
%\emline(10.02,42.97)(20.04,42.97)
\put(10.02,42.97){\line(1,0){10.02}}
%\end
\end{picture}

\medskip
Type A
\bigskip

\mbox{\rlap{\raisebox{5em}{$h=1$}}
%TexCad Options
%\grade{\off}
%\emlines{\off}
%\beziermacro{\off}
%\reduce{\on}
%\snapping{\off}
%\quality{2.00}
%\graddiff{0.01}
%\snapasp{1}
%\zoom{1.00}
\unitlength 0.25em \linethickness{0.4pt}
\begin{picture}(40.00,20.00)
%\emline(0.00,0.00)(30.00,0.00)
\put(0.00,0.00){\line(1,0){30.00}}
%\end
%\emline(30.00,0.00)(40.00,20.00)
\multiput(30.00,0.00)(0.12,0.24){84}{\line(0,1){0.24}}
%\end
%\emline(40.00,20.00)(10.00,20.00)
\put(40.00,20.00){\line(-1,0){30.00}}
%\end
%\emline(10.00,20.00)(0.00,0.00)
\multiput(10.00,20.00)(-0.12,-0.24){84}{\line(0,-1){0.24}}
%\end
%\emline(21.00,14.00)(21.00,10.00)
\put(21.00,14.00){\line(0,-1){4.00}}
%\end
%\emline(21.00,10.00)(15.00,5.00)
\multiput(21.00,10.00)(-0.14,-0.12){42}{\line(-1,0){0.14}}
%\end
\put(10.03,8.01){\makebox(0,0)[cc]{$\Bb_x$}}
%\emline(21.00,14.00)(13.00,19.00)
\multiput(21.00,14.00)(-0.19,0.12){42}{\line(-1,0){0.19}}
%\end
%\emline(15.00,5.00)(6.00,2.00)
\multiput(15.00,5.00)(-0.35,-0.12){26}{\line(-1,0){0.35}}
%\end
\put(16.00,11.50){\makebox(0,0)[cc]{$\Pp_x$}}
%\emline(21.00,14.00)(15.00,5.00)
\multiput(21.00,14.00)(-0.12,-0.18){50}{\line(0,-1){0.18}}
%\end
\put(28.00,8.00){\makebox(0,0)[ct]{$\eps_x$}}
\put(28.03,8.96){\circle*{0.58}} \put(21.00,10.02){\circle*{0.39}}
\put(21.00,13.97){\circle*{0.39}} \put(15.03,5.01){\circle*{0.39}}
\linethickness{0.05pt} \put(20.33,13.00){\line(0,-1){3.67}}
\put(19.67,12.00){\line(0,-1){3.00}}
\put(19.00,11.00){\line(0,-1){2.67}}
\put(18.33,10.00){\line(0,-1){2.33}}
\put(17.67,9.00){\line(0,-1){1.67}}
\put(17.00,8.00){\line(0,-1){1.33}}
\put(16.33,7.00){\line(0,-1){0.33}}
\end{picture}}
%TexCad Options
%\grade{\on}
%\emlines{\off}
%\beziermacro{\off}
%\reduce{\on}
%\snapping{\off}
%\quality{2.00}
%\graddiff{0.01}
%\snapasp{1}
%\zoom{1.00}
\unitlength 0.25em \linethickness{0.4pt}
\begin{picture}(40.00,20.00)
%\emline(0.00,0.00)(30.00,0.00)
\put(0.00,0.00){\line(1,0){30.00}}
%\end
%\emline(30.00,0.00)(40.00,20.00)
\multiput(30.00,0.00)(0.12,0.24){84}{\line(0,1){0.24}}
%\end
%\emline(40.00,20.00)(10.00,20.00)
\put(40.00,20.00){\line(-1,0){30.00}}
%\end
%\emline(10.00,20.00)(0.00,0.00)
\multiput(10.00,20.00)(-0.12,-0.24){84}{\line(0,-1){0.24}}
%\end
%\emline(24.00,9.00)(3.00,2.00)
\multiput(24.00,9.00)(-0.36,-0.12){59}{\line(-1,0){0.36}}
%\end
%\emline(24.00,9.00)(10.00,17.00)
\multiput(24.00,9.00)(-0.21,0.12){67}{\line(-1,0){0.21}}
%\end
\put(21.00,10.50){\makebox(0,0)[lb]{$\Pp_x$}}
\put(14.00,9.00){\makebox(0,0)[cc]{$\Bb_x$}}
\put(28.00,8.00){\makebox(0,0)[ct]{$\eps_x$}}
\put(28.03,8.96){\circle*{0.58}} \put(23.99,8.96){\circle*{0.58}}
\put(18.98,11.85){\circle*{0.58}} \linethickness{1.00pt}
%\emline(23.99,8.96)(18.98,11.85)
\multiput(23.99,8.96)(-0.20,0.12){25}{\line(-1,0){0.20}}
%\end
\end{picture}
%TexCad Options
%\grade{\off}
%\emlines{\off}
%\beziermacro{\off}
%\reduce{\on}
%\snapping{\off}
%\quality{2.00}
%\graddiff{0.01}
%\snapasp{1}
%\zoom{1.00}
\unitlength 0.25em \linethickness{0.4pt}
\begin{picture}(40.00,20.00)
%\emline(0.00,0.00)(30.00,0.00)
\put(0.00,0.00){\line(1,0){30.00}}
%\end
%\emline(30.00,0.00)(40.00,20.00)
\multiput(30.00,0.00)(0.12,0.24){84}{\line(0,1){0.24}}
%\end
%\emline(40.00,20.00)(10.00,20.00)
\put(40.00,20.00){\line(-1,0){30.00}}
%\end
%\emline(10.00,20.00)(0.00,0.00)
\multiput(10.00,20.00)(-0.12,-0.24){84}{\line(0,-1){0.24}}
%\end
%\emline(28.00,9.00)(3.00,1.00)
\multiput(28.00,9.00)(-0.37,-0.12){67}{\line(-1,0){0.37}}
%\end
%\emline(28.00,9.00)(11.00,19.00)
\multiput(28.00,9.00)(-0.20,0.12){84}{\line(-1,0){0.20}}
%\end
\put(17.00,10.00){\makebox(0,0)[cc]{$\Bb_x$}}
\put(28.00,8.00){\makebox(0,0)[ct]{$\eps_x$}}
\put(28.03,8.96){\circle*{0.58}}
\end{picture}\\*[2ex]
\mbox{\llap{\raisebox{2.5em}{$h=0$}}
%TexCad Options
%\grade{\off}
%\emlines{\off}
%\beziermacro{\off}
%\reduce{\on}
%\snapping{\off}
%\quality{2.00}
%\graddiff{0.01}
%\snapasp{1}
%\zoom{5.31}
\unitlength 0.25em \linethickness{0.4pt}
\begin{picture}(40.00,20.00)
%\emline(0.00,0.00)(30.00,0.00)
\put(0.00,0.00){\line(1,0){30.00}}
%\end
%\emline(30.00,0.00)(40.00,20.00)
\multiput(30.00,0.00)(0.12,0.24){84}{\line(0,1){0.24}}
%\end
%\emline(40.00,20.00)(10.00,20.00)
\put(40.00,20.00){\line(-1,0){30.00}}
%\end
%\emline(10.00,20.00)(0.00,0.00)
\multiput(10.00,20.00)(-0.12,-0.24){84}{\line(0,-1){0.24}}
%\end
%\emline(28.00,9.00)(3.00,1.00)
\multiput(28.00,9.00)(-0.37,-0.12){67}{\line(-1,0){0.37}}
%\end
%\emline(28.00,9.00)(11.00,19.00)
\multiput(28.00,9.00)(-0.20,0.12){84}{\line(-1,0){0.20}}
%\end
%\emline(18.00,13.00)(13.00,16.00)
\multiput(18.00,13.00)(-0.19,0.12){26}{\line(-1,0){0.19}}
%\end
%\emline(18.00,13.00)(18.00,9.00)
\put(18.00,13.00){\line(0,-1){4.00}}
%\end
%\emline(18.00,9.00)(12.00,4.00)
\multiput(18.00,9.00)(-0.14,-0.12){42}{\line(-1,0){0.14}}
%\end
%\emline(12.00,4.00)(7.00,7.00)
\multiput(12.00,4.00)(-0.19,0.12){26}{\line(-1,0){0.19}}
%\end
%\emline(7.00,7.00)(13.00,16.00)
\multiput(7.00,7.00)(0.12,0.18){50}{\line(0,1){0.18}}
%\end
\put(21.00,10.00){\makebox(0,0)[cc]{$\Kk$}}
\put(12.00,8.00){\makebox(0,0)[cc]{$\Pp$}}
\put(32.00,15.00){\makebox(0,0)[cc]{$\ES$}}
%\emline(12.04,4.05)(13.20,5.78)
\multiput(12.04,4.05)(0.12,0.17){10}{\line(0,1){0.17}}
%\end
%\emline(13.78,6.65)(14.93,8.38)
\multiput(13.78,6.65)(0.12,0.17){10}{\line(0,1){0.17}}
%\end
%\emline(15.51,9.25)(16.67,10.98)
\multiput(15.51,9.25)(0.12,0.17){10}{\line(0,1){0.17}}
%\end
%\emline(17.24,11.85)(18.02,13.01)
\multiput(17.24,11.85)(0.11,0.17){7}{\line(0,1){0.17}}
%\end
%\emline(13.01,15.99)(13.01,14.64)
\put(13.01,15.99){\line(0,-1){1.35}}
%\end
%\emline(13.01,13.39)(13.01,12.04)
\put(13.01,13.39){\line(0,-1){1.35}}
%\end
%\emline(13.01,12.04)(13.97,11.46)
\multiput(13.01,12.04)(0.19,-0.12){5}{\line(1,0){0.19}}
%\end
%\emline(14.93,10.89)(15.90,10.31)
\multiput(14.93,10.89)(0.19,-0.12){5}{\line(1,0){0.19}}
%\end
%\emline(16.86,9.73)(18.02,8.96)
\multiput(16.86,9.73)(0.17,-0.11){7}{\line(1,0){0.17}}
%\end
%\emline(13.01,12.04)(11.85,11.08)
\multiput(13.01,12.04)(-0.13,-0.11){9}{\line(-1,0){0.13}}
%\end
%\emline(11.27,10.60)(10.12,9.63)
\multiput(11.27,10.60)(-0.13,-0.11){9}{\line(-1,0){0.13}}
%\end
%\emline(9.54,9.15)(8.38,8.19)
\multiput(9.54,9.15)(-0.13,-0.11){9}{\line(-1,0){0.13}}
%\end
%\emline(7.80,7.71)(7.03,7.03)
\multiput(7.80,7.71)(-0.13,-0.11){6}{\line(-1,0){0.13}}
%\end
\linethickness{0.05pt} \put(16.99,13.57){\line(0,-1){5.35}}
\put(16.00,14.18){\line(0,-1){6.81}}
\put(15.01,14.77){\line(0,-1){8.21}}
\put(15.01,6.55){\line(0,0){0.00}}
\put(14.01,15.35){\line(0,-1){9.65}}
\put(13.02,15.96){\line(0,-1){11.08}}
\put(12.00,14.47){\line(0,-1){10.43}}
\put(11.00,12.98){\line(0,-1){8.36}}
\put(10.01,11.49){\line(0,-1){6.28}}
\put(9.01,9.97){\line(0,-1){4.18}}
\put(7.99,8.45){\line(0,-1){2.05}} \linethickness{0.8pt}
\multiput(18.00,13.00)(-0.19,0.12){39}{\line(-1,0){0.19}}
\put(18.00,13.00){\line(0,-1){4.00}}
\multiput(18.00,9.00)(-0.14,-0.12){42}{\line(-1,0){0.14}}
\multiput(12.00,3.80)(-0.37,-0.12){24}{\line(-1,0){0.37}}
\put(6.03,4.50){\makebox(0,0)[cc]{$\Bb$}}
\end{picture}}

\medskip
Type B
\end{center}
\end{figure}
\begin{definition}
For a hypercone $\Cc$ of type~{\typeB}, we define its
\emph{interior} $\intr\Cc=\bigcup_{x\in C}\intr\Cc_x\cup\intr\Kk$.

A \emph{face} of a hypercone $\Cc$ is a face $\Cc'$ of some
$\Cc_x$ such that $\Cc'\cap\Bb=\emptyset$.

A \emph{hyperface} of $\Cc$ is a hypercone $\Cc'=\Cc\cap\Ker\phi$,
where $\phi\in\Hyp^{*}$, $\langle\Cc,\phi\rangle\geq0$. The
hyperface $\Cc'$ is \emph{admissible} if any such $\phi$ is
admissible.

A hypercone is \emph{admissible} if all its hyperfaces of
type~{\typeB} are admissible.
\end{definition}
\begin{remark}
A hyperface $\Cc'\subseteq\Cc$ is of type~{\typeB} iff
$\Cc'\cap\Bb\neq\emptyset$. Indeed,
$\langle\eps_x+\Bb_x,\phi_x\rangle\geq0$ ($\forall x$) and
$\forall x:\ \Cc'_x\not\subseteq\ES\iff \forall x:\
\langle\eps_x+\Bb_x,\phi_x\rangle\ni0\iff
\sum\langle\eps_x+\Bb_x,\phi_x\rangle=
\sum\langle\eps_x,\phi_x\rangle+
\sum\langle\Bb_x,\phi_x\rangle=\langle\Bb,\phi\rangle\ni0$
\end{remark}

Properties of hypercones are similar to properties of convex
polyhedral cones. Let $\Cc$ be a hypercone. There is a separation
property:
\begin{lemma}\label{hc.sep}
\begin{roster}
\item\label{nsc} $q\notin\Cc\implies\exists\phi\in\Hyp^{*}:\
\langle\Cc,\phi\rangle\geq0,\ \langle q,\phi\rangle<0$
\item\label{sc} If $\Cc$ is strictly convex, then one may assume
that $\phi$ is admissible and
$\langle\Cc_x\setminus\{0\},\phi_x\rangle>0$ for any given finite
set of~$x$.
\end{roster}
\end{lemma}
\begin{proof}
\begin{roster}
\item[\ref{nsc}] If $q\in\EHS{y}$, then we construct a collection
of functionals $\phi_x$ on $\EHS{x}$ such that
$\phi_x|_{\ES}=\phi_y|_{\ES}$, $\langle q,\phi_y\rangle<0$,
$\langle\Cc_x,\phi_x\rangle\geq0$ for~$\forall x$ and the equality
is reached on $\Cc_x\setminus\{0\}$. Then
$\sum\langle\eps_x+\Bb_x,\phi_x\rangle=
\sum\langle\eps_x,\phi_x\rangle+\langle\Bb,\phi_y\rangle\geq0$ and
the equality is reached. But $\langle\Bb,\phi_y\rangle\geq0
\implies\sum\langle\eps_x,\phi_x\rangle\leq0$. It remains to
modify $\phi_x$ by a small variation (Lemma~\ref{small.var}) if
necessary. \item[\ref{sc}] In the proof of~\ref{nsc}, we may
assume $\langle\Kk\setminus\{0\},\phi_y\rangle>0\implies
{\langle\Bb,\phi_y\rangle>0}\implies
\sum\langle\eps_x,\phi_x\rangle<0$, and we may increase finitely
many $\phi_x$ to have $\langle\Cc_x,\phi_x\rangle\geq0$.
\qedhere\end{roster}
\end{proof}
This implies a dual characterization of a hypercone:
\begin{lemma}\label{dual}
For a (strictly convex) hypercone $\Cc=\Cc(Q)$, $q\in\Cc$ iff
$\langle Q,\phi\rangle\geq0\implies\langle q,\phi\rangle\geq0$ for
all (admissible)~$\phi$.
\end{lemma}
\begin{proof}
$\langle Q,\phi\rangle\geq0\implies\forall x\in C:\
\langle\eps_x+\Pp_x,\phi_x\rangle\geq0\implies
\langle\Pp,\phi\rangle\geq0\implies\langle\Cc,\phi\rangle\geq0$.
Lemma~\ref{hc.sep} completes the proof.
\end{proof}

For any $v\in\Cc$, there is a unique face or hyperface of
type~{\typeB} containing $v$ in its interior.
\begin{lemma}\label{face}
The face or (admissible) hyperface $\Cc'\subseteq\Cc$ such that
$v\in\intr\Cc'$ is the intersection of (admissible) hyperfaces of
$\Cc$ containing~$v$.
\end{lemma}
\begin{proof}
If $\langle\Cc,\phi\rangle\geq0$, $\langle v,\phi\rangle=0$, then
$\langle\Cc',\phi\rangle=0$. (If $\Cc'$ is a hyperface, then
$\langle\Kk',\phi\rangle=0\implies
\langle\eps_x+\Bb'_x,\phi_x\rangle=0\implies
\langle\Cc'_x,\phi_x\rangle=0$ for~$\forall x$.) Conversely, if
$\Cc'\subseteq\Cc_y$ is a face and $q\in\Cc\setminus\Cc'$, then we
construct an (admissible) functional $\phi$ such that
$\langle\Cc,\phi\rangle\geq0$, $\langle\Cc',\phi\rangle=0$,
$\langle q,\phi\rangle>0$ as follows. Take $\phi_y$ on $\EHS{y}$
such that $\langle\Cc_y,\phi_y\rangle\geq0$ and
$\Cc'=\Ker\phi_y\cap\Cc$. We may include $\phi_y$ in a collection
of functionals $\phi_x$ on $\EHS{x}$ such that
$\phi_x|_{\ES}=\phi_y|_{\ES}$, $\langle\Cc_x,\phi_x\rangle\geq0$
and the inequality is reached on $\Cc_x\setminus\{0\}$. But
$\Cc'\cap\Bb=\emptyset\implies\langle\Bb,\phi_y\rangle>0\implies
\sum\langle\eps_x,\phi_x\rangle<0$, as in
Lemma~\ref{hc.sep}\ref{sc}. Then we increase some $\phi_x$ to
obtain $\langle q,\phi\rangle>0$ and apply Lemma~\ref{small.var}.
\end{proof}

Now let $(\Ww,\Rr)$ be colored data from $\Cd$ and consider the
hypercone $\Cc=\Cc(\Ww,\Rr)$ generated by $\Ww\cup\res(\Rr)$.

Condition \reftag{C} means that $\Cc$ is strictly convex and
$\res(\Rr)\not\ni0$ (Lemma~\ref{hc.sep}\ref{sc}). We assume it in
the sequel.

Condition \reftag{W} means that the elements of $\Ww$ are exactly
the generators of those edges of $\Cc$ that do not
intersect~$\res(\Rr)$. (Indeed, \reftag{W} $\iff
w\notin\Cc(\Ww\setminus\{w\},\Rr)$.)

Condition \reftag{F} means that $\Cc$ is admissible. (This is
non-trivial, see \cite[Pr.4.1]{LVT}.)

Condition \reftag{V} means that $v\in\Cc$ (Lemma~\ref{dual}).

Conditions \reftag{V$'$} and \reftag{D$'$} say that $\Vv_Y$ and
$\Dd^B_Y$ consist of those elements of $\Ww$ and $\Rr$ which lie
in the face or hyperface (of type~{\typeB})
$\Cc_Y=\Cc(\Vv_Y,\Dd^B_Y)\subseteq\Cc$ such that $v\in\intr\Cc_Y$
(Lemma~\ref{face}).

Condition \reftag{S} says that $v\in\Vv\cap\intr\Cc_Y$
(Lemma~\ref{face}).

\begin{definition}
A \emph{colored hypercone} is a pair $(\Cc,\Rr)$, where
$\Rr\subseteq\Dd^B$, $\res(\Rr)\not\ni0$, and $\Cc$ is a strictly
convex hypercone generated by $\res(\Rr)$ and $\Ww\subseteq\Vv$.

A colored hypercone $(\Cc,\Rr)$ (of type~{\typeB}) is
\emph{supported} if $\intr\Cc\cap\Vv\neq\emptyset$.

A \emph{(hyper)face} of $(\Cc,\Rr)$ is a colored (hyper)cone
$(\Cc',\Rr')$, where $\Cc'$ is a (hyper)face of~$\Cc$ and
$\Rr'=\Rr\cap\res^{-1}(\Cc')$.

A \emph{colored hyperfan} is a finite set of supported colored
cones and hypercones of type~{\typeB} whose interiors are disjoint
inside~$\Vv$ and which is obtained as the set of all supported
(hyper)faces of finitely many colored hypercones.
\end{definition}
\begin{theorem}
\begin{roster}
\item $B$-charts are in bijection with colored hypercones
in~$\Hyp$. \item $G$-germs are in bijection with supported colored
cones and hypercones of type~{\typeB}. If $Y_1,Y_2\subseteq X$ are
$G$-subvarieties in a $G$-model, then $Y_1\adhereseq{Y_2}$ iff
$(\Cc_{Y_2},\Dd^B_{Y_2})$ is a (hyper)face of
$(\Cc_{Y_1},\Dd^B_{Y_1})$. \item $G$-models are in bijection with
colored hyperfans.
\end{roster}
\end{theorem}
Corollaries~\ref{sph.aff}--\ref{sph.comp} are easily generalized
to the case of complexity~$1$.

\begin{remark}
Let $\X$ be a $B$-chart defined by a colored hypercone
$(\Cc,\Rr)$. Then $\kk[\X]^B=\kk$ iff $\Cc$ is of type~{\typeB}:
if $f\in\kk[\X]^B$, then the respective $\phi\in\Hyp^{*}_{\ad}$
must be zero on~$\ES$ and non-negative on~$\Cc$. Thus we have two
types of $B$-charts:
\begin{enumerate}
\item[(\typeA)] $\kk[\X]^B\neq\kk$, or $\Cc$ is of type~{\typeA}.
\item[(\typeB)] $\kk[\X]^B=\kk$, or $\Cc$ is of type~{\typeB}.
\end{enumerate}

There are two types of $G$-germs:
\begin{enumerate}
\item[(\typeA)] $\Cc(\Vv_Y,\Dd^B_Y)$ is a colored cone.
\item[(\typeB)] $\Cc(\Vv_Y,\Dd^B_Y)$ is a colored hypercone.
\end{enumerate}
A $G$-germ is of type~{\typeA} iff $\Vv_Y,\Dd^B_Y$ are finite, and
of type~{\typeB} iff it has a minimal $B$-chart.
\end{remark}

\begin{example}\label{c_T=1}
Suppose that $G=B=T$ is a torus. We may assume (after factoring
out by the kernel of the action) that the stabilizer of general
position for any $T$-model is trivial.  Since a torus has no
non-trivial Galois cohomology (Hilbert's Theorem~90, see
\cite[2.6]{IT}), the birational type of the action is trivial,
i.e., any $T$-model is birationally isomorphic to $T\times C$. It
follows that $\ES=\Hom(\Ch(T),\QQ)$, $\Vv=\Hyp$,
$\Dd^B=\Dd^B(T)=\emptyset$.

A $T$-model is given by a set of cones and admissible hypercones
of type~{\typeB} with disjoint interiors which consists of all
faces and hyperfaces of type~{\typeB} of finitely many admissible
hypercones. (The word ``colored'' is needless here, since there
are no colors.) A $T$-chart $\X$ is of type~{\typeA}
(type~{\typeB}) iff $\X\by T\subset C$ is an open subset ($\X\by
T$ is a point).

If all germs of a $T$-model $X$ are of type~{\typeA}, then
quotient morphisms of its $T$-charts may be glued together into a
regular map $\pi:X\to C$ separating $T$-orbits of general
position. Such $T$-models were classified by Mumford in
\cite[Part~IV]{tor.emb} in the framework of the theory of toroidal
embeddings (for this theory see \cite[Part~II]{tor.emb}). The
hyperfan of $X$ is a union of fans $\Ff_x$ in $\EHS{x}$ having
common central part $\Ff=\{\Cc\in\Ff_x\mid\Cc\subseteq\ES\}$ and
such that $\Ff_x$ is a cylinder over $\Ff$ for $x\neq
x_1,\dots,x_s$ (finitely many exceptional points).

It is proved in \cite[Part~IV]{tor.emb} that $C$ is covered by
open neighborhoods $C_i$ of~$x_i$ such that $\pi^{-1}(C_i)\iso
C_i\times_{\AAA^1}X_i$, where $\nu_i:C_i\to\AAA^1$ are etale maps
such that $\nu_i^{-1}(0)=\{x_i\}$ and $X_i$ are toric
$(T\times\kk^{\times})$-varieties with fans $\Ff_{x_i}$ mapping
$\kk^{\times}$-equivariantly onto~$\AAA^1$.
\end{example}
\begin{remark}
The admissibility of a hypercone is essential for
condition~\reftag{F} as the following example~\cite{U-fin} shows.

Let $C$ be a smooth projective curve of genus $g>0$ and
$\delta_i=\sum n_{ix}x$, $i=1,2$, be divisors on $C$ having
infinite order in $\Pic C$ and such that $\deg\delta_1=0$,
$\deg\delta_2>2g-2$. Put $\Ll_i=\Lin[C]{\delta_i}$.

The total space $X$ of $\Ll_1^{*}\oplus\Ll_2^{*}$ is a
$T=(\kk^{\times})^2$-model, where the factors $\kk^{\times}$ act
on $\Ll_i^{*}$ by homotheties. There are the following $T$-germs
in~$X$: the divisors $D_i$, $D_x$ ($i=1,2$, $x\in C$), where $D_i$
is the total space of~$\Ll_j^{*}$, $\{i,j\}=\{1,2\}$, and $D_x$ is
the fiber of $X\to C$ over~$\{x\}$; $Y_{ix}=D_i\cap D_x$;
$C=D_1\cap D_2$; $\{x\}=D_1\cap D_2\cap D_x$.

Let $\ef{i}$ be a rational section of $\Ll_i$ such that $\divr
\ef{i}=\delta_i$. Then $\RG(X)$ is generated by the $T$-weights
$\omega_1,\omega_2$ of $\ef1,\ef2$. If $w_i,w_x$ are the
$T$-valuations corresponding to $D_i,D_x$, then $w_1=(0,1,0)$,
$w_2=(0,0,1)$, $w_x=(1,n_{1x},n_{2x})$ in the basis
$\eps_x,\omega^{\vee}_1,\omega^{\vee}_2$, where $\omega^{\vee}_i$
are the dual coweights to~$\omega_i$.

The algebra $\kk[X]$ is bigraded by the $T$-action:
$\kk[X]_{m,n}=\Ho^0(C,\Ll_1^m\otimes\Ll_2^n)\neq0$ iff $m\geq0,\
n>0$ or $m=n=0$. Hence $\RG_{+}(X)$ is not finitely generated, and
$\kk[X]$ as well. The reason is that $\kk[X]$ is defined by a
non-admissible hypercone $\Cc=\bigcup\Cc_{\{x\}}$. Indeed, its
hyperface $\Cc'=\bigcup\Cc_{Y_{2x}}$ is not admissible, because
$\Cc'=\Cc\cap\Ker\phi$, $\langle\omega^{\vee}_1,\phi\rangle=1$,
$\langle\omega^{\vee}_2,\phi\rangle=0$,
$\langle\eps_x,\phi\rangle=-n_{1x}$, and no multiple of
$-\delta_1=\sum\langle\eps_x,\phi\rangle x$ is a principal
divisor. See Figure~\ref{inf.gen}.
\begin{figure}[h]
\caption{}\label{inf.gen}
\begin{center}
%TexCad Options
%\grade{\off}
%\emlines{\off}
%\beziermacro{\on}
%\reduce{\on}
%\snapping{\off}
%\quality{2.00}
%\graddiff{0.01}
%\snapasp{1}
%\zoom{3.00}
\unitlength 0.50ex \linethickness{0.4pt}
\begin{picture}(53.20,31.58)
\put(0.00,0.00){\line(1,0){40.00}}
\put(40.00,0.00){\line(3,5){13.20}}
\put(53.20,22.00){\line(-1,0){8.20}}
\put(18.00,22.00){\line(-1,0){5.00}}
\put(13.00,22.00){\line(-3,-5){13.20}} \thicklines{}
\put(15.00,13.00){\line(2,-1){20.33}}
\put(29.63,31.58){\line(-4,-5){14.80}}
\put(15.00,13.00){\line(5,1){15.78}}
\put(46.00,19.00){\line(-5,-1){12.44}}
\put(26.00,7.44){\vector(2,-1){1.89}}
\put(36.11,17.00){\vector(4,1){2.78}}
\put(24.11,24.56){\vector(3,4){1.58}}
\put(10.78,5.00){\makebox(0,0)[cc]{$\ES$}}
\put(27.89,5.00){\makebox(0,0)[ct]{$w_1$}}
\put(38.00,16.00){\makebox(0,0)[ct]{$w_2$}}
\put(26.22,15.22){\circle*{1.11}}
\put(26.22,14.00){\makebox(0,0)[ct]{$\Pp$}}
\put(24.50,26.50){\makebox(0,0)[rb]{$w_x$}}
\put(27.00,21.00){\makebox(0,0)[cc]{$\Cc'_x$}}
\put(41.00,27.00){\makebox(0,0)[lb]{$\Cc_{\{x\}}$}} \thinlines{}
\linethickness{0.3pt}
%\emline(34.00,3.56)(29.00,30.56)
\multiput(34.00,3.56)(-0.12,0.64){42}{\line(0,1){0.64}}
%\end
%\emline(29.00,30.56)(45.00,18.89)
\multiput(29.00,30.56)(0.16,-0.12){98}{\line(1,0){0.16}}
%\end
%\emline(45.00,18.89)(34.00,3.56)
\multiput(45.00,18.89)(-0.12,-0.17){92}{\line(0,-1){0.17}}
%\end
\linethickness{0.05pt} \put(44.00,18.67){\line(0,1){0.89}}
\put(43.00,18.44){\line(0,1){1.89}}
\put(42.00,18.22){\line(0,1){2.78}}
\put(41.00,18.11){\line(0,1){3.56}}
\put(40.00,17.89){\line(0,1){4.56}}
\put(39.00,17.67){\line(0,1){5.56}}
\put(38.00,17.44){\line(0,1){6.44}}
\put(37.00,17.22){\line(0,1){7.33}}
\put(36.00,17.00){\line(0,1){8.33}}
\put(35.00,16.89){\line(0,1){9.11}}
\put(34.00,16.78){\line(0,1){10.00}}
\put(33.00,16.56){\line(0,1){10.89}}
\put(32.00,28.22){\line(0,-1){11.78}}
\put(31.00,16.44){\line(0,1){12.56}}
\put(30.00,29.78){\line(0,-1){13.78}}
\put(29.00,15.78){\line(0,1){14.67}}
\put(28.00,29.56){\line(0,-1){13.89}}
\put(27.00,15.44){\line(0,1){12.67}}
\put(26.00,27.11){\line(0,-1){12.00}}
\put(25.00,14.89){\line(0,1){10.78}}
\put(24.00,24.44){\line(0,-1){9.67}}
\put(23.00,14.56){\line(0,1){8.67}}
\put(22.00,22.00){\line(0,-1){7.56}}
\put(21.00,14.22){\line(0,1){6.56}}
\put(20.00,19.56){\line(0,-1){5.56}}
\put(19.00,13.78){\line(0,1){4.44}}
\put(18.00,17.22){\line(0,-1){3.56}}
\put(17.00,13.44){\line(0,1){2.33}}
\put(16.00,14.56){\line(0,-1){1.33}}
\end{picture}
\qquad\qquad
%TexCad Options
%\grade{\off}
%\emlines{\off}
%\beziermacro{\off}
%\reduce{\on}
%\snapping{\off}
%\quality{2.00}
%\graddiff{0.01}
%\snapasp{1}
%\zoom{1.00}
\unitlength 0.50ex \linethickness{0.4pt}
\begin{picture}(35.00,20.00)
%\vector(0.00,0.00)(0.00,20.00)
\put(0.00,20.00){\vector(0,1){0.2}}
\put(0.00,0.00){\line(0,1){20.00}}
%\end
%\vector(0.00,0.00)(35.00,0.00)
\put(35.00,0.00){\vector(1,0){0.2}}
\put(0.00,0.00){\line(1,0){35.00}}
%\end
\put(0.00,0.00){\circle*{1.50}} \put(5.00,0.00){\circle*{0.60}}
\put(10.00,0.00){\circle*{0.60}} \put(15.00,0.00){\circle*{0.60}}
\put(20.00,0.00){\circle*{0.60}} \put(25.00,0.00){\circle*{0.60}}
\put(0.00,5.00){\circle*{1.50}} \put(5.00,5.00){\circle*{1.50}}
\put(10.00,5.00){\circle*{1.50}} \put(15.00,5.00){\circle*{1.50}}
\put(20.00,5.00){\circle*{1.50}} \put(25.00,5.00){\circle*{1.50}}
\put(0.00,10.00){\circle*{1.50}} \put(0.00,15.00){\circle*{1.50}}
\put(5.00,10.00){\circle*{1.50}} \put(5.00,15.00){\circle*{1.50}}
\put(10.00,10.00){\circle*{1.50}}
\put(10.00,15.00){\circle*{1.50}}
\put(15.00,10.00){\circle*{1.50}}
\put(15.00,15.00){\circle*{1.50}}
\put(20.00,10.00){\circle*{1.50}}
\put(20.00,15.00){\circle*{1.50}}
\put(25.00,10.00){\circle*{1.50}}
\put(25.00,15.00){\circle*{1.50}} \put(30.00,0.00){\circle*{0.60}}
\put(30.00,5.00){\circle*{1.50}} \put(30.00,10.00){\circle*{1.50}}
\put(30.00,15.00){\circle*{1.50}}
\put(32.00,17.00){\makebox(0,0)[lb]{$\RG_{+}(X)$}}
\end{picture}
\qquad\qquad\strut
\end{center}
\end{figure}
\end{remark}

\begin{example}[$\SL_2$-embeddings]\label{SL2/e}
Suppose $G=\SL_2(\kk)$, $H=\{\1\}$. Then $\HS=\SL_2$ has
complexity one. Its embeddings were described in \cite[\S9]{LV}.

The elements of $G$ are matrices
%*
\begin{equation*}
g=\begin{pmatrix}
g_{11} & g_{12} \\
g_{21} & g_{22}
\end{pmatrix},\qquad
g_{11}g_{22}-g_{21}g_{12}=1,
\end{equation*}
%*
and $B$ consists of upper-triangular matrices ($g_{21}=0$). Let
$\omega$ be the fundamental weight: $\omega(g)=g_{11}$ for $g\in
B$.

All $B$-divisors in $\HS$ are regular. Their equations are the
nonzero elements of the two-dimensional subspace
$M=\kk[G]^{(B)}_{\omega}$ generated by $\eta_1(g)=g_{21}$,
$\eta_2(g)=g_{22}$.  Let $\eta_x=\alpha_1\eta_1+\alpha_2\eta_2$ be
an equation of $x\in\PP^1=\PP(M^{*})$.

The field $K^B=\kk(\PP^1)$ consists of rational functions in
$\eta_1,\eta_2$ of degree~$0$.

The group $\RG$ equals $\Ch(B)=\langle\omega\rangle\iso\ZZ$. We
may take $\ef{\omega}=\eta_{\infty}$, where $\infty\in\PP^1$ is a
certain fixed point.

The set of $G$-valuations is computed by the method of formal
curves~(\ref{form.curv}). First we determine $G$-valuations
corresponding to divisors with a dense orbit. Up to a multiple,
any such valuation is defined by the formula
$v_{x(t)}(f)=\ord_tf(gx(t))$, where $x(t)\in\SL_2(\kk\(t\))$ and
$g$~is the generic $\kk(\SL_2)$-point of~$\SL_2$. By the Iwasawa
decomposition (\ref{form.curv}), we may even assume $x(t)=\left(
\begin{smallmatrix}
t^m &  u(t) \\
 0  & t^{-m}
\end{smallmatrix}\right)$,
$u(t)\in\kk\(t\)$, $\ord_tu(t)=n\leq-m$.

The number
%*
\begin{equation*}
d=v_{x(t)}(\eta_x)=\ord_t\bigl((\alpha_1t^m+\alpha_2u(t))g_{21}+
\alpha_2t^{-m}g_{22}\bigr)
\end{equation*}
%*
is constant along $\PP^1$ except one~$x$, where it jumps, so that
$v_{x(t)}=(h,\ell)\in\Vv_x$. The following cases are possible:
%*
\begin{xxalignat}{1}
m\leq n &\implies
\begin{cases}
d=m,        & \ord_t(\alpha_1t^m+\alpha_2u(t))=m \\
d\in(m,-m], & \ord_t(\alpha_1t^m+\alpha_2u(t))>m
\end{cases}
\implies
\begin{cases}
h\in(0,-2m] \\
\ell=m\text{ (or $m+h$)}
\end{cases}
\\
m>n &\implies
\begin{cases}
d=n, & \alpha_2\ne0 \\
d=m, & \alpha_2=0
\end{cases}
\implies
\begin{cases}
h=m-n \\
\ell=n\text{ (or $n+h$)}
\end{cases}
\end{xxalignat}
%*
(Here $\ell=v_{x(t)}(\ef{\omega})$ increases by $h$ if the jump
occurs at $x=\infty$.)

In both cases, we obtain the subset in $\EHS{x}$ defined by the
inequalities $h>0$, $2\ell+h\leq0$ (or $2\ell-h\leq0$ if
$x=\infty$). Thus $\Vv_x$ is defined by $2\ell+h\leq0$ (or
$2\ell-h\leq0$) by~\ref{form.curv}. The colored equipment is
represented in Figure~\ref{hs(SL2)}. (Elements of
$\RG^{*}\times\ZZ_{+}\subset\EHS{x}$ are marked by dots.)
\begin{figure}[h]
\caption{}\label{hs(SL2)}
\begin{center}
%TexCad Options
%\grade{\off}
%\emlines{\off}
%\beziermacro{\off}
%\reduce{\on}
%\snapping{\off}
%\quality{2.00}
%\graddiff{0.01}
%\snapasp{1}
%\zoom{3.00}
\unitlength 0.70ex \linethickness{0.4pt}
\begin{picture}(30.00,25.00)
%\vector(0.00,5.00)(30.00,5.00)
\put(30.00,5.00){\vector(1,0){0.2}}
\put(0.00,5.00){\line(1,0){30.00}}
%\end
\put(5.00,5.00){\circle*{0.50}} \put(10.00,5.00){\circle*{0.50}}
\put(15.00,5.00){\circle*{0.50}} \put(20.00,5.00){\circle*{0.50}}
\put(25.00,5.00){\circle*{0.50}} \put(5.00,10.00){\circle*{0.50}}
\put(5.00,15.00){\circle*{0.50}} \put(5.00,20.00){\circle*{0.50}}
\put(10.00,10.00){\circle*{0.50}}
\put(10.00,15.00){\circle*{0.50}}
\put(10.00,20.00){\circle*{0.50}} \put(15.00,10.00){\circle{1.00}}
\put(15.00,15.00){\circle*{0.50}}
\put(15.00,20.00){\circle*{0.50}}
\put(20.00,10.00){\circle*{0.50}}
\put(20.00,15.00){\circle*{0.50}}
\put(20.00,20.00){\circle*{0.50}}
\put(25.00,10.00){\circle*{0.50}}
\put(25.00,15.00){\circle*{0.50}}
\put(25.00,20.00){\circle*{0.50}}
\put(4.00,13.00){\makebox(0,0)[rc]{$\Vv_x$}}
\put(16.00,9.89){\makebox(0,0)[lc]{$\Dd_x$}}
\put(16.00,25.00){\makebox(0,0)[lc]{$h$}}
\put(30.00,6.00){\makebox(0,0)[cb]{$\ell$}}
\put(15.00,3.00){\makebox(0,0)[ct]{\strut$x\neq\infty$}}
\put(15.00,5.00){\line(0,1){4.50}}
\put(15.00,10.50){\vector(0,1){14.50}}
\put(15.00,5.00){\line(-1,2){10.00}} \linethickness{0.05pt}
\put(0.00,7.00){\line(1,0){14.00}}
\put(0.00,9.00){\line(1,0){13.00}}
\put(0.00,11.00){\line(1,0){12.00}}
\put(0.00,13.00){\line(1,0){11.00}}
\put(0.00,15.00){\line(1,0){10.00}}
\put(0.00,17.00){\line(1,0){9.00}}
\put(0.00,19.00){\line(1,0){8.00}}
\put(0.00,21.00){\line(1,0){7.00}}
\put(0.00,23.00){\line(1,0){6.00}}
\end{picture}
\qquad\qquad
%TexCad Options
%\grade{\off}
%\emlines{\off}
%\beziermacro{\off}
%\reduce{\on}
%\snapping{\off}
%\quality{2.00}
%\graddiff{0.01}
%\snapasp{1}
%\zoom{3.00}
\unitlength 0.70ex \linethickness{0.4pt}
\begin{picture}(30.00,25.00)
%\vector(0.00,5.00)(30.00,5.00)
\put(30.00,5.00){\vector(1,0){0.2}}
\put(0.00,5.00){\line(1,0){30.00}}
%\end
\put(5.00,5.00){\circle*{0.50}} \put(10.00,5.00){\circle*{0.50}}
\put(15.00,5.00){\circle*{0.50}} \put(20.00,5.00){\circle*{0.50}}
\put(25.00,5.00){\circle*{0.50}} \put(5.00,10.00){\circle*{0.50}}
\put(5.00,15.00){\circle*{0.50}} \put(5.00,20.00){\circle*{0.50}}
\put(10.00,10.00){\circle*{0.50}}
\put(10.00,15.00){\circle*{0.50}}
\put(10.00,20.00){\circle*{0.50}} \put(20.00,10.00){\circle{1.00}}
\put(15.00,15.00){\circle*{0.50}}
\put(15.00,20.00){\circle*{0.50}}
\put(15.00,10.00){\circle*{0.50}}
\put(20.00,15.00){\circle*{0.50}}
\put(20.00,20.00){\circle*{0.50}}
\put(25.00,10.00){\circle*{0.50}}
\put(25.00,15.00){\circle*{0.50}}
\put(25.00,20.00){\circle*{0.50}}
\put(21.00,9.89){\makebox(0,0)[lc]{$\Dd_x$}}
\put(16.00,25.00){\makebox(0,0)[lc]{$h$}}
\put(30.00,6.00){\makebox(0,0)[cb]{$\ell$}}
\put(15.00,3.00){\makebox(0,0)[ct]{\strut$x=\infty$}}
\put(14.00,16.00){\makebox(0,0)[rb]{$\Vv_{\infty}$}}
\put(15.00,5.00){\vector(0,1){20.00}}
\put(15.00,5.00){\line(1,2){10.00}} \linethickness{0.05pt}
\put(0.00,7.00){\line(1,0){15.89}}
\put(0.00,9.00){\line(1,0){17.00}}
\put(0.00,11.00){\line(1,0){17.89}}
\put(0.00,13.00){\line(1,0){18.89}}
\put(0.00,15.00){\line(1,0){20.00}}
\put(0.00,17.00){\line(1,0){21.00}}
\put(0.00,19.00){\line(1,0){21.89}}
\put(0.00,21.00){\line(1,0){22.89}}
\put(0.00,23.00){\line(1,0){24.00}}
\end{picture}
\end{center}
\end{figure}

$G$-germs are given by colored cones or hypercones of
type~{\typeB} hatched in the figures below; their colors are
marked by bold dots. The notation for $G$-germs is taken from
\cite[\S9]{LV}.
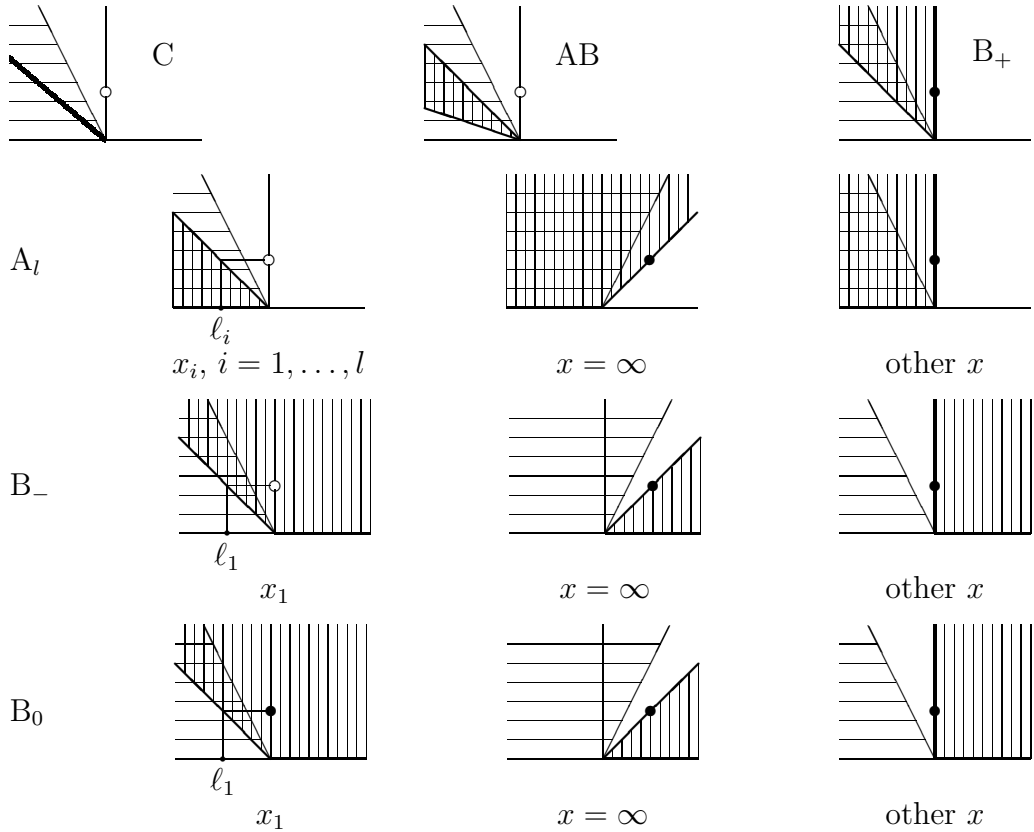
\begin{figure}[h]
\caption{$G$-germs of $\SL_2$-embeddings}\label{germs(SL2)}
\begin{center}
%TexCad Options
%\grade{\on}
%\emlines{\off}
%\beziermacro{\off}
%\reduce{\on}
%\snapping{\off}
%\quality{2.00}
%\graddiff{0.01}
%\snapasp{1}
%\zoom{4.00}
\unitlength 0.70ex \linethickness{0.4pt}
\begin{picture}(20.00,14.00)
\put(0.00,0.00){\line(1,0){20.00}} \put(10.00,5.00){\circle{1.00}}
\put(10.00,5.50){\line(0,1){8.50}}
\put(10.00,0.00){\line(0,1){4.50}}
\put(10.00,0.00){\line(-1,2){7.00}}
\put(16.00,9.00){\makebox(0,0)[cc]{$\text{C}$}}
\linethickness{0.05pt} \put(0.00,2.00){\line(1,0){9.00}}
\put(0.00,4.00){\line(1,0){8.00}}
\put(0.00,6.00){\line(1,0){7.00}}
\put(0.00,8.00){\line(1,0){5.89}}
\put(0.00,10.00){\line(1,0){4.89}}
\put(0.00,12.00){\line(1,0){4.00}} \thicklines{}
\linethickness{2.00pt}
%\emline(10.00,0.00)(0.00,8.56)
\multiput(10.00,0.00)(-0.14,0.12){72}{\line(-1,0){0.14}}
%\end
\end{picture}
\hfill
%TexCad Options
%\grade{\off}
%\emlines{\off}
%\beziermacro{\off}
%\reduce{\on}
%\snapping{\off}
%\quality{2.00}
%\graddiff{0.01}
%\snapasp{1}
%\zoom{4.00}
\unitlength 0.70ex \linethickness{0.4pt}
\begin{picture}(20.00,14.00)
\put(0.00,0.00){\line(1,0){20.00}} \put(10.00,5.00){\circle{1.00}}
\put(10.00,5.50){\line(0,1){8.50}}
\put(10.00,0.00){\line(0,1){4.50}}
\put(10.00,0.00){\line(-1,2){7.00}}
\put(16.00,9.00){\makebox(0,0)[cc]{$\text{AB}$}}
\linethickness{0.05pt} \put(0.00,2.00){\line(1,0){9.00}}
\put(0.00,4.00){\line(1,0){8.00}}
\put(0.00,6.00){\line(1,0){7.00}}
\put(0.00,8.00){\line(1,0){5.89}}
\put(0.00,10.00){\line(1,0){4.89}}
\put(0.00,12.00){\line(1,0){4.00}} \thicklines{}
\put(10.00,0.00){\line(-3,1){10.00}}
\put(10.00,0.00){\line(-1,1){10.00}} \thinlines{}
\linethickness{0.05pt} \put(0.00,10.00){\line(0,-1){6.58}}
\put(1.00,8.92){\line(0,-1){5.83}}
\put(2.00,8.00){\line(0,-1){5.33}}
\put(3.00,6.92){\line(0,-1){4.58}}
\put(4.00,6.00){\line(0,-1){4.00}}
\put(5.00,4.92){\line(0,-1){3.25}}
\put(6.00,4.00){\line(0,-1){2.67}}
\put(7.00,2.92){\line(0,-1){1.92}}
\put(8.00,2.00){\line(0,-1){1.33}}
\put(9.00,0.92){\line(0,-1){0.58}}
\end{picture}
\hfill
%TexCad Options
%\grade{\off}
%\emlines{\off}
%\beziermacro{\off}
%\reduce{\on}
%\snapping{\off}
%\quality{2.00}
%\graddiff{0.01}
%\snapasp{1}
%\zoom{4.00}
\unitlength 0.70ex \linethickness{0.4pt}
\begin{picture}(20.00,14.00)
\put(0.00,0.00){\line(1,0){20.00}}
\put(10.00,5.00){\circle*{1.00}}
\put(10.00,0.00){\line(-1,2){7.00}}
\put(16.00,9.00){\makebox(0,0)[cc]{$\text{B}_{+}$}}
\linethickness{0.05pt} \put(0.00,2.00){\line(1,0){9.00}}
\put(0.00,4.00){\line(1,0){8.00}}
\put(0.00,6.00){\line(1,0){7.00}}
\put(0.00,8.00){\line(1,0){5.89}}
\put(0.00,10.00){\line(1,0){4.89}}
\put(0.00,12.00){\line(1,0){4.00}} \thicklines{}
\put(10.00,0.00){\line(-1,1){10.00}}
\put(10.00,5.50){\line(0,1){8.50}}
\put(10.00,0.00){\line(0,1){4.50}} \thinlines{}
\linethickness{0.05pt} \put(9.00,14.00){\line(0,-1){13.00}}
\put(8.00,14.00){\line(0,-1){12.00}}
\put(7.00,14.00){\line(0,-1){11.08}}
\put(6.00,14.00){\line(0,-1){10.00}}
\put(5.00,14.00){\line(0,-1){9.00}}
\put(4.00,14.00){\line(0,-1){8.00}}
\put(3.00,14.00){\line(0,-1){7.00}}
\put(2.00,14.00){\line(0,-1){6.00}}
\put(1.00,14.00){\line(0,-1){5.00}}
\put(0.00,14.00){\line(0,-1){4.00}}
\end{picture}
\bigskip

\raisebox{7ex}{$\text{A}_l$}\hfill
%TexCad Options
%\grade{\off}
%\emlines{\off}
%\beziermacro{\off}
%\reduce{\on}
%\snapping{\off}
%\quality{2.00}
%\graddiff{0.01}
%\snapasp{1}
%\zoom{4.00}
\unitlength 0.70ex \linethickness{0.4pt}
\begin{picture}(20.00,20.00)
\put(10.00,6.00){\line(1,0){10.00}}
\put(10.00,11.00){\circle{1.00}}
\put(10.00,11.50){\line(0,1){8.50}}
\put(10.00,6.00){\line(0,1){4.50}}
\put(10.00,6.00){\line(-1,2){7.00}} \linethickness{0.05pt}
\put(0.00,8.00){\line(1,0){9.00}}
\put(0.00,10.00){\line(1,0){8.00}}
\put(0.00,12.00){\line(1,0){7.00}}
\put(0.00,14.00){\line(1,0){5.89}}
\put(0.00,16.00){\line(1,0){4.89}}
\put(0.00,18.00){\line(1,0){4.00}} \thicklines{}
\put(0.00,6.00){\line(1,0){10.00}}
\put(10.00,6.00){\line(-1,1){10.00}} \thinlines{}
\linethickness{0.05pt} \put(0.00,16.00){\line(0,-1){10.00}}
\put(1.00,14.92){\line(0,-1){8.92}}
\put(2.00,14.00){\line(0,-1){8.00}}
\put(3.00,12.92){\line(0,-1){6.92}}
\put(4.00,12.00){\line(0,-1){6.00}}
\put(5.00,10.92){\line(0,-1){4.92}}
\put(6.00,10.00){\line(0,-1){4.00}}
\put(7.00,8.92){\line(0,-1){2.92}}
\put(8.00,8.00){\line(0,-1){2.00}}
\put(9.00,6.92){\line(0,-1){0.92}} \linethickness{0.2pt}
\put(9.50,11.00){\line(-1,0){4.50}}
\put(5.00,11.00){\line(0,-1){5.00}}
\put(5.00,6.00){\circle*{0.50}}
\put(10.00,2.00){\makebox(0,0)[ct]{\strut$x_i$, $i=1,\dots,l$}}
\put(5.00,5.00){\makebox(0,0)[ct]{$\ell_i$}}
\end{picture}
\hfill
%TexCad Options
%\grade{\off}
%\emlines{\off}
%\beziermacro{\off}
%\reduce{\on}
%\snapping{\off}
%\quality{2.00}
%\graddiff{0.01}
%\snapasp{1}
%\zoom{4.00}
\unitlength 0.70ex \linethickness{0.4pt}
\begin{picture}(20.00,20.00)
\put(10.00,6.00){\line(1,0){10.00}}
\put(10.00,6.00){\line(0,1){14.00}}
\put(15.00,11.00){\circle*{1.00}} \linethickness{0.05pt}
\put(0.00,8.00){\line(1,0){11.00}}
\put(0.00,10.00){\line(1,0){12.00}}
\put(0.00,12.00){\line(1,0){13.00}}
\put(0.00,14.00){\line(1,0){14.00}}
\put(0.00,16.00){\line(1,0){15.00}}
\put(0.00,18.00){\line(1,0){16.00}} \thicklines{}
\put(0.00,6.00){\line(1,0){10.00}}
\put(10.00,6.00){\line(1,1){10.00}} \thinlines{}
\linethickness{0.05pt} \put(0.00,20.00){\line(0,-1){14.00}}
\put(1.00,20.00){\line(0,-1){14.00}}
\put(2.00,20.00){\line(0,-1){14.00}}
\put(3.00,20.00){\line(0,-1){14.00}}
\put(4.00,20.00){\line(0,-1){14.00}}
\put(5.00,20.00){\line(0,-1){14.00}}
\put(6.00,20.00){\line(0,-1){14.00}}
\put(7.00,20.00){\line(0,-1){14.00}}
\put(8.00,20.00){\line(0,-1){14.00}}
\put(9.00,20.00){\line(0,-1){14.00}}
\put(10.00,2.00){\makebox(0,0)[ct]{\strut$x=\infty$}}
\put(10.00,6.00){\line(1,2){7.00}}
\put(11.00,20.00){\line(0,-1){13.00}}
\put(12.00,20.00){\line(0,-1){12.00}}
\put(13.00,20.00){\line(0,-1){11.00}}
\put(14.00,20.00){\line(0,-1){10.00}}
\put(15.00,20.00){\line(0,-1){9.00}}
\put(16.00,20.00){\line(0,-1){8.00}}
\put(17.00,20.00){\line(0,-1){7.00}}
\put(18.00,20.00){\line(0,-1){6.00}}
\put(19.00,20.00){\line(0,-1){5.00}}
\end{picture}
\hfill
%TexCad Options
%\grade{\off}
%\emlines{\off}
%\beziermacro{\off}
%\reduce{\on}
%\snapping{\off}
%\quality{2.00}
%\graddiff{0.01}
%\snapasp{1}
%\zoom{4.00}
\unitlength 0.70ex \linethickness{0.4pt}
\begin{picture}(20.00,20.00)
\put(10.00,6.00){\line(1,0){10.00}}
\put(10.00,11.00){\circle*{1.00}}
\put(10.00,6.00){\line(-1,2){7.00}} \linethickness{0.05pt}
\put(0.00,8.00){\line(1,0){9.00}}
\put(0.00,10.00){\line(1,0){8.00}}
\put(0.00,12.00){\line(1,0){7.00}}
\put(0.00,14.00){\line(1,0){5.89}}
\put(0.00,16.00){\line(1,0){4.89}}
\put(0.00,18.00){\line(1,0){4.00}} \thicklines{}
\put(0.00,6.00){\line(1,0){10.00}}
\put(10.00,6.00){\line(0,1){14.00}} \thinlines{}
\linethickness{0.05pt} \put(0.00,20.00){\line(0,-1){14.00}}
\put(1.00,20.00){\line(0,-1){14.00}}
\put(2.00,20.00){\line(0,-1){14.00}}
\put(3.00,20.00){\line(0,-1){14.00}}
\put(4.00,20.00){\line(0,-1){14.00}}
\put(5.00,20.00){\line(0,-1){14.00}}
\put(6.00,20.00){\line(0,-1){14.00}}
\put(7.00,20.00){\line(0,-1){14.00}}
\put(8.00,20.00){\line(0,-1){14.00}}
\put(9.00,20.00){\line(0,-1){14.00}}
\put(10.00,2.00){\makebox(0,0)[ct]{\strut other~$x$}}
\end{picture}
\bigskip

\raisebox{7ex}{$\text{B}_{-}$}\hfill
%TexCad Options
%\grade{\off}
%\emlines{\off}
%\beziermacro{\off}
%\reduce{\on}
%\snapping{\off}
%\quality{2.00}
%\graddiff{0.01}
%\snapasp{1}
%\zoom{4.00}
\unitlength 0.70ex \linethickness{0.4pt}
\begin{picture}(20.00,20.00)
\put(0.00,6.00){\line(1,0){10.00}}
\put(10.00,11.00){\circle{1.00}}
\put(10.00,11.50){\line(0,1){8.50}}
\put(10.00,6.00){\line(0,1){4.50}}
\put(10.00,6.00){\line(-1,2){7.00}} \linethickness{0.05pt}
\put(0.00,8.00){\line(1,0){9.00}}
\put(0.00,10.00){\line(1,0){8.00}}
\put(0.00,12.00){\line(1,0){7.00}}
\put(0.00,14.00){\line(1,0){5.89}}
\put(0.00,16.00){\line(1,0){4.89}}
\put(0.00,18.00){\line(1,0){4.00}} \thicklines{}
\put(10.00,6.00){\line(1,0){10.00}}
\put(10.00,6.00){\line(-1,1){10.00}} \thinlines{}
\linethickness{0.2pt} \put(9.50,11.00){\line(-1,0){4.50}}
\put(5.00,11.00){\line(0,-1){5.00}}
\put(5.00,6.00){\circle*{0.50}} \linethickness{0.05pt}
\put(10.00,2.00){\makebox(0,0)[ct]{\strut$x_1$}}
\put(5.00,5.00){\makebox(0,0)[ct]{$\ell_1$}}
\put(20.00,20.00){\line(0,-1){14.00}}
\put(19.00,20.00){\line(0,-1){14.00}}
\put(18.00,20.00){\line(0,-1){14.00}}
\put(17.00,20.00){\line(0,-1){14.00}}
\put(16.00,20.00){\line(0,-1){14.00}}
\put(15.00,20.00){\line(0,-1){14.00}}
\put(14.00,20.00){\line(0,-1){14.00}}
\put(13.00,20.00){\line(0,-1){14.00}}
\put(12.00,20.00){\line(0,-1){14.00}}
\put(11.00,20.00){\line(0,-1){14.00}}
\put(9.00,20.00){\line(0,-1){13.00}}
\put(8.00,20.00){\line(0,-1){12.00}}
\put(7.00,20.00){\line(0,-1){11.00}}
\put(6.00,20.00){\line(0,-1){10.00}}
\put(5.00,20.00){\line(0,-1){9.00}}
\put(4.00,20.00){\line(0,-1){8.00}}
\put(3.00,20.00){\line(0,-1){7.00}}
\put(2.00,20.00){\line(0,-1){6.00}}
\put(1.00,20.00){\line(0,-1){5.00}}
\end{picture}
\hfill
%TexCad Options
%\grade{\off}
%\emlines{\off}
%\beziermacro{\off}
%\reduce{\on}
%\snapping{\off}
%\quality{2.00}
%\graddiff{0.01}
%\snapasp{1}
%\zoom{4.00}
\unitlength 0.70ex \linethickness{0.4pt}
\begin{picture}(20.00,20.00)
\put(0.00,6.00){\line(1,0){10.00}}
\put(10.00,6.00){\line(0,1){14.00}}
\put(15.00,11.00){\circle*{1.00}}
\put(10.00,6.00){\line(1,2){7.00}} \linethickness{0.05pt}
\put(0.00,8.00){\line(1,0){11.00}}
\put(0.00,10.00){\line(1,0){12.00}}
\put(0.00,12.00){\line(1,0){13.00}}
\put(0.00,14.00){\line(1,0){14.00}}
\put(0.00,16.00){\line(1,0){15.00}}
\put(0.00,18.00){\line(1,0){16.00}} \thicklines{}
\put(10.00,6.00){\line(1,0){10.00}}
\put(10.00,6.00){\line(1,1){10.00}} \thinlines{}
\linethickness{0.05pt} \put(20.00,16.00){\line(0,-1){10.00}}
\put(19.00,14.92){\line(0,-1){8.92}}
\put(18.00,14.00){\line(0,-1){8.00}}
\put(17.00,12.92){\line(0,-1){6.92}}
\put(16.00,12.00){\line(0,-1){6.00}}
\put(15.00,10.92){\line(0,-1){4.92}}
\put(14.00,10.00){\line(0,-1){4.00}}
\put(13.00,8.92){\line(0,-1){2.92}}
\put(12.00,8.00){\line(0,-1){2.00}}
\put(11.00,6.92){\line(0,-1){0.92}}
\put(10.00,2.00){\makebox(0,0)[ct]{\strut$x=\infty$}}
\end{picture}
\hfill
%TexCad Options
%\grade{\off}
%\emlines{\off}
%\beziermacro{\off}
%\reduce{\on}
%\snapping{\off}
%\quality{2.00}
%\graddiff{0.01}
%\snapasp{1}
%\zoom{4.00}
\unitlength 0.70ex \linethickness{0.4pt}
\begin{picture}(20.00,20.00)
\put(0.00,6.00){\line(1,0){10.00}}
\put(10.00,11.00){\circle*{1.00}}
\put(10.00,6.00){\line(-1,2){7.00}} \linethickness{0.05pt}
\put(0.00,8.00){\line(1,0){9.00}}
\put(0.00,10.00){\line(1,0){8.00}}
\put(0.00,12.00){\line(1,0){7.00}}
\put(0.00,14.00){\line(1,0){5.89}}
\put(0.00,16.00){\line(1,0){4.89}}
\put(0.00,18.00){\line(1,0){4.00}} \thicklines{}
\put(10.00,6.00){\line(1,0){10.00}}
\put(10.00,6.00){\line(0,1){14.00}} \thinlines{}
\linethickness{0.05pt} \put(11.00,20.00){\line(0,-1){14.00}}
\put(12.00,20.00){\line(0,-1){14.00}}
\put(13.00,20.00){\line(0,-1){14.00}}
\put(14.00,20.00){\line(0,-1){14.00}}
\put(15.00,20.00){\line(0,-1){14.00}}
\put(16.00,20.00){\line(0,-1){14.00}}
\put(17.00,20.00){\line(0,-1){14.00}}
\put(18.00,20.00){\line(0,-1){14.00}}
\put(19.00,20.00){\line(0,-1){14.00}}
\put(20.00,20.00){\line(0,-1){14.00}}
\put(10.00,2.00){\makebox(0,0)[ct]{\strut other~$x$}}
\end{picture}
\bigskip

\raisebox{7ex}{$\text{B}_0$}\hfill
%TexCad Options
%\grade{\off}
%\emlines{\off}
%\beziermacro{\off}
%\reduce{\on}
%\snapping{\off}
%\quality{2.00}
%\graddiff{0.01}
%\snapasp{1}
%\zoom{4.00}
\unitlength 0.70ex \linethickness{0.4pt}
\begin{picture}(20.00,20.00)
\put(0.00,6.00){\line(1,0){10.00}}
\put(10.00,11.00){\circle*{1.00}}
\put(10.00,11.50){\line(0,1){8.50}}
\put(10.00,6.00){\line(0,1){4.50}}
\put(10.00,6.00){\line(-1,2){7.00}} \linethickness{0.05pt}
\put(0.00,8.00){\line(1,0){9.00}}
\put(0.00,10.00){\line(1,0){8.00}}
\put(0.00,12.00){\line(1,0){7.00}}
\put(0.00,14.00){\line(1,0){5.89}}
\put(0.00,16.00){\line(1,0){4.89}}
\put(0.00,18.00){\line(1,0){4.00}} \thicklines{}
\put(10.00,6.00){\line(1,0){10.00}}
\put(10.00,6.00){\line(-1,1){10.00}} \thinlines{}
\linethickness{0.2pt} \put(9.50,11.00){\line(-1,0){4.50}}
\put(5.00,11.00){\line(0,-1){5.00}}
\put(5.00,6.00){\circle*{0.50}} \linethickness{0.05pt}
\put(10.00,2.00){\makebox(0,0)[ct]{\strut$x_1$}}
\put(5.00,5.00){\makebox(0,0)[ct]{$\ell_1$}}
\put(20.00,20.00){\line(0,-1){14.00}}
\put(19.00,20.00){\line(0,-1){14.00}}
\put(18.00,20.00){\line(0,-1){14.00}}
\put(17.00,20.00){\line(0,-1){14.00}}
\put(16.00,20.00){\line(0,-1){14.00}}
\put(15.00,20.00){\line(0,-1){14.00}}
\put(14.00,20.00){\line(0,-1){14.00}}
\put(13.00,20.00){\line(0,-1){14.00}}
\put(12.00,20.00){\line(0,-1){14.00}}
\put(11.00,20.00){\line(0,-1){14.00}}
\put(9.00,20.00){\line(0,-1){13.00}}
\put(8.00,20.00){\line(0,-1){12.00}}
\put(7.00,20.00){\line(0,-1){11.00}}
\put(6.00,20.00){\line(0,-1){10.00}}
\put(5.00,20.00){\line(0,-1){9.00}}
\put(4.00,20.00){\line(0,-1){8.00}}
\put(3.00,20.00){\line(0,-1){7.00}}
\put(2.00,20.00){\line(0,-1){6.00}}
\put(1.00,20.00){\line(0,-1){5.00}}
\end{picture}
\hfill
%TexCad Options
%\grade{\off}
%\emlines{\off}
%\beziermacro{\off}
%\reduce{\on}
%\snapping{\off}
%\quality{2.00}
%\graddiff{0.01}
%\snapasp{1}
%\zoom{4.00}
\unitlength 0.70ex \linethickness{0.4pt}
\begin{picture}(20.00,20.00)
\put(0.00,6.00){\line(1,0){10.00}}
\put(10.00,6.00){\line(0,1){14.00}}
\put(15.00,11.00){\circle*{1.00}}
\put(10.00,6.00){\line(1,2){7.00}} \linethickness{0.05pt}
\put(0.00,8.00){\line(1,0){11.00}}
\put(0.00,10.00){\line(1,0){12.00}}
\put(0.00,12.00){\line(1,0){13.00}}
\put(0.00,14.00){\line(1,0){14.00}}
\put(0.00,16.00){\line(1,0){15.00}}
\put(0.00,18.00){\line(1,0){16.00}} \thicklines{}
\put(10.00,6.00){\line(1,0){10.00}}
\put(10.00,6.00){\line(1,1){10.00}} \thinlines{}
\linethickness{0.05pt} \put(20.00,16.00){\line(0,-1){10.00}}
\put(19.00,14.92){\line(0,-1){8.92}}
\put(18.00,14.00){\line(0,-1){8.00}}
\put(17.00,12.92){\line(0,-1){6.92}}
\put(16.00,12.00){\line(0,-1){6.00}}
\put(15.00,10.92){\line(0,-1){4.92}}
\put(14.00,10.00){\line(0,-1){4.00}}
\put(13.00,8.92){\line(0,-1){2.92}}
\put(12.00,8.00){\line(0,-1){2.00}}
\put(11.00,6.92){\line(0,-1){0.92}}
\put(10.00,2.00){\makebox(0,0)[ct]{\strut$x=\infty$}}
\end{picture}
\hfill
%TexCad Options
%\grade{\off}
%\emlines{\off}
%\beziermacro{\off}
%\reduce{\on}
%\snapping{\off}
%\quality{2.00}
%\graddiff{0.01}
%\snapasp{1}
%\zoom{4.00}
\unitlength 0.70ex \linethickness{0.4pt}
\begin{picture}(20.00,20.00)
\put(0.00,6.00){\line(1,0){10.00}}
\put(10.00,11.00){\circle*{1.00}}
\put(10.00,6.00){\line(-1,2){7.00}} \linethickness{0.05pt}
\put(0.00,8.00){\line(1,0){9.00}}
\put(0.00,10.00){\line(1,0){8.00}}
\put(0.00,12.00){\line(1,0){7.00}}
\put(0.00,14.00){\line(1,0){5.89}}
\put(0.00,16.00){\line(1,0){4.89}}
\put(0.00,18.00){\line(1,0){4.00}} \thicklines{}
\put(10.00,6.00){\line(1,0){10.00}}
\put(10.00,6.00){\line(0,1){14.00}} \thinlines{}
\linethickness{0.05pt} \put(11.00,20.00){\line(0,-1){14.00}}
\put(12.00,20.00){\line(0,-1){14.00}}
\put(13.00,20.00){\line(0,-1){14.00}}
\put(14.00,20.00){\line(0,-1){14.00}}
\put(15.00,20.00){\line(0,-1){14.00}}
\put(16.00,20.00){\line(0,-1){14.00}}
\put(17.00,20.00){\line(0,-1){14.00}}
\put(18.00,20.00){\line(0,-1){14.00}}
\put(19.00,20.00){\line(0,-1){14.00}}
\put(20.00,20.00){\line(0,-1){14.00}}
\put(10.00,2.00){\makebox(0,0)[ct]{\strut other~$x$}}
\end{picture}
\bigskip
\end{center}
\end{figure}
Up to a ``change of coordinates'' (Remark~\ref{split}), we may
assume that for germs of types $\text{C},\text{AB},\text{B}_{+}$
the colored cone lies in $\EHS{x}$, $x\neq\infty$, and for germs
of types $\text{A}_l,\text{B}_{-},\text{B}_0$, $x_i\neq\infty$.
Moreover, for the hypercone to be strictly convex, we must have
$\sum\ell_i<-1$ for~$\text{A}_l$, $l\leq1$, and $\ell_1>-1$ for
$\text{B}_{-}$ and~$\text{B}_0$.

Affine $\SL_2$-embeddings correspond to minimal $B$-charts of
$G$-germs of type~$\text{B}_0$. They were first classified by
Popov~\cite{SL2.aff}. Embeddings of $\SL_2/H$, where $H$ is
finite, were classified in~\cite{SL2/fin}. Embeddings of $G/H$,
where $G$ has semisimple rank~$1$ and $H$ is finite, were
classified in~\cite[\S5]{LVT}.
\end{example}

\begin{example}[ordered triangles]\label{SL3/T}
Suppose $G=\SL_3(\kk)$, $H=T$ is the diagonal torus. Then
$\HS=G/H$ is the space of ordered triangles on a projective plane.
The standard Borel subgroup $B$ consists of upper-triangular
matrices $g=(g_{ij})$, $\det{g}=1$, $g_{ij}=0$ for $i>j$. As
usual, $\eps_i(g)=g_{ii}$ are the tautological weights of~$T$,
$\omega_1=\eps_1$, $\omega_2=\eps_1+\eps_2$ are the fundamental
weights, and $\alpha_1=\eps_1-\eps_2$, $\alpha_2=\eps_2-\eps_3$
are the simple roots. Denote by $\omega^{\vee}_i$ the fundamental
coweights, and let $\rho=\omega_1+\omega_2=\alpha_1+\alpha_2$.

The subregular $B$-divisors $D_i,\widetilde{D}_i$ are defined by
the $(B\times H)$-eigen\-func\-tions $\eta_i(g)=g_{3i}$,
$\widetilde{\eta}_i(g)=\left|
\begin{smallmatrix}
g_{2j} & g_{2k} \\
g_{3j} & g_{3k}
\end{smallmatrix}\right|$ of biweights $(\omega_2,\eps_i)$,
$(\omega_1,-\eps_i)$. $\widetilde{D}_i$~consists of triangles
whose \ordinal{$i$} side contains the $B$-fixed point in~$\PP^2$,
and $D_i$ consists of triangles whose \ordinal{$i$} vertex lies on
the $B$-fixed line.

The functions $\eta_i\widetilde{\eta}_i$ generate the
two-dimensional subspace $M=\kk[G]^{(B\times H)}_{(\rho,0)}$,
$\eta_1\widetilde{\eta}_1+\eta_2\widetilde{\eta}_2+
\eta_3\widetilde{\eta}_3=0$. Let $x_i\in\PP^1=\PP(M^{*})$ be the
points corresponding to $\eta_i\widetilde{\eta}_i$. The regular
$B$-divisors $D_x$, $x\neq x_1,x_2,x_3$, are defined by equations
$\eta_x=\alpha_1\eta_1\widetilde{\eta}_1+
\alpha_2\eta_2\widetilde{\eta}_2+\alpha_3\eta_3\widetilde{\eta}_3$.

The group $\RG=\langle\alpha_1,\alpha_2\rangle$ is the root
lattice, $\ef{\alpha_1}=
\widetilde{\eta}_1\widetilde{\eta}_2\widetilde{\eta}_3/
\eta_{\infty}$, $\ef{\alpha_2}=\eta_1\eta_2\eta_3/\eta_{\infty}$,
where $\infty\in\PP^1$ is a certain fixed point.

By~\ref{form.curv}, any $G$-valuation corresponding to a divisor
with dense $G$-orbit is proportional to $v_{x(t)}$, where
%*
\begin{equation*}
x(t)=\begin{pmatrix}
1 & t^m & u(t) \\
0 &  1  & t^n  \\
0 &  0  &  1
\end{pmatrix},
\end{equation*}
%*
and we may assume $m,n,r=\ord_tu(t)\leq0$. Computing the values
$v_{x(t)}(\eta_x)$ as in Example~\ref{SL2/e}, one finds that the
set of $G$-valuations $v=(h,\ell)\in\EHS{x}$ corresponding to
divisors with dense $G$-orbit is determined by the inequalities
$a_1,a_2\leq0\leq h$ (if $x=x_i$) or $a_1,a_2\leq-2h\leq0$ (if
$x=\infty$) or $a_1,a_2\leq-h\leq0$ (otherwise), and $h=0\implies
a_1\text{ or }a_2=0$, where
$\ell=a_1\omega^{\vee}_1+a_2\omega^{\vee}_2$. Hence $\Vv_x$ is
determined by the same inequalities without any restrictions for
$h=0$.

The colored equipment is represented in Figure~\ref{hs(trk)}. (The
intersections of $\Vv$ with the hyperplane sections $\ES= \{h=0\}$
and $\{h=1\}$ of $\EHS{x}$ are hatched.)
\begin{figure}[h]
\caption{}\label{hs(trk)}
\begin{center}
\raisebox{7ex}{$h=1$}\qquad
%TexCad Options
%\grade{\on}
%\emlines{\off}
%\beziermacro{\off}
%\reduce{\on}
%\snapping{\off}
%\quality{2.00}
%\graddiff{0.01}
%\snapasp{1}
%\zoom{4.00}
\unitlength 0.50ex \linethickness{0.4pt}
\begin{picture}(39.84,44.76)
\put(19.92,24.88){\oval(39.83,39.75)[]}
\put(20.00,25.00){\circle*{1.00}} \put(25.50,28.00){\circle{1.00}}
\put(25.50,22.00){\circle{1.00}} \thicklines{}
%\emline(20.00,25.00)(0.75,14.50)
\multiput(20.00,25.00)(-0.22,-0.12){88}{\line(-1,0){0.22}}
%\end
%\emline(20.00,25.00)(0.75,35.50)
\multiput(20.00,25.00)(-0.22,0.12){88}{\line(-1,0){0.22}}
%\end
%\emline(20.00,25.00)(25.08,27.75)
\multiput(20.00,25.00)(0.22,0.12){23}{\line(1,0){0.22}}
%\end
%\emline(20.00,25.00)(25.08,22.25)
\multiput(20.00,25.00)(0.22,-0.12){23}{\line(1,0){0.22}}
%\end
%\emline(25.75,28.42)(34.25,44.25)
\multiput(25.75,28.42)(0.12,0.22){71}{\line(0,1){0.22}}
%\end
%\emline(25.75,21.58)(34.25,5.75)
\multiput(25.75,21.58)(0.12,-0.22){71}{\line(0,-1){0.22}}
%\end
%\bezier{32}(23.42,23.17)(22.42,19.25)(26.75,19.58)
\qbezier(23.42,23.17)(22.42,19.25)(26.75,19.58)
%\put(23.42,23.17){\line(0,-1){2.03}}
%\multiput(23.31,21.13)(0.12,-0.15){8}{\line(0,-1){0.15}}
%\multiput(24.25,19.93)(0.83,-0.12){3}{\line(1,0){0.83}}
%\end
%\bezier{32}(23.42,26.83)(22.42,30.75)(26.75,30.42)
\qbezier(23.42,26.83)(22.42,30.75)(26.75,30.42)
%\put(23.42,26.83){\line(0,1){2.03}}
%\multiput(23.31,28.87)(0.12,0.15){8}{\line(0,1){0.15}}
%\multiput(24.25,30.07)(0.83,0.12){3}{\line(1,0){0.83}}
%\end
\thinlines{} \linethickness{0.05pt}
%\emline(0.00,35.00)(1.58,35.00)
\put(0.00,35.00){\line(1,0){1.58}}
%\end
%\emline(0.00,33.00)(5.25,33.00)
\put(0.00,33.00){\line(1,0){5.25}}
%\end
%\emline(0.00,31.00)(9.00,31.00)
\put(0.00,31.00){\line(1,0){9.00}}
%\end
%\emline(0.00,29.00)(12.58,29.00)
\put(0.00,29.00){\line(1,0){12.58}}
%\end
%\emline(0.00,27.00)(16.25,27.00)
\put(0.00,27.00){\line(1,0){16.25}}
%\end
%\emline(0.00,25.00)(20.00,25.00)
\put(0.00,25.00){\line(1,0){20.00}}
%\end
%\emline(0.00,15.00)(1.58,15.00)
\put(0.00,15.00){\line(1,0){1.58}}
%\end
%\emline(0.00,17.00)(5.25,17.00)
\put(0.00,17.00){\line(1,0){5.25}}
%\end
%\emline(0.00,19.00)(9.00,19.00)
\put(0.00,19.00){\line(1,0){9.00}}
%\end
%\emline(0.00,21.00)(12.58,21.00)
\put(0.00,21.00){\line(1,0){12.58}}
%\end
%\emline(0.00,23.00)(16.25,23.00)
\put(0.00,23.00){\line(1,0){16.25}}
%\end
\put(6.00,25.00){\makebox(0,0)[cc]{$Y$}}
\put(28.00,21.00){\makebox(0,0)[lc]{$\widetilde{D}_i$}}
\put(28.00,29.00){\makebox(0,0)[lc]{$D_i$}}
\put(20.00,3.00){\makebox(0,0)[ct]{\strut$x_i$, $i=1,2,3$}}
\put(6.00,16.00){\makebox(0,0)[lt]{$Y_i$}}
\put(6.00,34.00){\makebox(0,0)[lb]{$\widetilde{Y}_i$}}
\put(20.00,27.00){\makebox(0,0)[cb]{$W_i$}}
%\emline(1.00,14.67)(1.00,35.33)
\put(1.00,14.67){\line(0,1){20.66}}
%\end
%\emline(2.00,34.75)(2.00,15.17)
\put(2.00,34.75){\line(0,-1){19.58}}
%\end
%\emline(3.00,15.67)(3.00,34.33)
\put(3.00,15.67){\line(0,1){18.66}}
%\end
%\emline(4.00,33.67)(4.00,16.33)
\put(4.00,33.67){\line(0,-1){17.34}}
%\end
%\emline(5.00,16.83)(5.00,33.17)
\put(5.00,16.83){\line(0,1){16.34}}
%\end
%\emline(6.00,32.58)(6.00,17.33)
\put(6.00,32.58){\line(0,-1){15.25}}
%\end
%\emline(7.00,17.92)(7.00,32.08)
\put(7.00,17.92){\line(0,1){14.16}}
%\end
%\emline(8.00,31.50)(8.00,18.42)
\put(8.00,31.50){\line(0,-1){13.08}}
%\end
%\emline(9.00,19.00)(9.00,31.00)
\put(9.00,19.00){\line(0,1){12.00}}
%\end
%\emline(10.00,30.42)(10.00,19.58)
\put(10.00,30.42){\line(0,-1){10.84}}
%\end
%\emline(11.00,20.08)(11.00,29.92)
\put(11.00,20.08){\line(0,1){9.84}}
%\end
%\emline(12.00,29.33)(12.00,20.58)
\put(12.00,29.33){\line(0,-1){8.75}}
%\end
%\emline(13.00,21.17)(13.00,28.75)
\put(13.00,21.17){\line(0,1){7.58}}
%\end
%\emline(14.00,28.33)(14.00,21.67)
\put(14.00,28.33){\line(0,-1){6.66}}
%\end
%\emline(15.00,22.33)(15.00,27.67)
\put(15.00,22.33){\line(0,1){5.34}}
%\end
%\emline(16.00,27.17)(16.00,22.83)
\put(16.00,27.17){\line(0,-1){4.34}}
%\end
%\emline(17.00,23.42)(17.00,26.67)
\put(17.00,23.42){\line(0,1){3.25}}
%\end
%\emline(18.00,26.08)(18.00,23.92)
\put(18.00,26.08){\line(0,-1){2.16}}
%\end
%\emline(19.00,24.50)(19.00,25.50)
\put(19.00,24.50){\line(0,1){1.00}}
%\end
\put(20.00,40.00){\makebox(0,0)[cc]{$\Cc^{*}_{x_i}$}}
\put(20.00,10.00){\makebox(0,0)[cc]{$\Cc_{x_i}$}}
\end{picture}
\hfill
%TexCad Options
%\grade{\on}
%\emlines{\off}
%\beziermacro{\off}
%\reduce{\on}
%\snapping{\off}
%\quality{2.00}
%\graddiff{0.01}
%\snapasp{1}
%\zoom{4.00}
\unitlength 0.50ex \linethickness{0.4pt}
\begin{picture}(39.84,44.76)
\put(19.92,24.88){\oval(39.83,39.75)[]}
\put(13.00,25.00){\circle*{0.50}}
\put(18.50,28.00){\circle*{0.50}}
\put(18.50,22.00){\circle*{0.50}}
\put(29.50,28.00){\circle*{0.70}}
\put(29.50,22.00){\circle*{0.70}}
\put(35.00,25.00){\circle*{0.70}}
%\emline(13.00,25.00)(0.00,32.08)
\multiput(13.00,25.00)(-0.22,0.12){60}{\line(-1,0){0.22}}
%\end
%\emline(13.00,25.00)(0.00,17.92)
\multiput(13.00,25.00)(-0.22,-0.12){60}{\line(-1,0){0.22}}
%\end
\thicklines{}
%\emline(23.58,25.25)(0.00,37.83)
\multiput(23.58,25.25)(-0.22,0.12){105}{\line(-1,0){0.22}}
%\end
%\emline(24.25,25.42)(34.42,43.92)
\multiput(24.25,25.42)(0.12,0.22){85}{\line(0,1){0.22}}
%\end
%\emline(23.58,24.75)(0.00,12.17)
\multiput(23.58,24.75)(-0.22,-0.12){105}{\line(-1,0){0.22}}
%\end
%\emline(24.25,24.58)(34.42,6.08)
\multiput(24.25,24.58)(0.12,-0.22){85}{\line(0,-1){0.22}}
%\end
%\bezier{28}(22.00,26.00)(21.67,22.50)(25.08,23.00)
\qbezier(22.00,26.00)(21.67,22.50)(25.08,23.00)
%\multiput(22.00,26.00)(0.08,-0.66){3}{\line(0,-1){0.66}}
%\multiput(22.24,24.01)(0.32,-0.11){9}{\line(1,0){0.32}}
%\end
%\bezier{48}(21.00,23.42)(19.42,28.67)(25.58,27.92)
\qbezier(21.00,23.42)(19.42,28.67)(25.58,27.92)
%\multiput(21.00,23.42)(-0.11,0.64){3}{\line(0,1){0.64}}
%\multiput(20.68,25.35)(0.12,0.47){3}{\line(0,1){0.47}}
%\multiput(21.03,26.75)(0.13,0.11){8}{\line(1,0){0.13}}
%\multiput(22.05,27.64)(1.18,0.09){3}{\line(1,0){1.18}}
%\end
\thinlines{} \linethickness{0.2pt}
%\emline(39.83,27.67)(2.08,7.17)
\multiput(39.83,27.67)(-0.22,-0.12){171}{\line(-1,0){0.22}}
%\end
%\emline(39.83,22.33)(2.08,42.83)
\multiput(39.83,22.33)(-0.22,0.12){171}{\line(-1,0){0.22}}
%\end
\linethickness{0.05pt}
%\emline(0.00,29.00)(5.58,29.00)
\put(0.00,29.00){\line(1,0){5.58}}
%\end
%\emline(0.00,27.00)(9.25,27.00)
\put(0.00,27.00){\line(1,0){9.25}}
%\end
%\emline(0.00,25.00)(13.00,25.00)
\put(0.00,25.00){\line(1,0){13.00}}
%\end
%\emline(0.00,21.00)(5.58,21.00)
\put(0.00,21.00){\line(1,0){5.58}}
%\end
%\emline(0.00,23.00)(9.25,23.00)
\put(0.00,23.00){\line(1,0){9.25}}
%\end
\put(19.00,25.17){\makebox(0,0)[rc]{$Y$}}
\put(20.00,3.00){\makebox(0,0)[ct]{\strut$x=\infty$}}
\put(24.00,25.00){\circle*{1.00}}
\put(26.00,25.00){\makebox(0,0)[lc]{$D_{\infty}$}}
%\emline(0.00,31.00)(1.92,31.00)
\put(0.00,31.00){\line(1,0){1.92}}
%\end
%\emline(0.00,19.00)(1.92,19.00)
\put(0.00,19.00){\line(1,0){1.92}}
%\end
\put(20.00,40.00){\makebox(0,0)[cc]{$\Cc^{*}_{\infty}$}}
\put(20.00,10.00){\makebox(0,0)[cc]{$\Cc_{\infty}$}}
%\emline(4.92,14.67)(4.92,35.33)
\put(4.92,14.67){\line(0,1){20.66}}
%\end
%\emline(5.92,34.75)(5.92,15.17)
\put(5.92,34.75){\line(0,-1){19.58}}
%\end
%\emline(6.92,15.67)(6.92,34.33)
\put(6.92,15.67){\line(0,1){18.66}}
%\end
%\emline(7.92,33.67)(7.92,16.33)
\put(7.92,33.67){\line(0,-1){17.34}}
%\end
%\emline(8.92,16.83)(8.92,33.17)
\put(8.92,16.83){\line(0,1){16.34}}
%\end
%\emline(9.92,32.58)(9.92,17.33)
\put(9.92,32.58){\line(0,-1){15.25}}
%\end
%\emline(10.92,17.92)(10.92,32.08)
\put(10.92,17.92){\line(0,1){14.16}}
%\end
%\emline(11.92,31.50)(11.92,18.42)
\put(11.92,31.50){\line(0,-1){13.08}}
%\end
%\emline(12.92,19.00)(12.92,31.00)
\put(12.92,19.00){\line(0,1){12.00}}
%\end
%\emline(13.92,30.42)(13.92,19.58)
\put(13.92,30.42){\line(0,-1){10.84}}
%\end
%\emline(14.92,20.08)(14.92,29.92)
\put(14.92,20.08){\line(0,1){9.84}}
%\end
%\emline(15.92,29.33)(15.92,20.58)
\put(15.92,29.33){\line(0,-1){8.75}}
%\end
%\emline(16.92,21.17)(16.92,28.75)
\put(16.92,21.17){\line(0,1){7.58}}
%\end
%\emline(17.92,28.33)(17.92,21.67)
\put(17.92,28.33){\line(0,-1){6.66}}
%\end
%\emline(18.92,22.33)(18.92,27.67)
\put(18.92,22.33){\line(0,1){5.34}}
%\end
%\emline(19.92,27.17)(19.92,22.83)
\put(19.92,27.17){\line(0,-1){4.34}}
%\end
%\emline(20.92,23.42)(20.92,26.67)
\put(20.92,23.42){\line(0,1){3.25}}
%\end
%\emline(21.92,26.08)(21.92,23.92)
\put(21.92,26.08){\line(0,-1){2.16}}
%\end
%\emline(22.92,24.50)(22.92,25.50)
\put(22.92,24.50){\line(0,1){1.00}}
%\end
%\emline(1.00,37.25)(1.00,12.75)
\put(1.00,37.25){\line(0,-1){24.50}}
%\end
%\emline(2.00,13.17)(2.00,36.75)
\put(2.00,13.17){\line(0,1){23.58}}
%\end
%\emline(3.00,36.17)(3.00,13.75)
\put(3.00,36.17){\line(0,-1){22.42}}
%\end
%\emline(4.00,14.33)(4.00,35.67)
\put(4.00,14.33){\line(0,1){21.33}}
%\end
\end{picture}
\hfill
%TexCad Options
%\grade{\on}
%\emlines{\off}
%\beziermacro{\off}
%\reduce{\on}
%\snapping{\off}
%\quality{2.00}
%\graddiff{0.01}
%\snapasp{1}
%\zoom{4.00}
\unitlength 0.50ex \linethickness{0.4pt}
\begin{picture}(39.84,44.76)
\put(19.92,24.88){\oval(39.83,39.75)[]}
\put(20.00,25.00){\circle*{0.50}}
\put(25.50,28.00){\circle*{0.50}}
\put(25.50,22.00){\circle*{0.50}}
%\emline(20.00,25.00)(0.75,14.50)
\multiput(20.00,25.00)(-0.22,-0.12){88}{\line(-1,0){0.22}}
%\end
%\emline(20.00,25.00)(0.75,35.50)
\multiput(20.00,25.00)(-0.22,0.12){88}{\line(-1,0){0.22}}
%\end
\thicklines{}
%\emline(30.58,25.25)(1.00,41.25)
\multiput(30.58,25.25)(-0.22,0.12){134}{\line(-1,0){0.22}}
%\end
%\emline(30.58,24.75)(1.00,8.75)
\multiput(30.58,24.75)(-0.22,-0.12){134}{\line(-1,0){0.22}}
%\end
%\emline(31.25,25.42)(39.25,40.08)
\multiput(31.25,25.42)(0.12,0.22){67}{\line(0,1){0.22}}
%\end
%\emline(31.25,24.58)(39.25,9.92)
\multiput(31.25,24.58)(0.12,-0.22){67}{\line(0,-1){0.22}}
%\end
%\bezier{32}(29.00,26.08)(28.58,22.08)(32.17,22.75)
\qbezier(29.00,26.08)(28.58,22.08)(32.17,22.75)
%\multiput(29.00,26.08)(0.06,-1.02){2}{\line(0,-1){1.02}}
%\multiput(29.13,24.04)(0.11,-0.14){8}{\line(0,-1){0.14}}
%\multiput(30.04,22.90)(1.06,-0.08){2}{\line(1,0){1.06}}
%\end
%\bezier{48}(28.00,23.33)(27.00,29.50)(32.67,28.08)
\qbezier(28.00,23.33)(27.00,29.50)(32.67,28.08)
%\multiput(28.00,23.33)(-0.06,1.12){2}{\line(0,1){1.12}}
%\multiput(27.87,25.57)(0.11,0.40){4}{\line(0,1){0.40}}
%\multiput(28.32,27.15)(0.13,0.12){8}{\line(1,0){0.13}}
%\put(29.36,28.08){\line(1,0){3.31}}
%\end
\thinlines{} \linethickness{0.05pt}
%\emline(0.00,35.00)(1.58,35.00)
\put(0.00,35.00){\line(1,0){1.58}}
%\end
%\emline(0.00,33.00)(5.25,33.00)
\put(0.00,33.00){\line(1,0){5.25}}
%\end
%\emline(0.00,31.00)(9.00,31.00)
\put(0.00,31.00){\line(1,0){9.00}}
%\end
%\emline(0.00,29.00)(12.58,29.00)
\put(0.00,29.00){\line(1,0){12.58}}
%\end
%\emline(0.00,27.00)(16.25,27.00)
\put(0.00,27.00){\line(1,0){16.25}}
%\end
%\emline(0.00,25.00)(20.00,25.00)
\put(0.00,25.00){\line(1,0){20.00}}
%\end
%\emline(0.00,15.00)(1.58,15.00)
\put(0.00,15.00){\line(1,0){1.58}}
%\end
%\emline(0.00,17.00)(5.25,17.00)
\put(0.00,17.00){\line(1,0){5.25}}
%\end
%\emline(0.00,19.00)(9.00,19.00)
\put(0.00,19.00){\line(1,0){9.00}}
%\end
%\emline(0.00,21.00)(12.58,21.00)
\put(0.00,21.00){\line(1,0){12.58}}
%\end
%\emline(0.00,23.00)(16.25,23.00)
\put(0.00,23.00){\line(1,0){16.25}}
%\end
\put(6.00,25.00){\makebox(0,0)[cc]{$Y$}}
\put(20.00,3.00){\makebox(0,0)[ct]{\strut other $x$}}
\put(31.00,25.00){\circle*{1.00}}
\put(33.00,25.00){\makebox(0,0)[lc]{$D_x$}}
%\emline(1.00,8.75)(1.00,41.17)
\put(1.00,8.75){\line(0,1){32.42}}
%\end
%\emline(2.00,40.75)(2.00,9.33)
\put(2.00,40.75){\line(0,-1){31.42}}
%\end
%\emline(3.00,9.83)(3.00,40.17)
\put(3.00,9.83){\line(0,1){30.34}}
%\end
%\emline(4.00,39.58)(4.00,10.33)
\put(4.00,39.58){\line(0,-1){29.25}}
%\end
%\emline(5.00,10.92)(5.00,39.00)
\put(5.00,10.92){\line(0,1){28.08}}
%\end
%\emline(6.00,38.50)(6.00,11.42)
\put(6.00,38.50){\line(0,-1){27.08}}
%\end
%\emline(7.00,12.00)(7.00,38.00)
\put(7.00,12.00){\line(0,1){26.00}}
%\end
%\emline(8.00,37.42)(8.00,12.50)
\put(8.00,37.42){\line(0,-1){24.92}}
%\end
%\emline(9.00,13.08)(9.00,36.92)
\put(9.00,13.08){\line(0,1){23.84}}
%\end
%\emline(10.00,36.33)(10.00,13.67)
\put(10.00,36.33){\line(0,-1){22.66}}
%\end
%\emline(11.00,14.17)(11.00,35.83)
\put(11.00,14.17){\line(0,1){21.66}}
%\end
%\emline(12.00,35.25)(12.00,14.67)
\put(12.00,35.25){\line(0,-1){20.58}}
%\end
%\emline(13.00,15.25)(13.00,34.75)
\put(13.00,15.25){\line(0,1){19.50}}
%\end
%\emline(14.00,34.17)(14.00,15.75)
\put(14.00,34.17){\line(0,-1){18.42}}
%\end
%\emline(15.00,16.25)(15.00,33.67)
\put(15.00,16.25){\line(0,1){17.42}}
%\end
%\emline(16.00,33.08)(16.00,16.92)
\put(16.00,33.08){\line(0,-1){16.16}}
%\end
%\emline(17.00,17.42)(17.00,32.58)
\put(17.00,17.42){\line(0,1){15.16}}
%\end
%\emline(18.00,32.00)(18.00,17.92)
\put(18.00,32.00){\line(0,-1){14.08}}
%\end
%\emline(19.00,18.50)(19.00,31.50)
\put(19.00,18.50){\line(0,1){13.00}}
%\end
%\emline(20.00,31.00)(20.00,19.00)
\put(20.00,31.00){\line(0,-1){12.00}}
%\end
%\emline(21.00,19.58)(21.00,30.42)
\put(21.00,19.58){\line(0,1){10.84}}
%\end
%\emline(22.00,29.83)(22.00,20.08)
\put(22.00,29.83){\line(0,-1){9.75}}
%\end
%\emline(23.00,20.67)(23.00,29.33)
\put(23.00,20.67){\line(0,1){8.66}}
%\end
%\emline(24.00,28.75)(24.00,21.17)
\put(24.00,28.75){\line(0,-1){7.58}}
%\end
%\emline(25.00,21.75)(25.00,28.25)
\put(25.00,21.75){\line(0,1){6.50}}
%\end
%\emline(26.00,27.67)(26.00,22.33)
\put(26.00,27.67){\line(0,-1){5.34}}
%\end
%\emline(27.00,22.83)(27.00,27.17)
\put(27.00,22.83){\line(0,1){4.34}}
%\end
%\emline(28.00,26.58)(28.00,23.33)
\put(28.00,26.58){\line(0,-1){3.25}}
%\end
%\emline(29.00,23.92)(29.00,26.08)
\put(29.00,23.92){\line(0,1){2.16}}
%\end
%\emline(30.00,25.42)(30.00,24.50)
\put(30.00,25.42){\line(0,-1){0.92}}
%\end
\put(30.00,40.00){\makebox(0,0)[cc]{$\Cc^{*}_x$}}
\put(30.00,10.00){\makebox(0,0)[cc]{$\Cc_x$}}
\end{picture}
\bigskip

\raisebox{7ex}{$h=0$}\qquad \strut\hfill
%TexCad Options
%\grade{\on}
%\emlines{\off}
%\beziermacro{\off}
%\reduce{\on}
%\snapping{\off}
%\quality{2.00}
%\graddiff{0.01}
%\snapasp{1}
%\zoom{4.00}
\unitlength 0.50ex \linethickness{0.4pt}
\begin{picture}(39.84,39.76)
\put(19.92,19.88){\oval(39.83,39.75)[]}
\put(20.00,20.00){\circle*{0.50}}
\put(25.50,23.00){\circle*{0.50}}
\put(31.00,26.00){\circle*{0.50}}
\put(36.50,29.00){\circle*{0.50}}
\put(14.50,23.00){\circle*{0.50}} \put(9.00,26.00){\circle*{0.50}}
\put(3.50,29.00){\circle*{0.50}} \put(20.00,26.00){\circle*{0.50}}
\put(20.00,32.00){\circle*{0.50}}
\put(20.00,38.00){\circle*{0.50}}
\put(14.50,29.00){\circle*{0.50}} \put(9.00,32.00){\circle*{0.50}}
\put(14.50,35.00){\circle*{0.50}} \put(9.00,38.00){\circle*{0.50}}
\put(3.50,35.00){\circle*{0.50}} \put(25.50,29.00){\circle*{0.50}}
\put(31.00,32.00){\circle*{0.50}}
\put(25.50,35.00){\circle*{0.50}}
\put(31.00,38.00){\circle*{0.50}}
\put(36.50,35.00){\circle*{0.50}} \put(9.00,20.00){\circle*{0.50}}
\put(3.50,23.00){\circle*{0.50}} \put(31.00,20.00){\circle*{0.50}}
\put(36.50,23.00){\circle*{0.50}}
\put(25.50,17.00){\circle*{0.50}}
\put(31.00,14.00){\circle*{0.50}}
\put(36.50,11.00){\circle*{0.50}}
\put(14.50,17.00){\circle*{0.50}} \put(9.00,14.00){\circle*{0.50}}
\put(3.50,11.00){\circle*{0.50}} \put(20.00,14.00){\circle*{0.50}}
\put(20.00,8.00){\circle*{0.50}} \put(20.00,2.00){\circle*{0.50}}
\put(14.50,11.00){\circle*{0.50}} \put(9.00,8.00){\circle*{0.50}}
\put(14.50,5.00){\circle*{0.50}} \put(9.00,2.00){\circle*{0.50}}
\put(3.50,5.00){\circle*{0.50}} \put(25.50,11.00){\circle*{0.50}}
\put(31.00,8.00){\circle*{0.50}} \put(25.50,5.00){\circle*{0.50}}
\put(31.00,2.00){\circle*{0.50}} \put(36.50,5.00){\circle*{0.50}}
\put(3.50,17.00){\circle*{0.50}} \put(36.50,17.00){\circle*{0.50}}
%\vector(20.00,20.00)(25.50,23.00)
\put(25.50,23.00){\vector(2,1){0.2}}
\multiput(20.00,20.00)(0.21,0.12){26}{\line(1,0){0.21}}
%\end
%\vector(20.00,20.00)(25.50,17.00)
\put(25.50,17.00){\vector(2,-1){0.2}}
\multiput(20.00,20.00)(0.21,-0.12){26}{\line(1,0){0.21}}
%\end
\thicklines{}
%\emline(20.00,20.00)(0.75,9.50)
\multiput(20.00,20.00)(-0.22,-0.12){88}{\line(-1,0){0.22}}
%\end
%\emline(20.00,20.00)(0.75,30.50)
\multiput(20.00,20.00)(-0.22,0.12){88}{\line(-1,0){0.22}}
%\end
%\vector(20.00,20.00)(14.50,23.00)
\put(14.50,23.00){\vector(-2,1){0.2}}
\multiput(20.00,20.00)(-0.21,0.12){26}{\line(-1,0){0.21}}
%\end
%\vector(20.00,20.00)(14.50,17.00)
\put(14.50,17.00){\vector(-2,-1){0.2}}
\multiput(20.00,20.00)(-0.21,-0.12){26}{\line(-1,0){0.21}}
%\end
%\emline(20.00,20.00)(32.00,39.58)
\multiput(20.00,20.00)(0.12,0.20){100}{\line(0,1){0.20}}
%\end
%\emline(20.00,20.00)(32.00,0.42)
\multiput(20.00,20.00)(0.12,-0.20){100}{\line(0,-1){0.20}}
%\end
%\bezier{28}(18.00,21.08)(17.83,17.67)(21.33,17.83)
\qbezier(18.00,21.08)(17.83,17.67)(21.33,17.83)
%\multiput(18.00,21.08)(0.12,-0.66){3}{\line(0,-1){0.66}}
%\multiput(18.35,19.10)(0.27,-0.12){11}{\line(1,0){0.27}}
%\end
%\bezier{52}(17.00,18.42)(15.50,24.58)(22.00,23.33)
\qbezier(17.00,18.42)(15.50,24.58)(22.00,23.33)
%\multiput(17.00,18.42)(-0.09,0.70){3}{\line(0,1){0.70}}
%\multiput(16.72,20.51)(0.10,0.52){3}{\line(0,1){0.52}}
%\multiput(17.03,22.06)(0.11,0.13){8}{\line(0,1){0.13}}
%\multiput(17.93,23.06)(0.37,0.11){4}{\line(1,0){0.37}}
%\multiput(19.43,23.52)(1.29,-0.09){2}{\line(1,0){1.29}}
%\end
\thinlines{} \linethickness{0.2pt}
%\emline(0.00,30.00)(1.58,30.00)
\put(0.00,30.00){\line(1,0){1.58}}
%\end
%\emline(0.00,28.00)(5.25,28.00)
\put(0.00,28.00){\line(1,0){5.25}}
%\end
%\emline(0.00,26.00)(9.00,26.00)
\put(0.00,26.00){\line(1,0){9.00}}
%\end
%\emline(0.00,24.00)(12.58,24.00)
\put(0.00,24.00){\line(1,0){12.58}}
%\end
%\emline(0.00,22.00)(16.25,22.00)
\put(0.00,22.00){\line(1,0){16.25}}
%\end
%\emline(0.00,20.00)(20.00,20.00)
\put(0.00,20.00){\line(1,0){20.00}}
%\end
%\emline(0.00,10.00)(1.58,10.00)
\put(0.00,10.00){\line(1,0){1.58}}
%\end
%\emline(0.00,12.00)(5.25,12.00)
\put(0.00,12.00){\line(1,0){5.25}}
%\end
%\emline(0.00,14.00)(9.00,14.00)
\put(0.00,14.00){\line(1,0){9.00}}
%\end
%\emline(0.00,16.00)(12.58,16.00)
\put(0.00,16.00){\line(1,0){12.58}}
%\end
%\emline(0.00,18.00)(16.25,18.00)
\put(0.00,18.00){\line(1,0){16.25}}
%\end
\put(6.00,20.00){\makebox(0,0)[cc]{$Y$}}
\put(26.50,17.00){\makebox(0,0)[lc]{$\omega^{\vee}_1$}}
%\vector(20.00,20.00)(25.50,23.00)
\put(25.50,23.00){\vector(2,1){0.2}}
\multiput(20.00,20.00)(0.21,0.12){26}{\line(1,0){0.21}}
%\end
%\vector(20.00,20.00)(25.50,17.00)
\put(25.50,17.00){\vector(2,-1){0.2}}
\multiput(20.00,20.00)(0.21,-0.12){26}{\line(1,0){0.21}}
%\end
\put(15.00,16.00){\makebox(0,0)[lt]{$W$}}
%\emline(1.00,9.67)(1.00,30.33)
\put(1.00,9.67){\line(0,1){20.66}}
%\end
%\emline(2.00,29.75)(2.00,10.17)
\put(2.00,29.75){\line(0,-1){19.58}}
%\end
%\emline(3.00,10.67)(3.00,29.33)
\put(3.00,10.67){\line(0,1){18.66}}
%\end
%\emline(4.00,28.67)(4.00,11.33)
\put(4.00,28.67){\line(0,-1){17.34}}
%\end
%\emline(5.00,11.83)(5.00,28.17)
\put(5.00,11.83){\line(0,1){16.34}}
%\end
%\emline(6.00,27.58)(6.00,12.33)
\put(6.00,27.58){\line(0,-1){15.25}}
%\end
%\emline(7.00,12.92)(7.00,27.08)
\put(7.00,12.92){\line(0,1){14.16}}
%\end
%\emline(8.00,26.50)(8.00,13.42)
\put(8.00,26.50){\line(0,-1){13.08}}
%\end
%\emline(9.00,14.00)(9.00,26.00)
\put(9.00,14.00){\line(0,1){12.00}}
%\end
%\emline(10.00,25.42)(10.00,14.58)
\put(10.00,25.42){\line(0,-1){10.84}}
%\end
%\emline(11.00,15.08)(11.00,24.92)
\put(11.00,15.08){\line(0,1){9.84}}
%\end
%\emline(12.00,24.33)(12.00,15.58)
\put(12.00,24.33){\line(0,-1){8.75}}
%\end
%\emline(13.00,16.17)(13.00,23.75)
\put(13.00,16.17){\line(0,1){7.58}}
%\end
%\emline(14.00,23.33)(14.00,16.67)
\put(14.00,23.33){\line(0,-1){6.66}}
%\end
%\emline(15.00,17.33)(15.00,22.67)
\put(15.00,17.33){\line(0,1){5.34}}
%\end
%\emline(16.00,22.17)(16.00,17.83)
\put(16.00,22.17){\line(0,-1){4.34}}
%\end
%\emline(17.00,18.42)(17.00,21.67)
\put(17.00,18.42){\line(0,1){3.25}}
%\end
%\emline(18.00,21.08)(18.00,18.92)
\put(18.00,21.08){\line(0,-1){2.16}}
%\end
%\emline(19.00,19.50)(19.00,20.50)
\put(19.00,19.50){\line(0,1){1.00}}
%\end
\put(26.50,23.00){\makebox(0,0)[lc]{$\omega^{\vee}_2$}}
\put(15.00,24.00){\makebox(0,0)[lb]{$\widetilde{W}$}}
\put(4.00,10.00){\makebox(0,0)[lt]{$Y_i$}}
\put(4.00,30.00){\makebox(0,0)[lb]{$\widetilde{Y}_i$}}
\put(20.00,35.00){\makebox(0,0)[cc]{$\Cc^{*}$}}
\put(20.00,5.00){\makebox(0,0)[cc]{$\Cc$}}
\end{picture}
\hfill\strut
\bigskip
\end{center}
\end{figure}

The space of ordered triangles has three natural completions:
$(\PP^2)^3$, $({\PP^2}^{*})^3$, and
%*
\begin{xalignat*}{1}
X= \bigl\{(p_1,p_2,p_3,l_1,l_2,l_3)\bigm|p_j\in\PP^2,\
l_i\in{\PP^2}^{*},\ p_j\in l_i\text{ whenever }i\neq j\bigr\}.
\end{xalignat*}
%*
If $z_{kj}$ are the homogeneous coordinates of~$p_j$ in~$\PP^2$,
and $y_{ik}$ are the dual coordinates of~$l_i$ in~${\PP^2}^{*}$,
then $X$ is determined by 6~equations $(y\cdot z)_{ij}=0$, ($i\neq
j$) in $(\PP^2)^3\times({\PP^2}^{*})^3$. One verifies that the
Jacobian matrix is non-degenerate everywhere on $X\setminus Y$,
where $Y\subset X$ is given by the equations $p_1=p_2=p_3$,
$l_1=l_2=l_3$, and $\codim_XY=3$. Thence, by Serre's normality
criterion, $X$~is a normal complete intersection, smooth
outside~$Y$. It contains the following $G$-subvarieties of
degenerate triangles:
\begin{enumerate}
\item[$W_i$:] $p_j=p_k$ and $l_j=l_k$, $\{i,j,k\}=\{1,2,3\}$. (A
divisor.) \item[$\widetilde{W}$:] $p_1,p_2,p_3$ are collinear and
$l_1=l_2=l_3$. (The proper pullback of the divisor $\{\det z=0\}$
in~$(\PP^2)^3$.) \item[$W$:] $p_1=p_2=p_3$ and $l_1,l_2,l_3$ pass
through this point. (The proper pullback of the divisor $\{\det
y=0\}$ in~$({\PP^2}^{*})^3$.) \item[$\widetilde{Y}_i$:] $p_j=p_k$
and $l_1=l_2=l_3$ ($\codim=2$). \item[$Y_i$:] $l_j=l_k$ and
$p_1=p_2=p_3$ ($\codim=2$). \item[$Y$:] $p_1=p_2=p_3$ and
$l_1=l_2=l_3$ ($\codim=3$).
\end{enumerate}

Note that $\eta_i$ and $\widetilde{\eta}_i$ may be regarded as
certain homogeneous coordinates in the \ordinal{$i$} copy
of~$\PP^2$, resp.~${\PP^2}^{*}$, restricted to~$\HS$:
%*
\begin{equation*}
\left\{\begin{aligned}
\eta_i &=z_{3i}, \vphantom{\widetilde{\eta}_i} \\
\widetilde{\eta}_i &=
\begin{vmatrix}
z_{2j} & z_{2k} \\
z_{3j} & z_{3k}
\end{vmatrix},
\end{aligned}\right. \qquad
\text{or dually,} \qquad \left\{\begin{aligned} \eta_i &=
\begin{vmatrix}
y_{j1} & y_{j2} \\
y_{k1} & y_{k2}
\end{vmatrix}, \\
\widetilde{\eta}_i &=y_{i1},
\end{aligned}\right.
\end{equation*}
%*
and $\eta_{\infty}$ is a $3$-form in the matrix entries of $y$
or~$z$. Then
%*
\begin{align*}
\ef{\alpha_1}&=
\frac{\widetilde{\eta}_1(y)\widetilde{\eta}_2(y)\widetilde{\eta}_3(y)}
{\eta_{\infty}(y)}=
\frac{\widetilde{\eta}_1(z)\widetilde{\eta}_2(z)\widetilde{\eta}_3(z)}
{\eta_{\infty}(z)\det{z}} \\
\ef{\alpha_2}&=
\frac{\eta_1(y)\eta_2(y)\eta_3(y)}{\eta_{\infty}(y)\det{y}}=
\frac{\eta_1(z)\eta_2(z)\eta_3(z)}{\eta_{\infty}(z)}
\end{align*}
%*
One easily deduces that $\ef{\alpha_1},\eta_x$ ($\forall x\in C$)
are regular and do not vanish along~$W$, $\ef{\alpha_2},\eta_x$
along~$\widetilde{W}$, and $\ef{\alpha_1}$ (resp.~$\ef{\alpha_2}$)
has the \ordinal{1} order pole along~$\widetilde{W}$ (resp.~$W$).
Hence the $G$-valuations of $W,\widetilde{W},W_i$ are
$-\omega^{\vee}_2,-\omega^{\vee}_1,\eps_{x_i}$.

Since $X$ is complete and contains the minimal $G$-germ~$Y$ (the
closed orbit), we have $X=G\X$, where $\X$ is the minimal
$B$-chart of~$Y$ determined by the colored hypercone
$(\Cc_Y,\Dd^B_Y)$ of type~{\typeB} such that $\Cc_Y\supseteq\Vv$
and $\Vv_Y=\{-\omega^{\vee}_1,-\omega^{\vee}_2,
\eps_{x_1},\eps_{x_2},\eps_{x_3}\}$. It is easy to see from
Figure~\ref{hs(trk)} that there exists a unique such hypercone and
$\Dd^B_Y=\{D_x\mid x\neq x_1,x_2,x_3\}$. Its (hyper)faces
corresponding to various $G$-germs of~$X$ (including~$Y$) are
indicated in~Figure~\ref{hs(trk)} by the same letters.

A similar argument shows that $(\PP^2)^3$ is defined by the
colored hypercone $(\Cc,\{\widetilde{D}_i,D_x\mid i=1,2,3,\ x\neq
x_1,x_2,x_3\})$ and $({\PP^2}^{*})^3$ by
$(\Cc^{*},\{\widetilde{D}_i,D_x\mid i=1,2,3,\ x\neq
x_1,x_2,x_3\})$.

The space $\SL_3(\kk)/N(T)$ of unordered triangles and its
completion is studied in~\cite[\S9]{C-div}. The resolution of
singularities of $X$ was studied already by Schubert with
applications to enumerative geometry, see \cite{triang},
\cite[\S9]{C-div}, and~\ref{intersect}.
\end{example}

We say that a $G$-model $X$ \emph{is of type~{\typeA}} if it
contains no $G$-germs of type~{\typeB}, i.e., any $G$-orbit in $X$
is contained in finitely many $B$-stable divisors. For any~$X$,
there is a canonical proper birational morphism
$\nu:\widehat{X}\to X$ such that $\widehat{X}$ is of type~{\typeA}
and $\nu$ is isomorphic in codimension~$1$. (Just subdivide each
hypercone $\Cc$ from the hyperfan of $X$ by $\Kk=\Cc\cap\ES$.)

In characteristic zero, singularities of $G$-models of
type~{\typeA} are good.
\begin{theorem}
If $X$ is of type~{\typeA}, then all $G$-subvarieties $Y\subseteq
X$ are normal and have rational singularities.
\end{theorem}
\begin{proof}
By the local structure theorem, we may assume that $X$ is affine
and of type~{\typeA}, i.e., $\kk[X]^B\neq k$. Passing to the
categorical quotient by~$U$, we may assume that $G=B=T$. In the
notation of Example~\ref{c_T=1}, we may replace $X$ by $X_i$ and
assume that $X$ is an affine toric $(T\times\kk^{\times})$-variety
such that $X\by T\iso\AAA^1$ ($\kk^{\times}$-equivariantly). Then
each $T$-stable closed subvariety of~$X$ is either
$(T\times\kk^{\times})$-stable or lying in the fiber of the
quotient map $X\to\kk$ over a nonzero point, which is a toric
$T$-variety. Thus the question is reduced to the case of toric
varieties.
\end{proof}

\section{Divisors}
\label{div}

In the study of divisors on $G$-models, we may restrict our
attention to $B$-stable ones, by the following result.
\begin{proposition}
Let a connected solvable algebraic group $B$ act on a variety~$X$.
Then any Weil divisor on $X$ is rationally equivalent to a
$B$-stable one.
\end{proposition}
\begin{proof}
Replacing $X$ by $X^{\reg}$, we may assume that $X$ is smooth.
Replacing $X$ by a $B$-stable open subset, we may assume that $X$
is quasiprojective. Then any Weil divisor $\delta$ on $X$ is
Cartier. Furthermore, $\delta$~is the difference of two globally
generated divisors. Therefore we may assume that $\delta$ is
globally generated. The line bundle $\Lin{\delta}$ is
$B$-linearizable by Theorem~\ref{G-lin}, and the $B$-module
$\Ho^0(X,\Lin{\delta})$ contains a nonzero
$B$-eigensection~$\sigma$. The divisor $\divr\sigma$ is $B$-stable
and equivalent to~$\delta$.
\end{proof}
\begin{remark}
The proposition is true for any algebraic cycle, see
Theorem~\ref{B-cycle}.
\end{remark}

Our first aim is to describe Cartier divisors.

For any Cartier divisor $\delta$ on~$X$, we shall always equip the
respective line bundle with a $G$-linearization (see
Appendix~\ref{rat.mod&lin}).

\begin{lemma}[{\cite[2.2]{inv.mot}}]
\label{Cartier} Any prime divisor $D\subset X$ that does not
contain a $G$-orbit of $X$ is globally generated Cartier.
\end{lemma}
\begin{proof}
Let $\iota:X^{\reg}\embeds X$ be the inclusion of the subset of
smooth points. Then $D\cap X^{\reg}$ is Cartier on $X^{\reg}$, and
$D\cap X^{\reg}= \divr\eta$ for some
$\eta\in\Ho^0(X^{\reg},\Lin{D\cap X^{\reg}})$. As $X$ is normal,
$\Ll=\iota_{*}\Lin{D\cap X^{\reg}}$ is a trivial line bundle on
$X\setminus D$. As $G$ acts on~$\Ll$, the set of points where
$\Ll$ is not invertible is $G$-stable and contained in~$D$, hence
empty. Therefore $\Ll$ is a line bundle on~$X$.

If we regard $\eta$ as an element of $\Ho^0(X,\Ll)$, then
$D=\divr\eta$, because the equality holds on~$X^{\reg}$ and
$\codim_X(X\setminus X^{\reg})>1$. Furthermore, $\Ll$~is generated
by $g\eta$, $g\in G$, because the set of their common zeroes is
$\bigcap_{g\in G}gD=\emptyset$. The assertion follows.
\end{proof}

The following criterion says that a $B$-stable divisor is Cartier
iff it is determined by a local equation in a neighborhood of a
\emph{general} point of each $G$-subvariety.
\begin{theorem}\label{loc.prin}
Suppose $\delta$ is a $B$-stable divisor on~$X$. Then $\delta$ is
Cartier iff for any $G$-subvariety $Y\subseteq X$ there exists
$f_Y\in K^{(B)}$ such that each prime divisor $D\supseteq Y$
occurs in~$\delta$ with multiplicity~$v_D(f_Y)$.
\end{theorem}
\begin{proof}
Suppose $\delta$ is locally principal in general points of
$G$-subvarieties. Take any $G$-orbit $Y\subseteq X$ and a
$B$-chart $\X\subseteq X$ intersecting~$Y$. Replacing $\delta$ by
$\delta-\divr f_Y$, we may assume that no component of~$\delta$
contains~$Y$. Let $D_1,\dots,D_n$ be the components of~$\delta$
intersecting~$\X$ and $w_i=v_{D_i}$, $i=1,\dots,n$. We have either
$w_i\notin\Vv_Y$, or $D_i\notin\Dd^B_Y$ (depending on whether
$D_i$ is $G$-stable or not). By \reftag{V$'$} or~\reftag{D$'$},
$\exists f_i\in\kk[\X]^{(B)}:\ f_i|_Y\neq0,\ w_i(f_i)>0$. Now we
may replace $\X$ by its localization at $f_1\dots f_n$ and assume
that $\X$ intersects no component of~$\delta$. By
Lemma~\ref{Cartier}, $\delta$~is Cartier on~$G\X$, whence on~$X$.

Now suppose $\delta$ is Cartier, and $Y\subseteq X$ is a
$G$-subvariety. By Sumihiro's Theorem~\ref{Sumihiro}, there is an
open $G$-stable quasiprojective subvariety $X_0\subseteq X$
intersecting~$Y$. The restriction of $\delta$ to $X_0$ may be
represented as the difference of two globally generated divisors,
hence we may replace $X$ by $X_0$ and assume that $\delta$ is
globally generated.

It follows that the annihilator of $Y$ in $\Ho^0(X,\Lin{\delta})$
is a proper $G$-sub\-mod\-ule, whence there is a section
$\sigma\in\Ho^0(X,\Lin{\delta})^{(B)}$ such that $\sigma|_Y\neq0$.
Therefore $\delta$ is principal on the $B$-stable open subset
$X_{\sigma}$ intersecting~$Y$, and we may take for $f_Y$ the
equation of $\delta$ on $X_{\sigma}$.
\end{proof}
By Theorem~\ref{loc.prin}, a Cartier divisor on~$X$ is determined
by the following data:
\begin{enumerate}
\item\label{f_Y} a collection of rational $B$-eigenfunctions $f_Y$
given for each $G$-germ $Y\in{}_GX$ and such that
$w(f_{Y_1})=w(f_{Y_2})$, $v_D(f_{Y_1})=v_D(f_{Y_2})$, $\forall
w\in\Vv_{Y_1}\cap\Vv_{Y_2}$, $D\in\Dd^B_{Y_1}\cap\Dd^B_{Y_2}$.
\item\label{m_D} a collection of integers $m_D$,
$D\in\Dd^B\setminus\bigcup_{Y\subseteq X}\Dd^B_Y$, only finitely
many of them being nonzero ($m_D$ is the multiplicity of~$D$ in
the divisor).
\end{enumerate}
\begin{remark}
It suffices to specify the local equations $f_Y$ only for closed
$G$-orbits $Y\subseteq X$: if a $G$-subvariety $Y\subseteq X$
contains a closed orbit $Y_0$, then we may put $f_Y=f_{Y_0}$.
\end{remark}

When a Cartier divisor is replaced by a rationally equivalent one,
the local equations $f_Y$ are replaced by $f_Yf$ for some $f\in
K^{(B)}$, and $m_D$ are replaced by $m_D+v_D(f)$.

In the case of complexity $\leq1$, the data~\ref{f_Y}--\ref{m_D}
are retranslated to the language of convex geometry.

Consider first the spherical case. Each $f_Y$ defines a function
$\psi_Y$ on the cone~$\Cc_Y$, which is the restriction of a linear
functional $\lambda_Y\in\RG$. We may assume that $f_{Y_1}=f_{Y_2}$
if $\Cc_{Y_1}$ is a face of $\Cc_{Y_2}$, whence
$\psi_{Y_1}=\psi_{Y_2}|_{\Cc_{Y_1}}$. In particular,
$\psi_Y$~paste together in a piecewise linear function on
$\bigcup_{Y\subseteq X}(\Cc_Y\cap\Vv)$. A collection
$\psi=(\psi_Y)$ of functions $\psi_Y$ on $\Cc_Y$ with the above
properties is called an \emph{integral piecewise linear function}
on the colored fan $\Ff$ of $X$.

Note that generally $\psi$ is \emph{not} a well-defined function
on $\bigcup_{\Cc\in\Ff}\Cc$ as the following example shows.
\begin{example}[\cite{SL(3)/SL(2)}]
\label{PL-fun} Let $G=\SL_3(\kk)$, $H=\SL_2(\kk)=$ the common
stabilizer of $e_1\in\kk^3$, $x_1\in(\kk^3)^{*}$. Then $\HS=G/H$
is defined in $\kk^3\oplus(\kk^3)^{*}$ by an equation $\langle
v,v^{*}\rangle=1$. The $B$-divisors $D_1,D_2\in\Dd^B$ are defined
by the restrictions $\eta_1,\eta_2$ of linear $B$-eigenfunctions
on $\kk^3\oplus(\kk^3)^{*}$ of $B$-weights $\omega_1,\omega_2$.
Here $\eta_1,\eta_2\in K^{(B)}$ and
$\RG=\langle\omega_1,\omega_2\rangle$, whence
$\res(D_i)=\alpha^{\vee}_i$.

A one-dimensional torus acting on the summands of
$\kk^3\oplus(\kk^3)^{*}$ by the weights $\pm1$ commutes with $G$
and preserves $\HS$. The respective grading of $\kk(\HS)$
determines two $G$-valuations $v_{\pm}(f)=\pm\deg f_{\mp}$, where
$f_{\pm}$ is the highest/lowest degree term of $f\in\kk[\HS]$.
Since $\deg\eta_i=(-1)^{i-1}$, we have
$v_{\pm}=\pm\alpha^{\vee}_1\mp\alpha^{\vee}_2$. It follows easily
from Corollary~\ref{sph.qaff} that the colored space looks like in
Figure~\ref{SL(3)/SL(2)}.
\begin{figure}[h!]
\caption{}\label{SL(3)/SL(2)}
\begin{center}
%TexCad Options
%\grade{\on}
%\emlines{\off}
%\beziermacro{\off}
%\reduce{\on}
%\snapping{\off}
%\quality{2.00}
%\graddiff{0.01}
%\snapasp{1}
%\zoom{4.00}
\unitlength 0.50ex \linethickness{0.4pt}
\begin{picture}(32.00,40.00)
\thicklines{}
%\emline(20.00,20.00)(32.00,39.58)
\multiput(20.00,20.00)(0.12,0.20){100}{\line(0,1){0.20}}
%\end
%\emline(20.00,20.00)(32.00,0.42)
\multiput(20.00,20.00)(0.12,-0.20){100}{\line(0,-1){0.20}}
%\end
\put(26.08,30.00){\circle*{1.00}}
\put(28.00,30.00){\makebox(0,0)[lc]{$D_2$}}
%\bezier{112}(17.92,28.00)(30.75,27.00)(24.25,13.00)
\multiput(17.92,28.00)(0.71,-0.09){3}{\line(1,0){0.71}}
\multiput(20.06,27.72)(0.37,-0.10){5}{\line(1,0){0.37}}
\multiput(21.89,27.23)(0.25,-0.12){6}{\line(1,0){0.25}}
\multiput(23.41,26.53)(0.15,-0.11){8}{\line(1,0){0.15}}
\multiput(24.62,25.63)(0.11,-0.14){8}{\line(0,-1){0.14}}
\multiput(25.52,24.52)(0.12,-0.26){5}{\line(0,-1){0.26}}
\multiput(26.12,23.20)(0.10,-0.51){3}{\line(0,-1){0.51}}
\put(26.41,21.67){\line(0,-1){1.73}}
\multiput(26.39,19.94)(-0.11,-0.65){3}{\line(0,-1){0.65}}
\multiput(26.06,18.00)(-0.11,-0.36){6}{\line(0,-1){0.36}}
\multiput(25.42,15.85)(-0.12,-0.29){10}{\line(0,-1){0.29}}
%\end
\put(26.08,10.00){\circle*{1.00}}
\put(28.00,10.00){\makebox(0,0)[lc]{$D_1$}}
%\bezier{76}(18.17,13.00)(26.58,14.83)(23.08,25.00)
\multiput(18.17,13.00)(0.33,0.10){6}{\line(1,0){0.33}}
\multiput(20.18,13.63)(0.20,0.11){8}{\line(1,0){0.20}}
\multiput(21.77,14.54)(0.12,0.12){10}{\line(0,1){0.12}}
\multiput(22.95,15.74)(0.11,0.21){7}{\line(0,1){0.21}}
\multiput(23.72,17.24)(0.12,0.59){3}{\line(0,1){0.59}}
\put(24.08,19.02){\line(0,1){2.07}}
\multiput(24.03,21.09)(-0.12,0.49){8}{\line(0,1){0.49}}
%\end
%\bezier{68}(18.67,15.00)(12.00,20.00)(18.67,25.00)
\multiput(18.67,15.00)(-0.13,0.11){13}{\line(-1,0){0.13}}
\multiput(17.00,16.47)(-0.11,0.15){10}{\line(0,1){0.15}}
\multiput(15.90,17.94)(-0.10,0.29){5}{\line(0,1){0.29}}
\put(15.38,19.41){\line(0,1){1.47}}
\multiput(15.44,20.88)(0.11,0.25){6}{\line(0,1){0.25}}
\multiput(16.07,22.35)(0.12,0.12){22}{\line(0,1){0.12}}
%\end
%\emline(20.00,20.00)(14.92,39.75)
\multiput(20.00,20.00)(-0.12,0.46){43}{\line(0,1){0.46}}
%\end
%\emline(20.00,20.00)(14.92,0.25)
\multiput(20.00,20.00)(-0.12,-0.46){43}{\line(0,-1){0.46}}
%\end
\thinlines{}
%\emline(20.00,0.00)(20.00,40.00)
\put(20.00,0.00){\line(0,1){40.00}}
%\end
\linethickness{0.05pt}
%\emline(20.00,20.00)(0.00,20.00)
\put(20.00,20.00){\line(-1,0){20.00}}
%\end
%\emline(20.00,22.00)(0.00,22.00)
\put(20.00,22.00){\line(-1,0){20.00}}
%\end
%\emline(20.00,24.00)(0.00,24.00)
\put(20.00,24.00){\line(-1,0){20.00}}
%\end
%\emline(20.00,26.00)(0.00,26.00)
\put(20.00,26.00){\line(-1,0){20.00}}
%\end
%\emline(20.00,28.00)(0.00,28.00)
\put(20.00,28.00){\line(-1,0){20.00}}
%\end
%\emline(20.00,30.00)(0.00,30.00)
\put(20.00,30.00){\line(-1,0){20.00}}
%\end
%\emline(20.00,32.00)(0.00,32.00)
\put(20.00,32.00){\line(-1,0){20.00}}
%\end
%\emline(20.00,34.00)(0.00,34.00)
\put(20.00,34.00){\line(-1,0){20.00}}
%\end
%\emline(20.00,36.00)(0.00,36.00)
\put(20.00,36.00){\line(-1,0){20.00}}
%\end
%\emline(20.00,38.00)(0.00,38.00)
\put(20.00,38.00){\line(-1,0){20.00}}
%\end
%\emline(20.00,40.00)(0.00,40.00)
\put(20.00,40.00){\line(-1,0){20.00}}
%\end
%\emline(20.00,18.00)(0.00,18.00)
\put(20.00,18.00){\line(-1,0){20.00}}
%\end
%\emline(20.00,16.00)(0.00,16.00)
\put(20.00,16.00){\line(-1,0){20.00}}
%\end
%\emline(20.00,14.00)(0.00,14.00)
\put(20.00,14.00){\line(-1,0){20.00}}
%\end
%\emline(20.00,12.00)(0.00,12.00)
\put(20.00,12.00){\line(-1,0){20.00}}
%\end
%\emline(20.00,10.00)(0.00,10.00)
\put(20.00,10.00){\line(-1,0){20.00}}
%\end
%\emline(20.00,8.00)(0.00,8.00)
\put(20.00,8.00){\line(-1,0){20.00}}
%\end
%\emline(20.00,6.00)(0.00,6.00)
\put(20.00,6.00){\line(-1,0){20.00}}
%\end
%\emline(20.00,4.00)(0.00,4.00)
\put(20.00,4.00){\line(-1,0){20.00}}
%\end
%\emline(20.00,2.00)(0.00,2.00)
\put(20.00,2.00){\line(-1,0){20.00}}
%\end
%\emline(20.00,0.00)(0.00,0.00)
\put(20.00,0.00){\line(-1,0){20.00}}
%\end
\put(14.00,20.00){\makebox(0,0)[rc]{$\Cc_0$}}
\put(24.00,32.00){\makebox(0,0)[cb]{$\Cc_1$}}
\put(24.00,8.00){\makebox(0,0)[ct]{$\Cc_2$}}
\put(5.00,10.00){\makebox(0,0)[cc]{$\Vv$}}
\end{picture}
\end{center}
\end{figure}

Take a fan $\Ff$ determined by the cones $\Cc_0,\Cc_1,\Cc_2$ on
the figure. Since $\Cc_1\cap\Cc_2$ is a solid cone, a piecewise
linear function on $\Ff$ defines a function on
$\Cc_0\cup\Cc_1\cup\Cc_2=\ES$ iff it is linear.
\end{example}

Let $\PL(\Ff)$ be the group of all integral piecewise linear
functions on~$\Ff$, and $\LO(\Ff)$ be its subgroup of linear
functions $\psi=(\lambda|_{\Cc_Y})$, $\lambda\in\RG$.

The above discussion yields the following exact sequences:
%*
\begin{gather}
0\longrightarrow \ZZ\left(\Dd^B\setminus\bigcup_{Y\subseteq
X}\Dd^B_Y\right)
\longrightarrow\CaDiv(X)^B\longrightarrow\PL(\Ff)\longrightarrow0
\label{CaDiv}\\
\RG\cap\Ff^{\ann}\longrightarrow\PrDiv(X)^B\longrightarrow\LO(\Ff)
\longrightarrow0 \label{PrDiv}
\end{gather}
%*
where $\CaDiv(\cdot)$ and $\PrDiv(\cdot)$ denote the groups of
Cartier and principal divisors, respectively, and $\Ff^{\ann}$ is
the annihilator of the union of all cones in~$\Ff$.
\begin{theorem}[{\cite[3.1]{Pic(sph)}}] There is an exact
sequence
%*
\begin{equation*}
\RG\cap\Ff^{\ann}\longrightarrow
\ZZ\left(\Dd^B\setminus\bigcup_{Y\subseteq X}\Dd^B_Y\right)
\longrightarrow\Pic X\longrightarrow\PL(\Ff)/\LO(\Ff)
\longrightarrow0
\end{equation*}
%*
If $X$ contains a complete $G$-orbit, then $\Pic X$ is free
Abelian of finite rank.
\end{theorem}
\begin{proof}
The exact sequence is a consequence of
\eqref{CaDiv}--\eqref{PrDiv}. If $Y\subseteq X$ is a complete
$G$-orbit, then $\Ff^{\ann}\subseteq\Cc_Y^{\ann}=0$ by
Propositions \ref{cr(sph.orb)} and~\ref{r=0}. Then it is easy to
see that $\PL(\Ff)/\LO(\Ff)$ has finite rank and no torsion,
whence the second assertion.
\end{proof}

A spherical $G$-variety $X$ having only one closed orbit
$Y\subseteq X$ is called \emph{simple}. Its fan consists of all
supported colored faces of $(\Cc_Y,\Dd^B_Y)$.
\begin{corollary}
If $X$ is simple with the closed orbit $Y\subseteq X$, then there
is an exact sequence
%*
\begin{equation*}
\RG\cap\Cc_Y^{\ann}\longrightarrow
\ZZ(\Dd^B\setminus\Dd^B_Y)\longrightarrow\Pic X\longrightarrow0
\end{equation*}
%*
\end{corollary}
\begin{corollary}\label{Pic(s.c)}
If $X$ is simple and the closed orbit $Y\subseteq X$ is complete,
then $\Pic X=\ZZ(\Dd^B\setminus\Dd^B_Y)$ is free Abelian.
\end{corollary}
\begin{example}
If $X=G/P$ is a generalized flag variety, then $\Pic X$ is freely
generated by the Schubert divisors
$D_{\alpha}=\overline{B[w_Gr_{\alpha}]}$, $\alpha\in\Pi\setminus
I$, where $I\subseteq\Pi$ is the set of simple roots defining the
parabolic $P\supseteq B$.
\end{example}
\begin{example}[{\cite{Pic(sph)}}]
\label{conics} Let $G=\SL_3(\kk)$,
$H=N_G(\SO_3(\kk))=\SO_3(\kk)\times\ZZz_3$, where
$\ZZz_3=Z(\SL_3(\kk))$. Then $\HS=G/H$ is the space of conics
in~$\PP^2$. The coisotropy representation is the natural
representation of $\SO_3(\kk)$ in traceless symmetric matrices,
whence $H_{*}=\left\{\left(
\begin{smallmatrix}
\pm1 &   0  &   0 \\
  0  & \pm1 &   0 \\
  0  &   0  & \pm1
\end{smallmatrix}\right)\right\}$
is the Klein $4$-group, and
$\RG(\HS)=\langle2\alpha_1,2\alpha_2\rangle$, where $\alpha_i$ are
simple roots of $\SL_3(\kk)$. By~\eqref{2c+r(G/H)}, $\HS$~is
spherical.

We may consider $\HS$ as the projectivization of the open subset
of non-degenerate quadratic forms in $\Sym^2(\kk^3)^{*}$. The two
$(B\times H)$-eigenfunctions $\eta_1(q)=q_{11}$, $\eta_2(q)=
\begin{vmatrix}
q_{11} & q_{12} \\
q_{21} & q_{22}
\end{vmatrix}$ ($q\in\Sym^2(\kk^3)^{*}$) of biweights
$(2\omega_i,2i\eps)$, where $\omega_i$ are the fundamental weights
of~$G$ and $\eps$ is the weight of $\ZZz_3$ in~$\kk^3$, define the
two $B$-divisors $D_1,D_2$. They impose the conditions that a
conic passes through the $B$-fixed point, resp.\ is tangent to the
$B$-fixed line. Since $\ef{2\alpha_1}=\eta_1^2/\eta_2,\
\ef{2\alpha_2}=\eta_2^2/\eta_1\in K^{(B)}$, and their weights
$2\alpha_1,2\alpha_2$ generate~$\RG$, there are no other
$B$-divisors. (Indeed, if $\eta\in\kk[G]^{(B\times
H)}_{(\lambda,\chi)}$, then either $\chi=0$, $\lambda\in\RG$, or
$\chi=2i\eps$, $\lambda-2\omega_i\in\RG$, hence $\eta$ is
proportional to the product of $\eta_1,\eta_2$ and their
inverses.) Furthermore,
$\RG^{*}=\langle\omega^{\vee}_1/2,\omega^{\vee}_2/2\rangle$ and
$\res(D_i)=\alpha^{\vee}_i/2$.

The complement to $\HS$ in $(\PP^5)^{*}=\PP(\Sym^2(\kk^3)^{*})$ is
a $G$-stable prime divisor $\{\det q=0\}$, and the respective
$G$-valuation is $-\omega^{\vee}_2/2\in\ES$, because in
homogeneous coordinates $\ef{2\alpha_1}(q)=\eta_1(q)^2/\eta_2(q)$,
$\ef{2\alpha_2}(q)=\eta_2(q)^2/\eta_1(q)\det q$. The unique closed
orbit $Y=\{\rk q=1\}$ has the colored data
$\Vv_Y=\{-\omega^{\vee}_2/2\}$, $\Dd^B_Y=\{D_2\}$.

Similarly, we embed $\HS$ in $\PP^5=\PP(\Sym^2\kk^3)$ by the map
$q\to q^{\vee}$ (=the adjoint matrix of~$q$) sending a conic to
the dual one. Here the unique closed orbit $Y^{\vee}=\{\rk
q^{\vee}=1\}$ has the colored data
$\Vv_{Y^{\vee}}=\{-\omega^{\vee}_1/2\}$,
$\Dd^B_{Y^{\vee}}=\{D_1\}$.

Since $\PP^5,(\PP^5)^{*}$ are complete, the cones $\Cc_Y$ and
$\Cc_{Y^{\vee}}$ contain~$\Vv$, whence $\Vv$ is generated by
$-\omega^{\vee}_1/2,-\omega^{\vee}_2/2$. The colored equipment of
$\HS$ is represented at Figure~\ref{SL(3)/SO(3)}.
\begin{figure}[h!]
\caption{}\label{SL(3)/SO(3)}
\begin{center}
%TexCad Options
%\grade{\on}
%\emlines{\off}
%\beziermacro{\off}
%\reduce{\on}
%\snapping{\off}
%\quality{2.00}
%\graddiff{0.01}
%\snapasp{1}
%\zoom{4.00}
\unitlength 0.50ex \linethickness{0.4pt}
\begin{picture}(32.00,39.58)
%\vector(20.00,20.00)(30.00,25.50)
\put(30.00,25.50){\vector(2,1){0.2}}
\multiput(20.00,20.00)(0.22,0.12){46}{\line(1,0){0.22}}
%\end
%\vector(20.00,20.00)(30.00,14.50)
\put(30.00,14.50){\vector(2,-1){0.2}}
\multiput(20.00,20.00)(0.22,-0.12){46}{\line(1,0){0.22}}
%\end
\thicklines{}
%\vector(20.00,20.00)(10.00,25.50)
\put(10.00,25.50){\vector(-2,1){0.2}}
\multiput(20.00,20.00)(-0.22,0.12){46}{\line(-1,0){0.22}}
%\end
%\vector(20.00,20.00)(10.00,14.50)
\put(10.00,14.50){\vector(-2,-1){0.2}}
\multiput(20.00,20.00)(-0.22,-0.12){46}{\line(-1,0){0.22}}
%\end
%\emline(20.00,20.00)(0.75,9.50)
\multiput(20.00,20.00)(-0.22,-0.12){88}{\line(-1,0){0.22}}
%\end
%\emline(20.00,20.00)(0.75,30.50)
\multiput(20.00,20.00)(-0.22,0.12){88}{\line(-1,0){0.22}}
%\end
%\emline(20.00,20.00)(32.00,39.58)
\multiput(20.00,20.00)(0.12,0.20){100}{\line(0,1){0.20}}
%\end
%\emline(20.00,20.00)(32.00,0.42)
\multiput(20.00,20.00)(0.12,-0.20){100}{\line(0,-1){0.20}}
%\end
\put(6.00,20.00){\makebox(0,0)[cc]{$X$}}
\put(16.00,31.00){\makebox(0,0)[cb]{$(\PP^5)^{*}$}}
\put(16.00,9.00){\makebox(0,0)[ct]{$\PP^5$}}
%\bezier{116}(13.00,16.17)(10.92,30.92)(24.25,27.00)
\qbezier(13.00,16.17)(10.92,30.92)(24.25,27.00)
%\multiput(13.00,16.17)(-0.08,0.80){3}{\line(0,1){0.80}}
%\put(12.76,18.57){\line(0,1){2.13}}
%\multiput(12.74,20.70)(0.11,0.92){2}{\line(0,1){0.92}}
%\multiput(12.95,22.55)(0.11,0.39){4}{\line(0,1){0.39}}
%\multiput(13.40,24.12)(0.11,0.22){6}{\line(0,1){0.22}}
%\multiput(14.07,25.41)(0.11,0.13){8}{\line(0,1){0.13}}
%\multiput(14.97,26.43)(0.16,0.11){7}{\line(1,0){0.16}}
%\multiput(16.10,27.17)(0.34,0.12){4}{\line(1,0){0.34}}
%\multiput(17.46,27.63)(0.79,0.09){2}{\line(1,0){0.79}}
%\put(19.05,27.82){\line(1,0){1.82}}
%\multiput(20.87,27.73)(0.48,-0.10){7}{\line(1,0){0.48}}
%\end
%\bezier{84}(15.00,22.67)(13.17,12.33)(23.00,15.00)
\qbezier(15.00,22.67)(13.17,12.33)(23.00,15.00)
%\multiput(15.00,22.67)(-0.09,-0.76){3}{\line(0,-1){0.76}}
%\put(14.73,20.39){\line(0,-1){1.91}}
%\multiput(14.79,18.48)(0.10,-0.38){4}{\line(0,-1){0.38}}
%\multiput(15.18,16.94)(0.10,-0.17){7}{\line(0,-1){0.17}}
%\multiput(15.90,15.77)(0.15,-0.11){7}{\line(1,0){0.15}}
%\multiput(16.95,14.97)(0.35,-0.11){4}{\line(1,0){0.35}}
%\put(18.33,14.54){\line(1,0){1.71}}
%\multiput(20.05,14.47)(0.59,0.11){5}{\line(1,0){0.59}}
%\end
\thinlines{} \linethickness{0.05pt}
%\emline(0.00,30.00)(1.58,30.00)
\put(0.00,30.00){\line(1,0){1.58}}
%\end
%\emline(0.00,28.00)(5.25,28.00)
\put(0.00,28.00){\line(1,0){5.25}}
%\end
%\emline(0.00,26.00)(9.00,26.00)
\put(0.00,26.00){\line(1,0){9.00}}
%\end
%\emline(0.00,24.00)(12.58,24.00)
\put(0.00,24.00){\line(1,0){12.58}}
%\end
%\emline(0.00,22.00)(16.25,22.00)
\put(0.00,22.00){\line(1,0){16.25}}
%\end
%\emline(0.00,20.00)(20.00,20.00)
\put(0.00,20.00){\line(1,0){20.00}}
%\end
%\emline(0.00,10.00)(1.58,10.00)
\put(0.00,10.00){\line(1,0){1.58}}
%\end
%\emline(0.00,12.00)(5.25,12.00)
\put(0.00,12.00){\line(1,0){5.25}}
%\end
%\emline(0.00,14.00)(9.00,14.00)
\put(0.00,14.00){\line(1,0){9.00}}
%\end
%\emline(0.00,16.00)(12.58,16.00)
\put(0.00,16.00){\line(1,0){12.58}}
%\end
%\emline(0.00,18.00)(16.25,18.00)
\put(0.00,18.00){\line(1,0){16.25}}
%\end
%\emline(1.00,9.67)(1.00,30.33)
\put(1.00,9.67){\line(0,1){20.66}}
%\end
%\emline(2.00,29.75)(2.00,10.17)
\put(2.00,29.75){\line(0,-1){19.58}}
%\end
%\emline(3.00,10.67)(3.00,29.33)
\put(3.00,10.67){\line(0,1){18.66}}
%\end
%\emline(4.00,28.67)(4.00,11.33)
\put(4.00,28.67){\line(0,-1){17.34}}
%\end
%\emline(5.00,11.83)(5.00,28.17)
\put(5.00,11.83){\line(0,1){16.34}}
%\end
%\emline(6.00,27.58)(6.00,12.33)
\put(6.00,27.58){\line(0,-1){15.25}}
%\end
%\emline(7.00,12.92)(7.00,27.08)
\put(7.00,12.92){\line(0,1){14.16}}
%\end
%\emline(8.00,26.50)(8.00,13.42)
\put(8.00,26.50){\line(0,-1){13.08}}
%\end
%\emline(9.00,14.00)(9.00,26.00)
\put(9.00,14.00){\line(0,1){12.00}}
%\end
%\emline(10.00,25.42)(10.00,14.58)
\put(10.00,25.42){\line(0,-1){10.84}}
%\end
%\emline(11.00,15.08)(11.00,24.92)
\put(11.00,15.08){\line(0,1){9.84}}
%\end
%\emline(12.00,24.33)(12.00,15.58)
\put(12.00,24.33){\line(0,-1){8.75}}
%\end
%\emline(13.00,16.17)(13.00,23.75)
\put(13.00,16.17){\line(0,1){7.58}}
%\end
%\emline(14.00,23.33)(14.00,16.67)
\put(14.00,23.33){\line(0,-1){6.66}}
%\end
%\emline(15.00,17.33)(15.00,22.67)
\put(15.00,17.33){\line(0,1){5.34}}
%\end
%\emline(16.00,22.17)(16.00,17.83)
\put(16.00,22.17){\line(0,-1){4.34}}
%\end
%\emline(17.00,18.42)(17.00,21.67)
\put(17.00,18.42){\line(0,1){3.25}}
%\end
%\emline(18.00,21.08)(18.00,18.92)
\put(18.00,21.08){\line(0,-1){2.16}}
%\end
%\emline(19.00,19.50)(19.00,20.50)
\put(19.00,19.50){\line(0,1){1.00}}
%\end
\put(31.00,25.50){\makebox(0,0)[lc]{$\omega^{\vee}_2/2$}}
\put(30.00,36.67){\circle*{1.00}}
\put(31.50,36.33){\makebox(0,0)[lc]{$D_2$}}
\put(31.00,14.50){\makebox(0,0)[lc]{$\omega^{\vee}_1/2$}}
\put(30.00,3.33){\circle*{1.00}}
\put(31.50,3.67){\makebox(0,0)[lc]{$D_1$}}
\end{picture}
\end{center}
\end{figure}

The closure $X$ of the diagonal embedding
$\HS\embeds\PP^5\times(\PP^5)^{*}$ is called the \emph{space of
complete conics}. It is determined in $\PP^5\times(\PP^5)^{*}$ by
the equation ``$q\cdot q^{\vee}$ is s scalar matrix'', and this
implies by direct computations that $X$ is smooth. The unique
closed orbit $\widehat{Y}\subseteq X$ has the colored data
$\Vv_{\widehat{Y}}=\{-\omega^{\vee}_1/2,-\omega^{\vee}_2/2\}$,
$\Dd^B_{\widehat{Y}}=\emptyset$.

By Corollary~\ref{Pic(s.c)}, $\Pic\PP^5\iso\Pic(\PP^5)^{*}\iso\ZZ$
are freely generated by $D_1$, resp.~$D_2$, and $\Pic X\iso\ZZ^2$
is freely generated by $D_1,D_2$.
\end{example}
\begin{remark}
In fact, the space of smooth conics is a symmetric variety and the
space of complete conics is its ``wonderful completion'',
see~\ref{symmetric},~\ref{wonderful}.
\end{remark}

In the case of complexity~$1$, the description of Cartier divisors
is similar, but one should speak not only of cones, but also of
hypercones~$\Cc_Y$, and of admissible functionals $\lambda_Y$
which are integral on $\RG^{*}$ and such that
$\sum\langle\eps_x,\lambda_Y\rangle x$ is a principal divisor
on~$C$.

In particular, if $X$ is simple, then $\Pic X$ is generated by a
finite set $\Dd^B\setminus\Dd^B_Y$, where $Y\subseteq X$ is the
closed orbit, and Corollary~\ref{Pic(s.c)} is true.
\begin{example}
If $X$ is the space of complete triangles of Example~\ref{SL3/T},
then $\Pic X= \langle D_i,\widetilde{D}_i\mid
i=1,2,3\rangle\iso\ZZ^6$.
\end{example}

In characteristic zero, the $G$-module structure of the space of
global sections of a Cartier divisor is determined by the set of
$B$-eigensections. If $\sigma$ is a rational $B$-eigensection of
$\Lin{\delta}$ such that $\divr\sigma=\delta$, then
$\Ho^0(X,\Lin{\delta})^{(B)}=\{f\sigma\mid f\in K^{(B)},\ \divr
f+\delta\geq0\}$. The $B$-weight of $\eta=f\sigma\in
\Ho^0(X,\Lin{\delta})^{(B)}$ equals $\lambda+\pi(\delta)$, where
$\lambda$ is the $B$-weight of~$f$ and $\pi(\delta)$ is the
$B$-weight of~$\sigma$. The multiplicity of
$V(\lambda+\pi(\delta))$ in $\Ho^0(X,\Lin{\delta})$ equals
%*
\begin{equation}\label{mult(H^0)}
m_{\lambda}(\delta)=\dim\{f\in K^B\mid\divr
f+\divr\ef{\lambda}+\delta\geq0\}
\end{equation}
%*
The weight $\pi(\delta)$ is defined up to a character of~$G$ and
may be determined as follows. Consider a generic $G$-orbit
$Y\subseteq X$ and let $\widetilde{\delta}$ be the pull-back
on~$G$ of $\delta\cap Y$.  As $G$ is a factorial variety,
$\widetilde{\delta}$ is defined by an equation $F\in k(G)^{(B)}$.
Then $\pi(\delta)$ is the weight of~$F$.

In the case $c(X)\leq1$, the description of
$\Ho^0(X,\Lin{\delta})$ is given in the language of convex
geometry.

If $c(X)=0$, then the set of highest weights of
$\Ho^0(X,\Lin{\delta})$ is $\pi(\delta)+\Pp(\delta)\cap\RG$, where
%*
\begin{align}
\label{wt.pol} \Pp(\delta)&=\left\{\lambda\in\bigcap_{Y\subseteq
X} (-\lambda_Y+\Cc_Y^{\vee})\right.\left|\,\forall
D\in\Dd^B\setminus\bigcup_{Y\subseteq X}\Dd^B_Y:\
\langle\lambda,\res(D)\rangle+m_D\geq0\right\}
\end{align}
%*
and all highest weights occur with multiplicity~$1$.
\begin{example}
A Schubert divisor $D_{\alpha_i}\subseteq G/P_I$,
$\alpha_i\in\Pi\setminus I$, is defined by an equation $\langle
v,gv^{*}\rangle=0$, $v\in V(\omega_i^{*})$, $v^{*}\in
V(\omega_i)$. Hence $\pi(D_{\alpha_i})=\omega_i^{*}$. For
$\delta=\sum a_i D_{\alpha_i}$, we have $\Lin{\delta}= \ind{-\sum
a_i\omega_i}$, $\pi(\delta)=\sum a_i\omega_i^{*}$,
$\Pp(\delta)=\{0\}$, and $\Ho^0(G/P_I,\Lin{\delta})=V(\sum
a_i\omega_i^{*})$ (the Borel--Weil theorem, cf.~\ref{bundles}).
\end{example}
\begin{example}
Consider $X=\PP^{d-1}\times(\PP^{d-1})^{*}$ as a simple projective
embedding of a symmetric space
$\HS=\SL_d/\SG(\LO_1\times\LO_{d-1})$. Then $X\setminus\HS$ is a
homogeneous divisor consisting of all pairs $(x,y)$ such that the
point $x$ lies in the hyperplane~$y$. It is defined by an equation
$\sum x_iy_i=0$, where $x_1,\dots,x_d$ ($y_1,\dots,y_d$) are
projective coordinates on~$\PP^{d-1}$ (resp.~$(\PP^{d-1})^{*}$).
The two $B$-divisors $D,D'$ are defined by $B$-eigenfunctions
$y_1,x_d$ of biweights $(\omega_1,(d-1)\eps)$,
$(\omega_{d-1},(1-d)\eps)$, respectively, where $\eps$ generates
$\Ch(\SG(\LO_1\times\LO_{d-1}))$. One has
$\RG=\langle\omega_1+\omega_{d-1}\rangle\iso\ZZ$,
$\ef{\omega_1+\omega_{d-1}}(x,y)=x_dy_1/\sum x_iy_i$. It follows
easily that $\ES\iso\QQ\supset\Vv=\QQ_{-}$, $\res(D)=\res(D')=1$.

By Corollary~\ref{Pic(s.c)}, $\Pic(X)=\ZZ D\oplus\ZZ D'$, and for
$\delta=mD+nD'$ we have $\pi(\delta)=m\omega_1+n\omega_{d-1}$,
$\Pp(\delta)=\{\lambda=-k(\omega_1+\omega_{d-1})\mid
k,m-k,n-k\geq0\}$. On the other hand, $\Lin{\delta}=
\Lin{\PP^{d-1}}{n}\otimes\Lin{(\PP^{d-1})^{*}}{m}$, whence
$\Ho^0(X,\Lin{\delta})=\Sym^m\kk^d\otimes\Sym^n(\kk^d)^{*}=
V(m\omega_1)\otimes V(n\omega_{d-1})$. We obtain a decomposition
formula
%*
\begin{equation*}
V(m\omega_1)\otimes V(n\omega_{d-1})= \bigoplus\limits_{0\leq
k\leq\min(m,n)} V((m-k)\omega_1+(n-k)\omega_{d-1})
\end{equation*}
%*
\end{example}
For other applications to computing tensor product decompositions,
including Pieri formulae, see \cite[2.5]{Pic(sph)}.

Now assume $c(X)=1$. Put $\delta=\sum m_DD$ ($D$~runs through all
$B$-stable divisors on~$X$) and
$\res(D)=(h_d,\ell_D)\in\EHS{x_D}$, $x_d\in C$ (=the smooth
projective curve with the function field~$K^B$).

\begin{definition}
A \emph{pseudodivisor} on~$C$ is a formal linear combination
$\mu=\sum_{x\in C}m_x\cdot x$, where $m_x\in\RR\cup\{\pm\infty\}$,
and all but finitely many $m_x$ are~$0$. Put
$\Ho^0(C,\mu)=\{f\in\kk(C)\mid \divr{f}+\mu\geq0\}$. (Here we
assume $\forall c\in\RR:\ c+(\pm\infty)=\pm\infty$.)

If all $m_x\neq-\infty$, then $\Ho^0(C,\mu)$ is just the space of
global sections of the divisor $[\mu]=\sum[m_x]\cdot x$ on
$C\setminus\{x\mid m_x=+\infty\}$, otherwise $\Ho^0(C,\mu)=0$.
\end{definition}
Consider the pseudodivisor
%*
\begin{gather*}
\mu=\mu(\delta,\lambda)=\sum_{x\in C}\left(\min_{x_D=x}
\frac{\langle\lambda,\ell_D\rangle+m_D}{h_D}\right)x \\
\left(\text{Here we assume } \frac{c}0=
\begin{cases}
+\infty, & c\geq0\\
-\infty, & c<0
\end{cases}.\right)
\end{gather*}
%*
Since $\divr f=\sum h_Dv_{x_D}(f)\cdot D$, $\forall f\in K^B$ ,
and $\divr\ef{\lambda}=\sum\langle\lambda,\ell_D\rangle\cdot D$,
it follows from~\eqref{mult(H^0)} that
$m_{\lambda}(\delta)=h^0(\delta,\lambda):=\dim\Ho^0(C,\mu)$.

We have $h^0(\delta,\lambda)=0$ outside the polyhedral domain
%*
\begin{equation*}
\Pp(\delta)=\{\lambda\mid\langle\lambda,\ell_D\rangle\geq-m_D
\text{ for $\forall D$ such that $h_D=0$}\}
\subseteq\RG\otimes\RR.
\end{equation*}
%*
If there is $x\in C$ such that $x_D\neq x$ for all $D$ with
$h_D>0$, then $h^0(\delta,\lambda)=\infty$ for
$\forall\lambda\in\Pp(\delta)$, because in this case
$h^0(\delta,\lambda)$ is the dimension of the space of sections of
a divisor on an affine curve. Otherwise, by the Riemann--Roch
theorem,
%*
\begin{equation*}
\begin{split}
h^0(\delta,\lambda)&=\deg[\mu]-g+1+h^1(\delta,\lambda)\\
&=A(\delta,\lambda)-\sigma(\delta,\lambda)-g+1
+h^1(\delta,\lambda),
\end{split}
\end{equation*}
%*
where $g$ is the genus of~$C$, $h^1(\delta,\lambda)=\dim
\Ho^1(C,[\mu])$,
%*
\begin{equation*}
A(\delta,\lambda)=\sum_{x\in C}\left(\min_{x_D=x}
\frac{\langle\lambda,\ell_D\rangle+m_D}{h_D}\right)
\end{equation*}
%*
is a piecewise affine concave function of~$\lambda$, and
$\sigma(\delta,\lambda)$ is bounded non-negative for all
$\delta,\lambda$. Furthermore, as
$h^0(\delta,\lambda)\leq\deg[\mu]+1$ whenever $\deg[\mu]\geq0$
\cite[ex.IV.1.5]{AG}, we have $h^1(\delta,\lambda)\leq g$ if
$A(\delta,\lambda)\geq\sigma(\delta,\lambda)$. Note also that
$A(n\delta,n\lambda)=nA(\delta,\lambda)$.
\begin{question}
Maybe to consider an example $X=(\PP^{d-1})^3$ of complexity~$1$.
\end{question}

It follows that $h^0(\delta,\lambda)=0$ if $A(\delta,\lambda)<0$,
and $h^0(\delta,\lambda)=0$ differs from $A(\delta,\lambda)$ by a
globally bounded function whenever $A(\delta,\lambda)\geq0,$
$\lambda\in\Pp(\delta)$. This gives the asymptotic behaviors of
$h^0(\delta,\lambda)$ as $(\delta,\lambda)\to\infty$ in a fixed
direction.

Now we give criteria for a Cartier divisor to be globally
generated and ample.
\begin{theorem}\label{glob&ample}
Suppose $\delta$ is a Cartier divisor on~$X$ determined by the
data $\{f_Y\}$, $\{m_D\}$.
\begin{roster}
\item\label{glob} $\delta$~is globally generated iff local
equations $f_Y$ can be chosen in such a way that for any
$G$-subvariety $Y\subseteq X$ the following two conditions are
satisfied:
\begin{list}{}{}
\item[\reftag{a}] For any other $G$-subvariety $Y'\subseteq X$ and
each $B$-stable prime divisor $D\supseteq Y'$,\quad
$v_D(f_{Y})\leq v_D(f_{Y'})$. \item[\reftag{b}] $\forall
D\in\Dd^B\setminus\bigcup_{Y'\subseteq X}\Dd^B_{Y'}:\ v_D(f_Y)\leq
m_D$.
\end{list}
\item\label{ample} $\delta$ is ample iff, after replacing $\delta$
by a certain multiple, local equations $f_Y$ can be chosen in such
a way that, for any $G$-subvariety $Y\subseteq X$, there exists a
$B$-chart $\X$ of~$Y$ such that \reftag{a} and \reftag{b} are
satisfied and
\begin{list}{}{}
\item[\reftag{c}] the inequalities therein are strict iff
$D\cap\X=\emptyset$.
\end{list}
\end{roster}
\end{theorem}
\begin{proof}
\begin{roster}
\item[\ref{glob}] $\delta$~is globally generated iff for any
$G$-subvariety $Y\subseteq X$, there is $\eta\in
\Ho^0(X,\Lin{\delta})$ such that $\eta|_Y\neq0$. We may assume
$\eta$ to be a $B$-eigensection. This means that $\exists f\in
K^{(B)}:\ \divr{f}+\delta\geq0$, and no $D\supseteq Y$ occurs in
$\divr{f}+\delta$ with positive multiplicity. Replacing $f_Y$ by
$f^{-1}$ yields the conditions \reftag{a}--\reftag{b}. Conversely,
if \reftag{a}--\reftag{b} hold then $f=f_Y^{-1}$ yields the
desired global section. \item[\ref{ample}] Suppose $\delta$ is
ample. Replacing $\delta$ by a multiple, we may assume that
$\delta$ is very ample. Consider the $G$-equivariant projective
embedding $X\embeds\PP(M^{*})$ defined by a certain
finite-dimensional $G$-submodule
$M\subseteq\Ho^0(X,\Lin{\delta})$. Take a $G$-subvariety
$Y\subseteq X$. There exists a homogeneous $B$-eigenpolynomial in
homogeneous coordinates on $\PP(M^{*})$ (i.e., a section in
$\Ho^0(X,\Lin{\delta}^{\otimes N})^{(B)}$) that vanishes on
$\overline{X}\setminus X$ but not on~$Y$. Replacing $\delta$
by~$N\delta$, we may assume that
$\exists\eta\in\Ho^0(X,\Lin{\delta})^{(B)}:\
\eta|_{\overline{X}\setminus X}=0,\ \eta|_Y\neq0$. Then
$\X=X_{\eta}$ is a $B$-chart of~$Y$, and $\exists f\in K^{(B)}:\
\divr{f}+\delta=\divr\eta\geq0$. It remains to replace $f_Y$
by~$f^{-1}$.

Conversely, assume that the conditions \reftag{a}--\reftag{c}
hold. For any $G$-subvariety $Y\subseteq X$, there is a section
$\eta\in\Ho^0(X,\Lin{\delta})^{(B)}$ determined by $f_Y^{-1}$, and
$\X=X_{\eta}$ is a $B$-chart of~$Y$. We may pick finitely many
$B$-charts $\X_{\alpha}$ of this kind in such a way that
$G\X_{\alpha}$ cover~$X$. Let $\eta_{\alpha}\in
\Ho^0(X,\Lin{\delta})^{(B)}$ be the respective global sections.
Then
%*
\begin{equation*}
\kk[\X_{\alpha}]=\bigcup_{n\geq0}
\eta_{\alpha}^{-n}\Ho^0(X,\Lin{\delta}^{\otimes n})=\kk\left[
\frac{\sigma_{\alpha,1}}{\eta_{\alpha}^{n_{\alpha}}},\dots,
\frac{\sigma_{\alpha,s_{\alpha}}}{\eta_{\alpha}^{n_{\alpha}}}
\right]
\end{equation*}
%*
for some $n_{\alpha},s_{\alpha}\in\NN$, $\sigma_{\alpha,i}\in
\Ho^0(X,\Lin{\delta}^{\otimes n_{\alpha}})$.  Replacing $\delta$
by a multiple, we may assume $n_{\alpha}=1$.

Take the finite-dimensional $G$-submodule $M\subseteq
\Ho^0(X,\Lin{\delta})$ generated by
$\eta_{\alpha},\sigma_{\alpha,i}$, $\forall\alpha,i$. The
respective rational map $X\dasharrow\PP(M^{*})$ is $G$-equivariant
and defined on~$\X_{\alpha}$, hence everywhere. Moreover,
$\phi^{-1}(\PP(M^{*})_{\eta_{\alpha}})= \X_{\alpha}$ and
$\phi|_{\X_{\alpha}}$ is a closed embedding in
$\PP(M^{*})_{\eta_{\alpha}}$. Therefore $\phi$ is a locally closed
embedding and $\delta$ is very ample. \qedhere\end{roster}
\end{proof}
\begin{remark}
If $X$ is complete and $\delta$ is very ample, then the conditions
\reftag{a}--\reftag{c} hold for $\delta$ itself.
\end{remark}
\begin{corollary}\label{GX^0}
If $\X$ is a $B$-chart and $X=G\X$, then a divisor
$\sum_{D\subseteq X\setminus\X}m_DD$ is globally generated (ample)
iff all $m_D\geq0$ ($m_D>0$). In particular, $X$~is
quasiprojective.
\end{corollary}
\begin{corollary}
If $X$ is simple and $Y\subseteq X$ is the closed $G$-orbit, then
globally generated (ample) divisor classes in $\Pic X$ are those
$\delta=\sum_{D\in\Dd^B\setminus\Dd^B_Y}m_DD$ with $m_D\geq0$
($m_D>0$). In particular, any simple $G$-variety is
quasiprojective. (This also stems from Sumihiro's theorem.)
\end{corollary}

In the case of complexity $\leq1$, conditions
\reftag{a}--\reftag{b} mean that $\lambda_Y\leq\psi_{Y'}$ on
$\Cc_{Y'}$ and $\langle\lambda_Y,\res(D)\rangle\leq m_D$, and
\reftag{c} means that the inequalities therein are strict outside
$\Cc=\Cc(\Ww,\Rr)$ and~$\Rr$.

The description of globally generated and ample divisors on
spherical~$X$ in terms of piecewise linear functions is more
transparent if $X$ is complete (or all closed $G$-orbits
$Y\subseteq X$ are complete). Then maximal cones $\Cc_Y\in\Ff$ are
solid, and $\lambda_Y$ are determined by~$\psi_Y$.
\begin{definition}
A function $\psi\in\PL(\Ff)$ is (\emph{strictly}) \emph{convex} if
$\lambda_Y\leq\psi_{Y'}$ on $\Cc_{Y'}$ (resp.\
$\lambda_Y<\psi_{Y'}$ on $\Cc_{Y'}\setminus\Cc_Y$) for any two
maximal cones $\Cc_Y,\Cc_{Y'}\in\Ff$.
\end{definition}
\begin{corollary}\label{glob&ample(sph)}
If $X$ is complete (or all closed $G$-orbits in $X$ are complete)
and spherical, then $\delta$ is globally generated (ample) iff
$\psi$ is (strictly) convex on $\Ff$ and
$\langle\lambda_Y,\res(D)\rangle\leq m_D$ (resp.\ $<m_D$) for any
closed $G$-orbit $Y\subseteq X$ and $\forall
D\in\Dd^B\setminus\bigcup_{Y'\subseteq X}\Dd^B_{Y'}$.
\end{corollary}
\begin{proof}
It suffices to note that $Y$ has a unique $B$-chart $\X_Y$ given
by the colored cone $(\Cc_Y,\Dd^B_Y)$, and conditions
\reftag{a}--\reftag{c} are satisfied for $\delta$ iff they are
satisfied for its multiple.
\end{proof}
\begin{corollary}
On a complete spherical variety, every ample divisor is globally
generated.
\end{corollary}

The above results extend to the case of complexity~$1$, if all
closed $G$-orbits in $X$ are complete and of type~{\typeB}.

\begin{remark}
It follows from the proof of Theorem~\ref{glob&ample}\ref{ample}
that $\delta$ is very ample if $\kk[\X_{\alpha}]$ is generated by
$\eta_{\alpha}^{-1}\Ho^0(X,\Lin{\delta})$ for~$\forall\alpha$.
This may be effectively verified in some cases using
Lemma~\ref{G-qinv} for $R=\bigoplus_{n\geq0}
\Ho^0(X,\Lin{\delta}^{\otimes n})$, $S=\bigoplus_{n\geq0}
\bigl[\Ho^0(X,\Lin{\delta})\bigr]^n$, $\eta=\eta_{\alpha}$. For
example, if $X$ is complete and spherical, then it suffices to
verify that for each closed $G$-orbit $Y\subseteq X$ the
polyhedral domain $\lambda_Y+\Pp(\delta)$ contains the generators
of the semigroup $\Cc_Y^{\vee}\cap\RG$.
\end{remark}

\begin{example}
On a generalized flag variety $X=G/P$, globally generated (ample)
divisors are distinguished in the set of all $B$-divisors
$\delta=\sum a_iD_{\alpha_i}$ by the conditions $a_i\geq0$ (resp.\
$a_i>0$). Every ample divisor is very ample.
\end{example}
\begin{example}\label{non-proj}
The variety $X$ defined by the colored fan $\Ff$ from
Example~\ref{PL-fun} is complete, but not projective. Indeed,
since $\Cc_1\cap\Cc_2$ is a solid cone (Figure~\ref{SL(3)/SL(2)}),
a convex piecewise linear function on $\Ff$ is forced to be linear
on $\Cc_1\cup\Cc_2$, whence globally on~$\ES$. Hence there are no
non-principal globally generated divisors on~$X$.
\end{example}
\begin{remark}
If a fan $\Ff$ in a two dimensional colored space has no colors,
then the interiors of all cones in~$\Ff$ are disjoint and there
exists a strictly convex piecewise linear function on $\Ff$.
Therefore all toroidal spherical varieties (in particular, all
toric varieties) of rank~$2$ are quasiprojective. However, one can
construct a complete, but not projective, toric variety of
rank~$3$ \cite[p.71]{toric.intro}.
\end{remark}
\begin{example}
The same reasoning as in Example~\ref{non-proj} shows that an
$\SL_2$-embedding $X$ containing at least two $G$-germs of types
$\text{B}_{-}$, $\text{B}_0$ (Figure~\ref{germs(SL2)}) is not
quasiprojective. (Here $r(X)=1$, $\dim X=3$.) On the other hand,
if $X$ contains at most one $G$-germ of type $\text{B}_{-}$ or
$\text{B}_0$, then it is easy to construct a strictly convex
piecewise linear function on the hyperfan of~$X$, whence $X$ is
quasiprojective. (For smooth $X$, this was proved in
\cite[6.4]{SL(2)}.)
\end{example}

\section{Intersection theory}
\label{intersect}

Our basic reference in intersection theory is~\cite{int.theory}.
We begin our study of algebraic cycles on $G$-models with the
following general result reducing everything to $B$-stable cycles.
\begin{theorem}[\cite{int.sph}]\label{B-cycle}
Let a connected solvable algebraic group $B$ act on a variety~$X$.
Then the Chow group $A_d(X)$ is generated by the $B$-stable
$d$-cycles with the relations $[\divr f]=0$, where $f$ is a
rational $B$-eigenfunction on a $B$-stable $(d+1)$-subvariety
of~$X$.
\end{theorem}
\begin{proof}
Using the equivariant completion of $X$ and the equivariant Chow
lemma \cite[Th.1.3]{IT}, one reduces the assertion to the case of
projective $X$ by induction on $\dim X$ with the help of the
standard technique of exact sequences~\cite{int.sph}. The
projective case was handled by Vust \cite[6.1]{SL(2)} and Brion
\cite[1.3]{sph&Mori}. The idea is to consider the $B$-action on
the Chow variety $Z$ containing a given effective $d$-cycle~$z$.
Applying the Borel fixed point theorem, we find a $B$-stable cycle
$z_0\in\overline{Bz}$. An easy induction on $\dim B$ shows that
$z_0$ can be connected with $z$ by a sequence of rational curves,
whence is rationally equivalent to~$z$. The assertion on relations
is proved by a similar technique, see \cite[1.3]{sph&Mori} for
details.
\end{proof}
This theorem clarifies almost nothing in the structure of Chow
groups of general $G$-varieties, because the set of $B$-stable
cycles is almost as vast as the set of all cycles; however it is
very useful for $G$-varieties of complexity $\leq1$.

Assume $X$ is a unirational $G$-variety of complexity $\leq1$ or a
$B$-stable subvariety in it. (The assumption of unirationality is
needless in the spherical case, since $X$ has an open $B$-orbit.
If $c(X)=1$, then unirationality means that $K^B=\kk(\PP^1)$ for
$K=\kk(X)$.)
\begin{corollary}\label{f.g(Chow)}
$A_{*}(X)$ is finitely generated. If $U:X$ has finitely many
orbits, then $A_{*}(X)$ is freely generated by $U$-orbit closures.
\end{corollary}
\begin{proof}
If $c(X)=0$, then $B:X$ has finitely many orbits, whence
$A_{*}(X)$ is generated by $B$-orbit closures. If $c(X)=1$, then
by Theorem~\ref{cr(Y<X)}, each irreducible $B$-stable subvariety
$Y\subseteq X$ is either a $B$-orbit closure or the closure of a
one-parameter family of $B$-orbits. In the second case, $Y$~is one
of finitely many irreducible components of $X_k=\{x\in X\mid \dim
Bx\leq k\}$, $0\leq k<\dim X$, and it follows from
Lemma~\ref{B-eigenfun} that an open $B$-stable subset $\Y\subseteq
Y$ admits a geometric quotient $\Y/B$ which is a smooth rational
curve. Hence all $B$-orbits in $\Y$ are rationally equivalent and
each $B$-orbit, except finitely many of them, lies in one of $\Y$.
Therefore $A_{*}(X)$ is generated by finitely many $B$-orbit
closures and irreducible components of $X_k$.
\end{proof}
\begin{corollary}\label{A(comp)}
If $X$ is complete, then:
\begin{enumerate}
\item\label{A^+} The cone of effective cycles
$A_d^{+}(X)_{\QQ}\subseteq A_d(X)\otimes\QQ$ is a polyhedral cone
generated by the classes of rational subvarieties.
\item\label{alg=rat} Algebraic equivalence coincides with rational
equivalence of cycles on~$X$.
\end{enumerate}
\end{corollary}
\begin{proof}
\begin{roster}
\item[\ref{A^+}] Similar to the proof of Theorem~\ref{B-cycle}
using Corollary~\ref{f.g(Chow)}. \item[\ref{alg=rat}] The group of
cycles algebraically equivalent to~$0$ modulo rational equivalence
is divisible~\cite[19.1.2]{int.theory}. \qedhere\end{roster}
\end{proof}
\begin{corollary}\label{A=H}
If $X$ is smooth and complete (projective if $c(X)=1$), then the
cycle map $A_{*}(X)\to H_{*}(X)$ is an isomorphism of free Abelian
groups of finite rank. (Here $\kk=\CC$, but one may also consider
\'etale homology and Chow groups with corresponding coefficients
for arbitrary~$\kk$.)
\end{corollary}
\begin{proof}
If $c(X)=0$, then it is easy to deduce from Theorem~\ref{B-cycle}
that the K\"unneth map $A_{*}(X)\otimes A_{*}(Y)\to A_{*}(X\times
Y)$ is an isomorphism for~$\forall Y$, and the assertion follows
from the fact that $z=\sum(u_i\cdot z)v_i$ for $\forall z\in
A_{*}(X)$, where $\sum u_i\otimes v_i$ is the class of the
diagonal in~$X\times X$ \cite[\S3]{int.sph}. If $X$ is projective,
then one uses the Bia{\l}ynicki-Birula decomposition~\cite{BB}:
$X$~is covered by finitely many $B$-stable locally closed strata
$X_i$, where each $X_i$ is a vector bundle over a connected
component $X_i^T$ of $X^T$, and either $X_i^T=\pt$ or
$X_i^T=\PP^1$. This yields a cellular decomposition of~$X$, and we
conclude by~\cite[19.1.11]{int.theory}.
\end{proof}
\begin{remark}
The corollaries extend to an arbitrary variety $X$ with an action
of a connected solvable group $B$ having finitely many orbits.
\end{remark}
\begin{remark}
If $X$ is not unirational, then Corollaries~\ref{f.g(Chow)},
\ref{A(comp)}\ref{A^+} remain valid after replacing $A_{*}(X)$ by
the group $B_{*}(X)$ of cycles modulo algebraic equivalence.
\end{remark}
\begin{example}
If $X$ is a generalized flag variety or its Schubert subvariety,
then $A_{*}(X)\iso H_{*}(X)$ is freely generated by Schubert
subvarieties in~$X$.
\end{example}

Now we discuss intersection theory on varieties of
complexity~$\leq1$.

Let $X$ be a projective $G$-model of complexity $\leq1$, and
$\ch\kk=0$. A method to compute intersection numbers of Cartier
divisors on $X$ was introduced by Brion \cite[\S4]{Pic(sph)} in
the spherical case and generalized in \cite[\S8]{C-div} to the
case of complexity~$1$.

Put $\dim X=d$, $c(X)=c$ ($=0,1$), $r(X)=r$.

The N\'eron--Severi group $\NS(X)$ of Cartier divisors modulo
algebraic equivalence is finitely generated, and the intersection
form is a $d$-linear form on the finite-dimensional vector space
$\NS(X)_{\QQ}=\NS(X)\otimes\QQ$. This form is reconstructed via
polarization from the form $\delta\mapsto\deg_X\delta^d$ on
$\NS(X)_{\QQ}$ of degree~$d$. Moreover, each Cartier divisor
on~$X$ is a difference of two ample divisors, whence ample
divisors form an open solid convex cone in $\NS(X)_{\QQ}$, and the
intersection form is determined by values of $\deg\delta^d$ for
ample~$\delta$.

Retain the notation of~\ref{div}. Also put
$A(\delta,\lambda)\equiv1$, $\forall\lambda\in\ES^{*}$, if $c=0$,
and $\Pp_{+}(\delta)= \{\lambda\in\Pp(\delta)\mid
A(\delta,\lambda)\geq0\}$

\begin{theorem}\label{deg}
Suppose $\delta$ is an ample $B$-stable divisor on~$X$. Then
%*
\begin{align}
\label{dim=c+r+roots} d&=c+r+|\Delta_{+}^{\vee}\setminus
(\RG+\ZZ\pi(\delta))^{\ann}|+1, \qquad\text{and}\\
\label{deg=int}
\deg\delta^d&=d!\int\limits_{\pi(\delta)+\Pp_{+}(\delta)}
A(\delta,\lambda-\pi(\delta))
\prod_{\alpha^{\vee}\in\Delta_{+}^{\vee}\setminus
(\RG+\ZZ\pi(\delta))^{\ann}}
\frac{\langle\lambda,\alpha^{\vee}\rangle}
{\langle\rho,\alpha^{\vee}\rangle}\ d\lambda,
\end{align}
%*
where the Lebesgue measure on $\RG\otimes\RR$ is normalized so
that a fundamental parallelepiped of~$\RG$ has volume~$1$.
\end{theorem}
\begin{proof}
We have
%*
\begin{align*}
\dim\Ho^0(X,\Lin{\delta}^{\otimes n})&=\sum_{\lambda\in
n(\pi(\delta)+\Pp(\delta))\cap\RG}\dim V(\lambda)\cdot
m_{\lambda-n\pi(\delta)}(n\delta)\\
&=\sum_{\lambda\in (\pi(\delta)+\Pp(\delta))\cap\frac1n\RG}\dim
V(n\lambda)\cdot m_{n(\lambda-\pi(\delta))}(n\delta)\\
&=\sum_{\lambda\in(\pi(\delta)+\Pp_{+}(\delta))\cap\frac1n\RG}
\quad\prod_{\alpha^{\vee}\in\Delta_{+}^{\vee}}\left(1+
n\frac{\langle\lambda,\alpha^{\vee}\rangle}
{\langle\rho,\alpha^{\vee}\rangle}\right)
\qquad\text{if $c=0$,}\\
&\text{or}\sum_{\lambda\in(\pi(\delta)+\Pp_{+}(\delta))
\cap\frac1n\RG}\quad\prod_{\alpha^{\vee}\in\Delta_{+}^{\vee}}
\left(1+n\frac{\langle\lambda,\alpha^{\vee}\rangle}
{\langle\rho,\alpha^{\vee}\rangle}\right)
\bigl[nA(\delta,\lambda-\pi(\delta))\\
&\qquad\qquad-\sigma(n\delta,n(\lambda-\pi(\delta)))-g+1+
h^1(n\delta,n(\lambda-\pi(\delta)))\bigr]\\
\intertext{if $c=1$, using the Weyl dimension formula. In both
cases,} \dim\Ho^0(X,\Lin{\delta}^{\otimes n})&\sim
n^{c+r}\int\limits_{\pi(\delta)+\Pp_{+}(\delta)}
A(\delta,\lambda-\pi(\delta))
\prod_{\alpha^{\vee}\in\Delta_{+}^{\vee}\setminus
(\pi(\delta)+\Pp_{+}(\delta))^{\ann}}
n\frac{\langle\lambda,\alpha^{\vee}\rangle}
{\langle\rho,\alpha^{\vee}\rangle}\ d\lambda
\end{align*}
%*
On the other hand, the Euler characteristic
$\chi(\Lin{\delta}^{\otimes n})= \deg(\delta^d)n^d/d!+\dots$
equals $\dim \Ho^0(X,\Lin{\delta}^{\otimes n})$ for $n\gg0$. It
remains to note that $\Pp_{+}(\delta)$ generates $\RG\otimes\RR$,
because each rational $B$-eigenfunction on~$X$ is a quotient of
two $B$-eigensections of some $\Lin{\delta}^{\otimes n}$.
Therefore $(\pi(\delta)+\Pp_{+}(\delta))^{\ann}=
(\RG+\ZZ\pi(\delta))^{\ann}$.
\end{proof}
\begin{remark}
Formula~\eqref{dim=c+r+roots} may be proved using the local
structure theorem (Corollary~\ref{loc.str.gen}).
\end{remark}
\begin{remark}
Formula~\eqref{deg=int} is valid for globally generated~$\delta$,
because globally generated divisor classes lie on the boundary of
the cone of ample divisors in $\NS(X)_{\QQ}$ and the r.h.s.\
of~\eqref{deg=int} depends continuously on~$\delta$.
\end{remark}
\begin{remark}
The integral in the theorem can be easily computed using a
simplicial subdivision of the polyhedral domain $\Pp_{+}(\delta)$
and Brion's integration formula \cite[4.2, Rem.(ii)]{Pic(sph)}:
\begin{quote}
Suppose $F$ is a homogeneous polynomial of degree~$p$ on $\RR^r$,
and $[a_0,\dots,a_r]$ is a simplex with vertices $a_i\in\RR^r$.
Then
%*
\begin{align*}
\int\limits_{[a_0,\dots,a_r]}F(\lambda)\,d\lambda&=
\frac{r!\vol[a_0,\dots,a_r]}{(p+1)\dots(p+r)}
\Pi_rF(a_0,\dots,a_r)\\
\intertext{where}
\Pi_rF(a_0,\dots,a_r)&=\frac1{p!}\sum_{p_0+\dots+p_r=p}
\frac{\partial^pF(a_0t_0+\dots+a_rt_r)}{\partial{t_0}^{p_0}
\dots\partial{t_r}^{p_r}}
\end{align*}
%*
\end{quote}
\end{remark}
\begin{example}
For toric~$X$, $d=r$, $c=0$, and $\deg\delta^r=r!\vol\Pp(\delta)$
\cite[11.12.2]{toric.survey}.
\end{example}
\begin{example}
If $X=G/P_I$ is a generalized flag variety, then each ample
divisor $\delta$ defines an embedding $X\embeds V(\lambda)$,
$\lambda=\pi(\delta)^{*}$, and the degree of this embedding equals
%*
\begin{equation*}
|\Delta_{+}\setminus\Delta_{I,+}|!
\prod_{\alpha\in\Delta_{+}\setminus\Delta_{I,+}}
\frac{\langle\lambda,\alpha^{\vee}\rangle}
{\langle\rho,\alpha^{\vee}\rangle}
\end{equation*}
%*
In particular, the degree of the Pl\"ucker embedding
$\Gr_m(\kk^n)\embeds\PP(\E^m\kk^n)$ equals
%*
\begin{equation*}
[m(n-m)]!\frac{1!\dots(m-1)!}{(n-m)!\dots(n-1)!}
\end{equation*}
%*
(Schubert~\cite{deg(Plucker)}).
\end{example}
\begin{example}\label{int(conics)}
Consider the space $X$ of complete conics of Example~\ref{conics}.
Here $d=5$, $c=0$, $r=2$. If $\delta=a_1D_1+a_2D_2$ is an ample
divisor, then $\pi(\delta)=2a_1\omega_1+2a_2\omega_2$. Writing
$\lambda=-2x_1\alpha_1-2x_2\alpha_2$, we have
$d\lambda=dx_1\,dx_2$, and $\Pp(\delta)=\{\lambda\mid
x_1,x_2\geq0,\ 2x_1\leq x_2+a_1,\ 2x_2\leq x_1+a_2\}$ is a
quadrangle with the vertices
$\{0,-a_1\alpha_1,-a_2\alpha_2,-\pi(\delta)\}$.

We have $\pi(\delta)+\lambda=
(2a_1-4x_1+2x_2)\omega_1+(2a_2-4x_2+2x_1)\omega_2$, and
%*
\begin{equation*}
\begin{split}
\deg\delta^5&=5!\Bigl(\int_{\Pp(\delta)}
\tfrac{(2a_1-4x_1+2x_2)(2a_2-4x_2+2x_1)(2a_1+2a_2-2x_1-2x_2)}2\,
dx_1\,dx_2 \\
&=a_1^5+10a_1^4a_2+40a_1^3a_2^2+40a_1^2a_2^3+10a_1a_2^4+a_2^5
\end{split}
\end{equation*}
%*
Polarizing this $5$-form in $a_1,a_2$, we obtain the intersection
form on $\NS(X)=\langle D_1,D_2\rangle$: $\deg D_1^5=\deg
D_2^5=1$, $\deg D_1^4D_2=\deg D_1D_2^4=2$, $\deg D_1^3D_2^2=\deg
D_1^2D_2^3=4$ (Chasles~\cite{comp.con}).

This result can be applied to solving various enumerative problems
in the space $\HS$ of plane conics. For example, let us find the
number of conics tangent to 5~given conics in general position.
The set of conics tangent to a given one is a prime divisor
$D\subset\HS$. It is easy to see that (the closure of) $D$
intersects all $G$-orbits in $X$ properly. By Kleiman's
transversality theorem (see \cite[III.10.8]{AG} and below) five
general translates $g_iD$ ($g_i\in G$, $i=1,\dots,5$) are
transversal and intersect only inside~$\HS$. Thus the number we
are looking for equals $\deg_XD^5$.

Using local coordinates, one sees that the degree of (the closure
of) $D$ in $\PP^5$ or $(\PP^5)^{*}$ equals~$6$. (Take, e.g., a
parabola $\{y=x^2\}$. A conic $\{q(x,y)=0\}$ is tangent to this
parabola iff $q(x,x^2)=0$ and $\begin{vmatrix}
2x & \frac{\partial{q}}{\partial{x}}(x,x^2) \\
-1 & \frac{\partial{q}}{\partial{y}}(x,x^2)
\end{vmatrix}=0$
for some~$x$. The resultant of these two polynomials has degree
$7$ in the coefficients of~$q$. Cancelling it by the coefficient
at~$y^2$, we obtain the equation of $D$ of degree~$6$.) Since
$\deg_{\PP^5}D_1=\deg_{(\PP^5)^{*}}D_2=1$ and
$\deg_{\PP^5}D_2=\deg_{(\PP^5)^{*}}D_1=2$, one has $D\sim
2D_1+2D_2$ (on~$X$) and $\deg_XD^5=2^5(1+10+40+40+10+1)=3264$.
\end{example}

Spaces of conics and of quadrics in higher dimensions were studied
intensively from the origin of enumerative geometry
\cite{comp.con}, \cite{enum.geo}, \cite{enum.quad}. For a modern
approach, see \cite{comp.quad-1}, \cite{comp.quad-2},
\cite{comp.quad}, \cite{quadrics}, \cite{Pic(sph)}.

\begin{example}\label{int(triang)}
Let $X$ be a completion of the space $\HS$ of ordered triangles
from Example~\ref{SL3/T} with $d=6$, $c=1$, $r=2$. Consider an
ample divisor $\delta=a_1\widetilde{D}+a_2D$, where
$D=D_1+D_2+D_3$ imposes the condition that one of vertices of a
triangle lies on the $B$-stable line in~$\PP^2$ and
$\widetilde{D}=\widetilde{D}_1+\widetilde{D}_2+\widetilde{D}_3$
imposes the condition that a triangle passes through the $B$-fixed
point.

Writing $\lambda=-x_1\alpha_1-x_2\alpha_2$, we have
$d\lambda=dx_1\,dx_2$, $\Pp(\delta)=\{\lambda\mid x_1,x_2\geq0\}$,
$\pi(\delta)=3a_1\omega_1+ 3a_2\omega_2$, and
%*
\begin{equation*}
A(\delta,\lambda)=
\begin{cases}
A_0(\lambda)=x_1+x_2,      & x_i\leq a_i \\
A_1(\lambda)=3a_1-2x_1+x_2,& 0\leq x_1-a_1\geq x_2-a_2 \\
A_2(\lambda)=3a_2-2x_2+x_1, & 0\leq x_2-a_2\geq x_1-a_1
\end{cases}
\end{equation*}
%*
It follows that $\Pp_{+}(\delta)=\{\,\lambda\mid x_1,x_2\geq0;\
2x_1\leq x_2+3a_1;\ 2x_2\leq x_1+3a_2\,\}=
\Pp_0\cup\Pp_1\cup\Pp_2$, where $\Pp_i$ are quadrangles with
vertices
%*
\begin{gather*}
\{0,-a_1\alpha_1,-a_1\alpha_1-a_2\alpha_2,-a_2\alpha_2\},\\
\{-a_1\alpha_1,-\tfrac{3a_1}2\alpha_1,-\pi(\delta),
-a_1\alpha_1-a_2\alpha_2\},\\
\{-a_2\alpha_2,-\tfrac{3a_2}2\alpha_2,-\pi(\delta),
-a_1\alpha_1-a_2\alpha_2\},
\end{gather*}
%*
and $A(\delta,\lambda)=A_i(\lambda)$ on~$\Pp_i$. We have
$\pi(\delta)+\lambda=
(3a_1-2x_1+x_2)\omega_1+(3a_2-2x_2+x_1)\omega_2$, and
%*
\begin{equation*}
\begin{split}
\deg\delta^6&=6!\Bigl(\int_{\Pp_0}(x_1+x_2)\cdot
\tfrac{(3a_1-2x_1+x_2)(3a_2-2x_2+x_1)(3a_1+3a_2-x_1-x_2)}2\,
dx_1\,dx_2 \\
&\qquad+\int_{\Pp_1}(3a_1-2x_1+x_2)\cdot
\tfrac{(3a_1-2x_1+x_2)(3a_2-2x_2+x_1)(3a_1+3a_2-x_1-x_2)}2\,
dx_1\,dx_2 \\
&\qquad+\int_{\Pp_2}(3a_2-2x_2+x_1)\cdot
\tfrac{(3a_1-2x_1+x_2)(3a_2-2x_2+x_1)(3a_1+3a_2-x_1-x_2)}2\,
dx_1\,dx_2\Bigr) \\
&=90a_1^6+1080a_1^5a_2+4320a_1^4a_2^2+6840a_1^3a_2^3+
  4320a_1^2a_2^4+1080a_1a_2^5+90a_2^6
\end{split}
\end{equation*}
%*
It follows that $\deg D^6=\deg\widetilde{D}^6=90$, $\deg
D^5\widetilde{D}=\deg D\widetilde{D}^5=180$, $\deg
D^4\widetilde{D}^2=\deg D^2\widetilde{D}^4=288$, $\deg
D^3\widetilde{D}^3=342$.

Since $X$ has finitely many orbits and $D,\widetilde{D}$ intersect
all of them properly, it follows from Kleiman's transversality
theorem that any 6~general translates $\delta_i$ of
$D,\widetilde{D}$ are transversal and intersect only inside~$\HS$.
Thus the number of common points of $\delta_i$ in $\HS$ equals
$\deg_X(\delta_1\dots\delta_6)$. Dividing it by~6 (=the number of
ordered triangles corresponding to a given unordered triangle), we
obtain the number of triangles satisfying 6~conditions imposed
by~$\delta_i$. For example, there are
$\deg(D^3\widetilde{D}^3)/6=57$ triangles passing through $3$
given points in general position whose vertices lie on $3$ given
general lines.
\end{example}

Theorem~\ref{deg} was applied in \cite[\S10]{C-div} to computing
the degree of a closed $3$-dimensional orbit in any
$\SL_2$-module.

Brion~\cite[4.1]{geo(sph)} proved a formula similar
to~\eqref{deg=int} for the multiplicity of a spherical variety
along an orbit in it and deduced a criterion of smoothness for
spherical varieties \cite[4.2]{geo(sph)}.

For any complete variety $X$, there is a canonical pairing $\Pic
X\times A_1(X)\to\ZZ$ given by the degree of a line bundle
restricted to a curve in~$X$ (and pulled back to its
normalization). The following theorem is essentially due to
Brion~\cite{sph&Mori}.
\begin{question}
$\ch\kk=0$
\end{question}
\begin{theorem}
\begin{roster}
\item\label{Pic<A_1^*} If $X$ is a complete unirational $G$-model
of complexity $\leq1$, then $\Pic X\embeds
A_1(X)^{*}=\Hom(A_1(X),\ZZ)$ via the canonical pairing.
\item\label{Pic=A_1^*} If $X$ is complete and spherical, then
$\Pic X\isoto A_1(X)^{*}$. \item\label{duality} If in addition $X$
contains a unique closed $G$-orbit~$Y$, then $A_1(X)$ is
torsion-free, and the basis of $A_1(X)$ dual to the basis
$\Dd^B\setminus\Dd^B_Y$ of $\Pic X$ consists of (classes of)
irreducible $B$-stable curves. Moreover, these basic curves
generate the semigroup $A_1^{+}(X)$ of effective $1$-cycles.
\end{roster}
\end{theorem}
\begin{proof}
Using the equivariant Chow lemma and resolution of singularities,
we construct a proper birational $G$-morphism $\phi:\widehat{X}\to
X$, where $\widehat{X}$ is a smooth projective $G$-variety. In the
commutative diagram
%*
\begin{equation*}
\begin{CD}
    \Pic X      @>>>     A_1(X)^{*}       \\
     @VVV                    @VVV            \\
\Pic\widehat{X} @>>> A_1(\widehat{X})^{*}
\end{CD}
\end{equation*}
%*
the vertical arrows are injections, and the bottom arrow is an
isomorphism by Corollary~\ref{A=H} and by Poincar{\`e} duality,
whence~\ref{Pic<A_1^*}. Assertions~\ref{Pic<A_1^*},~\ref{duality}
are proved in \cite[\S3]{B-curves} using the description of
$B$-stable curves and their equivalences on spherical varieties
obtained in \cite{sph&Mori}, \cite{B-curves}.
\end{proof}
\begin{remark}[{\cite[1.6, 2.1]{sph&Mori}}]
On a spherical $G$-model~$X$, any line bundle $\Ll$ is
$G$-linearized and any $B$-stable curve $C$ is the closure of a
$1$-dimensional $B$-orbit. Let $\infty$ be a $B$-fixed point in
the normalization $\PP^1$ of~$C$ and
$0\in\AAA^1=\PP^1\setminus\{\infty\}$ be another $T$-fixed point.
Then $T$ acts on $\AAA^1\setminus\{0\}$ via a character
$\chi\neq0$. Let $x,y\in C$ be the images of $0,\infty$ under the
normalization map $\nu:\PP^1\to C$, and $\chi_x,\chi_y$ be the
weights of the $T$-action on $\Ll_x,\Ll_y$. Then $\chi_x-\chi_y$
is a multiple of~$\chi$, and
$\langle\Ll,C\rangle=\deg\nu^{*}\Ll|_C=(\chi_x-\chi_y)/\chi$.
\end{remark}
\begin{example}
For a generalized flag variety $X=G/P_I$, $\Pic X$ is freely
generated by Schubert divisors
$D_{\alpha_i}=\overline{B[w_Gs_i]}$, and $A_1(X)$ is freely
generated by Schubert curves $C_{\alpha_i}=\overline{B[s_i]}\iso
P_{\alpha_i}/B\iso\PP^1$. We have
$\Lin{D_{\alpha_i}}=G\itimes{P_I}\kk_{-\omega_i}$ and
$\Lin{D_{\alpha_i}}|_{C_{\alpha_j}}=
P_{\alpha_i}\itimes{B}\kk_{-\omega_i}=\Lin{1}$ if $i=j$ and
$\Lin{0}$ if $i\neq j$. Here the $T$-fixed points are
$[s_i],[\1]$, $\chi=\alpha_i$, and $\langle D_{\alpha_i},
C_{\alpha_j}\rangle=(\chi_{[s_i]}-\chi_{[\1]})/\chi=
(-s_i\omega_j+\omega_j)/\alpha_i=\delta_{ij}$. Hence the above
bases of $\Pic X$ and $A_1(X)$ are dual to each other.
\end{example}

Projective unirational normal varieties of complexity $\leq1$ are
well-behaved from the point of view of the Mori
theory~\cite{Mori}.
\begin{theorem}[\cite{sph&Mori}, \cite{Mori(c=1)}]
Suppose $X$ is a projective unirational $G$-model of complexity
$\leq1$. Then the cone $\NE(X)$ of effective $1$-cycles modulo
numerical equivalence is finitely generated by rational $B$-stable
curves and all its faces are contractible. If $X$ is
$\QQ$-factorial, then each contraction of an extremal ray of
$\NE(X)$ isomorphic in codimension~$1$ can be flipped, and every
sequence of directed flips terminates.
\end{theorem}

Explicit computations of Chow rings for some smooth completions of
classical homogeneous spaces were carried on by several authors.
Schubert, Pieri, Giambelli, A.~Borel, Kostant,
Bernstein--Gelfand--Gelfand, Demazure,
Lakshmibai--Musili--Seshadri et al contributed to computing Chow
(or cohomology) rings of generalized flag varieties.

Here and below we put $A^k(X)=A_{d-k}(X)$, $d=\dim X$.

Without loss of generality, assume that $G$ is semisimple simply
connected. Let $\Ch=\Ch(B)$ be the weight lattice of~$G$. Every
$\lambda\in\Ch$ defines an induced line bundle
$\ind{-\lambda}=G\itimes{B}\kk_{-\lambda}$ on~$G/B$, and this
gives rise to an isomorphism $\Ch\isoto\Pic G/B\iso A^1(G/B)$. Put
$S=\Sym^{\bullet}(\Ch\otimes\QQ)$.
\begin{theorem}[\cite{H(h.s)}, \cite{desing.Sch}]
\begin{roster}
\item\textup{\cite{H(h.s)}, \cite{desing.Sch}}
$A^{*}(G/B)_{\QQ}\iso S/SS^W_{+}$ (the quotient modulo the ideal
generated by $W$-invariants without constant term)
\item\textup{\cite[5.5]{BGG}} If $I\subseteq\Pi$, then
$A^{*}(G/P_I)_{\QQ}$ embeds in $A^{*}(G/B)_{\QQ}$ as
$S^{W_I}/S^{W_I}S^W_{+}$.
\end{roster}
\end{theorem}

Bernstein--Gelfand--Gelfand \cite{BGG} and
Demazure~\cite{desing.Sch} used divided difference operators to
introduce certain functionals $D_w$ on $S$ which represent
Schubert cells $S_w$ ($w\in W$) via the Poincar\`e duality. They
also found the basis of $S/SS^W_{+}$ dual to~$D_w$.

Chow rings of toric varieties were computed by Jurkiewicz and
Danilov, cf.~\cite[\S10]{toric.survey}. Namely, if $X$ is a smooth
complete toric variety, then $A^{*}(X)= S^{*}_{\ZZ}(\Pic X)/I$,
where the ideal $I$ is generated by monomials $[D_1]\dots[D_k]$
such that $D_i$ are $T$-stable prime divisors and
$D_1\cap\dots\cap D_k=\emptyset$.
\begin{question}
Prove it?
\end{question}

The above examples are spherical (see
also~\ref{monoids},~\ref{toroidal}).
\begin{question}
(Equivariant) Chow rings of smooth toroidal varieties.
\end{question}
In the case of complexity~$1$, Chow rings of complete
$\SL_2$-embeddings (cf.~Example~\ref{SL2/e}) were computed
in~\cite{SL(2)}. The space of complete triangles, which is a
desingularization of the space $X$ of Example~\ref{SL3/T}, was
studied in~\cite{triang}, and in particular, its Chow ring was
determined there.

Many enumerative problems arise on non-complete homogeneous
spaces. Given a homogeneous space $\HS=G/H$, typically a space of
geometric figures or tensors of certain type, one looks for the
number of points satisfying a number of conditions in general
position. The set of points satisfying a given condition is a
closed subvariety $Z\subset\HS$, and the configuration of
conditions $Z_1,\dots,Z_s\subset\HS$ is put in general position by
replacing $Z_i$ by their translates~$g_iZ_i$, where
$(g_1,\dots,g_s)\in G\times\dots\times G$ is a general $s$-tuple.
By Kleiman's transversality theorem \cite{trans},
\cite[III.10.8]{AG},
\begin{question}
Prove it?
\end{question}
the cycles $g_1Z_1,\dots,g_sZ_s$ intersect transversally in smooth
subvarieties of codimension $\sum\codim Z_i$, i.e.,
$g_1Z_1\cap\dots\cap g_sZ_s$ is empty if $\sum\codim Z_i>\dim\HS$
and finite if $\sum\codim Z_i=\dim\HS$, and the cardinality
$|g_1Z_1\cap\dots\cap g_sZ_s|$ is stable for general
$(g_1,\dots,g_s)$. Thus the natural intersection ring for the
enumerative geometry of $\HS$ is provided by the following
\begin{definition}[\cite{C^*(sym)}]
The \emph{intersection number} of irreducible subvarieties
$Z_1,\dots,Z_s$ whose codimensions sum up to~$\dim\HS$ is
$(Z_1\cdot\dots\cdot Z_s)_{\HS}=|g_1Z_1\cap\dots\cap g_sZ_s|$ for
all $(g_1,\dots,g_s)$ in a dense open subset of
$G\times\dots\times G$.

This defines a pairing between groups of cycles in $\HS$ of
complementary dimensions. The \emph{group of conditions}
$C^{*}(\HS)$ is the quotient of the group of all cycles modulo the
kernel of this pairing. Write $\[Z\]$ for the image in
$C^{*}(\HS)$ of a cycle~$Z$.
\end{definition}
\begin{theorem}[{\cite[6.3]{C^*(sym)}}]
\label{C^*(sph)} If $\HS$ is spherical, then $C^{*}(\HS)$ is a
graded ring w.r.t.\ the \emph{intersection product}
$\[Z\]\cdot\[Z'\]=\[gZ\cap g'Z'\]$, where $(g,g')\in G\times G$ is
a general pair. Furthermore, $C^{*}(\HS)=\varinjlim A^{*}(X)$ over
all smooth complete $G$-embeddings $X\embof\HS$.
\end{theorem}
\begin{proof}
The proof goes in several steps.
\begin{roster}
\item\label{prop.int} For any subvariety $Z\subset\HS$ and any
smooth complete $G$-embedding $X\embof\HS$, there is a smooth
complete $G$-embedding $X'\embof\HS$ dominating $X$ such that
$\overline{Z}$ intersects all $G$-orbits in $X'$ properly. First,
one constructs a smooth toroidal embedding $X'$ dominating $X$.
Then each $G$-orbit on $X'$ is a normal intersection of $G$-stable
prime divisors (Theorem~\ref{reg.var}), and one applies a general
result \cite[4.7]{C^*(sym)} that a cycle on a complete smooth
variety can be put in regular position w.r.t.\ a regular
configuration of hypersurfaces by blowing up several intersections
of pairs of these hypersurfaces. \item\label{int.prod} If
$Z,Z'\subset\HS$ intersect all $G$-orbits in $X\setminus\HS$
properly, then $[gZ\cap g'Z']=[Z]\cdot[Z']$ in $A^{*}(X)$ for
general $g,g'\in G$. Indeed, we may apply Kleiman's transversality
theorem to intersections of $Z,Z'$ with each of finitely many
$G$-orbits in~$X$ and deduce that $gZ$ and $g'Z'$ intersect
properly with each other and $gZ\cap g'Z'$ intersects
$X\setminus\HS$ properly. \item For any $z\in A^{*}(X)$, use the
Chow moving lemma to represent it as $z=\sum m_i[Z_i]$, where
$Z_i\subset X$ are closed subvarieties intersecting~$\HS$. For any
subvariety $Z'\subset\HS$ of complementary dimension, we may
assume by~\ref{prop.int} that $Z_i,Z'$ intersect all orbits in
$X\setminus\HS$ properly and deduce from~\ref{int.prod} that
$(\[Z'\],\sum m_i\[Z_i\])_{\HS}=\deg_X[Z']\cdot z$ depends only
on~$z$. Thus we have a well-defined map $A^{*}(X)\to C^{*}(\HS)$,
$z\mapsto\sum m_i\[Z_i\]$. \item This map gives rise to a
homomorphism $\varinjlim A^{*}(X)\to C^{*}(\HS)$
by~\ref{int.prod}. Its surjectivity is obvious, and injectivity
follows from the Poincar\`e duality on~$X$. \qedhere\end{roster}
\end{proof}

The ring of conditions $C^{*}(\HS)$ is also called \emph{Halphen
ring} in honor of \mbox{G.-H.~Halphen}, who used it in the
enumerative geometry of conics, see~\cite{C^*(con)}. If $\HS$ is a
torus, then $C^{*}(\HS)$ is MacMullen's polytope algebra
\cite{C^*(tor)}, \cite[3.3]{conv.pol}.
\begin{question}
Genaralize to spherical $\HS$ using cohomology of toroidal
completions.
\end{question}
The Halphen ring of the space of plane conics was computed
in~\cite{C^*(con)}.

Theorem~\ref{C^*(sph)} reflects an idea exploited already by
classics that in solving enumerative problems on~$\HS$, one has to
consider an appropriate completion $X\embof\HS$ with finitely many
orbits such that all ``conditions'' $Z_i$ under consideration
intersect all orbits in $X\setminus\HS$ properly. If $\sum\codim
Z_i=\dim\HS$, then for general~$g_i$, all intersection points of
$\bigcap g_iZ_i$ lie in~$\HS$, and the intersection number equals
$\deg_X\prod[Z_i]$. Applications of this idea can be found in
Examples~\ref{int(conics)}, \ref{int(triang)}.

Generalizing these examples, we describe a method to compute the
intersection number of $d$~divisors on a spherical homogeneous
space $\HS$ of dimension~$d$ in characteristic~$0$. Let $\delta$
be a divisor on~$\HS$. Replacing $\delta$ by a $G$-translate, we
may assume that $\delta$ contains no colors. Let $h\in K$ be an
equation of $\delta$ on the open $B$-orbit in~$\HS$, which is a
factorial variety.
\begin{definition}
The \emph{Newton polytope} of~$\delta$ is the set
%*
\begin{align*}
\Nn(\delta)&=\{\,\lambda\in\ES^{*}\mid\forall v\in\Vv:\ \langle
v,\lambda\rangle\geq v(h),\ \forall D\in\Dd^B:\
\langle\res(D),\lambda\rangle\geq v_D(h)\,\}
\end{align*}
%*
\end{definition}
\begin{remark}
We see below that $\Nn(\delta)$ is indeed a convex polytope
in~$\ES^{*}$. If $G=\HS=T$ is a torus, then $\Dd^B=\emptyset$,
$h=\sum c_i\lambda_i$, $c_i\in\kk^{\times}$, $\lambda_i\in\Ch(T)$,
and $v(h)=-\min\langle v,\lambda_i\rangle$. Thus
$\Nn(\delta)=-\conv\{\lambda_1,\dots,\lambda_s\}$ is a usual
Newton polytope.
\end{remark}

For any embedding $X\embof\HS$, we have $\divr h=\delta-\delta_X$
on~$X$, where $\delta_X=-\sum_iv_i(h)V_i- \sum_{D\in\Dd^B}v_D(h)D$
is a $B$-stable divisor and $V_i$ are $G$-stable prime divisors on
$X$ with valuations $v_i\in\Vv$.
\begin{theorem}\label{int.num}
$\Nn(\delta)=\bigcap_{X\embof\HS}\Pp(\delta_X)$. If $X$ is
complete and $\delta$ intersects all $G$-orbits in $X$ properly,
then $\Nn(\delta)=\Pp(\delta_X)$.
\end{theorem}
\begin{proof}
The first assertion is obvious from~\eqref{wt.pol}, since every
$G$-valuation corresponds to a divisor on some embedding of~$X$.
Suppose $X$ is complete and $\delta$ intersects all orbits
properly. Consider a $G$-linearized line bundle $\Ll=\Lin{\delta}
=\Lin{\delta_X}$. We may extend any $G$-valuation to a
$G$-valuation of $\bigoplus_{n\geq0}\Ho^0(X,\Ll^n)$
(Corollary~\ref{res&ext(G-val)}). If $\sigma$ is a rational
section of $\Ll$ such that $\divr\sigma=\delta$, then
$v(\sigma)=0$ for $\forall v\in\Vv$, because $v$ has a center on
$X$ and $\delta$ intersects it properly. Now
$\forall\lambda\in\Pp(\delta_X)\ \exists n\in\NN\ \exists\eta\in
\Ho^0(X,\Ll^n)^{(B)}_{n\lambda}
\implies\eta=c\ef{n\lambda}\sigma_X^n=c\ef{n\lambda}\sigma^n/h^n$,
where $c\in\kk^{\times}$, $\divr\sigma_X=\delta_X$. Then
$v(\eta)=v(\ef{n\lambda}/h^n)\geq0\implies\langle
v,\lambda\rangle\geq v_D(h)\implies
\Pp(\delta_X)\subseteq\Nn(\delta)$.
\end{proof}
\begin{corollary}[{\cite[4.2]{conv.pol}}]
\label{Bezout} For any effective divisor $\delta$ on~$\HS$,
%*
\begin{equation}
(\delta^d)_{\HS}=d!\int\limits_{\pi(\delta)+\Nn(\delta)}
\prod_{\alpha^{\vee}\in\Delta_{+}^{\vee}\setminus
(\RG+\ZZ\pi(\delta))^{\ann}}
\frac{\langle\lambda,\alpha^{\vee}\rangle}
{\langle\rho,\alpha^{\vee}\rangle}\ d\lambda,
\end{equation}
%*
where $\pi(\delta)=-\sum_{D\in\Dd^B}v_D(h)\pi(D)$
\end{corollary}
\begin{proof}
Follows from Theorems~\ref{C^*(sph)}, \ref{int.num},
and~\ref{deg}.
\end{proof}
\begin{remark}
In the toric case, the above formula transforms to
$(\delta^d)_{\HS}=d!\vol\Nn(\delta)$. Polarizing this formula, we
obtain a theorem of Bernstein~\cite{num(roots)} and
Kouchnirenko~\cite{Newton.pol}: for any effective divisors
$\delta_1,\dots,\delta_d$ on~$\HS$, the intersection number
$(\delta_1\dots\delta_d)_{\HS}$ is $d!$ times the mixed volume of
$\Nn(\delta_1),\dots,\Nn(\delta_d)$. In the general case, we have
a ``mixed integral'' instead.
\end{remark}
Corollary~\ref{Bezout} may be considered as a generalization of
the classical B\'ezout theorem.

\chapter{Invariant valuations}
\label{inv.val}

This chapter plays a significant, but auxiliary, role in the
general context of our survey. We investigate the set of
$G$-invariant valuations of the function field of a $G$-variety.
We have seen in Chapter~\ref{LV-theory} that $G$-valuations are of
importance in the embedding theory, because they provide a
material for constructing combinatorial objects (colored data)
that describe equivariant embeddings.

Remarkably, a $G$-valuation of a given $G$-field is uniquely
determined by its restriction to the multiplicative group of
$B$-eigenfunctions, the latter being a direct product of the
weight lattice and of the multiplicative group of $B$-invariant
functions. Thus a $G$-valuation is essentially a pair composed by
a linear functional on the weight lattice and by a valuation of
the field of $B$-invariants. Under these identifications, we prove
in~\ref{val.cones} that the set of $G$-valuations is a union of
convex polyhedral cones in certain half-spaces.

The common face of these valuation cones is formed by those
valuations, called central, that vanish on $B$-invariant
functions. The central valuation cone controls the situation
``over the field of $B$-invariant functions''. For instance, its
linear part determines the group of $G$-automorphisms acting
identically on $B$-invariants.

This cone has another remarkable property: it is a fundamental
chamber of a crystallographic reflection group called the little
Weyl group of a $G$-variety. This group is defined in~\ref{Weyl}
as the Galois group of a certain symplectic covering of the
cotangent bundle constructed in terms of the moment map. The
little Weyl group is linked with the central valuation cone via
the invariant collective motion on the cotangent variety, which is
studied in~\ref{inv.mot}.

For practical applications, we must be able to compute the set of
$G$-valuations. For central valuations, it suffices to know the
little Weyl group. In~\ref{form.curv} we describe the ``method of
formal curves'' for computing $G$-valuations on a homogeneous
space. Informally, one computes the order of functions at infinity
along a formal curve approaching to a boundary $G$-divisor.

Most of the results of this chapter are due to D.~Luna and
Th.~Vust, M.~Brion, F.~Pauer, and F.~Knop. We follow \cite{LV},
\cite{G-val}, \cite{inv.mot} in our exposition.

\section{$G$-valuations}
\label{G-vals}

An algebraic counterpart of a prime divisor on an algebraic
variety is the respective valuation of the field of rational
functions. Valuations obtained in this way are called geometric
(see Appendix~\ref{geom.val}). We consider invariant geometric
valuations.

Let $G$ be a connected algebraic group and $K$ a $G$-field, i.e.,
the function field of a $G$-variety.
\begin{definition}
A \emph{$G$-valuation} is a $G$-invariant geometric valuation
of~$K/\kk$. The set of $G$-valuations is denoted by $\Vv=\Vv(K)$.
\end{definition}

The following approximation result is due to Sumihiro.
\begin{proposition}[{\cite[\S4]{eq.compl}}]
\label{val->G-val} For any geometric valuation $v$ of $K$ there
exists a $G$-valuation $\overline{v}$ such that $\forall f\in K:\
\overline{v}(f)=v(gf)$ for general $g\in G$. If $A\subset K$ is a
rational $G$-algebra, then $\forall f\in A:\
\overline{v}(f)=\min_{g\in G}v(gf)$.
\end{proposition}
\begin{proof}
We may assume $v=v_D$ for a prime divisor $D$ on a model $X$
of~$K$. Then $v'=v_{G\times D}$ is a geometric valuation
of~$\kk(G\times X)$. It is clear that $\forall f\in\kk(G\times
X):\ v'(f)=v(f(g,\cdot))$ for general $g\in G$. The rational
action $G:X$ induces an embedding $\kk(X)\embeds\kk(G\times X)$.
It is easy to see that $v'|_{\kk(X)}=\overline{v}$ is the desired
$G$-valuation.

To prove the \ordinal{2} assertion, observe that $A_{(d)}=\{f\in
A\mid v(f)\geq d\}$ is a filtration of $A$ by linear subspaces,
$Gf$~is an algebraic variety and $Gf\cap A_{(d)}$ its closed
subvariety, $\forall d\in\QQ$.
\end{proof}
\begin{remark}
If $v$ has center $Y\subseteq X$, then $\overline{v}$ has center
$\overline{GY}$.
\end{remark}
\begin{example}
Let $G=\kk$ act rationally on the blow-up $X$ of $\AAA^2$ at~$0$
by translations along a fixed axis. In coordinates,
$u(x,y)=(x+u,y)$, $\forall u\in G$, $(x,y)\in\AAA^2$. The
valuation $v$ of $\kk(X)$ corresponding to the exceptional divisor
is given on $\kk[\AAA^2]=\kk[x,y]$ by the order of a polynomial
in~$x,y$ (i.e., the lowest degree of a homogeneous term) and has
center $0$ on~$\AAA^2$. But $\overline{v}(f)=\min_uv(f(x+u,y))$ is
the order of $f$ in~$y$, so that $\overline{v}=v_D$, where
$D=\{y=0\}$ is (the proper pullback of) the $x$-axis.
\end{example}

Together with Proposition~\ref{res&ext(val)}, Sumihiro's
approximation immediately implies
\begin{corollary}\label{res&ext(G-val)}
Let $K'\subseteq K$ be a $G$-subfield. The restriction of a
$G$-valuation of $K$ to $K'$ is a $G$-valuation, and any
$G$-valuation of $K'$ can be extended to a $G$-valuation of~$K$.
\end{corollary}

The next corollary is useful in applications.
\begin{corollary}\label{v(g*sec/sec)}
Let $X$ be a $G$-model of~$K$ and $\Ll$ a $G$-line bundle on~$X$.
Then $\forall\sigma,\eta\in\Ho^0(X,\Ll)$, $\eta\ne0$, $\forall
g\in G:\ v(\sigma/\eta)=v(g\sigma/\eta)$.
\end{corollary}
\begin{proof}
Consider a rational $G$-algebra
$R=\bigoplus_{n\geq0}\Ho^0(X,\Ll^n)$. Then $\Quot{R}=K'(\eta)$ is
a (purely transcendental) extension of a $G$-subfield $K'\subseteq
K$ (consisting of functions representable as ratio of sections of
some~$\Ll^n$). Now apply Corollary~\ref{res&ext(G-val)} to extend
$v$ to~$R$ and conclude by $v(\sigma/\eta)=
v(\sigma)-v(\eta)=v(g\sigma)-v(\eta)=v(g\sigma/\eta)$.
\end{proof}

A natural geometric characterization of $G$-valuations is given by
\begin{proposition}
Any $G$-valuation is proportional to $v_D$ for a $G$-stable prime
divisor $D$ on a normal $G$-model $X$ of~$K$.
\end{proposition}
\begin{proof}
Let $v\in\Vv$ and choose $f_1,\dots,f_s\in\Oo_v$ whose residues
generate~$\kk(v)$. Take a normal projective $G$-model $X$ of~$K$
and a $G$-line bundle $\Ll$ on~$X$ such that
$f_i=\sigma_i/\sigma_0$ for some $\sigma_0,\dots,\sigma_s\in
\Ho^0(X,\Ll)$. Let $M\subseteq\Ho^0(X,\Ll)$ be the $G$-submodule
generated by $\sigma_0,\dots,\sigma_s$. The respective rational
map $\phi:X\dasharrow\PP(M^{*})$ is $G$-equivariant. Replacing $X$
by the normalized closure of the graph of~$\phi$, we may assume
that $\phi$ is a $G$-morphism. Corollary~\ref{v(g*sec/sec)}
implies $v(M/\sigma_0)\geq0$, whence the center $Y\subseteq
X'=\phi(X)$ of $v|_{\kk(X')}$ intersects an affine chart
$X'_{\sigma_0}$, and $f_1,\dots,f_s\in\Oo_{X',Y}$. Therefore, if
$D$ is the center of $v$ on~$X$, then $f_1,\dots,f_s\in\Oo_{X,D}$,
whence $D$ is a divisor.
\end{proof}
Here is a relative version of this proposition.
\begin{proposition}\label{div->center}
Suppose a $G$-valuation $v$ has the center $Y$ on a $G$-model $X$
of~$K$. Then there exists a normal $G$-model $X'$ and a projective
morphism $\phi:X'\to X$ such that the center of $v$ on $X'$ is a
divisor $D'$ and $\phi(D')=Y$.
\end{proposition}
\begin{proof}
Take any projective $G$-model $X'$ such that the center of $v$ is
a divisor $D'\subset X'$. The rational map $\phi:X'\dasharrow X$
is defined on an open subset intersecting~$D'$, and
$\overline{\phi(D')}=Y$. Now we replace $X'$ by the normalized
closure $\widetilde{X}$ of the graph of $\phi$ in $X'\times X$.
Since $\widetilde{X}$ projects onto $X'$ isomorphically over the
domain of definition of~$\phi$, we can lift $D'$
to~$\widetilde{X}$.
\end{proof}

From now on, $G$~is a connected reductive group.

\begin{lemma}\label{v(G-alg)}
If $A\subset K$ is a rational $G$-algebra, then $\forall v\in\Vv,\
f\in A:\ v(f)=\min_{\widetilde{f}\in
(M^q)^{(B)}}v(\widetilde{f})/q$, where $M$ is a $G$-submodule
generated by~$f$, and $q$ is a sufficiently big power of the
characteristic exponent of~$\kk$.
\end{lemma}
\begin{proof}
As $M$ is generated by~$f$, we have $v(M)\geq v(f)$, whence
$v\left((M^q)^{(B)}\right)\geq qv(f)$. To prove that the equality
is reached, in characteristic zero ($\implies q=1$) it suffices to
note that $M$ is generated by~$M^{(B)}$ $\implies v(f)\in
v(M)\geq\min v\left(M^{(B)}\right)$. In the general case, this is
not true, and one has to consider powers of~$M$. We organize them
in a graded $G$-algebra $R=\bigoplus_{n\geq0}M^n$ and consider a
graded $G$-stable ideal $I=\bigoplus I_n\normin R$, $I_n=\{h\in
M^n\mid v(h)>nv(f)\}$.

As $M\not\subseteq I$, there exists $r\in M$ such that $0\ne
r\bmod I\in(R/I)^U$. By Lemma~\ref{(A/I)^H}, $(R/I)^U$ is a purely
inseparable finite extension of $R^U/I^U$. Hence $\exists h\in
I_q:\ \widetilde{f}=r^q+h\in R^U_q$, and
$v(\widetilde{f})=v(r^q)=qv(f)$.
\end{proof}
\begin{remark}\label{lin.red}
In characteristic zero or for $G=T$ Lemma~\ref{v(G-alg)} yields
$v(f)=\min_{\lambda\in\Ch_{+}}v(f_{(\lambda)})$, where
$f_{(\lambda)}$ is the projection of $f$ to the isotypic component
$A_{(\lambda)}$ of~$A$.
\end{remark}

Recall from \ref{B-charts} that $\Dd$ denotes the set of
non-$G$-stable prime divisors on (any) $G$-model of $K$ and
$K_B\subseteq K$ is the subalgebra of rational functions with
poles in~$\Dd^B$.

The following approximation lemma of Knop \cite[3.5]{G-val} allows
to simplify the study of $G$-valuations by restricting to
$B$-eigenfunctions.
\begin{lemma}\label{approx}
For any $G$-valuation $v\in\Vv$ and any rational function $f\in
K_B$ there exists a rational $B$-eigenfunction $\widetilde{f}\in
K^{(B)}$ such that:
%*
\begin{equation*}
\begin{cases}
v(\widetilde{f})=v(f^q)&\\
w(\widetilde{f})\geq w(f^q),&\forall w\in\Vv\\
v_D(\widetilde{f})\geq v_D(f^q),&\forall D\in\Dd^B
\end{cases}
\end{equation*}
%*
where $q$ is a sufficiently big power of the characteristic
exponent.
\end{lemma}
\begin{proof}
Let $X$ be a normal $G$-model and $\delta=\divr_{\infty}f$, the
divisor of poles on~$X$.
\begin{question}
Perhaps, we have to take the $B$-stable part of~$\divr f$.
\end{question}
Take $\eta\in\Ho^0(X,\Lin{\delta})$ such that $\divr\eta=\delta$
$\implies\sigma=f\eta\in \Ho^0(X,\Lin{\delta})$. Extend all
$G$-valuations to $R=\bigoplus_{n\geq0}\Ho^0(X,\Lin{n\delta})$ and
consider a $G$-submodule $M\subseteq\Ho^0(X,\Lin{\delta})$
generated by~$\sigma$. For any $B$-eigensection
$\widetilde\sigma\in(M^q)^{(B)}$, put
$\widetilde{f}=\widetilde\sigma/\eta^q$. Then
$v_D(\widetilde{f})\geq qv_D(f)$, $w(\widetilde{f})\geq qw(f)$,
and $v(\widetilde{f})=qv(f)$ for some
$\widetilde\sigma\in(M^q)^{(B)}$ by Lemma~\ref{v(G-alg)}.
\end{proof}
\begin{corollary}\label{restr}
$G$-valuations are determined uniquely by their restriction to
$K^{(B)}$.
\end{corollary}
\begin{proof}
As $\Quot{K_B}=K$, two distinct $v,w\in\Vv$ differ on some $f\in
K_B$, say $v(f)<w(f)$. Lemma~\ref{approx} yields
$v(\widetilde{f})<w(\widetilde{f})$ for some $\widetilde{f}\in
K^{(B)}$.
\end{proof}

\section{Valuation cones}
\label{val.cones}

We have seen in \ref{G-vals} that $G$-valuations of $K$ are
determined by their restriction to~$K^{(B)}$. In this section, we
give a geometric qualitative description of $\Vv$ in terms of this
restriction.

Let $\vv$ be a geometric valuation of~$K^B$. Factoring the exact
sequence \eqref{K^(B)} by $\Oo_{\vv}^{\times}$ yields an exact
sequence of lattices
%*
\begin{equation}
\label{K^(B)/v} 0\longrightarrow\ZZ_{\vv}\longrightarrow\RG_{\vv}
\longrightarrow\RG\longrightarrow0
\end{equation}
%*
where $\ZZ_{\vv}$ is the value group of~$\vv$. Passing to the dual
$\QQ$-vector spaces, we obtain
%*
\begin{equation}
\label{halfspace}
\begin{array}{ccccccccc}
0&\longleftarrow&\QQ_{\vv}&\longleftarrow&\ES_{\vv}
&\longleftarrow&\ES&\longleftarrow&0\\
&&|{\bigcup}&&|{\bigcup}&&\parallel&&\\
0&\longleftarrow&\QQ_{\vv,+}&\longleftarrow&\EHS{\vv}
&\longleftarrow&\ES&\longleftarrow&0\\
\end{array}
\end{equation}
%*
where $\QQ_{\vv}=\QQ$ and $\EHS{\vv}$ is the preimage of the
positive ray $\QQ_{\vv,+}$ for $\vv\ne0$, and
$\QQ_{0}=\QQ_{0,+}=0$, $\EHS{0}=\ES_0=\ES$.

\begin{definition}
The \emph{hyperspace} (of~$K$) is the union
$\Hyp=\bigcup_{\vv}\EHS{\vv}$, where $\vv$ runs over all geometric
valuations of $K^B$ considered up to proportionality. More
precisely, $\Hyp=\ES$ in the spherical case, and if $c(K)>0$, then
$\Hyp$ is the union of half-spaces $\EHS{\vv}$ (over all
$\vv\ne0$) glued together along their common boundary
hyperplane~$\ES$, called the \emph{center} of~$\Hyp$.
\begin{question}
Functoriality
\end{question}
\end{definition}

Since $\RG$ is a free Abelian group, the exact sequence
\eqref{K^(B)} splits. Any splitting of \eqref{K^(B)} gives rise to
simultaneous splittings of \eqref{K^(B)/v}, \eqref{halfspace},
$\forall\vv$. From time to time, we will fix such a splitting
$\ef{}:\RG\to K^{(B)}$, $\lambda\mapsto\ef{\lambda}$.

If $v$ is a geometric valuation of $K$ dominating $\vv$, then
$v|_{K^{(B)}}$ factors to a linear functional on $\RG_{\vv}$
non-negative on $\ZZ_{\vv,+}$, i.e., an element of~$\EHS{\vv}$.
Therefore $\Vv\embeds\Hyp$, and there is a restriction map
$\res:\Dd^B\to\Hyp$, which is in general not injective. Put
$\Vv_{\vv}=\Vv\cap\EHS{\vv}$ and
$\Dd^B_{\vv}=\res^{-1}(\EHS{\vv})$. We say that
$(\Hyp,\Vv,\Dd^B,\res)$ is the \emph{colored hyperspace}.

Our aim is to describe $\Vv_{\vv}$.

\begin{example}\label{G=T}
Assume that $G=B=T$ is a torus. Since every $T$-action has trivial
birational type, there exists a $T$-model $X=T/T_0\times C$, where
$T_0=\Ker(T:K)$ and $C=X/T$. We have $\RG=\Ch(T/T_0)$,
$K^T=\kk(C)$, $K_T=K^T[\RG]$. By Remark~\ref{lin.red}, there is
the only way to extend $v\in\Hyp$ to a $T$-valuation of~$K$: put
$v(f)=\min_{\lambda}v(f_{\lambda})$, $\forall f=\sum
f_{\lambda}\in K_T$, $f_{\lambda}\in K^{(T)}_{\lambda}$,
$\lambda\in\RG$.

To prove the multiplicative property, for $\forall f=\sum
f_{\lambda},\ g=\sum g_{\lambda}\in K_T$, choose $\gamma\in\ES$
such that $\min\langle\gamma,\lambda\rangle$ and
$\min\langle\gamma,\mu\rangle$ over all $\lambda$ with
$v(f_{\lambda})=\min$, resp.\ $\mu$ with $v(g_{\mu})=\min$, are
reached at only one point $\lambda_0$, resp.~$\mu_0$. Then
$fg=\sum f_{\lambda}g_{\mu}=
f_{\lambda_0}g_{\mu_0}+\sum_{(\lambda,\mu)\ne(\lambda_0,\mu_0)}
f_{\lambda}g_{\mu}$, and for any term of the \ordinal{2} sum we
have either $v(f_{\lambda}g_{\mu})>v(f_{\lambda_0}g_{\mu_0})$ or
$\langle\gamma,\lambda+\mu\rangle>
\langle\gamma,\lambda_0+\mu_0\rangle$. It follows that
$v(fg)=v(f_{\lambda_0}g_{\mu_0})=v(f_{\lambda_0})+v(g_{\mu_0})=
v(f)+v(g)$. Other properties of a valuation are obvious.

Finally, let $\vv=v|_{K^T}$ and consider a short exact subsequence
of~\eqref{K^(B)}:
%*
\begin{equation}
\label{K^(T)_0} 1\longrightarrow K^T_0\longrightarrow K^{(T)}_0
\longrightarrow\RG_0\longrightarrow0
\end{equation}
%*
where $K^{(T)}_0$ is the kernel of $v:K^{(T)}\to\QQ$, and
$K^T_0=\Oo_{\vv}^{\times}$. Note that any element of $K$ can be
written as $f=f_1/f_2$, $f_i\in K_T$, $v(f_2)=0$. It follows that
$\kk(v)$ is the fraction field of
$K_T\cap\Oo_v/K_T\cap\m_v\iso\kk(\vv)[\RG_0]$ $\implies \trdeg
K-\trdeg\kk(v)=\trdeg K^T+\rk\RG-\trdeg\kk(\vv)-\rk\RG_0=
\rk(K^T)^{\times}/K^T_0+\rk\RG/\RG_0=\rk K^{(T)}/K^{(T)}_0 \leq1$,
hence $v$ is geometric by Proposition~\ref{tr.deg.k(v)}.

We conclude that $\Vv=\Hyp$. By the way, we proved that every
$T$-invariant valuation of $K$ is geometric provided that its
restriction to $K^T$ is geometric.
\end{example}

The main result of this section is
\begin{theorem}\label{val.cone}
For any geometric valuation $\vv$ of~$K^B$, $\Vv_{\vv}$ is a
finitely generated solid convex cone in $\EHS{\vv}$.
\end{theorem}

We prove it in several steps.

\begin{lemma}\label{fin.div}
For any $G$-model $X$, there are only finitely many $B$-stable
prime divisors $D\subset X$ such that $v_D$ maps to~$\EHS{\vv}$.
\end{lemma}
\begin{proof}
Take a sufficiently small $B$-chart $\X\subseteq X$ such that a
geometric quotient $\pi:\X\to\X/B$ exists. Now if $v_D$ maps
to~$\EHS{\vv}$, then either $D$ is an irreducible component of
$X\setminus\X$ or $D=\pi^{-1}(D_0)$, where $D_0$ is the center of
$\vv$ on $\X/B$.
\end{proof}
\begin{corollary}\label{D^B_v}
$\Dd^B_{\vv}$ is finite.
\end{corollary}

In the study of $G$-valuations, it is helpful to consider their
centers on a sufficiently good projective $G$-model.

\begin{lemma}\label{toroid}
Let $P$ be the common stabilizer of all colors with a Levi
decomposition $P=L\semitimes\Ru{P}$, $L\supseteq T$. There exists
a projective $G$-model $X$, a $P$-stable open (not necessarily
affine) subset $\X\subseteq X$, and a $T$-stable closed subvariety
$Z\subseteq\X$ such that
\begin{question}
Relative version: $\forall X'$ $\exists$ projective morphism $X\to
X'$ \dots
\end{question}
\begin{enumerate}
\item\label{Xo/Pu} The action $\Ru{P}:\X$ is proper and has a
geometric quotient. \item\label{Xo=PZ} $\X=PZ$ and the natural
maps $\Ru{P}\times Z\to\X$, $Z\to\X/\Ru{P}$ are finite and
surjective. \item[\ref{Xo=PZ}$'$] In characteristic zero,
%*
\begin{equation*}
\Ru{P}\times Z=P\itimes{L}Z\isoto\X
\end{equation*}
%*
\item\label{L':triv} The $L'$-action on $\X/\Ru{P}$ is trivial.
\item\label{glob.chart} Every $G$-subvariety $Y\subset X$
intersects~$\X$, hence~$Z$. \item\label{tor.col.data} If $Y$ is
the center of $v\in\Vv_{\vv}$, then $\Dd^B_Y=\emptyset$,
$\Vv_Y\subset\EHS{\vv}$.
\end{enumerate}
\end{lemma}
\begin{proof}
Take any projective $G$-model $X$ and choose an ample $G$-line
bundle $\Ll$ and an eigensection $\sigma\in\Ho^0(X,\Ll)^{(B)}$
vanishing on sufficiently many colors such that
$G_{\langle\sigma\rangle}=P$. Put $M=\langle G\sigma\rangle$ and
take a lowest vector $u\in M^{*}$, $\langle\sigma,u\rangle\ne0$,
$G_{\langle u\rangle}=P^{-}$, so that $G\langle
u\rangle\subseteq\PP(M^{*})$ is the unique closed orbit.

There is a natural rational $G$-map $\phi:X\dasharrow\PP(M^{*})$.
Replacing $X$ by the normalized closure of its graph in
$X\times\PP(M^{*})$ makes $\phi$ regular. Put
$\X=X_{\sigma}=\phi^{-1}(\PP(M^{*})_{\sigma})$. Then \ref{Xo/Pu},
\ref{Xo=PZ}, \ref{Xo=PZ}$'$ follow from the local structure of
$\PP(M^{*})_{\sigma}$ (cf.~Lemma~\ref{loc.str.V}). Every $B\cap
L$-stable divisor on $\X/\Ru{P}$ is $L$-stable,
whence~\ref{L':triv}. Every closed $G$-orbit in $X$ maps onto
$G\langle u\rangle$, hence intersects~$\X$, which
yields~\ref{glob.chart}.

To prove~\ref{tor.col.data}, we modify the construction of~$X$.
First, we may choose $\sigma$ vanishing on $\forall
D\in\Dd^B_{\vv}$. Next, consider an affine model $C_0$ of $K^B$
such that $\vv$ has the center $D_0\subseteq C_0$ which is either
a prime divisor or the whole~$C_0$. Let
$\kk[C_0]=\kk[f_1,\dots,f_s]$. We may choose $\Ll$ and $\sigma$ so
that $f_i=\sigma_i/\sigma$, $\sigma_1,\dots,\sigma_s\in
\Ho^0(X,\Ll)^{(B)}$. Let $M'\subseteq\Ho^0(X,\Ll)$ be the
$G$-submodule generated by $\sigma,\sigma_1,\dots,\sigma_s$.

As above, we may assume that the natural rational map
$\phi':X\dasharrow\PP(M'^{*})$ is regular. Put $X'=\phi'(X)$.
Consider the composed map $\pi:\X\to\X'=X'_{\sigma}\to C_0$. By
Corollary~\ref{v(g*sec/sec)}, $v(M'/\sigma)\geq0$, whence the
center $Y'=\phi'(Y)\subseteq X'$ of $v|_{\kk(X')}$ intersects the
$B$-chart~$\X'$. Hence $\Y=Y\cap\X$ is non-empty and
$\overline{\pi(\Y)}\supseteq D_0$. It follows that
$\Vv_Y\sqcup\Dd^B_Y$ maps to~$\EHS{\vv}$.  But any
$D\in\Dd^B_{\vv}$ is contained in $\Zeros{\sigma}=X\setminus\X$,
thence $D\not\supseteq Y$, and we are done.
\end{proof}

\begin{proposition}[\cite{G-val}]\label{geom|K^B}
A $G$-invariant valuation of $K$ is geometric iff its restriction
to $K^B$ is geometric.
\end{proposition}
\newnum
\begin{proof}[Proof~\num~\textup{\cite[3.9, 4.4]{G-val}}\label{v->v|K^U}]
Let $v$ be a nonzero valuation of $K$ such that $v|_{K^B}$ is
geometric. Take a projective $G$-model $X$ as in
Lemma~\ref{toroid}. Then $v$ has the center $Y\subset X$ and
$\Y=Y\cap\X\ne\emptyset$. By Lemma~\ref{toroid}\ref{L':triv},
$\kk(\X/\Ru{P})=K^U$ and $\Y/\Ru{P}$ is the center of~$v|_{K^U}$.

Since $K^U$ is a $T$-field and $(K^U)^T=K^B$, it follows from
Example~\ref{G=T} that $v|_{K^U}$ is geometric. Now by
Proposition~\ref{div->center} there exists a projective birational
$L$-morphism $Z'\to\X/\Ru{P}$ such that the center of $v|_{K^U}$
on $Z'$ is a divisor $D'\subset Z'$. Consider a Cartesian square
%*
\begin{equation*}
\begin{CD}
 \X' @>>>   \X  \\
@VVV       @VVV \\
  Z' @>>> \X/\Ru{P}
\end{CD}
\end{equation*}
%*
where horizontal arrows are birational projective $P$-morphisms,
and vertical arrows are $\Ru{P}$-quotient maps. Therefore $v$ has
a center $D\subset\X'$, which is $P$-stable and maps onto~$D'$,
whence $D$ is the pull-back of~$D'$, i.e., a divisor. This means
that $v$ is geometric.
\end{proof}
\begin{proof}[Proof~\num\label{A(v)}]
Here we use the embedding theory of Chapter~\ref{LV-theory}.
\begin{question}
Choose one proof.
\end{question}
Assume $\vv=v|_{K^B}$. It is easy to construct an affine model
$C_0$ of $K^B$ containing a principal prime divisor $D_0=\divr(t)$
such that either $D_0$ is the center of~$\vv$, or $\vv=0$, $t=1$,
$D_0=\emptyset$, and $\oo{C}=C_0\setminus D_0=\X/B$ for a
(sufficiently small) $B$-chart~$\X$.

If $\Rr$ is the set of all $B$-stable prime divisors in~$\X$
(=preimages of prime divisors in~$\oo{C}$), then
$\{v\}\sqcup\Rr\in\Cd$  defines colored data, and we may consider
the respective Krull algebra $\Alg=\Alg(v,\Rr)$
(cf.~\ref{B-charts}). Recall that we need not to assume apriori
that all $G$-invariant valuations are geometric
(Remark~\ref{non-geom}). Clearly, $\Alg^U=\kk\bigl[\,f\in K^{(B)}
\bigm|\langle v,f\rangle,\langle\Rr,f\rangle\geq0\,\bigr]
\subseteq\kk[\oo{C}]\otimes\kk[\RG]$ is a subalgebra determined by
$v(f)\geq0$, whence $\Alg^U=\kk[C_0]\bigl[\,t^d\ef{\lambda}\bigm|
(d,\lambda)\in\RG_{\vv},\ \langle v,(d,\lambda)\rangle
\geq0\,\bigr]$.

The generating set of $\Alg^U$ over $\kk[C_0]$ forms a finitely
generated semigroup in $\RG_{\vv}$ consisting of lattice points in
the half-space $\{v\geq0\}$, whence condition~\reftag{F} holds
for~$\Alg$.

To prove \reftag{C}, we take $f=t^d\ef{\lambda}$ such that
$\langle v,(d,\lambda)\rangle>0$ and multiply $f$ by
$f_0\in\kk[C_0]$ vanishing on sufficiently many divisors in~$\Rr$.

Conversely, taking $f=t^d\ef{\lambda}$ such that $\langle
v,(d,\lambda)\rangle<0$ proves \reftag{W} for~$v$.

By Corollary~\ref{CFW}, $X_0=\Spec\Alg$ is a $B$-chart and $v$ is
the valuation of a $G$-stable divisor intersecting~$X_0$.
\end{proof}

\begin{remark}\label{aff.cone}
It is often helpful to assume that $K=\Quot{R}$, where $R$ is a
rational $G$-algebra. For instance, $R=\kk[X]$, where $X$ is a
(quasi)affine $G$-model of~$K$ (if any exists). The general case
is reduced to this special one by considering a projectively
normal $G$-model $X$ and taking the affine cone $\widehat{X}$
over~$X$. Then $\widehat{K}=\kk(\widehat{X})$ is a
$\widehat{G}$-field, where $\widehat{G}=G\times\kk^{\times}$ with
$\kk^{\times}$ acting by homotheties. Let us denote various
objects related to $\widehat{K}$ in the same way as for~$K$, but
equipped with a hat. We have short exact sequences
%*
\begin{equation*}
\begin{array}{ccccccccc}
1&\longrightarrow&(K^B)^{\times}&\longrightarrow&K^{(B)}
&\longrightarrow&\RG&\longrightarrow&0\\
&&\parallel&&|{\bigcap}&&\curvearrowleft&&\\
1&\longrightarrow&(\widehat{K}^{\widehat{B}})^{\times}&
\longrightarrow&\widehat{K}^{(\widehat{B})}&
\longrightarrow&\widehat{\RG}&\longrightarrow&0\\
\end{array}
\end{equation*}
%*
and dual sequences
%*
\begin{equation*}
\begin{array}{ccccccccc}
0&\longleftarrow&\QQ_{\vv,+}&\longleftarrow&\EHS{\vv}
&\longleftarrow&\ES&\longleftarrow&0\\
&&\parallel&&
{\uparrow}\makebox[0ex][r]{\raisebox{0.5ex}{$\uparrow$}}&&
{\uparrow}\makebox[0ex][r]{\raisebox{0.5ex}{$\uparrow$}}&&\\
0&\longleftarrow&\QQ_{\vv,+}&\longleftarrow&
\expandafter\widehat\EHS{\vv}&\longleftarrow&\widehat\ES&
\longleftarrow&0\\
\end{array}
\end{equation*}
%*

The set of colors $\widehat{\Dd}^{\widehat{B}}$ is identified
with~$\Dd^B$. The grading of $\kk[\widehat{X}]$ determines two
$G$-valuations $\pm v_0\in\widehat{\ES}$, which generate
$\Ker(\widehat{\ES}\to\ES)$. By Corollary~\ref{res&ext(G-val)},
$\widehat\Vv_{\vv}$ surjects onto $\Vv_{\vv}$ and even is the
preimage of $\Vv_{\vv}$ by Proposition~\ref{convex} below.
\end{remark}

\begin{definition}
Let $f_1,\dots,f_s\in R^{(B)}$ and $f\ne f_1\dots f_s$ be any
highest vector in $\langle Gf_1\rangle\dots\langle Gf_s\rangle$.
Then $f/f_1\dots f_s$ is called a \emph{tail vector} of $R$ and
its weight is called a \emph{tail weight} or just a \emph{tail}.
\begin{question}
Tails form a finitely generated semigroup (Alexeev--Brion, 2003).
\end{question}
Note that tails are negative linear combinations of simple roots.
In characteristic zero tails are the nonzero differences
$\mu-\lambda_1-\dots-\lambda_s$ over all highest weights $\mu$
occurring in the isotypic decomposition of $R_{(\lambda_1)}\dots
R_{(\lambda_s)}$, $\lambda_1,\dots,\lambda_s\in\RG_{+}$.
\end{definition}

Now we proceed in proving Theorem~\ref{val.cone}.

\begin{proposition}\label{convex}
$\Vv_{\vv}$~is a convex cone in~$\EHS{\vv}$.
\end{proposition}
\newnum
\begin{proof}[Proof~\num\label{tails}]
We may assume that $K=\Quot{R}$, where $R$ is a rational
$G$-algebra. The general case is reduced to this one by
considering the affine cone over a projectively normal $G$-model
as above. Then we prove for $\forall v\in\Hyp$ that $v\in\Vv$ iff
%*
\begin{equation*}
\text{$v$~is non-negative on all tail vectors of~$R$} \tag{T}
\end{equation*}
%*
Clearly, this condition is necessary, since $v(\langle
Gf_1\rangle\dots\langle Gf_s\rangle)\geq v(f_1\dots f_s)$,
$\forall v\in\Vv$, $f_1\dots f_s\in R$.

Conversely, assume that \reftag{T} is satisfied.
\begin{question}
Modify the \ordinal{\ref{tails}} proof (at least for $\ch\kk=0$)
using that $\gr R$ has no zero divisors.
\end{question}
It can be generalized as follows: let $f$ be any highest vector in
$\sum_i\langle Gf_{i1}\rangle\dots\langle Gf_{is}\rangle$,
$f_{ij}\in R^{(B)}$; then
$v(f)\geq\min_i\{v(f_{i1})+\dots+v(f_{is})\}$. Indeed, if $f=\sum
f_i$, where $f_i\in\langle Gf_{i1}\rangle\dots\langle
Gf_{is}\rangle$ are highest vectors of the same weight, then
$v(f)\geq\min v(f_i)\geq\min\{v(f_{i1})+\dots+v(f_{is})\}$
by~\reftag{T}. The general case is reduced to this one, because a
certain power $f^q$ belongs to the image of
$\Sym^q(\bigoplus\langle Gf_{i1}\rangle\dots\langle
Gf_{is}\rangle)^U$ by Corollary~\ref{(M/N)^H}.

Consider a rational $B$-algebra $A=\kk[\,gf/h\mid g\in G,\ f,h\in
R,\ {v(f/h)\geq0}\,]$. The ideal $I=(\,gf/h\mid v(f/h)>0\,)\normin
A$ is proper. Indeed, each $f\in I^{(B)}$ is a linear combination
of $(g_{i1}f_{i1}/h_1)\dots(g_{is}f_{is}/h_s)$, where
$v(f_{ij}/h_j)\geq0$ and $>$ occurs for~$\forall i$. By the above
$v(fh_1\dots h_s)\geq\min\{v(f_{i1})+\dots+v(f_{is})\}>
v(h_1)+\dots+v(h_s)\implies v(f)>0$, whence $v>0$ on $I^{(B)}$.

Take any valuation $v'$ non-negative on $A$ and positive on~$I$,
extend it to~$K$, and take the approximating
$G$-valuation~$\overline{v}$ (Proposition~\ref{val->G-val}). For
$\forall f\in R^{(B)},\ g\in G$ we have $v'(gf)\geq
v'(f)\implies\overline{v}(f)=v'(f)$. Now $\forall f,h\in R^{(B)}:\
v(f/h)\geq0\ (>0)\implies f/h\in A\ (\in
I)\implies\overline{v}(f/h)\geq0\ (>0)$. It follows that DVR's of
$v$ and $\overline{v}$ on $K^B$ coincide, hence
$\overline{v}\in\Vv_{\vv}$ provided $v\in\ES_{\vv}$, and
$\overline{v},v$ determine proportional linear functionals
on~$\RG_{\vv}$. Thus $v\in\Vv$.
\end{proof}
\begin{proof}[Proof~\num\label{A(v,w)}] It relies on the
embedding theory. Let $v_1,v_2\in\Vv_{\vv}$ be two
non-proportional vectors. It suffices to prove that
$c_1v_1+c_2v_2\in\Vv_{\vv}$, $\forall c_1,c_2\in\ZZ_{+}$.

Take an affine model $C_0$ of $K^B$ as in the \ordinal{\ref{A(v)}}
proof of Proposition~\ref{geom|K^B} and consider the algebra
$\Alg=\Alg(v_1,v_2,\Rr)$. Then
$\Alg^U\subseteq\kk[\oo{C}]\otimes\kk[\RG]$ is distinguished by
inequalities $v_i(f)\geq0$, whence
$\Alg^U=\kk[C_0]\bigl[\,t^d\ef{\lambda}\bigm|
(d,\lambda)\in\RG_{\vv},\ \langle v_i,(d,\lambda)\rangle \geq0,\
i=1,2\,\bigr]$.

The generating set of $\Alg^U$ over $\kk[C_0]$ forms a finitely
generated semigroup of lattice points in the dihedral cone
$\{v_1,v_2\geq0\}\subseteq\RG_{\vv}\otimes\QQ$, whence
condition~\reftag{F} holds for~$\Alg$. Conditions \reftag{C},
\reftag{W} are verified in the same way as in
Proposition~\ref{geom|K^B} by taking $f=t^d\ef{\lambda}$ with
$\langle v_i,(d,\lambda)\rangle,\langle v_j,(d,\lambda)\rangle>0$
or $\langle v_i,(d,\lambda)\rangle<0 \leq\langle
v_j,(d,\lambda)\rangle$, respectively, $\{i,j\}=\{1,2\}$. Thus by
Corollary~\ref{CFW}, $X_0=\Spec\Alg$ is a $B$-chart and $v_i$
correspond to $G$-stable divisors $D_i\subset X=GX_0$
intersecting~$X_0$.

We blow up the ideal sheaf $\Lin{-nc_2D_1}+\Lin{-nc_1D_2}$,
$n\gg0$, and prove that the exceptional divisor corresponds to
$c_1v_1+c_2v_2$.

The local structure of $X_0$ provided by
Proposition~\ref{loc.str.B} allows to replace $X_0$ by
$X_0/\Ru{P}$ and assume that $G=T$ and $X=X_0$ is an affine
$T$-model of~$K$. We may choose $f_i=t^{d_i}\ef{\lambda_i}$ such
that $v_i(f_j)=n\delta_{ij}$, $n\gg0$. The above ideal sheaf is
represented by a proper ideal
$I=(f_1^{c_2},f_2^{c_1})\normin\Alg$. (Indeed, it is easy to see
that $v_1+v_2>0$ on~$I^{(T)}$.)

The blow-up of $I$ is given in $X\times\PP^1$ by the equation
$[f_1^{c_2}:f_2^{c_1}]=[t_1:t_2]$, where $t_i$ are homogeneous
coordinates on~$\PP^1$. The exceptional divisor is given in the
open subset $\{t_2\ne0\}$ by the equation $f_2=0$.  Let $v_0$ be
the respective valuation. Up to a power, any $f\in\Alg^{(T)}$ is
represented as $f=f_0f_1^{k_1}f_2^{k_2}$, $v_i(f_0)=0$, $k_i\divby
c_j$, $\{i,j\}=\{1,2\}$. Then
$v_0(f)=v_0(f_0(t_1/t_2)^{k_1/c_2}f_2^{(c_1k_1+c_2k_2)/c_2})\sim
c_1k_1+c_2k_2$. It follows that $v_0\sim c_1v_1+c_2v_2$.
\end{proof}

\begin{proposition}\label{fin.gen}
The cone $\Vv_{\vv}$ is finitely generated.
\end{proposition}
\begin{proof}
Take a projective $G$-model $X$ as in Lemma~\ref{toroid}. Then any
$v\in\Vv_{\vv}$ has the center $Y$ on $X$, and by
Lemma~\ref{toroid}\ref{tor.col.data}, $\Dd^B_Y=\emptyset$,
$\Vv_Y\subset\Vv_{\vv}$. Condition \reftag{S} yields
$\forall\lambda\in\RG_{\vv}:\ \langle\Vv_Y,\lambda\rangle\geq0
\implies\langle v,\lambda\rangle\geq0$. Hence $v\in\QQ_{+}\Vv_Y$.
It remains to note that $\bigcup_Y\Vv_Y$ is finite by
Lemma~\ref{fin.div}.
\end{proof}

\begin{proof}[Proof of Theorem~\ref{val.cone}]
Due to Propositions \ref{convex}--\ref{fin.gen} and
Theorem~\ref{cent.cone}, it remains to prove that any geometric
valuation $\vv\ne0$ of $K^B$ extends to a $G$-valuation of~$K$.
For this, we modify the \ordinal{\ref{tails}} proof of
Proposition~\ref{convex}.
\begin{question}
Or use $C_0\embeds(\X')^{B^{-}}$ in the notation of
Lemma~\ref{toroid} (Knop's thesis, 4.1).
\end{question}

Namely, in the definition of $A$ we replace $v$ by $\vv$ and
assume $f/h\in K^B$. The respective ideal $I\normin A$ is still
proper. For otherwise, $1$~is a linear combination of
$(g_{i1}f_{i1}/h_1)\dots(g_{is}f_{is}/h_s)$, where
$\vv(f_{ij}/h_j)\geq0$ and $>$ occurs for~$\forall i$. But the
$T$-weights of all $T$-eigenvectors in $\langle Gf_{ij}\rangle$
except $f_{ij}$ are obtained from the weight of~$f_{ij}$ (=the
weight of~$h_j$) by subtracting simple roots. Hence we may assume
$g_{ij}=\1$, and $1$ is a linear combination of
$r_i=(f_{i1}/h_1)\dots(f_{is}/h_s)$, $\vv(r_i)>0$, a
contradiction.

Now reproducing the arguments of that proof yields a $G$-valuation
$\overline{v}$ such that $\overline{v}|_{K^B}=\vv$.
\end{proof}

Parabolic induction, which is helpful in various reduction
arguments, keeps the colored hyperspace ``almost'' unchanged.
Suppose $K$ is obtained from a $G_0$-field $K_0$ by parabolic
induction $G\supseteq Q\onto G_0$, i.e., a $G$-model $X$ of $K$ is
obtained from a $G_0$-model $X_0$ of $K_0$ by this procedure.
There is a natural projection $\pi:X=G\itimes{Q}X_0\to G/Q$. We
may assume $Q\supseteq B^{-}$. Then the colors of $G/Q$ are the
Schubert divisors $D_{\alpha}=\overline{B[r_{\alpha}]}$
($\alpha\in\Pi$, $r_{\alpha}\notin Q$). Let us denote various
objects related to $K_0$ in the same way as for~$K$, but with a
subscript~$0$.
\begin{proposition}\label{par.ind(hyp)}
There are natural identifications $\Hyp=\Hyp_0$, $\Vv=\Vv_0$, and
$\Dd^B=\Dd_0^{B_0}\sqcup\pi^{-1}(\Dd^B(G/Q))$ such that
$\res=\res_0$ on $\Dd_0^{B_0}$ and
$\res(\pi^{-1}(D_{\alpha}))\in\ES$ is the restriction of
$\alpha^{\vee}$ to~$\RG$.
\end{proposition}
\begin{proof}
Since $\pi^{-1}(B[e])\iso\Ru{Q^{-}}\times X_0$, the restriction to
the fiber identifies $K^{(B)}$ with~$K_0^{(B_0)}$. Therefore, the
hyperspaces are identified and the colors are extended by the
pullbacks of the Schubert divisors.

To prove $\Vv=\Vv_0$, we construct $X_0$ as in Lemma~\ref{toroid}.
Then $X$ satisfies the conditions of Lemma~\ref{toroid}, too. By
the condition~\ref{tor.col.data},
$\Vv\cap\EHS{\vv}=\Vv_0\cap\EHS{\vv}$.

Finally, $\pi^{-1}(D_{\alpha})$~is transversal to the subvariety
$X_{\alpha}=L_{\alpha}\itimes{B_{\alpha}^{-}}X_0\subseteq X$ along
$r_{\alpha}*X_0$, where $L_{\alpha}$ is the Levi subgroup of
$P_{\alpha}$ and $B_{\alpha}^{-}=B^{-}\cap
L_{\alpha}=TU_{-\alpha}$ acts on $X_0$ via~$T$. The variety
$X_{\alpha}$ is horospherical w.r.t.~$L_{\alpha}$ and intersects
generic $B$-orbits of~$X$. We easily deduce that generic $B$-orbit
closures in $X$ intersect $\pi^{-1}(D_{\alpha})$ and the
intersections cover a dense subset $Br_{\alpha}*X_0$. Hence all
$f\in K^B$ have order~$0$ along~$\pi^{-1}(D_{\alpha})$. To
determine $\res(\pi^{-1}(D_{\alpha}))$, it suffices to consider
the restriction of $K^{(B)}$ to a generic $L_{\alpha}$-orbit
in~$X_{\alpha}$, cf.~\ref{S-var}.
\end{proof}

\section{Central valuations}
\label{central}

$G$-valuations of $K$ that vanish on $(K^B)^{\times}$ are called
\emph{central}. They play a distinguished role among all
$G$-valuations. By \eqref{K^(B)} a central valuation restricted to
$K^{(B)}$ factors through a linear functional on~$\RG$. Thus
central valuations form a subset $\Zz\subset\ES=\Hom(\RG,\QQ)$.

If $G$ is linearly reductive (i.e., $\ch\kk=0$ or $G=T$) and
$K=\Quot{R}$ for a rational $G$-algebra $R$ with the isotypic
decomposition $R=\bigoplus_{\lambda\in\RG}R_{(\lambda)}$, then any
$v\in\Zz$ is constant on each isotypic component $R_{(\lambda)}$
(otherwise there would exist two highest vectors $f_1,f_2\in
R^{(B)}$ of the same highest weight $\lambda$ with $v(f_1)\ne
v(f_2)$) and $v|_{R_{(\lambda)}\setminus\{0\}}= \langle
v,\lambda\rangle$, $\forall\lambda\in\RG$.

\begin{theorem}\label{cent.cone}
$\Zz=\Vv\cap\ES$ is a solid convex polyhedral cone in $\ES$
containing the image of the antidominant Weyl chamber.
\end{theorem}
\begin{proof}
By Remark~\ref{aff.cone}, we reduce the question to the case
$K=\Quot{R}$, where $R$ is a rational $G$-algebra. Condition
\reftag{T} defining the subset $\Vv\subset\Hyp$ transforms under
restriction to $\ES$ to the following one: $v\in\Zz$ iff
%*
\begin{equation*}
\text{$v$~is non-negative on all tails of~$R$} \tag{T$_0$}
\end{equation*}
%*
Since tails are negative linear combinations of simple roots, we
see that \reftag{T$_0$} determines a convex cone containing the
image of the antidominant Weyl chamber, whence a solid cone. We
conclude by Proposition~\ref{fin.gen}.
\end{proof}
\begin{example}\label{cent.val(G)}
Let $G$ act on itself by right translations and $K=\kk(G)$. Here
$\RG=\Ch(T)$. Recall from Proposition~\ref{U(Mat)} that $\kk[G]$
is the union of subspaces $\Mat(V)$ of matrix entries over all
$G$-modules~$V$.

In characteristic zero, $\Mat(V(\lambda))=\kk[G]_{(\lambda)}$ are
the isotypic components of $\kk[G]$ and $\Mat(V(\lambda))\cdot
\Mat(V(\mu))=\bigoplus\Mat(V(\nu))$ over all simple submodules
$V(\nu)$ occurring in $V(\lambda)\otimes V(\mu)$ \eqref{mat}.
Generally, each highest vector $v\in V^{(B)}$ gives rise to
highest vectors $f_{\omega,v}\in\Mat(V)^{(B)}$, $\omega\in V^{*}$,
of the same highest weight.

If $v_{\lambda}\in V$ is a highest vector of regular highest
weight~$\lambda$, then $V$ contains $T$-eigenvectors
$v_{\lambda-\alpha}$ of weights $\lambda-\alpha$ for all simple
roots~$\alpha$, and $v_{\lambda}\otimes v_{\lambda-\alpha}-
v_{\lambda-\alpha}\otimes v_{\lambda}$ are highest vectors of
highest weights~$2\lambda-\alpha$ in $V\otimes V$. It follows that
all $-\alpha$ occur among tails of $\kk[G]$, whence $\Zz$ is the
antidominant Weyl chamber.
\end{example}

In characteristic zero, a much more precise information on the
structure of $\Zz$ can be obtained, see~\ref{Weyl}.

A special case of Theorem~\ref{cr(Y<X)} distinguishes central and
non-central $G$-valuations in terms of the complexity and the rank
of respective $G$-stable divisors.
\begin{proposition}\label{c&r(cent)}
A $G$-valuation $v\ne0$ is (non)central iff $c(\kk(v))=c(K)$,
$r(\kk(v))=r(K)-1$ ($c(\kk(v))=c(K)-1$, $r(\kk(v))=r(K)$,
respectively).
\end{proposition}
\begin{proof}
Let $Y\subset X$ be a $G$-stable divisor on a $G$-model of $K$
corresponding to~$v$. By Lemma~\ref{B-eigenfun},
$\kk(v)^U=\kk(Y)^U$ is a purely inseparable extension of the
residue field $\kk(v|_{K^U})$ of $\Oo_{X,Y}^U$, and similarly
for~$\kk(v)^B$. Thus by Proposition~\ref{c+r},
$c(\kk(v))+r(\kk(v))=\trdeg\kk(v)^U=\trdeg K^U-1=c(K)+r(K)-1$, and
$c(\kk(v))=c(K)$ iff $\kk(v|_{K^B})=K^B$ iff $v\in\Zz$.
\end{proof}

Now we explain the geometric meaning of the linear part of the
central valuation cone.

\begin{definition}
A $G$-automorphism of $K$ acting trivially on $K^B$ is called
\emph{central}. Denote by $\CA K=\Ker(\Aut_GK:K^B)$ the group of
central automorphisms.
\begin{question}
Require that central automorphisms act on $K^{(B)}$ by
homotheties. Is then $\CA K$ an algebraic group?
\end{question}
\end{definition}

\begin{theorem}\label{cent.aut}
There exists the largest connected algebraic subgroup
$S\subseteq\CA K$. It has the following properties:
\begin{enumerate}
\item\label{S:triv} $S$ acts on $\Vv,\Dd^B,{}_G\XX^{\norm}$
trivially. \item\label{S:norm} $S$ acts on every normal $G$-model
of $K$ regularly. \item\label{S<A} There is a canonical embedding
$S\embeds A=\Hom(\RG,\kk^{\times})$ via the action $S:K^{(B)}$, so
that $\Hom(\Ch(S),\QQ)\subseteq\Zz\cap-\Zz$. \item\label{dim S}
There exists a $G$-subfield $K'\subseteq K$ with $(K')^U=K^U$ and
with the same colored hyperspace as for $K$ such that $K$ is
purely inseparable over $K'$ and $\Hom(\Ch(S'),\QQ)=\Zz\cap-\Zz$
for the largest connected algebraic subgroup $S'\subseteq\CA K'$.
\end{enumerate}
\end{theorem}
\begin{proof}
Let $S\subseteq\CA K$ be any connected algebraic subgroup. Suppose
first that $K=\Quot{R}$, where $R$ is a rational $(G\times
S)$-subalgebra. Without loss of generality, we may assume in the
reasoning below that $R=\kk[X]$, where $X$ is a normal
(quasi)affine $G$-model of $K$ acted on by~$S$.

If $f\in R^{(B)}$, then $sf\in R^{(B)}$ has the same weight,
$\forall s\in S$, whence $sf/f\in K^B\subseteq K^S$. Hence $s\divr
f-\divr f$ is $S$-stable. (The divisors are considered on the
normalized projective closure of~$X$.) But on the other hand, this
divisor has no $S$-stable components, hence is zero, and
$sf\in\kk^{\times}f$. Therefore $R^{(B)}\subseteq R^{(S)}\implies
K^{(B)}\subseteq K^{(S)}$.

This yields a homomorphism $S\to A$. Let $S_0$ be its connected
kernel. Then $S_0$ acts on $R^U$ trivially. As $S_0$ commutes with
$G$, it acts on $R$ trivially. (In positive characteristic, this
stems from Lemma~\ref{A>GA^U}.) Hence $S_0=\{\1\}$ and $S$ is a
torus. Then every $s\in S$, $s\ne\1$, acts non-trivially
on~$K^{(B)}$: just take a ($G$-stable) eigenspace of $s$ in $R$ of
eigenvalue $\ne1$ and choose a $B$-eigenvector there. Thus
$S\embeds A$.

Since $S$-action multiplies $B$-eigenfunctions by scalars and
$G$-valuations are determined by their restriction to~$K^{(B)}$,
the action $S:\Vv$ is trivial. As any $D\in\Dd^B$ is a component
of $\divr f$, $f\in K^{(B)}\subseteq K^{(S)}$, $S$~fixes all
colors. Then $S$ fixes all $G$-germs by Proposition~\ref{germ},
whence~\ref{S:triv}. Assertion \ref{S:norm} stems
from~\ref{S:triv}.

Each one-parameter subgroup $\gamma\in\Hom(\Ch(S),\ZZ)$ defines a
$G$-stable grading of~$R$, which gives rise to a central valuation
(the order of the lowest homogeneous term) represented
by~$\gamma$: $f_{\lambda}\mapsto\langle\gamma,\lambda\rangle$,
$\forall f_{\lambda}\in R^{(B)}_{\lambda}$. This finally
yields~\ref{S<A}.

Furthermore, any $\tau\in\CA K$ commutes with~$\gamma$. Indeed, we
have two gradings of $R$ defined by $\gamma,\tau\gamma\tau^{-1}$.
They coincide on $R^{(B)}$, whence on the $G$-subalgebra of $R$
generated by~$R^U$, whence on~$R$. (The last implication is easily
deduced from Lemma~\ref{A>GA^U}.) Therefore any two subtori
$S_1,S_2\subseteq\CA K$ commute, and the natural homomorphism
$S_1\times S_2\to\CA K$ provides a larger subtorus. Since the
dimensions of subtori are restricted from above by
$\dim(\Zz\cap-\Zz)$, there exists the largest one.

We prove \ref{dim S} in characteristic zero referring to
\cite[8.2]{G-val} for $\ch\kk>0$. Every lattice vector
$\gamma\in\Zz\cap-\Zz$ defines a $G$-stable grading of $R$ such
that $R_{(\lambda)}$ are homogeneous of degree
$\langle\gamma,\lambda\rangle$. Since $\gamma$ vanishes on tails,
this grading respects multiplication and defines a 1-subtorus of
central automorphisms. These 1-subtori generate a subtorus
$S\subseteq A$ which is the connected common kernel of all tails,
and $\Hom(\Ch(S),\QQ)=\Zz\cap-\Zz$.

Finally, the general case is reduced to the above by taking a
projectively normal $G$-model $X$ acted on by $S$ and considering
the affine cone $\widehat{X}$ over~$X$. (Note that any central
subtorus action on $X$ lifts to $\widehat{X}$ if we consider a
sufficiently ample projective embedding of~$X$.) By
Remark~\ref{aff.cone}, $\Zz=\widehat{\Zz}/\QQ v_0$, and
$S=\widehat{S}/\kk^{\times}$, where the central valuation $v_0$ of
$\widehat{K}$ is defined by a 1-subtorus
$\kk^{\times}\subseteq\widehat{S}$ acting on $\widehat{X}$ by
homotheties. So all assertions on $\widehat{S}$ transfer to~$S$.
\end{proof}

In the generically transitive case, we can say more.
\begin{proposition}\label{cent.aut(G/H)}
If $K=\kk(G/H)$, then $\CA K$ is a quasitorus extended by a finite
$p$-group and $\dim\CA K=\dim(\Zz\cap-\Zz)$.
\end{proposition}
\begin{proof}
Since $\Aut_GK=\Aut_GG/H=N(H)/H$ is an algebraic group, $\CA K$~is
an algebraic group as well. Central automorphisms preserve generic
$B$-orbits in~$G/H$, whence there exist finitely many general
points $x_1,\dots,x_s\in G/H$ such that $\CA K\embeds
\Aut_BBx_1\times\dots\times\Aut_BBx_s$. The latter group is a
subquotient of $B\times\dots\times B$, which explains the
structure of~$\CA K$ in view of Theorem~\ref{cent.aut}\ref{S<A}.
By Theorem~\ref{cent.aut}\ref{dim S}, there exists a purely
inseparable $G$-map $G/H\onto\HS$ such that
$\dim\CA\HS=\dim(\Zz\cap-\Zz)$. But every $G$-automorphism of
$\HS$ lifts to $G/H$ by the universal property of quotients. Hence
$\dim\CA K$ is the same.
\end{proof}
\begin{remark}
If $\ch\kk=0$ and $G/H$ is quasiaffine, then $\CA G/H$ is
canonically embedded in $A$ as the common kernel of all tails of
$\kk[G/H]$ and acts on each $\kk[G/H]_{(\lambda)}$ by the
character $\lambda\in\RG_{+}(G/H)$.
\end{remark}
\begin{example}
In Example~\ref{cent.val(G)}, $\Aut_GG=G$ (acting by left
translations) and $\CA G=\Ker(G:G/B)=Z(G)$.
\end{example}
\begin{example}
If the orbit map $G\to\HS$ is not separable, then $\dim\CA\HS$ may
be smaller than ``the proper value''. For instance, let
$\ch\kk=2$, $G=\SL_2$, $X=\PP(\sgl_2)$.
\begin{question}
An example for arbitrary $\ch\kk$.
\end{question}
Then $X$ has an open orbit $\HS$ with stabilizer~$U$. By
Theorem~\ref{cent.aut}\ref{S:norm}, the central torus $S$ embeds
in $\Aut_GX$, but the latter group is trivial. Indeed, each
$G$-automorphism of $X$ lifts to an intertwining operator of
$\Ad\SL_2$, but it is easy to see that such an operator has to be
scalar. However $\kk(G/U)^B=\kk\implies\CA G/U=\Aut_GG/U=T
\implies\dim(\Zz\cap-\Zz)=1$, cf.\ Theorem~\ref{cent(hor)}.
\end{example}

We have seen in \ref{cotangent} that horospherical $G$-varieties
play an important role in studying general $G$-varieties.
Apparently they can be characterized in terms of central valuation
cones.
\begin{theorem}\label{cent(hor)}
A $G$-variety $X$ is horospherical iff $\Zz(X)=\ES(X)$.
\end{theorem}
\newnum
\begin{proof}[Proof~\num\/\label{no_tails}
\textup{($\ch\kk=0$) \cite[\S4]{contr}, \cite[\S5]{c&mod}}] This
proof goes back to Popov, who however considered tails of
coordinate algebras instead of central valuations,
cf.~Proposition~\ref{iso.mult}. The generically transitive case is
due to Pauer~\cite{hor.sub}.

By Remark~\ref{aff.cone}, we may assume $X$ to be (quasi)affine.
By~\reftag{T$_0$}, $\Zz=\ES$ iff $R=\kk[X]$ has no tails. However
one proves that a rational $G$-algebra $R$ has no tails iff
$R\iso(\kk[G/U^{-}]\otimes R^U)^T=
\bigoplus_{\lambda\in\Ch_{+}}\kk[G/U^{-}]_{(\lambda)}\otimes
R^{(B)}_{\lambda}$. Here $T=\Aut_GG/U^{-}$ acts on $G/U^{-}$ by
right translations, so that isotypic components
$\kk[G/U^{-}]_{(\lambda)}$ are at the same time $T$-eigenspaces of
weight~$-\lambda$ (cf.\ \eqref{k[G/U]}), and the isomorphism is
given by $g\ef{\lambda}\otimes f_{\lambda}\mapsto gf_{\lambda}$,
$\forall g\in G,\ f_{\lambda}\in R^{(B)}_{\lambda}$, where
$\ef{\lambda}\in\kk[G/U^{-}]^{(B)}_{\lambda}$,
$\ef{\lambda}(\1)=1$. In our situation, this implies that $X=(G\by
U^{-}\times X\by U)\by T$ is horospherical.

Conversely, if $R$ has tails, then tails do not vanish under
restriction of the isotypic decomposition of $R_{(\lambda)}\cdot
R_{(\mu)}$ to generic $G$-orbits. But the coordinate algebra of a
horospherical homogeneous space has no tails since isotypic
components of $\kk[G/U^{-}]$ are $T$-eigenspaces. Hence $X$ is not
horospherical.
\end{proof}
\begin{proof}[Proof~\num~\textup{\cite[8.5]{G-val}}\label{X/S}]
If $\Zz=\ES$, then by Theorem~\ref{cent.aut}\ref{dim S} we may
assume that $S=A$ and the geometric quotient $X/S$ exists. Then
$r(X/S)=0$ and by Propositions \ref{r(X)=r(G*)}, \ref{r=0}, orbits
of $G:X/S$ are projective homogeneous spaces.

Let $x\in X$ and $x\mapsto\bar{x}\in X/S$. We may assume
$G_{\bar{x}}\supseteq U$. Then $U$ preserves~$Sx$, and we have the
orbit map $U\to Ux\subseteq Sx$. As $U$ is an affine space (with
no non-constant invertible polynomials) and $Sx$ is a torus (whose
coordinate algebra is generated by invertibles), this map is
constant, whence $Ux=\{x\}$. Thus $X$ is horospherical.

Conversely, for horospherical $X$ put $Z=X^{U^{-}}$ and consider
the natural proper map $X'\to G\itimes{B^{-}}Z\onto X$. There are
natural maps $\ES(X')\onto\ES(X)$, $\Zz(X')\to\Zz(X)$. Restriction
of functions to $Z$ yields $\kk(X')^U=\kk(Z)$,
$\kk(X')^B=\kk(Z)^T$. Then any $T$-valuation of $\kk(Z)$ extends
to a $G$-valuation of~$\kk(X')$. Indeed, we may replace $Z$ by a
birationally isomorphic $B^{-}$-variety and assume that the
$T$-valuation corresponds to a $B^{-}$-stable divisor $D\subset Z$
($U^{-}$~acts trivially); then the desired $G$-valuation
corresponds to $D'={G\itimes{B^{-}}D}\subset X'$. By
Example~\ref{G=T}, $\Zz(X')=\ES(X')$, whence $\Zz(X)=\ES(X)$.
\end{proof}

We conclude this section by the description of $G$-valuations for
the residue field of a central valuation.
\begin{proposition}[{\cite[7.4]{G-val}}]\label{V(cent)}
Let $X$ be a $G$-model of~$K$, $D\subset X$ a $G$-stable prime
divisor with $v_D\in\Zz$, $X'$~the normal bundle of $X$ at~$D$.
Then:
\begin{enumerate}
\item\label{V(N_D)} $\Hyp(X')=\Hyp$ and $\Vv(X')=\Vv+\QQ{v_D}$;
\item\label{V(D)} $\Hyp(D)=\Hyp/\QQ{v_D}$ and $\Vv(D)$ is the
image of~$\Vv$.
\end{enumerate}
\end{proposition}
\begin{proof}
As usual, we may assume $K=\Quot{R}$, where $R$ is a rational
$G$-algebra. Then $K':=\kk(X')=\Quot\gr{R}$, where $\gr{R}$ is the
graded algebra associated with the filtration $R_{(d)}=\{f\in
R\mid v_D(f)\geq d\}$ of~$R$. Since $v_D$ is central, it is
constant on each $B$-eigenspace of~$R^U$, whence
$R^U\iso\gr(R^U)$. But $(\gr R)^U$ is a purely inseparable finite
extension of $\gr(R^U)$ by Corollary~\ref{(grA)^H}, hence
$(K')^U\supseteq K^U$ is a purely inseparable field extension.
This implies $\Hyp(X')=\Hyp$.

The $G$-invariant grading of $\gr{R}$ is defined by a central
1-torus acting on $\Vv(X')$ trivially by
Theorem~\ref{cent.aut}\ref{S:triv}. Hence it agrees with all
$G$-valuations. Thus $v\in\Vv(X')$ iff $v$ is non-negative at all
tail vectors of the form $\bar{f_0}/\bar{f_1}\dots\bar{f_s}$,
where $\bar{f_i}\in(\gr R)^{(B)}$ are homogeneous elements
represented by $f_i\in R$, and $v_D(f_0)=v_D(f_1)+\dots+v_D(f_s)$.
Replacing $\bar{f_i}$ by suitable powers, we may assume $f_i\in
R^{(B)}$. Thence $\Vv(X')$ is the set of all $v\in\Hyp$
non-negative on tail vectors of $R$ annihilated by~$v_D$, i.e.,
$\Vv(X')=\Vv+\QQ{v_D}$.

By Lemma~\ref{B-eigenfun}, $\kk(D)^U$~is a purely inseparable
extension of~$\kk(v_D|_{K^U})$, and $\kk(D)^B$ of~$K^B$. Hence
$\Hyp(D)=\Hyp/\QQ{v_D}$. Since $X'$ retracts onto~$D$,
$\Vv(D)=\Vv(X')/\QQ{v_D}$ by Corollary~\ref{res&ext(G-val)}.
\end{proof}

\section{Little Weyl group}
\label{Weyl}

In \ref{cotangent} we found out that important invariants of a
$G$-variety $X$ such as complexity, rank, and weight lattice,
which play an essential role in the embedding theory, are closely
related to the geometry of the cotangent bundle~$T^{*}X$. Knop
developed these observations further \cite{W_X}, \cite{inv.mot}
and described the cone of (central) $G$-valuations as a
fundamental chamber for the Galois group of a certain Galois
covering of~$T^{*}X$, called the little Weyl group of~$X$. As this
approach requires infinitesimal technique, we assume $\ch\kk=0$ in
this and the next section. We retain the notation and conventions
of~\ref{cotangent}.

The Galois covering of $T^{*}X$ is defined in terms of the moment
map. A disadvantage of the moment map $\Phi$ is that its image
$M_X$ can be non-normal and generic fibers can be reducible. A
remedy is to consider the ``Stein factorization'' of~$\Phi$. Let
$\widetilde{M}_X$ be the spectrum of the integral closure of
$\kk[M_X]$ (embedded via~$\Phi^{*}$) in~$\kk(T^{*}X)$. We may
assume that $X$ is smooth, whence $T^{*}X$ is smooth and normal,
and therefore $\kk[\widetilde{M}_X]\subseteq\kk[T^{*}X]$. It is
easy to see that $\kk(\widetilde{M}_X)$ is algebraically closed
in~$\kk(T^{*}X)$. Thus $\Phi$ decomposes into the product of a
finite morphism $\phi:\widetilde{M}_X\to M_X$ and the
\emph{normalized moment map}
$\widetilde\Phi:T^{*}X\to\widetilde{M}_X$ with irreducible generic
fibers. Set $\widetilde{L}_X=\widetilde{M}_X\by G$. We have the
quotient map $\widetilde\pi_G:\widetilde{M}_X\to\widetilde{L}_X$
and the natural finite morphism $\phi\by G:\widetilde{L}_X\to
L_X$.

The following result illustrates the role of the normalized moment
map in equivariant symplectic geometry.
\begin{proposition}\label{coll&inv}
The fields $\kk(\widetilde{M}_X)$ and $\kk(T^{*}X)^G$ are the
mutual centralizers of each other in $\kk(T^{*}X)$ w.r.t.\ the
Poisson bracket, and $\kk(\widetilde{L}_X)$ is the Poisson center
of both $\kk(\widetilde{M}_X)$ and $\kk(T^{*}X)^G$.
\end{proposition}
\begin{proof}
The field $\kk(M_X)$ is generated by Hamiltonians
$H_{\xi}=\Phi^{*}\xi$, $\xi\in\g$. Hence $f\in\kk(T^{*}X)$
Poisson-commutes with $\kk(M_X)$ iff $\{H_{\xi},f\}=\xi{f}=0$,
$\forall\xi\in\g$, i.e., $f$~is $G$-invariant. But then $f$ also
commutes with~$\kk(\widetilde{M}_X)$. Indeed, let $\mu_h$ be the
minimal polynomial of $h\in\kk(\widetilde{M}_X)$ over~$\kk(M_X)$.
Then $\{f,\mu_h(h)\}=\mu_h'(h)\{f,h\}=0\implies \{f,h\}=0$.
Therefore $\kk(T^{*}X)^G$ is the centralizer of $\kk(M_X)$ and
$\kk(\widetilde{M}_X)$.

Conversely, as generic orbits are separated by invariant
functions, $\g\alpha$~is the common kernel of $d_{\alpha}f$,
$f\in\kk(T^{*}X)^G$, for general $\alpha\in T^{*}X$. Hence $\Ker
d_{\alpha}\widetilde\Phi=\Ker d_{\alpha}\Phi=(\g\alpha)^{\sort}$
is generated by skew gradients of $f\in\kk(T^{*}X)^G$. It follows
that $h\in\kk(T^{*}X)$ commutes with $\kk(T^{*}X)^G$ iff $dh$
vanishes on $\Ker d_{\alpha}\widetilde\Phi=
T_{\alpha}\widetilde\Phi^{-1}\widetilde\Phi(\alpha)$ iff $h$ is
constant on $\widetilde\Phi^{-1}\widetilde\Phi(\alpha)$, because
generic fibers $\widetilde\Phi^{-1}\widetilde\Phi(\alpha)$ are
irreducible. Therefore $\kk(\widetilde{M}_X)$ is the centralizer
of~$\kk(T^{*}X)^G$.

Finally, $\kk(\widetilde{L}_X)=\kk(\widetilde{M}_X)^G=
\kk(\widetilde{M}_X)\cap\kk(T^{*}X)^G$, since quotient maps
$\pi_G$ and $\widetilde\pi_G$ separate generic orbits.
\end{proof}

Recall the local structure of an open subset of $X$ provided by
Corollary~\ref{loc.str.gen}: $\X\iso\Ru{P}\times Z$, where
$P=\Ru{P}\timessemi L$ is a parabolic, and the Levi subgroup $L$
acts on $Z$ with kernel $L_0\supseteq L'$, so that $Z\iso A\times
C$, $A=L/L_0$.

Generic $U$-orbits on $X$ coincide with generic $\Ru{P}$-orbits
and are of the form $\Ru{P}\times\{z\}$, $z\in Z$. Generic
$B$-orbits coincide with generic $P$-orbits and are of the form
$\Ru{P}\times A\times\{x\}$, $x\in C$. We have $T_xX=\Ru\p
x\oplus\ab x\oplus T_xC$, $\forall x\in C$.

Generic $U$- and $B$-orbits on $X$ form two foliations. Consider
the respective conormal bundles $\Uu\supseteq\Bb$. They are
$P$-vector bundles defined, e.g., over~$\X$. It follows from the
local structure that $\Uu\iso\Ru{P}\times T^{*}Z\iso\Ru{P}\times
A\times\ab^{*}\times T^{*}C$ and $\Bb\iso\Ru{P}\times A\times
T^{*}C$ over~$\X$. We have $\Uu/\Bb(x)= (\un x)^{\ann}/(\br
x)^{\ann}=(\br x/\un x)^{*}\iso
(\br/\un+\br_x)^{*}\subseteq\tr^{*}$. For $x\in Z$ we have
$\br_x=\br\cap\lv_0$, whence $\Uu/\Bb(x)\iso\ab^{*}$. As
$\Ad\Ru{P}$ acts on $\br/\un$ trivially, there is a canonical
isomorphism $\Uu/\Bb(x)\iso\ab^{*}$, $\forall x\in\X$. Therefore
$\Uu/\Bb\iso\ab^{*}\times\X$ is a trivial bundle over~$\X$, and we
have a canonical projection $\pi:\Uu\to\ab^{*}$, which is nothing
else, but the moment map for the $B$-action.

The bundle $\Uu/\Bb$ can be lifted (non-canonically) to $\Uu$
over~$\X$. Namely consider yet another foliation $\{g(\Ru{P}C)\mid
g\in P\}$ and let $\Aa$ be the respective conormal bundle. By the
local structure, $\Aa\iso\Ru{P}\times A\times\ab^{*}\times C$
over~$\X$, and $\Uu=\Aa\oplus\Bb$.

The isomorphism $\sigma:\ab^{*}\to\Aa(x)$, $x\in C$, defined by
the formula
%*
\begin{equation*}
\sigma(\lambda)=\begin{cases}
\lambda&\text{on $\ab x\iso\ab$}\\
0&\text{on $\Ru\p x\oplus T_xC$}
\end{cases}
\end{equation*}
%*
provides a section for~$\pi$. It depends on the choice of $x$ and
even more---we may replace $C$ by any subvariety in
$\X^{L_0}=\Ru{P}^{L_0}\times A\times C$ intersecting all
$P$-orbits transversally so that $x$ may be any point in $\X$ with
$P_x=L_0$ and $T_xC$ may be any ($L_0$-stable) complement to $\p
x$ in~$T_xX$.

Recall that $\ab$ embeds in $\lv$ as the orthocomplement
to~$\lv_0$. Consider the parabolic subgroup $Q\supseteq P$ having
the Levi decomposition $\Ru{Q}\timessemi M$,
$\Ru{Q}\subseteq\Ru{P}$, $M=Z_G(\ab)\supseteq L$.

\begin{lemma}\label{conorm&a^*}
There is a commutative square of dominant maps
%*
\begin{equation*}
\begin{CD}
   \Uu    @>{\Phi}>>            \ab
\text{\rlap{$\mathstrut\oplus\Ru\q\subseteq\g\iso\g^{*}$}} \\
@VV{\pi}V             @VV{\text{projection}}V \\
 \ab^{*}  @>{\sim}>>            \ab
\end{CD}
\end{equation*}
%*
\end{lemma}
\begin{proof}
Take $\alpha\in\Uu(x)$, $x\in\X$. Since all maps are
$P$-equivariant, we may assume $x\in C$ $\implies\Uu(x)\iso(\ab
x)^{*}\oplus T^{*}_xC$ and $\alpha=\sigma(\lambda)+\beta$ for some
$\lambda=\pi(\alpha)\in\ab^{*}$, $\beta\in T^{*}_xC=\Bb(x)$. Hence
$\langle\alpha,\xi x\rangle= \langle\lambda,\xi\rangle$,
$\forall\xi\in\p$
$\implies\Phi(\alpha)=\lambda\mod\p^{\ann}=\Ru\p$. Moreover,
$\Phi(\alpha)\in(\ab+\Ru\p)^{L_0}=\ab+\Ru\p^{L_0}\subseteq
\ab+\Ru\q$, because $\Ru\p^{L_0}\cap\m=\Ru\p^L=0$. Thus the square
is commutative. Finally, for general $\xi\in\ab$ we have
$\z(\xi)=\m\implies[\Ru\q,\xi]=\Ru\q\implies\xi+\Ru\q=\Ru{Q}\xi$
by Lemma~\ref{uni=>closed}. Therefore
$\overline{\Phi(\Uu)}=\ab+\Ru\q$.
\end{proof}
\begin{corollary}\label{M_X&a^*}
There is a commutative square
%*
\begin{equation*}
\begin{CD}
   \Uu    @>{\Phi}>>      M_X     \\
@VV{\pi}V             @VV{\pi_G}V \\
 \ab^{*}  @>{\pi_G}>>     L_X
\text{\rlap{$\mathstrut=\pi_G(\ab^{*})$}}
\end{CD}
\end{equation*}
%*
\end{corollary}

\begin{lemma}\label{norm(M_X)&a^*}
There exists a unique morphism $\psi:\ab^{*}\to\widetilde{L}_X$
making the following square commutative:
%*
\begin{equation*}
\begin{CD}
   \Uu    @>{\widetilde\Phi}>>    \widetilde{M}_X     \\
@VV{\pi}V                      @VV{\widetilde\pi_G}V \\
 \ab^{*}       @>{\psi}>>         \widetilde{L}_X
\end{CD}
\end{equation*}
%*
\end{lemma}
\begin{proof}
The uniqueness is evident. Take
$\psi=\widetilde\pi_G\widetilde\Phi\sigma$. The maps $\psi\pi$ and
$\widetilde\pi_G\widetilde\Phi$ coincide on $\sigma(\ab^{*})$, and
by Corollary~\ref{M_X&a^*} they map $\forall\alpha\in\Uu$ to one
and the same fiber of~$\phi\by G$. Thus for
$\forall\lambda\in\ab^{*}$ the irreducible subvariety
$\pi^{-1}(\lambda)\iso\Bb$ is mapped by
$\widetilde\pi_G\widetilde\Phi$ to the (finite) fiber of $\phi\by
G$ through~$\psi(\lambda)$, whence to~$\psi(\lambda)$.
\end{proof}

The normalization of $L_X=\pi_G(\ab^{*})$ equals
$\ab^{*}/W(\ab^{*})$, where $W(\ab^{*})=N_W(\ab^{*})/Z_W(\ab^{*})$
is the Weyl group of $\ab\subseteq\g$. By
Lemma~\ref{norm(M_X)&a^*}, there is a sequence of dominant finite
maps of normal varieties
$\ab^{*}\to\widetilde{L}_X\to\ab^{*}/W(\ab^{*})$. It follows from
the Galois theory that $\widetilde{L}_X\iso\ab^{*}/W_X$ for a
certain subgroup $W_X\subseteq W(\ab^{*})$ and the left arrow is
the quotient map.
\begin{definition}
The group $W_X$ is called the \emph{little Weyl group} of~$X$. It
is a subquotient of~$W$.
\end{definition}

By construction, $\widetilde{M}_X$, $\widetilde{L}_X$, and $W_X$
are $G$-birational invariants of~$X$. They are related to other
invariants such as the horospherical type~$S$.
\begin{proposition}[{\cite[6.4]{W_X}, \cite[7.3]{HC-hom}}]
$\widetilde{M}_{G/S}=
\widetilde{M}_X\times_{\widetilde{L}_X}\ab^{*}$ and
$\widetilde{M}_X=\widetilde{M}_{G/S}/W_X$.
\end{proposition}

As well as~$M_X$, $\widetilde{M}_X$~and $W_X$ are determined by a
generic $G$-orbit \cite[6.5.4]{W_X}. For functorial properties of
the normalized moment map and the little Weyl group, see
\cite[6.5]{W_X}. For geometric properties of the morphisms
$T^{*}X\to\widetilde{M}_X$, $\widetilde{M}_X\to\widetilde{L}_X$,
and $T^{*}X\to\widetilde{L}_X$, see \cite[7.4, 7.3, 6.6]{W_X},
\cite[\S\S5,7,9]{HC-hom}. $\widetilde{M}_X$~has rational
singularities \cite[4.3]{HC-hom}.

\begin{remark}\label{non-comm}
A non-commutative version of this theory was developed
in~\cite{HC-hom}. Here functions on $T^{*}X$ are replaced by
differential operators on $X$ and $\kk[\g^{*}]$ by~$\Uni\g$. The
analogue of $\kk[\widetilde{M}_X]$ consists of \emph{completely
regular} differential operators generated by velocity fields
locally on $X$ and ``at infinity''. Invariant completely regular
operators form a polynomial ring, which coincides with the center
of $\Dd(X)^G$ whenever $X$ is affine.
\begin{question}
It suffices to have $\kk(T^{*}X)^G=\Quot\kk[T^{*}X]^G$. Generally,
everything holds on the sheaf level.
\end{question}
This ring is isomorphic to $\kk[\rho+\ab^{*}]^{W_X}$
(``Harish-Chandra isomorphism''), where $\rho$ is half the sum of
positive roots and $W_X$ is naturally embedded in
$N_W(\rho+\ab^{*})$ being thus a subgroup, not only a subquotient,
of~$W$.
\end{remark}

\begin{example}
Take $X=G$ itself, with $G$ acting by left translations. Here
$T^{*}G\iso G\times\g^{*}$, and the moment map $\Phi$ is just the
coadjoint action map with irreducible fibers isomorphic to~$G$. We
have $A=T$, and
$\sigma:\tr^{*}\iso\tr\embeds\br\iso\un^{\ann}=\Uu(\1)$ is the
natural inclusion. Therefore $\widetilde{M}_G=M_G=\g^{*}$,
$\widetilde{L}_G=L_G=\g^{*}\by G\iso\tr^{*}/W$, and $W_G=W$.
\end{example}
\begin{example}
Let $X=G/T$, where $G$ is semisimple. Here $T^{*}X\iso
{G\itimes{T}(\un+\un^{-})}$, $A=\Ad_GT$, $M_X=\g^{*}$. The
subspace $\mathbf{e}+\g_{\mathbf{f}}\subset\un+\un^{-}$, where
$\mathbf{e}\in\un$, $\mathbf{f}\in\un^{-}$, $\mathbf{h}\in\tr$
form a principal $\sgl_2$-triple, is a cross-section for the
fibers of $\pi_G\Phi:T^{*}X\to\g^{*}\by G$. Indeed,
$\pi_G:\mathbf{e}+\g_{\mathbf{f}}\isoto\g^{*}\by G$
\cite[4.2]{Ad}. Hence
$\widetilde\pi_G\widetilde\Phi(\mathbf{e}+\g_{\mathbf{f}})$ is a
cross-section for the fibers of the finite map~$\phi\by G$. It
follows that $\widetilde{L}_{G/T}=L_{G/T}$, thence
$\widetilde{M}_{G/T}=M_{G/T}$ and $W_{G/T}=W$.
\end{example}
\begin{example}\label{W_(G/S)}
Consider a horospherical homogeneous space $X=G/S$. We have seen
in Theorem~\ref{hor.moment} that the moment map factors as
$\Phi=\overline\Phi\pi_A$, where
$\overline\Phi:G\itimes{P^{-}}\s^{\ann}\onto M_{G/S}$ is
generically finite proper and $\pi_A:G\itimes{S}\s^{\ann}\to
G\itimes{P^{-}}\s^{\ann}$ is the $A$-quotient map. It immediately
follows that
$\widetilde{M}_{G/S}=\Spec\kk[G\itimes{P^{-}}\s^{\ann}]$ and the
natural map $G\itimes{P^{-}}\s^{\ann}\to\widetilde{M}_{G/S}$ is a
resolution of singularities. The natural morphisms
$G\itimes{P^{-}}\s^{\ann}\to
G\itimes{P^{-}}\ab^{*}=G/P^{-}\times\ab^{*}\to\ab^{*}$ and
$\widetilde\pi_G:\widetilde{M}_{G/S}\to\widetilde{L}_{G/S}$ are
rational $G$-quotient maps. Indeed, $P^{-}$ acts on each fiber
$\lambda+\Ru{\p^{-}}$ of $\s^{\ann}\to\ab^{*}$ generically
transitively \cite[5.5]{Ad}, and fibers of $\widetilde\pi_G$ have
a dense orbit, because fibers of $\pi_G:M_{G/S}\to L_{G/S}$ do.
Passing to rational quotients, we see that
$\ab^{*}\to\widetilde{L}_{G/S}$ is birational, whence isomorphic.
Thus $W_{G/S}=\{\1\}$.
\begin{question}
An example with $W_X\neq W(\ab^{*})$, e.g.,
$\Sp_4/\kk^{\times}\times\Sp_2$ or $\SO_{2n+1}/\GL_n$; symmetric
case.
\end{question}
\end{example}
%\begin{example}[{\cite[Ex.4.3.3]{val.sph}}]
%Now let $X=G/H$, $G=\Sp_4$,
%$H=\kk^{\times}\times\Sp_2\subset\Sp_2\times\Sp_2\subset G$.
%The coisotropy representation of $H$ is
%$\h^{\ann}=\kk^1\oplus(\kk^1)^{*}\oplus\kk^2\oplus(\kk^2)^{*}$,
%where $\kk^1,\kk^2$ are the trivial and the tautological
%$\Sp_2$-module acted on by $\kk^{*}$ via the characters
%$2\eps_1,\eps_1$, respectively. Theorem~\ref{c&r(G/H)} yields
%$A=\Ad_GT$. The algebra
%$\kk[\widetilde{L}_{G/H}]\iso\kk[\h^{\ann}]^H$ is freely
%generated by two quadratic invariants $ab$, $\langle
%x,y\rangle$, where $(a,b,x,y)\in\h^{\ann}$. Hence
%$W_{G/H}\iso\ZZz_2^2$ is a subgroup of $W(\Sp_4)$ generated by
%the reflections along two orthogonal roots.
%\begin{question}
%long
%\end{question}
%This example shows that $W_X$ may be a proper subgroup
%of~$W(\ab^{*})$.
%\end{example}

The last example admits a conversion.
\begin{proposition}\label{W_hor}
$X$~is horospherical iff $W_X=\{\1\}$.
\end{proposition}
\begin{proof}
A horospherical variety of type $S$ is birationally $G$-isomorphic
to $G/S\times C$ by Proposition~\ref{loc.str.hor}. Therefore it
suffices to consider $X=G/S$, but then $W_X=\{\1\}$ by
Example~\ref{W_(G/S)}.

Conversely, suppose $W_X=\{\1\}$ and consider the morphism
$\Psi=\widetilde\pi_G\widetilde\Phi:T^{*}X\to\widetilde{L}_X=
\ab^{*}$. Then $\Psi^{*}$ embeds $\ab$ into the space of fiberwise
linear $G$-invariant functions on~$T^{*}X$, which restrict to
linear functions on $\sigma(\ab^{*})\iso\ab^{*}$.

Geometrically, $\Psi^{*}\ab$~is an Abelian subalgebra of
$G$-invariant vector fields on $X$ tangent to $G$-orbits.
Furthermore, $\Psi^{*}\ab=\pi^{*}\ab=0$ on $\Bb$ by
Lemma~\ref{norm(M_X)&a^*}, hence $\Psi^{*}\ab$ is tangent to
generic $P$-orbits. It follows that $\Psi^{*}\ab$ restricts
to~$Px$, $x\in C$, as an Abelian subalgebra in the algebra
$(\p/\p_x)^{P_x}=\ab\oplus\Ru\p^{L_0}$ of $P$-invariant vector
fields on~$Px$, and $\Psi^{*}\xi(x)=\xi{x}\mod\Ru\p^{L_0}x$,
$\forall\xi\in\ab$, whence $\Psi^{*}\ab$ projects onto~$\ab$. As
$\z(\ab)\cap\Ru\p^{L_0}=\Ru\p^L=0$, $\Psi^{*}\ab$~is conjugated to
$\ab$ by a unique $g_x\in\Ru{P}^{L_0}$. Moving each $x\in C$
by~$g_x$, we may assume $\Psi^{*}\xi(x)=\xi{x}$,
$\forall\xi\in\ab$, $x\in C$ (or $\forall x\in Z=AC$).

Therefore velocity fields of $A:Z$ extend to $G$-invariant vector
fields on~$X$. These vector fields can be integrated to an
$A$-action on $X$ by central automorphisms, which restricts to the
natural $A$-action on $\X=P\itimes{L}Z$ provided by $A:Z$. (The
induced action $A:T^{*}X$ integrates the invariant collective
motion, cf.~\ref{inv.mot}.)
\begin{question}
$\ab:T^{*}X$ is given by skew gradients of $\ab:X$, i.e., of
$\Psi^{*}\ab$.
\end{question}
We conclude by Theorems~\ref{cent.aut}--\ref{cent(hor)} that $X$
is horospherical.
\end{proof}

Here comes the main result linking the little Weyl group with
equivariant embeddings.
\begin{theorem}[{\cite[7.4]{inv.mot}}]\label{W_X&Z(X)}
The little Weyl group $W_X$ acts on $\ab^{*}$ as a
crystallographic reflection group preserving the lattice~$\RG(X)$,
and the central valuation cone $\Zz(X)$ is its fundamental chamber
in $\ES=\Hom(\RG(X),\QQ)$.
\end{theorem}
The proof relies on the integration of the invariant collective
motion in $T^{*}X$ and the study of the asymptotic behavior of its
projection to~$X$, see~\ref{inv.mot}. The description of
$\Zz=\Zz(X)$ as a fundamental chamber of a crystallographic
reflection group was first obtained by Brion \cite{W_X(sph)} in
the spherical case and generalized to arbitrary complexity in
\cite[\S9]{G-val}. In particular, $\Zz$~is a cosimplicial cone.

From this theorem, Knop derived the geometric look of all
valuation cones.
\begin{corollary}[{\cite[\S9]{G-val}}]\label{cosimp}
The cones $\Vv_{\vv}$ are cosimplicial.
\begin{question}
Prove it?
\end{question}
\end{corollary}

In proving Theorem~\ref{W_X&Z(X)}, we shall use its formal
consequence:
\begin{lemma}\label{W_X:cent.tor}
$W_X$~acts trivially on $\Zz\cap-\Zz$.
\end{lemma}
\begin{proof}
By Theorem~\ref{cent.aut}, there exists a torus $E$ of central
automorphisms such that $\e=(\Zz\cap-\Zz)\otimes\kk$ w.r.t.\ the
canonical embedding $\e\embeds\ab=\ES\otimes\kk$. Consider the
action $G^{+}=G\times E:X$ and indicate all objects related to
this action by the superscript~``+''. In particular, $A^{+}\iso A$
is the quotient of $A\times E$ modulo the antidiagonal copy
of~$E$, and $\Ch(A^{+})\subset\Ch(A)\oplus\Ch(E)$ is the graph of
the restriction homomorphism $\Ch(A)\to\Ch(E)$.

Obviously, $\Phi=\tau\Phi^{+}$, where $\tau:(\g^{+})^{*}\to\g^{*}$
is the canonical projection. It follows that
$\widetilde\Phi=\widetilde\tau\widetilde\Phi^{+}$ for a certain
morphism $\widetilde\tau:\widetilde{M}^{+}_X\to\widetilde{M}_X$.
The subalgebra $\kk[\widetilde{M}_X]$ is integrally closed
in~$\kk[\widetilde{M}^{+}_X]$, whence $\kk[\widetilde{L}_X]$ is
integrally closed in~$\kk[\widetilde{L}^{+}_X]$. On the other
hand, we have a commutative square
%*
\begin{equation*}
\begin{CD}
   (\ab^{+})^{*}         @>{\tau}>{\sim}>          \ab^{*}     \\
   @VV{\psi^{+}}V                                 @VV{\psi}V   \\
\widetilde{L}^{+}_X @>{\widetilde\tau\by{G}}>> \widetilde{L}_X
\end{CD}
\end{equation*}
%*
where $\psi,\psi^{+}$, and hence $\widetilde\tau\by{G}$ are finite
morphisms. Hence $\widetilde{L}^{+}_X\iso\widetilde{L}_X$ and
$W^{+}_X=W_X$.

It follows that $W_X$ preserves
$(\ab^{+})^{*}\subset\ab^{*}\oplus\e^{*}$ and acts trivially on
the \ordinal{2} summand
$\e^{*}\iso(\ab^{+})^{*}/(\ab^{+})^{*}\cap\ab^{*}\iso
\ab^{*}/(\ab^{+})^{*}\cap\ab^{*}$. Thus $W_X$ acts trivially on
$\e$ embedded in~$\ab$.
\end{proof}

\section{Invariant collective motion}
\label{inv.mot}

The skew gradients of functions in~$\kk[L_X]$
(or~$\kk[\widetilde{L}_X]$) pulled back to $T^{*}X$ generate an
Abelian flow of $G$-automorphisms preserving $G$-orbits, which is
called the invariant collective motion, see~\ref{cotangent}.
Restricted to a generic orbit $G\alpha\subset T^{*}X$, the
invariant collective motion gives rise to a connected Abelian
subgroup $A_{\alpha}=(G_{\Phi(\alpha)}/G_{\alpha})^0\subseteq
\Aut_GG\alpha$. It turns out that $A_{\alpha}\iso A$. However, in
general, this isomorphism cannot be made canonical in order to
produce an $A$-action on (an open subset of) $T^{*}X$ integrating
the invariant collective motion. This obstruction is overcome by
unfolding the cotangent variety by means of a Galois covering with
Galois group~$W_X$.

\begin{definition}
The fiber product
$\widehat{T}X=T^{*}X\times_{\widetilde{L}_X}\ab^{*}$ is called the
\emph{polarized cotangent bundle} of~$X$. Since generic fibers of
$T^{*}X\to\widetilde{L}_X$ are irreducible, $\widehat{T}X$~is
irreducible. Actually $\widehat{T}X$ is an irreducible component
of $T^{*}X\times_{L_X}\ab^{*}$, $W_X$~is its stabilizer in
$W(\ab^{*})$ acting on the set of components, and $\widehat{T}X\to
T^{*}X=\widehat{T}X/W_X$ is a rational Galois cover.
\end{definition}

Consider the \emph{principal stratum} $\ab^{\pr}\subseteq\ab^{*}$
obtained by removing all proper intersections with kernels of
coroots and with $W$-translates of $\ab^{*}$ in~$\tr^{*}$. The
group $W(\ab^{*})$ acts on $\ab^{\pr}$ freely. Put
$L_X^{\pr}=\pi_G(\ab^{\pr})=\ab^{\pr}/W(\ab^{*})$, the quotient
map being an {\'e}tale finite Galois covering. The preimages of
$L_X^{\pr}$ in various varieties under consideration will be
called \emph{principal strata} and denoted by the
superscript~``$\pr$''.

In particular, $\widehat{T}^{\pr}X\subset\widehat{T}X$ is a smooth
open subvariety (provided that $X$ is smooth) and the projection
$\widehat{T}^{\pr}X\onto T^{\pr}X\subseteq T^{*}X$ is an {\'e}tale
finite quotient map by~$W_X$. The $G$-invariant symplectic
structure on $T^{*}X$ is pulled back to~$\widehat{T}^{\pr}X$ so
that $\widehat{T}^{\pr}X\to T^{*}X\to M_X$ is the moment map.

The invariant collective motion on $\widehat{T}^{\pr}X$ is
generated by the skew gradients of Poisson-commuting functions
from~$\pi^{*}\ab$, where $\pi:\widehat{T}^{\pr}X\to\ab^{*}$ is the
other projection. These skew gradients constitute a commutative
$r$-dimensional subalgebra of Hamiltonian vector fields
($r=r(X)=\dim\ab$). Our aim is to show that these vector fields
are the velocity fields of a symplectic $A$-action so that $\pi$
is the respective moment map \cite[App.5]{SG}.
\begin{remark}
In particular, it will follow that the $W_X$-action on $\ab$ lifts
to $A$, so that $\widehat{T}^{\pr}X$ comes equipped with the
Poisson $G\times(W_X\semitimes A)$-action.
%Indeed, given such an $A$-action, $W_X$~permutes $A$-orbits.
%Since orbit maps $A\to A\hat\alpha$ are isomorphic for
%general $\hat\alpha\in\widehat{T}^{\pr}X$, each $w\in W_X$ gives
%rise to a family of automorphisms
%$A\iso A\hat\alpha\stackrel{w}\to Aw\hat\alpha\iso A$
%depending continuously on~$\hat\alpha$ and preserving the unity,
%hence a group automorphism of $A$ independent of $\hat\alpha$ by
%the rigidity of tori \cite[16.3]{Agr}. This yields the desired
%$W_X$-action.
\end{remark}

Following~\cite{inv.mot}, we shall restrict our considerations to
the symplectically stable case (Definition~\ref{symp.stab}) for
technical reasons.

\begin{proposition}
$G:T^{*}X$ is symplectically stable iff $T^{*}X=\overline{G\Uu}$.
\end{proposition}
\begin{proof}
In the notation of Lemma~\ref{conorm&a^*},
$\Phi(G\Uu^{\pr})=G(\ab^{\pr}+\Ru\q)=G\ab^{\pr}$. Hence density of
$G\Uu$ implies symplectic stability. Conversely, in the
symplectically stable case $\Ru{P^{-}}\alpha$ is transversal to
$\Uu$ for general $\alpha\in\Uu$. Indeed, we may assume
$\Phi(\alpha)\in\ab^{\pr}$, but then
$[\Ru{\p^{-}},\Phi(\alpha)]=\Ru{\p^{-}}$ is transversal to
$\overline{\Phi(\Uu)}=\ab+\Ru\p$. Therefore
$\dim\overline{\Ru{P^{-}}\Uu}=\dim\Ru{P}+\dim\Uu=\dim T^{*}X$.
\end{proof}

Suppose that the action $G:T^{*}X$ is symplectically stable. We
have observed in \ref{cotangent} that $M_X^{\pr}\iso
G\itimes{N(\ab)}\ab^{\pr}$. Then $\widetilde{M}_X^{\pr}\iso
M_X^{\pr}\times_{L_X}\widetilde{L}_X\iso G\itimes{N_X}\ab^{\pr}$,
where $N_X\subseteq N(\ab)$ is the extension of $W_X$ by
$Z(\ab)=L$. Hence $T^{\pr}X\iso G\itimes{N_X}\Sigma$ has a
structure of a homogeneous bundle over~$G/N_X$. Fibers of this
bundle are called \emph{cross-sections}. They are smooth and
irreducible, because generic fibers of $\widetilde\Phi$ are
irreducible. We may choose a canonical $N_X$-stable
cross-section~$\Sigma$, namely the unique cross-section in
$\Phi^{-1}(\ab^{\pr})$ intersecting~$\Uu$.
\begin{remark}
In fact, $\Uu\cap\Sigma$ is dense in~$\Sigma$. Indeed, $\Sigma\cap
T^{*}\X\subseteq\Uu$.
\end{remark}

\begin{lemma}
The kernel of $N_X:\Sigma$ is~$L_0$.
\end{lemma}
\begin{proof}
By Lemma~\ref{conorm&a^*}, $\Phi(\Uu^{\pr})=\ab^{\pr}+\Ru\p\iso
P\itimes{L}\ab^{\pr}\iso\Ru{P}\times\ab^{\pr}$. Hence
$\Uu^{\pr}=P\itimes{L}(\Uu\cap\Sigma)\iso
\Ru{P}\times(\Uu\cap\Sigma)$. On the other hand,
$\Uu|_{\X}=P\itimes{L}\Uu|_Z\iso\Ru{P}\times\Uu|_Z$, and all the
stabilizers of $L:\Uu|_Z\iso T^{*}Z$ are equal to~$L_0$. It
follows that generic stabilizers of $L:\Sigma$ are $P$-conjugate
to~$L_0$, hence coincide with~$L_0$.
\begin{question}
Just note that $\Uu\cap\Sigma$ and $\Uu|_Z$ are $L$-birationally
isomorphic.
\end{question}
\end{proof}
\begin{corollary}\label{W_X:A}
$W_X$~acts on $A=L/L_0$, i.e., preserves the character lattice
$\Ch(A)=\RG(X)\subset\ab^{*}$.
\end{corollary}
\begin{remark}\label{sgp(ss)}
The lemma implies Theorem~\ref{sgp(T^*)} in the symplectically
stable case. This observation was made in \cite[\S4]{inv.mot}.
\end{remark}

\begin{lemma}\label{int(coll.mot)}
The $A$-action integrates the invariant collective motion
on~$\Sigma$.
\end{lemma}
\begin{proof}
The skew gradient of $\Phi^{*}f$ ($f\in\kk(M_X)$) at $\alpha\in
T^{*}X$ equals $(d_{\Phi(\alpha)}f)\cdot\alpha$, where
$d_{\Phi(\alpha)}f$ is considered as an element of $\g^{**}=\g$ up
to a shift from $(T_{\Phi(\alpha)}M_X)^{\ann}$. Indeed, the skew
gradient of a function is determined by its linear portion at a
point, and for linear functions $f\in\g$ the assertion holds by
the definition of~$\Phi$.

If $\alpha\in\Sigma$, then $T_{\Phi(\alpha)}M_X=
\ab+[\g,\Phi(\alpha)]=\ab\oplus\Ru\p\oplus\Ru{\p^{-}}$,
$(T_{\Phi(\alpha)}M_X)^{\ann}=\lv_0$. The differentials
$d_{\Phi(\alpha)}f$ of $f\in\kk[L_X^{\pr}]$ generate the conormal
space of $G\cdot\Phi(\alpha)$ in $M_X$ at~$\Phi(\alpha)$, i.e.,
$[\g,\Phi(\alpha)]^{\ann}/(T_{\Phi(\alpha)}M_X)^{\ann}=
\lv/\lv_0=\ab$. Thus the invariant collective motion at $\alpha$
is $\ab_{\alpha}=\lv\alpha=\ab\alpha$.
\end{proof}

Translation by $G$ permutes cross-sections transitively and
extends the $A$-action to each cross-section. These actions
integrate the invariant collective motion, but in general, they do
not globalize to a regular $A$-action on the whole cotangent
bundle, due to non-trivial monodromy.

However, unfold $T^{\pr}X$ to $\widehat{T}^{\pr}X=
T^{\pr}X\times_{\widetilde{L}_X}\ab^{\pr}\iso
G\itimes{L}\widehat\Sigma$, where $\widehat\Sigma=
\{\hat\alpha=(\alpha,\Phi(\alpha))\mid\alpha\in\Sigma\}$. We
retain the name ``cross-sections'' for the fibers of this
homogeneous bundle, which are isomorphic to the cross-sections
in~$T^{\pr}X$. Now there is a natural $A$-action on
$\widehat{T}^{\pr}X$ provided by $A:\widehat\Sigma\iso\Sigma$,
which integrates the invariant collective motion
on~$\widehat{T}^{\pr}X$. The $W_X$-action on $\widehat{T}^{\pr}X$
is induced from the $N_X$-action on $G\times\widehat\Sigma$ given
by $n(g,\hat\alpha)=(gn^{-1},n\hat\alpha)$, $\forall n\in N_X$,
$g\in G$, $\hat\alpha\in\widehat\Sigma$. We sum up in the
following
\begin{theorem}[{\cite[4.1--4.2]{inv.mot}}]
There is a Poisson $G\times(W_X\semitimes A)$-action on
$\widehat{T}^{\pr}X$ with the moment map
$\Phi\times\pi:\widehat{T}^{\pr}X\to\g^{*}\oplus\ab^{*}$.
\end{theorem}
\begin{proof}
It remains only to explain why $\pi$ is the moment map for the
$A$-action. Take any $\hat\alpha\in\widehat\Sigma$ over
$\alpha\in\Sigma$. By (the proof of) Lemma~\ref{int(coll.mot)},
$\forall\xi\in\ab\ \exists f\in\kk[L_X^{\pr}]:\
d_{\Phi(\alpha)}f=\xi\mod\lv_0$. The skew gradient of $\pi^{*}\xi$
at $\hat\alpha$ is pulled back from that of $\Phi^{*}f$
at~$\alpha$, i.e., from~$\xi\alpha$, hence it equals
$\xi\hat\alpha$.
\begin{question}
Prove it directly in Lemma~\ref{int(coll.mot)}.
\end{question}
We conclude by $G$-equivariance.
\end{proof}

In particular, the orbit of the invariant collective motion
through $\hat\alpha\in\widehat{T}^{\pr}X$ over $\alpha\in
T^{\pr}X$ is $A\hat\alpha=G_{\Phi(\alpha)}\hat\alpha\iso
G_{\Phi(\alpha)}/G_{\alpha}$. For the purposes of the embedding
theory it is important to study the projections of these orbits to
$X$ and their boundaries.
\begin{definition}
A \emph{flat} in $X$ is
$F_{\alpha}=\pi_X(A\hat\alpha)=G_{\Phi(\alpha)}x$, where
$\alpha\in T^{\pr}_xX$, $\hat\alpha\in\widehat{T}^{\pr}X$ lies
over~$\alpha$, and $\pi_X:T^{*}X\to X$ is the canonical
projection. The composed map $A\to A\hat\alpha\to F_{\alpha}$ is
called the \emph{polarization} of the flat.
\begin{question}
Consider flats for symmetric spaces.
\end{question}
\end{definition}

For general $\alpha$ the polarization map is isomorphic: indeed,
w.l.o.g.\ assume $\alpha\in\Sigma\cap T^{*}_x\X$ $\implies
G_{\Phi(\alpha)}=L$, $G_{\Phi(\alpha)}\cap G_x=G_{\alpha}=L_0$
$\implies F_{\alpha}\iso A$. Generic flats are nothing else, but
$G$-translates of $L$- (or $A$-) orbits in~$Z$, under appropriate
choice of~$Z$.
\begin{question}
Better understanding of flats w.r.t.\ local structure.
\end{question}
Namely, by Lemma~\ref{conorm&a^*}, there is a commutative diagram
%*
\begin{equation}\label{loc.str.A}
\begin{CD}
\text{\llap{$\X\times\ab^{*}\iso\mathstrut$}}
   \Aa    @>{\Phi}>>  \ab
\text{\rlap{$\mathstrut\oplus\Ru\p$}} \\
@VV{\pi}V            @VVV \\
 \ab^{*}  @>{\sim}>>  \ab
\end{CD}
\end{equation}
%*
For $\forall\lambda\in\ab^{\pr}$ we have
$\lambda+\Ru\p=P\cdot\lambda\iso P/L$, hence $\X\iso
P\itimes{L}Z^{\lambda}$, where
$Z^{\lambda}=\pi_X(\Phi^{-1}(\lambda)\cap\Aa)$. Clearly, all
$L$-orbits in $Z=Z^{\lambda}$ are flats. On the other hand, for
$\forall\alpha\in\Sigma\cap T^{*}_x\X$ it is easy to construct a
subvariety $C\subset\X^{L_0}$ through $x$ intersecting all
$P$-orbits transversally such that $\alpha=0$ on $T_xC$, whence
$x\in Z^{\Phi(\alpha)}$.

The rigidity of torus actions implies that the closures of generic
flats are isomorphic.
\begin{proposition}[{\cite[\S6]{inv.mot}}]\label{cl(flats)}
The closures $\overline{F}_{\alpha}$ for general $\alpha\in
T^{\pr}X$ are $A$-isomorphic toric varieties, and the $W_X$-action
on $A\iso F_{\alpha}$ extends to~$\overline{F}_{\alpha}$.
\end{proposition}
\begin{proof}
We explain the affine case, the general case being reduced to this
one by standard techniques of invariant quasiprojective open
coverings and affine cones. Generic flats are $G$-translates of
generic $L$-orbits in~$\overline{\pi_X(\Sigma)}$. We may assume
that $X$ is embedded into a $G$-module. Since
$\overline{\pi_X(\Sigma)}$ is $N_X$-stable, the set of
eigenweights of $A=L/L_0$ in $\overline{\pi_X(\Sigma)}$ is
$W_X$-stable. For general $\alpha\in\Sigma$,
$\kk[\overline{F}_{\alpha}]$~is just the semigroup algebra
generated by these eigenweights.
\end{proof}

The following result is crucial for interdependence between flats
and central valuations. It partially describes the boundary of a
generic flat.
\begin{proposition}[{\cite[7.3]{inv.mot}}]\label{flat&div}
Let $D\subset X$ be a $G$-stable divisor with $v=v_D\in\Zz$. The
closure $\overline{F}_{\alpha}$ of a generic flat contains
$A$-stable prime divisors $D_{wv}\subseteq D$, $w\in W_X$, that
correspond to $wv$ regarded as $A$-valuations of~$\kk(A)$.
Furthermore, $\overline{F}_{\alpha}$~is smooth along~$D_{wv}$.
\end{proposition}
\begin{proof}
W.l.o.g.\ we may assume $\alpha\in\Sigma$. The $W_X$-action on
$\overline{F}_{\alpha}$ is given by $w:\overline{F}_{\alpha}\to
\overline{F}_{n\alpha}\to\overline{F}_{\alpha}$, where the left
arrow is the translation by $n\in N_X$ representing $w\in W_X$ and
the right arrow is the unique $A$-isomorphism mapping $n\alpha$
back to~$\alpha$. Since $D$ is $N_X$-stable, it suffices to prove
the assertion for $w=\1$.

Shrinking $\X$ if necessary, we find a $B$-chart $X_0$
intersecting $D$ such that $\X=X_0\setminus D$.
\begin{lemma}\label{mom(conorm)}
The morphism $\X\times\ab^{*}\to\ab+\Ru\p$ in \eqref{loc.str.A}
extends to $X_0\times\ab^{*}$.
\begin{question}
Replace $D$ by $Y$, $\Dd^B_Y=\emptyset$.
\end{question}
\end{lemma}
\begin{proof}
Trivializing sections of $\Aa\iso\X\times\ab^{*}$ corresponding to
$\lambda\in\RG$ are~$d\ef{\lambda}/\ef{\lambda}$, where
$\ef{\lambda}$ are $B$-eigenfunctions on $\X$ that are constant
on~$\Ru{P}C$. These sections extend to sections of $T^{*}X(\log
D)$ over~$X_0$, which trivialize
\begin{question}
Because any linear combination of $d\ef{\lambda}/\ef{\lambda}$
maps to the same combination of $\lambda\bmod\Ru\p$ under the
logarithmic moment map.
\end{question}
the subbundle $\Aa(\log D)= \overline\Aa\subseteq T^{*}X_0(\log
D)$. The moment map of $T^{*}X(\log D)$ restricted to $\Aa(\log
D)$ provides the desired extension.
\end{proof}

Consequently $X_0\iso P\itimes{L}Z_0$,
$Z_0=\overline{Z^{\Phi(\alpha)}}$, and $F_{\alpha}=Ax$ is a
generic $A$-orbit in~$Z_0$. The proposition stems from
\begin{lemma}
After possible shrinking of~$X_0$, $Z_0\iso F\times C$, where
$F=\overline{Ax}$ is the closure of a generic $A$-orbit in~$Z_0$.
\end{lemma}
\begin{proof}
Since $v$ is central and by Lemma~\ref{B-eigenfun}, the
restriction of functions identifies $\kk(D)^B\iso\kk(Z_0\cap D)^A$
with $\kk(X)^B\iso\kk(Z_0)^A\iso\kk(C)$. Hence removing
zeroes/poles of a $B$-invariant function preserves non-empty
intersection with~$D$. In particular, we may assume $\kk(Z_0\cap
D)^A=\Quot\kk[Z_0\cap D]^A$ and $\kk[Z_0\cap
D]^A\iso\kk[Z_0]^A\iso\kk[C]$ by shrinking~$X_0$. We have
$\kk[Z]=\kk[\,{\ef{\lambda}\mid\lambda\in\RG}\,]\otimes\kk[C]$ and
$\kk[Z_0]\subseteq\kk[\,{\ef{\lambda}\mid\lambda\in\RG_0}\,]
\otimes\kk[C]$, where $\RG_0$ is the weight semigroup of~$Z_0$.
There exist $h_{\lambda}\in\kk[C]$ such that
$h_{\lambda}\ef{\lambda}\in\kk[Z_0]$. Shrinking~$X_0$, we may
assume $\ef{\lambda}\in\kk[Z_0]$, $\forall\lambda\in\RG_0$
(because $\RG_0$ is finitely generated). Hence $Z_0\iso F\times
C$, where $F=\Spec\kk[\RG_0]$.
\end{proof}
\pushQED{}
\end{proof}

Now we explain how to deal with non-symplectically stable case.

We may assume $X$ to be quasiprojective. By $\widehat{X}$ denote
the cone over $X$ without the origin. In the notation of
Remark~\ref{aff.cone}, the $\widehat{G}$-action on
$T^{*}\widehat{X}$ is symplectically stable by
Proposition~\ref{qaff=>ss}.

The quotient space $T^{*}\widehat{X}/\kk^{\times}$ is a vector
bundle over $X$ containing $T^{*}X$ as a subbundle, the quotient
bundle being the trivial line bundle. The moment map for
$\widehat{G}:T^{*}\widehat{X}$ factors through
$T^{*}\widehat{X}/\kk^{\times}$, so that there is a commutative
diagram
%*
\begin{equation*}
\begin{array}{ccccc}
&&T^{*}\widehat{X}&& \\
&&\downarrow&& \\
T^{*}X&\subset&T^{*}\widehat{X}/\kk^{\times}&
\stackrel\Pi\longrightarrow&\kk \\
\downarrow\text{\rlap{$\Phi$}}&&
\downarrow\text{\rlap{$\widehat\Phi$}}&&\parallel \\
M_X&\subset&M_{\widehat{X}}&\longrightarrow&\kk \\
\end{array}
\end{equation*}
%*
Here $\Pi$ is induced by the evaluation at the Euler vector field
in~$\widehat{X}$, i.e., by the moment map for the
$\kk^{\times}$-action, and the lower right arrow is the projection
of $M_{\widehat{X}}\subseteq\hat\g^{*}=\g^{*}\oplus\kk$ to~$\kk$.
$T^{*}X$~and $M_X$ are the zero-fibers of the respective maps
to~$\kk$. Also, we have $\ab^{*}=\hat\ab^{*}\cap\g^{*}$. The
morphism $\widehat\Phi$ factors
through~$\widetilde{M}_{\widehat{X}}$, hence $\Phi$ factors
through the zero-fiber $M_X'$ of
$\widetilde{M}_{\widehat{X}}\to\kk$. As $M_X'\to M_X$ is finite,
there is a commutative diagram
%*
\begin{equation*}
\begin{array}{cccc@{}c}
T^{*}X&&\subset&&T^{*}\widehat{X}/\kk^{\times} \\
\downarrow\text{\rlap{$\widetilde\Phi$}}&&&&\downarrow \\
\widetilde{M}_X&\to&M_X'&\subset&\widetilde{M}_{\widehat{X}} \\
\end{array}
\end{equation*}
%*
Passing to quotients, we obtain
%*
\begin{equation*}
\begin{array}{ccc}
\ab^{*}&\subset&\hat\ab^{*} \\
\downarrow&&\downarrow \\
\widetilde{L}_X&\longrightarrow&\widetilde{L}_{\widehat{X}} \\
\end{array}
\end{equation*}
%*
whence $W_X\subseteq W_{\widehat{X}}$. Actually these groups
coincide by Theorem~\ref{W_X&Z(X)}. By Corollary~\ref{W_X:A},
$W_{\widehat{X}}$~preserves~$\RG(\widehat{X})$, whence $W_X$
preserves $\RG(X)=\RG(\widehat{X})\cap\ab^{*}$.

Instead of flats, one considers \emph{twisted flats} defined as
projectivizations of usual flats in~$\widehat{X}$. The above
results on flats and their closures in $\widehat{X}$ descend to
twisted flats in~$X$. If $T^{*}X$ is symplectically stable, then
$T^{\pr}X\subset T^{\pr}\widehat{X}/\kk^{\times}$, and flats are a
particular case of twisted flats.

\begin{example}
Let $G=\SL_2$ and $\widehat{X}\subset V(3)$ be the variety of
(nonzero) degenerate binary cubic forms (in the variables $x,y$).
Essentially $\widehat{G}=\GL_2$. The form $v=xy^2$ has the open
orbit $\widehat\HS\subset\widehat{X}$ and the stabilizer
$\widehat{H}=\left\{\left(\begin{smallmatrix}
t^2 &  0 \\
 0  & t^{-1}
\end{smallmatrix}\right)\mid t\in\kk^{\times}\right\}$
in~$\widehat{G}$. Passing to projectivizations, we obtain a
hypersurface $X\subset\PP(V(3))$ with the open orbit $\HS=G\langle
v\rangle$, and $G_{\langle v\rangle}=:H=
\left\{\left(\begin{smallmatrix}
t & 0 \\
0 & t^{-1}
\end{smallmatrix}\right)\mid t\in\kk^{\times}\right\}$. The
flats through $\langle v\rangle$ are the orbits of the isotropy
groups of $\h^{\ann}=\left\{\left(\begin{smallmatrix}
 0  & {*} \\
{*} &  0
\end{smallmatrix}\right)\right\}$ in~$G$, i.e., of
$L=\left\{\left(\begin{smallmatrix}
a & b \\
b & a
\end{smallmatrix}\right)\mid a^2-b^2=1\right\}$ and its
$H$-conjugates. However the twisted flats are the orbits of the
stabilizers in $G$ of non-degenerate matrices from
$\hat\h^{\ann}=\left\{\left(\begin{smallmatrix}
 c  & {*} \\
{*} & 2c
\end{smallmatrix}\right)\mid c\in\kk\right\}$, i.e., the orbits
of arbitrary $1$-tori in~$G$.

The boundary of the open orbit is a single orbit
$X\setminus\HS=G\langle v_0\rangle$, $v_0=y^3$, with the
stabilizer $G_{v_0}=B$. Put $Y=\PP(\langle v_0,v\rangle)$, a
$B$-stable subspace in~$X$. The natural bijective morphism
$\widetilde{X}=G\itimes{B}Y\to X$ is a desingularization. The
closures of generic twisted flats are isomorphic to $\PP^1$ and
intersect the boundary divisor $D=\widetilde{X}\setminus\HS$
transversally in two points permuted by the (little) Weyl group.
Indeed, it suffices to verify it for generic $T$-orbits in~$\HS$,
which is easy.
\end{example}

\begin{remark}
The fibers $T^cX=\Pi^{-1}(c)$, $c\in\kk\setminus\{0\}$, are called
\emph{twisted cotangent bundles} \cite[\S2]{D(G/H)}. They carry a
structure of affine bundles over $X$ associated with the vector
bundle~$T^{*}X$. Thus each $T^cX$ has a natural symplectic
structure. (This is a particular case of symplectic reduction
\cite[App.5]{SG} for the $\kk^{\times}$-action
on~$T^{*}\widehat{X}$.) The action $G:T^cX$ is Poisson and
symplectically stable: the moment map $\Phi^c$ is the composition
of $\widehat\Phi$ and the projection $\hat\g^{*}\to\g^{*}$, so
that $\overline{\Im\Phi^c}$ is identified with the fiber of
$M_{\widehat{X}}\to\kk$ over~$c$, and the symplectic stability
stems from that of~$T^{*}\widehat{X}$.

The whole theory can be developed for arbitrary $G$-varieties
replacing the usual cotangent bundle by its twisted analogue
\cite[\S9]{inv.mot}. If $X$ is quasiaffine, i.e., embedded in a
$G$-module~$V$, then $X\subset\PP(V\oplus\kk)$, $\widehat{X}\iso
X\times\kk^{\times}$, and $T^cX\iso T^{*}X$. Therefore in the
quasiaffine case the classical theory is included in the twisted
one.
\end{remark}

\begin{proof}[Proof of Theorem~\ref{W_X&Z(X)}]
We already know from the above that $W_X$ preserves $\RG(X)$ and
acts on~$\ES$. Let $W_X^{\#}\subset\GL(\ES)$ be the subgroup
generated by reflections at the walls of~$\Zz$. The first step is
to show that $W_X^{\#}\subseteq W_X$.

Choose a wall of $\Zz$ and a vector $v$ in its interior. We may
assume $v=v_D$ for a certain $G$-stable prime divisor $D\subset
X$. Consider the normal bundle $X'$ of $X$ at~$D$. By
Proposition~\ref{V(cent)}\ref{V(N_D)}, the central valuation cone
of $X'$ is a half-space $\Zz'=\Zz+\QQ{v}$. By
Theorem~\ref{cent(hor)}, $X'$~is not horospherical, whence
$W_{X'}\ne\{\1\}$ by Proposition~\ref{W_hor}. By
Lemma~\ref{W_X:cent.tor}, $W_{X'}$~acts trivially on
$\Zz'\cap-\Zz'$, whence $W_{X'}$ is generated by the reflection at
the chosen face of~$\Zz$.

On the other hand, $X'$~is deformed to~$X$, i.e., it is the
zero-fiber of the $(G\times\kk^{\times})$-equivariant flat family
$E\to\AAA^1$ with the other fibers isomorphic to $X$
\cite[5.1]{int.theory}. Since $E\setminus X'\iso
X\times\kk^{\times}$ and the moment map of $T^{*}(E\setminus X')$
factors through the projection onto~$T^{*}X$, we have
$\widetilde{M}_X=\widetilde{M}_E$. There is a commutative diagram
%*
\begin{equation*}
\begin{array}{ccccc}
T^{*}X'&\leftarrow&T^{*}E|_{X'}&\subset&T^{*}E \\
\downarrow&&&&\downarrow \\
M_{X'}&&=&&M_E\text{\rlap{$\mathstrut=M_X$}} \\
\end{array}
\end{equation*}
%*
As $\widetilde{M}_E\to M_E$ is finite and
$\kk[\widetilde{M}_{X'}]$ is integrally closed
in~$\kk[T^{*}E|_{X'}]$, there is a finite morphism
$\widetilde{M}_{X'}\to\widetilde{M}_E$, whence
$\ab^{*}\to\widetilde{L}_{X'}\to\widetilde{L}_E= \widetilde{L}_X$.
Thus $W_{X'}\subseteq W_E=W_X$.

At this point we may reduce the problem to the symplectically
stable case, because $\widehat\Zz=\Zz(\widehat{X})$ is the
preimage of $\Zz$ and $W_{\widehat{X}}^{\#}=W_X^{\#}\subseteq
W_X\subseteq W_{\widehat{X}}$.

It follows that $W_X^{\#}$ is a finite crystallographic reflection
group and $\Zz$ is a union of its fundamental chambers. To
conclude the proof, it remains to show that different vectors
$v_1,v_2\in\Zz$ cannot be $W_X$-equivalent.

Assume the converse, i.e., $v_2=wv_1$, $w\in W_X$. W.l.o.g.\ $X$
contains two $G$-stable prime divisors $D_1,D_2$ corresponding to
$v_1,v_2$. (Replace $X$ by the normalized closure of the graph of
the birational map $X_1\dasharrow X_2$, where $X_i$ is a complete
$G$-model of $K$ having a divisor with valuation~$v_i$.) Removing
$D_1\cap D_2$, we may assume that $D_1,D_2$ are disjoint. By
Proposition~\ref{flat&div}, the closure of a generic (twisted)
flat contains two $A$-stable prime divisors $D_{v_1},D_{v_2}$
lying both in $D_1$ and in~$D_2$, a contradiction.
\end{proof}

Proposition~\ref{flat&div}, together with Theorem~\ref{W_X&Z(X)},
leads to a description of the whole boundary of a generic flat (to
a certain extent).
\begin{definition}
A \emph{source} $Y\subset X$ is the center of a central valuation.
\end{definition}
\begin{proposition}[{\cite[7.6]{inv.mot}}]\label{cl(flat)}
Let $F_{\alpha}\subseteq X$ be a generic (twisted) flat. A vector
$v\in\ES$, regarded as an $A$-valuation of~$\kk(F_{\alpha})$, has
a center $F_v$ in $\overline{F}_{\alpha}$ iff the unique $v'\in
W_Xv\cap\Zz$ has a center $Y\subseteq X$. Furthermore,
$F_v\subseteq Y$.
\end{proposition}
\begin{proof}
Since $W_X$ acts on~$\overline{F}_{\alpha}$, we may assume $v=v'$.
Take a $G$-equivariant completion $\overline{X}\supseteq X$
\cite{eq.compl} and construct a proper birational $G$-morphism
$\phi:X'\to\overline{X}$ such that $X'$ contains a divisor $D$
with $v_D=v$ (Proposition~\ref{div->center}). By
Proposition~\ref{flat&div}, the center of $v$ on the closure
$\overline{F}\mathstrut'_{\alpha}$ of $F_{\alpha}$ in $X'$ is a
divisor $D_v\subseteq D$. Hence $\phi(D_v)\subseteq\phi(D)$ is the
center of $v$ on the closure $\dbloverline{F}_{\alpha}$ of
$F_{\alpha}$ in~$\overline{X}$. It intersects
$\overline{F}_{\alpha}$ (exactly in~$F_v$) iff $\phi(D)$
intersects $X$ (in~$Y$).
\end{proof}
\begin{remark}
$A$-valuations of $\kk(F_{\alpha})$ are determined by
one-parameter subgroups of~$A$, see Example~\ref{val<->1par}. The
open orbit in $F_v$ can be reached from $F_{\alpha}$ by taking the
limits of trajectories of the respective one-parameter subgroup.
Thus Proposition~\ref{cl(flat)} gives a full picture of the
asymptotic behavior of the invariant collective motion.
\end{remark}
\begin{corollary}\label{fin.sources}
There are finitely many sources in a $G$-variety, and the closure
of a generic (twisted) flat intersects all of them.
\end{corollary}
\begin{corollary}
Suppose a quasiaffine $G$-variety $X$ contains a proper source;
then $\kk[X]$ has a $G$-invariant non-negative grading induced by
a certain central one-parameter subgroup.
\end{corollary}
\begin{proof}
We may assume $X$ to be affine. Corollary~\ref{fin.sources} and
the assumptions imply that $\overline{F}_{\alpha}\neq F_{\alpha}$
is an affine toric variety acted on by~$W_X$. Its normalization is
determined by a strictly convex $W_X$-stable cone $\Cc\subset\ES$
(Example~\ref{toric}). Clearly, $\intr\Cc$~contains a
$W_X$-invariant vector $v\neq0$. Hence $v\in\Zz\cap-\Zz$ defines a
central one-parameter subgroup $\gamma$ acting on $X$ by
Theorem~\ref{cent.aut}, whence a $G$-invariant grading
of~$\kk[X]$. Generic $\gamma$-orbits in $X$ are contained in
generic flats and non-closed therein, because $\gamma$ contracts
$F_{\alpha}$ to the unique closed orbit
in~$\overline{F}_{\alpha}$. Hence the grading is non-negative.
\end{proof}

See \cite[\S\S8--9]{inv.mot} and \ref{toroidal} for a deeper
analysis of sources, flats, and their closures.

To any $G$-variety $X$ one can relate a root system
$\Delta_X\subset\RG(X)$, which is a birational invariant of the
$G$-action. Namely, let $\Pi_X^{\min}\subset\RG(X)$ be the set of
indivisible vectors generating the rays of the simplicial
cone~$-\Zz(X)^{\vee}$.  It is easy to deduce from
Theorem~\ref{W_X&Z(X)} that $\Delta_X^{\min}=W_X\Pi_X^{\min}$ is a
root system with base $\Pi_X^{\min}$ and the Weyl group~$W_X$,
called the \emph{minimal root system} of~$X$. It is a
generalization of the (reduced) root system of a symmetric variety
(see~\ref{symmetric} and Example~\ref{wonder(symm)}). The minimal
root system was defined by Brion \cite{W_X(sph)} for a spherical
variety and by Knop \cite{Aut&root} in the general case.
\begin{remark}
There are several natural root systems related to $G:X$ which
generate one and the same Weyl group~$W_X$ \cite[6.2, 6.4,
7.5]{Aut&root}; $\Delta_X^{\min}$~is the ``minimal'' one.
\begin{question}
Introduce also $\Delta^{\max}_X$.
\end{question}
\end{remark}
\begin{example}
If $X=G$ comes equipped with the $G$-action by left translations
and $G'$ is adjoint, then $\Delta_X^{\min}=\Delta_G$ by
Example~\ref{cent.val(G)}.  If $G'$ is not adjoint, then
$\Delta_X^{\min}$ may differ from~$\Delta_G$: this happens iff
some roots in $\Delta_G$ are divisible in~$\Ch(T)$. For
simple~$G$, $\Delta_G^{\min}=\Delta_G$ unless $G=\Sp_{2n}(\kk)$;
in the latter case $\Delta_G,\Delta_G^{\min}$ are of types
$\CCc_n,\BBb_n$, respectively.
\end{example}

An important result of Knop establishes a relation between the
minimal root system and central automorphisms.
\begin{theorem}[{\cite[6.4]{Aut&root}}]\label{root->cent}
A quasitorus
$S_X=\bigcap_{\alpha\in\Delta_X^{\min}}\Ker\alpha\subseteq A$ is
canonically embedded in~$\CA{X}$.
\begin{question}
Describe the whole $\CA{X}$ using $\Delta_X^{\max}$
\cite[6.5]{Aut&root}.
\end{question}
\end{theorem}
Note that $S_X^0$ is the largest connected algebraic subgroup of
$\CA{X}$ by Theorem~\ref{cent.aut}\ref{S<A}, but
Theorem~\ref{root->cent} is much more subtle.
\begin{proof}[Synopsis of a proof]
Standard reductions allow us to assume that $X$ is quasiaffine. It
is clear that $S_X\subseteq A^{W_X}$. The action
$A:\widehat{T}^{\pr}X$ descends to $A^{W_X}:T^{\pr}X$.

The most delicate part of the proof is to show that the action of
$S_X$ extends to $T^{*}X$ in codimension one. Knop shows that the
$A$-actions on the orbits of the invariant collective motion patch
together in an action on $T^{*}X$ of a smooth group scheme $\Ss$
over $\widetilde{L}_X$ with connected fibers, and
$A^{W_X}\subset\Ss(\kk(\widetilde{L}_X))$. Furthermore, $s\in
A^{W_X}$ induces a rational section of $\Ss\to\widetilde{L}_X$
which is defined in codimension one whenever $\alpha(s)=1$,
$\forall\alpha\in\Delta_X^{\min}$, whence the claim.

Now $S_X$ acts on an open subset $R\subseteq T^{*}X$ whose
complement has codimension $\geq2$, and this action commutes with
$G$ and with homotheties on the fibers. Hence $S_X$ acts by
$G$-automorphisms on $\PP(R)\subseteq\PP(T^{*}X)$. Since $X$ is
quasiaffine and generic fibers of $\PP(R)\to X$ have no
non-constant regular functions, we deduce that $S_X$ permutes the
fibers. This yields a birational action $S_X:X$ commuting with $G$
and preserving generic flats. The description of generic flats
shows that $S_X$ preserves $P$-orbits in $\X$, whence
$S_X\embeds\CA{X}$.
\end{proof}

\section{Formal curves}
\label{form.curv}

In the previous sections of this chapter we examined
$G$-valuations on arbitrary $G$-varieties. However our main
interest is in homogeneous varieties. In this section we take a
closer look at $\Vv=\Vv(\HS)$, where $\HS$ is a homogeneous space.

Namely we describe the subset $\Vv^1\subseteq\Vv$ consisting of
$G$-valuations $v$ such that $\kk(v)^G=\kk$. In geometric terms,
if $v$ is proportional to $v_D$ for a $G$-stable divisor $D$ on a
$G$-embedding $X\embof\HS$, then $v\in\Vv^1$ iff $D$ contains a
dense $G$-orbit.

The subset $\Vv^1$ is big enough. For instance, if $c(\HS)=0$,
then $\Vv^1=\Vv$, and if $c(\HS)=1$, then
$\Vv^1\supseteq\Vv\setminus\Zz$, because in this case
$c(\kk(v))=0$ by Proposition~\ref{c&r(cent)}. In general, any
$G$-valuation can be approximated by $v\in\Vv^1$ in a sense
\cite[4.11]{LV}.

In \cite{LV} Luna and Vust suggested to compute $v(f)$,
$v\in\Vv^1$, $f\in K=\kk(\HS)$, by restricting $f$ to a (formal)
curve in $\HS$ approaching to~$D$, in the above notation. More
precisely, take a smooth curve $\Theta\subseteq X$ meeting $D$
transversally in~$x_0$, $\overline{Gx_0}=D$. It is clear that
$v_D(f)$ equals the order of $f|_{g\Theta}$ at $gx_0$ for general
$g\in G$. More generally, take a germ of a curve
$\chi:\Theta\dasharrow\HS$ that converges to $x_0$ in~$X$, i.e.,
$\chi$~extends regularly to the base point $\theta_0\in\Theta$ and
$\chi(\theta_0)=x_0$, see Appendix~\ref{sch.pt}. Then
%*
\begin{equation}\label{curv->val}
v_D(f)\cdot\langle D,\Theta\rangle_{x_0}=
v_{\chi,\theta_0}(f):=v_{\theta_0}(\chi^{*}(gf)) \qquad\text{for
general $g\in G$,}
\end{equation}
%*
where $\langle D,\Theta\rangle_{x_0}$ is the local intersection
number \cite[Ch.7]{int.theory}.
\begin{theorem}\label{val<->curv}
For any germ of a curve
$(\chi:\Theta\dasharrow\HS,\theta_0\in\Theta)$,
Formula~\eqref{curv->val} defines a $G$-valuation
$v_{\chi,\theta_0}\in\Vv^1$, and every $v\in\Vv^1$ is proportional
to some~$v_{\chi,\theta_0}$. Furthermore, if $X\supseteq\HS$ is a
$G$-model of $K$ and $Y\subseteq X$ the center of~$v$, then the
germ converges in $X$ to $x_0\in Y$ such that $\overline{Gx_0}=Y$.
\end{theorem}
\begin{proof}
The $G$-action yields a rational dominant map
$\alpha:G\times\Theta\dasharrow\HS$, $(g,\theta)\mapsto
g\chi(\theta)$. By construction, $v_{\chi,\theta_0}$ is the
restriction of $v_{G\times\{\theta_0\}}\in\Vv^1(G\times\Theta)$
to~$K$, whence $v=v_{\chi,\theta_0}\in\Vv^1$. If $v$ has the
center $Y$ on~$X$, then $\alpha:G\times\Theta\dasharrow X$ is
regular along $G\times\{\theta_0\}$ and
$\overline{\alpha(G\times\{\theta_0\})}=Y$, whence $\chi$
converges to $x_0=\alpha(e,\theta_0)$ in the dense $G$-orbit
of~$Y$.
\end{proof}

For computations, it is more practical to adopt a more algebraic
point of view, namely to replace germs of curves by germs of
formal curves, i.e., by $\kk\(t\)$-points of~$\HS$, see
Appendix~\ref{sch.pt}.

Any germ of a curve $(\chi:\Theta\dasharrow\HS,\theta_0\in\Theta)$
defines a formal germ $x(t)\in\HS(\kk\(t\))$ if we replace
$\Theta$ by the formal neighborhood of~$\theta_0$. We have
%*
\begin{equation}\label{curv->form}
v_{\chi,\theta_0}(f)=v_{x(t)}(f):=\ord_tf(gx(t)) \qquad\text{for
generic~$g$,}
\end{equation}
%*
where ``generic'' means a sufficiently general point of~$G$
(depending on $f\in K$) or the generic $\kk(G)$-point of~$G$
(Example~\ref{gen.pt}).

The counterpart of Theorem~\ref{val<->curv} is
\begin{theorem}
For any $x(t)\in\HS(\kk\(t\))$, Formula~\eqref{curv->form} defines
a $G$-valuation $v_{x(t)}\in\Vv^1$, and every $v\in\Vv^1$ is
proportional to some~$v_{x(t)}$. Furthermore, if $X\supseteq\HS$
is a $G$-model of $K$ and $Y\subseteq X$ the center of~$v$, then
$x(t)\in X(\kk\[t\])$ and $Y=\overline{Gx(0)}$.
\end{theorem}

To prove this theorem it suffices to show that
$v_{x(t)}=v_{\chi,\theta_0}$ for a certain germ of a
curve~$(\chi,\theta_0)$. This stems from the two subsequent
lemmas.

\begin{lemma}[{\cite[4.4]{LV}}]
$\forall g(t)\in G(\kk\[t\]),\ x(t)\in\HS(\kk\(t\)):\
{v_{g(t)x(t)}=v_{x(t)}}$
\end{lemma}
\begin{proof}
The $G$-action on $x(t)$ yields $\kk(\HS)\embeds\kk(G)\(t\)$, so
that $v_{x(t)}$ coincides with $\ord_t$ w.r.t.\ this inclusion.
The lemma stems from the fact that $G(\kk\[t\])$ acts on
$\kk(G)\(t\)$ ``by right translations'' preserving~$\ord_t$.
\end{proof}

\begin{lemma}[{\cite[4.5]{LV}}]
Every germ of a formal curve in $\HS$ is $G(\kk\[t\])$-equivalent
to a formal germ induced by a germ of a curve.
\end{lemma}
\begin{proof}
Since $\HS$ is homogeneous, $G(\kk\[t\])$-orbits are open in
$\HS(\kk\(t\))$ in $t$-adic topology~\cite{adic}.
\begin{question}
Reduce to $\HS=G$ and prove directly.
\end{question}
Now the lemma stems from Theorem~\ref{approx(form)}.
\end{proof}

Germs of formal curves in $\HS$ are more accessible if they come
from formal germs in~$G$. Luckily, this is ``almost'' always the
case.
\begin{proposition}[{\cite[4.3]{LV}}]
For $\forall x(t)\in\HS(\kk\(t\))$ there exists $n\in\NN$ such
that $x(t^n)=g(t)\cdot\bp$ for some $g(t)\in G(\kk\(t\))$.
\end{proposition}
\begin{proof}
Consider the algebraic closure $\overline{\kk\(t\)}$
of~$\kk\(t\)$. The set $\HS(\overline{\kk\(t\)})$ is equipped with
a structure of an algebraic variety over $\overline{\kk\(t\)}$
with the transitive $G(\overline{\kk\(t\)})$-action. However
$\overline{\kk\(t\)}=\bigcup_{n=1}^{\infty}\kk\(\sqrt[n]t\)$,
whence $x(t)=g(\sqrt[n]t)\cdot\bp$ for some $g(\sqrt[n]t)\in
G(\kk\(\sqrt[n]t\))$.
\end{proof}

Note that $v_{x(t^n)}=n\cdot v_{x(t)}$. Thus we may describe
$\Vv^1$ in terms of germs of formal curves in~$G$, i.e., points
of~$G(\kk\(t\))$, considered up to left translations by
$G(\kk\[t\])$ and right translations by~$H(\kk\(t\))$. There is a
useful structural result shrinking the set of formal germs under
consideration:
\begin{Claim}[Iwasawa decomposition~\cite{Iwahori}]
$G(\kk\(t\))=G(\kk\[t\])\cdot\CoCh(T)\cdot U(\kk\(t\))$, where
$\CoCh(T)$ is regarded as a subset of~$T(\kk\(t\))$.
\end{Claim}
\begin{corollary}
Every $v\in\Vv^1$ is proportional to (the restriction
of)~$v_{g(t)}$, $g(t)\in\CoCh(T)\cdot U(\kk\(t\))$.
\end{corollary}
Let us mention a related useful result on the structure of
$G(\kk\(t\))$:
\begin{Claim}[Cartan decomposition~\cite{Iwahori}]
$G(\kk\(t\))=G(\kk\[t\])\cdot\CoCh(T)\cdot G(\kk\[t\])$.
\end{Claim}

\begin{example}
Suppose that $\HS=G/S$ is horospherical. We may assume $S\supseteq
U$; then $N(S)\supseteq B$ and $A:=\Aut_G\HS\iso N(S)/S=T/T\cap S$
is a torus. Since $\HS$ is spherical, $\Vv=\Vv^1$. Due to the
Iwasawa decomposition, every $v\in\Vv$ is proportional to
some~$v_{\gamma}$, $\gamma\in\CoCh(T)$. Let $\overline\gamma$ be
the image of $\gamma$ in~$\CoCh(A)$. By definition,
$v_{\gamma}(f)=\ord_{t=0}f(g\gamma(t)\bp)=
\ord_{t=0}f(\overline\gamma(t)\cdot g\bp)$ is the order of $f$
along generic trajectories of $\overline\gamma$ as $t\to0$. In
particular, $\Vv=\CoCh(A)\otimes\QQ$, cf.~Theorem~\ref{cent(hor)}.
\end{example}
\begin{example}\label{val<->1par}
Specifically, let $\HS=G=T$ be a torus. Every $T$-valuation of
$\kk(T)$ is proportional to~$v_{\gamma}$, $\gamma\in\CoCh(T)$,
where $v_{\gamma}$ is the order of a function restricted to
$s\gamma(t)$ as $t\to0$ for general $s\in T$. By
Theorem~\ref{val<->curv}, $v_{\gamma}$~has a center $Y$ on a toric
variety $X\supseteq T$ iff $\gamma(0):=\lim_{t\to0}\gamma(t)$
exists and belongs to the dense $T$-orbit in~$Y$. Thus the lattice
points in $\Vv_Y$ are exactly the one-parameter subgroups of $T$
converging to a point in the dense $T$-orbit of~$Y$. This is the
classical description of the fan of a toric variety
\cite{toric.book}, \cite[2.3]{toric.intro}.
\end{example}
Other examples can be found in~\ref{c=1}.

\chapter{Spherical varieties}
\label{spherical}

Although the theory developed in the previous chapters applies to
arbitrary homogeneous spaces of reductive groups, and even to more
general group actions, it acquires most complete and elegant form
for spherical homogeneous spaces and their equivariant embeddings,
called spherical varieties. A justification of the fact that
spherical homogeneous spaces are a significant mathematical object
is that they naturally arise in various fields, such as embedding
theory, representation theory, symplectic geometry, etc. In
\ref{sphericity} we collect various characterizations of spherical
spaces, the most important being: the existence of an open
$B$-orbit, the ``multiplicity free'' property for spaces of global
sections of line bundles, commutativity of invariant differential
operators and of invariant functions on the cotangent bundle
w.r.t.\ the Poisson bracket.

Then we examine most interesting classes of spherical homogeneous
spaces and spherical varieties in more details. Algebraic
symmetric spaces are considered in~\ref{symmetric}. We develop the
structure theory and classification of symmetric spaces, compute
the colored data required for description of their equivariant
embeddings, study $B$-orbits and (co)isotropy representation.
\ref{monoids}~is devoted to $(G\times G)$-equivariant embeddings
of a reductive group~$G$. A particular interest in this class is
explained, for example, by an observation that linear algebraic
monoids are nothing else but affine equivariant group embeddings.
Horospherical varieties of complexity~$0$ are classified and
studied in~\ref{S-var}.

Geometric structure of toroidal varieties, considered
in~\ref{toroidal}, is best understood among all spherical
varieties, since toroidal varieties are ``locally toric''. They
can be defined by several equivalent properties: their fans are
``colorless'', they are spherical and pseudo-free, the action
sheaf on a toroidal variety is the log-tangent sheaf w.r.t.\ a
$G$-stable divisor with normal crossings. An important property of
toroidal varieties is that they are rigid as $G$-varieties. The
so-called wonderful varieties are the most remarkable subclass of
toroidal varieties. They are canonical completions with nice
geometric properties of (certain) spherical homogeneous spaces.
The theory of wonderful varieties is developed in~\ref{wonderful}.
Applications include computation of the canonical divisor of a
spherical variety and Luna's conceptual approach to the
classification of spherical subgroups through the classification
of wonderful varieties.

The last \ref{Frob.split} is devoted to Frobenius splitting, a
technique for proving geometric and algebraic properties
(normality, rationality of singularities, cohomology vanishing,
etc) in positive characteristic. However, this technique can be
applied to zero characteristic using reduction ${\!}\bmod p$
provided that reduced varieties are Frobenius split. This works
for spherical varieties. As a consequence, one obtains vanishing
of higher cohomology of ample or nef line bundles on complete
spherical varieties, normality and rationality of singularities
for $G$-stable subvarieties, etc. Some of these results can be
proved by other methods, but Frobenius splitting provides a simple
uniform approach.

\section{Various characterizations of sphericity}
\label{sphericity}

Spherical homogeneous spaces can be considered from diverse
viewpoints: orbits and equivariant embeddings, representation
theory and multiplicities, symplectic geometry, harmonic analysis,
etc. The definition and some other implicit characterizations of
this remarkable class of homogeneous spaces are already scattered
in the text above. In this section, we review these issues and
introduce other important properties of homogeneous spaces which
are equivalent or closely related to sphericity.

As usual, $G$~is a connected reductive group, $\HS$~denotes a
homogeneous $G$-space with the base point~$\bp$, and $H=G_{\bp}$.
\begin{Claim}[Definition--Theorem]
A \emph{spherical homogeneous space}~$\HS$ (resp.\ a
\emph{spherical subgroup} $H\subseteq G$, a \emph{spherical
subalgebra} $\h\subseteq\g$, a \emph{spherical pair} $(G,H)$ or
$(\g,\h)$) can be defined by any one of the following properties:
\begin{enumerate}
\renewcommand{\theenumi}{\reftag{S\arabic{enumi}}}
\item\label{B-inv} $\kk(\HS)^B=\kk$. \item\label{open-B} $B$~has
an open orbit in~$\HS$. \item\label{open-H} $H$~has an open orbit
in~$G/B$. \item\label{b+gh=g} $\exists g\in G:\
\br+(\Ad{g})\h=\g$. \item\label{h+gb=g} There exists a Borel
subalgebra $\widetilde\br\subseteq\g$ such that
$\h+\widetilde\br=\g$. \item\label{fin-H} $H$~acts on $G/B$ with
finitely many orbits. \item\label{fin-G} For any $G$-variety $X$
and $\forall x\in X^H$, $\overline{Gx}$~contains finitely many
$G$-orbits. \item\label{fin-B} For any $G$-variety $X$ and
$\forall x\in X^H$, $\overline{Gx}$~contains finitely many
$B$-orbits.
\end{enumerate}
\end{Claim}
The term ``spherical homogeneous space'' is traced back to Brion,
Luna, and Vust \cite{BLV}, and ``spherical subgroup''
to~\cite{sph(s)}, though the notions themselves appeared much
earlier.
\begin{proof}
\begin{roster}
\item[\ref{B-inv}$\iff$\ref{open-B}] $B$-invariant functions
separate generic $B$-orbits \cite[2.3]{IT}.
\item[\ref{open-B}$\iff$\ref{open-H}] Both conditions are
equivalent to that $B\times H:G$ has an open orbit, where $B$ acts
by left and $H$ by right translations, or vice versa.
\item[\ref{b+gh=g} and \ref{h+gb=g}] are just reformulations of
\ref{open-B} and \ref{open-H} in terms of tangent spaces.
\item[\ref{open-B}$\implies$\ref{fin-B}]
$Gx$~satisfies~\ref{open-B}, too, and we conclude by
Corollary~\ref{B-fin}. \item[\ref{fin-B}$\implies$\ref{fin-G}]
Obvious. \item[\ref{fin-G}$\implies$\ref{open-B}] Stems from
Corollary~\ref{G-fin}. \item[\ref{fin-B}$\implies$\ref{fin-H}]
$B$~acts on $G/H$ with finitely many orbits, which are in
bijection with $(B\times H)$-orbits on $G$ and with $H$-orbits
on~$G/B$. \item[\ref{fin-H}$\implies$\ref{open-H}] Obvious.
\qedhere\end{roster}
\end{proof}
In particular, spherical spaces are characterized in the framework
of embedding theory as those having finitely many orbits in the
boundary of any equivariant embedding. The embedding theory of
spherical spaces is considered in~\ref{c=0}.

Another important characterization of spherical spaces is in terms
of representation theory, due to Kimelfeld and
Vinberg~\cite{sph&mult}. Recall from \ref{bundles} that the
multiplicity of a highest weight $\lambda$ in a $G$-module $M$ is
%*
\begin{equation*}
m_{\lambda}(M)=\dim\Hom_G(V(\lambda),M)=\dim M^{(B)}_{\lambda}
\end{equation*}
%*
In characteristic zero, $m_{\lambda}(M)$~is the multiplicity of
the simple $G$-module $V(\lambda)$ in the decomposition of~$M$. In
positive characteristic, $V(\lambda)$~denotes the respective Weyl
module. The module $M$ is said to be multiplicity free if all
multiplicities in $M$ are~$\leq1$.
\begin{theorem}
$\HS$~is spherical iff the following equivalent conditions hold:
\begin{roster}
\renewcommand{\theenumi}{\reftag{MF\arabic{enumi}}}
\item\label{P(V)^H} $\PP(V(\lambda))^H$~is finite
for~$\forall\lambda\in\Ch_{+}$. \item\label{V^(H)}
$\forall\lambda\in\Ch_{+},\ \chi\in\Ch(H):\ \dim
V(\lambda)^{(H)}_{\chi}\leq1$ \item\label{mult-free(L)} For any
$G$-line bundle $\Ll$ on~$\HS$, $\Ho^0(\HS,\Ll)$~is multiplicity
free.
\end{roster}
If $\HS$ is quasiaffine, then the last two conditions can be
weakened to
\begin{roster}
\renewcommand{\theenumi}{\reftag{MF\arabic{enumi}}}
\setcounter{enumi}{3} \item\label{V^H} $\forall\lambda\in\Ch_{+}:\
\dim V(\lambda)^H\leq1$ \item\label{mult-free(k[O])} $\kk[\HS]$~is
multiplicity free.
\end{roster}
\end{theorem}
The spaces satisfying these conditions are called
\emph{multiplicity free}.
\begin{proof}
\begin{roster}
\item[\ref{B-inv}$\iff$\ref{mult-free(L)}] If
$m_{\lambda}(\Ll)\geq2$, then there exist two non-proportional
sections $\sigma_0,\sigma_1\in\Ho^0(\HS,\Ll)^{(B)}_{\lambda}$.
Their ratio $f=\sigma_1/\sigma_0$ is a non-constant $B$-invariant
function. Conversely, any $f\in\kk(\HS)^B$ can be represented in
this way: the $G$-line bundle $\Ll$ together with the canonical
$B$-eigensection $\sigma_0$ is defined by a sufficiently big
multiple of $\divr_{\infty}f$ (cf.~Corollary~\ref{power-lin}).

Finally, if $\HS$ is quasiaffine, then we may take for $\Ll$ the
trivial bundle: for $\sigma_0$ take a sufficiently big power of
any $B$-eigenfunction in $\Ideal{D}\normin\kk[\HS]$,
$D=\Supp\divr_{\infty}f\subset\HS$. Hence
\ref{B-inv}$\iff$\ref{mult-free(k[O])}.

\item[\ref{P(V)^H}$\iff$\ref{V^(H)}] Stems from
$\PP(V(\lambda))^H=\PP\left(V(\lambda)^{(H)}\right)=
\bigsqcup_{\chi}\PP\left(V(\lambda)^{(H)}_{\chi}\right)$ (a finite
disjoint union).

\item[\ref{V^(H)}$\iff$\ref{mult-free(L)}] If $\HS=G/H$ is a
quotient space, then this is the Frobenius
reciprocity~\eqref{mult(G/H)}. Generally, there is a bijective
purely inseparable morphism $G/H\to\HS$ (Remark~\ref{insep}), and
$\HS$ is spherical iff $G/H$ is so. But we have already seen that
the sphericity is equivalent to~\ref{mult-free(L)}.

\item[\ref{V^H}$\iff$\ref{mult-free(k[O])}] is proved in the same
way. \qedhere\end{roster}
\end{proof}

The ``multiplicity free'' property leads to an interpretation of
sphericity in terms of automorphisms and group algebras associated
with~$G$. Since the complete reducibility of rational
representations is essential here, we assume $\ch\kk=0$ up to the
end of this section.

Recall from \ref{bundles} the algebraic versions of the group
algebra $\Aa(G)$ and the Hecke algebra $\Aa(\HS)$.

\begin{theorem}[\cite{weak.symm}, \cite{comm.hom}]
An affine homogeneous space $\HS$ is spherical iff either of the
following equivalent conditions is satisfied:
\begin{roster}
\renewcommand{\theenumi}{\reftag{GP\arabic{enumi}}}
\setcounter{enumi}{0} \item\label{Hecke} $\Aa(\HS)=\Aa(G)^{H\times
H}$ is commutative. \item\label{Hecke.fin} $\Aa(V)^{H\times H}$ is
commutative for all $G$-modules~$V$.
\end{roster}
\begin{roster}
\renewcommand{\theenumi}{\reftag{WS\arabic{enumi}}}
\item\label{Selberg} (Selberg condition) The $G$-action on $\HS$
extends to a cyclic extension $\widehat{G}=\langle G,s\rangle$ of
$G$ so that $(sx,sy)$ is $G$-equivalent to $(y,x)$ for general
$x,y\in\HS$. \item\label{Gelfand} (Gelfand condition) There exists
$\sigma\in\Aut{G}$ such that $\sigma(g)\in Hg^{-1}H$ for general
$g\in G$.
\end{roster}
\end{theorem}
The condition \ref{Hecke} is an algebraization of a similar
commutativity condition for the group algebra of a Lie group, see
\cite{symm=>MF}, \cite[I.2]{comm.hom}, and below. The condition
\ref{Gelfand} appeared in \cite{symm=>MF} and \ref{Selberg} was
first introduced by Selberg in the seminal paper on the trace
formula \cite{trace}, and by Akhiezer and Vinberg \cite{weak.symm}
in the context of algebraic geometry. The spaces satisfying
\ref{Selberg}--\ref{Gelfand} are called \emph{weakly symmetric}
and $(G,H)$ is said to be a \emph{Gelfand pair} if
\ref{Hecke}--\ref{Hecke.fin} hold.
\begin{proof}
\begin{roster}
\item[\ref{mult-free(k[O])}$\iff$\ref{Hecke}] Stems from Schur's
lemma.

\item[\ref{Hecke}$\iff$\ref{Hecke.fin}] Obvious.

\item[\ref{mult-free(k[O])}$\implies$\ref{Selberg}] There exists a
Weyl involution $\sigma\in\Aut{G}$, $\sigma(H)=H$
\cite{weak.symm}. There is a conceptual argument for symmetric
spaces and in general a case-by-case verification using the
classification from~\ref{c&r<=1}. Define $s\in\Aut\HS$ by
$s(g\bp)=\sigma(g)\bp$ and $\widehat{G}=G\timessemi\langle
s\rangle$ by $sgs^{-1}=\sigma(g)$. The $G$-action on
$\HS\times\HS$ is extended to $\widehat{G}$ by $s(x,y)=(sy,sx)$.

Consider the $(G\times G)$-isotypic decomposition
%*
\begin{gather*}
\kk[\HS\times\HS]=\bigoplus_{\lambda,\mu\in\RG_{+}(\HS)}
\kk[\HS\times\HS]_{(\lambda,\mu)},\qquad\text{where}\\
\kk[\HS\times\HS]_{(\lambda,\mu)}=
\kk[\HS]_{(\lambda)}\otimes\kk[\HS]_{(\mu)}\iso V(\lambda)\otimes
V(\mu)
\end{gather*}
%*
Clearly, $s$~twists the $G$-action by~$\sigma$, hence maps
$\kk[\HS\times\HS]_{(\lambda,\mu)}$ to
$\kk[\HS\times\HS]_{(\mu^{*},\lambda^{*})}$ and preserves each
summand of
%*
\begin{equation*}
\kk[\HS\times\HS]^G=\bigoplus_{\lambda\in\RG_{+}(\HS)}
\kk[\HS\times\HS]^G_{(\lambda,\lambda^{*})}
\end{equation*}
%*
However a $\widehat{G}$-invariant inner product
$(p,q)\to\average{(pq)}$ on $\kk[\HS]$ induces a nonzero pairing
between simple $G$-modules $\kk[\HS]_{(\lambda^{*})}$ and
$\kk[\HS]_{(\lambda)}$, whence by duality a
$\widehat{G}$-invariant function
in~$\kk[\HS\times\HS]_{(\lambda,\lambda^{*})}$, which spans
$\kk[\HS\times\HS]^G_{(\lambda,\lambda^{*})}$. It follows that $s$
acts trivially on~$\kk[\HS\times\HS]^G$.

But the action $G:\HS\times\HS$ is stable
(Theorem~\ref{stab(dbl)}), whence $s$ preserves generic $G$-orbits
in~$\HS\times\HS$, which is exactly the Selberg condition.

\item[\ref{Selberg}$\implies$\ref{Gelfand}] Multiplying $s$ by
$g\in G$ preserves the Selberg condition. Also, if
$(sx,sy)\sim(y,x)$, then the same is true for any $G$-equivalent
pair. Hence, w.l.o.g., $s\bp=\bp=x$. Define $\sigma\in\Aut{G}$ by
$\sigma(g)=sgs^{-1}$; then
$(s\bp,sg\bp)=(\bp,\sigma(g)\bp)\sim(g\bp,\bp)$ for general $g\in
G$. Hence $g'g\bp=\bp$, $g'\bp=\sigma(g)\bp$, i.e., $g'g=h\in H$,
$\sigma(g)=g'h'=hg^{-1}h'$ for some $h'\in H$.

\item[\ref{Gelfand}$\implies$\ref{Selberg}] The Gelfand condition
implies $\sigma(H)=H$, whence $s\in\Aut\HS$,
$s(g\bp)=\sigma(g)\bp$, is a well-defined automorphism. Put
$\widehat{G}=G\timessemi\langle s\rangle$, $sgs^{-1}=\sigma(g)$.
The Selberg condition is verified by reversing the previous
arguments.

\item[\ref{Gelfand}$\implies$\ref{Hecke}] The inversion map
$g\mapsto g^{-1}$ on $G$ extends to an involutive antiautomorphism
of~$\Aa(G)$. Its restriction to $\Aa(G)^{H\times H}$ coincides
with the automorphism induced by~$\sigma$. Hence $\Aa(G)^{H\times
H}$ is commutative. \qedhere\end{roster}
\end{proof}
\begin{remark}
Already in the quasiaffine case the classes of weakly symmetric
and spherical spaces are not contained in each other~\cite{Zor}.
\begin{question}
Provide an example?
\end{question}
\end{remark}

Now we characterize sphericity in terms of differential geometry.
%Since infinitesimal concepts behave well only in characteristic
%zero, we assume $\ch\kk=0$ up to the end of this section.

Recall from \ref{cotangent} that the action $G:T^{*}\HS$ is
Poisson w.r.t.\ the natural symplectic structure. Thus we have a
$G$-invariant Poisson bracket of functions on~$T^{*}\HS$.
Homogeneous functions on $T^{*}\HS$ are locally the symbols of
differential operators on~$\HS$, and the Poisson bracket is
induced by the commutator of differential operators.

The functions pulled back under the moment map
$\Phi:T^{*}\HS\to\g^{*}$ are called \emph{collective}. They
Poisson-commute with $G$-invariant functions on~$T^{*}\HS$
(Proposition~\ref{coll&inv}).

\begin{theorem}\label{WC&C}
$\HS$~is spherical iff the following equivalent conditions hold:
\begin{roster}
\renewcommand{\theenumi}{\reftag{WC\arabic{enumi}}}
\item\label{coiso} Generic orbits of $G:T^{*}\HS$ are coisotropic,
i.e., $\g\alpha\supseteq(\g\alpha)^{\sort}$ for general $\alpha\in
T^{*}\HS$. \item\label{weak.comm} $\kk(T^{*}\HS)^G$~is commutative
w.r.t.\ the Poisson bracket.
\renewcommand{\theenumi}{\reftag{CI}}
\item\label{coll.int} There exists a complete system of collective
functions in involution on~$T^{*}\HS$.
\end{roster}
If $\HS$ is affine, then these conditions are equivalent to
\begin{roster}
\renewcommand{\theenumi}{\reftag{Com}}
\item\label{comm} $\Dd(\HS)^G$~is commutative.
\end{roster}
\end{theorem}
The theorem goes back to Guillemin, Sternberg \cite{mult-free},
and Mikityuk \cite{sph(ss)}. The spaces satisfying
\ref{coiso}--\ref{weak.comm} are called \emph{weakly commutative}
and those satisfying \ref{comm} are said to be \emph{commutative}.
\begin{proof}
\begin{roster}
\item[\ref{open-B}$\iff$\ref{coiso}] By Theorem~\ref{c&r<=T^*},
$\cork T^{*}\HS=2c(\HS)$ is zero iff $\HS$ is spherical, and this
means exactly that generic orbits are coisotropic.

\item[\ref{coiso}$\iff$\ref{weak.comm}] Skew gradients of
$f\in\kk(T^{*}\HS)^G$ at a point $\alpha$ of general position span
$(\g\alpha)^{\sort}$. All $G$-invariant function Poisson-commute
iff their skew gradients are skew-orthogonal to each other, i.e.,
iff $(\g\alpha)^{\sort}$ is isotropic.

\item[\ref{Hecke}$\iff$\ref{comm}] If $\HS$ is quasiaffine, then
$\Dd(\HS)$ acts faithfully on $\kk[\HS]$ by linear endomorphisms.
Hence $\Dd(\HS)^G$ is a subalgebra in~$\Aa(\HS)$. It remains to
utilize the approximation of linear endomorphisms by differential
operators.

\begin{lemma}\label{end<-diff}
Let $X$ be an smooth affine $G$-variety.
\begin{roster}
\item\label{A<-D} For any linear operator $\phi:\kk[X]\to\kk[X]$
and any finite-dimensional subspace $M\subset\kk[X]$ there exists
$\partial\in\Dd(X)$ such that $\partial|_M=\phi|_M$.
\item\label{inv.D} If $\phi$ is $G$-equivariant, then one may
assume $\partial\in\Dd(X)^G$. \item\label{dbl.ann} Put $\Ii=\Ann
M\normin\Dd(X)$; then $\forall f\in\kk[X]:\Ii{f}=0\implies f\in
M$.
\end{roster}
\end{lemma}

We conclude by Lemma~\ref{end<-diff}\ref{inv.D} that $\Aa(\HS)$ is
commutative iff $\Dd(\HS)^G$ is so.

\begin{proof}[Proof of Lemma~\ref{end<-diff}]
\begin{roster}
\item[\ref{A<-D}] We deduce it from~\ref{dbl.ann}. Choose a basis
$f_1,\dots,f_n$ of~$M$. It suffices to construct
$\partial\in\Dd(X)$ such that $\partial f_i=0,\ \forall i<n,\
\partial f_n=1$. By \ref{dbl.ann} there exists
$\partial'\in\Ann(f_1,\dots,f_{n-1})$, $\partial'f_n\neq0$. As
$\kk[X]$ is a simple $\Dd(X)$-module \cite{NCA}
\begin{question}
Is it true for singular affine $X$?
\end{question}
we may find $\partial''\in\Dd(X)$, $\partial''(\partial'f_n)=1$
and put $\partial=\partial''\partial'$. \item[\ref{inv.D}]
W.l.o.g.~$M$ is $G$-stable. Assertion~\ref{A<-D} yields an
epimorphism of $G$-$\kk[X]$-modules $\Dd(X)\onto\Hom(M,\kk[X])$
given by restriction to~$M$. But taking $G$-invariants is an exact
functor.

\item[\ref{dbl.ann}] The assertion is trivial for $M=0$ and we
proceed by induction on $\dim M$. In the above notation, put
$\Ii'=\Ann(f_1,\dots,f_{n-1})$. For
$\forall\partial,\partial'\in\Ii'$ we have $(\partial
f_n)\partial'-(\partial'f_n)\partial\in\Ii$, whence
%*
\begin{equation}\label{dd'}
(\partial f_n)(\partial'f)=(\partial'f_n)(\partial f)
\end{equation}
%*
Taking $\partial'=\xi\partial$, $\xi\in\Vect{X}$, yields
$\xi(\partial f/\partial f_n)=0\implies\partial f/\partial
f_n=c_{\partial}=\const$. Substituting this in \eqref{dd'} yields
$c_{\partial'}=c_{\partial}=c$ (independent of~$\partial$). Thus
$\partial(f-cf_n)=0\implies f-cf_n\in\langle
f_1,\dots,f_{n-1}\rangle\implies f\in M$. \qedhere\end{roster}
\end{proof}

\item[\ref{comm}$\implies$\ref{weak.comm}] If $\HS$ is affine,
then $\gr\Dd(\HS)=\kk[T^{*}\HS]$. By complete reducibility,
$\kk[T^{*}\HS]^G=\gr\Dd(\HS)^G$ is Poisson-commutative. But the
$G$-action on $T^{*}\HS$ is stable (Remark~\ref{symp.stab=>stab}),
whence $\kk(T^{*}\HS)^G=\Quot\kk[T^{*}\HS]^G$ is
Poisson-commutative as well.
%
%%% A variant:
%
%Take an affine embedding $\HS\embeds X$ and consider the
%subalgebra $\Dd\subseteq\Dd(\HS)$ generated by vector fields
%on~$X$. Then $M=\Spec\gr\Dd$ is an affine embedding
%of~$T^{*}\HS$. By complete reducibility, $\kk[M]^G=\gr\Dd^G$ is
%Poisson-commutative. But the $G$-action on $M$ is stable
%(Remark~\ref{symp.stab=>stab}), whence
%$\kk(T^{*}\HS)^G=\Quot\kk[M]^G$ is Poisson-commutative as well.

\item[\ref{coll.int}$\iff$\ref{open-B}] This equivalence is due to
Mikityuk \cite{sph(ss)} (for affine~$\HS$).

A complete system of Poisson-commuting functions on $M_{\HS}$ can
be constructed by the method of argument shift~\cite{arg.shift}:
choose a regular semisimple element $\xi\in\g^{*}$ and consider
the derivatives $\partial_{\xi}^nf$ of all $f\in\kk[\g^{*}]^G$.
The functions $\partial_{\xi}^nf$ Poisson-commute and produce a
complete involutive system on $Gx\subset\g^{*}$ (for
general~$\xi$) whenever $\idx\g_x=\idx\g$, where
$\idx\g=d_G(\g^{*})$ \cite{comp.inv}. In the symplectically stable
case, general points $x\in M_{\HS}$ are semisimple and
$\idx\g_x=\idx\g=\rk\g$. Generally, the equality $\idx\g_x=\idx\g$
was conjectured by Elashvili and proved by Charbonnel \cite{index}
for $\forall x\in\g^{*}$.

Since symplectic leaves of the Poisson structure on $M_{\HS}$ are
$G$-orbits, there are $(d_G(M_{\HS})+\dim
M_{\HS})/2=\dim\HS-c(\HS)$ independent Poisson-commuting
collective functions. Thus a complete involutive system of
collective functions exists iff $c(\HS)=0$. \qedhere\end{roster}
\end{proof}

Since $T^{*}\HS=G\itimes{H}\h^{\ann}$, weak commutativity is
readily reformulated in terms of the coadjoint representation
\cite{sph(ss)}, \cite{c&r(hom)}, \cite[II.4.1]{comm.hom}.
\begin{theorem}
$(G,H)$ is a spherical pair iff general points
$\alpha\in\h^{\ann}$ satisfy any of the equivalent conditions:
\begin{roster}
\renewcommand{\theenumi}{\reftag{Ad\arabic{enumi}}}
\item\label{dim=2dim} $\dim G\alpha=2\dim H\alpha$
\item\label{Lagrangian} $H\alpha$~is a Lagrangian subvariety in
$G\alpha$ w.r.t.\ the Kirillov form. \item\label{Richardson}
(Richardson condition) $\g\alpha\cap\h^{\ann}=\h\alpha$
\end{roster}
\end{theorem}
The Richardson condition means that $G\alpha\cap\h^{\ann}$ is a
finite union of open $H$-orbits \cite[1.5]{IT}.
\begin{proof}
\begin{roster}
\item[\ref{coiso}$\iff$\ref{dim=2dim}] Recall that the moment map
$\Phi:G\itimes{H}\h^{\ann}\to\g^{*}$ is defined via replacing the
${*}$-action by the coadjoint action (Example~\ref{moment(G/H)}).
We have
%*
\begin{gather*}
d_G(T^{*}\HS)=\dim\HS-\dim H\alpha\qquad\text{and}\\
\df T^{*}\HS=\dim G_{\Phi(\1*\alpha)}/G_{\1*\alpha}=
\dim G_{\alpha}/H_{\alpha}\\
\intertext{Hence} \cork T^{*}\HS=d_G(T^{*}\HS)-\df T^{*}\HS= \dim
G\alpha-2\dim H\alpha
\end{gather*}
%*
\item[\ref{dim=2dim}$\iff$\ref{Lagrangian}] The Kirillov form
vanishes on~$\h\alpha$.
\item[\ref{Lagrangian}$\iff$\ref{Richardson}] Stems from
$(\g\alpha)\cap\h^{\ann}=(\h\alpha)^{\sort}$, the
skew-orthocomplement w.r.t.\ the Kirillov form.
\qedhere\end{roster}
\end{proof}

Invariant functions on cotangent bundles of spherical homogeneous
spaces have a nice structure.
\begin{proposition}[{\cite[7.2]{W_X}}]\label{coiso(sph)}
If $\HS=G/H$ is spherical, then
$\kk[T^{*}\HS]^G\iso\kk[\widetilde{L}_{\HS}]\iso\kk[\ab^{*}]^{W_{\HS}}$
is a polynomial algebra; there are similar isomorphisms for fields
of rational functions.
\end{proposition}
\begin{proof}
By Proposition~\ref{coll&inv} and \ref{weak.comm},
$\kk(T^{*}\HS)^G\iso\kk(\widetilde{L}_{\HS})\iso\kk(\ab^{*})^{W_{\HS}}$.
By Lemma~\ref{norm(M_X)&a^*},
$\widetilde\pi_G\widetilde\Phi:T^{*}\HS\to\widetilde{L}_{\HS}$ is
a surjective morphism of normal varieties. Therefore any
$f\in\kk(T^{*}\HS)^G$ having poles on $T^{*}\HS$ must have poles
on~$\widetilde{L}_{\HS}$, whence
$\kk[T^{*}\HS]^G=\kk[\widetilde{L}_{\HS}]$. The latter algebra is
polynomial for $W_{\HS}$ is generated by reflections
(Theorem~\ref{W_X&Z(X)}).
\end{proof}
In other words, invariants of the coisotropy representation form a
polynomial algebra $\kk[\h^{\ann}]^H\iso\kk[\ab^{*}]^{W_{G/H}}$
for any spherical pair $(G,H)$.
\begin{remark}
A similar assertion in the non-commutative setup was proved
in~\cite{HC-hom}. Namely, all invariant differential operators on
a spherical space $\HS$ are completely regular, whence
$\Dd(\HS)^G$ is a polynomial ring isomorphic to
$\kk[\rho+\ab^{*}]^{W_{\HS}}$ (see Remark~\ref{non-comm}). In
particular, every spherical homogeneous space is commutative.
\end{remark}

In our considerations $G$ was always assumed to be reductive.
However some of the concepts introduced above are reasonable even
for non-reductive $G$ assuming $H$ be reductive instead. Some of
the above results remain valid:
\begin{enumerate}
\item\label{WS=>GP} If $\HS=G/H$ is weakly symmetric, then $(G,H)$
is a Gelfand pair. \item $\HS$~is commutative iff $(G,H)$ is a
Gelfand pair. \item\label{C=>WC} A commutative space $\HS$ is
weakly commutative provided that $\kk[\h^{\ann}]^H$ separates
generic $H$-orbits in~$\h^{\ann}$.
\end{enumerate}
The above proofs work in this case: the functor $(\cdot)^G$ is
exact on global sections of $G$-sheaves on~$\HS$ since $\HS$ is
affine and $(\cdot)^H$ is exact on rational $H$-modules, and orbit
separation in~\ref{C=>WC} guarantees
$\kk(T^{*}\HS)^G=\Quot\kk[T^{*}\HS]^G$. The converse implication
in~\ref{WS=>GP} fails, the simplest counterexample being:
\begin{example}[\cite{not.weak.symm}]
Put $H=\Sp_{2n}(\kk)$, $G=H\semitimes N$, where $N=\exp\n$ is a
unipotent group associated with the Heisenberg type Lie algebra
$\n=(\kk^{2n}\oplus\kk^{2n})\oplus\kk^3$, the commutator in $\n$
being defined by the identification
$\E^2(\kk^{2n}\oplus\kk^{2n})^{\Sp_{2n}(\kk)}\iso\kk^3=\z(\n)$.
Then $(G,H)$ is a Gelfand pair, but $\HS$ is not weakly symmetric.
\end{example}
Also, the implication \ref{C=>WC} fails if the orbit separation is
violated. The reason is that there may be too few invariant
differential operators. For instance, in the previous example,
replace $H$ by $\kk^{*}$ acting on $\kk^{2n}$ via a character
$\chi\neq0$ and on $\kk^3$ via~$2\chi$. Then $\HS$ is not weakly
commutative while $\Dd(\HS)^G=\kk$.

The classes of weakly symmetric and (weakly) commutative
homogeneous spaces were first introduced and examined in
Riemannian geometry and harmonic analysis, see the
survey~\cite{comm.hom}. We shall review the analytic viewpoint
now.

Quitting a somewhat restrictive framework of algebraic varieties,
one may consider the above properties of homogeneous spaces in the
category of Lie group actions, making appropriate modifications in
formulations. For instance, instead of regular or rational
functions one considers arbitrary analytic or differentiable
functions. Some of these properties receive new interpretation in
terms of differential geometry, e.g., \ref{coll.int}~means that
invariant Hamiltonian dynamic systems on $T^{*}\HS$ are completely
integrable in the class of Noether integrals \cite{coll.int},
\cite{sph(ss)}.

The situation, where $H$ is a compact subgroup of a real Lie
group~$G$, i.e., $\HS=G/H$ is a Riemannian homogeneous space, has
attracted the main attention of researchers. Most of the above
results were originally obtained in this setting.

The properties \ref{V^H}, \ref{mult-free(k[O])} are naturally
reformulated here in the category of unitary representations of
$G$ replacing $\kk[\HS]$ by~$L^2(\HS)$. In \ref{Hecke} one
considers the algebra $\Aa(G)$ of complex measures with compact
support on~$G$. The conditions \ref{Selberg}, \ref{Gelfand} are
formulated for \emph{all} (not only general) points (which is
equivalent for compact~$H$); there is also an infinitesimal
characterization of weak symmetry \cite[I.1.5]{comm.hom}.

There are the following implications:
%*
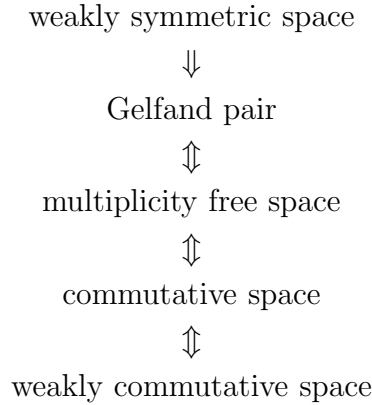
\begin{gather*}
\text{weakly symmetric space}\\
\Downarrow\\
\text{Gelfand pair}\\
\Updownarrow\\
\text{multiplicity free space}\\
\Updownarrow\\
\text{commutative space}\\
\Updownarrow\\
\text{weakly commutative space}
\end{gather*}
%*
The implication \ref{Gelfand}$\implies$\ref{Hecke} is due to
Gelfand~\cite{symm=>MF} and \ref{Hecke}$\iff$\ref{mult-free(k[O])}
was proved in~\cite{rep(grp)}. The equivalence
\ref{Hecke}$\iff$\ref{comm} is due to Helgason \cite[Ch.IV,
B13]{geom.anal} and Thomas~\cite{Thom}, for a proof see
\cite[I.2.5]{comm.hom}. The implication
\ref{comm}$\implies$\ref{weak.comm} is easy \cite[I.4.2]{comm.hom}
and the converse was recently proved by Rybnikov~\cite{WC=>C}.

A classification of commutative Riemannian homogeneous spaces was
obtained by Yakimova \cite{sph.indec}, \cite{GP} using partial
results of Vinberg \cite{Heisenberg} and the classification of
affine spherical spaces from~\ref{c&r<=1}.

An algebraic homogeneous space $\HS=G/H$ over $\kk=\CC$ may be
considered as a homogeneous manifold in the category of complex or
real Lie group actions. At the same time, if $(G,H)$ is defined
over~$\RR$, then $\HS$ has a real form $\HS(\RR)$ containing
$G(\RR)/H(\RR)$ as an open orbit (in classical topology). Thus
$G/H$ may be regarded as the complexification of $G(\RR)/H(\RR)$,
a homogeneous space of a real Lie group~$G(\RR)$.

It is easy to see that $G/H$ is commutative (resp.\ weakly
commutative, multiplicity free, weakly symmetric, satisfies
\ref{Hecke}, \ref{coll.int}) iff $G(\RR)/H(\RR)$ is so. In other
words, the above listed properties are stable under
complexification and passing to a real form.

This observation leads to the following criterion of sphericity,
which is a ``real form'' of Theorem~\ref{WC&C}.

By Chevalley's theorem, there exists a projective embedding
$\HS\subseteq\PP(V)$ for some $G$-module~$V$. Assume that $G$ is
reductive and $K\subset G$ is a compact real form. Then $V$ can be
endowed with a $K$-invariant Hermitian inner
product~$(\cdot|\cdot)$, which induces a K{\"a}hlerian metric on
$\PP(V)$ and on~$\HS$ (the \emph{Fubini--Studi metric}). The
imaginary part of this metric is a real symplectic form. The
action $K:\PP(V)$ is Poisson, the moment map
$\Phi:\PP(V)\to\ka^{*}$ being defined by the formula
%*
\begin{equation*}
\bigl\langle\Phi(\langle v\rangle),\xi\bigr\rangle=
\frac1{2i}\cdot\frac{(v|\xi{v})}{(v|v)}, \qquad\forall v\in V,\
\xi\in\ka
\end{equation*}
%*
\begin{theorem}[{\cite{moment}, \cite{mult.free(C)},
\cite[\S13]{spher.var}}] $\HS$~is spherical iff generic $K$-orbits
in $\HS$ are coisotropic w.r.t.\ the Fubini--Studi form, or
equivalently, $C^{\infty}(\HS)$~is Poisson-commutative.
\end{theorem}
\begin{proof}
First note that generic $K$-orbits in $\HS$ are coisotropic iff
%*
\begin{equation}\label{sph.ham}
d_K(\HS)=\df\HS=\rk K-\rk K_{*}
\end{equation}
%*
where $K_{*}$ is the stabilizer of general position for $K:\HS$.
The condition \eqref{sph.ham} does not depend on the symplectic
structure.

If $\HS$ is affine, then the assertion can be directly reduced to
Theorem~\ref{WC&C} by complexification. W.l.o.g.~$K\cap H$ is a
compact real form of~$H$. Using the Cartan decompositions
$G=K\cdot\exp i\ka$, $H=(K\cap H)\cdot\exp i(\ka\cap\h)$, we
obtain a $K$-diffeomorphism
%*
\begin{equation*}
\HS\iso K\itimes{K\cap H} i\ka/i(\ka\cap\h)\iso T^{*}(K/K\cap H)
\end{equation*}
%*
Complexifying the r.h.s.\ we obtain~$T^{*}\HS$.

In the general case, it is more convenient to apply the theory of
doubled actions~(\ref{cotangent}).

There exists a Weyl involution $\theta$ of $G$ commuting with the
Hermitian conjugation $g\mapsto g^{*}$. The mapping
$g\mapsto\overline{g}:=\theta(g^{*})^{-1}$ is a complex
conjugation on $G$ defining a split real form~$G(\RR)$. There
exists a $G(\RR)$-stable real form $V(\RR)\subset V$ such that
$(\cdot|\cdot)$ takes real values on~$V(\RR)$. The complex
conjugation on $V$, $\PP(V)$, or $G$ is defined by conjugating the
coordinates or matrix entries w.r.t.\ an orthonormal basis
in~$V(\RR)$.

It follows that the complex conjugate variety $\overline\HS$ is
naturally embedded in $\PP(V)$ as a $G$-orbit. Complexifying the
action $K:\HS$ we obtain the diagonal action
$G:\HS\times\overline\HS$,
$g(x,\overline{y})=(gx,\theta(g)\overline{y})$, $\forall g\in G,\
x,y\in\HS$. This action differs slightly from the doubled action,
but Theorems~\ref{sgp(dbl)}--\ref{c&r<=dbl} remain valid, together
with the proofs. Now it follows from
\eqref{2c+r<=dbl}--\eqref{r<=dbl} that $\HS$ is spherical iff
%*
\begin{equation*}
d_G(\HS\times\overline\HS)=\rk G-\rk G_{*}
\end{equation*}
%*
where $G_{*}=K_{*}(\CC)$ is the stabilizer of general position for
$G:\HS\times\overline\HS$. The latter condition coincides
with~\eqref{sph.ham}.
\end{proof}

\section{Symmetric spaces}
\label{symmetric}

The concept of a Riemannian symmetric space was introduced by
{\'E}.~Cartan \cite{Riemann.symm.1}, \cite{Riemann.symm.2}. A
(globally) symmetric space is defined as a connected Riemannian
manifold $\HS$ such that for $\forall x\in\HS$ there exists an
isometry $s_x$ of $\HS$ inverting the geodesics passing
through~$x$. Symmetric spaces form a very important class of
Riemannian spaces including all classical geometries. The theory
of Riemannian symmetric spaces is well developed,
see~\cite{symm.full}.

In particular, it is easy to see that a symmetric space $\HS$ is
homogeneous w.r.t.\ the identity component $G$ of the full
isometry group, so that $\HS=G/H$, where $H=G_{\bp}$ is the
stabilizer of a fixed base point. The geodesic symmetry
$s=s_{\bp}$ is an involutive automorphism of $\HS$
normalizing~$G$. It defines an involution $\sigma\in\Aut G$ by
$\sigma(g)=sgs^{-1}$. From the definition of a geodesic symmetry
one deduces that $(G^{\sigma})^0\subseteq H\subseteq G^{\sigma}$.
Furthermore, reducing $G$ to a smaller transitive isometry group
if necessary, one may assume that $\g$ is a reductive Lie algebra.
This leads to the following algebraic definition of a symmetric
space, which we accept in our treatment.
\begin{definition}
An (algebraic) \emph{symmetric space} is a homogeneous algebraic
variety $\HS=G/H$, where $G$ is a connected reductive group
equipped with a non-identical involution $\sigma\in\Aut G$, and
$(G^{\sigma})^0\subseteq H\subseteq G^{\sigma}$.
\end{definition}
Riemannian symmetric spaces are locally isomorphic to real forms
(with compact isotropy subgroups) of algebraic symmetric spaces
over~$\CC$.

It is reasonable to impose a restriction $\ch\kk\ne2$ on the
ground field.

\begin{remark}
If $G$ is semisimple simply connected, then $G^{\sigma}$ is
connected \cite[8.2]{Stein}, whence $H=G^{\sigma}$. On the other
side, if $G$ is adjoint, then $G^{\sigma}=N(H)$
\cite[2.2]{emb(symm)}.
\begin{question}
Prove it using~\eqref{H-fixed}.
\end{question}
\end{remark}

The differential of~$\sigma$, denoted by the same letter by abuse
of notation, induces a $\ZZz_2$-grading
%*
\begin{equation}\label{Z2-grad}
\g=\h\oplus\m
\end{equation}
%*
where $\h,\m$ are the $(\pm1)$-eigenspaces of~$\sigma$.

The subgroup $H$ is reductive \cite[\S8]{Stein}, thence $\HS$ is
an affine algebraic variety. More specifically, consider a
morphism $\tau:G\to G$, $\tau(g)=\sigma(g)g^{-1}$. Observe that
$\tau$ is the orbit map at $\1$ for the $G$-action on $G$ by
\emph{twisted conjugation}: $g\circ x=\sigma(g)xg^{-1}$. It is not
hard to prove the following result.
\begin{proposition}[{\cite[2.2]{inv}}]\label{symm<G}
$\tau(G)\iso G/G^{\sigma}$ is a connected component of $\{x\in
G\mid\sigma(x)=x^{-1}\}$.
\end{proposition}
\begin{example}
Let $G=\GL_n(\kk)$ and $\sigma$ be defined by
$\sigma(x)=(x^{\tran})^{-1}$. Then $G^{\sigma}=\Or_n(\kk)$ and
$\tau(G)=\{x\in G\mid\sigma(x)=x^{-1}\}$ is the set of
non-degenerate symmetric matrices, which is isomorphic to
$\GL_n(\kk)/\Or_n(\kk)$.

However, if $\sigma$ is an inner involution, i.e., the conjugation
by a matrix of order~$2$, then the set of matrices $x$ such that
$\sigma(x)=x^{-1}$ is disconnected. The connected components are
determined by the collection of eigenvalues of~$x$, which
are~$\pm1$.
\end{example}

The local and global structure of symmetric spaces is examined in
\cite{iso(symm)}, \cite{symm.full} (transcendental methods),
\cite{symm}, \cite{emb(symm)} ($\ch\kk=0$), \cite{inv(symm)},
\cite{inv}. We follow these sources in our exposition. The
starting point is an analysis of $\sigma$-stable tori.

\begin{lemma}\label{sigma-stab}
Every Borel subgroup $B\subseteq G$ contains a $\sigma$-stable
maximal torus~$T$.
\end{lemma}
\begin{proof}
The group $B\cap\sigma(B)$ is connected, solvable, and
$\sigma$-stable. By \cite[7.6]{Stein} it contains a
$\sigma$-stable maximal torus~$T$, which is a maximal torus
in~$G$, too.
\end{proof}
\begin{corollary}
Every $\sigma$-stable torus $S\subseteq G$ is contained in a
$\sigma$-stable maximal torus~$T$.
\end{corollary}
\begin{proof}
Put $T=Z\cdot T'$, where $Z$ is the connected center and $T'$ any
$\sigma$-stable maximal torus in the commutator subgroup
of~$Z_G(S)$.
\end{proof}

A $\sigma$-stable torus $T$ decomposes into an almost direct
product $T=T_0\cdot T_1$, where $T_0\subseteq H$ and $T_1$ is
\emph{$\sigma$-split}, which means that $\sigma$ acts on $T_1$ as
the inversion.

Let $\Delta$ denote the root system of $G$ w.r.t.~$T$ and
$\g_{\alpha}\subset\g$ the root subspace corresponding to
$\alpha\in\Delta$. One may choose root vectors
$e_{\alpha}\in\g_{\alpha}$ in such a way that
$e_{\alpha},e_{-\alpha},h_{\alpha}=[e_{\alpha},e_{-\alpha}]$ form
an $\sgl_2$-triple for $\forall\alpha\in\Delta$. Clearly,
$\sigma$~acts on $\Ch(T)$ leaving $\Delta$ stable. Choosing
$e_{\alpha}$ in a compatible way allows to subdivide all roots
into \emph{complex}, \emph{real}, and \emph{imaginary}
(\emph{compact} or \emph{non-compact}) ones, according to
Table~\ref{roots}.

\begin{table}[!h]
\caption{Root types w.r.t.\ an involution} \label{roots}
\begin{center}
\begin{tabular}{|c|c|c|c|c|}
\hline $\alpha$ & complex & real &
\multicolumn{2}{|c|}{imaginary} \\
\cline{4-5}
&&& compact & non-compact \\
\hline

$\sigma(\alpha)$ & $\neq\pm\alpha$ & $-\alpha$ &
$\alpha$ & $\alpha$ \\
\hline

$\sigma(e_{\alpha})$ & $e_{\sigma(\alpha)}$ & $e_{-\alpha}$ &
$e_{\alpha}$ & $-e_{\alpha}$ \\
\hline
\end{tabular}
\end{center}
\end{table}

We fix an inner product on $\Ch(T)\otimes\QQ$ invariant under the
Weyl group $W=N_G(T)/T$ and~$\sigma$.
\begin{question}
Use the language of coroots instead.
\end{question}
Then $\Ch(T)\otimes\QQ$ is identified with $\CoCh(T)\otimes\QQ$
and with the orthogonal sum of $\Ch(T_0)\otimes\QQ$ and
$\Ch(T_1)\otimes\QQ$. The coroots $\alpha^{\vee}\in\Delta^{\vee}$
($\alpha\in\Delta$) are identified with $2\alpha/(\alpha,\alpha)$.
Let $\langle\alpha|\beta\rangle=\langle\alpha^{\vee},\beta\rangle=
2(\alpha,\beta)/(\alpha,\alpha)$ denote the Cartan pairing
on~$\Ch(T)$ and
$r_{\alpha}(\beta)=\beta-\langle\alpha|\beta\rangle\alpha$ the
reflection of $\beta$ along~$\alpha$.

Two opposite classes of $\sigma$-stable maximal tori are of
particular importance in the theory of symmetric spaces.

\begin{lemma}
If $\dim T_0$ is maximal possible, then $T_0$ is a maximal torus
in~$H$ and $Z_G(T_0)=T$. Moreover, $T$~is contained in a
$\sigma$-stable Borel subgroup $B\subseteq G$ such that
$(B^{\sigma})^0$ is a Borel subgroup in~$H$.
\end{lemma}
\begin{proof}
If $Z_G(T_0)\neq T$, then the commutator subgroup $Z_G(T_0)'$ and
$(Z_G(T_0)')^{\sigma}$ have positive dimension. Hence $T_0$ can be
extended by a subtorus in $(Z_G(T_0)')^{\sigma}$, a contradiction.
Now choose a Borel subgroup of $H$ containing $T_0$ and extend it
to a Borel subgroup $B$ of~$G$. Then $B\supseteq T$. If $B$ were
not $\sigma$-stable, then there would exist a root
$\alpha\in\Delta^{+}$ such that $\sigma(\alpha)\in\Delta^{-}$.
Then $e_{\pm\alpha}+\sigma(e_{\pm\alpha})$ are opposite root
vectors in $\h$ outside the Borel subalgebra $\br^{\sigma}$, a
contradiction.
\end{proof}

In particular, if $T_0$ is maximal, then there are no real roots,
and $\sigma$ preserves the set $\Delta^{+}$ of positive roots
(w.r.t.~$B$) and induces a diagram involution $\overline\sigma$ of
the set $\Pi\subseteq\Delta^{+}$ of simple roots. If $G$ is of
simply connected type, then $\overline\sigma$ extends to an
automorphism of $G$ so that $\sigma=\overline\sigma\cdot\sigma_0$,
where $\sigma_0$ is an inner automorphism defined by an element
of~$T_0$.

Consider the set $\overline\Delta=\{\overline\alpha=\alpha|_{T_0}
\mid\alpha\in\Delta\}\subset\Ch(T_0)$.
\begin{question}
Perhaps, $\Delta'$ would look better.
\end{question}
Clearly, $\overline\Delta$~consists of the roots of $H$ w.r.t.\
$T_0$ and the nonzero weights of $T_0:\m$. The restrictions of
complex roots belong to both subsets, the eigenvectors being
$e_{\alpha}+\sigma(e_{\alpha})\in\h$,
$e_{\alpha}-\sigma(e_{\alpha})\in\m$, whereas (non-)compact roots
restrict to roots of~$H$ (resp.\ weights of~$\m$).
\begin{lemma}\label{Kac}
$\overline\Delta$~is a (possibly non-reduced) root system with
base~$\overline\Pi$. The simple roots of $H$ and the (nonzero)
lowest weights of $H:\m$ form an affine simple root
system~$\widetilde\Pi$, i.e.,
$\langle\overline\alpha|\overline\beta\rangle\in\ZZ_{-}$ for
distinct $\overline\alpha,\overline\beta\in\widetilde\Pi$.
\end{lemma}
\begin{proof}
Note that the restriction of $\alpha\in\Delta$ to $T_0$ is the
orthogonal projection to $\Ch(T_0)\otimes\QQ$, so that
$\overline\alpha=(\alpha+\sigma(\alpha))/2$. If $\alpha$ is
complex, then $\langle\alpha|\sigma(\alpha)\rangle=0\text{ or }-1$
(otherwise $\alpha-\sigma(\alpha)$ would be a real root), In the
second case, $2\overline\alpha=\alpha+\sigma(\alpha)$ is a
non-compact root with a root vector
$e_{\alpha+\sigma(\alpha)}=[e_{\alpha},\sigma(e_{\alpha})]$.

A direct computation shows that $\forall\alpha,\beta\in\Delta:
\langle\overline\alpha|\overline\beta\rangle\in\ZZ$ and the
reflections $r_{\overline\alpha}$ preserve~$\overline\Delta$, see
Table~\ref{Cartan}.
\begin{table}[!h]
\caption{Cartan numbers and reflections for restricted roots}
\label{Cartan}
\begin{center}
\renewcommand{\arraystretch}{1.2}
\begin{tabular}{|c|c|c|}
\hline Case & $\langle\overline\alpha|\overline\beta\rangle$ &
$r_{\overline\alpha}(\overline\beta)$ \\
\hline

$\alpha=\sigma(\alpha)$ & $\langle\alpha|\beta\rangle$ &
$\overline{r_{\alpha}(\beta)}$ \\
\hline

$\langle\alpha|\sigma(\alpha)\rangle=0$ &
$\langle\alpha|\beta\rangle+\langle\sigma(\alpha)|\beta\rangle$
& $\overline{r_{\alpha}r_{\sigma(\alpha)}(\beta)}$ \\
\hline

$\langle\alpha|\sigma(\alpha)\rangle=-1$ &
$2\langle\alpha|\beta\rangle+2\langle\sigma(\alpha)|\beta\rangle$
& $r_{2\overline{\alpha}}(\overline{\beta})=
\overline{r_{\alpha+\sigma(\alpha)}(\beta)}$ \\
\hline
\end{tabular}
\end{center}
\end{table}
Hence $\overline\Delta$ is a root system. Restricting to $T_0$ the
expression of $\alpha\in\Delta$ as a linear combination of $\Pi$
with integer coefficients of the same sign yields a similar
expression of $\overline\alpha$ in terms of~$\overline\Pi$. Since
$\overline\Pi$ is linearly independent, it is a base
of~$\overline\Delta$.

Note that $\overline\alpha=\overline\beta$ iff $\alpha=\beta$ or
$\sigma(\alpha)=\beta$. (Otherwise $\alpha-\beta$ or
$\sigma(\alpha)-\beta$ would be a real root, depending on whether
$\langle\alpha|\beta\rangle>0$ or
$\langle\sigma(\alpha)|\beta\rangle>0$.) Therefore the nonzero
weights occur in $\m$ with multiplicity~$1$.

To prove the second assertion, it suffices to consider the Cartan
numbers $\langle\overline\alpha|\overline\beta\rangle$ of lowest
weights of~$\m$. Assuming
$\langle\overline\alpha|\overline\beta\rangle>0$ yields w.l.o.g.\
$\langle\alpha|\beta\rangle>0$, whence
$\gamma=\alpha-\beta\in\Delta$,
$e_{\alpha}=[e_{\beta},e_{\gamma}]$. If $\beta$ is non-compact,
then ${[e_{\beta},e_{\gamma}+\sigma(e_{\gamma})]}=
e_{\alpha}-\sigma(e_{\alpha})$. If $\beta$ is complex, then
$\gamma$ is complex, too. (Indeed, $\overline\gamma$~is not longer
than $\overline\alpha,\overline\beta$, i.e., is shorter than any
root of~$\Delta$.) Then
${\beta+\sigma(\gamma),\sigma(\beta)+\gamma\notin\Delta}$, whence
$[e_{\beta}-\sigma(e_{\beta}),e_{\gamma}+\sigma(e_{\gamma})]=
e_{\alpha}-\sigma(e_{\alpha})$. In both cases, either
$\overline\beta$ or $\overline\alpha$ is not a lowest weight, a
contradiction.
\end{proof}
\begin{remark}
If some of the Cartan numbers of $\Delta$ vanish in~$\kk$, then
the previous arguments concerning lowest weights does not work.
The assertion on $\widetilde\Pi$ is true only if one interprets
lowest weights in the combinatorial sense as those weights of $\m$
which cannot be obtained from other weights by adding simple roots
of~$H$. However this happens only for $G=\GGg_2$ (if $\ch\kk=3$),
where the unique (up to conjugation) involution is easy to
describe by hand.
\end{remark}

In a usual way, the system $\widetilde\Pi$ together with the
respective Cartan numbers is encoded by an (affine) Dynkin
diagram. Marking the nodes corresponding to the simple roots of
$H$ by black, and those corresponding to the lowest weights of
$\m$ by white, one obtains the so-called \emph{Kac diagram} of the
involution~$\sigma$, or of the symmetric space~$\HS$. From the Kac
diagram one easily recovers $\h$ and (at least in characteristic
zero) the (co)isotropy representation $H^0:\m$.
\begin{example}\label{diag}
Let $H$ be diagonally embedded in $G=H\times H$, where $\sigma$
permutes the factors. Here the Kac diagram is the affine Dynkin
diagram of $H$ with the white nodes corresponding to the lowest
roots, e.g.:
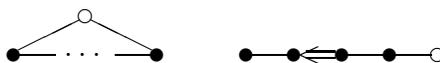
\begin{center}
%TeXCAD Picture [A0-K.pic]. Options:
%\grade{\off}
%\emlines{\off}
%\epic{\off}
%\beziermacro{\off}
%\reduce{\on}
%\snapping{\off}
%\quality{2.00}
%\graddiff{0.01}
%\snapasp{1}
%\zoom{1.00}
\unitlength .7ex % = .7pt
\linethickness{0.4pt}
\begin{picture}(15.67,5.67)(0,0)
\put(0,1){\circle*{1.33}} \put(15,1){\circle*{1.33}}
\put(.67,1){\line(1,0){3.67}} \put(10.67,1){\line(1,0){3.67}}
\put(7.5,1){\makebox(0,0)[cc]{$\dots$}} \put(7.5,5){\circle{1.33}}
\put(.5,1.5){\line(2,1){6.4}} \put(14.5,1.5){\line(-2,1){6.4}}
\end{picture}
\qquad
%TexCad Options
%\grade{\off}
%\emlines{\off}
%\beziermacro{\off}
%\reduce{\on}
%\snapping{\off}
%\quality{2.00}
%\graddiff{0.01}
%\snapasp{1}
%\zoom{1.00}
\unitlength 0.70ex \linethickness{0.4pt}
\begin{picture}(20.67,1.67)
\put(0.00,1.00){\circle*{1.33}} \put(5.00,1.00){\circle*{1.33}}
\put(10.00,1.00){\circle*{1.33}} \put(15.00,1.00){\circle*{1.33}}
\put(20.00,1.00){\circle{1.33}} \put(0.67,1.00){\line(1,0){3.67}}
\put(10.67,1.00){\line(1,0){3.67}}
\put(9.33,1.33){\line(-1,0){3.00}}
\put(9.33,0.67){\line(-1,0){3.00}}
\put(5.33,1.00){\makebox(0,0)[lc]{$<$}}
\put(15.67,1.00){\line(1,0){3.67}}
\end{picture}
\end{center}
\end{example}

Now consider an opposite class of $\sigma$-stable maximal tori.

\begin{lemma}\label{sigma-split}
There exist non-trivial $\sigma$-split tori.
\end{lemma}
\begin{proof}
In the converse case $\sigma$ acts identically on every
$\sigma$-stable torus. Lemma~\ref{sigma-stab} implies that all
Borel subgroups are $\sigma$-stable. Then all maximal tori are
$\sigma$-stable and even pointwise fixed, whence $\sigma$ is
identical.
\end{proof}

\begin{lemma}
If $T_1$ is a maximal $\sigma$-split torus, then $L=Z_G(T_1)$
decomposes into an almost direct product $L=L_0\cdot T_1$, where
$L_0=L\cap H$.
\end{lemma}
\begin{proof}
Clearly, $L$~and the commutator subgroup $L'$ are $\sigma$-stable.
If $L'\not\subseteq H$, then $T_1$ could be extended by a
non-trivial $\sigma$-split torus in $L'$ by
Lemma~\ref{sigma-split}, a contradiction. The assertion follows
from $L'\subseteq H$.
\end{proof}

Choose a general one-parameter subgroup $\gamma\in\CoCh(T_1)$ and
consider the associated parabolic subgroup $P=P(\gamma)$ with the
Lie algebra $\p=\tr\oplus
\bigoplus_{\langle\alpha,\gamma\rangle\geq0}\g_{\alpha}$. Clearly,
$L\subseteq P$ is a Levi subgroup and
$\Ru\p=\bigoplus_{\langle\alpha,\gamma\rangle>0}\g_{\alpha}$. Note
that $\sigma(P)=P^{-}$ (since
$\langle\sigma(\alpha),\gamma\rangle=-\langle\alpha,\gamma\rangle$,
$\forall\alpha\in\Delta$). In fact, all minimal parabolics having
this property are obtained as above \cite[1.2]{symm}. It follows
that $\h$ is spanned by $\lv_0$ and
$e_{\alpha}+\sigma(e_{\alpha})$ over all $\alpha\in\Delta$  such
that $\langle\alpha,\gamma\rangle\geq0$. This yields:
\begin{Claim}[Iwasawa decomposition]
$\g=\h\oplus\tr_1\oplus\Ru\p$
\end{Claim}

As a consequence, we obtain
\begin{theorem}\label{symm=>sph}
Symmetric spaces are spherical.
\end{theorem}

Indeed, choosing a Borel subgroup $B\subseteq P$, $B\supseteq T$
yields~\ref{h+gb=g}. There are many other ways to verify this
fact. For instance, it is easy to verify the Richardson
condition~\ref{Richardson}:  $\forall\xi\in\m\iso\h^{\ann}$ one
has $[\g,\xi]\cap\m=[\h,\xi]$ for $[\m,\xi]\subseteq\h$. One can
also check the Gelfand condition \ref{Gelfand} for elements in a
dense subset $\tau(G)H\subseteq G$: $g=xh,\ x\in\tau(G),\
h\in\h\implies \sigma(g)=x^{-1}h=hg^{-1}h$. The multiplicity free
property (for compact Riemannian symmetric spaces and unitary
representations) was established already by {\'E}.~Cartan
\cite[n$^{\circ}$17]{MF(symm)}.

The Iwasawa decomposition clarifies the local structure of a
symmetric space. Namely, $\HS$~contains a dense orbit
$P\cdot\bp\iso P/L_0\iso\Ru{P}\times A$, where $A=T/T\cap H$ is
the quotient of $T_1$ by an elementary Abelian $2$-group $T_1\cap
H$. We have $\ab\iso\tr_1$, $\RG(\HS)=\Ch(A)$, $r(\HS)=\dim\ab$.
The notation here agrees with Corollary~\ref{loc.str.gen}
and~\ref{horosph}.

\begin{lemma}\label{conj(split)}
All maximal $\sigma$-split tori are $H^0$-conjugate.
\end{lemma}
\begin{proof}
In the above notation, $PH$~is open in~$G$, whence the $H^0$-orbit
of $P$ is open in $G/P$. Since $P$ coincides with the normalizer
of the open $B$-orbit in~$\HS$, all such parabolics are
$G$-conjugate and therefore $H^0$-conjugate. Hence the Levi
subgroups $L=P\cap\sigma(P)$ and finally the maximal
$\sigma$-split tori $T_1=(Z(L)^0)_1$ are $H^0$-conjugate.
\end{proof}

If $T_1$ is maximal, then every imaginary root is compact and
$\sigma$ maps positive complex or real roots to negative ones.
Compact (simple) roots form (the base of) the root system of~$L$.

The endomorphism $\iota=-w_L\sigma$ of $\Ch(T)$ preserves
$\Delta^{+}$ and induces a diagram involution of the set $\Pi$ of
simple roots. (Here $w_L$ is the longest element in the Weyl group
of~$L$.) Since $w_Gw_L\sigma$ preserves $\Delta^{+}$ and differs
from $\sigma$ by an inner automorphism, it coincides with the
diagram automorphism~$\overline\sigma$, whence
$\iota(\lambda)=\overline\sigma(\lambda)^{*}$,
$\forall\lambda\in\Ch(T)$.

Consider the set $\Delta_{\HS}\subset\Ch(T_1)$ and the subset
$\Pi_{\HS}\subset\Delta_{\HS}$ consisting of the restrictions
$\overline\alpha=\alpha|_{T_1}$ of complex and real roots
$\alpha\in\Delta$ (resp. $\alpha\in\Pi$) to~$T_1$.
\begin{question}
Multiply by $2$ for compatibility with~\ref{inv.mot}.
\end{question}
\begin{lemma}
$\Delta_{\HS}$~is a (possibly non-reduced) root system with
base~$\Pi_{\HS}$, called the \emph{(little) root system} of the
symmetric space~$\HS$.
\end{lemma}
\begin{proof}
The proof is similar to that of Lemma~\ref{Kac}. The restriction
of $\alpha\in\Delta$ to $T_1$ is identified with the orthogonal
projection to $\Ch(T_1)\otimes\QQ$ given by
$\overline\alpha=(\alpha-\sigma(\alpha))/2$. We have
$\alpha+\sigma(\alpha)\notin\Delta$, $\forall\alpha\in\Delta$.
(Otherwise $\alpha+\sigma(\alpha)$ would be a non-compact root.)
The involution $\iota$ coincides with $-\sigma$ modulo the root
lattice of~$L$. One easily deduces that
$\overline\alpha=\overline\beta$ iff $\alpha=\beta$ or
$\iota(\alpha)=\beta$ for $\forall\alpha,\beta\in\Pi$ and that
$\Pi_{\HS}$ is linearly independent. Taking these remarks into
account, the proof repeats that of Lemma~\ref{Kac} with $\sigma$
replaced by~$-\sigma$.
\end{proof}

The Dynkin diagram of $\Pi$ with the ``compact'' nodes marked by
black and the remaining nodes by white, where the white nodes
transposed by $\iota$ are joined by two-headed arrows, is called
the \emph{Satake diagram} of the involution~$\sigma$, or of the
symmetric space~$\HS$. The Satake diagram encodes the embedding of
$\h$ into~$\g$. Besides, it contains information on the weight
lattice (semigroup) of the symmetric space (see
Propositions~\ref{wt.gr(symm)},~\ref{wt.sgr(symm)}).
\begin{example}
The Satake diagram of the symmetric space $\HS={H\times
H}/\diag{H}$ of Example~\ref{diag} consists of two Dynkin diagrams
of $H$, so that all nodes are white and each node of the
\ordinal{1} diagram is joined with the respective node of the
\ordinal{2} diagram, e.g.:
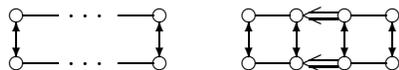
\begin{center}
%TeXCAD Picture [A0-S.pic]. Options:
%\grade{\off}
%\emlines{\off}
%\epic{\off}
%\beziermacro{\off}
%\reduce{\on}
%\snapping{\off}
%\quality{2.00}
%\graddiff{0.01}
%\snapasp{1}
%\zoom{1.00}
\unitlength .7ex % = .7pt
\linethickness{0.4pt}
\begin{picture}(15.67,6.67)(0,0)
\put(0,1){\circle{1.33}} \put(0,6){\circle{1.33}}
\put(15,1){\circle{1.33}} \put(15,6){\circle{1.33}}
\put(.67,1){\line(1,0){3.67}} \put(.67,6){\line(1,0){3.67}}
\put(7.5,1){\makebox(0,0)[cc]{$\dots$}}
\put(7.5,6){\makebox(0,0)[cc]{$\dots$}}
\put(10.67,1){\line(1,0){3.67}} \put(10.67,6){\line(1,0){3.67}}
\put(0,1.67){\line(0,1){3.67}}\put(0,5.45){\vector(0,1){0}}\put(0,1.55){\vector(0,-1){0}}
\put(15,1.67){\line(0,1){3.67}}\put(15,5.45){\vector(0,1){0}}\put(15,1.55){\vector(0,-1){0}}
\end{picture}
\qquad
%TeXCAD Picture [F0-S.pic]. Options:
%\grade{\off}
%\emlines{\off}
%\epic{\off}
%\beziermacro{\off}
%\reduce{\on}
%\snapping{\off}
%\quality{2.00}
%\graddiff{0.01}
%\snapasp{1}
%\zoom{1.00}
\unitlength .7ex % = .7pt
\linethickness{0.4pt}
\begin{picture}(15.67,6.67)(0,0)
\put(0,1){\circle{1.33}} \put(0,6){\circle{1.33}}
\put(5,1){\circle{1.33}} \put(5,6){\circle{1.33}}
\put(10,1){\circle{1.33}} \put(10,6){\circle{1.33}}
\put(15,1){\circle{1.33}} \put(15,6){\circle{1.33}}
\put(.67,1){\line(1,0){3.67}} \put(.67,6){\line(1,0){3.67}}
\put(10.67,1){\line(1,0){3.67}} \put(10.67,6){\line(1,0){3.67}}
\put(9.33,1.33){\line(-1,0){3}} \put(9.33,6.33){\line(-1,0){3}}
\put(9.33,.67){\line(-1,0){3}} \put(9.33,5.67){\line(-1,0){3}}
\put(5.33,1){\makebox(0,0)[lc]{$<$}}
\put(5.33,6){\makebox(0,0)[lc]{$<$}}
\put(0,1.67){\line(0,1){3.67}}\put(0,5.45){\vector(0,1){0}}\put(0,1.55){\vector(0,-1){0}}
\put(5,1.67){\line(0,1){3.67}}\put(5,5.45){\vector(0,1){0}}\put(5,1.55){\vector(0,-1){0}}
\put(10,1.67){\line(0,1){3.67}}\put(10,5.45){\vector(0,1){0}}\put(10,1.55){\vector(0,-1){0}}
\put(15,1.67){\line(0,1){3.67}}\put(15,5.45){\vector(0,1){0}}\put(15,1.55){\vector(0,-1){0}}
\end{picture}
\end{center}
\end{example}

The classification of symmetric spaces goes back to Cartan. To
describe it, first note that $\sigma$ preserves the connected
center and either preserves or transposes the simple factors
of~$G$. Hence every symmetric space is locally isomorphic to a
product of a torus $Z/Z\cap H$, of symmetric spaces ${H\times
H}/\diag{H}$ with $H$ simple, and of symmetric spaces of simple
groups.

Thus the classification reduces to simple~$G$. It can be obtained
using either Kac diagrams \cite[X.5]{symm.full}, \cite[Ch.3,
3.6--3.11]{LieG} or Satake diagrams \cite{class(inv)}, \cite[Ch.4,
4.1--4.3]{LieG}. For simple $G$ both Kac and Satake diagrams are
connected.

Further analysis shows that the underlying affine Dynkin diagram
for the Kac diagram of $\sigma$ depends only on the diagram
involution~$\overline\sigma$. This diagram is easily recovered
form the Dynkin diagram of $\Pi$ and from $\overline\sigma$ using
Table~\ref{Cartan}. Since the weight system of $T_0:\m$ is
symmetric, for each ``white'' root
$\overline\alpha\in\widetilde\Pi$ there exists a ``white'' root
$\overline{\alpha}_0$ and ``black'' roots
$\overline{\alpha}_1,\dots,\overline{\alpha}_r$ such that
$-\overline\alpha=
\overline{\alpha}_0+\overline{\alpha}_1+\dots+\overline{\alpha}_r$.
As $\widetilde\Pi$ is bound by a unique linear dependence, the
coefficients being positive integers, there exists either unique
``white'' root, with the coefficient $1$ or~$2$, or exactly two
``white'' roots, with the coefficients~$1$. The first possibility
occurs exactly for outer involutions, because in this case the
weight system contains the zero weight, while the other two
possibilities correspond to inner involutions. Using these
observations, it is easy to write down all possible Kac diagrams,
see Table~\ref{symm}.

On the other hand, all apriori possible Satake diagrams can also
be classified. One verifies that a Satake diagram cannot be one of
the following:
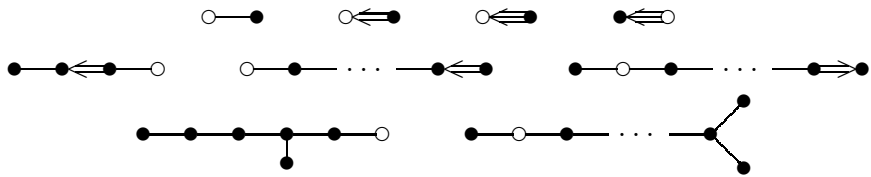
\begin{center}
%TeXCAD Picture [AV-S.pic]. Options:
%\grade{\off}
%\emlines{\off}
%\epic{\off}
%\beziermacro{\off}
%\reduce{\on}
%\snapping{\off}
%\quality{2.00}
%\graddiff{0.01}
%\snapasp{1}
%\zoom{10.0000}
\unitlength .7ex % = .7pt
\linethickness{0.4pt}
\begin{picture}(5.67,1.67)(0,0)
\put(0,1){\circle{1.33}} \put(5,1){\circle*{1.33}}
\put(.67,1){\line(1,0){3.67}}
\end{picture}
\qquad
%TeXCAD Picture [CIV-S.pic]. Options:
%\grade{\off}
%\emlines{\off}
%\epic{\off}
%\beziermacro{\off}
%\reduce{\on}
%\snapping{\off}
%\quality{2.00}
%\graddiff{0.01}
%\snapasp{1}
%\zoom{10.0000}
\unitlength .7ex % = .7pt
\linethickness{0.4pt}
\begin{picture}(5.67,1.67)(0,0)
\put(0,1){\circle{1.33}} \put(5,1){\circle*{1.33}}
\put(4.5,1.33){\line(-1,0){3.1}} \put(4.5,.67){\line(-1,0){3.1}}
\put(.33,1){\makebox(0,0)[lc]{$<$}}
\end{picture}
\qquad
%TeXCAD Picture [GI-S.pic]. Options:
%\grade{\off}
%\emlines{\off}
%\epic{\off}
%\beziermacro{\off}
%\reduce{\on}
%\snapping{\off}
%\quality{2.00}
%\graddiff{0.01}
%\snapasp{1}
%\zoom{10.0000}
\unitlength .7ex % = .7pt
\linethickness{0.4pt}
\begin{picture}(5.67,1.67)(0,0)
\put(0,1){\circle{1.33}} \put(5,1){\circle*{1.33}}
\put(.33,1){\makebox(0,0)[lc]{$<$}} \put(.67,1){\line(1,0){3.67}}
\put(5,1.5){\line(-1,0){3.33}} \put(5,.5){\line(-1,0){3.33}}
\end{picture}
\qquad
%TeXCAD Picture [GII-S.pic]. Options:
%\grade{\off}
%\emlines{\off}
%\epic{\off}
%\beziermacro{\off}
%\reduce{\on}
%\snapping{\off}
%\quality{2.00}
%\graddiff{0.01}
%\snapasp{1}
%\zoom{10.0000}
\unitlength .7ex % = .7pt
\linethickness{0.4pt}
\begin{picture}(5.67,1.67)(0,0)
\put(0,1){\circle*{1.33}} \put(5,1){\circle{1.33}}
\put(.33,1){\makebox(0,0)[lc]{$<$}} \put(.67,1){\line(1,0){3.67}}
\put(4.6,1.5){\line(-1,0){2.93}} \put(4.6,.5){\line(-1,0){2.93}}
\end{picture}
\\*[1ex]
%TeXCAD Picture [FIII-S.pic]. Options:
%\grade{\off}
%\emlines{\off}
%\epic{\off}
%\beziermacro{\off}
%\reduce{\on}
%\snapping{\off}
%\quality{2.00}
%\graddiff{0.01}
%\snapasp{1}
%\zoom{10.0000}
\unitlength .7ex % = .7pt
\linethickness{0.4pt}
\begin{picture}(15.67,1.67)(0,0)
\put(0,1){\circle*{1.33}} \put(5,1){\circle*{1.33}}
\put(10,1){\circle*{1.33}} \put(15,1){\circle{1.33}}
\put(.67,1){\line(1,0){3.67}} \put(10.67,1){\line(1,0){3.67}}
\put(9.5,1.33){\line(-1,0){3.1}} \put(9.5,.67){\line(-1,0){3.1}}
\put(5.33,1){\makebox(0,0)[lc]{$<$}}
\end{picture}
\qquad
%TeXCAD Picture [CV-S.pic]. Options:
%\grade{\off}
%\emlines{\off}
%\epic{\off}
%\beziermacro{\off}
%\reduce{\on}
%\snapping{\off}
%\quality{2.00}
%\graddiff{0.01}
%\snapasp{1}
%\zoom{10.0000}
\unitlength .7ex % = .7pt
\linethickness{0.4pt}
\begin{picture}(25.67,1.67)(0,0)
\put(0,1){\circle{1.33}} \put(5,1){\circle*{1.33}}
\put(20,1){\circle*{1.33}} \put(25,1){\circle*{1.33}}
\put(24.5,1.33){\line(-1,0){3.1}} \put(24.5,.67){\line(-1,0){3.1}}
\put(20.33,1){\makebox(0,0)[lc]{$<$}}
\put(.67,1){\line(1,0){3.67}} \put(5.67,1){\line(1,0){3.67}}
\put(15.67,1){\line(1,0){3.67}}
\put(12.5,1){\makebox(0,0)[cc]{$\dots$}}
\end{picture}
\qquad
%TeXCAD Picture [BI-S.pic]. Options:
%\grade{\off}
%\emlines{\off}
%\epic{\off}
%\beziermacro{\off}
%\reduce{\on}
%\snapping{\off}
%\quality{2.00}
%\graddiff{0.01}
%\snapasp{1}
%\zoom{10.0000}
\unitlength .7ex % = .7pt
\linethickness{0.4pt}
\begin{picture}(30.67,1.67)(0,0)
\put(5,1){\circle{1.33}} \put(10,1){\circle*{1.33}}
\put(0,1){\circle*{1.33}} \put(25,1){\circle*{1.33}}
\put(30,1){\circle*{1.33}} \put(25.57,1.33){\line(1,0){3}}
\put(25.57,.67){\line(1,0){3}}
\put(29.67,1){\makebox(0,0)[rc]{$>$}}
\put(5.67,1){\line(1,0){3.67}} \put(10.67,1){\line(1,0){3.67}}
\put(.67,1){\line(1,0){3.67}} \put(20.67,1){\line(1,0){3.67}}
\put(17.5,1){\makebox(0,0)[cc]{$\dots$}}
\end{picture}
\\*[1ex]
\raisebox{1.05ex}{%
%TeXCAD Picture [EX-S.pic]. Options:
%\grade{\off}
%\emlines{\off}
%\epic{\off}
%\beziermacro{\off}
%\reduce{\on}
%\snapping{\off}
%\quality{2.00}
%\graddiff{0.01}
%\snapasp{1}
%\zoom{10.0000}
\unitlength .7ex % = .7pt
\linethickness{0.4pt}
\begin{picture}(25.67,3.67)(0,0)
\put(0,3){\circle*{1.33}} \put(5,3){\circle*{1.33}}
\put(10,3){\circle*{1.33}} \put(15,3){\circle*{1.33}}
\put(15,0){\circle*{1.33}} \put(20,3){\circle*{1.33}}
\put(25,3){\circle{1.33}} \put(10.67,3){\line(1,0){3.67}}
\put(15.67,3){\line(1,0){3.67}} \put(20.67,3){\line(1,0){3.67}}
\put(15,2.33){\line(0,-1){1.67}} \put(5.67,3){\line(1,0){3.67}}
\put(.67,3){\line(1,0){3.67}}
\end{picture}}
\qquad
%TeXCAD Picture [DVI-S.pic]. Options:
%\grade{\off}
%\emlines{\off}
%\epic{\off}
%\beziermacro{\off}
%\reduce{\on}
%\snapping{\off}
%\quality{2.00}
%\graddiff{0.01}
%\snapasp{1}
%\zoom{10.0000}
\unitlength .7ex % = .7pt
\linethickness{0.4pt}
\begin{picture}(29.16,8.66)(0,0)
\put(5,4.5){\circle{1.33}} \put(10,4.5){\circle*{1.33}}
\put(0,4.5){\circle*{1.33}} \put(25,4.5){\circle*{1.33}}
\put(5.67,4.5){\line(1,0){3.67}} \put(10.67,4.5){\line(1,0){3.67}}
\put(.67,4.5){\line(1,0){3.67}} \put(20.67,4.5){\line(1,0){3.67}}
\put(17.5,4.5){\makebox(0,0)[cc]{$\dots$}}
\put(28.5,1){\circle*{1.33}} \put(28.5,8){\circle*{1.33}}
\multiput(25.4,4.1)(.116087,-.116087){23}{\line(1,0){.116087}}
\multiput(25.4,4.9)(.116087,.116087){23}{\line(1,0){.116087}}
\end{picture}
\end{center}
In all cases except the last two, the sum of all simple roots
would be a complex root $\alpha$ such that
$\alpha+\sigma(\alpha)\in\Delta$, a contradiction. In the
remaining two cases, $\sigma$~would be an inner involution
represented by an element $s\in S=Z_G(L_0)^0$. The group $S$ is a
simple $\SL_2$-subgroup corresponding to the highest root
$\delta\in\Delta$ and $T_1$ is a maximal torus in~$S$. Replacing
$T_1$ by another maximal torus containing $s$ one obtains
$\delta(s)=-1$. However the unique $\alpha\in\Pi$ such that
$\alpha(s)=-1$ occurs in the decomposition of $\delta$ with an
even coefficient, a contradiction.

By a \emph{fragment} of a Satake diagram we mean a $\iota$-stable
subdiagram such that no one of its nodes is joined with a black
node outside the fragment. A fragment is the Satake diagram of a
Levi subgroup in $G$. It follows that a Satake diagram cannot
contain the above listed fragments. Also, if a Satake diagram
contains a fragment
%TeXCAD Picture [A.pic]. Options:
%\grade{\off}
%\emlines{\off}
%\epic{\off}
%\beziermacro{\off}
%\reduce{\on}
%\snapping{\off}
%\quality{2.00}
%\graddiff{0.01}
%\snapasp{1}
%\zoom{10.0000}
\unitlength .7ex % = .7pt
\linethickness{0.4pt}
\begin{picture}(15.66,1.67)(0,0)
\put(0,1){\circle*{1.33}} \put(15,1){\circle*{1.33}}
\put(.67,1){\line(1,0){3.67}} \put(10.67,1){\line(1,0){3.67}}
\put(7.5,1){\makebox(0,0)[cc]{$\dots$}}
\end{picture}
of length $>1$, then there are no other black nodes and $\iota$ is
non-trivial. Having this in mind, it is easy to write down all
possible Satake diagrams, see Table~\ref{symm}.

Both Kac and Satake diagrams uniquely determine the
involution~$\sigma$. All apriori possible diagrams are realized
for simply connected~$G$. It follows that symmetric spaces of
simple groups are classified, up to a local isomorphism, by Kac or
Satake diagrams.

The classification is presented in Table~\ref{symm}. The column
``$\sigma$'' describes the involution for classical~$G$ in matrix
terms. Here
%*
\begin{equation*}
I_{n,m}=

\end{equation*}
%*
is the matrix of a standard symplectic form fixed
by~$\Sp_{2n}(\kk)$, where $E_k$ is the unit $(k\times k)$-matrix.

\begin{table}[!h]
\caption{Symmetric spaces of simple groups} \label{symm}
\begin{center}
\scriptsize\makebox[0pt]{\begin{tabular}{|c|c|c|c|c|c|c|} \hline
\No & $G$ & $H$ & $\sigma$ & Kac diagram & Satake diagram
& $\Delta_{G/H}$ \\
\hline

\end{tabular}}
\end{center}
\end{table}

\begin{example}
Let us describe the symmetric spaces of $G=\SL_n(\kk)$. Take the
standard Borel subgroup of upper-triangular matrices $B\subset G$
and the standard diagonal torus $T\subset B$. By
$\eps_1,\dots,\eps_n$ denote the weights of the tautological
representation in~$\kk^n$ (i.e.,~the diagonal entries of~$T$).

If $\sigma$ is inner, then $\overline\Delta=\Delta$ and the Dynkin
diagram of $\widetilde\Pi$ is the following one:
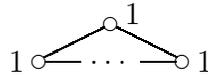
\begin{center}
%TeXCAD Picture [A0.pic]. Options:
%\grade{\off}
%\emlines{\off}
%\epic{\off}
%\beziermacro{\off}
%\reduce{\on}
%\snapping{\off}
%\quality{2.00}
%\graddiff{0.01}
%\snapasp{1}
%\zoom{10.0000}
\unitlength .7ex % = .7pt
\linethickness{0.4pt}
\begin{picture}(18,5)(0,0)
\put(1.5,1){\circle{1.33}} \put(16.5,1){\circle{1.33}}
\put(2.17,1){\line(1,0){3.67}} \put(12.17,1){\line(1,0){3.67}}
\put(9,1){\makebox(0,0)[cc]{$\dots$}} \put(9,5){\circle{1.33}}
\multiput(2,1.5)(.237037,.118519){27}{\line(1,0){.237037}}
\multiput(16,1.5)(-.237037,.118519){27}{\line(-1,0){.237037}}
\put(18,1){\makebox(0,0)[lc]{$1$}}
\put(0,1){\makebox(0,0)[r]{$1$}}
\put(10.5,5){\makebox(0,0)[lb]{$1$}}
\end{picture}
\end{center}
The coefficients of the unique linear dependence on
$\widetilde\Pi$ are indicated at the diagram. It follows that
there are exactly two white nodes in the Kac diagram. The
involution $\iota$ is non-trivial, whence there is at most one
black fragment in the Satake diagram, which is located in the
middle. Thus we obtain \No\,\ref{SL/GLSL} of Table~\ref{symm}.

The involution $\sigma$ is the conjugation by an element of order
$2$ in~$\GL_n(\kk)$. In a certain basis, $\sigma(g)=I_{n,m}\cdot
g\cdot I_{n,m}$. Then $T_0=T$, $H=\SG(\LO_m\times\LO_{n-m})$ is
embedded in $G$ by the two diagonal blocks, the simple roots being
$\eps_i-\eps_{i+1}$, $1\leq i<n$, $i\neq m$, and
$\m=\kk^m\otimes(\kk^{n-m})^{*} \oplus(\kk^m)^{*}\otimes\kk^{n-m}$
is embedded in $\g$ by the two antidiagonal blocks, the lowest
weights of the summands being $\eps_m-\eps_{m+1}$,
$\eps_n-\eps_1$, in accordance with the Kac diagram.

In another basis, $\sigma(g)=J_{n,m}\cdot g\cdot J_{n,m}$, where
%*
\begin{equation*}
J_{n,m}=
\begin{array}{r@{}|@{}c@{}@{}c@{}@{}c@{}|@{}l}
\multicolumn{4}{@{}r@{}} {\overbrace{\hphantom{
\begin{array}{|c@{}c@{}c}
\text{\large\strut$0$}   & & 1 \\[-1ex]
  &                \Dots++ &   \\[-1.5ex]
1 & & \smash{\text{\large$0$}} \\
\end{array}}}^m} & \\[-0.5ex]
\cline{2-4} & \smash{\raisebox{-1em}{\text{\huge$0$}}} &&
\begin{array}{|c@{}c@{}c}
\text{\large\strut$0$}   & & 1 \\[-1ex]
  &                \Dots++ &   \\[-1.5ex]
1 & & \smash{\text{\large$0$}} \\
\end{array} &
\left.
\begin{array}{@{}c@{}}
\text{\strut}\\[-1ex]
\\[-1.5ex]
\\
\end{array}
\right\}\scriptstyle{m} \\
\cline{3-4} &&
\begin{array}{|c@{}c@{}c|}
1 & & \text{\large\strut$0$}   \\[-1ex]
  & \Dots+-                &   \\[-1.5ex]
\smash{\text{\large$0$}} & & 1 \\
\end{array}
&& \\
\cline{2-3} \scriptstyle{m}\left\{
\begin{array}{@{}c@{}}
\text{\strut}\\[-1ex]
\\[-1.5ex]
\\
\end{array}
\right. &
\begin{array}{c@{}c@{}c|}
\text{\large\strut$0$} &   & 1 \\[-1ex]
  &                \Dots++ &   \\[-1.5ex]
1 & & \smash{\text{\large$0$}} \\
\end{array}
&& \smash{\text{\huge$0$}} & \\
\cline{2-4} \multicolumn{1}{@{}r@{}}{} & \multicolumn{4}{@{}l@{}}
{\raisebox{1.5ex}{$\underbrace{\hphantom{
\begin{array}{c@{}c@{}c|}
\text{\large\strut$0$}   & & 1 \\[-1ex]
  &                \Dots++ &   \\[-1.5ex]
1 & & \smash{\text{\large$0$}} \\
\end{array}}}_m$}}
\end{array}
\end{equation*}
%*
Now $T_1=
\{t=\diag(t_1,\dots,t_m,1,\dots,1,t_m^{-1},\dots,t_1^{-1})\}$ is a
maximal $\sigma$-split torus and the (compact) imaginary roots are
$\eps_i-\eps_j$, $m<i\neq j\leq n-m$, in accordance with the
Satake diagram. The little root system $\Delta_{\HS}$ consists of
the nonzero restrictions $\overline\eps_i-\overline\eps_j$, $1\leq
i,j\leq n$, i.e.,\ of $\pm\overline\eps_i\pm\overline\eps_j$,
$\pm2\overline\eps_i$, and $\pm\overline\eps_i$ unless $m=n/2$,
$1\leq i\neq j\leq m$. Thus $\Delta_{\HS}$ is of type $\BBb\CCc_m$
or $\CCc_{n/2}$.

If $\sigma$ is outer, then $\overline\sigma(\eps_i)=-\eps_{n+1-i}$
and $(T^{\overline\sigma})^0=
\{t=\diag(t_1,t_2,\dots,t_2^{-1},t_1^{-1})\}$. Restricting the
roots to this subtorus, we see that $\overline\Delta$ consists of
$\pm\eps_i'\pm\eps_j'$, $\pm2\eps_i'$, and $\pm\eps_i'$ for
odd~$n$, where $\eps_i'$ are the restrictions of~$\eps_i$, $1\leq
i\leq n/2$. The Dynkin diagram of $\widetilde\Pi$ has one of the
following forms:
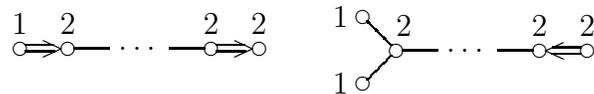
\begin{center}
\raisebox{1.9ex}{%
%TeXCAD Picture [Ao.pic]. Options:
%\grade{\off}
%\emlines{\off}
%\epic{\off}
%\beziermacro{\off}
%\reduce{\on}
%\snapping{\off}
%\quality{2.00}
%\graddiff{0.01}
%\snapasp{1}
%\zoom{10.0000}
\unitlength .7ex % = .7pt
\linethickness{0.4pt}
\begin{picture}(25.66,2.5)(0,0)
\put(0,1){\circle{1.33}} \put(5,1){\circle{1.33}}
\put(20,1){\circle{1.33}} \put(25,1){\circle{1.33}}
\put(20.5,1.33){\line(1,0){3}} \put(.5,1.33){\line(1,0){3}}
\put(20.5,.67){\line(1,0){3.1}} \put(.5,.67){\line(1,0){3.1}}
\put(24.67,1){\makebox(0,0)[rc]{$>$}}
\put(4.67,1){\makebox(0,0)[rc]{$>$}}
\put(5.67,1){\line(1,0){3.67}} \put(15.67,1){\line(1,0){3.67}}
\put(12.5,1){\makebox(0,0)[cc]{$\dots$}}
\put(0,2.5){\makebox(0,0)[cb]{$1$}}
\put(5,2.5){\makebox(0,0)[cb]{$2$}}
\put(20,2.5){\makebox(0,0)[cb]{$2$}}
\put(25,2.5){\makebox(0,0)[cb]{$2$}}
\end{picture}}
\qquad
%TeXCAD Picture [Ae.pic]. Options:
%\grade{\off}
%\emlines{\off}
%\epic{\off}
%\beziermacro{\off}
%\reduce{\on}
%\snapping{\off}
%\quality{2.00}
%\graddiff{0.01}
%\snapasp{1}
%\zoom{10.0000}
\unitlength .7ex % = .7pt
\linethickness{0.4pt}
\begin{picture}(25.67,7.67)(0,0)
\put(5,3.5){\circle{1.33}} \put(1.5,7){\circle{1.33}}
\put(1.5,0){\circle{1.33}} \put(20,3.5){\circle{1.33}}
\put(25,3.5){\circle{1.33}} \put(24.4,3.83){\line(-1,0){3}}
\put(24.4,3.17){\line(-1,0){3.1}}
\put(20.33,3.5){\makebox(0,0)[lc]{$<$}}
\put(5.67,3.5){\line(1,0){3.67}} \put(15.67,3.5){\line(1,0){3.67}}
\put(12.5,3.5){\makebox(0,0)[cc]{$\dots$}}
\multiput(4.5,4)(-.119048,.119048){22}{\line(0,1){.119048}}
\multiput(4.5,3)(-.119048,-.119048){22}{\line(0,-1){.119048}}
\put(0,0){\makebox(0,0)[rc]{$1$}}
\put(0,7){\makebox(0,0)[rc]{$1$}}
\put(5,5){\makebox(0,0)[lb]{$2$}}
\put(20,5){\makebox(0,0)[cb]{$2$}}
\put(25,5){\makebox(0,0)[cb]{$2$}}
\end{picture}
\end{center}
depending on whether $n$ is odd or even. Therefore the Kac diagram
has a unique white node, namely an extreme one.

The involution $\iota$ is trivial, whence either all nodes of the
Satake diagram are white or the black nodes are isolated from each
other and alternate with the white ones, the extreme nodes being
black. (Otherwise, there would exist an inadmissible fragment
%TeXCAD Picture [AV-S.pic]. Options:
%\grade{\off}
%\emlines{\off}
%\epic{\off}
%\beziermacro{\off}
%\reduce{\on}
%\snapping{\off}
%\quality{2.00}
%\graddiff{0.01}
%\snapasp{1}
%\zoom{10.0000}
\unitlength .7ex % = .7pt
\linethickness{0.4pt}
\begin{picture}(5.67,1.67)(0,0)
\put(0,1){\circle{1.33}} \put(5,1){\circle*{1.33}}
\put(.67,1){\line(1,0){3.67}}
\end{picture}.)
Thus we obtain \Nos\,\ref{SL/Sp}--\ref{SL/SO} of Table~\ref{symm}.

Any outer involution has the form $\sigma(g)=(g^{*})^{-1}$, where
${*}$ denotes the conjugation w.r.t.\ a non-degenerate
(skew-)symmetric bilinear form on~$\kk^n$. In the symmetric case,
choosing an orthonormal basis yields $\sigma(g)=(g^{\tran})^{-1}$,
whence $T_1=T$ is a maximal $\sigma$-split torus and
$\Delta_{\HS}=\Delta$. In a hyperbolic basis,
$\sigma(g)=(g^{\stran})^{-1}$, where $\stran$ denotes the
transposition w.r.t.\ the secondary diagonal. Then
$T_0=(T^{\overline\sigma})^0$ is a maximal torus in
$H=\SO_n(\kk)$. The space $\m$ consists of traceless symmetric
matrices, and the lowest weight is $-2\eps_1'$.

In the skew-symmetric case, choosing an appropriately ordered
symplectic basis yields
$\sigma(g)=I_{n,n/2}(g^{\stran})^{-1}I_{n,n/2}$. Here $T_0$ is a
maximal torus in $H=\Sp_n(\kk)$ and
$T_1=\{t=\diag(t_1,t_2,\dots,t_2,t_1)\mid t_1\dots,t_{n/2}=1\}$ is
a maximal $\sigma$-split torus. The roots of $H$ are
$\pm\eps_i'\pm\eps_j'$, $\pm2\eps_i'$, $1\leq i\neq j\leq n/2$,
and the lowest weight of $\m$ is $-\eps_1'-\eps_2'$. The compact
roots are $\eps_i-\eps_{n+1-i}$ ($1\leq i\leq n$), and
$\Delta_{\HS}$ consists of $\overline\eps_i-\overline\eps_j$,
$1\leq i\neq j\leq n/2$, thus having the type~$\AAa_{n/2-1}$.
\end{example}

From now on we assume that $T_1$ is a maximal $\sigma$-split
torus.

Consider the Weyl group $W_{\HS}$ of the little root
system~$\Delta_{\HS}$.
\begin{proposition}
$W_{\HS}\iso N_{H^0}(T_1)/Z_{H^0}(T_1)\iso N_G(T_1)/Z_G(T_1)$
\end{proposition}
\begin{proof}
First we prove that each element of $W_{\HS}$ is induced by an
element of $N_{H^0}(T_1)$. It suffices to consider a root
reflection~$r_{\overline\alpha}$. Let
$T_1^{\overline\alpha}\subseteq T_1$ be the connected kernel
of~$\overline\alpha$. Replacing $G$ by
$Z_G(T_1^{\overline\alpha})$ we may assume that
$W_{\HS}=\{\1,r_{\overline\alpha}\}$. The same argument as in
Lemma~\ref{conj(split)} shows that $P^{-}=\sigma(P)=hPh^{-1}$ for
some $h\in H^0$. It follows that $h\in N_{H^0}(L)=N_{H^0}(T_1)$
acts on $\Ch(T_1)$ as~$r_{\overline\alpha}$.

On the other hand, $N_G(T_1)$~acts on $T_1$ as a subgroup of the
``big'' Weyl group $W=N_G(T)/T$. Indeed, any $g\in N_G(T_1)$
normalizes $L=Z_G(T_1)$ and may be replaced by another element in
$gL$ normalizing~$T$. Since the Weyl chambers of $W_{\HS}$ in
$\Ch(T_1)\otimes\QQ$ are the intersections of Weyl chambers of $W$
with $\Ch(T_1)\otimes\QQ$, the orbits of $N_G(T_1)/Z_G(T_1)$
intersect them in single points. Thus $N_G(T_1)/Z_G(T_1)$ cannot
be bigger than~$W_{\HS}$. This concludes the proof.
\end{proof}

Since $\HS$ is spherical, there are finitely many $B$-orbits
in~$\HS$ (Corollary~\ref{B-fin}). Their structure plays an
important role in some geometric problems and, for $\kk=\CC$, in
the representation theory of the real reductive Lie group $G(\RR)$
acting on the Riemannian symmetric space~$\HS(\RR)$, the
non-compact real form of~$\HS$~\cite{Vog}. The classification and
the adherence relation for $B$-orbits were described in
\cite{inv}, \cite{RS}, \cite{RS.survey} (cf.~Example~\ref{RS}). We
explain the basic classification result under the assumption
$H=G^{\sigma}$. This is not an essential restriction
\cite[1.1(b)]{RS.survey}.

By Proposition~\ref{symm<G}, $\HS$~is identified with~$\tau(G)$,
where $G$ (and~$B$) acts by twisted conjugation.
\begin{proposition}\label{B&T-orb}
The (twisted) $B$-orbits in $\tau(G)\iso\HS$ intersect $N_G(T)$ in
$T$-orbits. Thus $\B(\HS)$ is in bijective correspondence with the
set of twisted $T$-orbits in $N(T)\cap\tau(G)$.
\end{proposition}
\begin{proof}
Consider a $B$-orbit $Bg\bp\subseteq\HS$. By
Lemma~\ref{sigma-stab}, replacing $g$ by $bg$, $b\in B$, one may
assume that $g^{-1}Tg$ is a $\sigma$-stable maximal torus in
$g^{-1}Bg$. This holds iff $\tau(g)\in N(T)$. One the other hand,
taking another point $g'\bp\in Bg\bp$, $g'=bgh$, $b\in B$, $h\in
H$, we have $\tau(g')=\sigma(b)\tau(g)b^{-1}\in N(T)$ iff
$\tau(g')=\sigma(t)\tau(g)t^{-1}$, where $b=tu$, $t\in T$, $u\in
U$, by standard properties of the Bruhat decomposition
\cite[28.4]{Agr}.
\end{proof}

There is a natural map $\B(\HS)\to W$, $Bg\bp\mapsto w$, where
$\sigma(B)wB$ is the unique Bruhat cell containing the respective
$B$-orbit $\tau(BgH)$. By Proposition~\ref{B&T-orb},
$\tau(BgH)\cap N(T)\subseteq wT$. This map plays an important role
in the study of $B$-orbits \cite{RS}, \cite{RS.survey}. Its image
is contained in the set of \emph{twisted involutions} $\{w\in
W\mid\sigma(w)=w^{-1}\}$, but in general is neither injective nor
surjective onto this set.

\begin{example}
Let $G=\GL_n(\kk)$, $\sigma(g)=(g^{\tran})^{-1}$, $H=\Or_n(\kk)$.
Then $\tau(G)$ is the set of non-degenerate symmetric matrices,
viewed as quadratic forms on~$\kk^n$. The group $B$ of
upper-triangular matrices acts on $\tau(G)$ by base changes
preserving the standard flag in~$\kk^n$. It is an easy exercise in
linear algebra that for any inner product on $\kk^n$ one can
choose a basis $e_1,\dots,e_n$ compatible with a given flag and
having the property that for any $i$ there is a unique $j$ such
that $(e_i,e_j)=1$ and $(e_i,e_k)=0$, $\forall k\neq j$. The
matrix of the quadratic form in this basis is the permutation
matrix of the involution transposing $i$ and~$j$. It lies
in~$N(T)$ (where $T$ is the diagonal torus) and is uniquely
determined by the $B$-orbit of the quadratic form. Thus $\B(\HS)$
is in bijection with the set of involutions in $W=\Sym_n$.
\end{example}

Now we describe the colored equipment of a symmetric space,
according to~\cite{emb(symm)}.

The weight lattice of a symmetric space is read off the Satake
diagram, at least up to an finite extension. Let $Z=Z(G)^0$ and
$\omega_i$ be the fundamental weights corresponding to the simple
roots $\alpha_i\in\Pi$.
\begin{proposition}\label{wt.gr(symm)}
If $G$ is of simply connected type, then
%*
\begin{equation}\label{wts<=Satake}
\RG(\HS)=\Ch(Z/Z\cap H)\oplus \left\langle\,\widehat\omega_j,\
\omega_k+\omega_{\iota(k)}\bigm| j,k\,\right\rangle
\end{equation}
%*
where $j,k$ run over all $\iota$-fixed, resp.\ $\iota$-unstable,
white nodes of the Satake diagram, and $\widehat\omega_j=\omega_j$
or~$2\omega_j$, depending on whether the \ordinal{$j$} node is
adjacent to a black one or not. In the general case, $\RG(\HS)$~is
a sublattice of finite index in the r.h.s.\
of~\eqref{wts<=Satake}.
\end{proposition}
\begin{remark}\label{2(fund.wts)}
The weight lattice $\RG(\HS)=\Ch(T/T\cap H)=\Ch(T_1/T_1\cap H)$
injects into $\Ch(T_1)$ via restriction of characters from $T$
to~$T_1$. The space $\ES=\Hom(\RG(\HS),\QQ)$ is then identified
with $\CoCh(T_1)\otimes\QQ$. The second direct summand in the
r.h.s.\ of \eqref{wts<=Satake} is nothing else but the doubled
weight lattice $2(\ZZ\Delta_{\HS}^{\vee})^{*}$ of the little root
system~$\Delta_{\HS}$. Indeed, $\widehat\omega_j/2$ and
$(\omega_k+\omega_{\iota(k)})/2$ restrict to the fundamental
weights dual to the simple coroots
$\overline\alpha^{\vee}_j=\alpha^{\vee}_j-\sigma(\alpha^{\vee}_j)$
or $\alpha^{\vee}_j$ and
$\overline\alpha^{\vee}_k=\alpha^{\vee}_k-\sigma(\alpha^{\vee}_k)$.
\end{remark}
\begin{proof}
W.l.o.g.\ we may assume that $G$ is semisimple simply connected,
whence $H=G^{\sigma}$. The sublattice $\RG(\HS)\subseteq\Ch(T)$
consisting of the weights vanishing on~$T^{\sigma}$, i.e., of
$\mu-\sigma(\mu)$, $\mu\in\Ch(T)$, is contained in $\Ch(T/T_0)=
\{\lambda\in\Ch(T)\mid\sigma(\lambda)=-\lambda\}=
\langle\omega_j,\ \omega_k+\omega_{\iota(k)}\mid j,k\rangle$. The
latter lattice injects into $\Ch(T_1)$ so that $\RG(\HS)$ is
identified with~$2\Ch(T_1)$. It remains to prove that
$\Ch(T_1)=(\ZZ\Delta_{\HS}^{\vee})^{*}$ or, equivalently, that
$\CoCh(T_1)=\ZZ\Delta_{\HS}^{\vee}$ is the coroot lattice of the
little root system.

We have $\CoCh(T)=\ZZ\Delta^{\vee}$ and
$\CoCh(T_1)=\ZZ\Delta^{\vee}\cap\ES\supseteq\ZZ\Delta_{\HS}^{\vee}$.
The \emph{alcoves} (=fundamental polyhedra, see \cite[IV,
\S2]{Lie:4-6}) of the affine Weyl group
$W_{\text{aff}}(\Delta_{\HS}^{\vee})$ are the intersections of
$\ES$ with alcoves of $W_{\text{aff}}(\Delta_{\HS}^{\vee})$. Hence
each alcove of $W_{\text{aff}}(\Delta_{\HS}^{\vee})$ contains a
unique point from~$\CoCh(T_1)$. It follows that $\CoCh(T_1)$
coincides with $\ZZ\Delta_{\HS}^{\vee}$.
\end{proof}

Let $\DWC=\DWC(\Delta^{+})$ denote the dominant Weyl chamber of a
root system~$\Delta$ (w.r.t.~a chosen subset of positive
roots~$\Delta^{+}$). The weight semigroup $\RG_{+}(\HS)$ is
contained both in $\RG(\HS)$ and in~$\DWC(\Delta^{+})$. Note that
$\DWC(\Delta^{+})\cap\ES=\DWC(\Delta_{\HS}^{+})$.
\begin{proposition}\label{wt.sgr(symm)}
$\RG_{+}(\HS)=\RG(\HS)\cap\DWC(\Delta_{\HS}^{+})$
\end{proposition}
\begin{proof}
Since $\RG_{+}(\HS)$ is the semigroup of all lattice points in a
cone (see~\ref{c=0}), it suffices to prove that
$\QQ_{+}\RG_{+}(\HS)=\DWC(\Delta_{\HS}^{+})$. Take any dominant
$\lambda\in\RG(\HS)$. We prove that $2\lambda\in\RG_{+}(\HS)$.

First note that $\lambda=-\sigma(\lambda)$ is orthogonal to
compact roots, whence $\lambda$ is extended to $P$ and
$V^{*}(\lambda)=\Ind^G_P\kk_{-\lambda}$. Consider another dual
Weyl module obtained by twisting the $G$-action by~$\sigma$:
$V^{*}(\lambda)^{\sigma}=\Ind^G_{\sigma(P)}\kk_{-\sigma(\lambda)}
\iso V^{*}(\lambda^{*})$. We have the canonical $H$-equivariant
linear isomorphism $\omega:V^{*}(\lambda)^{\sigma}\isoto
V^{*}(\lambda)$. (If the dual Weyl modules are realized in
$\kk[G]$ as in Example~\ref{Ind}, then $\omega$ is just the
restriction of $\sigma$ acting on~$\kk[G]$.) In other words,
$\omega\in(V^{*}(\lambda)\otimes V(\lambda^{*}))^H$. Note that
$\omega$ maps a $T$-eigenvector of weight $\mu$ to an eigenvector
of weight $\sigma(\mu)$. Hence
%*
\begin{equation}\label{H-fixed}
\omega=v_{-\lambda}\otimes v'_{-\lambda}+
\sum_{\mu\neq\lambda}v_{\sigma(\mu)}\otimes v'_{-\mu}
\end{equation}
%*
where $v_{\chi},v'_{\chi}$ denote basic eigenvectors of weight
$\chi$ in $V^{*}(\lambda)$ and $V(\lambda^{*})$, respectively.
Applying the homomorphisms $V(\lambda^{*})\to V^{*}(\lambda)$,
$v'_{-\lambda}\mapsto v_{-\lambda}$, and $V^{*}(\lambda)\otimes
V^{*}(\lambda)\to V^{*}(2\lambda)$ (induced by multiplication
in~$\kk[G]$), we obtain a nonzero element $\overline\omega\in
V^{*}(2\lambda)^H$, whence $2\lambda\in\RG_{+}(\HS)$
by~\eqref{mult(G/H)}.
\end{proof}

Now we are ready to describe the colors and $G$-valuations of a
symmetric space.
\begin{theorem}\label{c.d.(symm)}
The $B$-divisors of a symmetric space $\HS$ are represented by the
vectors from $\frac12\Pi_{\HS}^{\vee}\subset\ES$ (where
$\Pi_{\HS}^{\vee}$ is the base of
$\Delta_{\HS}^{\vee}\subset\CoCh(T_1)$). The valuation cone $\Vv$
is the antidominant Weyl chamber of $\Delta_{\HS}^{\vee}$
in~$\ES$.
\end{theorem}
\begin{corollary}
$W_{\HS}$~is the little Weyl group of $\HS$ in the sense
of~\ref{Weyl}.
\begin{question}
and $\Delta_{\HS}^{\max}$ is the reduced root system
of~$2\Delta_{\HS}$.
\end{question}
\end{corollary}
\begin{proof}
W.l.o.g.\ $G$ is assumed to be of simply connected type. In the
notation of Remarks~\ref{D^B(G/H)},~\ref{D^B(sph)}, each
$f\in\kk[\HS]^{(B)}_{\lambda}$ is represented as
$f=\eta_1^{d_1}\dots\eta_s^{d_s}$, where the $\eta_i$ are
equations of the $B$-divisors $D_i\in\Dd^B$, $d_i\in\ZZ_{+}$, and
$\lambda=\sum d_i\lambda_i$, $\sum d_i\chi_i=0$, where
$(\lambda_i,\chi_i)$ are the biweights of~$\eta_i$.

In the notation of Proposition~\ref{wt.gr(symm)}, if
$\lambda=\widehat\omega_j$ or $\omega_k+\omega_{\iota(k)}$, then
$f=\eta_j$, or $\eta_j'\eta_j''$, or $\eta_k$, or
$\eta_k'\eta_k''$, where the biweights of
$\eta_j,\eta_j',\eta_j'',\eta_k,\eta_k',\eta_k''$ are
$(\widehat\omega_j,0)$, $(\omega_j,\chi_j)$, $(\omega_j,-\chi_j)$,
$(\omega_k+\omega_{\iota(k)},0)$, $(\omega_k,\chi_k)$,
$(\omega_{\iota(k)},-\chi_k)$, respectively, for some nonzero
$\chi_j,\chi_k\in\Ch(H)$. In particular, the respective
$B$-divisors $D_j,D_j',D_j'',D_k,D_k',D_k''$ are pairwise
distinct, and all $B$-divisors occur among them since these $f$'s
span the multiplicative semigroup
$\kk[\HS]^{(B)}/\kk[\HS]^{\times}$ by
Proposition~\ref{wt.sgr(symm)} and Remark~\ref{2(fund.wts)}. The
assertion on colors stems now from Remarks~\ref{D^B(sph)}
and~\ref{2(fund.wts)}.

Now we treat $G$-valuations. Take any $v=v_D\in\Vv$, where $D$ is
a $G$-stable prime divisor on a $G$-model $X$ of~$\kk(\HS)$. It
follows from the local structure theorem that
$F=\overline{T_1\bp}$ is an $N_H(T_1)$-stable subvariety of $X$
intersecting $D$ in the union of $T_1$-stable prime divisors
$D_{wv}$, $w\in W_{\HS}$, that correspond to $wv$ regarded as
$T_1$-valuations of~$\kk(T_1\bp)$
(cf.~Proposition~\ref{flat&div}). By Theorem~\ref{cent.cone}
$\Vv$~contains the antidominant Weyl chamber. It remains to show
as in the proof of Theorem~\ref{W_X&Z(X)} that different vectors
from $\Vv$ cannot be $W_{\HS}$-equivalent.
\end{proof}

The proof of Theorem~\ref{c.d.(symm)} shows that the map
$\res:\Dd^B\to\ES$ may be non-injective if $H$ is not semisimple.
There is a more precise description of $B$-divisors in the spirit
of Proposition~\ref{B&T-orb} \cite[5.4]{inv},
\cite[\S4]{wonder(symm)}.

It suffices to consider simple $G$.  Assume first that $H$ is
connected. For any $\overline\alpha^{\vee}\in\Pi_{\HS}^{\vee}$
there exist either a unique or exactly two $B$-divisors mapping to
$\overline\alpha^{\vee}/2$. They correspond to the twisted
$T$-orbits in $\tau(G)\cap r_{\alpha}T$ (for real~$\alpha$) or in
$\tau(G)\cap\left(r_{\sigma(\alpha)}r_{\alpha}T\cup
r_{\sigma(\iota(\alpha))}r_{\iota(\alpha)}T\right)$ (for
complex~$\alpha$).

If $H$ is semisimple, then such an orbit (and the respective
$B$-divisor~$D_{\alpha}$) is always unique. In particular,
$\iota(\alpha)=\alpha$ or
$\iota(\alpha)=-\sigma(\alpha)\perp\alpha$.

If $H$ is not semisimple (\emph{Hermitian case}), then inspection
of Table~\ref{symm} shows that $\dim Z(H)=1$ and $\Delta_{\HS}$ is
of type $\BBb\CCc_n$ or~$\CCc_n$. The $B$-divisor mapped to
$\overline\alpha^{\vee}/2$ is unique except for the case where
$\overline\alpha^{\vee}$ is the short simple coroot.

In the latter case, if $\alpha$ is complex, then $\tau(G)\cap
r_{\sigma(\alpha)}r_{\alpha}T$ and $\tau(G)\cap
r_{\sigma(\iota(\alpha))}r_{\iota(\alpha)}T$ are the twisted
$T$-orbits corresponding to the two $\Aut_G\HS$-stable
$B$-divisors $D_{\alpha},D_{\iota(\alpha)}$ mapped
to~$\overline\alpha^{\vee}/2$. Here $\Delta_{\HS}=\BBb\CCc_n$ and
$c(G/H')=0$.

If $\alpha$ is real, then $\tau(G)\cap r_{\alpha}T$ consists of
two twisted $T$-orbits corresponding to the two $B$-divisors
$D_{\alpha}',D_{\alpha}''$ mapped to~$\alpha^{\vee}/2$ and swapped
by~$\Aut_G\HS$. Here $\Delta_{\HS}=\CCc_n$ and $c(G/H')=1$.

For disconnected $H$ the divisors
$D_{\alpha}',D_{\alpha}''\in\Dd(G/H^0)^B$ may patch together into
a single divisor $D_{\alpha}\in\Dd(G/H)^B$.

The (co)isotropy representation $H:\m$ has nice
invariant-theoretic properties in characteristic zero. They were
examined by Kostant and Rallis~\cite{iso(symm)}. From now on
assume $\ch\kk=0$.

Semisimple elements in $\m$ are exactly those having closed
$H$-orbits, and the unique closed $H$-orbit in $\overline{H\xi}$
($\xi\in\m$) is~$H\xi_{\sms}$. Generic elements of $\m$ are
semisimple. One may deduce it from the fact that $T^{*}\HS$ is
symplectically stable (Proposition~\ref{qaff=>ss}) or prove
directly:
$\g=\lv\oplus[\g,\tr_1]\implies\m=\tr_1\oplus[\h,\tr_1]\implies
\m=\overline{H\tr_1}$. This argument also shows that $H$-invariant
functions on $\m$ are uniquely determined by their restrictions
to~$\tr_1$. A more precise result was obtained by Kostant and
Rallis.
\begin{proposition}[\cite{iso(symm)}]\label{iso(symm)}
Every semisimple $H$-orbit in $\m$ intersects $\tr_1$ in a
$W_{\HS}$-orbit. Restriction of functions yields an isomorphism
$\kk[\m]^H\iso\kk[\tr_1]^{W_{\HS}}$.
\end{proposition}
\begin{proof}
Every semisimple element $\xi\in\m$ is contained in the Lie
algebra of a maximal $\sigma$-split torus, hence by
Lemma~\ref{conj(split)}, $\xi'=(\Ad h)\xi\in\tr_1$ for some $h\in
H$. If $\xi\in\tr_1$, then $T_1,h^{-1}T_1h$ are two maximal
$\sigma$-split tori in~$Z_G(\xi)$. Again by
Lemma~\ref{conj(split)}, $zT_1z^{-1}=h^{-1}T_1h$ for some $z\in
Z_G(\xi)\cap H$, whence $h'=hz\in N_H(\tr_1)$, $\xi'=(\Ad
h')\xi\in W_{\HS}\xi$.

The second assertion is a particular case of
Proposition~\ref{coiso(sph)}. It suffices to observe that the
surjective birational morphism $\m\by{H}\to\tr_1/W_{\HS}$ of two
normal affine varieties has to be an isomorphism.
\end{proof}

Global analogues of these results for the $H$-action on~$\HS$ (in
any characteristic) were obtained by Richardson~\cite{inv(symm)}.

It is not by chance that the description of the valuation cone of
a symmetric space was obtained by the same reasoning as
in~\ref{inv.mot}.
\begin{proposition}[{\cite[\S6]{inv.mot}}]\label{flats(symm)}
Flats in $\HS$ are exactly the $G$-translates of~$T_1\cdot\bp$.
\end{proposition}
\begin{proof}
It suffices to consider flats $F_{\alpha}$, $\alpha\in
T^{\pr}_{\bp}\HS$. We have $T^{*}\HS=G\itimes{H}\m$,
$\alpha=\1*\xi$, $\xi=\Phi(\alpha)\in\m^{\pr}$. By
Proposition~\ref{iso(symm)}, $\xi\in(\Ad h)\tr_1^{\pr}$, $h\in
H^0$. It follows that $G_{\xi}=hLh^{-1}$, whence
$F_{\alpha}=hL\bp=hT_1\bp$.
\end{proof}

The $W_{\HS}$-action on the flat $T_1\cdot\bp$ comes
from~$N_H(T_1)$.

In the case $\kk=\CC$, flats in $\HS$ are ($G$-translates of) the
complexifications of maximal totally geodesic flat submanifolds in
a Riemannian symmetric space $\HS(\RR)$ which is a real form
of~$\HS$ \cite{symm.full}.

\section{Algebraic monoids and group embeddings}
\label{monoids}

Alike algebraic groups, defined by superposing the concepts of an
abstract group and an algebraic variety, it is quite natural to
consider \emph{algebraic semigroups}, i.e., algebraic varieties
equipped with an associative multiplication law which is a regular
map.
\begin{example}
All linear operators on a finite-dimensional vector space $V$ form
an algebraic semigroup $\LO(V)\iso\LO_n(\kk)$ ($n=\dim V$). The
operators (matrices) of rank $\leq r$ form a closed subsemigroup
$\LO^{(r)}(V)$ ($\LO^{(r)}_n(\kk)$), a particular example of a
determinantal variety.
\end{example}

However the category of \emph{all} algebraic semigroups is
immense. (For instance, every algebraic variety $X$ turns into an
algebraic semigroup being equipped with the ``zero''
multiplication $X\times X\to\{0\}$, where $0\in X$ is a fixed
element.) In order to make the theory really substantive, one has
to restrict the attention to algebraic semigroups not too far from
algebraic groups.
\begin{definition}
An \emph{algebraic monoid} is an algebraic semigroup with unit,
i.e., an algebraic variety $X$ equipped with a morphism
$\mu:X\times X\to X$, $\mu(x,y)=:x\cdot y$ (the
\emph{multiplication law}) and with a distinguished \emph{unity
element} $\1\in X$ such that $(x\cdot y)\cdot z=x\cdot(y\cdot z)$,
$e\cdot x=x\cdot e=x$, $\forall x,y,z\in X$.
\end{definition}

Let $G=G(X)$ denote the group of invertible elements in~$X$. The
following elementary result can be found, e.g.,
in~\cite[\S2]{monoids}.
\begin{proposition}
$G$~is open in~$X$.
\end{proposition}
\begin{proof}
Since the left translation $x\mapsto g\cdot x$ by an element $g\in
G$ is an automorphism of~$X$, it suffices to prove that $G$
contains an open subset of an irreducible component of~$X$.
W.l.o.g.~we may assume that $X$ is irreducible. Let $p_1,p_2$ be
the two projections of $\mu^{-1}(\1)\subset X\times X$ to~$X$. By
the fiber dimension theorem, every component of $\mu^{-1}(\1)$ has
dimension $\geq\dim X$, and $p_i^{-1}(\1)=(\1,\1)$. Hence $p_i$
are dominant maps and $G=p_1(\mu^{-1}(\1))\cap p_2(\mu^{-1}(\1))$
is a dense constructible set containing an open subset of~$X$.
\end{proof}
\begin{corollary}
$G$~is an algebraic group.
\end{corollary}

Those irreducible components of $X$ which do not intersect $G$ do
not ``feel the presence'' of $G$ and their behavior is beyond of
control. Therefore it is reasonable to restrict oneself to
algebraic monoids $X$ such that $G=G(X)$ is dense in~$X$. In this
case, left translations by $G$ permute the components of $X$
transitively and many questions are reduced to the case, where $X$
is irreducible.

Monoids of this kind form an interesting category of algebraic
structures closely related to algebraic groups (e.g., they arise
as the closures of linear algebraic groups in the spaces of linear
operators). The theory of algebraic monoids was created in major
part during the last 25 years by M.~S.~Putcha, L.~E.~Renner,
E.~B.~Vinberg, A.~Rittatore, et al. The interested reader may
consult a detailed survey \cite{monoids.survey} of the theory from
the origin up to latest developments. In this section, we discuss
algebraic monoids from the viewpoint of equivariant embeddings. A
link between these two theories is provided by the following
result.
\begin{theorem}[{\cite[\S2]{monoids}}]\label{mon<=>emb}
\begin{enumerate}
\item\label{mon=>emb} Any algebraic monoid $X$ is a ${G\times
G}$-equivariant embedding of $G=G(X)$, where the factors of
$G\times G$ act by left/right multiplication, having a unique
closed ${G\times G}$-orbit. \item\label{emb=>mon} Conversely, any
\emph{affine} ${G\times G}$-equivariant embedding $X\embof G$
carries a structure of algebraic monoid with $G(X)=G$.
\end{enumerate}
\end{theorem}
\begin{proof}
\begin{roster}
\item[\ref{mon=>emb}] One has only to prove the uniqueness of a
closed orbit $Y\subseteq X$. Note that $X\cdot Y\cdot X=
\overline{G\cdot Y\cdot G}=Y$, i.e., $Y$~is a (two-sided)
\emph{ideal} in~$X$. For any other ideal $Y'\subseteq X$ we have
$Y\cdot Y'\subseteq Y\implies Y=Y\cdot Y'\subseteq Y'$. Thus $Y$
is the smallest ideal, called the \emph{kernel} of~$X$.
\item[\ref{emb=>mon}] The actions of the left and right copy of
$G\times G$ on $X$ define coactions $\kk[X]\to\kk[G]\otimes\kk[X]$
and $\kk[X]\to\kk[X]\otimes\kk[G]$, which are the restrictions to
$\kk[X]\subseteq\kk[G]$ of the comultiplication
$\kk[G]\to\kk[G]\otimes\kk[G]$. Hence the image of $\kk[X]$ lies
in $(\kk[G]\otimes\kk[X])\cap(\kk[X]\otimes\kk[G])=
\kk[X]\otimes\kk[X]$, and we have a comultiplication in~$\kk[X]$.
Now $G$ is open in $X=\overline{G}$ and consists of invertibles.
For any invertible $x\in X$, we have $xG\cap G\ne\emptyset$, hence
$x\in G$. \qedhere\end{roster}
\end{proof}
\begin{remark}
Assertion~\ref{emb=>mon} was first proved for reductive $G$ by
Vinberg \cite{red.sgr} in a different way.
\end{remark}

Among general algebraic groups, affine (=linear) ones occupy a
privileged position due to their most rich and interesting
structure. The same holds for algebraic monoids. We provide two
results confirming this observation.
\begin{theorem}[{\cite[\S4]{ab.var}}]
Complete irreducible algebraic monoids are just Abelian varieties.
\end{theorem}
\begin{theorem}[\cite{mon=>aff}]\label{mon=>aff}
An algebraic monoid $X$ is affine provided that $G(X)$ is affine.
\end{theorem}
This theorem was proved by Renner \cite{qaff.mon=>aff} for
quasiaffine $X$ using some structure theory and Rittatore
\cite{mon=>aff} reduced the general case to the quasiaffine one by
considering total spaces of certain line bundles over~$X$.

A theorem of Barsotti \cite{grp.var} and Rosenlicht
\cite{Agr.basic} says that every algebraic group has a unique
affine normal subgroup such that the quotient group is an Abelian
variety. It is an interesting unsolved problem to obtain an
analogous structure result for algebraic monoids.

\begin{theorem}\label{aff=>lin}
Any affine algebraic monoid $X$ admits a closed homomorphic
embedding $X\embeds\LO(V)$. Furthermore, $G(X)=X\cap\GL(V)$.
\end{theorem}
The proof is essentially the same as that of a similar result for
algebraic groups \cite[8.6]{Agr}. Thus the adjectives ``affine''
and ``linear'' are synonyms for algebraic monoids, alike for
algebraic groups.

In the notation of Theorem~\ref{aff=>lin}, the space of matrix
entries $\Mat(V)$ generates $\kk[X]\subseteq\kk[G]$. Generally,
$\kk[X]\supset\Mat(V)$ iff the representation $G:V$ is extendible
to~$X$. It follows from Theorem~\ref{aff=>lin} and \eqref{mat}
that
%*
\begin{equation}\label{k[mon]}
\kk[X]=\bigcup\Mat(V)
\end{equation}
%*
over all $G$-modules $V$ that are $X$-modules
(cf.~Proposition~\ref{U(Mat)}).

\begin{example}
By Theorem~\ref{mon<=>emb}\ref{emb=>mon}, every affine toric
variety $X$ carries a natural structure of algebraic monoid
extending the multiplication in the open torus~$T$. By
Theorem~\ref{aff=>lin}, $X$~is the closure of $T$ in $\LO(V)$ for
some faithful representation $T:V$, i.e., a closed submonoid in
the monoid of all diagonal matrices in some~$\LO_n(\kk)$. The
coordinate algebra $\kk[X]$ is the semigroup algebra of the
semigroup $\Sigma\subseteq\Ch(T)$ consisting of all characters
$T\to\kk^{\times}$ extendible to~$X$. Conversely, every finitely
generated semigroup $\Sigma\ni0$ such that $\ZZ\Sigma=\Ch(T)$
defines a toric monoid $X\supseteq T$.
\end{example}

The classification and structure theory for algebraic monoids is
most well developed in the case, where the group of invertibles is
reductive.
\begin{definition}
An irreducible algebraic monoid $X$ is called \emph{reductive} if
$G=G(X)$ is a connected reductive group.
\end{definition}
In the sequel we consider only reductive monoids, thus returning
to the general convention of our survey that $G$ is a connected
reductive group. By Theorems~\ref{mon<=>emb},~\ref{mon=>aff},
reductive monoids are nothing else but ${G\times G}$-equivariant
affine embeddings of~$G$. They were classified by Vinberg
\cite{red.sgr} in characteristic zero. Rittatore \cite{monoids}
extended this classification to arbitrary characteristic using the
embedding theory of spherical homogeneous spaces.

Considered as a homogeneous space under $G\times G$ acting by
left/right multiplication, $G$~is a symmetric space
(Example~\ref{diag}). All $\sigma$-stable maximal tori of $G\times
G$ are of the form $T\times T$, where $T$ is a maximal torus
in~$G$. The maximal $\sigma$-split tori are $({T\times
T})_1=\{(t^{-1},t)\mid t\in T\}$. Choose a Borel subgroup
$B\supseteq T$ of~$G$.  Then $B^{-}\times B$ is a Borel subgroup
in $G\times G$ containing $T\times T$ and $\sigma({B^{-}\times
B})=B\times B^{-}$ is the opposite Borel subgroup.

The weight lattice $\RG=\Ch({T\times
T}/\diag{T})=\{(-\lambda,\lambda)\mid \lambda\in\Ch(T)\}$ is
identified with $\Ch(T)$ and the little root system
with~$\frac12\Delta$. The eigenfunctions
$\ef{\lambda}\in\kk(G)^{(B^{-}\times B)}$ ($\lambda\in\Ch(T)$) are
defined on the ``big'' open cell $U^{-}\times T\times U\subseteq
G$ by the formula $\ef{\lambda}(u^{-}tu)=\lambda(t)$. For
$\lambda\in\Ch_{+}$ they are matrix entries:
$\ef{\lambda}(g)=\langle v_{-\lambda},gv_{\lambda}\rangle$, where
$v_{\lambda}\in V$, $v_{-\lambda}\in V^{*}$ are highest,
resp.~lowest, vectors of weights~$\pm\lambda$.

By Theorem~\ref{c.d.(symm)}, the valuation cone $\Vv$ is
identified with the antidominant Weyl chamber in
$\CoCh(T)\otimes\QQ$ (this can also be deduced from
Example~\ref{cent.val(G)}) and the colors are represented by the
simple coroots
$\alpha_1^{\vee},\dots,\alpha_l^{\vee}\in\Pi^{\vee}$. The
respective $B$-divisors are $D_i=\overline{B^{-}r_{\alpha_i}B}$.
Indeed, the equation of $D_i$ in $\kk[\widetilde{G}]$
is~$\ef{\omega_i}$, where $\omega_i$ denote the fundamental
weights.

The results of \ref{form.curv} (in particular, the Cartan
decomposition) imply that all ${G\times G}$-valuations are
proportional to $v=v_{\gamma}$, $\gamma\in\CoCh(T)$.  Since
$v_{w\gamma}=v_{\gamma}$, $\forall w\in W=N_G(T)/T$,
w.l.o.g.~$\gamma\in\Vv$. Then a direct computation shows
$v(\ef{\lambda})=\langle\gamma,\lambda\rangle$,
$\forall\lambda\in\Ch_{+}$, whence $v$ is identified with~$\gamma$
(as a vector in the valuation cone).

Now Corollary~\ref{sph.aff} yields
\begin{theorem}\label{norm.mon}
Normal reductive monoids $X$ are in bijection with strictly convex
cones $\Cc=\Cc(X)\subset\CoCh(T)\otimes\QQ$ generated by all
simple coroots and finitely many antidominant vectors.
\end{theorem}
\begin{remark}
The normality assumption is not so restrictive, because the
multiplication on $X$ lifts to its normalization $\widetilde{X}$
turning it into a monoid with the same group of invertibles.
\end{remark}
\begin{corollary}
There are no non-trivial monoids with semisimple group of
invertibles.
\end{corollary}
\begin{corollary}[{\cite{Put}, \cite[Pr.9]{monoids}}]
\label{mon->0} Every normal reductive monoid has the structure
$X=(X_0\times G_1)/Z$, where $X_0$ is a monoid with zero, and $Z$
is a finite central subgroup in $G(X_0)\times G_1$ not
intersecting the factors.
\end{corollary}
\begin{proof}
Identify $\Ch(T)\otimes\QQ$ with $\ES=\CoCh(T)\otimes\QQ$ via a
$W$-invariant inner product. Consider an orthogonal decomposition
$\ES=\ES_0\oplus\ES_1$, where $\ES_0=\langle\Cc\cap\Vv\rangle$,
$\ES_1=(\Cc\cap\Vv)^{\perp}$. It is easy to see that each root is
contained in one of the~$\ES_i$. Then $G=G_0\cdot G_1=(G_0\times
G_1)/Z$, where $G_i$ are the connected normal subgroups with
$\CoCh(T\cap G_i) =\ES_i\cap\CoCh(T)$. Take a reductive monoid
$X_0\supseteq G_0$ defined by $\Cc_0=\Cc\cap\ES_0$. Since
$\intr\Cc_0$ intersects $\Vv(G_0)=\Vv\cap\ES_0$, the kernel of
$X_0$ is a complete variety, hence a single point~$0$, the zero
element w.r.t.\ the multiplication on~$X_0$. Now $X$ coincides
with $(X_0\times G_1)/Z$, because both monoids have the same
colored data.
\end{proof}

This classification can be made more transparent via coordinate
algebras and representations. Recall from \ref{c=0} that
$\kk[X]^{U^{-}\times U}=\kk[\Cc^{\vee}\cap\Ch(T)]$. The algebra
$\kk[X]$ itself is given by~\eqref{k[mon]}. It remains to
determine which representations of $G$ extend to~$X$.
\begin{proposition}\label{rep(mon)}
The following conditions are equivalent:
\begin{enumerate}
\item\label{ext} The representation $G:V$ is extendible to $X$.
\item\label{ht-wts} The highest weights of the simple factors of
$V$ are in~$\Cc^{\vee}$. \item\label{T-wts} All dominant
$T$-weights of $V$ are in~$\Cc^{\vee}$.
\end{enumerate}
\end{proposition}
\begin{proof}
\begin{roster}
\item[\ref{ext}$\implies$\ref{ht-wts}] Choose a $G$-stable
filtration of $V$ with simple factors and consider the associated
graded $G$-module~$\gr V$. If $G:V$ extends to~$X$, then $\gr V$
is an $X$-module. Hence $\ef{\lambda}\in\kk[X]$ whenever $\lambda$
is a highest weight of~$\gr V$.
\item[\ref{ht-wts}$\iff$\ref{T-wts}] All $T$-weights of $V$ are
obtained from the highest weights of simple factors by subtracting
positive roots. The structure of $\Cc$ implies that all dominant
vectors obtained this way from $\lambda\in\Cc^{\vee}$ belong
to~$\Cc^{\vee}$.

\item[\ref{T-wts}$\implies$\ref{ext}] Assume that
$\kk[X]\not\supset\Mat(V)$. Choose $f\in\Mat(V)$ representing a
nonzero ${B^{-}\times B}$-eigenvector $\mod\kk[X]$. Then by
Corollary~\ref{(M/N)^H}, $f^q=\ef{\lambda}\mod\kk[X]$ for some
$\lambda\notin\Cc^{\vee}$. It follows that $\lambda/q$ is a
$T$-weight of $V$ outside~$\Cc^{\vee}$. \qedhere\end{roster}
\end{proof}
\begin{corollary}\label{C^}
If $X\subseteq\LO(V)$, then $\Cc^{\vee}=\Kk(V)\cap\DWC$, where
$\Kk(V)$ denotes the convex cone spanned by the $T$-weights
of~$V$.
\end{corollary}
\begin{proof}
The proposition implies $\Cc^{\vee}\supseteq\Kk(V)\cap\DWC$. On
the other hand, all ${T\times T}$-weights of $\kk[X]$ are of the
form $(-\lambda,\mu)$, $\lambda,\mu\in\Kk(V)$, whence
$\Cc^{\vee}\subseteq\Kk(V)$.
\end{proof}

In characteristic zero, Proposition~\ref{rep(mon)} together with
\eqref{k[mon]} yields
%*
\begin{equation}\label{k[red.norm]}
\kk[X]=\bigoplus_{\lambda\in\Cc^{\vee}\cap\Ch(T)} \Mat(V(\lambda))
\end{equation}
%*
(cf.~Theorem~\ref{reg.rep} and~\eqref{k[G]}). In positive
characteristic, $\kk[X]$~has a ``good'' filtration with factors
$V^{*}(\lambda)\otimes V^{*}(\lambda^{*})$ \cite{Dot},
\cite[\S4]{F-split.red}, \cite[Cor.9.9]{monoids.survey}.

The embedding theory provides a combinatorial encoding for
${G\times G}$-orbits in~$X$, which reflects the adherence
relation. This description can be made more explicit using the
following
\begin{proposition}\label{cl(T)}
Suppose $X\embof G$ is an equivariant embedding. Then
$F=\overline{T}$ intersects each ${G\times G}$-orbit $Y\subset X$
in finitely many $T$-orbits permuted transitively by~$W$. Exactly
one of these orbits $F_Y\subseteq F\cap Y$ satisfies
$\intr\Cc_{F_Y}\cap\Vv\neq\emptyset$; then
$\Cc_{F_Y}=W({\Cc_Y\cap\Vv})\cap\Cc_Y$.
\begin{question}
Perhaps $\Cc_Y\cap\Vv\subseteq\Cc_{F_Y}\subseteq
W({\Cc_Y\cap\Vv})$?
\end{question}
\end{proposition}
\begin{remark}
Since $T$ is a flat of~$G$ (Proposition~\ref{flats(symm)}), some
of the assertions stem from the results of~\ref{inv.mot}. However,
the proposition here is more precise. In particular, it completely
determines the fan of~$F$.
\end{remark}
\begin{proof}
Take any $v\in\Ss_Y$; then $v=v_{\gamma}$,
$\gamma\in\CoCh(T)\cap\Vv$,
$\exists\lim_{t\to0}\gamma(t)=\gamma(0)\in Y$.
\begin{question}
Use toroidal resolution instead?
\end{question}
The associated parabolic subgroup $P=P(\gamma)$ contains~$B^{-}$.
Consider the Levi decomposition $P=L\Ru{P}$, $L\supseteq T$. One
verifies that $({G\times G})_{\gamma(0)}\supseteq
(\Ru{P^{-}}\times\Ru{P})\cdot\diag{L}$. It easily follows that
$({B^{-}\times B})\gamma(0)=\Y$ is the open ${B^{-}\times
B}$-orbit in $Y$ and $F_Y:=T\gamma(0)=\Y^{\diag{T}}$ is the unique
$T$-orbit in $F$ intersecting~$\Y$.

In view of Example~\ref{val<->1par}, this implies
$\intr\Cc_{F_Y}\supseteq(\intr\Cc_Y)\cap\Vv$. On the other hand,
each $T$-orbit in $F\cap Y$ is accessed by a one-parameter
subgroup $\gamma\in\CoCh(T)$, $\gamma(0)\in Y$. Taking $w\in W$
such that $w\gamma\in\Vv$ yields $w(T\gamma(0))=F_Y$. All
assertions of the proposition are deduced from these observations.
\end{proof}

Now suppose $X\subseteq\LO(V)$ and denote $\Kk=\Kk(V)$.
\begin{theorem}\label{orb(mon)}
The ${G\times G}$-orbits in $X$ are in bijection with the faces of
$\Kk$ whose interiors intersect~$\DWC$. The orbit $Y$
corresponding to a face $\Ff$ is represented by the
$T$-equivariant projector $\1_{\Ff}$ of $V$ onto the sum of
$T$-eigenspaces of weights in~$\Ff$. The cone $\Cc_Y$ is dual to
the corner cone of $\Kk\cap\DWC$ at the face $\Ff\cap\DWC$, and
$\Dd^B_Y$ consists of simple coroots orthogonal to~$\Ff$.
\end{theorem}
\begin{proof}
A complete set of $T$-orbit representatives in $F=\overline{T}$ is
formed by the limits of one-parameter subgroups, i.e., by the
$\1_{\Ff}$'s over \emph{all} faces $\Ff$ of~$\Kk$. The respective
cones in the fan of $F$ are the dual faces
$\Ff^{*}=\Kk^{\vee}\cap\Ff^{\ann}$ of $\Kk^{\vee}=W(\Cc\cap\Vv)$.
By Proposition~\ref{cl(T)}, the orbits $Y$ are bijectively
represented by those $\1_{\Ff}$ which satisfy
$\intr\Ff^{*}\cap\Vv\neq\emptyset$. This happens iff $\Ff^{*}$
lies on a face of $\Cc$ of the same dimension (namely on~$\Cc_Y$)
or, equivalently, $\Ff$~contains a face of
$\Cc^{\vee}=\Kk\cap\DWC$ of the same dimension
(namely~$\Cc_Y^{*}=\Ff\cap\DWC$), i.e.,
$\intr\Ff\cap\DWC\neq\emptyset$. The assertion on $(\Cc_Y,\Dd_Y)$
stems from the description of a dual face.
\end{proof}

\begin{example}\label{Mat(n)}
Let $G=\GL_n(\kk)$ and $X=\LO_n(\kk)$. For~$B$ and $T$ take the
standard Borel subgroup of upper-triangular matrices and diagonal
torus, respectively. We have
$\Ch(T)=\langle\eps_1,\dots,\eps_n\rangle$, where the $\eps_i$ are
the diagonal matrix entries of~$T$. We identify $\Ch(T)$ with
$\CoCh(T)$ via the inner product such that the $\eps_i$ form an
orthonormal basis. Let $(k_1,\dots,k_n)$ denote the coordinates on
$\Ch(T)\otimes\QQ$ w.r.t.\ this basis. The Weyl group $W=S_n$
permutes them.

The weights $\lambda_i=\eps_1+\dots+\eps_i$ span~$\Ch(T)$ and
$\ef{\lambda_i}\in\kk[X]$ are the upper-left corner $i$-minors of
a matrix. Put $D_i=\{x\in X\mid\ef{\lambda_i}(x)=0\}$. Then
$\Dd^B=\{D_1,\dots,D_{n-1}\}$, $D_i$~are represented by
$\alpha_i=\eps_i-\eps_{i+1}$, $\forall i<n$, and $D_n$ is the
unique $G$-stable prime divisor, $v_{D_n}=\eps_n$.

Therefore $\Cc=\{k_1+\dots+k_i\geq0,\ i=1,\dots,n\}$ is the cone
spanned by ${\eps_i-\eps_{i+1}},\eps_n$, and
$\Cc^{\vee}=\{k_1\geq\dots\geq k_n\geq0\}$ is spanned
by~$\lambda_i$. The lattice vectors of $\Cc^{\vee}$ are exactly
the dominant weights of polynomial representations
(cf.~Proposition~\ref{rep(mon)}).  The lattice vectors of
$\Kk=W\Cc^{\vee}=\{k_1,\dots,k_n\geq0\}$ are all polynomial
weights of~$T$.

The ${G\times G}$-orbits in $X$ are $Y_r=\{x\in X\mid\rk x=r\}$.
Clearly, $\Dd^B_{Y_r}=\{D_i\mid r<i<n\}$ and $\Cc_{Y_r}$ is a face
of $\Cc$ cut off by the equations $k_1=\dots=k_r=0$. The dual face
$\Cc_{Y_r}^{*}$ of $\Cc^{\vee}$ is the dominant part of the face
$\Ff_r=\{k_i\geq0=k_j\mid i\leq r<j\}\subseteq\Kk$, and all faces
of $\Kk$ whose interiors intersect $\DWC=\{k_1\geq\dots\geq k_n\}$
are obtained this way. Clearly, the respective projectors
$\1_{\Ff_r}=\diag(1,\dots,1,0,\dots,0)$ are the ${G\times
G}$-orbit representatives, and the representatives of all
$T$-orbits in $\overline{T}$ are obtained from $\1_{\Ff_r}$ by the
$W$-action.
\end{example}

In characteristic zero, it is possible to classify (to a certain
extent) arbitrary (not necessarily normal) reductive
monoids~\cite{red.sgr} via their coordinate algebras
alike~\eqref{k[red.norm]}. The question is to describe finitely
generated ${G\times G}$-stable subalgebras of $\kk[G]$ with the
quotient field~$\kk(G)$. They are of the form
%*
\begin{equation}\label{k[red]}
\kk[X]=\bigoplus_{\lambda\in\Sigma}\Mat(V(\lambda))
\end{equation}
%*
where $\Sigma$ is a finitely generated subsemigroup of $\Ch_{+}$
such that $\ZZ\Sigma=\Ch(T)$ and the r.h.s.\ of \eqref{k[red]}
remains closed under multiplication, i.e., all highest weights of
$V(\lambda)\otimes V(\mu)$ belong to $\Sigma$ whenever
$\lambda,\mu\in\Sigma$. Such a semigroup $\Sigma$ is called
\emph{perfect}.

\begin{definition}
We say that $\lambda_1,\dots,\lambda_m$ \emph{$G$-generate}
$\Sigma$ if $\Sigma$ consists of all highest weights of
$G$-modules $V(\lambda_1)^{\otimes k_1}\otimes\dots\otimes
V(\lambda_m)^{\otimes k_m}$, $k_1,\dots,k_m\in\ZZ_{+}$. (In
particular any generating set $G$-generates~$\Sigma$.)
\end{definition}
\begin{example}
In Example~\ref{Mat(n)}, $\Sigma=\Cc^{\vee}\cap\Ch(T)$ is
generated by $\lambda_1,\dots,\lambda_n$ and $G$-generated
by~$\lambda_1$.
\end{example}
It is easy to see that $X\embeds\LO(V)$ iff the highest weights
$\lambda_1,\dots,\lambda_m$ of $G:V$ $G$-generate~$\Sigma$. The
highest weight theory implies that $\Kk=\Kk(V)$ is the $W$-span of
%*
\begin{equation}\label{wt.cone}
\Kk\cap\DWC=
(\QQ_{+}\{\lambda_1,\dots,\lambda_m,-\alpha_1,\dots,-\alpha_l\})
\cap\DWC
\end{equation}
%*
Theorem~\ref{orb(mon)} generalizes to this context.

By Theorem~\ref{U-inv}\ref{X<=>X/U}, $X$~is normal iff
$\kk[X]^{U^{-}\times U}=\kk[\Sigma]$ is integrally closed iff
$\Sigma$ is the semigroup of all lattice vectors in a polyhedral
cone. In general, taking the integral closure yields
%*
\begin{equation*}
\QQ_{+}\Sigma=\Cc^{\vee}=\Kk\cap\DWC
\end{equation*}
%*
where $\Cc=\Cc(\widetilde{X})$. Here is a representation-theoretic
interpretation: each dominant vector in $\Kk$ eventually occurs as
a highest weight in a tensor power of~$V$, see
{\cite[\S2]{grp.comp}} for a direct proof.

Given~$G:V$, the above normality condition for $X\subseteq\LO(V)$
is generally not easy to verify, because the reconstruction of
$\Sigma$ from $\{\lambda_1,\dots,\lambda_m\}$ requires decomposing
tensor products of arbitrary $G$-modules. Of course, there is no
problem if $\lambda_i$ already generate $\Kk\cap\Ch_{+}$---a
sufficient condition for normality. Here is an effective necessary
condition:
\begin{proposition}[{\cite{Ren}, \cite[Th.5.4(b)]{monoids.survey}}]
If $X$ is normal, then $F=\overline{T}$ is normal, i.e., the
$T$-weights of $V$ generate $\Kk\cap\Ch(T)$.
\end{proposition}
\begin{proof}
We can increase $V$ by adding new highest weights $\lambda_i$ so
that $\lambda_1,\dots,\lambda_m$ will generate
$\Sigma=\Kk\cap\Ch_{+}$. (This operation does not change $X$
and~$F$.) Then $W\{\lambda_1,\dots,\lambda_m\}$ generates
$\Kk\cap\Ch(T)$, i.e., $\kk[F]=\kk[\Kk\cap\Ch(T)]$ is integrally
closed.
\end{proof}

If $V=V(\lambda)$ is irreducible, then the center of $G$ acts by
homotheties, whence $G=\kk^{\times}\cdot G_0$, where $G_0$ is
semisimple, $\Ch(T)\subseteq\ZZ\oplus\Ch({T\cap G_0})$ is a
cofinite sublattice, and $\lambda=(1,\lambda_0)$. Recently de
Concini showed that $\Kk(V(\lambda))\cap\Ch_{+}$ is $G$-generated
by the $T$-dominant weights of~$V(\lambda)$ \cite{norm.sgr}.
However $\Sigma$ contains no $T$-weights of $V(\lambda)$
except~$\lambda$. It follows that $X$ is normal iff $\lambda_0$ is
a minuscule weight for~$G_0$ \cite{norm.sgr},
\cite[\S12]{grp.comp}.

It turns out that Example~\ref{Mat(n)} is essentially the unique
non-trivial example of a smooth reductive monoid.
\begin{theorem}[{cf.~\cite{Ren}, \cite[\S11]{grp.comp}}]
Smooth reductive monoids are of the form
$X=(G_0\times\LO_{n_1}(\kk)\times\dots\times\LO_{n_s}(\kk))/Z$,
where $Z\subset G_0\times\GL_{n_1}(\kk)\times\dots\times\GL_{n_s}$
is a finite central subgroup not intersecting
$\GL_{n_1}(\kk)\times\dots\times\GL_{n_s}$.
\end{theorem}
\begin{proof}
By Corollary~\ref{mon->0}, $X=G_0\itimes{Z}X_0$, where $X_0$ has
the zero element. Thus it suffices to consider monoids with zero.
We explain how to handle this case in characteristic zero.
\begin{question}
Prove in arbitrary characteristic.
\end{question}

Assume $X\subseteq\LO(V)$. There exists a coweight
$\gamma\in\intr\Cc\cap\Vv$, $\gamma\perp\Delta$. It defines a
one-parameter subgroup $\gamma(t)\in Z(G)$ contracting $V$ to~$0$
(as $t\to0$). The algebra $\Aa=\Aa(V)$ spanned by $X$ in $\LO(V)$
is semisimple, i.e., a product of matrix algebras, and $T_0X$ is
an ideal in~$\Aa$. As $X$ is smooth and the multiplication by
$\gamma(t)$ contracts $X$ to~$0$, the equivariant projection $X\to
T_0X$ is an isomorphism.
\end{proof}

We conclude this section by a discussion of arbitrary (not
necessarily affine) equivariant embeddings of~$G$. For simplicity,
we assume $\ch\kk=0$.

In the same way as a faithful linear representation $G:V$ defines
a reductive monoid $\overline{G}\subseteq\LO(V)$, a faithful
projective representation $G:\PP(V)$ (arising from a linear
representation of a finite cover of $G$ in~$V$) defines a
projective completion $X=\overline{G}\subseteq\PP(\LO(V))$. These
group completions are studied in~\cite{grp.comp}. There are two
main tool to reduce their study to reductive monoids.

First, the cone $\widehat{X}\subseteq\LO(V)$ over $X$ is a
reductive monoid whose group of invertibles $\widehat{G}$ is the
extension of $G$ by homotheties. Conversely, any such monoid gives
rise to a projective completion. This allows to transfer some of
the above results to projective group completions. For instance,
Theorem~\ref{orb(mon)} transfers verbatim if we only replace the
weight cone $\Kk(V)$ by the weight polytope $\Pp=\Pp(V)$ (=the
convex hull of the $T$-weights of~$V$), see \cite[\S9]{grp.comp}
for details.

Another approach, suitable for local study, is to use the local
structure theorem. By the above, closed $(G\times G)$-orbits
$Y\subset X$ correspond to the dominant vertices $\lambda\in\Pp$,
and the representatives are $y=\langle v_{\lambda}\otimes
v_{-\lambda}\rangle$, where $v_{\lambda}\in V$, $v_{-\lambda}\in
V^{*}$ are highest, resp.~lowest, vectors of weights~$\pm\lambda$,
$\langle v_{\lambda},v_{-\lambda}\rangle\neq0$. Consider the
parabolic $P=P(\lambda)$ and its Levi decomposition $P=L\Ru{P}$,
$L\supseteq T$. Then $V_0=\langle v_{-\lambda}\rangle^{\ann}$ is
an $L$-stable complement to $\langle v_{\lambda}\rangle$ in~$V$.
Put $\X=X\setminus\Zeros{\ef{\lambda}}$.
\begin{lemma}
$\X\iso\Ru{P^{-}}\times Z\times\Ru{P}$, where
$Z\iso\overline{L}\subseteq\LO(V_0\otimes\kk_{-\lambda})$ is a
reductive monoid with the zero element~$y$.
\end{lemma}
\begin{proof}
Applying Lemma~\ref{loc.str.V} to ${G\times G}:\LO(V)=V\otimes
V^{*}$, passing to projectivization and intersecting with~$X$, we
obtain a neighborhood of the desired structure with
$Z=X\cap\PP\bigl(\kk^{\times}(v_{\lambda}\otimes
v_{-\lambda})+E_0\bigr)$, where
%*
\begin{equation*}
E_0=(\g\times\g)(v_{-\lambda}\otimes v_{\lambda})^{\ann}=(\g
v_{-\lambda}\otimes v_{\lambda}+v_{-\lambda}\otimes\g
v_{\lambda})^{\ann} \supseteq V_0\otimes V_0^{*}=\LO(V_0)
\end{equation*}
%*
Hence $Z=\overline{L}\subseteq\PP\bigl(\kk^{\times}
(v_{\lambda}\otimes v_{-\lambda})\oplus\LO(V_0)\bigr)
\iso\LO(V_0\otimes\kk_{-\lambda})$.
\end{proof}

The monoids $Z$ are transversal slices to the closed orbits
in~$X$. They can be used to study the local geometry of~$X$. For
instance, one can derive criteria for normality and smoothness
\cite[\S\S10,11]{grp.comp}.

\begin{example}
Take $G=\Sp_4(\kk)$, with the simple roots
$\alpha_1=\eps_1-\eps_2$, $\alpha_2=2\eps_2$, and the fundamental
weights $\omega_1=\eps_1$, $\omega_2=\eps_1+\eps_2$,
$\pm\eps_i$~being the weights of the tautological representation
$\Sp_4(\kk):\kk^4$. Let $\lambda_1=3\omega_1$,
$\lambda_2=2\omega_2$ be the highest weights of~$V$. The weight
polytope $\Pp$ is depicted in Figure~\ref{Sp(4)}(a), the highest
weights are indicated by bold dots.
\begin{figure}[!h]
\caption{A projective completion of $\Sp_4(\kk)$}\label{Sp(4)}
\begin{center}
\begin{tabular}{c@{\qquad}c}
%TexCad Options
%\grade{\off}
%\emlines{\off}
%\beziermacro{\off}
%\reduce{\on}
%\snapping{\off}
%\quality{2.00}
%\graddiff{0.01}
%\snapasp{1}
%\zoom{2.00}
\unitlength 0.40ex \linethickness{0.4pt}
\begin{picture}(80,65)
\put(70.00,60.00){\line(-1,-1){60.00}}
\put(0.00,30.00){\line(1,0){80}}
\put(40.00,30.00){\vector(0,1){20.00}}
\put(40.00,30.00){\vector(1,-1){10.00}}
\put(70.00,30.00){\line(-1,2){10.00}}
\put(40.00,60.00){\line(2,-1){20.00}}
\put(10.00,30.00){\line(1,2){10.00}}
\put(40.00,60.00){\line(-2,-1){20.00}}
\put(70.00,30.00){\line(-1,-2){10.00}}
\put(40.00,0.00){\line(2,1){20.00}}
\put(10.00,30.00){\line(1,-2){10.00}}
\put(40.00,0.00){\line(-2,1){20.00}}
\put(40.00,30.00){\circle*{1.00}}
\put(50.00,30.00){\circle*{1.00}}
\put(60.00,30.00){\circle*{1.00}}
\put(30.00,30.00){\circle*{1.00}}
\put(20.00,30.00){\circle*{1.00}}
\put(10.00,30.00){\circle*{1.00}}
\put(40.00,40.00){\circle*{1.00}}
\put(40.00,50.00){\circle*{1.00}}
\put(40.00,60.00){\circle*{1.00}}
\put(50.00,40.00){\circle*{1.00}}
\put(50.00,50.00){\circle*{1.00}}
\put(30.00,40.00){\circle*{1.00}}
\put(30.00,50.00){\circle*{1.00}}
\put(20.00,40.00){\circle*{1.00}}
\put(20.00,50.00){\circle*{1.00}}
\put(40.00,20.00){\circle*{1.00}}
\put(40.00,10.00){\circle*{1.00}}
\put(50.00,20.00){\circle*{1.00}}
\put(50.00,10.00){\circle*{1.00}}
\put(60.00,20.00){\circle*{1.00}}
\put(60.00,10.00){\circle*{1.00}}
\put(30.00,20.00){\circle*{1.00}}
\put(30.00,10.00){\circle*{1.00}}
\put(20.00,20.00){\circle*{1.00}}
\put(20.00,10.00){\circle*{1.00}} \put(40.00,0.00){\circle*{1.00}}
\put(70.00,30.00){\circle*{2.00}}
\put(60.00,40.00){\circle*{1.00}}
\put(60.00,50.00){\circle*{2.00}}
\put(40.00,51.00){\makebox(0,0)[cb]{$\alpha_2$}}
\put(51.00,19.00){\makebox(0,0)[lt]{$\alpha_1$}}
\put(50.00,31.00){\makebox(0,0)[cb]{$\omega_1$}}
\put(71.00,31.00){\makebox(0,0)[lb]{$\lambda_1$}}
\put(62.00,50.00){\makebox(0,0)[lc]{$\lambda_2$}}
\put(75.00,43.00){\makebox(0,0)[cc]{$\DWC$}}
\put(51.00,41.00){\makebox(0,0)[rb]{$\omega_2$}}
\end{picture} &
%TeXCAD Picture [Sp4-fan.pic]. Options:
%\grade{\off}
%\emlines{\off}
%\epic{\off}
%\beziermacro{\off}
%\reduce{\on}
%\snapping{\off}
%\quality{2.00}
%\graddiff{0.01}
%\snapasp{1}
%\zoom{10.0000}
\unitlength .4ex % = .4pt
\linethickness{0.4pt}
\begin{picture}(80,65)(0,0)
\thicklines \put(40,30){\line(-1,-1){30}}
\put(0,30){\line(1,0){40}} \put(40,30){\vector(1,-1){10}}
\put(40,30){\vector(0,1){10}} \put(5,12){\makebox(0,0)[cc]{$\Vv$}}
\put(40,40){\circle*{2}} \put(50,20){\circle*{2}}
\put(42,41){\makebox(0,0)[lb]{$\alpha^{\vee}_2$}}
\put(53,19){\makebox(0,0)[lb]{$\alpha^{\vee}_1$}} \thinlines
\put(70,60){\line(-1,-1){30}} \put(40,30){\line(1,0){40}}
%\emline(70,30)(60,50)
\multiput(70,30)(-.1190476,.2380952){84}{\line(0,1){.2380952}}
%\end
%\emline(40,60)(60,50)
\multiput(40,60)(.2380952,-.1190476){84}{\line(1,0){.2380952}}
%\end
%\emline(10,30)(20,50)
\multiput(10,30)(.1190476,.2380952){84}{\line(0,1){.2380952}}
%\end
%\emline(40,60)(20,50)
\multiput(40,60)(-.2380952,-.1190476){84}{\line(-1,0){.2380952}}
%\end
%\emline(70,30)(60,10)
\multiput(70,30)(-.1190476,-.2380952){84}{\line(0,-1){.2380952}}
%\end
%\emline(40,0)(60,10)
\multiput(40,0)(.2380952,.1190476){84}{\line(1,0){.2380952}}
%\end
%\emline(10,30)(20,10)
\multiput(10,30)(.1190476,-.2380952){84}{\line(0,-1){.2380952}}
%\end
%\emline(40,0)(20,10)
\multiput(40,0)(-.2380952,.1190476){84}{\line(-1,0){.2380952}}
%\end
%\emline(40,30)(10,15)
\multiput(40,30)(-.24,-.12){125}{\line(-1,0){.24}}
%\end
\put(40,30){\line(0,1){35}}
%\emline(40,30)(65,5)
\multiput(40,30)(.11961722,-.11961722){209}{\line(0,-1){.11961722}}
%\end
%\qbezier(30,25)(26,39.6)(40,45)
\put(30,25){\line(0,1){.241}} \put(29.94,25.24){\line(0,1){.239}}
\put(29.87,25.48){\line(0,1){.238}}
\put(29.81,25.72){\line(0,1){.237}}
\put(29.76,25.96){\line(0,1){.236}}
\put(29.7,26.19){\line(0,1){.234}}
\put(29.65,26.43){\line(0,1){.233}}
\put(29.6,26.66){\line(0,1){.232}}
\put(29.55,26.89){\line(0,1){.231}}
\put(29.5,27.12){\line(0,1){.229}}
\put(29.46,27.35){\line(0,1){.228}}
\put(29.42,27.58){\line(0,1){.227}}
\put(29.38,27.81){\line(0,1){.226}}
\put(29.35,28.03){\line(0,1){.224}}
\put(29.32,28.26){\line(0,1){.223}}
\put(29.28,28.48){\line(0,1){.222}}
\put(29.26,28.7){\line(0,1){.221}}
\put(29.23,28.92){\line(0,1){.219}}
\put(29.21,29.14){\line(0,1){.218}}
\put(29.19,29.36){\line(0,1){.217}}
\put(29.17,29.58){\line(0,1){.216}}
\put(29.15,29.79){\line(0,1){.214}}
\put(29.14,30){\line(0,1){2.274}}
\put(29.16,32.28){\line(0,1){.199}}
\put(29.17,32.48){\line(0,1){.198}}
\put(29.19,32.68){\line(0,1){.197}}
\put(29.21,32.87){\line(0,1){.195}}
\put(29.24,33.07){\line(0,1){.194}}
\put(29.26,33.26){\line(0,1){.193}}
\put(29.29,33.46){\line(0,1){.192}}
\put(29.32,33.65){\line(0,1){.19}}
\put(29.36,33.84){\line(0,1){.189}}
\put(29.39,34.03){\line(0,1){.188}}
\put(29.43,34.21){\line(0,1){.187}}
\put(29.47,34.4){\line(0,1){.185}}
\put(29.51,34.59){\line(0,1){.184}}
\put(29.56,34.77){\line(0,1){.183}}
\put(29.61,34.95){\line(0,1){.182}}
\put(29.66,35.14){\line(0,1){.18}}
\put(29.71,35.32){\line(0,1){.179}}
\put(29.77,35.5){\line(0,1){.178}}
\put(29.83,35.67){\line(0,1){.177}}
\put(29.89,35.85){\line(0,1){.175}}
\put(29.95,36.02){\line(0,1){.174}}
\put(30.01,36.2){\line(0,1){.173}}
\put(30.08,36.37){\line(0,1){.172}}
\put(30.15,36.54){\line(0,1){.17}}
\put(30.23,36.71){\line(0,1){.169}}
\put(30.3,36.88){\line(0,1){.168}}
\put(30.38,37.05){\line(0,1){.167}}
\put(30.46,37.22){\line(0,1){.165}}
\put(30.54,37.38){\line(0,1){.164}}
\put(30.63,37.55){\line(0,1){.163}}
\put(30.71,37.71){\line(0,1){.162}}
\put(30.8,37.87){\line(0,1){.16}}
\put(30.9,38.03){\line(0,1){.159}}
\put(30.99,38.19){\line(0,1){.158}}
\put(31.09,38.35){\line(0,1){.156}}
\put(31.19,38.5){\line(0,1){.155}}
\put(31.29,38.66){\line(0,1){.154}}
\put(31.4,38.81){\line(0,1){.153}}
\put(31.5,38.97){\line(0,1){.151}}
\put(31.61,39.12){\line(0,1){.15}}
\put(31.73,39.27){\line(0,1){.149}}
\put(31.84,39.42){\line(0,1){.148}}
\put(31.96,39.56){\line(0,1){.146}}
\multiput(32.08,39.71)(.061,.0726){2}{\line(0,1){.0726}}
\multiput(32.2,39.86)(.0622,.072){2}{\line(0,1){.072}}
\multiput(32.32,40)(.0635,.0713){2}{\line(0,1){.0713}}
\multiput(32.45,40.14)(.0647,.0707){2}{\line(0,1){.0707}}
\multiput(32.58,40.28)(.0659,.0701){2}{\line(0,1){.0701}}
\multiput(32.71,40.42)(.0671,.0694){2}{\line(0,1){.0694}}
\multiput(32.85,40.56)(.0684,.0688){2}{\line(0,1){.0688}}
\multiput(32.98,40.7)(.0696,.0682){2}{\line(1,0){.0696}}
\multiput(33.12,40.84)(.0708,.0676){2}{\line(1,0){.0708}}
\multiput(33.26,40.97)(.0721,.0669){2}{\line(1,0){.0721}}
\multiput(33.41,41.11)(.0733,.0663){2}{\line(1,0){.0733}}
\multiput(33.55,41.24)(.0745,.0657){2}{\line(1,0){.0745}}
\multiput(33.7,41.37)(.0757,.0651){2}{\line(1,0){.0757}}
\multiput(33.85,41.5)(.077,.0644){2}{\line(1,0){.077}}
\multiput(34.01,41.63)(.0782,.0638){2}{\line(1,0){.0782}}
\multiput(34.16,41.76)(.0794,.0632){2}{\line(1,0){.0794}}
\multiput(34.32,41.88)(.0807,.0625){2}{\line(1,0){.0807}}
\multiput(34.48,42.01)(.0819,.0619){2}{\line(1,0){.0819}}
\multiput(34.65,42.13)(.0831,.0613){2}{\line(1,0){.0831}}
\multiput(34.81,42.25)(.0844,.0607){2}{\line(1,0){.0844}}
\multiput(34.98,42.38)(.0856,.06){2}{\line(1,0){.0856}}
\put(35.15,42.5){\line(1,0){.174}}
\put(35.33,42.61){\line(1,0){.176}}
\put(35.5,42.73){\line(1,0){.179}}
\put(35.68,42.85){\line(1,0){.181}}
\put(35.86,42.96){\line(1,0){.183}}
\put(36.05,43.08){\line(1,0){.186}}
\put(36.23,43.19){\line(1,0){.188}}
\put(36.42,43.3){\line(1,0){.191}}
\put(36.61,43.41){\line(1,0){.193}}
\put(36.81,43.52){\line(1,0){.196}}
\put(37,43.63){\line(1,0){.198}}
\put(37.2,43.73){\line(1,0){.201}}
\put(37.4,43.84){\line(1,0){.203}}
\put(37.6,43.94){\line(1,0){.206}}
\put(37.81,44.04){\line(1,0){.208}}
\put(38.02,44.15){\line(1,0){.211}}
\put(38.23,44.25){\line(1,0){.213}}
\put(38.44,44.34){\line(1,0){.215}}
\put(38.66,44.44){\line(1,0){.218}}
\put(38.87,44.54){\line(1,0){.22}}
\put(39.09,44.63){\line(1,0){.223}}
\multiput(39.32,44.73)(.2265,.0918){2}{\line(1,0){.2265}}
\put(39.77,44.91){\line(1,0){.23}}
%\end
%\qbezier(26,23)(35.5,11)(53,17)
\multiput(26,23)(.0703,-.0881){2}{\line(0,-1){.0881}}
\multiput(26.14,22.82)(.0708,-.0871){2}{\line(0,-1){.0871}}
\multiput(26.28,22.65)(.0712,-.0861){2}{\line(0,-1){.0861}}
\multiput(26.42,22.48)(.0716,-.0851){2}{\line(0,-1){.0851}}
\multiput(26.57,22.31)(.0721,-.0841){2}{\line(0,-1){.0841}}
\multiput(26.71,22.14)(.0725,-.0832){2}{\line(0,-1){.0832}}
\multiput(26.86,21.97)(.0729,-.0822){2}{\line(0,-1){.0822}}
\multiput(27,21.81)(.0734,-.0812){2}{\line(0,-1){.0812}}
\multiput(27.15,21.65)(.0738,-.0802){2}{\line(0,-1){.0802}}
\multiput(27.3,21.49)(.0743,-.0792){2}{\line(0,-1){.0792}}
\multiput(27.45,21.33)(.0747,-.0783){2}{\line(0,-1){.0783}}
\multiput(27.6,21.17)(.0751,-.0773){2}{\line(0,-1){.0773}}
\multiput(27.75,21.02)(.0756,-.0763){2}{\line(0,-1){.0763}}
\multiput(27.9,20.86)(.076,-.0753){2}{\line(1,0){.076}}
\multiput(28.05,20.71)(.0764,-.0743){2}{\line(1,0){.0764}}
\multiput(28.2,20.56)(.0769,-.0734){2}{\line(1,0){.0769}}
\multiput(28.36,20.42)(.0773,-.0724){2}{\line(1,0){.0773}}
\multiput(28.51,20.27)(.0777,-.0714){2}{\line(1,0){.0777}}
\multiput(28.67,20.13)(.0782,-.0704){2}{\line(1,0){.0782}}
\multiput(28.82,19.99)(.0786,-.0694){2}{\line(1,0){.0786}}
\multiput(28.98,19.85)(.079,-.0685){2}{\line(1,0){.079}}
\multiput(29.14,19.71)(.0795,-.0675){2}{\line(1,0){.0795}}
\multiput(29.3,19.58)(.0799,-.0665){2}{\line(1,0){.0799}}
\multiput(29.46,19.44)(.0804,-.0655){2}{\line(1,0){.0804}}
\multiput(29.62,19.31)(.0808,-.0645){2}{\line(1,0){.0808}}
\multiput(29.78,19.18)(.0812,-.0636){2}{\line(1,0){.0812}}
\multiput(29.94,19.06)(.0817,-.0626){2}{\line(1,0){.0817}}
\multiput(30.1,18.93)(.0821,-.0616){2}{\line(1,0){.0821}}
\multiput(30.27,18.81)(.0825,-.0606){2}{\line(1,0){.0825}}
\put(30.43,18.69){\line(1,0){.166}}
\put(30.6,18.57){\line(1,0){.167}}
\put(30.77,18.45){\line(1,0){.168}}
\put(30.93,18.34){\line(1,0){.169}}
\put(31.1,18.22){\line(1,0){.169}}
\put(31.27,18.11){\line(1,0){.17}}
\put(31.44,18){\line(1,0){.171}}
\put(31.61,17.89){\line(1,0){.172}}
\put(31.78,17.79){\line(1,0){.173}}
\put(31.96,17.69){\line(1,0){.174}}
\put(32.13,17.58){\line(1,0){.175}}
\put(32.31,17.48){\line(1,0){.176}}
\put(32.48,17.39){\line(1,0){.176}}
\put(32.66,17.29){\line(1,0){.177}}
\put(32.84,17.2){\line(1,0){.178}}
\put(33.01,17.1){\line(1,0){.179}}
\put(33.19,17.01){\line(1,0){.18}}
\put(33.37,16.93){\line(1,0){.181}}
\put(33.55,16.84){\line(1,0){.182}}
\put(33.73,16.76){\line(1,0){.182}}
\put(33.92,16.67){\line(1,0){.183}}
\put(34.1,16.59){\line(1,0){.184}}
\put(34.28,16.52){\line(1,0){.185}}
\put(34.47,16.44){\line(1,0){.186}}
\put(34.66,16.37){\line(1,0){.187}}
\put(34.84,16.29){\line(1,0){.188}}
\put(35.03,16.22){\line(1,0){.189}}
\put(35.22,16.16){\line(1,0){.189}}
\put(35.41,16.09){\line(1,0){.19}}
\put(35.6,16.02){\line(1,0){.191}}
\put(35.79,15.96){\line(1,0){.192}}
\put(35.98,15.9){\line(1,0){.193}}
\put(36.17,15.84){\line(1,0){.194}}
\put(36.37,15.79){\line(1,0){.195}}
\put(36.56,15.73){\line(1,0){.196}}
\put(36.76,15.68){\line(1,0){.196}}
\put(36.96,15.63){\line(1,0){.197}}
\put(37.15,15.58){\line(1,0){.198}}
\put(37.35,15.53){\line(1,0){.199}}
\put(37.55,15.49){\line(1,0){.2}}
\put(37.75,15.45){\line(1,0){.201}}
\put(37.95,15.41){\line(1,0){.202}}
\put(38.15,15.37){\line(1,0){.203}}
\put(38.35,15.33){\line(1,0){.203}}
\put(38.56,15.29){\line(1,0){.204}}
\put(38.76,15.26){\line(1,0){.205}}
\put(38.97,15.23){\line(1,0){.206}}
\put(39.17,15.2){\line(1,0){.415}}
\put(39.59,15.15){\line(1,0){.418}}
\put(40.01,15.1){\line(1,0){.422}}
\put(40.43,15.07){\line(1,0){.425}}
\put(40.85,15.04){\line(1,0){.429}}
\put(41.28,15.02){\line(1,0){.432}}
\put(41.71,15.01){\line(1,0){1.317}}
\put(43.03,15.01){\line(1,0){.446}}
\put(43.48,15.03){\line(1,0){.45}}
\put(43.93,15.06){\line(1,0){.453}}
\put(44.38,15.09){\line(1,0){.456}}
\put(44.84,15.13){\line(1,0){.46}}
\put(45.3,15.18){\line(1,0){.463}}
\put(45.76,15.24){\line(1,0){.467}}
\put(46.23,15.31){\line(1,0){.47}}
\put(46.7,15.38){\line(1,0){.474}}
\put(47.17,15.46){\line(1,0){.477}}
\put(47.65,15.55){\line(1,0){.481}}
\put(48.13,15.65){\line(1,0){.484}}
\put(48.61,15.75){\line(1,0){.488}}
\multiput(49.1,15.86)(.2457,.0601){2}{\line(1,0){.2457}}
\multiput(49.59,15.98)(.2474,.0641){2}{\line(1,0){.2474}}
\multiput(50.09,16.11)(.2492,.068){2}{\line(1,0){.2492}}
\multiput(50.59,16.25)(.2509,.0719){2}{\line(1,0){.2509}}
\multiput(51.09,16.39)(.2526,.0758){2}{\line(1,0){.2526}}
\multiput(51.59,16.54)(.2548,.0807){3}{\line(1,0){.2548}}
\multiput(52.36,16.78)(.3215,.1076){2}{\line(1,0){.3215}}
%\end
\put(25,39.5){\makebox(0,0)[cc]{$\Cc_{Y_1}$}}
\put(35.5,11){\makebox(0,0)[cc]{$\Cc_{Y_2}$}}
\end{picture} \\[1ex]
(a) Weight polytope & (b) Colored fan \\[2ex]
%TeXCAD Picture [SL2T1.pic]. Options:
%\grade{\off}
%\emlines{\off}
%\epic{\off}
%\beziermacro{\off}
%\reduce{\on}
%\snapping{\off}
%\quality{2.00}
%\graddiff{0.01}
%\snapasp{1}
%\zoom{10.0000}
\unitlength .4ex % = .4pt
\linethickness{0.4pt}
\begin{picture}(72,40)(0,0)
\put(6,25){\circle*{1}} \put(6,35){\circle*{1}}
\put(16,15){\circle*{1}} \put(16,25){\circle*{1}}
\put(36,5){\circle*{2}} \put(36,15){\circle*{1}}
\put(26,35){\circle*{1}} \put(16,35){\circle*{1}}
\put(26,25){\circle*{1}} \put(16,25){\circle*{1}}
\put(26,15){\circle*{1}} \put(16,15){\circle*{1}}
\put(6,35){\circle*{1}} \put(36,5){\circle*{2}}
\put(46,15){\circle*{2}} \put(56,15){\circle*{2}}
\put(36,5){\line(0,1){35}} \put(46,25){\circle*{2}}
\put(46,35){\circle*{2}} \put(56,25){\circle*{2}}
\put(56,35){\circle*{2}} \put(66,25){\circle*{2}}
\put(66,35){\circle*{2}} \put(36,25){\circle*{2}}
\put(36,35){\circle*{2}} \put(36,5){\vector(1,0){20}}
\put(56,3){\makebox(0,0)[ct]{$\alpha_2$}}
\put(36,5){\line(2,1){36}} \put(36,5){\line(-2,1){36}}
\end{picture} &
%TeXCAD Picture [GL2.pic]. Options:
%\grade{\off}
%\emlines{\off}
%\epic{\off}
%\beziermacro{\off}
%\reduce{\on}
%\snapping{\off}
%\quality{2.00}
%\graddiff{0.01}
%\snapasp{1}
%\zoom{10.0000}
\unitlength .4ex % = .4pt
\linethickness{0.4pt}
\begin{picture}(84,40)(0,0)
\put(12,35){\circle*{1}} \put(22,25){\circle*{1}}
\put(2,25){\circle*{1}} \put(42,5){\circle*{2}}
\put(32,35){\circle*{1}} \put(22,25){\circle*{1}}
\put(2,25){\circle*{1}} \put(32,15){\circle*{1}}
\put(52,15){\circle*{1}} \put(12,15){\circle*{1}}
\put(12,35){\circle*{1}} \put(42,5){\circle*{2}}
\put(72,15){\circle*{2}} \put(42,5){\line(0,1){35}}
\put(52,35){\circle*{2}} \put(62,25){\circle*{2}}
\put(82,25){\circle*{2}} \put(72,35){\circle*{2}}
\put(42,25){\circle*{2}} \put(42,5){\vector(1,0){20}}
\put(62,3){\makebox(0,0)[ct]{$\alpha_1$}}
\put(42,5){\line(3,1){42}} \put(42,5){\line(-3,1){42}}
\end{picture} \\[1ex]
\multicolumn{2}{c}{(c) Slice semigroups}
\end{tabular}
\end{center}
\end{figure}
There are two closed orbits $Y_1,Y_2\subset X$. The respective
Levi subgroups are $L_1=\SL_2(\kk)\times\kk^{\times}$ and
$L_2=\GL_2(\kk)$, with the simple roots $\alpha_2$ and $\alpha_1$,
respectively.

Consider the slice monoids $Z_i$ for~$Y_i$. The weight semigroups
of $F_i=\overline{T}$ (the closure in~$Z_i$) are plotted by dots
in Figure~\ref{Sp(4)}(c), the bold dots corresponding to the
weight semigroups $\Sigma_i$ of~$Z_i$. (They are easily computed
using the Clebsch--Gordan formula.) We can now see that $F_i$ are
normal, but $Z_i$ are not, i.e., $X$~is non-normal along
$Y_1,Y_2$. However, if we increase $V$ by adding two highest
weights $\lambda_2=2\omega_1$, $\lambda_3=\omega_1+\omega_2$, then
$X$ becomes normal. Its colored fan is depicted in
Figure~\ref{Sp(4)}(b).
\end{example}

The projective completions of adjoint simple groups in projective
linear operators on fundamental and adjoint representation spaces
were studied in detail in \cite[\S12]{grp.comp}. In particular,
the orbital decomposition was described, and normal and smooth
completions were identified.

\begin{example}
Suppose $G=\SO_{2l+1}(\kk)$, and $V=V(\omega_i)$ is a fundamental
representation. We have a unique closed orbit $Y\subset X$. If
$i<l$, then
$L\not\iso\GL_{n_1}(\kk)\times\dots\times\GL_{n_s}(\kk)$, hence
$Z$ and $X$ are singular. But for $i=l$ (the spinor
representation), $L\iso\GL_l(\kk)$ and
$V(\omega_l)\otimes\kk_{-\omega_l}$ is $L$-isomorphic to
$\E^{\bullet}\kk^l$. It follows that $Z\iso\LO_l(\kk)$, whence $X$
is smooth.
\end{example}

\begin{example}
Suppose that all vertices of $\Pp$ are regular weights. Then the
slice monoids $Z$ are toric and their weight semigroups $\Sigma$
are generated by the weights $\mu-\lambda$, where $\mu$ runs over
all $T$-weights of~$V$. The variety $X$ is toroidal, and normal
(smooth) iff each $\Sigma$ consists of all lattice vectors in the
corner cone of $\Pp$ at~$\lambda$ (resp.\ $\Sigma$ is generated by
linearly independent weights).

In particular, if $V=V(\lambda)$ is a simple module of regular
highest weight, then $\Sigma=\ZZ(-\Pi)$, whence $X$ is smooth.
This is a particular case of a wonderful completion,
see~\ref{wonderful}.

A interesting model for the wonderful completion of $G$ in terms
of Hilbert schemes was proposed by Brion~\cite{wonder<-Hilb}.
Namely, given a generalized flag variety $M=G/Q$, he proves that
the closure $X=\overline{(G\times G)[\diag{M}]}$ in the Hilbert
scheme (or the Chow variety) of $M\times M$ is isomorphic to the
wonderful completion.
\begin{question}
Generalize to arbitrary symmetric spaces via embedding into~$G$.
\end{question}
If $G=(\Aut M)^0$ (e.g., if $Q=B$), then $X$ is an irreducible
component of the Hilbert scheme (the Chow variety). All fibers of
the universal family over $X$ are reduced and Cohen--Macaulay
(even Gorenstein if $Q=B$).

Toroidal and wonderful group completions were studied intensively
in the framework of the general theory of toroidal and wonderful
varieties (see~\ref{toroidal}--\ref{wonderful}) and by their own.
De Concini and Procesi \cite{H(grp.comp)} and Strickland
\cite{H_G(grp.comp)} computed ordinary and equivariant rational
cohomology of smooth toroidal completions over $\kk=\CC$ (see also
\cite{H(reg.emb)}, \cite{H_G(wonder)}). Brion \cite{inf.Bruhat}
carried out a purely algebraic treatment of these results
replacing cohomology by (equivariant) Chow rings.

The basis of the Chow ring $A(X)$ of a smooth toroidal completion
$X=\overline{G}$ is given by the closures of the
Bia{\l}ynicki-Birula cells~\cite{BB}, which are isomorphic to
affine spaces and intersect ${G\times G}$-orbits in ${B^{-}\times
B}$-orbits \cite[2.3]{struct.loc}. The latter were described in
\cite[2.1]{inf.Bruhat}. The ${B^{-}\times B}$-orbit closures in
$X$ are smooth in codimension~$1$, but singular in codimension~$2$
(apart from trivial exceptions arising from $G=\PSL_2(\kk)$)
\cite[\S2]{inf.Bruhat}. For wonderful~$X$, the
Bia{\l}ynicki-Birula cells are described in
\cite[3.3]{inf.Bruhat}, and those intersecting $G$ (=the closures
in $X$ of ${B^{-}\times B}$-orbits in~$G$) are normal and
Cohen--Macaulay~\cite{large.Schubert}. Geometry of ${B^{-}\times
B}$-orbit closures in $X$ was studied in \cite{IH(B*B-orb)},
\cite{BWB(grp.comp)}.
\end{example}

The class of reductive group embeddings is not closed under
degenerations. Alexeev and Brion \cite{stab.red.aff},
\cite{stab.red.proj} introduced a more general class of
(\emph{stable}) \emph{reductive varieties} closed under flat
degenerations with irreducible (resp.\ reduced) fibers. Affine
(stable) reductive varieties may be defined as affine spherical
${G\times G}$-varieties $X$ such that $\RG(X)=\RG(G)\cap\Ss$ for
some subspace $\Ss\subseteq\RG(G)\otimes\QQ$ (resp.\ as seminormal
connected unions of reductive varieties); projective (stable)
reductive varieties are the projectivizations of affine ones.
Affine reductive varieties provide examples of algebraic
semigroups without unit.

Alexeev and Brion gave a combinatorial classification and
described the orbital decomposition for stable reductive varieties
in the spirit of Theorems~\ref{norm.mon},~\ref{orb(mon)}. They
constructed moduli spaces for affine stable reductive varieties
embedded in a ${G\times G}$-module and for stable reductive pairs,
i.e., projective stable reductive varieties with a distinguished
effective ample divisor containing no ${G\times G}$-orbit.

An interesting family of reductive varieties was introduced by
Vinberg~\cite{red.sgr}. Consider the group $\widehat{G}=(G\times
T)/Z$, where $Z=\{(t^{-1},t)\mid t\in Z(G)\}$. The cone
$\Cc\subset\ES(\widehat{G})$ spanned by (the projections to $\ES$
of) $(\alpha_i^{\vee},0)$ and $(-\gamma,\gamma)$,
$\gamma\in\Vv(G)$, defines a normal reductive monoid $\Env{G}$,
called the \emph{enveloping semigroup} of~$G$, with group of
invertibles~$\widehat{G}$. The projection
$\ES(\widehat{G})\to\ES(T/Z(G))$ maps $\Cc$ onto~$\DWC$. Hence by
Theorem~\ref{sph->sph} we have an equivariant map
$\pi_G:\Env{G}\to\AAA^l$, where ${G\times G}$ acts on $\AAA^l$
trivially and $T$ acts with the weights $\alpha_1,\dots,\alpha_l$.

The algebra $\kk[\Env{G}]=\bigoplus_{\chi\in\lambda+\ZZ_{+}\Pi}
\Mat(V(\lambda))\otimes\kk\chi$ is a free module over
$\kk[\AAA^l]=\kk[\ZZ_{+}\Pi]$ and $\kk[\Env{G}]^{U^{-}\times U}=
\kk[\Ch_{+}]\otimes\kk[\ZZ_{+}\Pi]$, i.e., all schematic fibers of
$\pi_G$ have the same algebra of ${U^{-}\times
U}$-invariants~$\kk[\Ch_{+}]$. Hence $\pi_G$ is flat and all its
fibers are reduced and irreducible by
Theorem~\ref{U-inv}\ref{A<=>A^U}, i.e., $\Env{G}$~is the total
space of a family of reductive varieties. (In fact,
$\AAA^l=(\Env{G})\by(G\times G)$ and $\pi_G$ is the categorical
quotient map.)

It is easy to see that the fibers of $\pi_G$ over points with
nonzero coordinates are isomorphic to~$G$. Degenerate fibers are
obtained from $G$ by a deformation of the multiplication law
in~$\kk[G]$. In particular, the ``most degenerate'' fiber
$\As{G}:=\pi_G^{-1}(0)$, called the \emph{asymptotic semigroup}
of~$G$, is just the horospherical contraction of~$G$
(see~\ref{horosph}). In a sense, the asymptotic semigroup reflects
the behavior of $G$ at infinity.

The enveloping semigroup is used in \cite[7.5]{stab.red.aff} to
construct families of affine reductive varieties with given
generic fiber~$X$: $\Env{X}=({\Env{G}\times X})\by{G}$, where $G$
acts as $\{\1\}\times\diag{G}\times\{\1\}\subset G\times G\times
G\times G$, so that
$\kk[\Env{X}]=\bigoplus_{\chi\in\lambda+\ZZ_{+}\Pi}
\kk[X]_{(\lambda)}\otimes\kk\chi\subseteq\kk[{X\times T}]$. The
map $\pi_G$ induces a flat morphism $\pi_X:\Env{X}\to\AAA^l$ with
reduced and irreducible fibers.

It was proved in \cite[7.6]{stab.red.aff} that $\pi_X$ is a
locally universal family of reductive varieties with generic
fiber~$X$, i.e., every flat family of affine reductive varieties
with reduced fibers over irreducible base is locally a pullback
of~$\pi_X$. The universal property for enveloping semigroups was
already noticed in~\cite{red.sgr}.

\begin{example}
Let us describe the enveloping semigroup of $G=\SL_n(\kk)$, using
the notation of Example~\ref{Mat(n)}. Here $\RG_{+}(\Env{G})$ is
generated by $(\omega_i,\omega_i)$, $(0,\alpha_i)$,
$i=1,\dots,n-1$. Recall that $\omega_i=\eps_1+\dots+\eps_i$ is the
highest weight of $\E^i\kk^n$. Thus $\Env\SL_n(\kk)$ is the
closure in $\LO(V)$ of the image of ${\SL_n(\kk)\times T}$ acting
on $V=\E^{\bullet}\kk^n\oplus\kk^{n-1}$, where $\SL_n(\kk)$ acts
on $\E^{\bullet}\kk^n$ in a natural way, and $T$ acts on
$\E^k\kk^n$ by the weight $\eps_1+\dots+\eps_k$ and on $\kk^{n-1}$
by the weights $\eps_i-\eps_{i+1}$. In other words, the image of
${\SL_n(\kk)\times T}$ consists of tuples of the form
%*
\begin{equation*}
(t_1g,\dots,t_1\cdots t_k\E^kg,\dots, t_1/t_2,\dots,t_{n-1}/t_n)
\end{equation*}
%*
where $g\in\SL_n(\kk)$, $t=\diag(t_1,\dots,t_n)\in T$ ($t_1\cdots
t_n=1$). It follows that
%*
\begin{multline*}
\Env\SL_n(\kk)=\biggl\{(a_1,\dots,a_k,\dots,z_1,\dots,z_{n-1})
\biggm|a_k\in\LO(\E^k\kk^n),\ z_i\in\kk,\\
a_k\wedge a_l=a_{k+l}
\prod\limits_{\substack{i=1,\dots,k\\j=1,\dots,l}}z_{i+j-1},\
a_n=1\biggr\}
\end{multline*}
%*
In particular, $\Env\SL_2(\kk)=\LO_2(\kk)$ and $\As\SL_2(\kk)$ is
the subsemigroup of degenerate matrices.
\end{example}

\section{S-varieties}
\label{S-var}

Horospherical varieties of complexity $0$ form another class of
spherical varieties whose structure and embedding theory is
understood better than in the general case.

\begin{definition}
An \emph{S-variety} is an equivariant embedding of a horospherical
homogeneous space $\HS=G/S$.
\end{definition}

This terminology is due to Popov and Vinberg \cite{S-var}, though
they considered only the affine case. General S-varieties were
studied by Pauer \cite{norm(G/U)}, \cite{smooth(G/U)} in the case,
where $S$ is a maximal unipotent subgroup of~$G$.

S-varieties are spherical. We shall examine them from the
viewpoint of the Luna--Vust theory. In order to apply it, we have
to describe the colored space $\ES=\ES(\HS)$.

It is convenient to assume $S\supseteq U^{-}$; then
$S=L_0\semitimes\Ru{P^{-}}$ for a certain parabolic $P\supseteq B$
with the Levi subgroup $L\supseteq L_0\supseteq L'$ and the
unipotent radical $\Ru{P}$ (Lemma~\ref{hor.sub}). We may assume
that $L\supseteq T$. Put $T_0=T\cap L_0$.

We have $\RG(\HS)=\Ch(A)$, where $A=P^{-}/S\iso L/L_0\iso T/T_0$.
By Theorem~\ref{cent(hor)}, $\Vv(\HS)=\ES$. The space
$\ES=\CoCh(A)\otimes\QQ$ may be identified with the
orthocomplement of $\CoCh(T_0)\otimes\QQ$ in $\CoCh(T)\otimes\QQ$.
It follows from the Bruhat decomposition that the $B$-divisors on
$\ES$ are of the form $D_{\alpha}=\overline{Br_{\alpha}\bp}$,
$\alpha\in\Pi\setminus\Pi_0$, where $\Pi_0\subseteq\Pi$ is the
simple root system of~$L$. An argument similar to that in
\ref{monoids} shows that $D_{\alpha}$ maps to
$\overline{\alpha^{\vee}}$, the image of $\alpha^{\vee}$ under the
projection $\CoCh(T)\to\CoCh(A)$.

Theorem~\ref{sph.emb}\ref{col.fans} says that normal S-varieties
are classified by colored fans in~$\ES$, each fan consisting of
finitely many colored cones $(\Cc_i,\Rr_i)$, so that the cones
$\Cc_i$ form a polyhedral fan in~$\ES$, $\Rr_i\subseteq\Dd^B$, and
each $\Cc_i\setminus\{0\}$ contains all $\overline{\alpha^{\vee}}$
such that $D_{\alpha}\in\Rr_i$. The colored cones in a fan
correspond to the $G$-orbits $Y_i$ in the respective
S-variety~$X$, and $X$ is covered by simple open $S$-subvarieties
$X_i=\{x\in X\mid\overline{Gx}\supseteq Y_i\}$.

The following result ``globalizing'' Theorem~\ref{loc.str.sph} is
a nice example of how the combinatorial embedding theory of
\ref{c=0} helps to clarify the geometric structure of S-varieties.
For any $G$-orbit $Y\subseteq X$ let
$P(Y)=P[\Dd^B\setminus\Dd^B_Y]$ be the normalizer of generic
$B$-orbits in~$Y$ and $S(Y)\subseteq P(Y)$ the normalizer of
generic $U$-orbits, so that $S(Y)^{-}$ is the stabilizer of $G:Y$
(see~\ref{horosph}). The Levi subgroup $L(Y)\subseteq P(Y)$
containing $T$ has the simple root system
$\Pi_0\cup\{\alpha\in\Pi\mid D_{\alpha}\in\Dd^B_Y\}$, and
$S(Y)=L(Y)_0\semitimes\Ru{P(Y)}$, where the Levi subgroup $L(Y)_0$
is intermediate between $L(Y)$ and $L(Y)'$ and is in fact the
common kernel of all characters in $\RG(Y)$ or in
$\Ch(A)\cap\Cc_Y^{\ann}$.
\begin{theorem}[{cf.~\cite[5.4]{norm(G/U)}}]\label{loc.str.S}
Let $X$ be a simple normal S-variety with the unique closed
$G$-orbit $Y\subseteq X$.
\begin{enumerate}
\item There exists a $P(Y)$-stable affine closed subvariety
$Z\subseteq X$ such that $\Ru{P(Y)^{-}}$ acts on $Z$ trivially and
$X\iso G\itimes{P(Y)^{-}}Z$. \item There exists an $S(Y)$-stable
closed subvariety $Z_0\subseteq Z$ with a fixed point such that
$Z\iso {P(Y)^{-}\itimes{S(Y)^{-}}Z_0}\iso{L(Y)\itimes{L(Y)_0}Z_0}$
and $X\iso{G\itimes{S(Y)^{-}}Z_0}$. \item\label{hor.slices} The
varieties $Z$ and $Z_0$ are equivariant affine embeddings of
${L(Y)/L(Y)\cap S}$ and ${L(Y)_0/L(Y)_0\cap S}$ whose weight
lattices are $\Ch(A)$ and ${\Ch(A)/\Ch(A)\cap\Cc_Y^{\ann}}$,
colored spaces are $\ES$ and $\ES_0:=\langle\Cc_Y\rangle$, and
colored cones coincide with $(\Cc_Y,\Dd^B_Y)$.
\end{enumerate}
\end{theorem}
\begin{proof}
The idea of the proof is to construct normal affine S-varieties
$Z$ and $Z_0$ with the colored data as in~\ref{hor.slices} and
then to verify that the colored data of $L(Y)\itimes{L(Y)_0}Z_0$
coincide with those of~$Z$ and the colored data of
$G\itimes{P(Y)^{-}}Z$ with those of~$X$. In each case both
varieties under consideration are simple normal embeddings of one
and the same homogeneous space. The restriction of
$B$-eigenfunctions to the fiber of each homogeneous bundle above
preserves the orders along $B$-stable divisors. It follows that
the colored cones of both varieties coincide with the colored cone
of the fiber, whence the varieties are isomorphic. Note that $Z_0$
contains a fixed point since it is determined by a colored cone of
full dimension.
\end{proof}

The theorem shows that the local geometry of (normal) S-varieties
is completely reduced to the affine case (even to affine
S-varieties with a fixed point). Affine S-varieties were studied
in \cite{S-var} in characteristic~$0$ and in \cite[\S17]{HS} in
arbitrary characteristic.

First note that $\HS$ is quasiaffine iff all
$\overline{\alpha^{\vee}}$ are nonzero (whenever
$\alpha\in\Pi\setminus\Pi_0$) and generate a strictly convex cone
in~$\ES$ (Corollary~\ref{sph.qaff}). This holds iff there exists a
dominant weight $\lambda$ such that
$\langle\lambda,\Pi^{\vee}\setminus\Pi_0^{\vee}\rangle>0$ and
$\lambda|_{T_0}=1$, i.e, iff $S$ is regularly embedded in the
stabilizer of a highest weight vector of weight~$\lambda$
(cf.~Theorem~\ref{subpar}).

\begin{theorem}\label{norm.aff.S-var}
Let $X$ be the normal affine S-variety determined by a colored
cone $(\Cc,\Dd^B)$. Then
%*
\begin{equation*}
\kk[X]\iso
\bigoplus_{\lambda\in\Ch(A)\cap\Cc^{\vee}}V^{*}(\lambda^{*})
\subseteq\kk[G/S]=
\bigoplus_{\lambda\in\Ch(A)\cap\DWC}V^{*}(\lambda^{*})
\end{equation*}
%*
If the semigroup $\Ch(A)\cap\Cc^{\vee}$ is generated by dominant
weights $\lambda_1,\dots,\lambda_m$, then
$X\iso\overline{Gv}\subseteq V(\lambda_1^{*})\oplus\dots\oplus
V(\lambda_m^{*})$, where $v=v_{-\lambda_1}+\dots+v_{-\lambda_m}$
is the sum of respective lowest weight vectors.
\end{theorem}
\begin{proof}
Observe that
$R=\bigoplus_{\lambda\in\Ch(A)\cap\Cc^{\vee}}V^{*}(\lambda^{*})$
is the largest subalgebra of $\kk[G/S]$ with the given algebra of
$U$-invariants $R^U=\kk[X]^U\iso\kk[\Ch(A)\cap\Cc^{\vee}]$. Hence
$R\supseteq\kk[X]\supseteq\langle G\cdot R^U\rangle$, and the
extension is integral by Lemma~\ref{A>GA^U}. Now $R=\kk[X]$ since
$\kk[X]$ is integrally closed.

It is easy to see that $\overline{Gv}$ is an affine embedding of
$\HS$ such that $\kk[\overline{Gv}]$ is generated by
$V^{*}(\lambda_1^{*})\oplus\dots\oplus
V^{*}(\lambda_m^{*})\subset\kk[\HS]$. By Lemma~\ref{V(l)*V(m)},
$\kk[\overline{Gv}]=R$.
\end{proof}

Every (even non-normal) affine S-variety with the open orbit $\HS$
is realized in a $G$-module $V$ as $X=\overline{Gv}$, $v\in V^S$.
We may assume $V=\langle Gv\rangle$ and decompose
$v=v_{-\lambda_1}+\dots+v_{-\lambda_m}$, where $v_{-\lambda_i}$
are lowest vectors of certain antidominant weights~$-\lambda_i$.

In characteristic zero, $V\iso V(\lambda_1^{*})\oplus\dots\oplus
V(\lambda_m^{*})$ and the same arguments as in the proof of
Theorem~\ref{norm.aff.S-var} show that
$\kk[X]=\bigoplus_{\lambda\in\Sigma}V^{*}(\lambda^{*})$, where
$\Sigma$ is the semigroup generated by
$\lambda_1,\dots,\lambda_m$, and the dual Weyl modules
$V^{*}(\lambda^{*})\iso\kk[G]^{(B)}_{\lambda^{*}}\iso V(\lambda)$
are the (simple) $G$-isotypic components of~$\kk[X]$. It is easy
to see that $Gv\iso\HS$ iff $\lambda_1,\dots,\lambda_m$
span~$\Ch(A)$. Thus we obtain the following
\begin{proposition}[{\cite[3.1, 3.4]{S-var}}]
In the case $\ch\kk=0$, affine S-varieties $X$ with the open orbit
$\HS$ bijectively correspond to finitely generated semigroups
$\Sigma$ of dominant weights spanning~$\Ch(A)$, via
$\Sigma=\RG_{+}(X)$. The variety $X$ is normal iff the semigroup
$\Sigma$ is saturated, i.e., $\Sigma=\QQ_{+}\Sigma\cap\Ch(A)$.
Moreover, the saturation
$\widetilde\Sigma=\QQ_{+}\Sigma\cap\Ch(A)$ of $\Sigma$ corresponds
to the normalization $\widetilde{X}$ of~$X$.
\end{proposition}

$G$-orbits in an affine S-variety $X=\overline{Gv}\subseteq V$
have a transparent description ``dual'' to that in
Theorem~\ref{sph.emb}\ref{orb(sph)}.
\begin{proposition}[{\cite[Th.8]{S-var}}]
The orbits in $X$ are in bijection with the faces of
$\Cc^{\vee}=\QQ_{+}\lambda_1+\dots+\QQ_{+}\lambda_m$. The orbit
corresponding to a face $\Ff$ is represented by
$v_{\Ff}=\sum_{\lambda_i\in\Ff}v_{\lambda_i}$. The adherence of
orbits agrees with the inclusion of faces.
\end{proposition}
\begin{proof}
We have $X=G\overline{Tv}$ since $\overline{Tv}$ is $B$-stable.
The $T$-orbits in $\overline{Tv}$ are represented by $v_{\Ff}$
over all faces $\Ff\subseteq\Cc^{\vee}$, and the adherence of
orbits agrees with the inclusion of faces. On the other hand, it
is easy to see that the $U^{-}$-fixed point set in each $G$-orbit
of $X$ is a $T$-orbit, hence distinct $v_{\Ff}$ represent distinct
$G$-orbits.
\end{proof}

In characteristic zero, one can describe the defining equations of
$X$ in~$V$. Let $c=\sum\xi_i\xi_i^{*}\in\Uni\g$ be the Casimir
element w.r.t.\ to a $G$-invariant inner product on~$\g$,
$\xi_i,\xi_i^{*}$ being mutually dual bases. It is well known that
$c$ acts on $V(\lambda^{*})$ by a scalar
$c(\lambda)=(\lambda+2\rho,\lambda)$. Note that $c(\lambda)$
depends on $\lambda$ monotonously w.r.t.\ the partial order
induced by positive roots: if $\lambda=\mu+\sum k_i\alpha_i$,
$k_i\ge0$, then $c(\lambda)=c(\mu)+\sum
k_i\bigl((\lambda+2\rho,\alpha_i)+(\alpha_i,\mu)\bigr)\ge c(\mu)$,
and the inequality is strict, except for $\lambda=\mu$. The
following result is due to Kostant:
\begin{proposition}[\cite{LT}]\label{eqs(S-var)}
If $\ch\kk=0$ and $\lambda_1,\dots,\lambda_m$ are linearly
independent, then $\Ideal{X}\normin\kk[V]$ is generated by the
relations
%*
\begin{equation*}
c(x_i\otimes x_j)=
(\lambda_i+\lambda_j+2\rho,\lambda_i+\lambda_j)(x_i\otimes x_j),
\qquad i,j=1,\dots,m,
\end{equation*}
%*
where $x_k$ denotes the projection of $x\in V$
to~$V(\lambda_k^{*})$.
\end{proposition}
\begin{proof}
The algebra
$\kk[V]=\bigoplus_{k_1,\dots,k_m}\Sym^{k_1}V(\lambda_1)\otimes
\dots\otimes\Sym^{k_m}V(\lambda_m)$ is multigraded and $\Ideal{X}$
is a multihomogeneous ideal. The structure of $\kk[X]$ implies
that each homogeneous component $\Ideal{X}_{k_1,\dots,k_m}$ is the
kernel of the natural map
$\Sym^{k_1}V(\lambda_1)\otimes\dots\otimes
\Sym^{k_m}V(\lambda_m)\to V(k_1\lambda_1+\dots+k_m\lambda_m)$.

Consider a series of linear endomorphisms $\pi=c-c(\sum
k_i\lambda_i)\id$ of the subspaces
$\Sym^{k_1,\dots,k_m}V=\Sym^{k_1}V(\lambda_1^{*})
\otimes\dots\otimes\Sym^{k_m}V(\lambda_m^{*})\subset\Sym^{\bullet}V$.
Note that $\Ker\pi\iso V(\sum k_i\lambda_i^{*})$ is the highest
irreducible component of $\Sym^{k_1,\dots,k_m}V$, annihilated
by~$\Ideal{X}_{k_1,\dots,k_m}$, and
$\Im\pi\iso\Ideal{X}_{k_1,\dots,k_m}^{*}$ is the complementary
$G$-module.

It follows that $\Ideal{X}$ is spanned by the coordinate functions
of all $\pi(x_1^{k_1}\cdots x_m^{k_m})$. An easy calculation shows
that
%*
\begin{multline*}
\pi(x_1^{k_1}\cdots x_m^{k_m})=
\sum_i\frac{k_i(k_i-1)}2\pi(x_i^2)x_1^{k_1}\cdots
x_i^{k_i-2}\cdots x_m^{k_m}\\
+\sum_{i<j}k_ik_j\pi(x_ix_j)x_1^{k_1}\cdots x_i^{k_i-1}\cdots
x_j^{k_j-1}\cdots x_m^{k_m}
\end{multline*}
%*
Thus $\Ideal{X}$ is generated by the relations $\pi(x_ix_j)=0$,
$i,j=1,\dots,m$.
\end{proof}

If the generators of $\RG_{+}(X)$ are not linearly independent,
one has to extend the defining equations of $X$ by those arising
from the linear dependencies between the~$\lambda_i$'s,
see~\cite{Smi}.

The results of \ref{div} allow to compute the divisor class group
of a normal affine S-variety~$X$. Every Weil divisor is rationally
equivalent to a $B$-stable one $\delta=\sum
m_{\alpha}D_{\alpha}+\sum m_iY_i$, where $Y_i$ are the $G$-stable
prime divisors corresponding to the generators $v_i$ of the rays
of~$\Cc$ containing no colors. The divisor $\delta$ is principal
iff $m_{\alpha}=\langle\lambda,\alpha^{\vee}\rangle$ and
$m_i=\langle\lambda,v_i\rangle$ for a certain $\lambda\in\Ch(A)$.
This yields a finite presentation for~$\Cl{X}$. In particular, we
have
\begin{proposition}
An affine S-variety $X$ is factorial iff $\RG_{+}(X)$ is generated
by weights
$\lambda_1\dots,\lambda_s,\pm\lambda_{s+1},\dots,\pm\lambda_r$
($s\leq r$), where the $\lambda_i$'s are linearly independent and
the projection $\Ch(T)\to\Ch(T\cap G')$ maps them to distinct
fundamental weights or to~$0$.
\end{proposition}
For semisimple~$G$, we conclude that factorial S-varieties are
those corresponding to weight semigroups $\Sigma$ generated by
some of the fundamental weights \cite[Th.11]{S-var}.

The simplest class of affine S-varieties is formed by
\emph{HV-varieties}, i.e., cones of highest (or lowest) vectors
$X=\overline{Gv_{-\lambda}}$, $v_{-\lambda}\in V^{(B^{-})}$,
see~\ref{double.cone}. Particular examples are quadratic cones or
Grassmann cones of decomposable polyvectors. The above results on
affine S-varieties imply Proposition~\ref{HV}, which describes
basic properties of HV-varieties. It follows from
Proposition~\ref{eqs(S-var)} that an HV-cone is defined by
quadratic equations in the ambient simple $G$-module. For a
Grassmann cone we recover the Pl\"ucker relations between the
coordinates of a polyvector.

Now we describe smooth S-varieties in characteristic zero. By
Theorem~\ref{loc.str.S}, the problem is reduced to affine
S-varieties with a fixed point, which are nothing else but
$G$-modules with a dense orbit of a $U$-fixed vector.
\begin{lemma}
If a $G$-module $V$ is an S-variety, then $V=V_0\oplus
V_1\oplus\dots\oplus V_s$ so that $Z=Z(G)^0$ acts on $V_0$ with
linearly independent weights and each $V_i$ ($i>0$) is a simple
submodule acted on non-trivially by a unique simple factor
$G_i\subseteq G$, $G_i\iso\SL(V_i)$ or~$\Sp(V_i)$.
\end{lemma}
\begin{proof}
Since $Z$ has a dense orbit in $V_0=V^{G'}$, it acts with linearly
independent weights. If $G_i$ acts non-trivially on two simple
submodules $V_i,V_j$, and $v_i\in V_i^{(B)}$, $v_j\in
V_j^{(B^{-})}$, then the stabilizer of $v_i+v_j$ is not
horospherical, i.e., $V$~is not an S-variety. Therefore we may
assume that $V$ is irreducible and each simple factor of $G$ acts
non-trivially.

Then $G$ acts transitively on $\PP(V)$, which implies
$G'\iso\SL(V)$ or~$\Sp(V)$ \cite{decomp}. Indeed, we have $V=\br
v_{-\lambda}$, where $v_{-\lambda}\in V$ is a lowest vector. Hence
there exists a unique root $\delta$ such that
$e_{\delta}v_{-\lambda}=v_{\lambda^{*}}$ is a highest vector. One
easily deduces that the root system of $G$ is indecomposable and
$\delta$ is the highest root, so that $\delta=\lambda+\lambda^{*}$
is the sum of two dominant weights, whence the assertion.
\end{proof}
The colored data of such a $G$-module $V$ are easy to write down.
Namely $\Pi\setminus\Pi_0=\{\alpha_1,\dots,\alpha_s\}$, where
$\alpha_i$ are the first simple roots in some components of $\Pi$
having the type $\AAa_l$ or~$\CCc_l$. The weight lattice $\Ch(A)$
is spanned by linearly independent weights
$\lambda_1,\dots,\lambda_r$, where $\lambda_1,\dots,\lambda_s$ are
the highest weights of~$V_i^{*}$, which project to the fundamental
weights $\omega_i$ corresponding to~$\alpha_i$, and
$\lambda_{s+1},\dots,\lambda_r$ are the weights of~$V_0^{*}$,
which are orthogonal to~$\Pi$. The cone $\Cc$ is spanned by the
basis
$\overline{\alpha_1^{\vee}},\dots,\overline{\alpha_s^{\vee}},
v_{s+1},\dots,v_r$ of $\CoCh(A)$ dual to
$\lambda_1,\dots,\lambda_r$. Using Theorem~\ref{loc.str.S} we
derive the description of colored data of arbitrary smooth
S-varieties:
\begin{theorem}[{cf.~\cite[3.5]{smooth(G/U)}}]
An $S$-variety $X$ is smooth iff all colored cones
$(\Cc_Y,\Dd^B_Y)$ in the colored fan of $X$ satisfy the following
properties:
\begin{enumerate}
\item\label{loc.fact} $\Cc_Y$~is generated by a part of a basis
of~$\CoCh(A)$, and all $\overline{\alpha^{\vee}}$ such that
$D_{\alpha}\in\Dd^B_Y$ are among the generators. \item The simple
roots $\alpha$ such that $D_{\alpha}\in\Dd^B_Y$ are isolated from
each other at the Dynkin diagram of~$G$, and each $\alpha$ is
connected with at most one component $\Pi_{\alpha}$ of~$\Pi_0$;
moreover, $\{\alpha\}\cup\Pi_{\alpha}$ has the type $\AAa_l$
or~$\CCc_l$, $\alpha$~being the first simple root therein.
\end{enumerate}
\end{theorem}
The condition \ref{loc.fact} is equivalent to the local
factoriality of~$X$.

\section{Toroidal embeddings}
\label{toroidal}

In this section we assume $\ch\kk=0$. Recall that a
$G$-equivariant normal embedding $X$ of a spherical homogeneous
space $\HS=G/H$ is said to be \emph{toroidal} if
$\Dd^B_Y=\emptyset$ for each $G$-orbit $Y\subseteq X$. Toroidal
embeddings are defined by fans in~$\Vv$ and $G$-morphisms between
them correspond to subdivisions of these fans in the same way as
in toric geometry \cite{toric.intro}. There is a more direct
relation between toroidal and toric varieties. Put $P=P(\HS)$,
with the Levi decomposition $P=L\Ru{P}$ and other notation
from~\ref{horosph}.
\begin{theorem}[{\cite[3.4]{val.sph}, \cite[2.4]{spher.full}}]
\label{loc.str.tor} A toroidal embedding $X\embof\HS$ is covered
by $G$-translates of an open $P$-stable subset
%*
\begin{equation*}
\X=X\setminus\bigcup_{D\in\Dd^B}D\iso P\itimes{L}Z\iso
\Ru{P}\times Z
\end{equation*}
%*
where $Z$ is a locally closed $L$-stable subvariety pointwise
fixed by~$L_0$. The variety $Z$ is a toric embedding of $A=L/L_0$
defined by the same fan as~$X$, and the $G$-orbits in $X$
intersect $Z$ in $A$-orbits.
\end{theorem}
\begin{proof}
The problem is easily reduced to the case, where $X$ contains a
unique closed orbit $Y$ with $\Cc_Y=\Vv$, $\Dd^B_Y=\emptyset$.
Such toroidal embeddings, called wonderful, are discussed
in~\ref{wonderful}. Indeed, consider another spherical homogeneous
space $\overline{\HS}=G/N(H)$. Then
$\overline{\Vv}=\Vv(\overline{\HS})=\Vv/(\Vv\cap-\Vv)$ is strictly
convex, whence there exists a wonderful embedding
$\overline{X}\embof\overline{\HS}$. The canonical map
$\phi:\HS\to\overline{\HS}$ extends to $X\to\overline{X}$ by
Theorem~\ref{sph->sph}. We have $P(\overline{\HS})=P$, and $\X,Z$
are the preimages of the respective subvarieties defined
for~$\overline{X}$. The assertion on fans and orbits is easy,
cf.~Remark~\ref{c.d.(slice)}.

For wonderful~$X$ one applies the local structure theorem in a
neighborhood of~$Y$: by Theorem~\ref{loc.str.sph}, $\X=\X_Y\iso
P\itimes{L}Z$, where $Z$ is toric since $Z\cap\HS$ is a single
$A$-orbit.
\end{proof}

It follows that toroidal varieties are locally toric and have at
worst Abelian quotient singularities. They inherit many nice
geometric properties from toric varieties. On the other hand, each
spherical variety is the image of a toroidal one by a proper
birational equivariant map: to obtain this toroidal covering
variety, just remove all colors from the fan. This universality of
toroidal varieties can be used to derive some properties of
spherical varieties from the toroidal case.

A toroidal variety is smooth iff all cones of its fan are
simplicial and generated by a part of a basis of~$\RG(\HS)^{*}$:
for toric varieties this is deduced from the description of the
coordinate algebra \cite[2.1]{toric.intro}
(cf.~Example~\ref{toric}) and the general case follows by
Theorem~\ref{loc.str.tor}. For a singular toroidal variety one may
construct an equivariant desingularization by subdividing its fan,
cf.~\cite[2.6]{toric.intro}.

Every (smooth) toroidal variety admits an equivariant (smooth)
completion, which is defined by adding new cones to the fan in
order to cover all of~$\Vv(\HS)$. Smooth complete toroidal
varieties have other interesting characterizations.
\begin{theorem}[\cite{reg.var}]\label{reg.var}
For a smooth $G$-variety $X$ consider the following conditions:
\begin{enumerate}
\item\label{tor} $X$~is toroidal. \item\label{reg.BCP} There is a
dense open orbit $\HS\subseteq X$ such that
$\partial{X}=X\setminus\HS$ is a divisor with normal crossings,
each orbit $Gx\subset X$ is locally the intersection of several
components of~$\partial{X}$, and $G_x$ has a dense orbit in
$T_xX/\g{x}$. \item\label{reg.Gin} There is a $G$-stable divisor
$D\subset X$ with normal crossings such that $\Gg_X=\Tt_X(-\log
D)$. \item\label{pseudo-free} $X$~is spherical and pseudo-free.
\end{enumerate}
Then \ref{pseudo-free}$\implies$\ref{tor}$\implies$%
\ref{reg.BCP}$\iff$\ref{reg.Gin}. If $X$ is complete or spherical,
then all conditions are equivalent.
\end{theorem}
$G$-varieties satisfying the condition~\ref{reg.BCP},
resp.~\ref{reg.Gin}, are known as \emph{regular} in the sense of
Bifet--de~Concini--Procesi \cite{H(reg.emb)}, resp.~of
Ginzburg~\cite{Gin}.
\begin{proof}
\begin{roster}
\item[\ref{tor}$\implies$\ref{reg.BCP}\&\ref{reg.Gin}]
Theorem~\ref{loc.str.tor} reduces the problem to smooth toric
varieties. The latter are covered by invariant affine open charts
of the form $X=\AAA^s\times({\AAA^1\setminus0})^{r-s}$, where
$(\kk^{\times})^r$ acts in the natural way, so that
$D=\partial{X}$ is the union of coordinate hyperplanes
$\{x_i=0\}$, $X$~is isomorphic to the normal bundle of the closed
orbit, and $\Tt_X(-\log D)$ is a free sheaf spanned by velocity
fields $x_1\partial_1,\dots,x_n\partial_n$
($\partial_i:=\partial/\partial x_i$).

\item[\ref{reg.BCP}$\iff$\ref{reg.Gin}] First observe that
$\HS=X\setminus D$ is a single $G$-orbit iff $\Gg_{X\setminus
D}=\Tt_{X\setminus D}$. Now consider a neighborhood of any $x\in
D$. Due to local nature of the
conditions~\ref{reg.BCP},~\ref{reg.Gin}, we may assume that all
components $D_1,\dots,D_k$ of $D$ contain~$x$. Choose local
parameters $x_1,\dots,x_n$ at $x$ such that $D_i$ are locally
defined by the equations $x_i=0$. Let
$\partial_1,\dots,\partial_n$ denote the vector fields dual to
$dx_1,\dots,dx_n$. Then $\Tt_X(-\log D)$ is locally generated by
$x_1\partial_1,\dots,x_k\partial_k,\partial_{k+1},\dots,\partial_n$.

Let $Y=D_1\cap\dots\cap D_k$ and
$\pi:N=\Spec\Sym^{\bullet}(\Ii_Y/\Ii_Y^2)\to Y$ be the normal
bundle. There is a natural embedding $\pi^{*}\Tt_X(-\log D)|_Y
\embeds\Tt_N$: each vector field in $\Tt_X(-\log D)$ preserves
$\Ii_Y$ whence induces a derivation
of~$\Sym^{\bullet}(\Ii_Y/\Ii_Y^2)$. The image of
$\pi^{*}\Tt_X(-\log D)|_X$  is $\Tt_N(-\log\bigcup N_i)$, where
$N_i$ are the normal bundles to $Y$ in~$D_i$. Indeed,
$\bar{x}_i=x_i\bmod\Ii_Y^2$ ($i\leq k$), $\bar{x}_j=\pi^{*}x_j|_Y$
($j>k$) are local parameters on~$N$ and $x_i\partial_i,\partial_j$
induce the derivations $\bar{x}_i\bar\partial_i,\bar\partial_j$.
Note that $N=\bigoplus L_i$, where $L_i=\bigcap_{j\neq i}N_j$ are
$G$-stable line subbundles. Hence the $G_x$-action on
$T_xX/T_xY=N(x)=\bigoplus L_i(x)$ is diagonalizable.

Condition \ref{reg.BCP} implies that $Gx$ is open in~$Y$ and the
weights of $G_x:L_i(x)$ are linearly independent. This yields
velocity fields $\bar{x}_i\bar\partial_i$ on~$N(x)$ and in
transversal directions, which locally generate $\Tt_N(-\log\bigcup
N_i)$. Therefore $\Tt_X(-\log D)|_Y$ is generated by velocity
fields. By Nakayama's lemma, $\Tt_X(-\log D)=\Gg_X$ in a
neighborhood of~$x$.

Conversely, \ref{reg.Gin}~implies $\Gg_Y=\Tt_Y$ and
$\Gg_N=\Tt_N(-\log\bigcup N_i)$. Hence $Gx$ is open in~$Y$ and
$N|_{Gx}=G\itimes{G_x}T_xX/\g{x}$ has an open $G$-orbit. Thus
$T_xX/\g{x}$ contains an open $G_x$-orbit.

\item[\ref{reg.BCP}\&\ref{reg.Gin}$\implies$\ref{pseudo-free}]
Since $\Tt_X(-\log\partial{X})$ is a vector bundle, the
implication is trivial provided that $X$ is spherical. It remains
to prove that $X$ is spherical if it is complete.

A closed orbit $Y\subseteq X$ intersects a $B$-chart $\X\iso
P\itimes{L}Z$, where $L\subseteq P=P(Y)$ is the Levi subgroup and
$Z$ is an $L$-stable affine subvariety intersecting $Y$ in a
single point~$z$. Since the maximal torus $T\subseteq L\subseteq
G_z=P^{-}$ acts on $T_zZ\iso T_zX/\g{z}$ with linearly independent
weights, $Z\iso T_zZ$ contains an open $T$-orbit, whence $\X$ has
an open $B$-orbit.

\item[\ref{pseudo-free}$\implies$\ref{tor}] There is a morphism
$X\to\Gr(\g)$, $x\mapsto[\h_x]$, extending the map
$x\mapsto[\g_x]$ on~$\HS$. If $X$ is not toroidal, then there
exists a $G$-orbit $Y\subset X$ contained in a $B$-divisor
$D\subset X$. Then we have $\br+\h_{gy}=\br+(\Ad{g})\h_y\neq\g$,
$\forall y\in Y,\ g\in G$, i.e., $\h_y$~is not spherical.

To obtain a contradiction, it suffices to prove that all $\h_x$
are spherical subalgebras. Passing to a toroidal variety mapping
onto~$X$, one may assume that $X$~itself is toroidal. Consider the
normal bundle $N$ to $Y=Gx$. Since $\Gg_N=\pi^{*}\Gg_X|_Y$,
$\h_x$~is the stabilizer subalgebra of general position for $G:N$.
But $N$ is spherical, because the minimal $B$-chart $\X$ of $Y$ is
$P$-isomorphic to~$N|_{Y\cap\X}$. \qedhere\end{roster}
\end{proof}

Toric varieties and generalized flag varieties form two
``extreme'' classes of toroidal varieties. A number of geometric
and cohomological results generalize from these particular cases
to general toroidal varieties. A powerful vanishing theorem was
proved by Bien and Brion (1991) and refined by Knop (1992).
\begin{theorem}[\cite{reg.var}]\label{H(log.tan)}
If $X$ is a smooth complete toroidal variety, then
$\Ho^i(X,\Sym^{\bullet}\Tt_X(-\log\partial{X}))=0$, $\forall i>0$.
\end{theorem}
For flag varieties, this result is due to Elkik (vanishing of
higher cohomology of the tangent sheaf was proved already by Bott
in 1957). In fact, Bien and Brion proved a twisted version of
Theorem~\ref{H(log.tan)} \cite[3.2]{reg.var}:
$\Ho^i(X,\Ll\otimes\Sym^{\bullet}\Tt_X(-\log\partial{X}))=0$ for
all $i>0$ and any globally generated line bundle $\Ll$ on~$X$,
under a technical condition that the stabilizer $H$ of $\HS$ is
parabolic in a reductive subgroup of~$G$. (Generally, higher
cohomology of globally generated line bundles vanishes on every
complete spherical variety, see Corollary~\ref{H(nef)}.)

In view of Theorem~\ref{reg.var}, Theorem~\ref{H(log.tan)} stems
from a more general vanishing result of Knop:
\begin{theorem}[{\cite[4.1]{HC-hom}}]\label{H(act.sheaf)}
If $X$ is a pseudo-free equivariant completion of a homogeneous
space~$\HS$, then $\Ho^i(X,\Uu^m_X)=\Ho^i(X,\Sym^m\Gg_X)=0$,
$\forall i>0,\ m\geq0$.
\end{theorem}
\begin{proof}[Synopsis of a proof]
The assertions on $\Uu_X$ are reduced to those on
$\Sym^{\bullet}\Gg_X=\gr\Uu_X$. Since $\pi_X:T^{\g}X\to X$ is an
affine morphism, the Leray spectral sequence reduces the question
to proving $\Ho^i\left(T^{\g}X,\Oo_{T^{\g}X}\right)=0$. The
localized moment map $\overline\Phi:T^{\g}X\to M_X$ factors
through $\widetilde\Phi:T^{\g}X\to\widetilde{M}_X$. As
$\widetilde{M}_X$ is affine,
$\Ho^i\left(T^{\g}X,\Oo_{T^{\g}X}\right)=
\Ho^0\bigl(\widetilde{M}_X,\Rf^i\widetilde\Phi_{*}\Oo_{T^{\g}X}\bigr)$,
and it remains to prove $\Rf^i\widetilde\Phi_{*}\Oo_{T^{\g}X}=0$.
Here one applies to $\widetilde\Phi$ a version of Koll\'ar's
vanishing theorem \cite[4.2]{HC-hom}:
\begin{quote}
If $Y$ is smooth, $Z$~has rational singularities, and $\phi:Y\to
Z$ is a proper morphism with connected generic fibers $F$, which
satisfy $\Ho^i(F,\Oo_F)=0$, $\forall i>0$, then
$\Rf^i\phi_{*}\Oo_Y=0$ for all $i>0$.
\end{quote}
It remains to verify the conditions. The morphism $\widetilde\Phi$
is proper by Example~\ref{comp.moment}. The variety
$\widetilde{M}_X$ has rational singularities by
\cite[4.3]{HC-hom}. To show vanishing of the higher cohomology
of~$\Oo_F$, it suffices to prove that $F$ is
unirational~\cite{fund(unirat)}. Here one may assume $X=\HS$,
$\widetilde\Phi:T^{*}\HS\to\widetilde{M}_{\HS}$. Unirationality of
the fibers of the moment map is the heart of the proof
\cite[\S5]{HC-hom}.
\end{proof}
There is a relative version of Theorem~\ref{H(act.sheaf)}
asserting $\Rf^i\psi_{*}\Uu^m_X=\Rf^i\psi_{*}\Sym^m\Gg_X=0$,
$\forall i>0$, for a proper morphism $\psi:X\to Y$ separating
generic orbits, where $X$ is pseudo-free and $Y$ has rational
singularities.

The vanishing theorems of Bien--Brion and Knop have a number of
important consequences. For instance, on a pseudo-free completion
$X$ of $\HS$ the symbol map $\gr\Ho^0(X,\Uu_X)\to\kk[T^{\g}X]$ is
surjective. In the toroidal case,
$\Ho^1(X,\Tt_X(-\log\partial{X}))=0$ implies that the pair
$(X,\partial{X})$ is locally rigid, by deformation theory of
Kodaira--Spencer~\cite{deform}. Using this observation, Alexeev
and Brion proved Luna's conjecture on rigidity of spherical
subgroups.
\begin{theorem}[{\cite[\S3]{rigid.sph}}]\label{sph.s.g.p}
For any (irreducible) $G$-variety with spherical (generic) orbits,
the stabilizers of points in general position are conjugate.
\end{theorem}
\begin{proof}
Let $\Xx$ be a $G$-variety with spherical orbits. Passing to an
open subset, we may assume that $\Xx$ is smooth quasiprojective
and there exists a smooth $G$-invariant morphism $\pi:\Xx\to Z$
whose fibers contain dense orbits. Regarding $\Xx$ as a family of
spherical $G$-orbit closures, we may replace it by a birationally
isomorphic family of smooth projective toroidal varieties.

Indeed, there is a locally closed $G$-embedding of $\Xx$
into~$\PP(V)$, and therefore into $\PP(V)\times Z$, for some
$G$-module~$V$. Replacing $\Xx$ by its closure and taking a
pseudo-free desingularization, we may assume that $\Xx$ is
pseudo-free and $\pi$ is a projective morphism. By
Theorem~\ref{reg.var}, the fibers of $\pi$ are smooth projective
toroidal varieties. Shrinking $Z$ if necessary, we obtain that the
$G$-orbits of non-maximal dimension in $\Xx$ form a divisor with
normal crossings $\partial\Xx=\Dd_1\cup\dots\cup\Dd_k$ whose
components $\Dd_i$ are smooth over~$Z$.

Morally, an equivariant version of Kodaira--Spencer theory should
imply that all fibers of $\pi$ are $G$-isomorphic, which should
complete the proof. An alternative argument uses nested Hilbert
schemes~\cite{Hilb}.

Let $X$ be any fiber of~$\pi$, with $\partial{X}=D_1\cup\dots\cup
D_k$, $D_i=\Dd_i\cap X$. Applying a suitable Veronese map, we
satisfy a technical condition that the restriction map
$V^{*}\to\Ho^0(X,\Lin{1})$ is surjective.

The nested Hilbert scheme $\Hilb$ parametrizes tuples
$(Y,Y_1,\dots,Y_k)$ of projective subvarieties in~$\PP(V)$ having
the same Hilbert polynomials as $X,D_1,\dots,D_k$. The varieties
$\Xx,\Dd_1,\dots,\Dd_k$ are obtained as the pullbacks under
$Z\to\Hilb$ of the universal families
$\Yy,\Yy_1,\dots,\Yy_k\to\Hilb$. The groups $\GL(V)$ and $G$ act
on $\Hilb$ in a natural way, so that $\Hilb^G$ parametrizes tuples
of $G$-subvarieties. Since the centralizer $\GL(V)^G$ of $G$ maps
$G$-subvarieties to $G$-isomorphic ones, it suffices to prove that
the $\GL(V)^G$-orbit of $(X,D_1,\dots,D_k)$ is open in~$\Hilb^G$.

This is done by considering tangent spaces. Let
$\Nn_Z,\Nn_{Z/Z_i}$ denote the normal bundles to $Z$ in~$\PP(V)$,
resp.~to $Z_i$ in~$Z$. Then
$T_{(X,D_1,\dots,D_k)}\Hilb=\Ho^0(X,\Nn)$, where
$\Nn\subset\Nn_X\oplus\Nn_{D_1}\oplus\dots\oplus\Nn_{D_k}$ is
formed by tuples $(\xi,\xi_1,\dots,\xi_k)$ of normal vector fields
such that $\xi|_{D_i}=\xi_i\bmod\Nn_{X/D_i}$, $i=1,\dots,k$.
(These vector fields define infinitesimal deformations of
$X,D_1,\dots,D_k$, so that the deformation of $D_i$ is determined
by the deformation of $X$ modulo a deformation inside~$X$.) There
are exact sequences
%*
\begin{gather*}
0\longrightarrow\Tt_X(-\log\partial{X})\longrightarrow
\Tt_{\PP(V)}|_X\longrightarrow\Nn\longrightarrow0 \\
0\longrightarrow\Oo_X\longrightarrow V\otimes\Lin[X]{1}
\longrightarrow\Tt_{\PP(V)}|_X\longrightarrow0
\end{gather*}
%*
Taking cohomology yields
%*
\begin{gather*}
\Ho^0\left(X,\Tt_{\PP(V)}\right)\longrightarrow
T_{(X,D_1,\dots,D_k)}\Hilb
\longrightarrow\Ho^1(X,\Tt_X(-\log\partial{X}))=0 \\
V\otimes\Ho^0(X,\Lin{1})\longrightarrow\Ho^0\left(X,\Tt_{\PP(V)}\right)
\longrightarrow\Ho^1(X,\Oo_X)=0
\end{gather*}
%*
(The first cohomologies vanish by Theorem~\ref{H(log.tan)}
and~\cite{fund(unirat)}, since $X$ is a smooth projective rational
variety.) Hence the differential of the orbit map
%*
\begin{equation*}
\gl(V)\iso V\otimes V^{*}\longrightarrow V\otimes\Ho^0(X,\Lin{1})
\longrightarrow\Ho^0\left(X,\Tt_{\PP(V)}\right)\longrightarrow
T_{(X,D_1,\dots,D_k)}\Hilb
\end{equation*}
%*
is surjective. By linear reductivity of~$G$, the composite map
%*
\begin{equation*}
\gl(V)^G\longrightarrow T_{(X,D_1,\dots,D_k)}(\Hilb^G)\subseteq
\left(T_{(X,D_1,\dots,D_k)}\Hilb\right)^G
\end{equation*}
%*
is surjective as well. Hence $(X,D_1\dots,D_k)$ is a smooth point
of $\Hilb^G$ and $\GL(V)^G(X,D_1\dots,D_k)$ is open.
\end{proof}

Cohomology rings of smooth complete toroidal varieties (over
$\kk=\CC$) were computed by
Bifet--de~Concini--Procesi~\cite{H(reg.emb)}, see also
\cite{H_G(wonder)} for toroidal completions of symmetric spaces.
By Corollary~\ref{A=H}, cohomology coincides with the Chow ring in
this situation. The most powerful approach is through equivariant
cohomology or equivariant intersection theory of Edidin--Graham,
see~\cite{H_G&A_G}. In particular, Chow (or cohomology) rings of
smooth (complete) toric varieties and flag varieties are easily
computed in this way \cite[I.4]{H(reg.emb)},
\cite[2,~3]{conv.pol}, \cite{H_G&A_G}, cf.~\ref{intersect}.

The local structure of toroidal varieties can be refined in order
to obtain a full description for the closures of generic flats.
\begin{proposition}[{\cite[8.3]{inv.mot}}]
The closure of a generic twisted flat in a toroidal variety $X$ is
a normal toric variety whose fan is the $W_X$-span of the fan
of~$X$.
\end{proposition}
\begin{proof}
It suffices to choose the toric slice $Z$ in
Theorem~\ref{loc.str.tor} in such a way that the open $A$-orbit in
$Z$ is a generic (twisted) flat $F_{\alpha}$. Then
$\overline{Z}=\overline{F}_{\alpha}$ (the closure in~$X$), so that
Theorem~\ref{loc.str.tor} and Proposition~\ref{cl(flats)} imply
the claim.

If $X$ is smooth and $T^{*}X$ is symplectically stable, then the
conormal bundle to generic $U$-orbits extends to a trivial
subbundle $\X\times\ab^{*}\embeds T\X(\log\partial{X})$, the
trivializing sections being $d\ef{\lambda}/\ef{\lambda}$,
$\lambda\in\RG$. The logarithmic moment map restricts to
$\Phi:\X\times\ab^{*}\to\ab\oplus\Ru\p$,
cf.~Lemma~\ref{mom(conorm)}. It follows that $\X\iso
P\itimes{L}Z$, where $Z=\pi_X\Phi^{-1}(\lambda)$,
$\lambda\in\ab^{\pr}$, and $F_{\alpha}$ is the open $L$-orbit in
$Z$ for $\forall\alpha\in\Phi^{-1}(\lambda)\cap T^{*}\HS$.

If $X$ is singular, then it admits a toroidal resolution of
singularities $\nu:X'\to X$. Then
$\X':=\nu^{-1}(\X)=X'\setminus\bigcup_{D\in\Dd^B}D\iso
P\itimes{L}Z'$ and $Z'\supseteq F_{\alpha}$. The map
$\Phi:\X'\times\ab^{*}\to\ab\oplus\Ru\p$ descends to~$\X$, because
$\kk[\X']=\kk[\X]$. Thus one may put $Z=\nu(Z')$.

If $T^{*}X$ is not symplectically stable, then passing to affine
cones and back to projectivizations yields $Z$ such that the open
$L$-orbit in $Z$ is a twisted flat.
\end{proof}
\begin{example}
If $X$ is a toroidal $G\times G$-embedding of $G$, then $T$ is a
flat and $F=\overline{T}$ is a toric variety whose fan is the
$W$-span of the fan of~$X$ (in the antidominant Weyl chamber),
cf.~Proposition~\ref{cl(T)}. For instance, if
$X=\overline{G}\subseteq\PP(\LO(V))$ for a faithful projective
representation $G:\PP(V)$ with regular highest weights, then the
fan of $F$ is formed by the duals to the corner cones of the
weight polytope $\Pp(V)$, and the fan of $X$ is its antidominant
part (see~\ref{monoids}).
\end{example}
\begin{example}
Consider the variety of complete conics
$X\subset\PP(\Sym^2(\kk^3)^{*})\times\PP(\Sym^2\kk^3)$ from
Example~\ref{conics}. The set $F=\{([q],[q^{\vee}])\mid q\text{
diagonal, }\det q\neq0\}$ is a flat. Using the Segre embedding
$\PP(\Sym^2(\kk^3)^{*})\times\PP(\Sym^2\kk^3)\embeds
\PP(\Sym^2(\kk^3)^{*}\otimes\Sym^2\kk^3)$ and observing that the
$T$-weights occurring in the weight decomposition of $q\otimes
q^{\vee}$ are $2(\eps_i-\eps_j)$, we conclude that the fan of
$\overline{F}$ is the set of all Weyl chambers of $G=\SL_3(\kk)$
together with their faces, while the fan of $X$ consists of the
antidominant Weyl chamber and its faces.
\end{example}

\section{Wonderful varieties}
\label{wonderful}

In the study of a homogeneous space $\HS$ it is useful to consider
its equivariant completions. The reason is that properties of
$\HS$ and of related objects (subvarieties and their intersection,
functions, line bundles and their sections, etc) often become
apparent ``at infinity'', and equivariant completions of $\HS$
take into account the points at infinity. Also, complete varieties
behave better than non-complete ones from various points of view
(e.g., in intersection theory).

Among all equivariant completions of a spherical homogeneous space
$\HS$ one distinguishes two opposite classes. Toroidal completions
have nice geometry (see~\ref{toroidal}) and a universal property:
each equivariant completion of $\HS$ is dominated by a toroidal
one. On the other hand, simple completions of~$\HS$ (i.e., those
having a unique closed orbit) are the most ``economical'' ones:
their boundaries are ``small''. Simple completions exist iff the
valuation cone $\Vv$ is strictly convex.

These two classes intersect in a unique element, called the
wonderful completion.
\begin{definition}
A spherical subgroup $H\subseteq G$ is called \emph{sober} if
$N(H)/H$ is finite or, equivalently, if $\Vv(G/H)$ is strictly
convex.

The \emph{wonderful embedding} of $\HS=G/H$ is the unique toroidal
simple complete $G$-embedding $X\embof\HS$, defined by the colored
cone $(\Vv,\emptyset)$, provided that $H$ is sober.
\end{definition}
The wonderful embedding has a universal property: for any toroidal
completion $X'\embof\HS$ and any simple completion $X''\embof\HS$,
there exist unique proper birational $G$-morphisms $X'\to X\to
X''$ extending the identity map on~$\HS$.

Wonderful embeddings were first introduced by De Concini and
Procesi \cite{comp.symm} for symmetric spaces. Their remarkable
properties were studied by many researchers (see below) mainly in
characteristic zero, though some results in special cases, e.g.,
for symmetric spaces~\cite{wonder(symm)}, are obtained in
arbitrary characteristic. For simplicity, we assume $\ch\kk=0$
from now on.

Every spherical subgroup $H\subseteq G$ is contained in the
smallest sober overgroup $H\cdot N(H)^0$. This stems, e.g., from
the following useful lemma.
\begin{lemma}\label{N(sph)}
If $H\subseteq G$ is a spherical subgroup, then
$N(H)=N(\overline{H})$ for any intermediate subgroup
$\overline{H}$ between $H$ and~$N(H)$.
\end{lemma}
\begin{proof}
As $N(H)/H$ is Abelian, we have $N(H)\subseteq N(\overline{H})$.
In particular, $N(H)=N(H^0)$. To prove the converse inclusion, we
may assume w.l.o.g.\ that $H$ is connected and $\br+\h=\g$. Then
the right multiplication by $N(\overline{H})$ preserves
$BH=B\overline{H}$, the unique open $(B\times H)$-orbit in~$G$.
Hence the $N(\overline{H})$-action on $\kk(G)$ by right
translations of an argument preserves~$\kk[G]^{(B\times H)}$ (=the
set of regular functions on $G$ invertible on~$BH$). Since this
action commutes with the $G$-action by left translations, it
preserves~$\kk[G]^{(H)}$, whence~$\kk(G/H)$, too. Hence
$N(\overline{H})$ acts on $G/H$ by $G$-automorphisms, i.e., is
contained in~$N(H)$.
\end{proof}

Now let $H\subseteq G$ be a sober subgroup and $X$ the wonderful
embedding of $\HS=G/H$. The local structure theorem reveals the
orbit structure and local geometry of~$X$: by
Theorem~\ref{loc.str.tor} there are an affine open chart
$\X=X\setminus\bigcup_{D\in\Dd^B}D$ and a closed subvariety
$Z\subset\X$ such that $\X$ is stable under $P=P(\HS)$, the Levi
subgroup $L\subset P$ leaves $Z$ stable and acts on it via the
quotient torus $A=L/L_0$, $\X\iso P\itimes{L}Z\iso\Ru{P}\times Z$,
and each $G$-orbit of $X$ intersects $Z$ in an $A$-orbit. Actually
$\X$ is the unique $B$-chart of $X$ intersecting all $G$-orbits.

The affine toric variety $Z$ is defined by the cone~$\Vv$, so that
$\kk[Z]=\kk[{\Vv^{\vee}\cap\RG}]$, where $\RG=\RG(\HS)=\Ch(A)$.
The orbits (of $A:Z$ or of $G:X$) are in an order-reversing
bijection with the faces of~$\Vv$, and each orbit closure is the
intersection of invariant divisors containing the orbit. If $\Vv$
is generated by a basis of $\RG^{*}$, then $Z\iso\AAA^r$ with the
natural action of $A\iso(\kk^{\times})^r$; the eigenweight set for
$A:Z$ is~$\Pi_{\HS}^{\min}$. Generally, since $\Vv$ is simplicial
(Theorem~\ref{W_X&Z(X)}), one deduces that $Z\iso\AAA^r/\Gamma$
with the natural action of $A\iso(\kk^{\times})^r/\Gamma$, where
$\Gamma\iso\RG^{*}/N$ is the common kernel of all $\lambda\in\RG$
in $(\kk^{\times})^r=N\otimes\kk^{\times}$, the sublattice
$N\subseteq\RG^{*}$ being spanned by the indivisible generators of
the rays of~$\Vv$.

In particular, $X$~is smooth iff $\Vv$ is generated by a basis of
$\RG^{*}$ iff $\RG=\ZZ\Delta_{\HS}^{\min}$. It is a delicate
problem to characterize the (sober) spherical subgroups
$H\subseteq G$ such that the wonderful embedding $X\embof\HS=G/H$
is smooth.

Note that $N(H)/H=\Aut_G\HS$ acts on a finite set~$\Dd^B$.
\begin{definition}
A spherical subgroup $H\subseteq G$ is called \emph{very sober} if
$N(H)/H$ acts on $\Dd^B$ effectively. (In particular, $H$~is
sober, because $(N(H)/H)^0$ leaves $\Dd^B$ pointwise fixed.) The
\emph{very sober hull} of $H$ is the kernel $\overline{H}$ of
$N(H):\Dd^B$.
\end{definition}
\begin{remark}
It is easy to deduce from Lemma~\ref{N(sph)} that $\overline{H}$
is the smallest very sober subgroup of $G$ containing~$H$. The
colored space $\overline\ES=\ES(G/\overline{H})$ is identified
with $\ES/(\Vv\cap-\Vv)$, the valuation cone is
$\overline\Vv=\Vv/(\Vv\cap-\Vv)$, and the set of colors
$\overline\Dd^B$ is identified with $\Dd^B$ via pullback.

Observe that $\overline{H}$ is the kernel of $N(H):\Ch(H)$
\cite[7.4]{Aut&root}. Indeed, (a multiple of) each $B$-stable
divisor $\delta$ on $\HS$ is defined by an equation
$\eta\in\kk(G)^{(B\times H)}_{(\lambda,\chi)}$, and each
$\chi\in\Ch(H)$ arises in this way (because every $G$-line bundle
$\ind[G/H]{\chi}$ has a rational $B$-eigensection). The right
multiplication by $n\in N(H)$ maps $\eta$ to
$\eta'\in\kk(G)^{(B\times H)}_{(\lambda,\chi')}$, the equation of
$\delta'=n(\delta)$, where $\chi'(h)=\chi(n^{-1}hn)$. Since
$\kk(\HS)^B=\kk$, we have $\chi'=\chi\iff\eta'/\eta=\const\iff
\delta'=\delta$.

In particular, $\overline{H}\supseteq Z_G(H)$.
\end{remark}

\begin{theorem}[{\cite[7.6, 7.2]{Aut&root}}]\label{sober=>smooth}
If $H$ is very sober, then the wonderful embedding $X\embof G/H$
is smooth. In particular, $X$~is smooth if $N(H)=H$.
\end{theorem}
\begin{remark}
If all simple factors of $G$ are isomorphic to~$\PSL_{n_i}$, then
very soberness is also a necessary condition for $X$ be smooth
\cite[7.1]{sph(A)}. This is not true in general:
$S^{n-1}=\SO_n/\SO_{n-1}$ and $\SL_4/\Sp_4$ are symmetric spaces
of rank~$1$, hence their wonderful embeddings are smooth
(Proposition~\ref{emb(r=1)}), while $\overline{\SO_{n-1}}=\Or_n$,
$\overline{\Sp_4}=\Sp_4\cdot Z(\SL_4)$.
\end{remark}
\begin{proof}
By Theorem~\ref{root->cent},
$S_{\HS}=\bigcap_{\alpha\in\Delta_{\HS}^{\min}}\Ker\alpha\embeds
\Aut_G\HS=N(H)/H$. It suffices to show that $S_{\HS}$ fixes all
colors; then $S_{\HS}=\{\1\}$, i.e.,
$\Delta_{\HS}^{\min}$~spans~$\RG$.

Take any $D\in\Dd^B$. Replacing $D$ by a multiple, we may assume
that $\Lin{D}$ is $G$-linearized. Consider the total space
$\widehat\HS=\widehat{G}/\widehat{H}$ of~$\Lin{D}^{\times}$, where
$\widehat{G}=G\times\kk^{\times}$, cf.~Remark~\ref{aff.cone}.
Using the notation of Remark~\ref{aff.cone}, we have
%*
\begin{equation*}
0\longrightarrow\RG\longrightarrow\widehat\RG
\longrightarrow\ZZ\longrightarrow0,
\end{equation*}
%*
$\Vv=\widehat\Vv/(\widehat\Vv\cap-\widehat\Vv)$, and
$\Delta_{\widehat\HS}^{\min}=\Delta_{\HS}^{\min}$. Therefore
$S_{\HS}=S_{\widehat\HS}/\kk^{\times}$.

However the pullback $\widehat{D}\subset\widehat\HS$ of $D$ is
principal. Since $S_{\widehat\HS}$ multiplies the equation of
$\widehat{D}$ by scalars, it leaves $\widehat{D}$ stable, whence
$S_X$ leaves $D$ stable.
\end{proof}

If $N(H)=H$, then $\HS\iso G[\h]$, the orbit of $\h$ in
$\Gr_k(\g)$, $k=\dim\h$. The closure
$X(\h)=\overline{G[\h]}\subseteq\Gr_k(\g)$ is called the
\emph{Demazure embedding} of~$\HS$.
\begin{proposition}[{\cite[1.4]{W_X(sph)}}]\label{Demazure}
The wonderful embedding $X$ is the normalization of~$X(\h)$.
\end{proposition}
\begin{proof}
The decomposition $\g=\Ru\p\oplus\ab\oplus\h$ yields
$\h=\lv_0\oplus\langle e_{-\alpha}+\xi_{\alpha}\mid
\alpha\in\Delta^{+}\setminus\Delta_L^{+}\rangle$, where
$\xi_{\alpha}\in\Ru\p\oplus\ab$ is the projection of
$-e_{-\alpha}$ along~$\h$. Hence
%*
\begin{equation*}
\widehat\h=\widehat{\lv_0}\wedge
\bigwedge_{\alpha\in\Delta^{+}\setminus\Delta_L^{+}}
(e_{-\alpha}+\xi_{\alpha})=\widehat\s+ \text{terms of higher
$T$-weights}
\end{equation*}
%*
where $\widehat\q\in\E^{\bullet}\g$ denotes a generator of
$[\q]\in\Gr(\g)$, $\s=\lv_0\oplus\Ru{\p^{-}}$, and the weights of
other terms differ from that of $\widehat\s$ by
$\sum(\alpha_i+\beta_i)$,
$\alpha_i,\beta_i\in\Delta^{+}\setminus\Delta_L^{+}$ or
$\beta_i=0$.

Let $Z(\h)$ be the closure of $T[\h]$ in the affine chart defined
by non-vanishing of the highest weight covector dual
to~$\widehat\s$. It is an affine toric variety with the fixed
point~$[\s]$. Thus $Y=G[\s]\subset X(\h)$ is a closed orbit. The
local structure theorem in a neighborhood of $[\s]$ provides a
$B$-chart $\X(\h)\subset X(\h)$, $\X(\h)\iso\Ru{P}\times Z(\h)$.
Note that for any $[\q]\in Z(\h)\setminus T[\h]$ the subalgebra
$\q$ is transversal to $\Ru\p\oplus\ab$ while
$\n(\q)\cap\ab\neq0$, whence $\dim G_{[\q]}>\dim H$. It follows
that $\X(\h)$ intersects no colors, i.e., $X(\h)$ is toroidal in a
neighborhood of~$Y$.

On the other hand, every smooth toroidal embedding of $\HS$ maps
to $X(\h)$ by Theorem~\ref{reg.var}.
%Brion shows \cite[1.3]{W_X(sph)} that every
%$G$-orbit in $X(\h)$ contains a point $[\q]$ such that
%$\g=\Ru\p\oplus\ab\oplus\q$. As above,
%$\overline{T[\q]}\ni[\s]$. Hence $Y$ is the unique closed
%orbit in~$X(\h)$.
It follows that $\widetilde{X(\h)}$ is simple, whence wonderful.
\end{proof}
It is an open question whether the Demazure embedding is always
smooth.

If $N(H)\neq H$, then the normalization of $X(\h)$ is the
wonderful embedding of~$G/N(H)$.

\begin{example}\label{wonder(symm)}
Let $G$ be an adjoint semisimple group and $H=G^{\sigma}$ a
symmetric subgroup. Here $N(H)=H$. We have
$\RG(\HS)=\Ch(T/T^{\sigma})=\{\mu-\sigma(\mu)\mid\mu\in\Ch(T)\}$,
where $T$ is a $\sigma$-stable maximal torus such that $T_1$ is a
maximal $\sigma$-split torus. Hence $\RG(\HS)$ is the root lattice
of~$2\Delta_{\HS}$. Since $\Vv(\HS)$ is the antidominant Weyl
chamber of $\Delta_{\HS}^{\vee}$ in $\RG(\HS)^{*}\otimes\QQ$ (by
Theorem~\ref{c.d.(symm)}), $\Delta_{\HS}^{\min}$~is the reduced
root system associated with~$2\Delta_{\HS}$. It follows that the
wonderful completion $X$ is smooth in this case.

Wonderful completions of symmetric spaces were studied in
\cite{comp.symm}, \cite{wonder(symm)}. In particular, a geometric
realization for a wonderful completion as an embedded projective
variety was constructed. Let $\lambda$ be a dominant weight of
$\widetilde{G}$ such that $\sigma(\lambda)=-\lambda$ and
$\lambda\in\intr\DWC(\Delta_{\HS}^{+})$. There exists a unique (up
to proportionality) $\widetilde{G}^{\sigma}$-fixed vector $v'\in
V^{*}(\lambda)$. Then
$X'=\overline{G[v']}\subseteq\PP(V^{*}(\lambda))$ is the wonderful
embedding of $G[v']\iso\HS$.

Indeed, a natural closed embedding
$\PP(V^{*}(\lambda))\embeds\PP(V^{*}(2\lambda))$ (given by the
multiplication $V^{*}(\lambda)\otimes V^{*}(\lambda)\to
V^{*}(2\lambda)$ in~$\kk[\widetilde{G}]$) identifies $X'$ with
$X''=\overline{G[v'']}$, where $v''\in V^{*}(2\lambda)$ is a
unique $\widetilde{G}^{\sigma}$-fixed vector. As $X''$ is a simple
projective embedding of~$G[v'']$, the natural map $\HS\to G[v'']$
extends to $X\to X''$. On the other hand, the homomorphism
$V^{*}(\lambda)\otimes V^{*}(\lambda)\to V^{*}(2\lambda)$ maps
$\omega$ to~$v''$, where $\omega$ is defined by~\eqref{H-fixed}.
Let $Z''$ be the closure of $T[v'']$ in the affine chart of
$\PP(V^{*}(2\lambda))$ defined by non-vanishing of the highest
weight covector of weight~$2\lambda$. From~\eqref{H-fixed} it is
easy to deduce that $Z''\iso\AAA^r$ is acted on by $T$ via the
eigenweight set $\Pi_{\HS}^{\min}$ and the closed orbit
$G[v_{-2\lambda}]$ is transversal to $Z''$ at~$[v_{-2\lambda}]$.
Hence $Z\isoto Z''$, $\X\isoto PZ''\iso\Ru{P}\times Z''$, and
finally $X\isoto X''\iso X'$. (A similar reasoning shows
$X\iso\overline{G[\omega]}$. A slight refinement carries over the
construction to positive characteristic~\cite{wonder(symm)}.)

Another model for the wonderful completion is the Demazure
embedding. First note that $\h=\lv_0\oplus\langle
e_{\alpha}+e_{\sigma(\alpha)}\mid
\alpha\in\Delta^{+}\setminus\Delta_L^{+}\rangle$. Arguing as in
the proof of Proposition~\ref{Demazure}, we see that
$Z(\h)=\overline{T[\h]}\iso\AAA^r$ is acted on by $T$ with the
eigenweights $\alpha-\sigma(\alpha)$, $\alpha\in\Pi$, and
$Y=G[\s]$ is transversal to $Z(\h)$ at~$[\s]$. This yields
$\Cc_Y=\Vv$. Now the Luna--Vust theory together with the
description of the colored data for symmetric spaces implies that
$X(\h)$ is wonderful. The varieties $X(\h)$ were first considered
by Demazure in the case, where $G=\PSL_n(\kk)$ and $H$ is the
projective orthogonal or symplectic group~\cite{Dem}.
\end{example}

Using the Demazure embedding, Brion computed the canonical class
of any spherical variety.
\begin{proposition}[{\cite[1.6]{W_X(sph)}}]\label{can(sph)}
Suppose $X$ is a spherical variety with the open orbit $\HS\iso
G/H$. Consider the $G$-morphism $\phi:\HS\to\Gr_k(\g)$,
$\phi(\bp)=[\h]$, $k=\dim H$. Then a canonical divisor of $X$ is
%*
\begin{equation*}
K_X=-\sum_iD_i-\overline{\phi^{*}\Hh}
=-\sum_iD_i-\sum_{D\in\Dd^B}m_DD
\end{equation*}
%*
where $D_i$ runs over all $G$-stable prime divisors in~$X$,
$\Hh$~is a hyperplane section of $X(\h)$ in~$\PP(\E^k\g)$, and
$m_D\in\NN$.
\end{proposition}
Explicit formul{\ae} for $m_D$ are given in \cite[4.2]{B-curves}.
\begin{proof}
Removing all $G$-orbits of codimension $>1$, we may assume that
$X$ is smooth and toroidal. Then by Theorem~\ref{reg.var},
$\phi$~extends to~$X$, and we have an exact sequence
%*
\begin{equation*}
0\longrightarrow\phi^{*}\Ee\longrightarrow\Oo_X\otimes\g
\longrightarrow\Tt_X(-\log\partial{X})\longrightarrow0
\end{equation*}
%*
where $\Ee$ is the tautological vector bundle on~$\Gr_k(\g)$.
Taking the top exterior powers yields
$\omega_X\otimes\Lin[X]{\partial{X}}\iso\E^k\phi^{*}\Ee
=\Lin[X]{-\phi^{*}\Hh}$, whence the first expression for~$K_X$. If
$\Hh$ is defined by a covector in $(\E^k\g^{*})^{(B)}$ dual
to~$\widehat\s$, then $\X(\h)=X(\h)\setminus\Hh$ intersects all
$G$-orbits in open $B$-orbits, whence $\phi^{*}\Hh=\sum m_DD$ with
$m_D>0$ for $\forall D\in\Dd^B$.
\end{proof}
Using the characterization of ample divisors on complete spherical
varieties (Corollary~\ref{glob&ample(sph)}), one deduces that
certain smooth wonderful embeddings (e.g., flag varieties,
wonderful completions of symmetric spaces, of affine spherical
spaces of rank~$1$) are Fano varieties (i.e., anticanonical
divisor is ample).

Smooth wonderful embeddings can be characterized intrinsically by
the configuration of $G$-orbits.
\begin{theorem}[\cite{wonder=>sph}]
A smooth complete $G$-variety $X$ is a wonderful embedding of a
spherical homogeneous space iff it satisfies the following
conditions:
\begin{enumerate}
\item\label{open.orb} $X$~contains a dense open orbit~$\HS$.
\item\label{norm.cross} $X\setminus\HS$ is a divisor with normal
crossings, i.e., its components $D_1,\dots,D_r$ are smooth and
intersect transversally. \item\label{orb(wonder)} For each tuple
$1\leq i_1<\dots<i_k\leq r$, the set $D_{i_1}\cap\dots\cap
D_{i_k}\setminus\bigcup_{i\neq i_1,\dots,i_k}D_i$ is a single
$G$-orbit. (In particular, it is non-empty.)
\end{enumerate}
\end{theorem}
$G$-varieties satisfying the conditions of the theorem are called
\emph{wonderful varieties}.
\begin{proof}[Sketch of a proof]
Smooth wonderful embeddings obviously satisfy the conditions
\ref{open.orb}--\ref{orb(wonder)}, as a particular case of
Theorem~\ref{reg.var}: the toric slice $Z\iso\AAA^r$ is
transversal to all orbits and the $G$-stable prime divisors
intersect it in the coordinate hyperplanes.

To prove the converse, consider the local structure of $X$ in a
neighborhood of the closed orbit $Y$ which is provided by an
embedding of $X$ into a projective space. Let $P=L\Ru{P}$ be a
Levi decomposition of $P=P(Y)$. There is a $B$-chart
$\X\iso\Ru{P}\times Z$ such that $Z$ is a smooth $L$-stable
locally closed subvariety intersecting $Y$ transversally at the
unique $P^{-}$-fixed point~$z$. It is easy to see that a general
dominant one-parameter subgroup $\gamma\in\CoCh(Z(L))$ contracts
$\X$ to~$z$. Hence $Z$ is $L$-isomorphic to~$T_zZ$.

Consider the wonderful subvarieties $X_i=\bigcap_{j\neq i}D_j$ of
rank~$1$ and let $\{\lambda_i\}$ be the $T$-weights of~$T_z(Z\cap
X_i)$, $i=1,\dots,r$. Since $T_zZ=\bigoplus T_z(Z\cap X_i)$, it
suffices to prove that $\lambda_1,\dots,\lambda_r$ are linearly
independent.

The latter is reduced to the cases $r=1$~or~$2$. Indeed, if we
already know that $X_i$ and $X_{ij}=\bigcap_{k\neq i,j}D_k$ are
wonderful embeddings of spherical spaces, then
$\Pi_{X_i}^{\min}=\{\lambda_i\}$ and
$\Pi_{X_{ij}}^{\min}=\{\lambda_i,\lambda_j\}$. Thus the
$\lambda_i$'s are positive linear combinations of positive roots
located at obtuse angles to each other. This implies the linear
independence.

The case $r=1$ stems from Proposition~\ref{emb(sph,r=1)}.

The case $r=2$ can be reduced to $G=\SL_2$. Indeed, assuming that
$\lambda_1,\lambda_2$ are proportional, we see that $c(X)=r(X)=1$.
By Proposition~\ref{r=1}, $\HS$~is obtained from a $3$-dimensional
homogeneous $\SL_2$-space by parabolic induction. Let us describe
the colored hypercone $(\Cc_Y,\Dd^B_Y)$.

Since $T_zZ$ is contracted to $0$ by~$\gamma$, we have
$\RG(X)=\ZZ\lambda$, $\lambda_i=h_i\lambda$, where $\lambda$ is
dominant and $h_1,h_2$ are coprime positive integers. W.l.o.g.\
$\ell_1h_1-\ell_2h_2=1$ for some $\ell_1,\ell_2\in\NN$. Consider
$T_z(Z\cap X_i)$ as coordinate axes in $T_zZ\iso Z$ and extend the
respective coordinates to $f_1,f_2\in\kk(X)^{(B)}$.  Then we may
put $\ef{\lambda}=f_2^{\ell_2}/f_1^{\ell_1}$, and
$\kk(X)^B=\kk\left(f_2^{h_1}/f_1^{h_2}\right)$. We have the
following picture for $(\Cc_Y,\Dd^B_Y)$ ($f_2^{h_1}/f_1^{h_2}$~is
regarded as affine coordinate on $\PP^1$, colors in $\Dd_Y^B$ are
marked by bold dots):
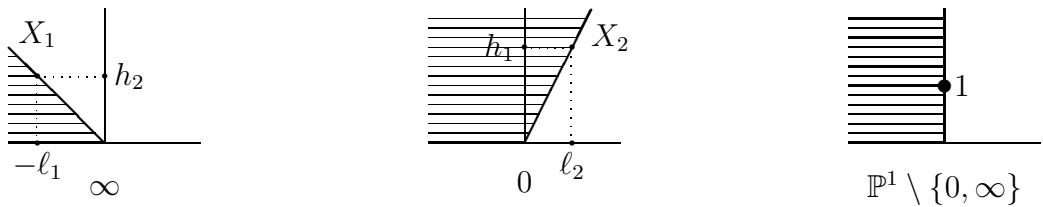
\begin{center}
%TeXCAD Picture [wonder-1.pic]. Options:
%\grade{\off}
%\emlines{\off}
%\epic{\off}
%\beziermacro{\off}
%\reduce{\on}
%\snapping{\off}
%\quality{2.00}
%\graddiff{0.01}
%\snapasp{1}
%\zoom{10.0000}
\unitlength .7ex % = .7pt
\linethickness{0.4pt}
\begin{picture}(20,20)(0,0)
\put(10,6){\line(1,0){10}} \put(10,2){\makebox(0,0)[ct]{$\infty$}}
\thicklines{} \put(0,6){\line(1,0){10}}
\put(10.00,6.00){\line(-1,1){10.00}}
\put(1,16){\makebox(0,0)[lb]{$X_1$}} \thinlines{}
\put(10,6){\line(0,1){14}} \put(3,13){\circle*{.5}}
\put(10,13){\circle*{.5}} \put(3,6){\circle*{.5}}
\put(3,5){\makebox(0,0)[ct]{$-\ell_1$}}
\put(11,13){\makebox(0,0)[lc]{$h_2$}}
%\dottedline(3,13)(10,13)
\multiput(2.86,12.86)(.875,0){9}{{\rule{.2pt}{.2pt}}}
%\end
%\dottedline(3,13)(3,6)
\multiput(2.86,12.86)(0,-.875){9}{{\rule{.2pt}{.2pt}}}
%\end
\linethickness{0.05pt} \put(9,7){\line(-1,0){9}}
\put(8,8){\line(-1,0){8}} \put(7,9){\line(-1,0){7}}
\put(6,10){\line(-1,0){6}} \put(5,11){\line(-1,0){5}}
\put(4,12){\line(-1,0){4}} \put(3,13){\line(-1,0){3}}
\put(2,14){\line(-1,0){2}} \put(1,15){\line(-1,0){1}}
\end{picture}
\hfill
%TeXCAD Picture [wonder-2.pic]. Options:
%\grade{\off}
%\emlines{\off}
%\epic{\off}
%\beziermacro{\off}
%\reduce{\on}
%\snapping{\off}
%\quality{2.00}
%\graddiff{0.01}
%\snapasp{1}
%\zoom{10.0000}
\unitlength .7ex % = .7pt
\linethickness{0.4pt}
\begin{picture}(20,20)(0,0)
\put(10,6){\line(1,0){10}} \put(10,6){\line(0,1){14}}
\put(10,3){\makebox(0,0)[ct]{$0$}} \thicklines{}
\put(10,6){\line(1,2){7}} \put(0,6){\line(1,0){10}}
\put(17,17){\makebox(0,0)[lc]{$X_2$}} \thinlines
\put(15,6){\circle*{.5}} \put(15,16){\circle*{.5}}
\put(10,16){\circle*{.5}} \put(15,5){\makebox(0,0)[ct]{$\ell_2$}}
\put(9,16){\makebox(0,0)[rc]{$h_1$}}
%\dottedline(10,16)(15,16)
\multiput(9.86,15.86)(.8333,0){7}{{\rule{.2pt}{.2pt}}}
%\end
%\dottedline(15,6)(15,16)
\multiput(14.86,5.86)(0,.9091){12}{{\rule{.2pt}{.2pt}}}
%\end
\linethickness{0.05pt} \put(10.5,7){\line(-1,0){10.5}}
\put(11,8){\line(-1,0){11}} \put(11.5,9){\line(-1,0){11.5}}
\put(12,10){\line(-1,0){12}} \put(12.5,11){\line(-1,0){12.5}}
\put(13,12){\line(-1,0){13}} \put(13.5,13){\line(-1,0){13.5}}
\put(14,14){\line(-1,0){14}} \put(14.5,15){\line(-1,0){14.5}}
\put(15,16){\line(-1,0){15}} \put(15.5,17){\line(-1,0){15.5}}
\put(16,18){\line(-1,0){16}} \put(16.5,19){\line(-1,0){16.5}}
\end{picture}
\hfill
%TeXCAD Picture [wonder-pgp.pic]. Options:
%\grade{\off}
%\emlines{\off}
%\epic{\off}
%\beziermacro{\off}
%\reduce{\on}
%\snapping{\off}
%\quality{2.00}
%\graddiff{0.01}
%\snapasp{1}
%\zoom{10.0000}
\unitlength .7ex % = .7pt
\linethickness{0.4pt}
\begin{picture}(20,20)(0,0)
\put(10,6){\line(1,0){10}}
\put(10,3){\makebox(0,0)[ct]{$\PP^1\setminus\{0,\infty\}$}}
\thicklines{} \put(10,6){\line(0,1){14}} \put(0,6){\line(1,0){10}}
\put(10,12){\circle*{1.3}} \put(11,12){\makebox(0,0)[lc]{$1$}}
\linethickness{0.05pt} \thinlines{} \put(0,7){\line(1,0){10}}
\put(0,8){\line(1,0){10}} \put(0,9){\line(1,0){10}}
\put(0,10){\line(1,0){10}} \put(0,11){\line(1,0){10}}
\put(0,12){\line(1,0){10}} \put(0,13){\line(1,0){10}}
\put(0,14){\line(1,0){10}} \put(0,15){\line(1,0){10}}
\put(0,16){\line(1,0){10}} \put(0,17){\line(1,0){10}}
\put(0,18){\line(1,0){10}} \put(0,19){\line(1,0){10}}
\end{picture}
\end{center}

Since $\Dd^B_Y$ does not contain central colors, Propositions
\ref{par.ind(hyp)} and~\ref{par.ind(c.d.)} imply that $X$ is
induced from a wonderful $\SL_2$-variety. However it is easy to
see (e.g., from the classification in~\cite[\S5]{LVT}) that there
exist no $\SL_2$-germs with the colored data as above. (Luna uses
different arguments in~\cite{wonder=>sph}.)
\end{proof}

Wonderful varieties play a distinguished role in the study of
spherical homogeneous spaces, because they are the canonical
completions of these spaces having nice geometric properties. To a
certain extent this role is analogous to that of (generalized)
flag varieties in the theory of reductive groups. For symmetric
spaces this was already observed by de Concini and Procesi
\cite{comp.symm}. For general spherical spaces this principle was
developed by Brion, Knop, Luna, et al \cite{val.sph},
\cite{W_X(sph)}, \cite{Aut&root}, \cite{sph.solv},
\cite{big.cells}, \cite{sph(A)}.

In particular, wonderful varieties are applied to classification
of spherical subgroups. The strategy, proposed by Luna, is to
reduce the classification to very sober subgroups, which are
stabilizers of general position for wonderful varieties, and then
to classify the wonderful varieties.

By Theorem~\ref{sph.s.g.p}, there are no continuous families of
non-conjugate spherical subgroups, and even more:
\begin{proposition}
There are finitely many conjugacy classes of sober spherical
subgroups $H\subseteq G$.
\end{proposition}
\begin{proof}
Sober spherical subalgebras of dimension~$k$ form a locally closed
$G$-subvariety in~$\Gr_k(\g)$. Indeed, the set of spherical
subalgebras is open in the variety of $k$-dimensional Lie
subalgebras, and sober subalgebras are those having
$k$-dimensional orbits. Theorem~\ref{sph.s.g.p} implies that this
variety is a finite union of locally closed strata such that all
orbits in each stratum have the same stabilizer. But the isotropy
subalgebras are nothing else but the points of the strata. Hence
each stratum is a single orbit, i.e., there are finitely many
sober subalgebras, up to conjugation. As for subgroups, there are
finitely many ways to extend $H^0$ by a (finite) subgroup
in~$N(H^0)/H^0$.  \end{proof} Note that finiteness fails for
non-sober spherical subgroups:  $H^0$~can be extended by countably
many quasitori in~$N(H^0)/H^0$.

These results create an evidence that spherical subgroups should
be classified by some discrete invariants. Such invariants were
suggested by Luna, under the names of \emph{spherical systems} and
\emph{spherical homogeneous data} (Definition~\ref{sph.data}).
They are defined in terms of roots and weights of $G$ and
wonderful $G$-varieties of rank~$1$.

For spherical homogeneous spaces of rank~$1$, wonderful embeddings
are always smooth. Indeed, they are normal $G$-varieties
consisting of two $G$-orbits---a dense one and another of
codimension~$1$. Furthermore, spherical homogeneous spaces of
rank~$1$ are characterized by existence of a completion by
homogeneous divisors.
\begin{proposition}[\cite{comp.div}, \cite{r(sph)=1}]
\label{emb(sph,r=1)} The following conditions are equivalent:
\begin{enumerate}
\item\label{sph,r=1} $\HS=G/H$~is a spherical homogeneous space of
rank~$1$. \item\label{comp.hom} There exists a smooth complete
embedding $X\embof\HS$ such that $X\setminus\HS$ is a union of
$G$-orbits of codimension~$1$.
\end{enumerate}
Moreover, if $\HS$ is horospherical, then $X\setminus\HS$ consists
of two orbits and $X\iso G\itimes{Q}\PP^1$, where $Q\subseteq G$
is a parabolic acting on $\PP^1$ via a character. Otherwise
$X\setminus\HS$ is a single orbit and $X$ is a wonderful embedding
of~$\HS$.
\end{proposition}
\begin{proof}
The implication \ref{sph,r=1}$\implies$\ref{comp.hom} and the
properties of $X$ easily stem from the Luna--Vust theory: the
colored space $\ES$ is a line, whence there exists a unique smooth
complete toroidal embedding~$X$, which is obtained by adding two
homogeneous divisors (corresponding to the two rays of~$\ES$) if
$\Vv=\ES$ and is wonderful if $\Vv$ is a ray.

To prove \ref{comp.hom}$\implies$\ref{sph,r=1}, we consider the
local structure of $X$ in a neighborhood of a closed orbit $Y$.
Let $P=L\Ru{P}$ be a Levi decomposition of $P=P(Y)$. There is a
$B$-chart $\X\iso\Ru{P}\times Z$ such that $Z$ is an $L$-stable
affine curve intersecting $Y$ transversally at the unique
$P^{-}$-fixed point. Note that $T\subseteq L$ cannot fix $Z$
pointwise for otherwise $\HS^T$ would be infinite, which is
impossible. Hence $T:Z$ has an open orbit, whence~\ref{sph,r=1}.
\end{proof}
\begin{remark}\label{emb(r=1)}
A similar reasoning proves an embedding characterization of
arbitrary rank~$1$ spaces, due to Panyushev~\cite{r=1}: $r(\HS)=1$
iff there exists a complete embedding $X\embof\HS$ such that
$X\setminus\HS$ is a divisor consisting of closed $G$-orbits. Here
$Z$ is an affine $L$-stable subvariety with a pointwise $L$-fixed
divisor $Z\setminus\HS$ (provided that $Y$ is a generic closed
orbit), which readily implies that generic orbits of $L:Z$ are
one-dimensional, whence $r(X)=r(Z)=1$. On the other hand, it is
easy to construct a desired embedding $X$ for a homogeneous space
$\HS$ parabolically induced from $\SL_2$ modulo a finite subgroup,
cf.~Proposition~\ref{r=1}.
\end{remark}

Spherical homogeneous spaces $G/H$ of rank~$1$ were classified by
Akhiezer \cite{comp.div} and Brion~\cite{r(sph)=1}. It is easy to
derive the classification from a regular embedding of $H$ into a
parabolic $Q\subseteq G$. In the notation of
Theorem~\ref{c&r(nqaff)} we have an alternative: either
$r(M/K)=1$, $r_{M_{*}}(\Ru{Q}/\Ru{H})=0$, or vice versa.

In the first case $\Ru{H}=\Ru{Q}$, i.e., $G/H$~is parabolically
induced from an affine spherical homogeneous rank~$1$ space~$M/K$.
Except the trivial case $M/K\iso\kk^{\times}$ (where $H$ is
horospherical), $K$~is sober in $M$ and $H$ in~$G$.

In the second case $M=K=M_{*}$, and $\Ru{Q}/\Ru{H}\iso\Ru\q/\Ru\h$
is an $M$-module such that $(\Ru\q/\Ru\h)\setminus\{0\}$ is a
single $M$-orbit. Indeed, $\kk[\Ru\q/\Ru\h]^{U(M)}$ is generated
by one $B(M)$-eigenfunction, namely a highest weight covector in
$(\Ru\q/\Ru\h)^{*}$, whence $\Ru\q/\Ru\h$ is an HV-cone. Therefore
$M$ acts on $\Ru\q/\Ru\h\iso\kk^n$ as $\GL_n(\kk)$ or
$\kk^{\times}\cdot\Sp_n(\kk)$ and the highest weight of
$\Ru\q/\Ru\h$ is a negative simple root.

We deduce that every spherical homogeneous space of rank~$1$ is
either horospherical or parabolically induced from a
\emph{primitive} rank~$1$ space $\HS=G/H$ with $G$ semisimple and
$H$ sober. Primitive spaces are of the two types:
\begin{enumerate}
\item $H$~is reductive. \item $H$~is regularly embedded in a
maximal parabolic $Q\subseteq G$ which shares a Levi subgroup $M$
with $H$ and $\Ru\q/\Ru\h$ is a simple $M$-module of type
$\GL_n(\kk):\kk^n$ or $\kk^{\times}\cdot\Sp_n(\kk):\kk^n$
generated by a simple root vector.
\end{enumerate}
Primitive spherical homogeneous spaces of rank~$1$ are listed in
Table~\ref{wonder(r=1)}. Those of the \ordinal{1} type are easy to
classify, e.g., by inspection of Tables \ref{sph(s)},
\ref{sph(ss)}, and~\ref{symm}. We indicate the embedding $H\embeds
G$ by referring to Table~\ref{symm}~(\ref{sph(s)}) in the
(non-)symmetric case. Primitive spaces of the \ordinal{2} type are
classified by choosing a Dynkin diagram and its node corresponding
to a short simple root $\alpha$ which is adjacent to an extreme
node of the remaining diagram, the latter being of type $\AAa_l$
or~$\CCc_l$. The diagrams are presented in the column ``$H\embeds
G$'', with the white node corresponding to~$\alpha$.

The wonderful embeddings of spherical homogeneous space of
rank~$1$ are parabolically induced from those of primitive spaces.
The latter are easy to describe. For type~1 the construction of
Example~\ref{wonder(symm)} works whenever $N(H)=H$: the wonderful
embedding of $G/H$ is realized as
$X=\overline{G[v]}\subseteq\PP(V(\lambda))$, where $v\in
V(\lambda)^{(\widetilde{H})}$,
$\lambda\in\RG_{+}(\widetilde{G}/\widetilde{H}^0)$. If $N(H)\neq
H$, then $X$ is the projective closure of $Gv$ in
$\PP(V(\lambda)\oplus\kk)$. The simple minimal root of $X$ is the
generator of~$\RG_{+}(G/H)$.

For type~2 the wonderful embedding is
$X=G\itimes{Q}\PP(\Ru{Q}/\Ru{H}\oplus\kk)$. Indeed, $\Ru{Q}$~acts
on the $M$-module $\Ru{Q}/\Ru{H}$ by affine translations, whence
the projective closure of $\Ru{Q}/\Ru{H}$ consists of two
$Q$-orbits---the affine part and the hyperplane at infinity. Here
$\Pi_{X}^{\min}=\{w_M\alpha\}$.

Also the Demazure embedding is wonderful in all cases where
$N(H)=H$.

\begin{table}[!h]
\caption{Wonderful varieties of rank~$1$} \label{wonder(r=1)}
\begin{center}
\small\makebox[0pt]{\begin{tabular}{|c|c|c|c|c|c|} \hline \No &
$G$ & $H$ & $H\embeds G$ & $\Pi_{G/H}^{\min}$ &
Wonderful embedding \\
\hline

&&&&& $X=\{(x:t)\mid\det{x}=t^2\}$ \\
\tabitem & $\SL_2\times\SL_2$ & $\SL_2$ & diagonal &
$\omega+\omega'$ &
$\subset\PP(\LO_2\oplus\kk)$ \\
\cline{1-3}\cline{5-6}

\tabitem & $\PSL_2\times\PSL_2$ & $\PSL_2$ && $2\omega+2\omega'$
& $\PP(\LO_2)$ \\
\hline

\tabitem & $\SL_n$ & $\GL_{n-1}$ & symmetric \No\,\ref{SL/GLSL} &
$\omega_1+\omega_{n-1}$
& $\PP^n\times(\PP^n)^{*}$ \\
&&&
%TeXCAD Picture [SLGL.pic]. Options:
%\grade{\off}
%\emlines{\off}
%\epic{\off}
%\beziermacro{\off}
%\reduce{\on}
%\snapping{\off}
%\quality{2.00}
%\graddiff{0.01}
%\snapasp{1}
%\zoom{10.0000}
\unitlength .7ex % = .7pt
\linethickness{0.4pt}
\begin{picture}(20.66,1.67)(0,0)
\put(5,1){\circle*{1.33}} \put(0,1){\circle{1.33}}
\put(20,1){\circle*{1.33}} \put(5.67,1){\line(1,0){3.67}}
\put(.67,1){\line(1,0){3.67}}
\put(12.5,1){\makebox(0,0)[cc]{$\dots$}}
\put(15.67,1){\line(1,0){3.67}}
\end{picture} && \\
\hline

\tabitem & $\PSL_2$ & $\PO_2$ & symmetric \No\,\ref{SL/SO} &
$4\omega_1$ & $\PP(\sgl_2)$ \\
\hline

\tabitem & $\Sp_{2n}$ & $\Sp_2\times\Sp_{2n-2}$ & symmetric
\No\,\ref{Sp/SpSp} & $\omega_2$ &
$\Gr_2(\kk^{2n})$ \\
\hline

\tabitem & $\Sp_{2n}$ & $B(\Sp_2)\times\Sp_{2n-2}$ &
%TeXCAD Picture [CIII-S.pic]. Options:
%\grade{\off}
%\emlines{\off}
%\epic{\off}
%\beziermacro{\off}
%\reduce{\on}
%\snapping{\off}
%\quality{2.00}
%\graddiff{0.01}
%\snapasp{1}
%\zoom{10.0000}
\unitlength .7ex % = .7pt
\linethickness{0.4pt}
\begin{picture}(25.67,3)(0,0)
\put(0,1){\circle{1.33}} \put(5,1){\circle*{1.33}}
\put(20,1){\circle*{1.33}} \put(25,1){\circle*{1.33}}
\put(24.5,1.33){\line(-1,0){3.1}} \put(24.5,.67){\line(-1,0){3.2}}
\put(20.33,1){\makebox(0,0)[lc]{$<$}}
\put(5.67,1){\line(1,0){3.67}} \put(.67,1){\line(1,0){3.67}}
\put(15.67,1){\line(1,0){3.67}}
\put(12.5,1){\makebox(0,0)[cc]{$\dots$}}
\end{picture}
& $\omega_2$ & $\Fl_{1,2}(\kk^{2n})$ \\
\hline

&&&&& $X=\{(x:t)\mid(x,x)=t^2\}$ \\
\tabitem & $\SO_n$ & $\SO_{n-1}$ & symmetric & $\omega_1$ &
$\subset\PP^n$ \\
\cline{1-3}\cline{5-6}

\tabitem & $\SO_n$ & $\SG(\Or_1\times\Or_{n-1})$ &
\No\,\ref{SO/SOSO} & $2\omega_1$ & $\PP^{n-1}$ \\
\hline

&&&&& $X=\{(V_1,V_2)\mid V_1\subset V_1^{\perp}\}$ \\
\tabitem & $\SO_{2n+1}$ & $\GL_n\semitimes\E^2\kk^n$ &
%TeXCAD Picture [SOGL.pic]. Options:
%\grade{\off}
%\emlines{\off}
%\epic{\off}
%\beziermacro{\off}
%\reduce{\on}
%\snapping{\off}
%\quality{2.00}
%\graddiff{0.01}
%\snapasp{1}
%\zoom{10.0000}
\unitlength .7ex % = .7pt
\linethickness{0.4pt}
\begin{picture}(20.66,1.67)(0,0)
\put(0,1){\circle*{1.33}} \put(15,1){\circle*{1.33}}
\put(20,1){\circle{1.33}} \put(15.57,1.33){\line(1,0){3}}
\put(15.57,.67){\line(1,0){3}}
\put(19.67,1){\makebox(0,0)[rc]{$>$}}
\put(.67,1){\line(1,0){3.67}} \put(10.67,1){\line(1,0){3.67}}
\put(7.5,1){\makebox(0,0)[cc]{$\dots$}}
\end{picture}
& $\omega_1$ & $\subset\Fl_{n,2n}(\kk^{2n+1})$ \\
\hline

&&&&& $X=\{(x:t)\mid(x,x)=t^2\}$ \\
\tabitem & $\Spin_7$ & $\GGg_2$ & non-symmetric & $\omega_3$ &
$\subset\PP(V(\omega_3)\oplus\kk)$ \\
\cline{1-3}\cline{5-6}

\tabitem & $\SO_7$ & $\GGg_2$ & \No\,\ref{B3/G2} & $2\omega_3$ &
$\PP(V(\omega_3))$ \\
\hline

\tabitem & $\FFf_4$ & $\BBb_4$ & symmetric \No\,\ref{F4/B4} &
$\omega_1$ & \\
\hline

&&&&& $X=\{(x:t)\mid(x,x)=t^2\}$ \\
\tabitem & $\GGg_2$ & $\SL_3$ & non-symmetric & $\omega_1$ &
$\subset\PP(V(\omega_1)\oplus\kk)$ \\
\cline{1-3}\cline{5-6}

\tabitem & $\GGg_2$ & $N(\SL_3)$ & \No\,\ref{G2/A2} & $2\omega_1$
& $\PP(V(\omega_1))$ \\
\hline

\tabitem & $\GGg_2$ &
$\GL_2\semitimes(\kk\oplus\kk^2)\otimes\E^2\kk^2$ &
%TeXCAD Picture [GI-S.pic]. Options:
%\grade{\off}
%\emlines{\off}
%\epic{\off}
%\beziermacro{\off}
%\reduce{\on}
%\snapping{\off}
%\quality{2.00}
%\graddiff{0.01}
%\snapasp{1}
%\zoom{10.0000}
\unitlength .7ex % = .7pt
\linethickness{0.4pt}
\begin{picture}(5.67,1.67)(0,0)
\put(0,1){\circle{1.33}} \put(5,1){\circle*{1.33}}
\put(.33,1){\makebox(0,0)[lc]{$<$}} \put(.67,1){\line(1,0){3.67}}
\put(5,1.5){\line(-1,0){3.33}} \put(5,.5){\line(-1,0){3.33}}
\end{picture}
& $\omega_2-\omega_1$ & \\
\hline

\end{tabular}}
\end{center}
\end{table}

Simple minimal roots of arbitrary wonderful $G$-varieties of
rank~$1$ are called \emph{spherical roots} of~$G$. They are
non-negative linear combinations of~$\Pi_G$. Let $\Sigma_G$ denote
the set of all spherical roots. It is a finite set, which is easy
to find from the classification of wonderful varieties of
rank~$1$.

Spherical roots of reductive groups of simply connected type are
listed in Table~\ref{sph.roots}. Namely, $\lambda\in\Sigma_G$ iff
it is a spherical root of a simple factor, or a product of two
simple factors, indicated in the \ordinal{1} column of the table.
For each spherical root~$\lambda$, we indicate the Dynkin diagram
of the simple roots occurring in the decomposition of $\lambda$
with positive coefficients. The numbering of the simple roots
$\alpha_i$ is according to~\cite{Sem}, and $\alpha_i,\alpha'_j$
denote simple roots of different simple factors.
\begin{table}[!h]
\caption{Spherical roots} \label{sph.roots}
\begin{center}
\normalsize\makebox[0pt]{\begin{tabular}{|c|c|} \hline
$G$ & $\Sigma_G$ \\
\hline

$\AAa_l$ & $\alpha_i+\dots+\alpha_j$ ($\AAa_{j-i+1}$, $i\leq
j$),\quad $2\alpha_i$ ($\AAa_1$),\\ & $\alpha_i+\alpha_j$
($\AAa_1\times\AAa_1$, $i\leq j-2$),\quad $(\alpha_1+\alpha_3)/2$
($\AAa_1\times\AAa_1$, $l=3$),\\ &
$\alpha_{i-1}+2\alpha_i+\alpha_{i+1}$ ($\DDd_3$, $1<i<l$),\quad
$(\alpha_1+2\alpha_2+\alpha_3)/2$
($\DDd_3$, $l=3$) \\
\hline

& $\alpha_i+\dots+\alpha_j$ ($\AAa_{j-i+1}$, $i\leq j<l$),\quad
$2\alpha_i$ ($\AAa_1$),\\ & $\alpha_i+\alpha_j$
($\AAa_1\times\AAa_1$, $i\leq j-2$),\quad $(\alpha_1+\alpha_3)/2$
($\AAa_1\times\AAa_1$, $l=3,4$),\\
$\BBb_l$ & $\alpha_{i-1}+2\alpha_i+\alpha_{i+1}$ ($\DDd_3$,
$1<i<l-1$),\quad $(\alpha_1+2\alpha_2+\alpha_3)/2$ ($\DDd_3$,
$l=4$),\\ & $\alpha_i+\dots+\alpha_l$ ($\BBb_{l-i+1}$,
$i<l$),\quad $2\alpha_i+\dots+2\alpha_l$ ($\BBb_{l-i+1}$,
$i<l$),\\ & $\alpha_{l-2}+2\alpha_{l-1}+3\alpha_l$
($\BBb_3$),\quad $(\alpha_1+2\alpha_2+3\alpha_3)/2$
($\BBb_3$, $l=3$) \\
\hline

& $\alpha_i+\dots+\alpha_j$ ($\AAa_{j-i+1}$, $i\leq
j<l$),\quad $2\alpha_i$ ($\AAa_1$),\\
$\CCc_l$ & $\alpha_i+\alpha_j$ ($\AAa_1\times\AAa_1$, $i\leq
j-2$),\quad $\alpha_{i-1}+2\alpha_i+\alpha_{i+1}$ ($\DDd_3$,
$1<i<l-1$),\\ &
$\alpha_i+2\alpha_{i+1}+\dots+2\alpha_{l-1}+\alpha_l$
($\CCc_{l-i+1}$, $i<l$),\quad
$2\alpha_{l-1}+2\alpha_l$ ($\CCc_2$) \\
\hline

& $\alpha_{i_1}+\dots+\alpha_{i_k}$ ($\AAa_k$, $k\geq1$),\quad
$2\alpha_i$ ($\AAa_1$),\quad $\alpha_i+\alpha_j$
($\AAa_1\times\AAa_1$),\\ &
$2\alpha_{i_1}+\alpha_{i_2}+\alpha_{i_3}$ ($\DDd_3$),\quad
$(2\alpha_{i_1}+\alpha_{i_2}+\alpha_{i_3})/2$
($\DDd_3$, $l=4$),\\
$\DDd_l$ & $2\alpha_i+\dots+2\alpha_{l-2}+\alpha_{l-1}+\alpha_l$
($\DDd_{l-i+1}$, $i<l-1$),\\ &
$\alpha_i+\dots+\alpha_{l-2}+(\alpha_{l-1}+\alpha_l)/2$
($\DDd_{l-i+1}$, $i<l-1$),\quad $(\alpha_{l-1}+\alpha_l)/2$
($\AAa_1\times\AAa_1$),\\ & $(\alpha_1+\alpha_3)/2$
($\AAa_1\times\AAa_1$, $l=4$),\quad $(\alpha_1+\alpha_4)/2$
($\AAa_1\times\AAa_1$, $l=4$) \\
\hline

$\EEe_l$ & $\alpha_{i_1}+\dots+\alpha_{i_k}$ ($\AAa_k$,
$k\geq1$),\quad $2\alpha_i$ ($\AAa_1$),\quad $\alpha_i+\alpha_j$
($\AAa_1\times\AAa_1$),\\ &
$2\alpha_{i_1}+\dots+2\alpha_{i_{k-2}}+\alpha_{i_{k-1}}+\alpha_{i_k}$
($\DDd_k$, $k\geq3$) \\
\hline

& $\alpha_i$ ($\AAa_1$),\quad $2\alpha_i$ ($\AAa_1$),\quad
$\alpha_i+\alpha_j$ ($\AAa_1\times\AAa_1$),\quad
$\alpha_i+\alpha_{i+1}$
($\AAa_2$, $i\neq2$),\\
$\FFf_4$ & $\alpha_2+\alpha_3$ ($\CCc_2$),\quad
$2\alpha_2+2\alpha_3$ ($\CCc_2$),\quad
$\alpha_1+2\alpha_2+\alpha_3$ ($\CCc_3$),\\ &
$\alpha_2+\alpha_3+\alpha_4$ ($\BBb_3$),\quad
$2\alpha_2+2\alpha_3+2\alpha_4$ ($\BBb_3$),\quad
$3\alpha_2+2\alpha_3+\alpha_4$ ($\BBb_3$),\\ &
$2\alpha_1+3\alpha_2+2\alpha_3+\alpha_4$ ($\FFf_4$) \\
\hline

$\GGg_2$ & $\alpha_i$ ($\AAa_1$),\quad $2\alpha_i$
($\AAa_1$),\quad $\alpha_1+\alpha_2$ ($\GGg_2$),\quad
$2\alpha_1+\alpha_2$ ($\GGg_2$),\quad
$4\alpha_1+2\alpha_2$ ($\GGg_2$) \\
\hline

$\XXx_l\times\YYy_m$ & $\alpha_i+\alpha'_j$
($\AAa_1\times\AAa_1$),\quad $(\alpha_l+\alpha'_m)/2$
($\AAa_1\times\AAa_1$,
$\XXx=\YYy=\CCc$, $l,m\geq1$) \\
\hline

\end{tabular}}
\end{center}
\end{table}
For arbitrary~$G$, $\Sigma_G$~is obtained from
$\Sigma_{\widetilde{G}}$ by removing the spherical roots that are
not in the weight lattice of~$G$. Note that if
$\lambda,\mu\in\Sigma_G$ are proportional, then $\lambda=2\mu$ or
$\mu=2\lambda$, and also if
$\lambda\in\Sigma_G\setminus\ZZ\Delta_G$, then
$2\lambda\in\Sigma_G\cap\ZZ\Delta_G$.

More generally, two-orbit complete (normal) $G$-varieties were
classified by Cupit-Foutou~\cite{2-orb} and
Smirnov~\cite{q-closed}. All of them are spherical.

Wonderful varieties of rank~$2$ were classified by
Wasserman~\cite{wonder-2}.

For arbitrary wonderful varieties, many questions can be reduced
to the case of rank $\leq2$ via the procedure of localization
\cite{big.cells}, \cite[3.2]{sph(A)}.

Given a wonderful variety $X$ with the open $G$-orbit~$\HS$, there
is a bijection $D_i\leftrightarrow\lambda_i$ ($i=1,\dots,r$)
between the component set of $\partial{X}$ and $\Pi_X^{\min}$.
Namely $\lambda_i$ is orthogonal to the facet of $\Vv$
complementary to the ray which corresponds to~$D_i$. Also,
$\lambda_i$~is the $T$-weight of $T_zX/T_zD_i$ at the unique
$B^{-}$-fixed point~$z$.

For any subset $\Sigma\subset\Pi_X^{\min}$, put
$X^{\Sigma}=\bigcap_{\lambda_i\notin\Sigma}D_i$, the
\emph{localization} of $X$ at~$\Sigma$. It is a wonderful variety
with $\Pi_{X^{\Sigma}}^{\min}=\Sigma$, and all colors in
$\Dd(X^{\Sigma})^B$ are obtained as irreducible components of
$\overline{D}\cap X^{\Sigma}$, $D\in\Dd^B$. (To see the latter,
observe that every $B$-divisor on $X^{\Sigma}$ is contained in the
zeroes of a $B$-eigenform in projective coordinates, which extends
to $X$ by complete reducibility of $G$-modules.) In particular,
the above considered wonderful subvarieties $X_i,X_{ij}$ of ranks
$1,2$ are the localizations of $X$ at $\{\lambda_i\}$,
$\{\lambda_i,\lambda_j\}$, respectively.

Another kind of localization is defined by choosing a subset
$I\subset\Pi$. Let $P_I$ be the respective standard parabolic
in~$G$, with the standard Levi subgroup~$L_I$, and $T_I=Z(L_I)^0$.
Denote by $Z^I,\X^I,X^I$ the sets of $T_I$-fixed points in $Z$,
$\X$, and $X_I:=P_I\X$, respectively.
\begin{lemma}\label{loc(wonder)}
\begin{enumerate}
\item The contraction by a general dominant one-parameter subgroup
$\gamma\in\CoCh(T_I)$ gives a $P_I$-equivariant retraction
$\pi_I:X_I\to X^I$, $\pi_I(x)=\lim_{t\to0}\gamma(t)x$.
\item\label{X^I} $X^I$~is a wonderful $L_I$-variety with
$P(X^I)=P\cap L_I$, $\Pi_{X^I}^{\min}=\Pi_X^{\min}\cap\langle
I\rangle$, and the colors of $X^I$ are in bijection, given by the
pull-back along~$\pi_I$, with the $P_I$-unstable colors of $X$.
\item $\X^I\iso(\Ru{P}\cap L_I)\times Z^I$ is the $(B\cap
L_I)$-chart of $X^I$ intersecting all orbits.
\item\label{triv.act} $\Radu{P_I^{-}}$~fixes $X^I$ pointwise.
\end{enumerate}
\end{lemma}
\begin{proof}
It is obvious that $\pi_I$ contracts $\X\iso\Ru{P}\times Z$ onto
$\X^I\iso(\Ru{P}\cap L_I)\times Z^I$, while the conjugation by
$\gamma(t)$ contracts $P_I$ to~$L_I$. Hence $\pi_I$ extends to a
retraction of $X_I$ onto $X^I=L_I\X^I$, and $\pi_I^{-1}(\X^I)=\X$
since $X_I\setminus\X$ is closed and $\gamma$-stable. Thus the
$P_I$-unstable $B$-divisors on $X$ intersect $X_I$ and are the
pull-backs of the $(B\cap L_I)$-divisors on~$X^I$.  The structure
of $\X^I$ readily implies the remaining assertions on $X^I$
in~\ref{X^I}: both $G$- and $B$-orbits intersect $\X^I$ in the
orbits of $(\Ru{P}\cap L_I)T$, and $\Pi_X^{\min}\cap\langle
I\rangle$ is the set of $T$-weights of~$Z^I$.

Since $\Radu{P_I^{-}}$-orbits are connected, it suffices to prove
in~\ref{triv.act} that $\Radu{P_I^{-}}x\cap\X=\{x\}$, $\forall
x\in\X^I$. If $gx\in\X$ for some $g\in\Radu{P_I^{-}}$, then
$\gamma(t)gx=\gamma(t)g\gamma(t)^{-1}x\to x$ as $t\to\infty$,
whence $gx=x$, because $\gamma(t)$ contracts $\X$ to $\X^I$ as
$t\to0$.
\end{proof}
The wonderful variety $X^I$ is called the \emph{localization} of
$X$ at~$I$. It is easy to see that $X^I\subseteq X^{\Sigma}$,
where $\Sigma=\Pi_X^{\min}\cap\langle I\rangle$. Conversely, for
any $\Sigma\subseteq\Pi_X^{\min}$ one has $X^{\Sigma}\iso
G\itimes{Q}X^I$, where $I$ is the set of simple roots occurring in
the decompositions of $\lambda_i\in\Pi_X^{\min}$ and $Q$ is the
parabolic generated by $P_I^{-}$ and~$P^{-}$.

It is helpful to extend localization to an arbitrary spherical
homogeneous space~$\HS=G/H$ using an arbitrary smooth complete
toroidal embedding $X\embof\HS$ instead of the wonderful one. For
any subset $\Sigma\subseteq\Pi_{\HS}^{\min}$ one can find an orbit
$\HS'\subset X$ such that $\Cc_{\HS'}$ is a solid subcone in the
face $\Vv\cap\Sigma^{\perp}$ of $\Vv=\Vv(\HS)$. Then
$X^{\Sigma}:=\overline{\HS'}$ is a wonderful subvariety of $X$ and
$\Pi_{X^{\Sigma}}^{\min}=\Sigma$. In particular, the localization
at a single $\lambda\in\Pi_{\HS}^{\min}$ yields
$\Pi_{\HS}^{\min}\subseteq\Sigma_G$.

Also for any $I\subset\Pi$ such that $\RG(\HS)\not\subset\langle
I\rangle$ one can find $X\embof\HS$ and a general dominant
one-parameter subgroup $\gamma\in\CoCh(T_I)$ which is contained in
a unique solid cone $\Cc_Y$ in the fan of~$X$. (It suffices to
take care lest the image of $\langle I\rangle^{\ann}$ in
$\ES(\HS)$ should lie in a hyperplane which separates two
neighboring solid cones in the fan.) Starting with $\X=\X_Y$, one
defines $X^I$ as above and generalizes Lemma~\ref{loc(wonder)}. If
$\RG(\HS)\subset\langle I\rangle$, then $X^I$ can be defined for
any toroidal embedding $X\embof\HS$ using
$\X=X\setminus\bigcup_{D\in\Dd^B}D$. Lemma~\ref{loc(wonder)}
extends to this setup except that $X^I$ may be no longer
wonderful.

In particular, the localization at a single root $\alpha\in\Pi$
yields a smooth complete subvariety $X^{\alpha}$ of rank~$1$ acted
on by $S_{\alpha}=L_{\alpha}'\iso\SL_2(\kk)\text{ or
}\PSL_2(\kk)$. The classification of complete varieties of
rank~$1$, together with Lemma~\ref{loc(wonder)}, allows to
subdivide all simple roots in 4~types:
\begin{enumerate}
\item[\reftag{$p$}] $\alpha\in\Pi_L$. Here $X^{\alpha}$ is a point
and $P_{\alpha}$ leaves all colors stable. \item[\reftag{$b$}]
$\alpha\notin\QQ_{+}\Pi_{\HS}^{\min}\cup\Pi_L$. If $X^{\alpha}$ is
wonderful, then $r(X^{\alpha})=0$ whence
$X^{\alpha}=S_{\alpha}/B\cap S_{\alpha}\iso\PP^1$; otherwise
$r(X^{\alpha})=1$ and $X^{\alpha}\iso S_{\alpha}\itimes{B\cap
S_{\alpha}}\PP^1$, where $B\cap S_{\alpha}$ acts on $\PP^1$ via a
character. There is a unique $P_{\alpha}$-unstable color
$D_{\alpha}=\overline{\pi_{\alpha}^{-1}([\1])}\text{ or
}\overline{\pi_{\alpha}^{-1}(\1\itimes{}\PP^1)}$.
\item[\reftag{$a$}] $\alpha\in\Pi_{\HS}^{\min}$. Here
$r(X^{\alpha})=1$ and $X^{\alpha}\iso\PP^1\times\PP^1$. There are
two $P_{\alpha}$-unstable colors
$D_{\alpha}^{+}=\overline{\pi_{\alpha}^{-1}(\PP^1\times[\1])}$ and
$D_{\alpha}^{-}=\overline{\pi_{\alpha}^{-1}([\1]\times\PP^1)}$.
\item[\reftag{$a'$}] $2\alpha\in\Pi_{\HS}^{\min}$. Here
$r(X^{\alpha})=1$ and $X^{\alpha}\iso\PP^2=\PP(\s_{\alpha})$.
There is a unique $P_{\alpha}$-unstable color
$D_{\alpha}'=\overline{\pi_{\alpha}^{-1}(\PP(\br\cap\s_{\alpha}))}$.
\end{enumerate}

The type of a color $D\in\Dd^B$ is defined as the type of
$\alpha\in\Pi$ such that $P_{\alpha}$ moves~$D$. Using the
localization at $\{\alpha,\beta\}\subseteq\Pi$ and the
classification of wonderful varieties of rank $\leq2$, one
verifies that, as a rule, each $D\in\Dd^B$ is moved by a
unique~$P_{\alpha}$, with the following exceptions:
$D_{\alpha}=D_{\beta}$ iff $\alpha,\beta$ are pairwise orthogonal
simple roots of type~$b$ such that
$\alpha+\beta\in\Pi_{\HS}^{\min}\sqcup2\Pi_{\HS}^{\min}$; two sets
$\{D_{\alpha}^{\pm}\}$ and $\{D_{\beta}^{\pm}\}$ may intersect in
one color for distinct $\alpha,\beta$ of type~$a$. In particular,
each color belongs to exactly one type. We obtain disjoint
partitions $\Pi=\Pi^{a}\sqcup\Pi^{a'}\sqcup\Pi^{b}\sqcup\Pi^{p}$,
$\Dd^B=\Dd^{a}\sqcup\Dd^{a'}\sqcup\Dd^{b}$ according to the types
of simple roots and colors.

\begin{lemma}\label{colors}
For $\forall\lambda\in\RG(\HS)$ we have
%*
\begin{align*}
\langle D_{\alpha}^{+},\lambda\rangle+ \langle
D_{\alpha}^{-},\lambda\rangle&=
\langle\alpha^{\vee},\lambda\rangle,\qquad
\forall\alpha\in\Pi^{a}\\
\langle D_{\alpha}',\lambda\rangle&=
\langle\tfrac{\alpha^{\vee}}2,\lambda\rangle,\qquad
\forall\alpha\in\Pi^{a'}\\
\langle D_{\alpha},\lambda\rangle&=
\langle\alpha^{\vee},\lambda\rangle,\qquad
\forall\alpha\in\Pi^{b}\\
\end{align*}
%*
\end{lemma}
\begin{proof}
Let $Y^{\alpha}\iso S_{\alpha}/(B^{-}\cap S_{\alpha})\iso\PP^1$ be
a closed $S_{\alpha}$-orbit in~$X^{\alpha}$. Namely $Y^{\alpha}$
is the diagonal of $X^{\alpha}\iso\PP^1\times\PP^1$ in type~$a$, a
conic in $X^{\alpha}\iso\PP^2$ in type~$a'$, and a section of
$X^{\alpha}\to\PP^1$ in type~$b$. Put
$\delta_{\lambda}=\sum_{D\in\Dd^B}\langle D,\lambda\rangle D$.
From the description of $P_{\alpha}$-stable and unstable colors,
we readily derive $\langle
Y^{\alpha},\delta_{\lambda}\rangle=\langle
D_{\alpha}^{+},\lambda\rangle+\langle
D_{\alpha}^{-},\lambda\rangle$, $2\langle
D_{\alpha}',\lambda\rangle$, or $\langle
D_{\alpha},\lambda\rangle$, depending on the type of~$\alpha$.

On the other hand, $\delta_{\lambda}\sim-\sum\langle
v_i,\lambda\rangle D_i$, where $D_i$ runs over all $G$-stable
prime divisors in~$X$ and $v_i\in\Vv$ is the corresponding
$G$-valuation. Since $\Lin{D_i}|_{D_i}$ is the normal bundle
of~$D_i$, the fiber of $\Lin{D_i}$ at the $B^{-}$-fixed point
$z\in Y^{\alpha}$ is $T_zX/T_zD_i$ for $\forall D_i\supseteq
Y^{\alpha}$. Note that the $T$-weights $\lambda_i$ of these fibers
form the basis of $-\Cc_Y^{\vee}$ dual to $-v_i$'s, where $Y=Gz$
is the closed $G$-orbit in $X$ containing~$Y^{\alpha}$. Hence
$\Lin{\delta_{\lambda}}|_{Y^{\alpha}}=\ind{-\sum_{D_i\supseteq
Y}\langle v_i,\lambda\rangle\lambda_i}=\ind{\lambda}$ and $\langle
Y^{\alpha},\delta_{\lambda}\rangle=\deg\ind{\lambda}=
\langle\alpha^{\vee},\lambda\rangle$. The lemma follows.
\end{proof}

These results show that $\Dd^{a'},\Dd^{b}$ as abstract sets and
their representation in $\ES(\HS)$ are determined by $\Pi^{p}$
and~$\Pi_{\HS}^{\min}$. The colors of type~$a$, together with the
weight lattice, the parabolic $P$, and the simple minimal roots,
form a collection of combinatorial invariants supposed to identify
$\HS$ up to isomorphism. Namely
$(\RG(\HS),\Pi^{p},\Pi_{\HS}^{\min},\Dd^{a})$ is a homogeneous
spherical datum in the sense of the following
\begin{definition}[{\cite[\S2]{sph(A)}}]
\label{sph.data} A \emph{homogeneous spherical datum} is a
collection $(\RG,\Pi^{p},\Sigma,\Dd^{a})$, where $\RG$ is a
sublattice in~$\Ch(T)$, $\Pi^{p}\subseteq\Pi_G$,
$\Sigma\subseteq\Sigma_G\cap\RG$ is a linearly independent set
consisting of indivisible vectors in~$\RG$, and $\Dd^{a}$ is a
finite set equipped with a map $\res:\Dd^{a}\to\RG^{*}$, which
satisfies the following axioms:
\begin{enumerate}
\item[\reftag{A1}] $\langle\res(D),\lambda\rangle\leq1$, $\forall
D\in\Dd^{a},\ \lambda\in\Sigma$, and the equality is reached iff
$\lambda=\alpha\in\Sigma\cap\Pi$ and $D=D_{\alpha}^{\pm}$, where
$D_{\alpha}^{+},D_{\alpha}^{-}\in\Dd^{a}$ are two distinct
elements depending on~$\alpha$. \item[\reftag{A2}]
$\res(D_{\alpha}^{+})+\res(D_{\alpha}^{-})=\alpha^{\vee}$ on~$\RG$
for $\forall\alpha\in\Sigma\cap\Pi$. \item[\reftag{A3}]
$\Dd^{a}=\{D_{\alpha}^{\pm}\mid\alpha\in\Sigma\cap\Pi\}$
\item[\reftag{$\Sigma$1}] If $\alpha\in\Pi\cap\frac12\Sigma$, then
$\langle\alpha^{\vee},\RG\rangle\subseteq2\ZZ$ and
$\langle\alpha^{\vee},\Sigma\setminus\{2\alpha\}\rangle\leq0$.
\item[\reftag{$\Sigma$2}] If $\alpha,\beta\in\Pi$,
$\alpha\perp\beta$, and $\alpha+\beta\in\Sigma\sqcup2\Sigma$, then
$\alpha^{\vee}=\beta^{\vee}$ on~$\RG$. \item[\reftag{S}]
$\langle\alpha^{\vee},\RG\rangle=0$, $\forall\alpha\in\Pi^{p}$,
and the pair $(\lambda,\Pi^p)$ comes from a wonderful variety of
rank~$1$ for any $\lambda\in\Sigma$.
\end{enumerate}
A \emph{spherical system} is a triple $(\Pi^{p},\Sigma,\Dd^{a})$
satisfying the above axioms with $\RG=\ZZ\Sigma$.
\end{definition}

The homogeneous spherical datum of the open orbit in a wonderful
variety amounts to its spherical system. It is easy to see that
there are finitely many spherical systems for given~$G$.

For the homogeneous spherical datum of~$\HS$, most of the axioms
\reftag{A1}--\reftag{A3}, \reftag{$\Sigma$1}--\reftag{$\Sigma$2},
\reftag{S} are verified using the above results together with some
additional general arguments. For instance, the inequality in
\reftag{$\Sigma$1} stems from the fact that
$\Sigma=\Pi_{\HS}^{\min}$ is a base of a root
system~$\Delta_{\HS}^{\min}$. On the other hand, each axiom
involves at most two simple or spherical roots, like the axioms of
classical root systems. Thus the localizations at one or two
simple or spherical roots reduce the verification to wonderful
varieties of rank $\leq2$.

Actually the list of axioms was obtain by inspecting the
classification of wonderful varieties of rank $\leq2$, which leads
to the following conclusion: spherical systems (homogeneous data)
with $|\Sigma|\leq2$ bijectively correspond to wonderful varieties
of rank $\leq2$ (resp.\ to spherical homogeneous spaces $\HS=G/H$
with $r(G/N(H))\leq2$). It is tempting to extend this
combinatorial classification to arbitrary wonderful varieties and
spherical spaces. Luna succeeded to fulfil this program in the
case, where all simple factors of $G$ are locally isomorphic
to~$\SL_{n_i}$.
\begin{theorem}[\cite{sph(A)}]\label{sph(A)}
Suppose $G$ is a reductive group with all simple factors of
type~$\AAa$; then there are natural bijections:
%*
\begin{align}
\left\{\text{\begin{tabular}{c} spherical homogeneous \\
$G$-spaces
\end{tabular}}\right\}&\longleftrightarrow
\left\{\text{\begin{tabular}{c} homogeneous spherical \\ data
for~$G$
\end{tabular}}\right\} \label{sph.h.s.(A)}\\
\{\text{wonderful $G$-varieties}\}&\longleftrightarrow
\{\text{spherical systems for~$G$}\} \label{wonder(A)}
\end{align}
%*
\end{theorem}
Recently this result was generalized to the groups with the simple
factors of types $\AAa$ and~$\DDd$ \cite{sph(D)} or $\AAa$
and~$\CCc$ (with some technical restrictions) \cite{sph(C)}.

Actually \eqref{sph.h.s.(A)} was proved by Luna for any $G$
provided that $G/Z(G)$ satisfies~\eqref{wonder(A)}. The basic idea
is to replace $\HS=G/H$ by $\overline\HS=G/\overline{H}$. This
passage preserves the types of simple roots and colors, and
$\Pi_{\overline\HS}^{\min}$ is obtained from $\Pi_{\HS}^{\min}$ by
a dilation: some $\lambda\in\Pi_{\HS}^{\min}\setminus\QQ_{+}\Pi$
are replaced by~$2\lambda$. It is not hard to prove that spherical
subgroups $H$ with a fixed very sober hull $\overline{H}$
bijectively correspond to homogeneous spherical data
$(\RG,\Pi^{p},\Sigma,\Dd^{a})$ such that
$(\Pi^{p},\Pi_{\overline\HS}^{\min},\Dd^{a})$ is the spherical
system of $\overline\HS$, $\RG\supset\Pi_{\overline\HS}^{\min}$,
and $\Sigma$ is obtained from $\Pi_{\overline\HS}^{\min}$ by
replacing $\lambda\in\Pi_{\overline\HS}^{\min}\setminus\QQ_{+}\Pi$
by $\lambda/2$ whenever $\lambda/2\in\RG$ \cite[\S6]{sph(A)}.

If \eqref{wonder(A)} holds for the adjoint group of~$G$, then
$\Pi_{\overline\HS}^{\min}$ is obtained from $\Pi_{\HS}^{\min}$ by
the ``maximal possible'' dilation: every
$\lambda\in\Pi_{\HS}^{\min}\setminus\QQ_{+}\Pi$ such that
$2\lambda\in\Sigma_G$ is replaced by~$2\lambda$
\cite[7.1]{sph(A)}. It follows that the spherical homogeneous
datum of $\HS$ determines the spherical system of $\overline\HS$
in a pure combinatorial way. Conversely, this spherical system
together with $\RG$ determines $\overline\HS$ and $\HS$ by the
above, which proves~\eqref{sph.h.s.(A)}.

The proof of \eqref{wonder(A)} for adjoint~$G$ is much more
involved for lacking of a uniform conceptual argument. Let us
explain the general scheme. The first stage is to prove that
certain geometric operations on wonderful varieties (localization,
parabolic induction, direct product, etc) are expressed in a pure
combinatorial language of spherical systems. Every spherical
system is obtained by these combinatorial operations from a list
of \emph{primitive} systems, which are explicitly classified. For
primitive spherical systems, the existence and uniqueness of a
geometric realization is proved case by case.

Another situation, where the assertions of Theorem~\ref{sph(A)}
remain valid, is the classification of solvable spherical
subgroups (of arbitrary~$G$) \cite{sph.solv}. More precisely, the
bijections \eqref{sph.h.s.(A)}--\eqref{wonder(A)} hold for
spherical homogeneous spaces with stabilizers contained in a Borel
subgroup of~$G$. Spherical data arising here satisfy
$\Sigma=\Pi^{a}$, $\Pi^{a'}=\Pi^{p}=\emptyset$,
$D_{\alpha}^{-}\neq D_{\beta}^{\pm}$
($\forall\alpha,\beta\in\Pi^{a}$, $\alpha\neq\beta$), and
$-(\QQ_{+}\Sigma)^{\vee}+\sum\QQ_{+}\res(D_{\alpha}^{+})=
\RG\otimes\QQ$.
\begin{question}
$H\subseteq B^{-}\implies D_{\alpha},D_{\alpha}^{-}$ are lifted
from $G/B^{-}$; $\langle D_{\alpha}^{+},\beta\rangle<0\implies
\langle D_{\beta}^{+},\alpha\rangle=0$
\end{question}

\section{Frobenius splitting}
\label{Frob.split}

Frobenius splitting is a powerful tool of modern algebraic
geometry which allows to prove various geometric and cohomological
results by reduction to positive characteristic. This notion was
introduced by Mehta and Ramanathan \cite{Schubert.split} in their
study of Schubert varieties.

Let $X$ be an algebraic variety over an algebraically closed field
$\kk$ of characteristic $p>0$. The Frobenius endomorphism
$f\mapsto f^p$ of $\Oo_X$ gives rise to the Frobenius morphism
$F:X^{1/p}\to X$, where $X^{1/p}=X$ as ringed spaces but the
$\kk$-algebra structure on $\Oo_{X^{1/p}}$ is defined as
$c*f=c^pf$, $\forall c\in\kk$.

If $X$ is a subvariety in $\AAA^n$ or $\PP^n$, then $X^{1/p}$ is,
too. The defining equations of $X^{1/p}$ are obtained from those
of $X$ by replacing all coefficients by their \ordinal{$p$}
powers. The Frobenius morphism $F$ is given by raising all
coordinates to the power~$p$.

The Frobenius endomorphism may be regarded as an injection of
$\Oo_X$-modules $\Oo_X\embeds F_{*}\Oo_{X^{1/p}}$, where
$F_{*}\Oo_{X^{1/p}}=\Oo_X$ is endowed with another $\Oo_X$-module
structure: $f*h=f^ph$ for any local sections $f$ of $\Oo_X$ and
$h$ of~$F_{*}\Oo_{X^{1/p}}$.

\begin{definition}
The variety $X$ is said to be \emph{Frobenius split} if the
Frobenius homomorphism has an $\Oo_X$-linear left inverse
$\sigma:F_{*}\Oo_{X^{1/p}}\to\Oo_X$, called a \emph{Frobenius
splitting}. In other words, $\sigma$~is a $\ZZz_p$-linear
endomorphism of $\Oo_X$ such that $\sigma(1)=1$ and
$\sigma(f^ph)=f\sigma(h)$.

For any subvariety $Y\subset X$ one has
$\sigma(\Ii_Y)\supseteq\Ii_Y$, because $\Ii_Y\supseteq\Ii_Y^p$.
The splitting $\sigma$ is \emph{compatible} with $Y$ if
$\sigma(\Ii_Y)=\Ii_Y$. Clearly, a compatible splitting induces a
splitting of~$Y$.

More generally, let $\delta$ be an effective Cartier divisor
on~$X$, with the canonical section
$\eta_{\delta}\in\Ho^0(X,\Lin{\delta})$,
$\divr\eta_{\delta}=\delta$. We say that $X$ is \emph{Frobenius
split relative to $\delta$} if there exists an $\Oo_X$-module
homomorphism, called a \emph{$\delta$-splitting},
$\sigma_{\delta}:F_{*}\Lin[X^{1/p}]{\delta}\to\Oo_X$ such that
$\sigma(h)=\sigma_{\delta}(h\eta_{\delta})$ is a Frobenius
splitting, or equivalently, $\sigma_{\delta}(\eta_{\delta})=1$ and
$\sigma_{\delta}(f^p\eta)=f\sigma_{\delta}(\eta)$ for any local
section $\eta$ of~$\Lin{\delta}$. The $\delta$-splitting
$\sigma_{\delta}$ is \emph{compatible} with $Y$ if $\Supp\delta$
contains no component of~$Y$ (i.e., $\delta$~restricts to a
divisor on~$Y$) and $\sigma$ is compatible with~$Y$. Then
$\sigma_{\delta}$ induces a $(\delta\cap Y)$-splitting of $Y$.
\end{definition}

For a systematic treatment of Frobenius splitting and its
applications, we refer to a monograph of Brion and
Kumar~\cite{F-split}. Here we recall some of its most important
properties.

Clearly, a Frobenius splitting of $X$ (compatible with~$Y$,
relative to~$\delta$) restricts to a splitting of every open
subvariety $U\subset X$ (compatible with~${Y\cap U}$, relative
to~$\delta\cap U$). Conversely, if $X$ is normal and
$\codim(X\setminus U)>1$, then any splitting of $U$ extends
to~$X$. In applications it is often helpful to consider
$U=X^{\reg}$.

If $\phi:X\to Y$ is a morphism such that $\phi_{*}\Oo_X=\Oo_Y$,
then a Frobenius splitting of $X$ descends to a splitting of~$Y$.
If the splitting of $X$ is compatible with $Z\subset X$, then the
splitting of $Y$ is compatible with~$\overline{\phi(Z)}$. For
instance, one obtains a splitting of a normal variety $X$ from
that of its desingularization.

It is not hard to prove that Frobenius split varieties are
\emph{weakly normal}, i.e., every bijective finite birational map
onto a Frobenius split variety has to be isomorphic
\cite[1.2.5]{F-split}.

\begin{proposition}\label{H(split)}
\begin{enumerate}
\item\label{vanish} Suppose $X$ is a Frobenius split projective
variety; then $\Ho^i(X,\Ll)=0$ for any ample line bundle $\Ll$ on
$X$ and $\forall i>0$. \item\label{surject} If $Y\subset X$ is a
compatibly split subvariety, then the restriction map
$\Ho^0(X,\Ll)\to\Ho^0(Y,\Ll)$ is surjective.
\item\label{rel.ample} If the splittings above are relative to an
ample divisor~$\delta$, then the assertions of
\ref{vanish}--\ref{surject} hold for any numerically effective
line bundle. \item\label{relative} There are relative versions of
assertions \ref{vanish}--\ref{rel.ample} for a proper morphism
$\phi:X\to Z$ stating that $\Rf^i\phi_{*}\Ll=0$ and
$\phi_{*}\Ll\to\phi_{*}(\iota_{*}\iota^{*}\Ll)$ is surjective
under the same assumptions, with $\iota:Y\embeds X$.
\end{enumerate}
\end{proposition}
\begin{proof}
The idea is to embed the cohomology of $\Ll$ as a direct summand
in the cohomology of a sufficiently big power of~$\Ll$. Namely the
canonical homomorphism $\Ll\to
F_{*}F^{*}\Ll=\Ll\otimes_{\Oo_X}F_{*}\Oo_{X^{1/p}}$ has a left
inverse $\id\otimes\sigma$, whence $\Ll$ is a direct summand
in~$F_{*}F^{*}\Ll$. Taking the cohomology yields a split injection
%*
\begin{equation*}
\Ho^i(X,\Ll)\embeds\Ho^i(X,F_{*}F^{*}\Ll)\iso
\Ho^i(X^{1/p},F^{*}\Ll)\iso\Ho^i(X,\Ll^{\otimes p}), \qquad\forall
i\geq0
\end{equation*}
%*
(The right isomorphism is only $\ZZz_p$-linear.) Iterating this
procedure yields a split $\ZZz_p$-linear injection
$\Ho^i(X,\Ll)\embeds\Ho^i(X,\Ll^{\otimes p^k})$ compatible with
the restriction to~$Y$. Thus the assertions \ref{vanish} and
\ref{surject} are reduced to the case of the line bundle
$\Ll^{\otimes p^k}$, $k\gg0$, where the Serre theorem applies
\cite[III.5.3]{AG}.

Similar reasoning applies to \ref{rel.ample} making use of a split
injection $\Ho^i(X,\Ll)\embeds\Ho^i(X,{\Ll^{\otimes p}
\otimes\Lin{\delta}})$ together with ampleness of~$\Ll^{\otimes
p}\otimes\Lin{\delta}$. The relative assertions are proved by the
same arguments.
\end{proof}

Among other cohomology vanishing results for Frobenius split
varieties we mention the extension of the Kodaira vanishing
theorem \cite[1.2.10(i)]{F-split}: if $X$ is smooth projective and
Frobenius split, then $\Ho^i(X,\Ll\otimes\omega_X)=0$ for ample
$\Ll$ and $i>0$.

Now we reformulate the notion of Frobenius splitting for smooth
varieties in terms of differential forms.

The de Rham derivation of $\Omega^{\bullet}_X$ may be considered
as an $\Oo_X$-linear derivation
of~$F_{*}\Omega^{\bullet}_{X^{1/p}}$. Let $\Hh^k_X$ denote the
respective cohomology sheaves. It is easy to check that
$f\mapsto[f^{p-1}df]$ is a $\kk$-derivation of $\Oo_X$ taking
values in~$\Hh^1_X$ (where $[\,\cdot\,]$ denotes the de Rham
cohomology class). By the universal property of K\"ahler
differentials, it induces a homomorphism of graded
$\Oo_X$-algebras
%*
\begin{equation*}
c:\Omega^{\bullet}_X\to\Hh^{\bullet}_X,\qquad
c(f_0df_1\wedge\dots\wedge df_k)=[f_0^p(f_1\cdots f_k)^{p-1}
df_1\wedge\dots\wedge df_k],
\end{equation*}
%*
called the \emph{Cartier operator}. Cartier proved that $c$ is an
isomorphism for smooth~$X$. (Using local coordinates, the proof is
reduced to the case $X=\AAA^n$, where the verification is
straightforward \cite[1.3.4]{F-split}.)

Now suppose that $X$ is smooth. Then we have the \emph{trace map}
%*
\begin{equation*}
\tau:{F_{*}\omega_{X^{1/p}}\to\omega_X},\qquad
\tau(\omega)=c^{-1}[\omega]
\end{equation*}
%*
In local coordinates $x_1,\dots,x_n$, the trace map can be
characterized as the unique $\Oo_X$-linear map taking $(x_1\cdots
x_n)^{p-1} dx_1\wedge\dots\wedge dx_n\mapsto dx_1\wedge\dots\wedge
dx_n$ and $x_1^{k_1}\cdots x_n^{k_n} dx_1\wedge\dots\wedge dx_n
\mapsto0$ unless $k_1\equiv\dots\equiv k_n\equiv p-1\pmod p$.

Using the trace map, it is easy to establish an isomorphism
%*
\begin{equation*}
\Homsh(F_{*}\Oo_{X^{1/p}},\Oo_X)\iso
F_{*}\omega_{X^{1/p}}^{1-p},\qquad
\sigma\leftrightarrow\widehat\sigma
\end{equation*}
%*
such that $\sigma(h)\omega=\tau(h\omega^{\otimes p}
\otimes\widehat\sigma)$ for any local sections $h$ of
$F_{*}\Oo_{X^{1/p}}$ and $\omega$ of~$\omega_X$. Similarly, for
any divisor $\delta$ on $X$ we have
%*
\begin{equation*}
\Homsh(F_{*}\Lin[X^{1/p}]{\delta},\Oo_X)\iso
F_{*}\omega_{X^{1/p}}^{1-p}(-\delta)
\end{equation*}
%*
This leads to the following conclusion.
\begin{proposition}[{\cite[1.3.8, 1.4.10]{F-split}}]
\label{diff.split} Suppose that $X$ is smooth and irreducible.
Then $\sigma\in\Hom(F_{*}\Oo_{X^{1/p}},\Oo_X)$ is a splitting of
$X$ iff the Taylor expansion of $\widehat\sigma$ at some (whence
any) $x\in X$ has the form
%*
\begin{equation*}
\left((x_1\cdots x_n)^{p-1}+\sum c_{k_1,\dots,k_n} x_1^{k_1}\cdots
x_n^{k_n}\right)
(\partial_1\wedge\dots\wedge\partial_n)^{\otimes(p-1)}
\end{equation*}
%*
where the sum is taken over all multiindices $(k_1,\dots,k_n)$
such that $\exists k_i\not\equiv p-1\pmod p$. (Here $x_i$ denote
local coordinates and $\partial_i$ the vector fields dual
to~$dx_i$.) If $X$ is complete, then it suffices to have
%*
\begin{equation*}
\widehat\sigma=((x_1\cdots x_n)^{p-1}+\cdots)
(\partial_1\wedge\dots\wedge\partial_n)^{\otimes(p-1)}
\end{equation*}
%*
The splitting $\sigma$ is relative to any effective divisor
$\delta\leq\divr\widehat\sigma$.
\end{proposition}

By abuse of language, we shall say that $\widehat\sigma$ splits
$X$ if $\sigma$ does. Also, $X$~is said to be \emph{split by a
\ordinal{$(p-1)$} power} if $\alpha^{\otimes(p-1)}$ splits $X$ for
some $\alpha\in\Ho^0(X,\omega_X^{-1})$. This splitting is
compatible with~$\VV(\alpha)$. For instance, a smooth complete
variety $X$ is split by the \ordinal{$(p-1)$} power of $\alpha$ if
the divisor of $\alpha$ in a neighborhood of some $x\in X$ is a
union of $n=\dim X$ smooth prime divisors intersecting
transversally at~$x$.
\begin{example}
Every smooth toric variety $X$ is Frobenius split by a
\ordinal{$({p-1})$} power compatibly with~$\partial{X}$. For
complete~$X$, this stems from the structure of its canonical
divisor, given by Proposition~\ref{can(sph)} (which extends to
positive characteristic in the toric case). The general case
follows by passing to a smooth toric completion. Now toric
resolution of singularities readily implies that all normal toric
varieties are Frobenius split compatibly with their invariant
subvarieties.
\end{example}
\begin{example}[{\cite{F-split(G/P)}, \cite[Ch.2--3]{F-split}}]
Generalized flag varieties are Frobenius split by a
\ordinal{$(p-1)$} power. For $X=G/B$, $\omega_X^{-1}=\ind{-2\rho}$
and the splitting is provided by $\alpha=f_{\rho}\cdot
f_{-\rho}\in V^{*}(2\rho)$, where $f_{\pm\rho}\in V^{*}(\rho)$ are
$T$-weight vectors of weights~$\pm\rho$.

Moreover, this splitting is compatible with all Schubert
subvarieties $S_w=\overline{B[w]}\subset X$, $w\in W$. Using the
weak normality of $S_w$ and the Bott--Samelson resolution of
singularities
%*
\begin{gather*}
\phi:\check{S}=\check{S}_{\alpha_1,\dots,\alpha_l}:=
P_{\alpha_1}\itimes{B}\cdots\itimes{B}P_{\alpha_l}/B\to S_w\\
w=r_{\alpha_1}\cdots r_{\alpha_l},\quad \alpha_i\in\Pi,\quad
l=\dim S_w
\end{gather*}
%*
with connected fibers and $\Rf^i\phi_{*}\Oo_{\check{S}}=0$,
$\forall i>0$, one deduces that $S_w$ are normal (Demazure,
Seshadri) and have rational singularities (Andersen, Ramanathan).
These properties descend to Schubert subvarieties in~$G/P$,
$\forall P\supset B$.
\end{example}

Splitting by a \ordinal{$(p-1)$} power has further important
consequences. For instance, the Grauert--Riemenschneider theorem
extends to this situation, due to Mehta--van der Kallen
\cite{MvdK}:
\begin{quote}
If $\phi:X\to Y$ is a proper birational morphism, $X$~is smooth
and split by~$\alpha^{\otimes(p-1)}$ such that $\phi$ is
isomorphic on $X\setminus\VV(\alpha)$, then
$\Rf^i\phi_{*}\omega_X=0$, $\forall i>0$.
\end{quote}

Although the concept of Frobenius splitting is defined in
characteristic $p>0$, it successfully applies to algebraic
varieties in characteristic zero via reduction~${\!}\bmod p$.

Namely let $X$ be an algebraic variety over an algebraically
closed field $\kk$ of characteristic~$0$. One can find a finitely
generated subring $R\subset\kk$ such that $X$ is defined over~$R$,
i.e., is obtained from an $R$-scheme $\Xx$ by extension of
scalars. One may assume that $\Xx$ is flat over~$R$, after
replacing $R$ by a localization.
\begin{question}
Is flatness local in $\Spec R$? First replace $X$ by a completion?
\end{question}
For any maximal ideal $\p\normin R$ we have $R/\p\iso\FF_{p^k}$.
The variety $X_{\p}$ obtained from the fiber $\Xx_{\p}$ of
$\Xx\to\Spec R$ over $\p$ by an extension of scalars
$\FF_{p^k}\to\FF_{p^{\infty}}$ is called a \emph{reduction
${\!}\bmod p$} of~$X$ and sometimes denoted simply by~$X_p$ (by
abuse of notation).

Reductions ${\!}\bmod p$ exist and share geometric properties of
$X$ (affinity, projectivity, completeness, smoothness, normality,
rationality of singularities, etc) for all sufficiently large~$p$.
Conversely, a local geometric property of open type (e.g.,
smoothness, normality, rationality of singularities) holds for $X$
if it holds for $X_p$ whenever $p\gg0$. Replacing $R$ by an
appropriate localization, one may always assume that a given
finite collection of algebraic and geometric objects on~$X$
(subvarieties, line bundles, coherent sheaves, morphisms, etc) is
defined over~$R$, whence specializes to $X_p$ for $p\gg0$;
coherent sheaves may be supposed to be flat over~$R$.

Cohomological applications of reduction ${\!}\bmod p$ are based on
the semicontinuity theorem \cite[III.12.8]{AG}, which may be
reformulated in our setup as follows:
\begin{quote}
If $X$ is complete and $\Ff$ is a coherent sheaf on~$X$, then
$\dim\Ho^i(X,\Ff)=\dim\Ho^i(X_p,\Ff_p)$ for all $p\gg0$.
\end{quote}
This implies, for instance, that the assertions of
Proposition~\ref{H(split)} hold in characteristic zero provided
that $X_p$ are Frobenius split for $p\gg0$. This is the case,
e.g., for Fano varieties. Another case, which is important in the
framework of this chapter, are spherical varieties.

\begin{theorem}[\cite{F-split.sph}]
If $X$ is a spherical $G$-variety in characteristic~$0$, then
$X_p$ is Frobenius split by a \ordinal{$(p-1)$} power compatibly
with all $G$-subvarieties and relative to any given $B$-stable
effective divisor, for $p\gg0$.
\end{theorem}
\begin{proof}
Using an equivariant completion of $X$ and its toroidal
desingularization, we may assume that $X$ is smooth, complete, and
toroidal. Consider the natural morphism $\phi:X\to X(\h)$, where
$\h$ is a generic isotropy subalgebra for $G:X$. By
Proposition~\ref{can(sph)},
$\omega_X^{-1}=\Lin{\partial{X}+\phi^{*}\Hh}$, where $\Hh$ is a
hyperplane section of~$X(\h)$.

The restriction of $\Lin{\partial{X}}$ to a closed $G$-orbit
$Y\subset X$ is the top exterior power of the normal bundle
to~$Y$, whence $\omega_Y^{-1}=
\omega_X^{-1}|_Y\otimes\Lin{-\partial{X}}|_Y=\Lin{\phi^{*}\Hh}|_Y$.
Since $Y$ is a generalized flag variety, $Y_p$~is split by the
\ordinal{$(p-1)$} power of (the reduction ${\!}\bmod p$ of) some
$\alpha_Y\in\Ho^0(Y,\omega_Y^{-1})$. The $G$-module
$\Ho^0(Y,\omega_Y^{-1})$ being irreducible and $\Lin{\phi^{*}\Hh}$
globally generated, the restriction map
$\Ho^0(X,\Lin{\phi^{*}\Hh})\to\Ho^0(Y,\omega_Y^{-1})$ is
surjective and $\alpha_Y$ extends to
$\alpha_0\in\Ho^0(X,\Lin{\phi^{*}\Hh})$.

We have $\partial{X}=D_1\cup\dots\cup D_k$, where $D_i$ runs over
all $G$-stable prime divisors of~$Y$. It is easy to see from
Proposition~\ref{diff.split} that
$\alpha=\alpha_0\otimes\alpha_1\otimes\dots\otimes\alpha_k$
provides a splitting for~$X_p$, where
$\alpha_i\in\Ho^0(X,\Lin{D_i})$, $\divr\alpha_i=D_i$. Moreover,
this splitting is compatible with all $(D_i)_p$ and therefore with
all $G$-subvarieties in~$X_p$, because the latter are unions of
transversal intersections of some~$(D_i)_p$.

Finally, for any $B$-stable effective divisor $\delta$ we have
$\delta\leq(1-p)K_X$ for $p\gg0$, by Proposition~\ref{can(sph)}.
Hence the splitting is relative to $\delta_p$ by
Proposition~\ref{diff.split}.
\end{proof}
It is worth noting that not all spherical varieties in positive
characteristic are Frobenius split. Counterexamples are provided
by some complete homogeneous spaces with non-reduced isotropy
group subschemes~\cite{split.comp}.

Frobenius splitting of spherical varieties provides short and
conceptual proofs for a number of important geometric and
cohomological properties. In particular, Theorem~\ref{norm&rat}
can be deduced in the following way.

Consider a resolution of singularities $\psi:X'\to X$, where $X'$
is toroidal and quasiprojective. Choose an ample $B$-stable
effective divisor $\delta$ on~$X'$; then $X'_p$ is split relative
to~$\delta_p$ for $p\gg0$. By semicontinuity and
Proposition~\ref{H(split)}\ref{relative}, applied to the trivial
line bundle over~$X'_p$, $\Rf^i\psi_{*}\Oo_{X'}=0$, whence $X$ has
rational singularities. By the same reason,
$\Oo_X=\psi_{*}\Oo_{X'}$~surjects onto $\psi_{*}\Oo_{Y'}$ for any
irreducible closed $G$-subvariety $Y'\subset X'$, whence
$\psi_{*}\Oo_{Y'}=\Oo_Y$ for $Y=\psi(Y')$. Since $Y'$ is smooth,
$Y$~is normal and has rational singularities by the above.

For any line bundle $\Ll$ on~$X$ denote $\Ll'=\psi^{*}\Ll$. The
Leray spectral sequence
%*
\begin{equation*}
\Ho^{i+j}(X',\Ll')\Longleftarrow\Ho^i(X,\Rf^j\psi_{*}\Ll')=
\Ho^i(X,\Ll\otimes\Rf^j\psi_{*}\Oo_{X'})
\end{equation*}
%*
degenerates to $\Ho^i(X',\Ll')=\Ho^i(X,\Ll)$, $\forall i\geq0$.
The same holds for direct images instead of cohomology. Together
with Proposition~\ref{H(split)}, applied to $X'_p$ and~$\Ll'_p$,
this proves the following
\begin{corollary}\label{H(nef)}
Suppose $\ch\kk=0$. If $X$ is a complete spherical $G$-variety,
$Y\subset X$ a $G$-subvariety, and $\Ll$ a numerically effective
line bundle on~$X$, then $\Ho^i(X,\Ll)=0$, $\forall i>0$, and the
restriction map $\Ho^0(X,\Ll)\to\Ho^0(Y,\Ll)$ is surjective. More
generally, if $X$ is spherical and $\phi:X\to Z$ is a proper
morphism, then $\Rf^i\phi_{*}\Ll=0$, $\forall i>0$, and
$\phi_{*}\Ll\to\phi_{*}(\iota_{*}\iota^{*}\Ll)$ is surjective,
where $\iota:Y\embeds X$.
\end{corollary}
See \cite{H(sph)}, \cite{B-curves} for other proofs.

More precise results on Frobenius splitting of spherical varieties
and their subvarieties (usually $G$- or $B$-orbit closures) are
obtained in special cases.

As noted above, generalized flag varieties are Frobenius split
compatibly with their Schubert subvarieties, and the latter have
rational singularities in positive, hence any (by semicontinuity),
characteristic.

Equivariant normal embeddings of $G$ (\ref{monoids}) are Frobenius
split compatibly with their $(G\times G)$-subvarieties, in all
positive characteristics. For wonderful completions of adjoint
semisimple groups, this was established by
Strickland~\cite{F-split.wonder}. The general case is due to
Rittatore~\cite{F-split.red}, see also \cite[Ch.6]{F-split}. This
implies that normal reductive group embeddings have rational
singularities (in particular, they are Cohen--Macaulay) and that
the coordinate algebras of normal reductive monoids have ``good''
filtration \cite[\S4]{F-split.red}, \cite[6.2.13]{F-split}.

Brion and Polo proved that the closures of the Schubert cells in
wonderful completions of adjoint semisimple groups (called
\emph{large Schubert varieties}) are compatibly split and deduced
that they are normal and Cohen--Macaulay \cite{large.Schubert}.

De Concini and Springer proved that wonderful embeddings of
symmetric spaces for adjoint $G$ are Frobenius split compatibly
with their $G$-subvarieties \cite[5.9]{wonder(symm)} in odd
characteristics. However this splitting is not always compatible
with $B$-orbit closures; in fact, the latter may be neither normal
nor Cohen--Macaulay \cite{orb.closures}.

See \cite{orb.closures} for a detailed study of $B$-orbits in
spherical varieties and their closures. This is an area of active
current research, with many open questions.

\appendix
\chapter*{\appendixname}
\addcontentsline{toc}{chapter}{\appendixname}
\markboth{\uppercase{\appendixname}}{} \makeatletter
\setcounter{section}{0}
\renewcommand{\thesection}{A\arabic{section}}
\renewcommand{\p@section}{}
\makeatother

\section{Rational modules and linearization}
\label{rat.mod&lin}

Rational modules are representations of algebraic groups in the
category of algebraic varieties.

\begin{definition}
Let $G$ be a linear algebraic group. A finite-dimensional
$G$-module $M$ is called \emph{rational} if the representation map
$R:G\to\GL(M)$ is a homomorphism of algebraic groups. The
terminology is explained by observing that for
$G\subseteq\GL_n(\kk)$ the matrix entries of $R(g)$ are rational
functions in the matrix entries of $g\in G$ (the denominator being
a power of~$\det g$). Generally, a \emph{rational $G$-module} is a
union of finite-dimensional rational submodules.

A $G$-algebra $A$ is said to be \emph{rational} if it is a
rational $G$-module and $G$ acts on $A$ by algebra automorphisms.

If a rational $G$-module $M$ is at the same time an $A$-module and
$g(am)=(ga)(gm)$ for $\forall g\in G, a\in A, m\in M$, then $M$ is
called a \emph{rational $G$-$A$-module}.
\end{definition}

By $\Mor(X,M)$ denote the set of all morphisms of an algebraic
variety $X$ to a vector space~$M$. (If $\dim M=\infty$, then a
morphism $X\to M$ is by definition a morphism to a
finite-dimensional subspace of~$M$.) It is a free $\kk[X]$-module:
$\Mor(X,M)\iso\kk[X]\otimes M$. If $X$ is a $G$-variety and $M$ is
a rational $G$-module, then $\Mor(X,M)$ is a rational
$G$-$\kk[X]$-module.

The $\kk[X]^G$-submodule $\Mor_G(X,M)\iso(\kk[X]\otimes M)^G$ of
equivariant morphisms is called the \emph{module of covariants} on
$X$ with values in~$M$. If $G$ is reductive and $\dim M<\infty$,
\begin{question}
and $X$ is affine
\end{question}
then $\Mor_G(X,M)$ is finite over $\kk[X]^G$ \cite[3.12]{IT}.

Another source for infinite-dimensional rational $G$-modules are
functions or global sections of sheaves on $G$-varieties.

Let $X$ be a $G$-variety, $\alpha,\pi_X:G\times X\to X$ the action
morphism and the projection, and $\Ff$ a quasicoherent sheaf
on~$X$.

\begin{definition}
A \emph{$G$-linearization} of $\Ff$ is an isomorphism
$\widehat\alpha:\pi_X^{*}\Ff\isoto\alpha^{*}\Ff$ inducing a
$G$-action on the set of local sections of $\Ff$ via isomorphisms
$\widehat\alpha|_{g\times X}:\Ff(U)\isoto\Ff(gU)$ over all $g\in
G$, $U$~open in~$X$.

A \emph{$G$-sheaf} is a quasicoherent sheaf equipped with a
$G$-linearization.
\end{definition}

\begin{theorem}[\cite{G-sheaves}]\label{rat.act}
Given a $G$-variety $X$ and a $G$-sheaf $\Ff$ on~$X$, $\kk[X]$~is
a rational $G$-algebra and $\Ho^i(X,\Ff)$ are rational
$G$-$\kk[X]$-modules.
\begin{question}
Give a proof (\cite[III.9.3]{AG}, \cite[2.5]{G-lin})
\end{question}
\end{theorem}

If $\Ff$ is the sheaf of sections of a vector bundle $F\to X$,
then its $G$-linearization is given by a fiberwise linear action
$G:F$ compatible with the projection onto~$X$.

By abuse of language we often make no terminological difference
between vector bundles and the respective locally free sheaves of
sections since they determine each other.

An important problem is to construct $G$-linearizations for line
bundles on $G$-varieties. Its treatment goes back to Mumford. Here
we follow~\cite{G-lin}.

Assume $G$ is connected.

\begin{theorem}[{\cite[2.4]{G-lin}}]\label{G-lin}
If $G$ is factorial, i.e., $\Pic G=0$, then any line bundle $\Ll$
on a normal $G$-variety $X$ is $G$-linearizable.
\end{theorem}

We say that an algebraic group $\widetilde{G}$ is a
\emph{universal cover} of $G$ if
$\widetilde{G}/\Radu{\widetilde{G}}$ is a product of a torus and a
simply connected semisimple group, and there is an epimorphism
$\widetilde{G}\to G$ with finite kernel. Every connected group has
a universal cover: it is well known for reductive groups
\cite[\S\S32,33]{Agr}, \cite[XXIII]{gr.sch}, and generally we may
put $\widetilde{G}=G\times_{G_{\red}}\widetilde{G_{\red}}$, where
$G_{\red}=G/\Radu{G}$.

\begin{corollary}\label{G^-lin}
Any line bundle $\Ll$ on $X$ is $\widetilde{G}$-linearizable.
\end{corollary}
Indeed, $\widetilde{G}$~is factorial.
\begin{corollary}\label{power-lin}
A certain power $\Ll^{\otimes d}$ of $\Ll$ is $G$-linearizable.
\end{corollary}
For $d$ one may take the degree of the universal covering or the
order of $\Pic G$ \cite[2.4]{G-lin}.

The existence of a $G$-linearization has fundamental consequences
in the local description of $G$-varieties, due to Sumihiro:
\begin{theorem}[\cite{eq.compl}, {\cite[\S1]{G-lin}}]
\label{Sumihiro} Let $G$ be a connected group acting on a normal
variety~$X$. Then any point $x\in X$ has an open $G$-stable
neighborhood $U$ which admits a locally closed $G$-equivariant
embedding ${U\embeds\PP(V)}$ for some $G$-module~$V$.
\end{theorem}
\begin{proof}
Take an affine neighborhood $U_0\ni x$. The complement
$D=X\setminus U_0$ may support no effective Cartier divisor.
However if we remove $\bigcap_{g\in G}gD$ from $X$, then any
effective Weil divisor with support $D$ becomes base point free,
hence Cartier (cf.~Lemma~\ref{Cartier}).

Take $\sigma_0\in\Ho^0(X,\Ll)$ such that $U_0=X_{\sigma_0}$. Then
%*
\begin{equation*}
\kk[U_0]=\bigcup_{d\geq0}\Ho^0(X,\Ll^{\otimes d})/\sigma_0^d
=\kk\left[\frac{\sigma_1}{\sigma_0^{d_1}},\dots,
\frac{\sigma_m}{\sigma_0^{d_m}}\right]
\end{equation*}
%*
for some $\sigma_i\in\Ho^0(X,\Ll^{\otimes d_i})$, $d_i\in\NN$.
Replacing $\Ll$ by a power, we may assume it to be a $G$-bundle
and all $d_i=1$. Include $\sigma_0,\dots,\sigma_m$ in a
finite-dimensional $G$-submodule $M\subseteq\Ho^0(X,\Ll)$. The
induced rational map $X\dasharrow\PP(V)$, $V=M^{*}$, is a locally
closed embedding on $U=GU_0$.
\end{proof}
\begin{remark}
If $X$ is itself quasiprojective, then one may take $U=X$. Indeed,
a certain power of an ample line bundle on $X$ is
$G$-linearizable, and we can find a finite-dimensional $G$-stable
space of sections inducing a projective embedding of~$X$.
\end{remark}

\section{Invariant theory}
\label{inv.th}

Let $G$ be a linear algebraic group and $A$ a rational
$G$-algebra. The subject of algebraic invariant theory is the
structure of the subalgebra $A^G$ of $G$-invariant elements.

A geometric view on the subject is to consider an affine
$G$-variety $X=\Spec A$, provided that $A$ is finitely generated.
(Note that each rational $G$-algebra is a union of finitely
generated $G$-stable subalgebras.) If $A^G$ is finitely generated,
too, then one may consider $X\by G:=\Spec A^G$ and the natural
dominant morphism $\pi=\pi_G:X\to X\by G$. The variety $X\by G$,
considered together with~$\pi$, is called the \emph{categorical
quotient} of $G:X$, because it is the universal object in the
category of $G$-invariant morphisms from $X$ to affine varieties.
This means that every morphism $\phi:X\to Y$ (with affine~$Y$)
which is constant on $G$-orbits fits into a unique commutative
triangle:
\begin{center}
%TeXCAD Picture [cat-quot.pic]. Options:
%\grade{\off}
%\emlines{\off}
%\epic{\off}
%\beziermacro{\off}
%\reduce{\on}
%\snapping{\off}
%\quality{8.00}
%\graddiff{0.01}
%\snapasp{1}
%\zoom{10.0000}
\unitlength 1ex % = 1pt
\linethickness{0.4pt}
\begin{picture}(13,11.5)(0,0)
\put(2,1){\makebox(0,0)[cc]{$X\by G$}}
\put(2,11.5){\makebox(0,0)[cc]{$X$}} \put(2,9.5){\vector(0,-1){6}}
\put(3,10){\vector(3,-1){9}} \put(3,3){\vector(3,1){9}}
\put(13,6.5){\makebox(0,0)[lc]{$Y$}}
\put(8,9){\makebox(0,0)[cb]{$\phi$}}
\put(8,4){\makebox(0,0)[ct]{$\overline\phi$}}
\put(1,6.5){\makebox(0,0)[rc]{$\pi$}}
\end{picture}
\end{center}

Geometric properties of $X\by G$ and $\pi_G$ translate into
algebraic properties of $A^G$ and of its embedding into~$A$, and
vice versa. The case of reductive $G$ is considered by Geometric
Invariant Theory (GIT). We collect basic results on invariants and
quotients of affine varieties by reductive groups in the following
theorem.
\begin{Claim}[Main Theorem of GIT]
Let $G$ be a reductive group and $A$ a rational $G$-algebra.
\begin{enumerate}
\item If $A$ is finitely generated, then $A^G$ is so.
\end{enumerate}
Under this assumption, put $X=\Spec A$. Then $\pi:X\to X\by G$ is
well defined and has the following properties:
\begin{enumerate}
\setcounter{enumi}{1} \item $\pi$~is surjective and maps closed
$G$-stable subsets of $X$ to closed subsets of~$X\by G$. \item
$X\by G$ carries the quotient topology w.r.t.~$\pi$, and
$\Oo_{X\by G}=\pi_{*}\Oo_X^G$. \item If $Z_1,Z_2\subset X$ are
disjoint closed $G$-stable subsets, then
$\pi(Z_1)\cap\pi(Z_2)=\emptyset$. In particular, each fiber of
$\pi$ contains a unique closed orbit.
\end{enumerate}
\end{Claim}
Thus $X\by G$ may be regarded as the ``variety of closed orbits''
for $G:X$. It is not hard to show that $\pi_G:X\to X\by G$ is the
categorical quotient in the category of \emph{all} algebraic
varieties.

Finite generation of $G$-invariants goes back to Hilbert and Weyl
(in characteristic zero), the general case is due to Nagata and
Haboush. Other assertions are due to Mumford.

If $G$ is \emph{linearly reductive}, i.e., all rational
$G$-modules are completely reducible (e.g., $\ch\kk=0$ or $G$ is a
torus), then the proof is considerably simplified \cite[3.4,
4.4]{IT} by using the $G$-$A^G$-module decomposition $A=A^G\oplus
A_G$, where $A_G$ is the sum of all nontrivial irreducible
$G$-submodules. The respective projection
%*
\begin{equation*}
A\mapsto A^G,\qquad f\mapsto f^{\natural}
\end{equation*}
%*
is known as the \emph{Reynolds operator}. For finite $G$ and
$\ch\kk=0$, it is just the group averaging:
%*
\begin{equation*}
\average{f}=\frac1{|G|}\sum_{g\in G}gf,\qquad\forall f\in A
\end{equation*}
%*
For a complex reductive group~$G$, the Reynolds operator may be
defined by averaging over a compact real form of~$G$.

The proof in positive characteristic may be found in
\cite[App.~1A,~1C]{GIT}.

For non-reductive groups the situation is not so nice---even
finite generation of invariants fails due to famous Nagata's
counterexample \cite{Nagata.ex} and results of
Popov~\cite{Hilb-14}. However for subgroups of reductive groups
acting on algebras or affine varieties, there are positive results
on finite generation and the structure of invariant algebras and
categorical quotients.

\begin{lemma}\label{(A/I)^H}
Let $G$ be a reductive group, $H\subseteq G$ be an algebraic
subgroup, $A$~be a rational $G$-algebra, and $I\normin A$ be a
$G$-stable ideal. Then $(A/I)^H$ is a purely inseparable finite
extension of~$A^H/I^H$.
\end{lemma}
The lemma is obvious for $\ch\kk=0$, since $A/I$ lifts to a
$G$-submodule of~$A$, whence $(A/I)^H=A^H/I^H$. The proof for
$H=G$ may be found in \cite[Lemma~A.1.2]{GIT} and the general case
follows by the transfer principle (Remark~\ref{transfer}).
\begin{corollary}\label{(M/N)^H}
Let $M$ be a $G$-module and $N\subset M$ a $G$-submodule. For any
$\overline{m}\in(M/N)^{(H)}$ there exist $q=p^n$ and
$m\in(\Sym^qM)^{(H)}$ such that
$m\mapsto\overline{m}^q\in\Sym^q(M/N)$.
\end{corollary}
\begin{proof}
Just apply Lemma~\ref{(A/I)^H} to $A=\Sym^{\bullet}M$, $I=AN$,
replacing $H$ by the common kernel $H_0$ of all $\chi\in\Ch(H)$,
and use the fact that $H/H_0$ is diagonalizable.
\end{proof}
\begin{corollary}\label{(grA)^H}
If $A$ carries a $G$-stable filtration, then $(\gr A)^H$ is a
purely inseparable extension of~$\gr(A^H)$.
\end{corollary}
\begin{proof}
Each homogeneous component of $\gr A$ has the form $M/N$, where
$M,N$ are two successive members of the filtration. It remains to
apply Corollary~\ref{(M/N)^H} to $M$ and~$N$.
\end{proof}

The most important case is $H=U$, a maximal unipotent subgroup
in~$G$. $U$-invariants of rational $G$-algebras were studied by
Hadzhiev, Vust, Popov (in characteristic zero), Donkin, Grosshans
(in arbitrary characteristic), et al. We refer to \cite{HS} for
systematic exposition of the theory.
\begin{lemma}[{\cite[14.3]{HS}}]\label{A>GA^U}
A rational $G$-algebra is integral over its subalgebra $\langle
G\cdot A^U\rangle$.
\end{lemma}
The proof relies on Lemma~\ref{(A/I)^H}. In characteristic zero,
Lemma~\ref{A>GA^U} is trivial, since $A=\langle G\cdot A^U\rangle$
by complete reducibility of $G$-modules and highest weight theory.

The fundamental importance of $U$-invariants is explained by the
fact that $A$ and $A^U$ share a number of properties.
\begin{theorem}\label{U-inv}
Let $G$ be a connected reductive group and $A$ a rational
$G$-algebra.
\begin{enumerate}
\item\label{A<=>A^U} A rational $G$-algebra $A$~is finitely
generated (resp.\ has no nilpotents, is an integral domain) iff
$A^U$ is so. \item In particular, for any affine $G$-variety~$X$,
the categorical quotient $\pi_U:X\to X\by U$ is well defined.
\item\label{X<=>X/U} In characteristic zero, $X$~is normal (has
rational singularities) iff $X\by U$ is so.
\end{enumerate}
\end{theorem}
In characteristic zero, finite generation of $U$-invariants is due
to Hadzhiev; other assertions were partially proved by Brion
(nilpotents and zero divisors), Vust (normality), Kraft
(rationality of singularities), and by Popov in full
generality~\cite{contr}. The theorem was extended to arbitrary
characteristic by Grosshans. We shall give an outline of the proof
scattered in~\cite{HS}.

Finite generation of $A^U\iso(\kk[G/U]\otimes A)^G$ (see
Remark~\ref{transfer}) stems from that of $A$ and of~$\kk[G/U]$.
The latter is proved by a representation-theoretic argument
(Lemma~\ref{V(l)*V(m)}) or by providing an explicit embedding of
$G/U$ into a $G$-module, with the boundary of codimension~$2$
\cite[5.6]{HS} (cf.~Theorem~\ref{norm.aff.S-var}).

The other assertions are proved using horospherical contraction
(cf.~\ref{horosph}).

The algebra $A$ is endowed with a $G$-stable increasing filtration
$A^{(n)}$ such that $\gr A$ has an integral extension
$S=(\kk[G/U^{-}]\otimes A^U)^T$ and $A^U\iso\gr(A^U)=(\gr A)^U\iso
S^U$, where $T$ acts on $G/U^{-}$ by right translations. In
characteristic zero this filtration is described in~\ref{horosph}
(for $A=\kk[X]$) and the general case is considered in
\cite[\S15]{HS}.

Now finite generation of $A^U$ implies that of~$S$, hence of~$\gr
A$ (both are finite modules over $\langle G\cdot S^U\rangle$), and
finally of $A$ by a standard argument. Moreover, the algebra
$R=\bigoplus_{n=0}^{\infty}A^{(n)}t^n\subseteq A[t]$ is finitely
generated, too (because $R^U\iso A^U[t]$) \cite[16.5]{HS}.

The remaining assertions may be proved for finitely generated $A$
and $X=\Spec A$. As in~\ref{horosph}, put $E=\Spec R$ and consider
the natural $G\times\kk^{\times}$-morphism $\delta:E\to\AAA^1$
with the zero fiber $X_0=\Spec\gr A$ and other fibers isomorphic
to~$X$. Note that $\kk^{\times}$ contracts $E$ to~$X_0$ (i.e.,
$\forall x\in E\ \exists\lim_{t\to0}t\cdot x\in X_0$), because the
grading on $R$ is non-negative.

Since $\delta$ is flat, the set of $x\in E$ such that the
schematic fiber $\delta^{-1}(\delta(x))$ has a given local
property of open type (e.g., is reduced, irreducible, normal, or
has rational singularity) at $x$ is open in~$E$. The complementary
closed subset of $E$ is $\kk^{*}$-stable, whence it intersects
$X_0$ whenever it is non-empty. It follows that $X$ has the
property of open type whenever $X_0$ has this property.

If an affine $\kk$-scheme $Z$ of finite type is reduced (resp.\
irreducible or normal), then $Z\by H$ is so for any algebraic
group $H$ acting on~$Z$. (Only normality requires some
explanation.) In particular, these properties are inherited by
$X\by U$ from~$X$.

Conversely, if $X\by U$ has one of these properties, then $\Spec
S=(G\by U^{-}\times X\by U)\by T$ and $X_0$ have it, too. (For
normality, we use the isomorphism $S\iso\gr A$ in characteristic
zero.) By the above, $X$~has this property. More elementary (and
lengthy) arguments are given in \cite[\S18]{HS}.

The same reasoning works for rational singularities in
characteristic zero, using the facts that $G\by U$ has rational
singularities \cite{collaps} and that rational singularities are
preserved by products \cite{rat.sing} and categorical quotients
modulo reductive groups \cite{rat.quot} (this guarantees that
$X\by U=(G\by U\times X)\by G$ has rational singularities provided
that $X$ is so).

\section{Geometric valuations}
\label{geom.val}

Let $K$ be a function field, i.e., a finitely generated field
extension of~$\kk$. By a valuation $v$ of $K$ we always mean a
discrete $\QQ$-valued valuation of $K/\kk$, i.e., assume the
following properties:
\begin{enumerate}
\item $v:K^{\times}\to\QQ$, $v(0)=\infty$ \item
$v(K^{\times})\iso\ZZ$ \item $v(\kk^{\times})=0$ \item
$v(fg)=v(f)+v(g)$ \item $v(f+g)\geq\min(v(f),v(g))$
\end{enumerate}
\begin{remark}
If $v$ is defined only on a $\kk$-algebra $A$ with $\Quot{A}=K$,
then it is extended to $K$ in a unique way by putting
$v(f/g)=v(f)-v(g)$, $f,g\in A$.
\end{remark}

Our main source in the valuation theory is \cite[Ch.4, App.2]{ZS}.

\begin{definition}
A valuation $v$ is called \emph{geometric} if there exists a
normal variety $X$ with $\kk(X)=K$ (a \emph{model} of~$K$) and a
prime divisor $D\subset X$ such that $v(f)=c\cdot v_D(f)$,
$\forall f\in K$, for some $c\in\QQ_{+}$. Here $v_D(f)$ is the
order of $f$ along~$D$.
\end{definition}

To any valuation corresponds a (\emph{discrete}) \emph{valuation
ring} (\emph{DVR}) $\Oo_v=\{f\in K\mid v(f)\geq0\}$, which is a
local ring with the maximal ideal $\m_v=\{f\in K\mid v(f)>0\}$ and
quotient field~$K$. The \emph{residue field} of $v$ is
$\kk(v)=\Oo_v/\m_v$.
\begin{example}
If $v$ is geometric, then $\Oo_v=\Oo_{X,D}$, $\kk(v)=\kk(D)$.
\end{example}
\begin{Claim}[Properties]
\begin{enumerate}
\item $\Oo_v$ is a maximal subring of~$K$. \item $\Oo_v$
determines $v$ up to proportionality.
\end{enumerate}
\end{Claim}

\begin{definition}
Let $X$ be a model of~$K$. A closed irreducible subvariety
$Y\subseteq X$ is the \emph{center} of $v$ on~$X$ if $\Oo_v$
dominates $\Oo_{X,Y}$ (i.e., $\Oo_v\supseteq\Oo_{X,Y}$,
$\m_v\supseteq\m_{X,Y}$, which implies $\kk(v)\supseteq\kk(Y)$).
\end{definition}
\begin{example}
A prime divisor $D\subset X$ is the center of the respective
geometric valuation.
\end{example}
If $\phi:X\to X'$ is a dominant morphism and $v$ has the center
$Y\subseteq X$, then the restriction $v'$ of $v$ to $K'=\kk(X')$
has the center $Y'=\overline{\phi(Y)}\subseteq X'$.

\begin{description}
\item[Valuative criterion of separation:] $X$ is separated iff any
(geometric) valuation has at most one center on~$X$
\item[Valuative criterion of properness:] The map $\phi:X\to X'$
is proper iff any (geometric) valuation of $K$ has the center
on~$X$ provided that its restriction to $K'$ has the center
on~$X'$. \item[Valuative criterion of completeness:] $X$ is
complete iff any (geometric) valuation has the center on~$X$.
\end{description}

\begin{proposition}
If $X$ is affine, then $v$ has the center $Y\subseteq X$ iff
$v|_{\kk[X]}\geq0$, and then $\Ideal{Y}=\kk[X]\cap\m_v$.
\end{proposition}

\begin{proposition}\label{tr.deg.k(v)}
A valuation $v\ne0$ is geometric iff $\trdeg\kk(v)=\trdeg K-1$.
\end{proposition}
\begin{proof}
Assume that $\trdeg K=n$ and that the residues of
$f_1,\dots,f_{n-1}\in\Oo_v$ form a transcendence base of
$\kk(v)/\kk$. Take a nonzero $f_n\in\m_v$; then $f_1,\dots,f_n$
are easily seen to be a transcendence base of $K/\kk$. Consider an
affine variety $X$ such that $\kk[X]$ is the integral closure of
$\kk[f_1,\dots,f_n]$ in~$K$. It is easy to show that
$v|_{\kk[X]}\geq0$, whence $v$ has the center $D\subset X$ and
$f_1,\dots,f_{n-1}\in\kk[D]$ are algebraically independent. Hence
$D$ is a prime divisor, and $\Oo_v=\Oo_{X,D}$ implies $v=v_D$ up
to a multiple. The converse implication is obvious.
\end{proof}

\begin{proposition}\label{res&ext(val)}
Let $\kk\subseteq K'\subseteq K$ be a subfield.
\begin{enumerate}
\item\label{res(val)} If $v$ is a geometric valuation of~$K$, then
$v'=v|_{K'}$ is geometric. \item\label{ext(val)} Any geometric
valuation $v'$ of $K'$ extends to a geometric valuation $v$ of
$K$.
\end{enumerate}
\end{proposition}
\begin{proof}
\begin{roster}
\item[\ref{res(val)}] Take $f_1,\dots,f_k\in\Oo_v$ whose residues
form a transcendence base of $\kk(v)/\kk(v')$. They are
algebraically independent over~$K'$ (otherwise one can take an
algebraic dependence of $f_1,\dots,f_k$ over $\Oo_{v'}$ with at
least one coefficient not in~$\m_{v'}$, and pass to residues
obtaining a contradiction). Hence
$\trdeg\kk(v')=\trdeg\kk(v)-\trdeg_{\kk(v')}\kk(v)\geq \trdeg
K-1-\trdeg_{K'}K=\trdeg K'-1$, and we conclude by
Proposition~\ref{tr.deg.k(v)}. \item[\ref{ext(val)}] Take a
complete normal variety $X'$ with a prime divisor $D'\subset X'$
such that $v'$ is proportional to~$v_{D'}$. We may construct a
complete normal variety $X$ with $\kk(X)=K$ mapping onto~$X'$:
take any complete model $X$ of $K$ and replace it by the
normalization of the closure of the graph of the rational map
$X\dasharrow X'$. Let $D\subset X$ be a component of the preimage
of~$D'$ mapping onto~$D'$. Then we may take $v=v_D$ up to a
multiple. \qedhere\end{roster}
\end{proof}

\section{Schematic points}
\label{sch.pt}

Given a $\kk$-scheme~$X$, it is often instructive to consider the
respective representable functor associating with any $\kk$-scheme
$S$ the set $X(S)$ of $\kk$-mor\-phisms $S\to X$, called
\emph{$S$-points} of~$X$. If $S=\Spec{A}$ is affine, then
$S$-points are called \emph{$A$-points} and the notation
$X(A):=X(S)$ is used.
\begin{example}\label{aff(alg)}
If $X\subseteq\AAA^n$ is an (embedded) affine scheme of finite
type, then its $A$-point is given by an algebra homomorphism
$\kk[X]=\kk[t_1,\dots,t_n]/\Ideal{X}\to A$, i.e., by an $n$-tuple
$x=(x_1,\dots,x_n)\in A^n$ satisfying the defining equations
of~$X$. A similar description works for quasiaffine schemes.
\end{example}

We require a closer look at this notion in case, where $X$ is an
algebraic variety over $\kk$ and $A$ is a local $\kk$-algebra with
the maximal ideal~$\m$. Given $\chi\in X(A)$, the closed point of
$\Spec A$ is mapped by $\chi$ to the generic point of an
irreducible subvariety $Y\subset X$ called the \emph{center}
of~$\chi$. If $\X\subseteq X$ is an affine chart meeting~$Y$, then
$\chi\in\X(A)$. Thus $X(A)=\bigcup\X(A)$ over all affine open
subsets $\X\subseteq X$. From the algebraic viewpoint, an
$A$-point of $X$ is given by an irreducible subvariety $Y\subseteq
X$ and a local algebra homomorphism $\Oo_{X,Y}\to A$, $\m_Y\to\m$,
or by a homomorphism $\kk[\X]\to A$, where $\X\subseteq X$ is an
affine chart (intersecting~$Y$).
\begin{example}\label{gen.pt}
The \emph{generic point} of an irreducible variety $X$ over
$\kk(X)$ has the center~$X$, and $\Oo_{X,X}\to\kk(X)$ is the
identity map. Informally, the coordinates of the generic point are
indeterminates bound only by relations that hold identically
on~$X$.
\end{example}
\begin{example}
If $v$ is a valuation of $\kk(X)$ with center $Y\subseteq X$, then
the inclusion $\Oo_{X,Y}\subseteq\Oo_v$ yields an $\Oo_v$-point of
$X$ with center~$Y$.
\end{example}
\begin{example}
Any $A$-point of a quasiprojective scheme $X\subseteq\PP^n$ is at
the same time an $A$-point of $X\cap\AAA^n$ for a certain affine
chart $\AAA^n\subseteq\PP^n$. In view of Example~\ref{aff(alg)},
$A$-points of $X$ are identified with tuples $x=(x_0:\dots:x_n)$,
$x_i\in A$, considered up to proportionality, satisfying the
defining equations of~$X$, and such that at least one $x_i$ is
invertible.
\end{example}

It is quite common in algebraic geometry to consider the case,
where $A$ is a field. For applications in~\ref{form.curv}, we
consider points over the function field of an algebraic curve or
its formal analogue.
\begin{definition}
A \emph{germ of a curve} in $X$ is a pair $(\chi,\theta_0)$, where
$\chi\in X(\kk(\Theta))$, $\Theta$~is a smooth projective curve,
and $\theta_0\in\Theta$. In other words, a germ of a curve is
given by a rational map from a curve to $X$ and a fixed base point
on the curve.

The germ is said to be \emph{convergent} if $\chi\in
X(\Oo_{\Theta,\theta_0})$, i.e., the rational map
$\chi:\Theta\dasharrow X$ is regular at~$\theta_0$. The point
$x_0=\chi(\theta_0)$ is the \emph{limit} of the germ.
\end{definition}

There is a formal analytic analogue of this notion.
\begin{definition}
A \emph{germ of a formal curve} in $X$ is a $\kk\(t\)$-point
of~$X$. A $\kk\[t\]$-point is called a \emph{convergent formal
germ}, and its center $x_0\in X$ is the \emph{limit} of the formal
germ.
\end{definition}

It is natural to think of a formal germ as of a ``parametrized
formal analytic curve'' $x(t)$ in~$X$. In local coordinates,
$x(t)$~is a tuple of Laurent series satisfying the defining
equations. If $x(t)$ converges, then its coordinates are power
series, and their constant terms are the coordinates of the limit
$x_0=:x(0)=\lim_{t\to0}x(t)$.

With any germ of a curve $(\theta_0\in\Theta\dasharrow X)$ one can
associate a formal germ via the inclusions
$\Oo_{\Theta,\theta_0}\subset\widehat\Oo_{\Theta,\theta_0}\iso
\kk\[t\]$, $\kk(\Theta)\subset\kk\(t\)$, depending on the choice
of a formal uniformizing parameter
$t\in\widehat\Oo_{\Theta,\theta_0}$.
\begin{proposition}
A formal germ is induced by a germ of a curve iff its center has
dimension~$\leq1$.
\end{proposition}
\begin{proof}
The ``only if'' direction and the case, where the center is a
point, are clear. Suppose the center of a formal germ is a curve
$C\subseteq X$. Then $\kk(C)\embeds\kk\(t\)$. Choose any
$f\in\kk(C)$, $\ord_tf=k>0$, and consider $s\in\kk\[t\]$, $s^k=f$.
Then $\kk(C)(s)=\kk(\Theta)$ is a function field of a smooth
projective curve~$\Theta$, and
$\kk(\Theta)\cap\kk\[t\]=\Oo_{\Theta,\theta_0}$ for a certain
$\theta_0\in\Theta$, so that
$\widehat\Oo_{\Theta,\theta_0}=\kk\[s\]=\kk\[t\]$.
\end{proof}

There is a \emph{$t$-adic topology} on $X(\kk\(t\))$ thinner than
the Zariski topology~\cite{adic}. For $X=\AAA^n$, a basic $t$-adic
neighborhood of $x(t)=(x_1(t),\dots,x_n(t))$ consists of all
$y(t)=(y_1(t),\dots,y_n(t))$ such that $\ord_t(y_i(t)-x_i(t))\geq
N$, $\forall i=1,\dots,n$, where $N\in\NN$. The $t$-adic topology
on arbitrary varieties is induced from that on affine spaces using
affine charts.

An important approximation result is due to Artin:
\begin{theorem}[\cite{approx}]\label{approx(form)}
The set of formal germs induced by germs of curves is dense in
$X(\kk\(t\))$ w.r.t.\ the $t$-adic topology.
\end{theorem}

\cleardoublepage
\addcontentsline{toc}{chapter}{\bibname}
\footnotesize
\bibliography{mrabbrev,embed}
\bibliographystyle{house}
\end{document}